\documentclass[11pt]{amsart}

\usepackage{amssymb} 
\usepackage[latin1]{inputenc}
\usepackage[matrix,arrow,curve]{xy}
\usepackage[mathscr]{eucal}
\usepackage{epic,eepic}
\usepackage{bbm}
\usepackage{manfnt}

\makeindex
\newtheorem{theorem}{Theorem}[subsection]
\newtheorem{lemma}[theorem]{Lemma}
\newtheorem{prop}[theorem]{Proposition}

\newtheorem{cor}[theorem]{Corollary}

\newtheorem{conj}[theorem]{Conjecture}
\theoremstyle{definition}
\newtheorem{definition}[theorem]{Definition}
\newtheorem{remark}[theorem]{Remark}
\newtheorem{example}[theorem]{Example}
\newtheorem{exemple}[theorem]{Example}

\newtheorem{notation}[theorem]{Notation}

\newenvironment{pf}
{\medskip\noindent {\it Proof --- \ }}
{\hfill\nobreak $\Box$ \par\bigbreak}

\newcommand{\RR}{{\mathcal R}}

\newcommand{\Hom}{\text{Hom}}

\newcommand{\p}{{\mathfrak P}}
\newcommand{\Zp}{{\mathbb Z}_p }
\newcommand{\Qp}{{\mathbb Q}_p }
\newcommand{\F}{ \mathbb F} 
\newcommand{\C}{{ \mathbb C  }}
\newcommand{\R}{{ \mathbb R  }}
\newcommand{\Q}{{ \mathbb Q } }
\newcommand{\Z}{{ \mathbb Z  }}
\newcommand{\N}{{ \mathbb N  }}
\newcommand{\got}{\mathfrak}

\newcommand{\HH}{{\mathcal H}}
\newcommand{\PP}{{\mathcal P}}

\newcommand{\m}{m}
\newcommand{\Ker}{{\text{Ker}\,}}

\newcommand{\im}{{\text{Im}\,}}

\renewcommand{\Im}{{\text{Im}\,}}

\renewcommand{\ker}{{\text{Ker}\,}}
\newcommand{\Gal}{{\mathrm{Gal}\,}}

\newcommand{\rhobar}{{\bar \rho}}

\newcommand{\A}{\mathbb A}

\newcommand{\anneau}{{ \mathcal O}}
\newcommand{\MM}{{\mathcal M}}
\newcommand{\M}{{\mathcal M}}

\newcommand{\Aut}{{\text{Aut}\,}}
\newcommand{\End}{{\text{End}}}
\newcommand{\Id}{{\text{Id}}}

\newcommand{\U}{{\text{U}}}
\newcommand{\SU}{{\text{SU}}}

\newcommand{\Gl}{{\text {GL}}}

\newcommand{\GL}{{\text {GL}}}

\newcommand{\SL}{{\text {SL}}}
\newcommand{\Sl}{{\text {SL}}}

\newcommand{\opp}{\rm opp}

\newcommand{\spec}{{\text{Spec\,}}}

\newcommand{\vareps}{{\varepsilon}}

\newcommand{\pr}{{\text{pr}}}

\newcommand{\rec}{{\text{rec}}}

\newcommand{\Sc}{{\mathcal{S}}}

\newcommand{\G}{{\mathfrak g}}

\newcommand{\Ind}{{\mathrm{Ind}}}
\newcommand{\ad}{{\text{ad}}}

\newcommand{\Ext}{{\text{Ext}}}
\newcommand{\Tor}{{\text{Tor}}}

\newcommand{\Res}{{\text{Res}}}
\newcommand{\ses}{{\text{ss}}}

\newcommand{\diag}{{\text{diag}}}
\newcommand{\Frob}{{\text{Frob\,}}}

\newcommand{\ord}{{\text{ord}}}
\newcommand{\tr}{{\text{tr\,}}}

\newcommand{\rhob}{{\bar \rho}}
\newcommand{\sd}{\ltimes}

\newcommand{\Rb}{\overline{R}}
\newcommand{\Tb}{\overline{T}}
\newcommand{\rad}{{\mathrm{rad}}}
\newcommand{\ps}{\par \smallskip}

\newcommand{\surjection}{\rightarrow}
\newcommand{\rig}{{\text{rig}}}

\newcommand{\cris}{{\text{crys}}}

\renewcommand{\hat}{\widehat}

\newcommand{\isomo}{\overset{\sim}{\longrightarrow}}
\newcommand{\cal}{\mathcal}
\newcommand{\Sym}{{\rm Symm}}
\newcommand{\Bun}{{B^{\text{univ}}}}
\newcommand{\rot}{{\rm rot}}
\newcommand{\saut}{${}^{}$ \ps}

\newcommand{\Bc}{B_{\rm crys}}
\newcommand{\Bdr}{B_{\rm dR}}
\newcommand{\Ro}{\cal{R}}
\newcommand{\W}{W_{\Q_p}}
\newcommand{\Fil}{{\rm Fil}}
\newcommand{\Dc}{D_{\rm crys}}
\newcommand{\fg}{(\varphi,\Gamma)}
\newcommand{\Dr}{D_{\rm rig}}
\newcommand{\Ddr}{D_{\rm dR}}
\newcommand{\Ddrf}{\cal{D}_{\rm dR}}
\newcommand{\Dcf}{\cal{D}_{\rm crys}}
\newcommand{\Ds}{D_{\rm Sen}}
\newcommand{\Dsf}{\cal{D}_{\rm Sen}}
\newcommand{\Ref}{\cal{F}}
\newcommand{\bdelta}{\overline{\delta}}
\newcommand{\fr}{{\rm Frac}}
\newcommand{\Gp}{{G_p}}
\newcommand{\Hp}{{H_p}}
\newcommand{\rk}{{\rm rk}}
\newcommand{\Qpb}{{\overline{\Q}_p}}
\renewcommand{\Gal}{{\rm Gal}}

\newcommand{\Tc}{{\cal{T}}}
\newcommand{\WD}{\cal{X}}
\newcommand{\dwk}{{{d}\choose{k}}}
\newcommand{\selr}{{\rm Sel}(\rho)}

\newcommand{\Dcp}{D_{\rm crys}^+}
\renewcommand{\Ds}{D_{\rm Sen}}
\renewcommand{\Ddr}{D_{\rm DR}}
\renewcommand{\Dc}{D_{\rm crys}}
\newcommand{\Bcp}{B_{\rm crys}^+}
\newcommand{\OO}{\mathcal{O}}
\newcommand{\Wk}{{\Lambda^k}}
\newcommand{\MB}{\overline{\MM}}

\newcommand{\Qb}{\overline{\Q}}
\newcommand{\AAA}{\mathbb{A}}
\newcommand{\wt}{\widehat{\otimes}}

\newcommand{\disc}{{\rm disc}}

\newcommand{\LG}{{}^L\!G}

\renewcommand{\G}{G}
\renewcommand{\rec}{{\rm rec}} 
\renewcommand{\W}{{\rm W}}
\newcommand{\Irr}{{\rm Irr}}
\newcommand{\Cc}{\mathcal{C}_c}
\newcommand{\ATL}{\mathcal{A}_p}
\renewcommand{\ss}{\mathrm ss}
\newcommand{\unit}{\text{unit}}
\newcommand{\cusp}{\text{cusp}}
\newcommand{\temp}{\text{temp}}
\newcommand{\bs}{{\bf{S}}}

\newcommand{\Ad}{\text{Ad}}
\newcommand{\kum}{\text{kum}}

\newcommand{\GAFP}{G(\AAA_f^p)}
\newcommand{\FT}{\widetilde{F}}
\newcommand{\Mod}{{\rm Mod}}
\newcommand{\Fc}{{\cal F}}
\newcommand{\K}{\cal K}
\newcommand{\Kb}{\overline{\cal K}}
\newcommand{\kbx}{\overline{k(x)}}
\newcommand{\rhog}{{\rho}^{\rm gen}}

\newcommand{\nb}{{\bar n}}
\newcommand{\Nb}{\overline{N}}
\newcommand{\Ng}{N^{\rm gen}}
\newcommand{\rb}{\overline{r}}
\newcommand{\rg}{r^{\rm gen}}
\newcommand{\XX}{X}
\newcommand{\dkappa}{\partial\kappa}
\newcommand{\WW}{\mathcal W}
\newcommand{\WG}{\WW \times \mathbb{G}_m}

\begin{document}

\baselineskip 14pt

\bibliographystyle{style} 

\title[$p$-adic families of Galois representations]{$p$-adic families of Galois representations and higher rank Selmer groups}

\author[J.~Bella\"iche]{Jo\"el Bella\"iche}
\email{jbellaic@math.columbia.edu}
\address{509 Math Building, Columbia University \\
         2990 Broadway \\
         New York, NY 10027}

\author[G.~Chenevier]{Ga\"etan Chenevier}
\email{Gaetan.Chenevier@normalesup.org}
\address{Laboratoire Analyse G\'eom\'etrie et Applications \\ 
Institut Galil\'ee \\ Université Paris 13 \\ 
99 Av. J.-B. Cl\'ement \\
93430 Villetaneuse \\ France}
\maketitle

% Commande concernant les tables des matieres de chaque chapitre...
%\dosecttoc
%\setcounter{secttocdepth}{3}
%\adjuststc
%\renewcommand{\stctitle}{}
%\nostcrule
\setcounter{tocdepth}{2}
\tableofcontents
\newpage 

${}^{}$

\par \bigskip

${}^{}$

\par \bigskip

${}^{}$

\par \bigskip

${}^{}$

\par \bigskip

{{ \large \bf \center{Ce livre est d\'edi\'e \`a la m\'emoire de Serge Bella\"\i che, fr\`ere et ami.}}}

\newpage

\section*{Introduction}\label{introgenerale}

This book\footnote{During the elaboration and writing of this book, 
Jo\"el Bella\"\i che was supported by the NSF grant DMS-05-01023. 
Ga\"etan Chenevier would like to thank the I.H.E.S. for their hospitality during the last three months of this work.}
takes place in a now thirty years long trend of researches, 
initiated by Ribet (\cite{Ribet}) 
aiming at constructing "arithmetically interesting" non trivial extensions 
between global Galois representations 
(either on finite $p^\infty$-torsion modules, 
or on $p$-adic vector spaces) or, as we shall say, 
non-zero elements of Selmer groups, by studying congruences or variations of automorphic forms. 
 As far as we know, despite of its great successes 
(to name one : the proof of Iwasawa's main 
conjecture for totally real fields by Wiles \cite{wiles}), 
this current of research has never established, in any case, 
the existence of {\it two linearly independent} 
elements in a Selmer group -- 
although well-established conjectures predict that sometimes such elements should exist.\footnote{By a very different approach, let us mention here that the parity theorem of Nekovar \cite{Nekovar} shows, in the {\it sign +1} case and for $p$-ordinary modular forms, that the rank of the Selmer group is bigger than $2$ if nonzero.}
The final aim of the book is, focusing on the characteristic zero case, 
to understand the conditions under which, by this kind of method, 
existence of two or more independent elements in a Selmer space could be
proved.

To be somewhat more precise, 
let $G$ be reductive group over a number field. 
We assume, to fix ideas, 
that the existence of the $p$-adic rigid analytic 
eigenvariety $\cal E$ of $G$, as well as the existence and basic properties
of the Galois representations attached 
to algebraic automorphic forms of $G$ are 
known\footnote{Besides $\Gl_2$ over a totally real field and its forms, 
the main examples in the short term of such $G$ would be 
suitable unitary groups and suitable forms of $\rm{GSp}_4$. Concerning
unitary groups in $m\geq 4$ variables, it is one of the goals of 
the book project of the GRFA seminar~\cite{GRFAbook} to construct the
expected Galois representations, which makes the assumption relevant. All our applications to unitary
eigenvarieties for such groups (hence to Selmer groups) will be conditional to
their work. However, thanks to Rogawski's work and \cite{ZFPMS}, everything concerning $U(3)$ will be unconditional.}. Thus $\cal E$ carries a
family of $p$-adic Galois representations. Our main result has the form of a numerical relation between the dimension of the tangent space at suitable points 
$x \in \cal E$ and the dimension of the part of Selmer groups of 
components of $\ad \rho_x$ that are "seen by $\cal E$", where 
$\rho_x$ is the Galois representation carried by $\cal E$ at the point $x$.  
 
Such a result can be used both way : if the Selmer groups are known, and small, it can be used to study the geometry of 
$\cal E$ at $x$, for example (see~\cite{kis}, \cite{lissite}) to prove its smoothness. On the other direction, it can be used to get a lower bound
on the dimension of some interesting Selmer groups, 
lower bound that depends on the dimension of the tangent space of $\cal E$ at $x$.
An especially interesting case is the case of unitary groups with $n+2$ variables, 
and of some particular points $x  \in \cal E$ attached to non-tempered automorphic
forms\footnote{In general, their existence
is predicted by Arthur's conjectures and known in some cases.}. These forms were
already used in \cite{Bellaichethese} for a unitary group with three variables, and later for
$\rm{GSp}_4$ in \cite{SU}, and for $\U(3)$ again in \cite{BC}. At those points, the Galois representation $\rho_x$ is a sum of an irreducible $n$-dimensional 
representation\footnote{Say, of the absolute Galois group of a quadratic imaginary
field.}
$\rho$, the trivial character, and the cyclotomic character $\chi$. 
The representations $\rho$ we could get this way are, at least conjecturally,
all irreducible $n$-dimensional representations satisfying some auto-duality condition, 
and such that the order of vanishing of $L(\rho,s)$ at the center of its functional equation
is odd. Our result then gives a lower bound on the dimension of the Selmer group of
$\rho$. Let us call $\selr$ this Selmer group\footnote{Precisely, this is the
group usually denoted by $H^1_f(E,\rho)$.}. This lower bound implies,
in any case, that $\selr$ is non zero 
(which is predicted by the Bloch-Kato conjecture), and if $\cal E$ is non smooth at $x$,
that the dimension of $\selr$ is at least $2$. 

This first result (the non-triviality of $\selr$, proved in chapter 8) extends to any dimension $n$ a previous work of the
authors \cite{BC} in which we proved that $\selr\neq 0$ in
the case $n=1$, i.e. $G=\U(3)$, and the work of Skinner-Urban \cite{SU} in
the case $n=2$ and $\rho$ ordinary. Moreover, 
the techniques developed in this paper shed 
also many lights back on those works. For example, the arguments in
\cite{BC} to produce a non trivial element in $\selr$ involved some 
arbitrary choice of a germ of irreducible curve 
of $\cal E$ at $x$, and it was not clear in which way the
resulting element depended on them. With our new method, we do not have to
make these choices and we construct directly a canonical subvectorspace of
$\selr$. 

In order to prove our second, main, result (the lower bound on $\dim (\selr)$, see chapter 9) 
we study the reducibility loci of the family of Galois representations on
$\cal E$. An originality of the present work is that we focus on points $x \in \cal{E}$ at which the
Galois deformation at $p$ is as non trivial as possible (we call some of
them {\it anti-ordinary points})\footnote{A bit more precisely, among the
(finite number of) points of $x \in \cal E$ having the same Galois
representation $\rho_x$, we choose one which is {\it refined} in quite a special way.}. 
We discovered that at these points, the local Galois deformation is highly
irreducible, that is not only generically irreducible\footnote{Note that although we do
assume in the applications of this paper to Selmer groups the existence 
of Galois representations attached to
algebraic automorphic forms on $\U(m)$ with $m\geq 4$, we do not assume
that the expected ones are irreducible, but instead our arguments prove this
irreducibility for some of
them.}, but even irreducible on every proper artinian thickening on the
point $x$ inside $\cal E$ (recall however that $\rho_x=1\oplus 
\chi\oplus \rho$ is
reducible). In other words, the reducibility locus of the 
family is schematically equal to the point $x$. It should be pointed out here that the situation is quite orthogonal to the one of Iwasawa's main conjecture
(see~\cite{mazurwiles}, \cite{wiles}), 
for which there is a big known part in the reducibility locus at $x$ (the
Eisenstein part), and this locus cannot be controlled {\it a priori}. In our
case, this fact turns out to be one of
the main ingredients in order to get some geometric control on the size of
the subspace we construct in $\selr$, and it is maybe the main reason
why our points $x$ are quite susceptible to produce independent elements in
$\selr$.
The question of whether we should (or not) expect to construct the full Selmer group
of $\rho$ at $x$ remains a very interesting mystery, whatever the answer
may be. Although it might not be easy to decide it even in explicit examples
(say with $L(\rho,s)$ vanishing at order $>1$ at its center), 
our geometric criterion reduces this question to some computations of spaces of
$p$-adic automorphic forms on explicit definite unitary groups, which should
be feasible. Last but not least, we hope
that it could be possible to relate the geometry of $\cal E$ at
$x$ (which is built on spaces of $p$-adic automorphic forms) to the $L$-function (or
rather a $p$-adic $L$-function) of $\rho$, so that our results could be used in
order to prove the "lower bound of Selmer group" part of the Bloch-Kato conjecture.
However, this is beyond the scope of this book.

\bigskip

The four first chapters form a detailed study of $p$-adic families of
Galois representations, especially near reducible points, 
and how their behavior is related to Selmer groups. There are no references to
"automorphic forms" in them, as opposite to the following chapters 5 to 9
which are devoted to the applications to eigenvarieties. In what
follows, we describe very briefly the contents of each of the different
chapters by focusing on the way they fit in the general theme of the
book. As they contain quite a number of results of independent
interest, we invite the reader to consult then their respective introductions for more details. 

When we deal with families of representations 
$(\rho_x)_{x \in X}$ of a group $G$ (or an algebra) over a "geometric" space
$X$, there are two natural notions to consider. The most obvious one is the
data of a "big" representation of $G$ of a locally free sheaf of $\anneau_X$-modules whose
evaluations at each $x \in X$ are the $\rho_x$. Another one, visibly weaker, 
is the data of a "trace map" $G \longrightarrow \anneau(X)$ whose evaluation at
each $x \in X$ is
$\tr(\rho_x)$; these abstract traces are then called {\it pseudocharacters}
(or pseudorepresentations). As a typical example, when we are interested in the space of all representations (up to isomorphisms) 
of a given group say, we usually only get a pseudocharacter on that full space. This is 
what happens also for the family of Galois representations on the eigenvarieties. 
When all the $\rho_x$ are irreducible, the two definitions turn out to be essentially the same, 
but the links between them are much more 
subtle around a reducible $\rho_x$ and they are related to the extensions between the
irreducible constituents of $\rho_x$, what we are interested in. 

Thus our first chapter is a general study of 
pseudocharacters $T$ over a henselian local ring $A$ (having in view that the local rings 
of a rigid analytic space are henselian). 
There is no mention of a Galois group in all this chapter, and those results can be applied to any 
group or algebra. Most of our work is based on the assumption 
that the residual pseudocharacter $\bar T$ (that is, the pseudocharacter 
one gets after tensorizing $T$ by the residue field of $A$) 
is without multiplicity, so it may be reducible, which is fundamental, but
all its components appear only once. Under this hypothesis, 
we prove a precise structure theorem for $T$, describe the groups of 
extensions between the constituents of $\bar T$ we can construct from $T$, 
and define and characterize the {\it reducibility loci} of $T$ (intuitively the subscheme of $\spec A$
where $T$ has a given reducibility structure). We also discuss conditions
under which $T$ is, or cannot be, the trace of a true representation. This chapter provides the 
framework of many of our subsequent results. 
        
In the second chapter we study infinitesimal (that is, artinian) families 
of $p$-adic local Galois representations, 
and their Fontaine's theory, having in mind to characterize abstractly those coming from eigenvarieties.
A key role is played by the theory of $\fg$-modules over the Robba ring and 
Colmez' notion of {\it trianguline representation} \cite{colmeztri}. 
We generalize some results of Colmez to any dimension and with artinian
coefficients, giving in particular a fairly complete description 
of the trianguline deformation space of a {\it non critical} trianguline
representation (of any rank). For the applications to eigenvarieties, 
we also give a criterion for an infinitesimal family to be trianguline in terms of crystalline periods.

In the third chapter, we generalize a recent result of 
Kisin in \cite{kis} on the analytic continuation of crystalline periods in a family of local Galois
representations. This result was proved there for the strong definition of
families, namely for true representations of $\Gal(\Qpb/\Qp)$ on a locally free 
$\anneau$-module, and we prove it more generally for any torsion free
coherent $\anneau$-module. Our main technical tool is a method of descent by
blow-up of crystalline periods (which turns out to be rather general) and a
reduction to Kisin's case by a flatification argument. 

In the fourth chapter, we give our working definition of "$p$-adic families of
refined Galois representation", motivated by the families carried by eigenvarieties, and we apply to them the 
results of chapters $2$ and $3$. In particular, we are able in favorable cases to understand 
their reducibility loci in terms of the Hodge-Tate-Sen weight maps, and to
prove that they are infinitesimally trianguline. 

In the fifth chapter we discuss our main motivating conjecture relating the 
dimension of Selmer groups of geometric semi-simple Galois representations 
to the order of the zeros of their $L$-functions at integers. We are mainly
interested in ``one half'' of this conjecture, namely to
give a lower bound on the dimension of the Selmer groups, as well as in a very special case of it that we call the {\it sign conjecture}.
As was explained in \cite{Bellaichethese}, 
an important feature of the method we use is that we need as an input some results (supposedly simpler) about upper bounds of other Selmer groups. For the sign conjecture, we only need the vanishing of ${\rm Sel}(\chi)$ (for a quadratic imaginary field) which is elementary. However, we need more upper bounds results for our second main theorem, and we cannot prove all of them in general. Thus we formulate as hypotheses the results we shall need, which will appear as assumptions in the results of chapter 9. Using results of Kisin and Kato, we are able to prove all that we need in most cases when $n=2$, and in all the cases for $n=1$.

The sixth chapter contains all the results we need about the unitary groups, their automorphic forms, and the 
Galois representations attached to them, and 
includes a detailled discussion of Langlands and Arthur's conjectures. 
In particular, we formulate there the two hypotheses (AC($\pi$)) and (Rep($m$)) that we use in chapters 8 and 9. 

In the seventh chapter, we introduce and study in details the eigenvarieties of definite unitary groups and we prove the basic 
properties of the (sometimes conjectural) family of Galois representations that they carry. We essentially rely on one of us' thesis \cite{Ch}), and actually go a bit further on several respects. It furnishes a lot of interesting examples where all the concepts studied in this book occur, and provides also an important tool for the applications to Selmer groups. The first half of this chapter only concerns eigenvarieties and may be read independently, whereas the second one depends on chapters 1 to 5.

Finally, in chapters 8 and 9 we prove our main results, and we refer to those chapters for precise statements.

As a parallel goal, we made considerable efforts all along the redaction of this book to develop concepts and techniques adapted to a proper study of eigenvarieties. We hope that the reader will enjoy playing with them as much as we did.
\medskip

The leitfaden is as follows :
$$
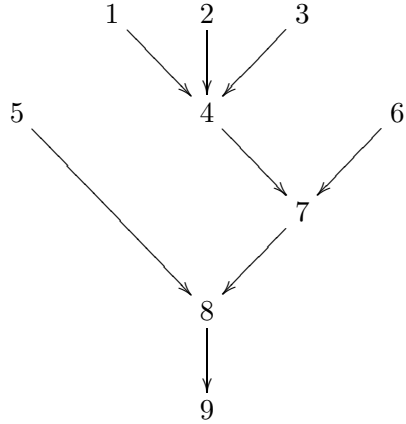
\xymatrix{ & 1\ar[dr] & 2\ar[d] & 3\ar[dl] & \\ 5 \ar[ddrr] & & 4\ar[dr] & & 6\ar[dl] \\ & & & 7 \ar[dl] & \\ & & 8 \ar[d] & & \\ & & 9 & &}$$ 
\medskip

It is a pleasure to thank  Laurent Berger, Brian Conrad, Laurent Clozel, Pierre Colmez, Jean-Francois Dat, Guy Henniart, Michael Harris, Kiran Kedlaya, Colette Moeglin and Claudio Procesi for many useful conversations. 

 Last but not least, the authors would like to thank Cl\'ementine, Sarah and \hbox{Valeria}: this book is also theirs.

\newpage
 
\section{Pseudocharacters, representations and extensions}  
\label{pseudocharacters}  
\subsection{Introduction}
This section is devoted to the local (in the sense of the \'etale topology) study of pseudocharacters $T$ satisfying a residual multiplicity freeness hypothesis. Two of our main objectives are to determine when those pseudocharacters  
come from a true representation and to prove the optimal generalization of  ``Ribet's lemma'' for them.  
   
Let us precise our main notations and hypotheses. Throughout this section, we will work with a pseudocharacter $T : R \longrightarrow A$ of dimension $d$,   
where $A$ is a {\it local henselian} commutative ring of residue field $k$ where $d!$ is invertible  
and $R$ a (not necessarily commutative)   
$A$-algebra\footnote{In the applications, $R$ will be the group algebra   
$A[G]$ where $G$ a group, especially a Galois group. However, it is important to keep this degree of generality as most of the statements concerning pseudocharacters are ring theoretic.}.   
To formulate our residual hypothesis,   
we assume\footnote{Since it is allowed to replace $k$ by a separable extension,   
this assumption is actually not a restriction.} that $T \otimes k : R \otimes k \longrightarrow k$ is the   
sum of $r$ pseudocharacters of the form $\tr \rhob_i$   
where the $\rhob_i$'s are absolutely irreducible   
representations of $R \otimes k$ defined over $k$.    
Our {\it residually multiplicity free} hypothesis is   
that the $\rhob_i$'s are two by two non isomorphic.   
In this context, ``Ribet's lemma'' amounts to determine   
how much we can deduce about the existence of non-trivial extensions   
between the representations $\rhob_i$ from the existence and   
irreducibility properties of $T$.   
Before explaining our work and results in more details,  
let us recall the history of those two interrelated themes :   
pseudocharacters and the generalizations of ``Ribet's lemma''.  
  
We begin with the original Ribet's lemma (\cite[Prop. 2.1]{Ribet}).  
Ribet's hypotheses are that $d=r=2$, and that $A$ is a complete discrete valuation   
ring. He works with a representation $\rho : R \longrightarrow M_2(A)$, but   
that is no real supplementary restriction   
since every pseudocharacter over a strictly complete discrete valuation   
ring is the trace of a true   
representation\footnote{We leave the proof of this assertion to the   
interested reader (use the fact that the Brauer group of any finite extension of $K$ is trivial, e.g. by \cite[XII \S2, especially exercise 2]{corpslocaux}).}.   
Ribet proves that if $\rho \otimes K$  
($K$ being the fraction field of $A$) is irreducible, then  
 a non-trivial extension of $\rhob_1$ by $\rhob_2$   
(resp. of $\rhob_2$ by $\rhob_1$) arises as a subquotient of $\rho$.   
This seminal result suggests numerous generalizations : we may   
wish to weaken the hypotheses on  the dimension $d$, the number of  
residual factors $r$, the ring $A$, and for more general $A$,  
to work with general pseudocharacters instead of representations.   
We may also wonder if we can obtain, under suitable hypotheses,  
extensions between deformations $\rho_1$ and $\rho_2$ of $\rhob_1$ and   
$\rhob_2$ over some suitable artinian quotient of $A$, not only over $k$.  
  
A big step forward is made in the papers by Mazur-Wiles and Wiles  
(\cite{mazurwiles}, \cite{wiles}) on Iwasawa's main conjecture.   
As their works is the primary source of inspiration for this section, let us explain it with some details  
(our exposition owes much to~\cite{harderpink}; see also   
\cite[\S2]{lissite}). They still assume $d=r=2$, but the ring $A$   
now is any finite flat reduced local $A_0$-algebra $A$,   
where $A_0$ is a complete discrete valuation ring.   
Though the notion of pseudocharacter at that time was still to be defined, their formulation amount to consider a   
pseudocharacter (not necessarily coming from a representation) $T: R \longrightarrow A$, where  
$R$ is the group algebra of a global Galois group. The pseudocharacter is supposed to be {\it odd} which implies our multiplicity free hypothesis. They introduce an ideal $I$ of $A$, which turns out to be the smallest   
ideal of $A$ such that $T \otimes A/I$ is   
the sum of two characters $\rho_1,\rho_2 :\  R \longrightarrow A/I$ deforming   
respectively $\rhob_1$ and $\rhob_2$.   
Assuming that $I$ has cofinite length $l$,   
their result is the construction of a finite $A/I$-module of length at least $l$ in    
$\Ext^1_{R/IR}(\rho_1,\rho_2)$. We note that it is not possible to determine the precise structure of this module, so we do   
not know if their method constructs,   
for example, $l$ independent extensions over $k$   
of $\rhob_2$ by $\rhob_1$ or, on   
the contrary, one ``big'' extension of $\rho_2$ by $\rho_1$ over the   
artinian ring $A/I$, that would generate a free $A/I$-module in   
$\Ext^1_{R/IR}(\rho_1,\rho_2)$.  
  
The notion of {\it pseudocharacter}   
was introduced soon after by Wiles in dimension 2 (\cite{wiles}), and by Taylor in full generality   
(\cite{Tay}), under the name of {\it  
pseudorepresentation}.   
Besides their elementary properties, the main   
questions that has been studied until now is   
whether they arise as the trace of a true   
representation.   
Taylor showed in 1990, relying on earlier   
results by Procesi,   
that the answer is always yes in the case   
where $A$ is an algebraically closed field   
of characteristic zero;   
this result was extended, with a different   
method, to  
any algebraically closed field (of   
characteristic prime to $d!$)   
by Rouquier. The question was settled   
affirmatively in 1996   
for any local henselian ring $A$, in the case   
where the residual   
pseudocharacter $\bar T$ is {\it absolutely   
irreducible},   
independently by Rouquier (\cite{Rou}) and   
Nyssen (\cite{Nys}).  
  
We now return to the progresses on Ribet's   
lemma.  
  
Urban's work (\cite{urban}) deals with the   
question of obtaining, using notations of the   
paragraph describing Mazur-Wiles   
modules, a {\it free} $A/I$-module of   
extensions of $\rho_1$ by $\rho_2$. His   
hypotheses are as follows : the dimension   
$d$ is arbitrary, but the number $r$ of residual factors is still $2$. The ring $A/I$ is an arbitrary artinian local ring,  
and the pseudocharacter $T$ is (over $A/I$) equal to $\tr \rho_1 + \tr \rho_2$, but he also assumes that $T$ comes from a   
true representation $\rho$ (at least over $A/I$), which moreover is modulo the maximal ideal of $A$ a non-trivial extension   
of $\rhob_1$ by $\rhob_2$. Then he proves that   
$\rho$ is indeed a non trivial extension of $\rho_1$ by $\rho_2$. Thus he obtains a more precise result than Mazur and Wiles, but with the much stronger assumption that his pseudocharacter comes from a representation that gives already  
the searched extension modulo the maximal ideal.  
Our work (see~\S\ref{exampler2}) will actually show that the   
possibility of producing a free $A/I$-module of extensions as he does   
depends fundamentally on that hypotheses, which is very hard to   
check in practice excepted when $A$ is a discrete valuation ring,   
or when $T$ allows to construct only one extension of $\rhob_1$ by  
$\rhob_2$.

One of us studied (\cite{bellaiche}) the case of an arbitrary number of residual factors $r$ (and   
an arbitrary $d$) but as Ribet with $A$ a complete discrete valuation ring. The main feature   
here is that the optimal result about extensions becomes more combinatorially involved.   
Assuming that $\rho$ is generically irreducible, we can say nothing  
about the vanishing of an individual space of extensions   
$\Ext^1_{R\otimes k}(\rhob_i,\rhob_j)$. What we can say is that there are enough   
couples $(i,j)$ in $\{1,\dots,r\}^2$ with non-zero  
$\Ext^1_{R\otimes k}(\rhob_i,\rhob_j)$ for the graph  
drawn by the oriented edges $(i,j)$ to be connected.    
This result was soon after extended to deal with extensions over   
$A/I$ assuming the residual multiplicity one hypothesis, in a joint work with P. Graftieaux in \cite{BG}.   
The combinatorial description of extensions we will obtain here is reminiscent of the results of that work.  
  
Let us conclude those historical remarks by noting that two basic questions are not answered by all the results mentioned above : about Ribet's lemma, is it possible to find reasonable hypotheses so that two independent extensions  
of $\rhob_1$ by $\rhob_2$ over $k$ exist ? About pseudocharacters (over a strictly local henselian ring $A$, say), for which conceptual reasons a pseudocharacter may not be the trace of a true representation ?  
  
In this work, we will obtain the most general form of Ribet's lemma (for any $A$ and $T$, and implying all the ones above) as well as a satisfactory answer to both questions above and more. Indeed we will derive a precise structure theorem for residually multiplicity free pseudocharacters, and using this result we are able   
to understand precisely and to provide links (some expected, others   
rather surprising) between the questions of   
when does a pseudocharacter come from a representation,  
how many extensions it defines, and how its (ir)reducibility behaves with respect to changing the ring $A$ by a quotient of it.  
  
We now explain our work, roughly following the order of the   
subparts of this section.   
  
The first subpart \S\ref{appendix} deals with generalities on   
pseudocharacters. There $A$ is not local henselian but can be any commutative ring.   
Though this part is obviously influenced by~\cite{Rou}, we have tried   
to make it self-contained, partly for the convenience of the reader,   
and partly because we anyway needed to improve and generalize most   
of the arguments of Rouquier.  We begin by recalling Rouquier's definition of a pseudocharacter of dimension   
$d$. We then introduce the notion of  
 {\it Cayley-Hamilton} pseudocharacter $T$ : it means that every $x$ in $R$ is killed by its ``characteristic polynomial''   
whose coefficients are computed from the  $T(x^i)$, $i=1,\dots,d$.   
This notion is weaker than the notion of faithfulness that was used by Taylor and Rouquier,   
but it is stable by many operations, and this fact   
allows us to give more general statements with often simpler proofs. This notion is also closely related with the {\it  
Cayley-Hamilton trace algebras} studied by Procesi (see \cite{Proc2}).   
Every $A$-algebra $R$ with a pseudocharacter   
$T$ has a bunch of quotients on which $T$ factors and becomes Cayley-Hamilton, the smallest of those being the unique faithful quotient $R/\ker T$. We also prove results concerning idempotents,   
and the radical of an algebra with a Cayley-Hamilton pseudocharacters, that will be useful in our analysis of residually   
multiplicity free pseudocharacters.   
Finally, we define and study the notion of Schur functors of a pseudocharacter.  
  
In \S\ref{lesgma} and \S\ref{structuremult1}, we study the structure of the residually multiplicity free pseudo\-characters   
over a local henselian ring $A$. We introduce the notion of {\it generalized matrix algebra}, or briefly {\it GMA}, over $A$.   
Basically, a GMA over $A$ is an $A$-algebra whose elements are square matrices (say, of size $d$) but where we allow the non diagonal entries to be elements of arbitrary $A$-modules instead of $A$ - say the $(i,j)$-entries are elements of the given $A$-module $A_{i,j}$. Of course, to define the multiplications of such matrices, we need to suppose given some morphisms  
$A_{i,j} \otimes_A A_{j,k} \longrightarrow A_{i,k}$   
satisfying suitable rules. The result motivating the introduction of   
GMA is our main structure result (proved in \S\ref{structuremult1}), namely : if $T: R \longrightarrow A$ is a residually multiplicity free pseudocharacter, then every Cayley-Hamilton quotient of $R$ is a GMA. Conversely, we prove that the trace function on any GMA  
is a Cayley-Hamilton pseudocharacter, which is residually multiplicity free if we assume that $A_{i,j}A_{j,i} \subset m$ (the maximal ideal   
of $A$) for every $i \neq j$, which provides us with many non trivial examples of such pseudocharacters.   
This result is a consequence of the main theorem of our study   
of GMA's (\S\ref{lesgma}) which states that any GMA over $A$ can be embedded, compatibly with the traces function,   
in an algebra $M_d(B)$ for some explicit commutative $A$-algebra $B$. Those two results take place in the long-studied topic of embedding an  abstract algebra in a matrix algebra. It should be compared to a result of Procesi (\cite{Proc2}) on embeddings of trace algebras in matrix algebras : our results deal with less general algebras $R$, but with more general $A$, since we avoid the characteristic zero   
hypothesis that was fundamental in Procesi's invariants theory methods.   
  
In \S\ref{redlocext}, we get the dividends of our rather abstract work on the structure of residually multiplicity free pseudocharacters. Firstly, for such a pseudocharacter, and for every partition of $\{1,\dots,r\}$ of cardinality $k$, we prove that there exists a greatest subscheme of $\spec{A}$ on which $T$ is a sum of $k$   
pseudocharacters, each of which being residually the sum of $\tr \rhob_i$  
for $i$ belonging to an element of the partition.   
We also show that this decomposition of $T$ as a sum of $k$ such pseudocharacters is unique,   
and that those subschemes of $\spec(A)$  
do not change if $R$ is changed into a quotient on which $T$ factors. They are called the {\it reducibility loci}\footnote{We  
stress the reader that we could not define those loci without the   
assumption of residually multiplicity one (see \cite{BCnonramf}).}   
attached to the given partition, and they will become one of our main object of study in section~\ref{families}. Moreover, if $S$ is any Cayley-Hamilton quotient of $R$, hence a GMA defined by some $A$-modules $A_{i,j}$'s, we give a simple and explicit   
description of the ideals of the reducibility loci in terms of the $A_{i,j}$.  
  
Secondly, we construct submodules (explicitly described in   
terms of the modules $A_{i,j}$) of the extensions modules $\Ext_R(\rho_j,\rho_i)$. This is our version of  
``Ribet's lemma'', as it provides a link between non-trivial extensions of $\rho_j$ by $\rho_i$ and the irreducibility properties of $T$   
encoded in its reducibility loci, and we show that it is in any reasonable   
sense optimal.   
  
Nevertheless, and though its simplicity, this result may not seem  
perfectly satisfactory, as it involves the unknown modules $A_{i,j}$'s.   
It may seem desirable to get a more direct link  
between the module of extensions we can construct   
and the reducibility ideals, solving out the modules   
$A_{i,j}$. However, this is actually a very complicated task,   
that has probably no nice answer in general, as it involves in the same   
times combinatorial and ring-theoretical difficulties :   
for the combinatorial difficulties, and how  
they can be solved (at a high price in terms of simplicity of statements)  
in a context that is ring-theoretically trivial   
(namely $A$ a discrete valuation ring),   
we refer the reader to \cite{BG} ; for the ring-theoretical difficulties  
in a context that is combinatorially trivial   
we refer the reader to \S\ref{exampler2}.   
In this subpart, we make explicit in the simple case $r=2$ the subtle relations our results implies between, for a given   
pseudocharacter $T$, the modules of extensions that $T$ allows to construct, the existence of a representation whose trace is $T$, the reducibility ideal of $T$ and the ring-theoretic properties of $A$. We also give some criteria for our method to construct several   
{\it independent} extensions. Finally, let us say that the final sections of this paper will   
show that our version of Ribet's lemma, as stated, can actually easily be used in practice.  
   
In \S\ref{representations}, we determine the local henselian rings $A$ on which   
{\it every} residually multiplicity free pseudocharacter comes from a representation.   
The answer is surprisingly simple, if we restrict ourselves to noetherian $A$.   
Those $A$'s are exactly the unique factorization domains. The proof relies on our structure result and its converse.   
  
Finally, in \S\ref{pseudosymmetry}, we study pseudocharacters having a property of symmetry   
of order two (for example, autodual pseudocharacter). It is natural to expect to   
retrieve this symmetry on the modules of extensions we have constructed, and this is what this   
subsection elucidates. Our main tool is a (not so easy) lemma about lifting idempotents ``compatibly   
with an automorphism or an anti-automorphism of order two'' which may be of independent interest. \ps  
  
It is a pleasure to acknowledge the influence of all the persons mentioned in the historical part of this introduction.  
Especially important to us have been the papers and surveys of Procesi, as well as a few but illuminating discussions   
with him, either at Rome, the ENS, or by email.

\subsection{Some preliminaries on pseudocharacters}  
\label{appendix}  
\subsubsection{Definitions} \label{definitionspseudo}Let $A$\index{A@$A$, a commutative ring|(} be a commutative ring\footnote{In all the paper, {\it rings} and {\it algebras} are associative and have a unit, and a ring homomorphism preserves the unit.} and $R$ an  
$A$-algebra (not necessarily commutative). \index{RR@$R$, an $A$-algebra} Let us recall the  
definition of an {\it $A$-valued pseudocharacter on $R$} introduced by R. Taylor in \cite[\S1]{Tay}.   
Let $T: R \longrightarrow A$ \index{T@$T: R \longrightarrow A$, a   
pseudocharacter}   
be an $A$-linear map which is central, that is such that $T(xy)=T(yx)$ for all $x, \, y \in R$. For  
each integer $n\geq 1$, define a map $S_n(T): R^n \longrightarrow A$ by  
\newcommand{\ST}{T^\sigma}  
$$S_n(T)(x)=\sum_{\sigma \in {\mathfrak{S}}_n} \varepsilon(\sigma) \ST(x),$$  
where $\ST: R^n \longrightarrow A$ is defined as follows. Let  
$x=(x_1,...,x_n) \in R^n$. If $\sigma$ is a cycle, say $(j_1,...,j_m)$, then  
set $\ST(x)=T(x_{j_1}\cdots x_{j_m})$, which is well defined. In  
general, we let $\ST(x)=\prod_{i=1}^r T^{\sigma_i}(x)$, where  
$\sigma=\prod_{i=1}^r \sigma_i$ is the decomposition in cycles of  
the permutation $\sigma$. We set $S_0(T):=1$. \par  
The central function $T$ is called a {\it pseudocharacter} on $R$ if there exists  
an integer $n$ such that $S_{n+1}(T)=0$, and such that $n!$ is invertible in $R$.   
The smallest such $n$ is then called the dimension of $T$, and it satisfies $T(1)=n$ (see  
Lemma \ref{dimpseudo} (2))\footnote{The definition of a pseudocharacter of dimension $n$ used here looks   
slightly more restrictive than the one introduced in \cite{Tay} or \cite{Rou},   
as we assume that $n!$ is invertible. This assumption on $n!$ is first  
crucial to express the Cayley-Hamilton theorem from the trace, which is a  
basic link between pseudocharacters and true representations, and also to avoid a   
strange behavior of the dimension of pseudocharacters with base change. Note that Taylor's theorem concerns the case where $A$ is a  
field of characteristic $0$, hence $n!$ is invertible. Moreover, Lemma 2.13 of \cite{Rou} does not hold when $d!$ is not invertible, hence this hypothesis should be added {\it a priori} in the hypothesis of Lemma 4.1 and Theorem 5.1 there (see also remark \ref{fauterouquier}).}.   
\par  
These notions apply in the special case where $R=A[G]$ for some group (or monoid) $G$. In this case,  $T$ is uniquely determined by the data of its restriction to $G$ (central, and satisfying $S_{n+1}(T)=0$ on $G^{n+1}$). \par  
If $T:R \longrightarrow A$   
is an $A$-valued pseudocharacter on $R$ of dimension $d$ and if $A'$ is a commutative $A$-algebra, then the induced linear map $T \otimes A' :\, R \otimes A' \longrightarrow A'$ is an $A'$-valued pseudocharacter on $R$ of dimension $d$.   
  
\subsubsection{Main example.} \label{mainexample}   
 Let $V:=A^d$ and $\rho: R \longrightarrow \End_A(V)$   
be a morphism of $A$-algebras. For each $n\geq 1$, $V^{\otimes_A n}$ carries an $A$-linear  
representation of $\mathfrak{S}_n$ and a diagonal action of the underlying multiplicative monoid   
of $R^n$. If $e=\sum_{\sigma \in {\mathfrak{S}}_n}\varepsilon(\sigma)\sigma \in A[{\mathfrak{S}}_n]$, then a computation\footnote{For instance, reducing to the case where $R=\End_A(V)$ and   
using the polarization identity for symmetric multi-linear forms, it suffices to check it when $x=(y,y,\cdots,y)$ (see also \cite[\S1.1]{Proc1}, \cite[prop. 3.1]{Rou}).}  
shows that for $x \in R^n$, $$\tr (x e | V^{\otimes_A n})=S_n(\tr(\rho))(x).$$   
As $e$ acts by $0$ on $V^{\otimes_A n}$ if $n>d$, the central function  
$T:=\tr(\rho)$ is a pseudocharacter of dimension $d$ (assuming that   
$d$ is invertible in $A$). Moreover, when $\rho$ is an isomorphism, an easy computation using standard matrices shows that $T$ is the unique $A$-valued pseudocharacter of dimension $d$ of $R=M_d(A)$. By faithfully flat descent,   
these results hold also when $\End_A(V)$ is replaced by any Azumaya algebra of  
rank $d^2$ over $A$, and when $\tr$ is its reduced trace. \ps  
  
	It is known that a pseudocharacter $T: R \longrightarrow A$ arises from a representation of $R$ in an Azumaya algebra as above, essentially uniquely, under some suitable irreducibility hypothesis: when $A$ is a field and   
$T$ is absolutely irreducible (\cite[thm. 1.2]{Tay} in characteristic $0$, relying on the work of Procesi in \cite{Proc1}, and \cite[thm. 4.2]{Rou} in any characteristic by a different proof, see these references for the definition of absolutely irreducible), and more generally when $T \otimes A/m$ is absolutely irreducible for all $m \in {\rm Specmax}(A)$   
(\cite[thm. 5.1]{Rou}, \cite{Nys} when $A$ is local henselian). One main goal of this section is to study the general case where $A$ is local and $T \otimes A/m$ is reducible and satisfies   
a multiplicity one hypothesis (def. \ref{defmult1}, \S\ref{representations})).   
  
\subsubsection{The Cayley-Hamilton identity and Cayley-Hamilton  
pseudocharacters.}\label{operations}   
Let $T: R \longrightarrow A$ be a pseudocharacter of dimension $d$. For $x \in R$, let  
$$P_{x,T}(X):=X^d+\sum_{k=1}^{d} \frac{(-1)^k}{k!} S_k(T)(x,\dots,x) X^{d-k} \in A[X].$$  
In the example given in \S\ref{mainexample}, $P_{x,T}(X)$ is the usual characteristic   
polynomial of $x$. We will say that $T$ is {\it Cayley-Hamilton} if it satisfies the {\it Cayley-Hamilton identity}, that is if   
$$\text{ for all }x \in R,\ \ P_{x,T}(x)=0.$$  
 In this case, $R$ is integral over $A$.   
The algebra $R$ equipped with $T$ is then a {\it Cayley-Hamilton algebra} in the sense of C. Procesi \cite[def. 2.6]{Proc2}. \ps \newcommand{\CH}{{\rm CH}}  
An important observation is that for a general pseudocharacter   
$T: R \longrightarrow A$ of dimension $d$,   
the map $R \longrightarrow R$, $x\mapsto P_{x,T}(x)$, is the evaluation at   
$(x,\dots,x)$ of a $d$-linear symmetric map $\CH(T): R^d \longrightarrow R$, explicitly given by:  
  
$$\CH(T)(x_1,\dots,x_d):=\frac{(-1)^d}{d!}\underset{I, \sigma}{\sum}(-1)^{|I|}S_{d-|I|}(T)(\{x_i, i \notin I\})x_{\sigma(1)}\cdots x_{\sigma(|I|)},$$  
where $I$ is a subset of $\{1,\dots,d\}$ and $\sigma$ a bijection from   
$\{1,\dots,|I|\}$ to $I$. A first consequence of the polarization identity (\cite[Alg., Chap. I, \S8, prop. 2]{Bki}, applied to the ring $\Sym^d_A(R)$)   
is that $T$ is Cayley-Hamilton if and only if $\CH(T)=0$. In particular, if $T$ is Cayley-Hamilton then   
for any $A$-algebra $A'$, $T \otimes A'$ is also Cayley-Hamilton. \ps  
  
In the same way, we see that for $x_1,\dots,x_{d+1} \in R$, we have \begin{equation}\label{pseudocayley}  
S_{d+1}(T)(x_1,\dots,x_{d+1})=d!T(\CH(T)(x_1,\dots,x_d)x_{d+1}),\end{equation}   
hence a good way to think about the identity $S_{d+1}(T)=0$   
defining a pseudocharacter is to see it as a   
polarized, $A$-valued, form of the Cayley-Hamilton identity.\ps  
  
\subsubsection{Faithful pseudocharacters, the kernel of a pseudocharacter.} We recall that the   
{\it kernel} of $T$ is the two-sided ideal $\Ker T$ of $R$  
defined by $$\ker T:=\{x \in R,\  \forall y \in R, T(xy)=0\}$$   
\index{kert@$\ker T$, the kernel of the pseudocharacter $T$}  
$T$ is said to be {\it faithful} when $\ker{T}=0$. If $R \longrightarrow S$ is a surjective morphism of   
$A$-algebras whose kernel is included in $\Ker T$, then $T$ factors uniquely as a pseudocharacter   
$T_S: S\longrightarrow A$, which is   
still of dimension $d$, and which will be often denoted by $T$. In   
particular, $T$ induces a faithful   
pseudocharacter on $R/\ker T$.\ps  
  
	Assume that $T$ is faithful,   
then $T$ is Cayley-Hamilton by the formula (\ref{pseudocayley}) above (see also \cite[lemme 2.12]{Rou}).  
 More generally, let $(T_i)_{i=1}^r$ be a family of pseudocharacters $R \longrightarrow A$ such   
that the integer $d:=\dim T_1+ \dots + \dim T_r$ is invertible in $A$. Then $T:=\sum_i T_i$  
 is a pseudocharacter of dimension $d$, and for all $x \in R$, $$P_{x,T}=\prod_{i=1}^rP_{x,T_i}$$  
 (we may assume that $r=2$, in which case it follows from \cite[lemme 2.2]{Rou}).  
 As a consequence, $P_{x,T}(x) \in (\Ker T_1)(\Ker T_2) \cdots (\Ker T_r) \subset \bigcap_i \Ker T_i$, hence  
$T: R/(\cap_i \Ker T_i) \longrightarrow A$ is Cayley-Hamilton. The following lemma is obvious from the formula of $P_{x,T}(X)$, but useful.  
  
\begin{lemma}\label{kerrad} Let $T: R \longrightarrow A$ be a Cayley-Hamilton pseudocharacter of dimension $d$, then for each $x \in \Ker T$ we have $x^d=0$. In particular $\Ker T$ is a nil ideal, and is contained in the Jacobson ideal of $R$.  
\end{lemma}  
  
\begin{remark}	If $A'$ is an $A$-algebra and $T$ is faithful, it is not true in general that $T \otimes A'$ is still faithful. Although we will not need it in what follows, let us mention that this is however the case when $A'$ is projective as an $A$-module (so e.g. when $A$ is a field),   
or when $A'$ is flat over $A$ and either $R$ of finite type over $A$ (see \cite[prop. 2.11]{Rou}) or $A$ is noetherian (mimic the proof {\it loc. cit.} and use that $A^X$ is flat over $A$ for any set $X$). \end{remark}  
\subsubsection{Cayley-Hamilton quotients}  
\label{CHquotient}  
\begin{definition}   
Let $T : R \longrightarrow A$ be a pseudocharacter of dimension $d$. Then a quotient $S$ of $R$ by a two-sided   
ideal of $R$ which is included in $\ker T$, and such that the induced pseudocharacter $T: S \longrightarrow A$ is Cayley-Hamilton, is called a   
{\it Cayley-Hamilton quotient} of $(R,T)$.  
\end{definition}  
\index{S@$S$, a Cayley-Hamilton quotient of $(R,T)$}  
\begin{example}  
\begin{itemize}   
\item[(i)] $R/\ker T$ is a Cayley-Hamilton quotient of $(R,T)$.  
\item[(ii)] Let $I$ be the two-sided ideal of $R$ generated by the elements $P_{x,T}(x)$ for all $x \in R$. Then  
$S_0:=R/I$ is a Cayley-Hamilton quotient of $(R,T)$. Indeed, $I \subset \ker T$ by~(\ref{pseudocayley})   
and $T$ is obviously Cayley-Hamilton on $S_0$.    
\item[(iii)] Let $B$ be a commutative $A$-algebra and $\rho : R \longrightarrow M_d(B)$ be a representation such that   
$\tr \circ \rho=T$. Then $\rho(R)$ is a Cayley-Hamilton quotient on $T$. Indeed, $\ker \rho$ is obviously  
included in $\ker T$ and $T$ is Cayley-Hamilton on $\rho(R)$ by the usual Cayley-Hamilton theorem.  
\end{itemize}  
\end{example}  
  
The Cayley-Hamilton quotients of $(R,T)$ form in a natural way a category : morphisms are $A$-algebra morphisms   
which are compatible with the morphism from $R$. Thus in that category   
any morphism  $S_1 \longrightarrow S_2$ is surjective, and   
has kernel $\ker T_{S_1}$ which is a nil ideal by Lemma \ref{kerrad}. Note that $S_0$ is the initial object and $R/\ker T$ the final object of that category.

\subsubsection{Two useful lemmas on pseudocharacters.} \label{dimension} Let $T:R \longrightarrow A$ be a pseudocharacter of dimension $d$. Recall that an element $e \in R$ is said to be idempotent if $e^2=e$. The subset $eRe \subset R$ is then  
an $A$-algebra whose unit element is $e$.  
  
\begin{lemma} \label{dimpseudo} Assume that $\spec(A)$ is connected. \begin{itemize}  
\item[(1)] For each idempotent $e \in R$, $T(e)$ is an integer less   
or equal to $d$.   
\item[(2)] We have $T(1)=d$. Moreover, if $A'$ is any $A$-algebra,   
the pseudocharacter $T \otimes A'$ has dimension $d$.  
\item[(3)] If $e \in R$ is an idempotent, the restriction $T_e$ of $T$ to the $A$-algebra $eRe$ is a pseudocharacter of dimension $T(e)$.  
\item[(4)] If $T$ is Cayley-Hamilton (resp. faithful), then so is $T_e$.  
\item[(5)] Assume that $T$ is Cayley-Hamilton.    
If $e_1, \dots, e_r$ is a family of (nonzero) orthogonal idempotents of $R$, then $r \leq d$.  
Moreover, if $T(e)=0$ for some idempotent $e$ of $R$, then $e=0$.  
\end{itemize}  
\end{lemma}  
\begin{pf}  
Let us prove (1). By definition of $S_{d+1}(T)$ and \cite[cor. 3.2]{Rou},   
\begin{equation} S_{d+1}(T)(e,e,\dots,e)=\sum_{\sigma\in\got{S}_{d+1}}\varepsilon(\sigma)T(e)^{|\sigma|}=T(e)(T(e)-1)\cdots(T(e)-d)=0 \end{equation}  
\noindent in $A$, where $|\sigma|$ is the number of cycles of $\sigma$.   
The discriminant of the split polynomial $X(X-1)\cdots (X-d) \in A[X]$ is $d!$ , hence is invertible in $A$. As $\spec(A)$ is connected, we get that $T(e)=i$ for some $i\leq d$. This proves (1). \par   
To prove (2), apply (1) to  
$e=1$. We see that $T(1)=i$ is an integer less than $d$. But following the  
proof of \cite[prop. 2.4]{Rou}, there is $x \in A-\{0\}$ such that  
$x(T(1)-d)=0$, Then $x(i-d)=0$, and because $i-d$ is invertible if non zero,  
we must have $i=d=T(1)$. In particular,   
$S_d(T)(1,1,\dots,1)=T(1)(T(1)-1)\dots(T(1)-d+1)=d!$ is invertible, hence $S_d(T  
\otimes A')(1,\dots,1)$ is non zero, which proves (2).\par Let  
$T_e:=T_{|eRe}: eRe \longrightarrow A$. For all $n$, we have  
$S_n(T_e)=S_n(T)_{|(eRe)^{n+1}}$, so that $T_e$ is a pseudocharacter of  
dimension $\leq d$. As $e$ is the unit of $eRe$, and $T(e)!$ is invertible  
in $A$ by (1), part (2) implies that $\dim T_e=T(e)$. \par If $x \in eRe$ and  
$y \in R$, then $T(xy)=T(exey)=T(xeye)=T_e(xeye)$, hence $T_e$ is faithful  
if $T$ is.  Assume now that $T$ is Cayley-Hamilton and fix $x \in eRe$. Let us compute  
$$e \CH(T)(x,\dots,x,(1-e),\dots,(1-e))$$ where $x$ appears  
$r:=T(e)$ times. As $x(1-e)=e(1-e)=0$, we see that the only nonvanishing  
terms defining the sum above are the ones with $(I,\sigma)$   
satisfying $|I|\leq r$ and $\sigma(\{1,\dots,|I|\}) \subset \{1,\dots,r\}$.  
For such a term, it follows from \cite[lemme 2.5]{Rou} that  
$$S_{d-|I|}(T)(\{x_i, i \notin  
I\})=S_{r-|I|}(T)(x,\dots,x)S_{d-r}(T)(1-e,\dots,1-e).$$   
As we have seen in proving part (2) above, and by (3),   
$S_{d-r}(T)(1-e,\dots,1-e)=S_{d-r}(T_{1-e})(1-e,\dots,1-e)=(d-r)!$ is  
invertible. We proved that:  
$$e\CH(T)(x,\dots,x,(1-e),\dots,(1-e))=\frac{(d-r)!^2}{d!}\CH(T_e)(x,\dots,x)e,$$  
hence $T_e$ is Cayley-Hamilton if $T$ is.  
\par  
Let us prove (5), we assume that $T$ is Cayley-Hamilton. Let $e$ be an idempotent of $R$. If $e$ satisfies $T(e)=0$, then we see that $P_{e,T}(X)=X^d$, hence $e^d=e=0$ by the  
Cayley-Hamilton identity. As a consequence of (1), if $e$ is  
nonzero then $T(e)$ is invertible. Assume now by contradiction that  
$e_1, \dots, e_{d+1}$ is a family of orthogonal nonzero idempotents of $R$. Then  
we get that $S_{d+1}(e_1,\dots,e_{d+1})=T(e_1)\cdots T(e_{d+1})$, which has  
to be invertible and zero, a contradiction. \end{pf}   
  
\begin{remark} \label{fauterouquier} Lemma 2.14 of \cite{Rou} is obviously incorrect  
as stated, and must be replaced by the part (5) of the above lemma (it is used in the proofs of Lemma 4.1 and Theorem 5.1 there).    
\end{remark}  
  
\medskip  
We conclude by computing the Jacobson radical of $R$ when $T$ is  
Cayley-Hamilton.  In what follows, $A$ is a   
local ring with maximal ideal $\m$ and residue field $k:=A/\m$. We will denote by $\bar R$ the $k$-algebra $R\otimes_A k= R/\m R$, and by   
$\bar T$ the pseudocharacter $T \otimes k : \bar R \longrightarrow k$.   
  
\begin{lemma} \label{radical}  
Assume that $T$ is Cayley-Hamilton. Then the kernel of the canonical surjection  
$R \longrightarrow  \Rb/ \ker \Tb$ is the Jacobson radical $\rad(R)$ of $R$.  
\end{lemma}  
\begin{pf}  
Let $J$ denote the kernel above, it is a two-sided ideal of $R$. By \cite[lemma 4.1]{Rou} (see precisely the sixth paragraph of the proof there), $\Rb/(\ker \Tb)$  
 is a semisimple $k$-algebra, hence $\rad(R) \subset J$. \par  
Let $x \in J$ ; we will show that $1+x \in R^*$. We  
have $T(xy) \in \m, \forall y \in R$, hence $T(x^i) \in \m$  
for all $i$, so that by the Cayley-Hamilton identity $x^d \in m(A[x])$. Let us consider the  
commutative finite $A$-algebra $B:=A[x]$. Then $B$ is local with maximal ideal  
$(\m,x)$, as $B/\m B$ is. As a consequence, $1+x$ is invertible in $B$, hence in  
$R$. As $J$ is a two-sided ideal of $R$ such that $1+J \subset R^*$, we have  
$J \subset \rad(R)$.   
\end{pf}  
  
\subsubsection{Tensor operations on pseudocharacters}  
  
\label{tenspseudo} In this section we assume  that $A$ is a $\Q$-algebra.   
All the tensor products involved below are assumed to be over $A$. \ps  
Let $R$ be an $A$-algebra, $T: R \longrightarrow A$ be a pseudocharacter of dimension $d$, and $m$ a positive integer.  
We define $T^{\otimes m}:R^{\otimes m} \longrightarrow A$ as the $A$-linear form that satisfies   
\begin{eqnarray} \label{puretens}  
T(x_1 \otimes \dots \otimes x_m)=T(x_1)  
\dots T(x_m).  
\end{eqnarray}  
Let us denote by $R^{\otimes m}[\got{S}_m]$  the twisted group algebra of $\got{S}_m$ over $R^{\otimes m}$ satisfying $$\sigma\cdot  x_1\otimes \cdots \otimes x_m = x_{\sigma(1)}\otimes \cdots \otimes x_{\sigma(m)}\cdot \sigma.$$ We can then extend $T^{\otimes m}$ to an $A$-linear map $R^{\otimes m}[\got{S}_m] \longrightarrow A$ by setting $$U(x_1\otimes \cdots \otimes x_m\cdot \sigma):=T^{\sigma}(x_1,\cdots,x_m)$$ (see \S\ref{definitionspseudo}, note that this map coincides with $T^{\otimes m}$ on the subalgebra $R^{\otimes m}$).   
  
\begin{prop}\label{tensprodpseudo} $T^{\otimes m}$ and $U$ are both pseudocharacters of dimension $d^m$.   
\end{prop}  
\begin{pf} By \cite{Proc2}, there is a commutative $A$-algebra $B$, with $A \subset B$,  
and a morphism $\rho : R \longrightarrow M_d(B)=\End_B(B^d)$ of $A$-algebras  
such that $\tr \rho (x)=T(x) 1_B$ for every $x \in R$. Let $\rho^{\otimes m} : R^{\otimes m} \longrightarrow \End_B((B^d)^{\otimes m})=M_{d^m}(B)$ be the $m$th   
tensor power of $\rho$. The equality $\tr \rho^{\otimes m}(z)=T^{\otimes m} (z) 1_{B}$ follows from~(\ref{puretens}) for pure tensors $z \in R^{\otimes m}$ and then by $A$-linearity for all $z$. We deduce $\tr \rho^{\otimes m} = T$ as   
$A \subset B$.   
Thus $T^{\otimes m}$ is a pseudocharacter, being the trace of a   
representation.  
  
We can extend the morphism $\rho^{\otimes m}: R^{\otimes m}   
\longrightarrow \End_B((B^d)^{\otimes m})$ into a morphism   
$\rho' : R^{\otimes m}[\got S_m] \longrightarrow \End_B((B^d)^{\otimes m})$  
by letting $\got S_m$ act by permutations on the $m$ tensor components of $(B^d)^{\otimes m}$. It is an easy computation to check  
that the trace of $\rho'$ is $U$. So $U$ is a pseudocharacter.  
  \end{pf}

\begin{remark} It should be true that more generally,   
if for $i=1,2$, $T_i: R_i \longrightarrow A$ is a pseudocharacter of dimension $d_i$, and if $T: R_1\otimes R_2 \longrightarrow A$   
is the $A$-linear map defined by $$T(x_1\otimes x_2)=T_1(x_1)T_2(x_2),$$ then $T$ is a pseudocharacter of dimension $d_1d_2$.   
It is probably possible to deduce directly the formula $S_{d_1d_2+1}(T)=0$ from the formulas $S_{d_i+1}(T_i)=0$, $i=1,2$,   
but we have not written down a proof\footnote{Note also that  the proof of the proposition above would break   
on the fact that if $\rho_i : R_i \longrightarrow M_d(B_i)$ are   
representations of trace $T_i$ given by \cite{Proc2}, it does not follow from $A \subset B_i$ that the map $A   
\longrightarrow B_1 \otimes_A B_2$ is injective, so that we cannot find a representation whose trace is $T$, but only a representation whose trace coincides with $T$ after reduction to the image of $A$ in $B_1 \otimes_A B_2$. However this line of reasoning would imply the result that $T$ is a pseudocharacter in two cases : if $A$ is reduced, because in that case, we can take $B_1=B_2$ equal to the product   
of algebraic closures of residue fields of all points of $\spec(A)$, and  
$\rho_i:R_i \longrightarrow M_d(B)$ be the "diagonal" representation deduced from $T_i$ ; and if $A$ is local henselian, $T_i$  
residually multiplicity free (see \S\ref{defmult1}), since in this case we may use Proposition~\ref{univ} to produce representations   
$\rho_i : R \longrightarrow M_d(B_i)$ of trace $T_i$ such that $A$ is a {\it direct factor} of $B_i$, so that  
we know that $A \subset B_1 \otimes_A B_2$.}.    
\end{remark}    
  
To conclude this paragraph, we give an application of the preceding proposition to the construction of the Schur functors of a given   
pseudocharacter in the case when $R:=A[G]$ with $G$ a group or a monoid. \par  
Let $T: A[G] \longrightarrow A$ be pseudocharacter and let $m\geq 1$ be an integer. There is a natural $A$-algebra   
embedding $$\iota_m: A[G] \longrightarrow R^{\otimes m}[\got{S}_m]=A[G^m\rtimes \got{S}_m]$$   
extending the diagonal map $G \longrightarrow G^m$. Let $e\in \Q[\got{S}_m]$ be any central idempotent. As the image of $\iota_m$ commutes with $\got{S}_m$, the map $$T^e: A[G] \longrightarrow A, \, \, x\mapsto U(\iota_m(x)e),$$ is a pseudocharacter by Proposition \ref{tensprodpseudo}. \par\ps  
  
\begin{remark}\label{remtensprod}\begin{itemize}  
\item[(i)] In the special case when $e=\frac{2}{m!}(\sum_{\sigma \in \got{S}_m} \varepsilon(\sigma)\sigma)$, then we set as usual $\Lambda^m(T):=T^e$. Note that $T^e(g)=\frac{2}{m!}S_m(T)(g,\cdots,g)$ for $g\in G$.  
\item[(ii)] It follows easily from the definitions that when $T$ (resp. $T_1$ and $T_2$) is the trace of a representation $G \longrightarrow \GL(V)$ (resp. of some representations $V_1$ and $V_2$), then $T^e$ (resp. $T_1T_2$) is the trace of the representation of $G$ on $e(V^{\otimes m})$ (resp. on $V_1\otimes V_2$).   
\end{itemize}  
\end{remark}  
\subsection{Generalized matrix algebras} \label{lesgma} ${}^{}$\ps

\noindent Let $d_1, \dots, d_r$ be nonzero positive integers, and   
$d:=d_1+\dots+d_r$.  
  
\subsubsection{Definitions, notations and examples}  
  
\begin{definition} \label{abstractGMA} Let $A$ be a commutative ring and $R$ an $A$-algebra. We  
will say that $R$ is a {\it generalized  matrix algebra (GMA)}   
of type $(d_1,\dots,d_r)$ if $R$ is  
equipped with:  
\begin{itemize}  
\item[(i)] a family orthogonal idempotents $e_1,\dots,e_r$ of sum $1$,  
\index{ei@$e_i$, idempotent of a GMA $R$ or $S$,   
part of the {\it data of idempotents} $\cal{E}$}  
\item[(ii)] for each $i$, an  
$A$-algebra isomorphism $\psi_i: e_iRe_i \longrightarrow M_{d_i}(A)$,  
\end{itemize}  
such that the {\it trace} map $T: R \longrightarrow A$, defined by  
$T(x):=\sum_{i=1}^r \tr(\psi_i(e_ixe_i))$, satisfies $T(xy)=T(yx)$ for all $x,  
y \in R$. We will call $\cal{E}=\{e_i,\psi_i, i=1,\dots, r\}$ the   
{\it data of idempotents} of $R$. \index{E@$\cal E$, the data of idempotents of the GMA $R$ or $S$}  
\end{definition}  
  
\begin{remark}  
If $R$ is a GMA as above, then $R$ equipped with the map $T(\bullet)1_R$   
is a {\it trace algebra} in the sense of   
Procesi \cite{Proc2}.  
\end{remark}  
  
\begin{notation} If $(R,\cal{E})$ is GMA as above, we shall often use the following notations.   
For $1\leq i\leq r, 1 \leq k, l \leq d_i$,   
there is a unique element $E_i^{k,l} \in e_iRe_i$ such that   
$\psi_i(E_i^{k,l})$  
is the elementary matrix of $M_{d_i}(A)$ with unique nonzero coefficient at  
 row $k$ and column $l$. These elements satisfy the usual relations    
$$E_i^{k,l}E_{i'}^{k',l'}=\delta_{i,i'}\delta_{l,k'}E_i^{k,l'},$$  
$e_i=\sum_{1\leq k \leq d_i} E_i^{k,k}$, and $AE_i^{1,1}$ is free of rank one over $A$. Clearly, the data of the $E_i^{k,l}$ satisfying these last three conditions is equivalent to condition (ii) in the definition of $R$. For each $i$, we set also $E_i:=E_i^{1,1}$.  
\end{notation}  
  
\begin{exemple} \label{gma}   
Let $A$ be a commutative ring, and $B$ \index{B@$B$, a commutative $A$-algebra} be a commutative $A$-algebra.   
Let $A_{i,j}$,   
$1\leq i, j\leq  
r$, be a family of $A$-submodules of $B$ satisfying the following properties:  
\begin{equation} \label{aij}    
{\rm For\, \,  all\, \, } i,\, j,\, k, \hspace{1cm}A_{i,i}=A, \hspace{1cm} A_{i,j} A_{j,k} \subset A_{i,k}  
\end{equation}  
Then the following $A$-submodule $R$ of $M_d(B)$  
\begin{equation} \label{forme}  
\left( \begin{array}{cccc} M_{d_1}(A_{1,1}) & M_{d_1,d_2}(A_{1,2}) & \hdots &   
M_{d_1,d_r}(A_{1,r}) \\  
M_{d_2,d_1}(A_{2,1}) & M_{d_2}(A_{2,2}) & \hdots & M_{d_2,d_r}(A_{2,r}) \\  
\vdots & \vdots & \ddots & \vdots \\  
M_{d_r,d_1}(A_{r,1}) & M_{d_r,d_2}(A_{r,2}) & \hdots & M_{d_r}(A_{r,r}) \end{array} \right)  
\end{equation}  
is an $A$-subalgebra.   
Let $e_i \in M_d(B)$ be the matrix which is the identity in the $i^{th}$ diagonal block  
(of size $d_i$) and $0$ elsewhere. As $A_{i,i}=A$, $e_i$ belongs to $R$, and  
in $R$ we have a decomposition in orthogonal idempotents   
$$1=e_1+e_2+\dots+e_r.$$  
We also have canonical isomorphisms $\psi_i:\ e_iRe_i \isomo M_{d_i}(A).$   
Hence $R$ together with $\{e_i,\psi_i,i=1,\dots,r\}$ is a GMA, and the   
trace $T$ is the restriction of the trace of $M_d(B)$. Note that assuming $d!$ invertible in $A$,    
\S\ref{mainexample} shows that $T$ is a pseudo-character of dimension $d$   
over $R$, which is Cayley-Hamilton (see \S\ref{operations}).  
  
The GMA $R$ is called the {\it standard GMA of type $(d_1,\dots,d_r)$   
associated to the $A$-submodules $A_{i,j}$ of $B$}.  
\end{exemple}  
\index{Aij@$A_{i,j}$, a structural module of a GMA $R$ or $S$,   
seen as an $A$-submodule of a commutative $A$-algebra $B$}   
\subsubsection{Structure of a GMA}  
\label{structureGMA}  
\label{structuregma}  
Let $R$ be  GMA of type $(d_1,\dots,d_r)$.   
We will attach to it a canonical family of $A$-modules   
$\cal{A}_{i,j} \subset R$, $1\leq i,j \leq r$, as follows. Set  
$$\cal{A}_{i,j}:=E_iRE_j.$$  
For each triple $1\leq i, j, k\leq r$, we have   
$$\cal{A}_{i,j}\cal{A}_{j,k} \subset \cal{A}_{i,k}$$ in $R$, hence the product in $R$   
induces a map   
$$\varphi_{i,j,k}: \cal{A}_{i,j}\otimes \cal{A}_{j,k} \longrightarrow  
\cal{A}_{i,k}.$$ Moreover, $T$ induces an $A$-linear isomorphism   
$$\cal{A}_{i,i} \isomo  
A.$$ \ps \smallskip  
\index{Aij@$\cal{A}_{i,j}$, a structural module of a GMA $R$ or $S$,   
seen as an abstract $A$-module}  
By Morita equivalence, the map induced by the   
product\footnote{All the tensor products below are assumed to be over $A$.} of $R$ $$e_iRE_i \otimes \cal{A}_{i,j} \otimes E_jRe_j \longrightarrow e_iRe_j$$   
is an isomorphism of $e_iRe_i \otimes e_jR^{\rm opp}e_j$-modules. In particular, with the help of $\psi_i$ and $\psi_j$, we get a canonical identification   
$$e_iRe_j = M_{d_i,d_j}(\cal{A}_{i,j}),$$ as a module over   
$e_iRe_i \otimes e_jR^{\rm opp}e_j=M_{d_i}(A)\otimes M_{d_j}(A)^{\rm opp}$. Moreover, in terms of these identifications, the natural map induced by the product in $R$, $e_iRe_j \otimes e_jRe_k \longrightarrow e_iRe_k$, is the map $M_{d_i,d_j}(\cal{A}_{i,j})\otimes M_{d_j,d_k}(\cal{A}_{j,k}) \longrightarrow M_{d_i,d_k}(\cal{A}_{i,k})$ induced by $\varphi_{i,j,k}$.\ps \smallskip  
  
To summarize all of this, there is a canonical isomorphism of $A$-algebra  
\begin{equation} \label{forme2}  
R \simeq \left( \begin{array}{cccc} M_{d_1}(\cal{A}_{1,1}) & M_{d_1,d_2}(\cal{A}_{1,2}) & \hdots &   
M_{d_1,d_r}(\cal{A}_{1,r}) \\  
M_{d_2,d_1}(\cal{A}_{2,1}) & M_{d_2}(\cal{A}_{2,2}) & \hdots & M_{d_2,d_r}(\cal{A}_{2,r}) \\  
\vdots & \vdots & \ddots & \vdots \\  
M_{d_r,d_1}(\cal{A}_{r,1}) & M_{d_r,d_2}(\cal{A}_{r,2}) & \hdots &   
M_{d_r}(\cal{A}_{r,r}) \end{array} \right),  
\end{equation}  
where the right hand side is a notation for the  
algebra that is $\bigoplus_{i,j} M_{d_i,d_j}(\cal{A}_{i,j})$  
as an $A$-module, and whose product is defined by the usual matrix product   
formula, using the $\varphi_{i,j,k}$'s to multiply entries. Moreover, we have canonical isomorphisms $\cal{A}_{i,i} \isomo A$. By an abuse of language, we will often write this precise isomorphism as an equality $\cal{A}_{i,i}=A$. \ps  
  
\noindent Let us consider the following sets of conditions on the   
$\varphi_{i,j,k}$'s:\ps  
(UNIT) For all $i$, $\cal{A}_{i,i}=A$ and for all $i, j$, $\varphi_{i,i,j}: A \otimes \cal{A}_{i,j}\longrightarrow \cal{A}_{i,j}$ (resp. $\varphi_{i,j,j}: \cal{A}_{i,j} \otimes A \longrightarrow \cal{A}_{i,j}$) is the $A$-module structure of $\cal{A}_{i,j}$. \ps   
(ASSO) For all $i,j,k,l$, the two natural maps $\cal{A}_{i,j}\otimes \cal{A}_{j,k}\otimes \cal{A}_{k,l} \longrightarrow \cal{A}_{i,l}$ coincide.\ps  
(COM) For all $i, j$ and for all $x \in \cal{A}_{i,j}$, $y\in \cal{A}_{j,i}$, we have
$\varphi_{i,j,i}(x\otimes y)=\varphi_{j,i,j}(y\otimes x)$. \ps  
  
\begin{lemma} \label{lesaij}  
The $\varphi_{i,j,k}$'s above satisfy the conditions (UNIT), (ASSO) and (COM). The $\varphi_{i,j,i}$'s are nondegenerate if and only if $T: R\otimes R \longrightarrow A$, $x\otimes y \mapsto T(xy)$, is nondegenerate.  
\end{lemma}  
\begin{pf} First, (ASSO) follows from the associativity of the product in $R$.  
To check (UNIT), we must show that for all $i, j$, and for all $x, y\in R$,  
then $E_ixE_iyE_j=T(E_ixE_i)E_iyE_j$ and $E_ixE_jyE_j=T(E_jyE_j)E_ixE_j$.  
As $T(R)=A$ is in the center of $R$, it suffices to check that for  
all $i$, and for all $x\in R$, $$E_ixE_i=T(E_ixE_i)E_i,$$  
but this is obvious. The property (COM) holds as $T(xy)=T(yx)$ for all   
$x, y\in R$,\ps  
Note that if $x \in R$ and $i \neq j$, $T(e_ixe_j)=T(e_je_ix)=0$. Hence for $x \in \cal{A}_{i,j}$ and $y \in \cal{A}_{i',j'}$ with $i'\neq j$ or  
$j' \neq i$, we have $T(xy)=\sum_{i=1}^r T(e_ixye_i)=0$.   
\end{pf}  
  
Reciprocally, if we have a family of $A$-modules $\cal{A}_{i,j}$, $1\leq i, j \leq r$, equipped with $A$-linear maps $\varphi_{i,j,k}: \cal{A}_{i,j} \otimes \cal{A}_{j,k} \longrightarrow \cal{A}_{i,k}$ satisfying (UNIT), (ASSO) and (COM), then we leave as an exercise to the reader to check that $R:=\oplus_{i,j} M_{d_i,d_j}(\cal{A}_{i,j})$ has a unique structure of $GMA$ of type $(d_1,\dots,d_r)$ such that for all $i, j$, $E_iRE_j=\cal{A}_{i,j}$.   
  
\subsubsection{Representations of a GMA}\label{reprGMA}   
  
\par  
  
If $R$ is an $A$-algebra, we will call {\it representation} of   
$R$ any morphism of $A$-algebras $\rho : R \longrightarrow M_n(B)$,   
where $B$ is a commutative $A$-algebra.   
If $R$ is equipped with a central function $T:R \longrightarrow A$,   
we will say that $\rho$ is a {\it trace representation} if    
$\tr \circ \rho(x)=T(x) 1_B$ for any $x \in R$. \ps  
Let $(R,\cal{E})$ be a GMA of type $(d_1,\dots, d_r)$.  
We will be interested by the trace representations of $R$, and especially  
by those that are compatible with the structure $\cal{E}$, as follows:  
  
\begin{definition} \label{adapted} Let $B$ be a commutative $A$-algebra. A representation $\rho: R \longrightarrow M_d(B)$ \index{Zrho@$\rho$, a morphism $R \longrightarrow M_d(B)$, of trace $T$ (often with $B=K$, or $B=A$)} is said to be {\it adapted to $\cal{E}$} if its restriction to the $A$-subalgebra   
$\oplus_{i=1}^r e_iRe_i$ is the composite of the representation $\oplus_{i=1}^r  
\psi_i$ by the natural "diagonal" map   
$M_{d_1}(A) \oplus \cdots \oplus M_{d_r}(A) \longrightarrow M_d(B)$.    
\end{definition}  
Obviously, an adapted representation is a trace representation. In the other   
direction we have :  
\begin{lemma} \label{traceadapted}  
 Let $B$ be a commutative $A$-algebra and   
$\rho: R \longrightarrow M_d(B)$ be a trace representation.   
Then there is a commutative ring $C$ containing $B$ and a $P\in \GL_d(C)$ such that $P \rho P^{-1} : \  R \longrightarrow M_d(C)$ is adapted to $\cal{E}$. Moreover, if every finite type projective $B$-module  
is free, then  we can take $C=B$.   
\end{lemma}  
\begin{pf} As $\tr \circ \rho = T$, the $\rho(E_i^{k,k})$'s form an orthogonal family of $d$ idempotents of trace $1$ of $M_d(B)$ whose sum is $1$. As a consequence, in the $B$-module decomposition  
$$B^d=\oplus_{i,k} \rho(E_i^{k,k})(B^d),$$  
the modules $\rho(E_i^{k,k})(B^d)$'s are projective, hence become free (of rank $1$) over a suitable ring $C$ containing $B$ (and of course we can take $C=B$ if those modules are already free). We now define a $C$-basis $f_1, \dots, f_d$ of $C^d$ as follows. For each $1\leq i \leq r$, choose first $g_i$ a $C$-basis $\rho(E_i^{1,1})(C^d)$. Then for $1\leq k\leq d_i$,  $$f_{d_1+\dots + d_{i-1}+k}:=\rho(E_i^{k,1})(g_i)$$ is a $C$-basis of $\rho(E_1^{i,i})(C^d)$. By construction, in this new $C$-basis, $\rho$ is adapted to $\cal{E}$.   
\end{pf}

Let us call $G$ the natural covariant functor from commutative $A$-algebras to sets such that for a commutative $A$-algebra $B$, $G(B)$ is the set of   
representations $\rho: R \longrightarrow M_d(B)$ adapted to $\cal{E}$. \ps  
  
Let $B$ be a commutative algebra and $\rho \in G(B)$. By a slight abuse of   
language we set $E_i:=\rho(E_i)\in M_d(B)$. By definition, for each $i, j$,  
$\rho(E_iRE_j)=E_i\rho(R)E_j$, hence it falls into the $B$-module of   
matrices whose  
coefficients are $0$ everywhere, except on line $d_1+\dots+d_{i-1}+1$ and   
row $d_1+\dots+d_{j-1}+1$. We get this way an $A$-linear map  
$f_{i,j}: \cal{A}_{i,j} \longrightarrow B$, whose image is an $A$-submodule   
of $B$ which we denote by $A_{i,j}$. Hence  
\begin{prop} \label{imageadapted} The subalgebra $\rho(R)$ of $M_d(B)$ is   
the   
standard GMA of type $(d_1,\dots,d_r)$ associated to the $A$-submodules $A_{i,j}$ of $B$ (see example~\ref{gma}).   
\end{prop}  
  
Moreover, the $f_{i,j}$'s have the two following properties: \begin{itemize}  
\item [(i)] $f_{i,i}$ is the structural map $A \longrightarrow B$,  
\item [(ii)] the product $\cdot: B \otimes B \longrightarrow B$ induces the   
$\varphi_{i,j,k}$'s, i.e. $$\forall i,j,k,\, \, \, f_{i,k} \circ \varphi_{i,j,k} = f_{i,j}\cdot f_{j,k}.$$  
\end{itemize}  
This leads us to introduce the following new functor.   
If $B$ is a commutative $A$-algebra,   
let $F(B)$ be the set of $(f_{i,j})_{1\leq i, j\leq r}$, where $f_{i,j}: \cal{A}_{i,j} \longrightarrow B$ is an $A$-linear map, satisfying conditions (i) and (ii) above. It is easy to check that $F$ is a covariant functor from commutative $A$-algebras to sets. In the discussion above, we attached to each   
$\rho \in G(B)$ an element $f_{\rho}=(f_{i,j}) \in F(B)$.  
  
\begin{prop} \label{FetG}$\rho \mapsto f_{\rho}$ induces an isomorphism of functors $G \isomo F$. Both those functors are representable by a commutative $A$-algebra $\Bun$.  
\end{prop}  
  
\begin{pf}Let $B$ be a   
commutative $A$-algebra and $f:=(f_{i,j}) \in F(B)$.   
Then $f$ induces coefficient-wise a natural map   
$$\rho_f: R=\oplus_{i,j} M_{d_i,d_j}(\cal{A}_{i,j}) \longrightarrow   
\oplus_{i,j}M_{d_i,d_j}(B)=M_d(B).$$  
This map is by definition a morphism of $A$-algebra which is adapted to $  
\cal{E}$. We get this way a morphism $F \longrightarrow G$   
which is obviously an inverse of $\rho \mapsto f_{\rho}$ constructed above.  
  
To prove the second assertion, it suffices to prove that $F$ is representable.  
If $M$ is an $A$-module, we will denote by $\Sym(M):=\oplus_{k\geq 0} \Sym^k(M)$ the symmetric $A$-algebra of $M$. We set $$\cal{B}:=\Sym(\bigoplus_{i\neq j}\cal{A}_{i,j}).$$ Let $J$ be the ideal of $\cal{B}$ generated by all the elements of the form $b\otimes c - \varphi(b\otimes c)$, where $b\in \cal{A}_{i,j}, c\in \cal{A}_{j,k}$ and $\varphi=\varphi_{i,j,k}$, for all $i,\  j$ and $k$ in $\{1,\dots,r\}$. It is obvious that  
$\Bun:=\cal{B}/I$, equipped with the canonical element $(f_{i,j}: \cal{A}_{i,j} \rightarrow \Bun)_{i,j} \in F(\Bun)$, is the universal object we are looking for.  
\end{pf}  
  
\subsubsection{An embedding problem.} \label{embedding}  
  
It is a natural question to ask when a trace algebra $(R,T)$ has   
an {\it injective} trace representation of dimension $d$,   
that is, when it can be embedded trace compatibly in a   
matrix algebra over a commutative ring.  A beautiful theorem of Procesi \cite{Proc2} gives a very satisfactory answer when $A$ is a $\Q$-algebra : $(R,T)$ has an {\it injective} trace representation of dimension $d$ if and only if $T$ satisfies the $d$-th Cayley-Hamilton identity (see~\cite{Proc2} and \S\ref{operations}). \ps  
Assume that $(R,\cal E)$ is a GMA. Then we may ask two natural questions :  
\begin{itemize}   
\item[(1)] Is there an injective $d$-dimensional trace representation of $R$ ?   
\item[(2)] Is there an injective $d$-dimensional adapted representation of $R$ ?   
\end{itemize}  
Actually, it turns out that those questions are equivalent. Indeed, if $\rho : R \longrightarrow M_d(B)$ is an injective trace representation, then Lemma~\ref{traceadapted} gives an   
injective adapted representation $R \rightarrow M_d(C)$ for some ring $C \supset B$. By elementary reasoning, question (2) is equivalent to the following questions (3) and (4).  
\begin{itemize}  
\item[(3)] Is the universal adapted representation $\rho:\ R \longrightarrow M_d(\Bun)$ injective ?  
\item[(4)] Are the universal maps $f_{i,j}: \cal{A}_{i,j} \longrightarrow \Bun$ injectives ?  
\end{itemize}  
  
For a GMA for which we know {\it a priori} that $T$ is   
a Cayley-Hamilton pseudocharacter of dimension $d$ (residually   
multiplicity free Cayley-Hamilton pseudocharacters over local henselian rings are examples of such a situation - see \S\ref{structuremult1}), Procesi's result gives a positive answer to question (1), hence to  
questions (2) to (4) as well, in the case where $A$ is a $\Q$-algebra. We shall give below a positive answer in the general case to those questions. As a consequence, by Proposition~\ref{imageadapted}, any GMA is isomorphic to some standard GMA of Example~\ref{gma}, and its trace is a Cayley-Hamilton   
pseudocharacter of dimension $d$. Note that it does not seem much easier to prove first this last fact. \ps  
  
This result (the positive answer to questions (1) to (4)) will be used in its full generality only   
in the proof of the Theorem~\ref{factoriel} (and here only for $r=2$), and  
also to prove the converse of Theorem~\ref{structure} (i)   
(see Example~\ref{importantexample}). In particular, it is not needed for the Galois theoretic applications of the following sections. However,   
we shall use several times this result in a special case    
(see \S\ref{reducedcase} below)  
where there is a much simpler proof, and where more precise results are available.   
Hence, for the commodity of the reader, we  
first give the proof in this special case.

\subsubsection{Solution of the embedding problem in the reduced and nondegenerate case}  
\label{reducedcase}  
  
Let $I=\{1,\dots,r\}$ and assume that we are given a family of $A$-modules   
$\cal{A}_{i,j}$, $i, j \in I$, and for each $i, j, k$ in $I$ an $A$-linear   
map\footnote{All the tensor products below are assumed to be over $A$.}  
$$\varphi_{i,j,k}: \cal{A}_{i,j}\otimes \cal{A}_{j,k} \longrightarrow   
\cal{A}_{i,k},$$ which  
satisfy (UNIT), (ASSO) and (COM).   
  We denote by $F$ again the functor from   
commutative $A$-algebras to sets which is associated to this data, as defined in \S\ref{reprGMA}.  
  
\begin{lemma}\label{specialcaselemma} \begin{itemize}  
\item[(i)] Assume that the $\cal{A}_{i,j}$'s are free of rank $1$ over $A$, and  
that the $\varphi_{i,j,k}$ are isomorphisms. Then there is a $(f_{i,j}) \in F(A)$ such that the $f_{i,j}$'s are  
isomorphisms.  
\item[(ii)] The relation $i \sim j$ if, and only if, $\cal{A}_{i,j}$ is free of rank one and $\varphi_{i,j,i}$ is an   
isomorphism, is an equivalence relation on $I$. Moreover, if $i \sim j \sim k$, then $\varphi_{i,j,k}$ is an isomorphism.  
\end{itemize}  
\end{lemma}  
\begin{pf}  
We show first (i). Let $e_{i,j}$ be an $A$-basis of $\cal{A}_{i,j}$. As $\varphi_{i,j,k}$ is an isomorphism, there exists a unique  
$\lambda_{i,j,k} \in A^*$ such that $\varphi_{i,j,k}(e_{i,j}\otimes  
e_{j,k})=\lambda_{i,j,k}e_{i,k}$. Let us fix some $i_0 \in I$. For all $i, j$, set  
$\mu_{i,j}:=\lambda_{i,i_0,j}$. We claim that the $A$-linear isomorphisms $f_{i,j}:= \cal{A}_{i,j}  
\longrightarrow A$ defined by $f_{i,j}(e_{i,j})=\mu_{i,j}$ satisfy  
$(f_{i,j}) \in F(A)$. It suffices to check that for all $i, j, k$, we have  
$\mu_{i,j}\mu_{j,k}=\lambda_{i,j,k}\mu_{i,k}$. But this is the hypothesis  
(ASSO) applied to $i, i_0, j$ and $k$. \par  
Let us show (ii). By (UNIT) we have $i \sim i$, and by   
(COM) $i \sim j$ implies $j \sim i$.   
If $i \sim j$ and $j \sim k$ we claim that $\varphi_{i,j,k}$   
and is an isomorphism. It will imply that $\cal{A}_{i,k}$ and $\cal{A}_{k,i}$ is free of rank $1$ over $A$, and that $\varphi_{i,k,i}$ is an isomorphism by (ASSO), hence $i \sim k$. Using (UNIT)  
and (ASSO), we check easily\footnote{If $(x,y,z,t)\in (\cal{A}_{i,j}\times\cal{A}_{j,k} \times  
\cal{A}_{k,j}\times \cal{A}_{j,i})$, using (ASSO), (ASSO) again, and (UNIT), we have with the obvious notations:  
$(xy)(zt)=x(y(zt))=x((yz)t)=(yz)(xt)$. In general, to check this kind  
of identities with values in some $\cal{A}_{k,l}$, it suffices to do it   
in the GMA of type $(1,1,\dots,1)$ defined by the $\cal{A}_{i,j}$, which  
might be a bit easier (e.g. in the proof of Proposition \ref{univ}).} the equality   
of linear maps $$\varphi_{i,k,i}\circ  
(\varphi_{i,j,k}\otimes\varphi_{k,j,i})=\varphi_{j,k,j}\cdot\varphi_{i,j,i}:\, \, \, \cal{A}_{i,j}\otimes\cal{A}_{j,k} \otimes \cal{A}_{k,j}\otimes \cal{A}_{j,i} \longrightarrow A.$$ As $i\sim j$ and $j\sim k$, it implies that $\varphi_{i,j,k}$ is injective. The surjectivity of $\varphi_{i,j,k}$ comes from the fact that the natural map $$\cal{A}_{i,k} \otimes \cal{A}_{k,j} \otimes \cal{A}_{j,k} \longrightarrow \cal{A}_{i,k}$$   
\noindent is an isomorphism (as $j \sim k$) whose image is contained in $\Im(\varphi_{i,j,k})$ by (ASSO).   
\end{pf}   

	Before stating the main proposition of this subsection, we need to recall some definitions from commutative algebra. If $A$ is a 
commutative ring, recall that the {\it total fraction ring} of $A$ is the fraction ring ${\rm Frac}(A):=S^{-1}A$ where $S \subset A$ 
is the multiplicative subset of nonzerodivisors of $A$, that is $f \in S$ if and only if the map 
$g \mapsto gf, \, \, A \rightarrow A,$ is injective. We check at once that the natural map $A \rightarrow S^{-1}A$
is injective and flat, and that each nonzerodivisor of $S^{-1}A$ is invertible. 
Of course, $S^{-1}A$ is the fraction field of $A$ if $A$ is a domain. 
\begin{prop}\label{propfinitefrac} Assume $A$ is reduced. The following properties are equivalent: \begin{itemize}
\item[(i)] $A$ has a finite number of minimal prime ideals,
\item[(ii)] $A$ embeds into a finite product of fields,
\item[(iii)] $S^{-1}A$ is a finite product of fields.
\end{itemize} If they are satisfied, $S^{-1}A = \prod_{P} A_P$ where the product is over the finite set of minimal prime ideals of $A$.
\end{prop}

\begin{pf} It is clear that (i) is equivalent to (ii). Note that
${\rm Spec}(S^{-1}A) \subset {\rm Spec}(A)$ is the subset of prime ideals that do not meet $S$. For $P$ any
minimal prime ideal of $A$, remark that the image of $f$ in $A_P={\rm Frac}(A/P)$ is not a zero divisor
of this latter ring by flatness of $A_P$ over $A$, so
$S\cap P = \emptyset$. In particular, $A$ and $S^{-1}A$ have the same minimal prime ideals, and (iii)
implies (i). Moreover, if $A$ has a finite number of minimal prime ideals, say $P_1,\dots,P_r$, then we
have an injection
$$A \longrightarrow \prod_{i=1\dots r} A_{P_i},$$ so \begin{equation}\label{finitenumbprimmin} S=A\backslash \left(P_1\cup
\cdots \cup P_r\right).\end{equation}
Assume now that (i) holds, we will show (iii) as well as the last assertion of the statement. As $A$ and
$S^{-1}A$ have the same minimal prime ideals, we may assume that $S^{-1}A=A$, {\it i.e.} that each
nonzerodivisor of $A$ is invertible. By (\ref{finitenumbprimmin}), we get that for each maximal ideal $m$
of $A$, $m \subset \cup_{i=1\dots r} P_i$. By \cite[Chap. II, \S1.1, Prop. 2]{BouAc}, 
this implies that each $P_i$ is maximal, hence $$A \isomo
\prod_{i=1\dots r} A_{P_i}$$ and we are done.
\end{pf}

An $A$-module $M$ is said to be {\it torsion free}
if the multiplication by each $f \in S$ on $M$ is injective, {\it i.e.} if $M$ factors through an $S^{-1}A$-module. 
An $A$-submodule $M$ of $S^{-1}A$ is said to be a {\it fractional ideal of $S^{-1}A$} if $fM \subset A$ for some $f\in A$ which is not a zerodivisor.
Assume that $A$ is reduced and that $S^{-1}A=\prod_s K_s$ is a finite product of fields.  
Note that if $A_s=\im(A \longrightarrow K_s)$, then $\prod_s A_s$ is a fractional ideal of $K$. 
As a consequence, $M \subset K$ is a fractional ideal if, and only if, for each $s$, $\Im(M \longrightarrow K_s)$ 
is a fractional ideal of $K_s$. We will often denote by $K$ the total fraction ring $S^{-1}A$.

\begin{prop}  
\label{specialcase}  
Assume that $A$ is reduced and that its total fraction ring $K$ is a finite product of fields. Assume moreover that the maps $\varphi_{i,j,i}: \cal{A}_{i,j} \otimes \cal{A}_{j,i} \longrightarrow A$ are nondegenerate\footnote{That is, that the induced maps $\cal{A}_{i,j} \longrightarrow \Hom_A(\cal{A}_{j,i},A)$ are injective.}.  
\ps  
Then there exists   
$(f_{i,j}) \in F(K)$ such that each $f_{i,j}: \cal{A}_{i,j} \longrightarrow  
K$ is an injection whose image is a fractional ideal of $K$.  
Moreover, if $A=K$ is a field, the relation $i\sim j$ if, and only if,  
$\cal{A}_{i,j} \neq 0$ coincides with the one of Lemma  
\ref{specialcaselemma}.   
\end{prop}  
\begin{pf}  
Write $K=\prod_s K_s$ as a finite product of fields. As $\cal{A}_{i,j}$ embeds into $\Hom_A(\cal{A}_{j,i},A)$ by assumption, it is torsion free over $A$, hence embeds into $\cal{A}_{i,j}\otimes K$. As $A \rightarrow K$ is an injection into a fraction ring, we check easily that $\varphi_{i,j,i} \otimes K$ is 
again nondegenerate\footnote{If $A$ is any commutative ring with total fraction ring $S^{-1}A$, and $M$ any $A$-module (not necessarily of finite type), 
then the natural map $S^{-1}\Hom_A(M,A) \rightarrow \Hom_{S^{-1}A}(S^{-1}M,S^{-1}A)$ is injective.}, hence so are the $\varphi_{i,j,i} \otimes K_s$'s.   
By (ASSO) applied to $i, j, i, j$, and by (COM) and (UNIT), we have:  
$$\forall x, x' \in \cal{A}_{i,j}, \forall y\in \cal{A}_{j,i}, \, \, \, \, \, \varphi_{i,j,i}(x',y)x=\varphi_{i,j,i}(x,y)x'.$$   
hence $\cal{A}_{i,j} \otimes K_s$ has $K_s$-dimension $\leq 1$ and $\cal{A}_{i,j}$ is isomorphic to a fractional ideal of $K$. It remains only to construct the injections $f_{i,j}$ of the statement. By what we have just seen, we can assume that $A=K$ is a field,   
and in this case each $\cal{A}_{i,j}$ is either $0$ or one dimensional over $K$, and the $\varphi_{i,j,i}$'s are  
nondegenerate, hence isomorphisms. \par   
For $i, j\in I$, say $i \sim j$ if $\cal{A}_{i,j} \neq 0$.   
As the $\varphi_{i,j,i}$ are isomorphisms, this relation coincides with the one defined in Lemma \ref{specialcaselemma} (ii). On each equivalence class of the relation $\sim$, we define some $f_{i,j}$'s by Lemma~\ref{specialcaselemma} (i), and we set $f_{i,j}:=0$ if $i \not \sim j$.  
\end{pf}  
  
\subsubsection{Solution of the embedding problem in the general case}\label{generalcase} Same notations as in \S\ref{reducedcase}. We recall that $\Bun$ is the universal $A$-algebra representing $F$ (see Proposition~\ref{FetG}).  
  
\begin{prop}\label{univ} The universal maps $f_{i,j}: \cal{A}_{i,j} \longrightarrow B$ are $A$-split injections.\end{prop}  
\begin{pf} We use the notations of the proof of Proposition~\ref{FetG}.  
Recall that $I=\{1,\dots,r\}$ and set   
$\Omega:=\{(i,j),\, i,\, j \in I,\, i\neq j\}$;   
if $x=(i',j') \in \Omega$ we will write $i(x):=i'$ and $j(x):=j'$.\ps  
If $\gamma=(x_1,\dots, x_s)$ is a sequence of elements of $\Omega$ such that   
for all $k\in \{1,\dots, s-1\}$   
we have $j(x_k)=i(x_{k+1})$, then we will say that $\gamma$ is a {\it path} from $i(x_1)$ to $j(x_s)$, and we will set $\cal{A}_\gamma:=\cal{A}_{i(x_1),j(x_1)}\otimes \cdots \otimes \cal{A}_{i(x_s),j(x_s)}$. If moreover $i(x_1)=j(x_s)$, we will say that $\gamma$ is a {\it cycle}. In this case, $\rot(\gamma):=(x_s,x_1,\dots,x_{s-1})$ is again a cycle. Let $i, j\in I$, $\gamma$ a path from $i$ to $j$, and $c_1,\dots,c_n$ a sequence of cycles   
(which can be empty). We will call the sequence of paths $\Gamma=(c_1,\dots,c_n,\gamma)$ an {\it extended path} from $i$ to $j$.  
If $\Gamma$ is such a sequence and $(i',j')\in \Omega$,   
we denote by $\Gamma_{i',j'}$ the total number of times   
that $(i',j')$ appears in the $c_k$'s or in $\gamma$.  
 It will be convenient to identify $\N^{\Omega}$ with the set of   
oriented graphs\footnote{In an oriented graph, we authorize   
multiple edges between two vertices $i$ and $j$, with $i \neq j$, but we   
do not authorize edges from a vertex to itself.}   
with set of vertices $I$, by associating to   
$\tau=(\tau_{i,j})_{(i,j) \in \Omega}$ the graph with $\tau_{i,j}$  
 edges from $i$ to $j$. If $\Gamma$ is an extended path from $i$ to $j$,  
 we shall say that $\tau(\Gamma):=(\Gamma_{i',j'})\in \N^{\Omega}$ is  
 the {\it underlying graph} of $\Gamma$. \par  
Let $\deg: \N^\Omega \longrightarrow \Z^I$ be the map such that, for   
$\tau \in \N^\Omega, i\in I$, $\deg(\tau)_i$ is the number of arrows in $\tau$ arriving at $i$ minus the number of arrows departing from $i$. If $(i,j) \in \Omega$, let $\tau(i,j)$ be the graph with a unique arrow, which goes from $i$ to $j$. If $i\in I$, set $\tau(i,i)=0$. The following lemma is easily checked.  
  
\begin{lemma} \label{graph}Let $i, j \in I$.   
\begin{itemize}  
\item[(i)] If $\Gamma$ is an extended path from $i$ to $j$, then $\deg(\tau(\Gamma))=\deg(\tau(i,j))$.  
\item[(ii)] If $\tau$ is a graph such that   
$\deg(\tau)=\deg(\tau(i,j))$, then $\tau=\tau(\Gamma)$ for some extended  
 path $\Gamma$ from $i$ to $j$.   
If moreover $\tau_{i',j'}\neq 0$ and $\tau_{j',k'}\neq 0$ for   
some $i', j', k' \in I$, then we can assume that the   
sequence $\Gamma$ has a path containing $((i',j'),(j',k'))$ as a subpath.  
\end{itemize}  
\end{lemma}   
   
By (ASSO), for each path $\gamma$ from $i$ to $j$, we have a canonical {\it contraction} map $\varphi_\gamma: \cal{A}_\gamma \longrightarrow \cal{A}_{i,j}$. If $\gamma$ is a cycle, $\varphi_\gamma$ goes from $\cal{A}_{\gamma}$ to $A$ by (UNIT), and the assumption (COM) implies that $\varphi_{\rot(\gamma)}=\varphi_{\gamma}\circ \rot$, where $\rot: \cal{A}_{\rot(\gamma)} \longrightarrow \cal{A}_\gamma$ is the canonical circular map. We claim now that the following property holds:\ps  
(SYM) For any cycle $c$ having some $(i',j') \in \Omega$ in common with some path $\gamma'$, the map $\varphi_{c} \otimes {\rm id}: \cal{A}_c\otimes \cal{A}_{\gamma'} \longrightarrow \cal{A}_{\gamma'}$ is symmetric in that two $\cal{A}_{i',j'}$'s.\ps  
	Indeed, by the rotation property we can assume that $c$ begins with $(i',j')$, and by (ASSO) and (UNIT) that $\gamma'=(i',j')$. By (ASSO) and (UNIT) again, we can assume then that $c=((i',j'),(j',i'))$, in which case it is an easy consequences of (ASSO) (applied with $i,j,i,j$), (UNIT) and (COM).\ps  
	Fix $i, j \in I$. Let $\Gamma=(c_1,\dots,c_n,\gamma)$ be an extended path from $i$ to $j$. We can consider the following $A$-linear map $\varphi_{\Gamma}: \cal{A}_{c_1} \otimes \cdots \otimes \cal{A}_{c_n} \otimes \cal{A}_\gamma\longrightarrow \cal{A}_{i,j}$, $$\left(\bigotimes_{k=1}^n x_k\right)\otimes y  \mapsto \left(\prod_{k=1}^n \varphi_{c_k}(x_k)\right)\varphi_\gamma(y).$$  
\newcommand{\bphi}{\overline{\varphi}}\par  
\noindent By the property (SYM), $\varphi_{\Gamma}$ factors canonically through a map  
$$\bphi_{\Gamma}: \bigotimes_{(k,l)\in \Omega} \Sym^{\Gamma_{k,l}}(\cal{A}_{k,l}) \longrightarrow \cal{A}_{i,j}.$$  
\noindent It is clear that: \begin{itemize}   
\item[(i)] for any permutation $\sigma \in \got{S}_n$, $\bphi_{(c_{\sigma(1)},\dots,c_{\sigma(n)},\gamma)}=\bphi_{(c_1,\dots,c_n,\gamma)}$,   
\item[(ii)] as the $\varphi_{c_k}$'s are invariant under rotation, $\bphi_{(\rot(c_1),\dots,c_n,\gamma)}=\bphi_{(c_1,\dots,c_n,\gamma)}$.   
\end{itemize}  
\noindent Let $\gamma=(x_1,\dots,x_s)$ be a path from $i$ to $j$ and $c=(y_1,\dots,y_{s'})$ be a cycle. We will say that $\gamma$ and $c$ are {\it linked at $i' \in I$} if there exists $x_k \in \gamma$ and $y_{k'} \in c$ with same origin, that is such that $i(x_k)=i(y_{k'})=i'$. Then can consider the path $\gamma\cup c:=(x_1,\dots,x_{k-1},y_{k'},\dots,y_{s'},y_1,\dots,y_{k'-1},x_k,\dots,x_s)$, which still goes from $i$ to $j$. Then we see that $\bphi_{\gamma\cup c}=\bphi_{c,\gamma}$, and it does not depend in particular on the $i'$ such that $\gamma$ and $c$ are linked at $i'$. As a consequence, going back to the notation of the paragraph above, we have:   
\begin{itemize}\item[(iii)] if $\gamma$ and $c_1$ (resp. $c_1$ and $c_2$) are linked, then $\bphi_{(c_1,\dots,c_n,\gamma)}=\bphi_{(c_2,\dots,c_n,\gamma\cup c_1)}$ (resp. $\bphi_{(c_1,c_2,\dots,c_n,\gamma)}=\bphi_{(c_1 \cup c_2,c_3, \dots,c_n,\gamma )}$).   
\end{itemize}  
	Let now $\Gamma'$ be another extended path from $i$ to $j$. Then using several times the "moves" (i), (ii) and (iii), we check at once that $\bphi_{\Gamma}=\bphi_{\Gamma'}$. Let $\tau \in \N^\Omega$ satisfies $\deg(\tau)=\deg(\tau(i,j))$. By Lemma \ref{graph} (ii), we can choose an extended path $\Gamma$ from $i$ to $j$ with   
underlying graph $\tau$, and define  
$$\bphi_{\tau}:=\bphi_{\Gamma}, \, \, \, \,  \bigotimes_{(k,l)\in \Omega} \Sym^{\tau_{k,l}}(\cal{A}_{k,l}) \longrightarrow \cal{A}_{i,j},$$  
which does not depend on $\Gamma$ (whose associated graph is $\tau$) by what we said above.\ps  
	Let us finish the proof of the proposition. The $A$-algebras $\cal{B}$ and is naturally graded by the additive monoid $\N^{\Omega}$. We have $\cal{B}=\oplus_{\tau \in \N^\Omega} \cal{B}_\tau$, where   
$\cal{B}_{\tau}=\bigotimes_{i\neq j} \Sym^{\tau_{i,j}}(\cal{A}_{i,j})$.  
The map $\deg: \cal{G}\longrightarrow \Z^I$ is additive, hence we get a   
$\Z^I$-graduation\footnote{Actually, it is even graded by the subgroup of $\Z^I$ whose elements $(n_i)$ satisfy $\sum_i n_i=0$.} on $\cal{B}$.  
Obviously, if $n\in \Z^I$, then $\cal{B}_n=\bigoplus_{\tau\in \N^\Omega, \deg(\tau)=n} \cal{B}_\tau$. For this latter graduation, the ideal $J \subset \cal{B}$ is homogeneous, hence $\Bun$ is also graded by $\Z^I$. \ps  
	Fix now $i, j \in I$, and let $n:=\deg(\tau(i,j))$. If $\deg(\tau)=n$, we constructed above a map $\bphi_{\tau}: \cal{B}_{\tau} \longrightarrow \cal{A}_{i,j}$. By summing all of them we get an $A$-linear map:  
$$\bphi_n: \cal{B}_n \longrightarrow \cal{A}_{i,j}.$$  
We claim that $\bphi_n(I_n)=0$. Assuming that, $\bphi_n$ factors through a map $$\psi_n: (\Bun)_n \longrightarrow \cal{A}_{i,j}.$$ Let $f_{i,j}: \cal{A}_{i,j} \longrightarrow (\Bun)_n$ denote the canonical map. Then by construction, $\psi_n\circ f_{i,j}=\bphi_{\tau(i,j)}$ is the identity map. It concludes the proof. \ps  
	Let us check the claim. Let $b \in \cal{A}_{i',j'}$, $c \in A_{j',k'}$ and $\varphi=\varphi_{i',j',k'}$, for some $(i',j'),\, (j',k') \in \Omega$. By $A$-linearity, is suffices to show that $\bphi_n$ vanishes on the elements of the form $x=f\otimes(b\otimes c - \varphi(b \otimes c))$, where $f$ is in $\cal{B}_{\tau}$ for some graph $\tau$ satisfying $\deg(\tau+\tau(i',k'))=n$. By Lemma \ref{graph} (ii), we can find an extended path $\Gamma$ from $i$ to $j$ with underlying graph $\tau+\tau(i',j')+\tau(j',k')$, such that some path $\gamma'$ of $\Gamma$ contains $((i',j'),(j',k'))$ as a subpath.   
Let $\Gamma'$ be the extended path from $i$ to $j$ obtained from $\Gamma$ by replacing $\gamma'=(\cdots,(i',j'),(j',k'),\cdots)$ by $(\cdots,(i',k'),\cdots)$. By construction, $\bphi_{\Gamma}(f\otimes c\otimes b)=\bphi_{\Gamma'}(f \otimes \varphi(b \otimes c))$, hence $\bphi(x)=0$. \end{pf}  
  
\begin{remark}\label{remarkuniv} When $r=2$, a slight modification of the above proof shows that the $A$-linear map $A\oplus \bigoplus_{n\geq 1}\left(\Sym^n(\cal{A}_{1,2})\oplus\Sym^n(\cal{A}_{2,1})\right) \longrightarrow \Bun$, induced by $f_{1,2}$ and $f_{2,1}$, is an isomorphism. This describes $\Bun$ completely in this case.  
%\item[(ii)] More generally, if we add in the axioms that for all $i, j, k$, the two natural maps $\cal{A}_{i,j} \otimes \cal{A}_{i,j} \otimes \cal{A}_{j,k} \longrightarrow A_{i,k}$ coincide\footnote{This condition holds for instance if the all the $\varphi_{i,j,i}$'s are nondegenerate, which will be satisfied in many interesting cases: see Lemma \ref{lesaij}.}, then we can modify the proof and show that for all $i, j$, $\Sym (f_{i,j}): \Sym(\cal{A}_{i,j}) \longrightarrow B$ is a split injection. If $|I| \leq 3$, it implies that $B=A\oplus \bigoplus_{i \neq j, n\geq 1} \Sym^n(\cal{A}_{i,j})$.  
%\item[(iii)] (Change of base ring) Let $A'$ be a commutative $A$-algebra, and set $\cal{A}_{i,j}':=\cal{A}_{i,j}\otimes_A A'$ and  
%$\varphi_{i,j,k}':=\varphi_{i,j,k}\otimes_A A'$. Then the new universal problem is represented by $B':=B \otimes_A A'$ and the $f'_{i,j}:=f_{i,j} \otimes_A A'$.  
\end{remark}  
  
As we have noted in \S\ref{embedding}, we have :  
  
\begin{cor} \label{coruniv} If $(R,\cal{E})$ is a GMA of type $(d_1,\dots,d_r)$, and if $d!$ is invertible in $A$ (where $d=d_1+\dots+d_r$), then the trace $T$ of $R$ is a Cayley-Hamilton pseudocharacter of dimension $d$.  
\end{cor}

\subsection{Residually multiplicity-free pseudocharacters}\label{structuremult1}  
  
\subsubsection{Definition}\label{hypomfree} In all this section, $A$\index{A@$A$, a local henselian ring|(} is a local henselian ring (see \cite{Ray}), $m$   
\index{m@$m$, the maximal ideal of the local ring $A$}   
is the maximal ideal of $A$, and $k:=A/m$\index{k@$k$, the residue field of $A$}. Let $R$ be an $A$-algebra and   
$T :\ R \longrightarrow A$ be a pseudocharacter of dimension $d$\index{d@$d$, the dimension of $T$}.   
Let $\bar R := R \otimes_A k$ and   
$\bar T :=T \otimes_A k: \bar R \longrightarrow k$\index{Ta@$\bar T$, the reduction of $T$ mod $m$} be the   
reductions mod $\m$ of $R$ and $T$.  
  
\begin{definition} \label{defmult1}We say that $T$ is residually multiplicity   
free if there are representations   
$\rhob_i: R \longrightarrow M_{d_i}(k)$, $i=1,\dots, r$,   
which are absolutely irreducible \index{ZRhoib@$\rhob_i$,   
a representation of $R$ over $k$ of dimension $d_i$,   
component of $\bar T$} and pairwise nonisomorphic, such that   
$\bar T=\sum_{i=1}^r \tr \rhob_i$.  
\end{definition}  
\index{di@$d_i$, the dimension of $\rhob_i$} \index{r@$r$, the number   
of factors of $\Tb$, i.e. the number of $\rhob_i$'s}  
We set $d_i:= \dim \rhob_i$, we have $\sum_{i=1}^r d_i = d$. \ps  
   
\begin{example} \label{importantexample}  
Let us give an important example. Let $(R,\cal{E})$ be a GMA   
(\S\ref{abstractGMA}), then its trace $T: R \longrightarrow A$ is a   
Cayley-Hamilton pseudocharacter by Corollary~\ref{coruniv}. We use the notations of \S\ref{structuregma}. Assume moreover that for all $i \neq j$,   
we have $$T(\cal{A}_{i,j}\cal{A}_{j,i}) \subset m.$$  
Now, for each $i$, let $\rhob_i: R \longrightarrow M_{d_s}(k)$, $r \mapsto (\psi_i(e_ire_i) \bmod m)$. Then we see easily   
that the $\rhob_i$'s are pairwise non isomorphic surjective   
representations\footnote{Note that the maps $f_{i,j}: \cal{A}_{i,j} \longrightarrow k$, defined to be $0$ if $i\neq j$, and $A \overset{\rm can}{\longrightarrow} k$ if $i=j$, define an element of $F(k)$.}, and that $T=\sum_{s=1}^r \tr \rhob_i$, hence $T$ is residually multiplicity free. The main result of this section will show that this example is the general case.   
\end{example}  
  
\subsubsection{Lifting idempotents.}  \label{lifting} Let $A$, $R$ and $T$ be as in \S\ref{hypomfree}, and assume that $T$ is residually multiplicity free. In particular, we have some representations $\rhob_i: R \longrightarrow M_{d_i}(k)$ as in definition \ref{defmult1}.   
\begin{lemma} \label{idempotents} Suppose $T$ Cayley-Hamilton.   
There are orthogonal idempotents $e_1,\dots,e_r$ in   
$R$ such   
that   
\begin{itemize}  
\item[(1)] $\sum_{i=1}^r e_i=1$.  
\item[(2)] For each $i$, $T(e_i)=d_i$   
\item[(3)] For all $x \in R$, we have $T(e_i x e_i) \equiv \tr \rhob_i(x) \pmod{\m}$  
\item[(4)] If $i \neq j$, $T(e_ixe_jye_i) \in m$ for any $x,y \in R$.  
\item[(5)] There is an $A$-algebra isomorphism $\psi_i: e_iRe_i \longrightarrow M_{d_i}(A)$ lifting $(\rhob_i)_{|e_iRe_i}: e_iRe_i \longrightarrow M_{d_i}(k)$, and such that for all $x \in e_iRe_i$, $T(x)=\tr(\psi_i(x))$.  
\end{itemize}   
Moreover, if $e'_1,\dots,e'_r$ is another family of orthogonal idempotents of $R$ satisfying (3), then there exists $x \in 1+\rad(R)$ such that for all $i$, $e'_i=xe_ix^{-1}$.   
\end{lemma}  
  
\begin{pf} Let $\rhob : \bar R \longrightarrow M_d(k)$ be the product of the $\rhob_i$'s. Because the $\rhob_i$'s are pairwise distinct,  
the image of $\rhob$ is $\prod_{i=1}^r M_{d_i}(k)$.   
As $\rhob$ is semisimple, \cite[thm. 1.1]{Tay}  
implies that $\ker \rhob = \ker \bar T$. We have the following diagram  
  
\begin{equation} \label{isom}  
\xymatrix{ \bar R / \Ker \bar T \ar[rr]^{\rhobar=\prod_{i=1}^r \rhob_i}   
\ar[rd]_{\bar T} &  & \prod_{i=1}^r M_{d_i}(k) \ar[dl]^{tr} \\  
&k& }  
\end{equation}  
which commutes by assumption on $T$, and whose first row is an isomorphism. Let us call $\epsilon_i$, for $l=1,\dots,r$, the central idempotents of $\bar R / \Ker \bar T$ corresponding to the unit of $M_{d_i}(k)$ in this decomposition.  
  
By the Cayley-Hamilton identity, and following \cite[chap. III,\S~4, exercice 5(b)]{BouAc}\footnote{The statement is that if $A$ is an henselian local ring, $R$ an $A$-algebra which is integral over $A$, and $I$ a two-sided ideal of $R$, then any family of orthogonal idempotents of $R/I$ lifts to $R$. Note that is stated there with $R$ a finite $A$-algebra, but the same proof holds in the integral case.}, there exists a family of orthogonal idempotents $e_i \in R$,   
$i=1,\dots,r $, with $e_i$ lifting the $\epsilon_i$. The element   
$1-\sum_{i=1}^r e_i$ is an idempotent which is in the radical of $R$   
by Lemma~\ref{radical}, hence it is $0$, which proves (1). By   
Lemma~\ref{dimpseudo}(1) $T(e_i)$ is an integer less than $d$, and because   
$\overline {T(e_i)} = \bar T (\epsilon_i)=d_i$, we have $T(e_i)=d_i$,  
which is (2). \par  
The assertion (3) follows from the diagram~(\ref{isom}). In order to prove (5) it suffices to show that the image of $e_ixe_jye_i$   
is zero in $\bar R/\ker \bar T$. But this image is $\epsilon_i \bar x \epsilon_j \bar y \epsilon_i$ which is zero   
by the diagram~(\ref{isom}), and we are done.\par  
Now consider the restriction $T_i$ of $T$ to the subalgebra $e_i R e_i$  
(with unit  element  $e_i$) of $R$. By Lemma~\ref{dimpseudo}(3), $T_i$ is  
a pseudocharacter of dimension $d_i=T(e_i)$, faithful if $T$ is. By (3), $T_i$ is moreover residually absolutely irreducible. If we had assumed $T$ faithful, we could have applied \cite[thm. 5.1 or   
cor. 5.2]{Rou} and get (5). As we assume only $T$ Cayley-Hamilton, we have to argue a bit more. By Lemma \ref{dimpseudo} (4), $T_i$ is Cayley-Hamilton, hence we may assume that $r=1$, and we have to prove that $R=M_d(A)$. By Lemma \ref{radical} and \cite[chap. III,\S~4, exercice 5(c)]{BouAc}, we can lift the basic matrices of $R/\Ker(\rhob_1)=M_{d_1}(k)$, i.e. find elements $(E^{k,l})_{1\leq k, l \leq d}$ in $R$ satisfying the relations $E^{k,l}E^{k',l'}=\delta_{l,k'} E^{k,l'}$. By Lemma \ref{dimpseudo} (1), for each $k \in \{1,\dots,d\}$ we have $T(E^{k,k})=1$. By Lemma \ref{dimpseudo} (4), $T_k: E^{k,k}R E^{k,k}\longrightarrow A$ is Cayley-Hamilton of dimension $1$, hence $T_k$ is an isomorphism and $E^{k,k}RE^{k,k}=A E^{k,k}$ is free of rank $1$ over $A$. Now, if $x \in E^{k,k}R E^{l,l}$, then $$x=E^{k,l}(E^{l,k}x)=E^{k,l}(T(E^{l,k}x)E^{l,l})=T(E^{l,k}x)E^{k,l} \in A E^{k,l},$$  
hence $R=\sum_{k,l} A E^{k,l}$. This concludes the proof of (5) (we even showed that Rouquier's   
Theorem 5.1 holds when faithful is replaced by Cayley-Hamilton). \ps  
To prove the last assertion, note first that the hypothesis on the $e'_i$ means that $\overline{e'_i}=\varepsilon_i$, hence by the work above the properties (1) to (5) hold also for the $e'_i$'s. As $e_iRe_i\simeq M_{d_i}(A)$ is a local ring, the Krull-Schmidt-Azumaya  Theorem \cite[thm. (6.12)]{CR} (see the remark there, \cite[prop. 6.6]{CR} and \cite[chap. 6, exercise 14]{CR}), there exists an $x\in R^*$ such that for each $i$, $xe_ix^{-1}=e'_i$. Up to conjugation by an element in $\sum_i (e_iRe_i)^*$, we may assume that $x\in 1+\rad(R)$.\end{pf}  
\subsubsection{The structure theorem} Let $A$, $R$, $T$ be as in \S\ref{lifting}.   
  
\begin{theorem}\label{structure}   
\begin{itemize}  
\item[(i)] Let $S$ be a Cayley-Hamilton quotient of $(R,T)$.  
  
Then there is a data $\cal{E}=\{e_i,\psi_i, 1\leq i \leq r\}$ on $S$   
for which $S$ is a GMA and such that for each $i$,   
$\psi_i \otimes k=(\rhob_i)_{|e_iSe_i}$. Two such data on $S$ are conjugate under $S^*$. Every such data defines  
$A$-submodules $\cal{A}_{i,j}$ of $S$ that satisfy   
$$\cal{A}_{i,j} \cal{A}_{j,k} \subset \cal{A}_{i,k},\ \ T: \cal{A}_{i,i}\isomo A,\ \   
T(\cal{A}_{i,j}\cal{A}_{j,i}) \subset m$$   
and   
$$S \simeq \left( \begin{array}{cccc} M_{d_1}(\cal{A}_{1,1}) & M_{d_1,d_2}(\cal{A}_{1,2}) & \hdots &   
M_{d_1,d_r}(\cal{A}_{1,r}) \\  
M_{d_2,d_1}(\cal{A}_{2,1}) & M_{d_2}(\cal{A}_{2,2}) & \hdots & M_{d_2,d_r}(\cal{A}_{2,r}) \\  
\vdots & \vdots & \ddots & \vdots \\  
M_{d_r,d_1}(\cal{A}_{r,1}) & M_{d_r,d_2}(\cal{A}_{r,2}) & \hdots &   
M_{d_r}(\cal{A}_{r,r}) \end{array} \right)$$  
  
\item[(ii)]  \index{Zrho@$\rho$, a morphism $R \longrightarrow M_d(B)$, of trace $T$ (often with $B=K$, or $B=A$)}  
Assume that $A$ is reduced, and that its total fraction ring $K$ is a finite product of fields. \index{K@$K$, when $A$ is reduced, its total fraction ring}  
Take $S=R/\ker T$. Choose a data $\cal{E}$ on $S$ as in (i). Then there exists an adapted injective representation  
$\rho:  S \longrightarrow M_d(K)$ whose image has the form   
$$\left( \begin{array}{cccc} M_{d_1}({A}_{1,1}) & M_{d_1,d_2}({A}_{1,2}) & \hdots &   
M_{d_1,d_r}({A}_{1,r}) \\  
M_{d_2,d_1}({A}_{2,1}) & M_{d_2}({A}_{2,2}) & \hdots & M_{d_2,d_r}({A}_{2,r}) \\  
\vdots & \vdots & \ddots & \vdots \\  
M_{d_r,d_1}({A}_{r,1}) & M_{d_r,d_2}({A}_{r,2}) & \hdots &   
M_{d_r}({A}_{r,r}) \end{array} \right)$$ where the $A_{i,j}$ are fractional  
ideals of $K$ that satisfy  
 $${A}_{i,j} A_{j,k} \subset {A}_{i,k},\ \ {A}_{i,i}=A,\ \   
{A}_{i,j}{A}_{j,i} \subset m.$$   
Moreover the $A_{i,j}$'s are isomorphic to the $\cal{A}_{i,j}$'s of part (i), in such a way that the map   
$A_{i,j} \otimes_A A_{j,k} \longrightarrow A_{i,k}$ given by the product in $K$ and the map   
$\cal{A}_{i,j} \otimes_A \cal{A}_{j,k} \longrightarrow \cal{A}_{i,k}$ given by the product in $R$ coincide.   
\item[(iii)] Let $P\in {\rm Spec}(A)$, $L:={\rm Frac}(A/P)$, and assume that   
$T\otimes L$ is irreducible\footnote{This means that  
$T\otimes L$ is not the sum of two $L$-valued pseudocharacters on  
$S\otimes L$.}. If $S$ is any Cayley-Hamilton  
quotient of $(R,T)$, then $S \otimes L$ is trace isomorphic to $M_d(L)$.  
In particular, $T \otimes L$ is faithful and absolutely irreducible.  
\end{itemize}  
\end{theorem}  
\begin{pf} As $S$ is Cayley-Hamilton, Lemma \ref{idempotents} gives us a data $\cal{E}=\{e_i,\psi_i,1\leq i \leq r\}$ satisfying (i).\par  
Assume now moreover that $A$ is as in (ii), and set $S:=R/\ker T$. Since $T$ is faithful on $S$, Lemma~\ref{lesaij}  
proves that the $\varphi_{i,j,i}$'s are nondegenerate. Then Proposition~\ref{specialcase} gives us a family of   
injections $f_{ij} : \cal{A}_{i,j} \longrightarrow L$, $(f_{i,j}) \in F(L)$ whose image are fractional ideals. Set $A_{i,j}:=f_{i,j}(\cal{A}_{i,j})$.   
By Proposition~\ref{FetG}, $(f_{i,j})$ defines an adapted representation $\rho :S \longrightarrow M_d(L)$   
that satisfies (ii).\par  
Let us prove (iii). Note that $A/P$ is still local henselian and that $S\otimes A/P$  
is a Cayley-Hamilton quotient of $(R \otimes A/P, T\otimes A/P)$, hence we may assume that  
$A$ is a domain and that $P=0$. In this case, we check at once that the natural map   
$(\Ker T)\otimes L \longrightarrow \Ker (T\otimes L)$ is an isomorphism. By  
this and by (i) applied to $T: S/\Ker T \longrightarrow A$, we see that   
$S':=(S\otimes L)/(\Ker T\otimes L)$ is a GMA of type $(d_1,\dots,d_r)$ over $L$  
whose trace $T \otimes L$ is faithful. As $T \otimes L$   
is irreducible by assumption, Proposition \ref{specialcase} implies that  
$S'$ is trace isomorphic to $M_d(L)$, as the equivalence relation there  may only have one class.   
Let us consider now the surjective map $$\psi: S \otimes L \longrightarrow  
(S\otimes L)/(\Ker T\otimes L)\isomo M_d(L).$$   
By Lemma \ref{kerrad}, its kernel is in $\rad(S\otimes L)$. By an argument  
already given in part (5) of Lemma \ref{idempotents} (using the lifting of  
the $E^{k,l}$'s of $M_d(L)$ to $S\otimes L$, and checking that they span $S  
\otimes L$ by Lemma \ref{dimpseudo} (1) and (4)), $\psi$ is an  
isomorphism, which concludes the proof.  
\end{pf}  
\begin{remark} \label{remthmstru}If $A$ is reduced and noetherian, it satisfies the conditions of (ii), hence the $A_{i,j}$'s and $R/\ker T$ are finite type torsion free $A$-modules.  
\end{remark}

\subsection{Reducibility loci and $\Ext$-groups}\label{redlocext}  
  
\subsubsection{Reducibility loci} \label{assredloc}Let $A$ be an henselian local ring, $R$ an $A$-algebra and $T: R \longrightarrow A$ a   residually multiplicity free pseudocharacter of dimension $d$. We shall use the notations of \S\ref{defmult1}.   
  
\begin{prop} \label{redloc}Let $\PP=(\PP_1,\dots,\PP_s)$ be a partition of   
$\{1,\dots,r\}$. \index{P@$\PP$, a partition of $\{1,\dots,r\}$} There exists an ideal $I_\PP$ of $A$  
such that for each ideal $J$ of $A$, the following property   
holds if and only if $I_\PP \subset J$ :\ps  
\index{IP@$I_\PP$, the reducibility ideal attached to the partition $\PP$}  
\noindent $(\text{dec}_{\PP})$ There exists pseudocharacters   
$T_1, \dots, T_s : R \otimes A/J    
\longrightarrow A/J$ such that   
\begin{itemize}  
\item[(i)] $T \otimes A/J =   
\sum_{l=1}^s T_l$,  
\item[(ii)] for each $l \in \{1,\dots,s\}$,  
 $T_l \otimes k = \sum_{i \in \PP_l} \tr\rhob_i.$  
\end{itemize}\ps  
  
\noindent If this property holds, then the $T_l$'s are uniquely determined and satisfy $\Ker T_l \subset \Ker (T\otimes A/J)$.   
  
Moreover, if $S$ is a Cayley-Hamilton quotient of $(R,T)$  
then, using the notations of Theorem~\ref{structure}, we have (for any   
choice of the data $\cal{E}$ on $S$)  
$$I_\PP = \sum_{\stackrel{(i,j)}{i,j \text{ are not in the same } \PP_l }} T(\cal{A}_{i,j}\cal{A}_{j,i})$$  
\end{prop}  
  
\begin{pf}   
	Let $S$ be a Cayley-Hamilton quotient of $(R,T)$.   
We can then chose a GMA data $\cal{E}$ for $S$ as in Theorem \ref{structure} (i),   
and consider the structural modules $\cal{A}_{i,j}=E_iSE_j$. We set  
$$I_{\PP}(T,S,\cal{E}):=\sum_{i,j \text{ are not in the same } \PP_l} T(\cal{A}_{i,j}\cal{A}_{j,i}).$$   
By Theorem \ref{structure} (i), $I_{\PP}(T,S,\cal{E})$ does not depend on the choice of the data   
$\cal{E}$ used to define it. We claim that it does not depend on $S$. Indeed, we check at once that the image of $\cal{E}$ under the surjective homomorphism $\psi: S \longrightarrow R/\Ker T$ is a data of idempotents for $R/\Ker T$ (and even that $\psi$ is an isomorphism on $\oplus_{i=1}^r e_iSe_i$). As $T \circ \psi = T$, we have that  $$T(\psi(\cal{A}_{i,j})\psi(\cal{A}_{j,i}))=T(\cal{A}_{i,j}\cal{A}_{j,i}),$$  
which proves the claim. We can now set without ambiguity $I_{\PP}:=I_{\PP}(T)$.   
As a first consequence of all of this, we see that if $J\subset A$ is an ideal, then $I_{\PP}(T\otimes A/J)$ is the   
image in $A/J$ of $I_{\PP}(T)$.   
\ps  
  
To prove the proposition we are reduced to show the following statement:   
\begin{flushleft} {\it $T: R \longrightarrow A$ satisfies $I_{\PP}=0$ if and only if we can write $T=T_1+\dots+T_s$ as a sum of pseudocharacters satisfying assumption (ii) in $(\text{dec}_{\PP})$.}  
\end{flushleft}  \ps  
Let us prove first the "only if" part of the statement above.   
Let $S=R/\Ker T$ and fix a GMA data $\cal{E}$ as in Theorem \ref{structure}  
(i). Set   
\begin{equation}\label{lesfl} f_l := \sum_{i \in \PP_l} e_i \in S,\end{equation} then $1=f_1+\dots+f_s$ is a decomposition in orthogonal idempotents. In this setting, the condition $I_{\PP}=0$ means that for each $l$,  
\begin{equation}\label{meaningIP} T(f_lS(1-f_l)Sf_l)=0. \end{equation}   
As a consequence, the two-sided ideal $f_lS(1-f_l)Sf_l$ of the ring  
$f_lSf_l$ is included in the kernel of the pseudocharacter  
$T_{f_l}=T_{|f_lRf_l}: f_lRf_l \longrightarrow A$ (see Lemma \ref{dimpseudo}  
(3)). The map $T_l: R \longrightarrow A$ defined by $T_l(x):=T(f_lxf_l)$ is then the composite of the  
$A$-algebra homomorphism   
\begin{equation}\label{formuletl}S \longrightarrow  
f_lSf_l/(f_lS(1-f_l)Sf_l), \, \, \, x \mapsto f_lxf_l+f_lS(1-f_l)Sf_l,  
\end{equation}  
by $T_{f_l}$, hence it is a pseudocharacter. As $1=f_1+\dots+f_s$, we have $T=T_1+\dots+T_s$, and the $T_l$'s satisfy (ii) of $(\text{dec}_{\PP})$ by Lemma \ref{idempotents} (3), hence we are done. In particular, we have shown that $I_{\PP}$ always satisfies $(\text{dec}_{\PP})$.\ps  
Let us prove now the "if" part of the statement. Let $K=\bigcap_i \Ker T_i$, by assumption $K \subset \Ker T$. By \S\ref{operations}, $T: R/K \longrightarrow A$ is Cayley-Hamilton, hence we can choose a data $\cal{E}$ for $S:=R/K$ and consider again the $f_l \in S$'s defined from the $e_i$'s as in formula (\ref{lesfl}) above. To check that $I_{\PP}=0$, it suffices to check that $I_{\PP}(T,S,\cal{E})=0$ or, which is the same, that $T(f_lSf_{l'}Sf_l)=0$ for $l \neq l'$. As $T=T_1+\dots+T_s$, it suffices to show that for all $x \in S$, $T_l(f_{l'}x)=0$ if $l \neq l'$. But if $l \neq l'$, $T_l(f_l')$ is in the maximal ideal $m$ by assumption (ii) of $(\text{dec}_{\PP})$ and Lemma \ref{idempotents} (3). By Lemma \ref{dimpseudo} (1), it implies that $T_l(f_{l'})=0$. By Lemma \ref{dimpseudo} (5), we get $f_{l'} \in \Ker T_l$, what we wanted.  
  
 In particular, we proved that for all $x \in S$, $T_l(x)=T(f_lx)$. As a   
consequence, $\Ker T_l  \subset \Ker T$, $K=\Ker T$, $S=R/\Ker T$, and the   
$T_l$'s are  
unique. \end{pf}   
  
\begin{definition} \label{defredloc}  
We call $I_\PP$ the {\it reducibility ideal of $T$ for the partition $\PP$}. We call the closed   
subscheme $\spec(A/I_\PP)$ of $\spec A$ the  {\it reducibility locus of $T$ for the partition $\PP$}.  
When $\PP$ is the total partition $\{\{1\},\{2\},\dots,\{r\}\}$, we call $I_\PP$ the {\it total} reducibility ideal and $\spec(A/I_\PP)$ the total reducibility locus of $T$.  
\end{definition}  
  
Note that $I_{\PP} \subset I_{\PP'}$ if $\PP'$ is a finer partition than $\PP$.

\subsubsection{The representation $\rho_i$.} We keep the assumptions of \S 
\ref{assredloc}, and we assume now that $\{i\} \in \PP$. Then for each ideal $J$ containing   
$I_\PP$, there is by Proposition \ref{redloc} a unique   
pseudocharacter $T_i : R \otimes A/J \longrightarrow A/J$ with   
$T_i \otimes k = \tr \rhob_i$ and $T = T_i + T'$ with   
$T' \otimes k = \sum_{j \neq i} \tr \rhob_i$. If   
$J \subset J'$, the pseudocharacter  
 $T_i : R\otimes A/J' \longrightarrow J'$ is just $T_i \otimes_{R/J} R/J'$, hence it is not dangerous   
to forget the ideal $J$ in the notation. As $\rhob_i$ is irreducible, we know that there is a  
(surjective, unique up to conjugation) representation $\rho_i : R/JR \longrightarrow M_{d_i}(A/J)$ of trace $T_i$  
which reduces to $\rhob_i$ modulo $m$.   
  
\begin{definition} \label{lesrhoi} If $\{i\} \in \PP$ and $J \supset I_\PP$, we let $\rho_i: R/JR \longrightarrow M_{d_i}(A/J)$ be the surjective representation defined above.  
\end{definition}   
  
As usual, by a slight abuse of notation, we will denote also by $\rho_i$ the $R$-module $(A/J)^{d_i}$ on  
which $R$ acts via $\rho_i$. It will be useful for the next section to collect here the following  
facts which are easy consequences of the proof of Proposition \ref{redloc}:  
  
\begin{lemma} \label{conseq}Let $S$ be a Cayley-Hamilton quotient of $(R,T)$,   
$\PP$ a partition of $\{1,\dots,r\}$ such that $\{i\} \in \PP$ and $J\supset I_\PP$.   
\begin{itemize}  
\item[(i)] If $j\neq i$, $e_i(S/J)e_j(S/J)e_i=0$.   
\item[(ii)] The canonical projection $$a_{i,i}: S/JS \longrightarrow  
e_i(S/JS)e_i \simeq M_{d_i}(A/J), \, \, \, x \mapsto e_ixe_i,$$ is an $A/J$-algebra  
homomorphism and satisfies $T \circ a_{i,i} =  
T_i$. As a consequence, $\rho_i$ factors through $S/JS$, $a_{i,i} \simeq  
\rho_i$, and $\rho_i(e_k)=\delta_{i,k}{\rm id}$.  
\item[(iii)] Assume moreover that $\{j\}\in \PP$ for some $j\neq i$, then we  
have $$a_{i,j}(xy)-(a_{i,i}(x)a_{i,j}(y)+a_{i,j}(x)a_{j,j}(y)) \in \sum_{k\neq i,  
j} e_i(S/J)e_k(S/J)e_j, \, \, \, \forall x, y \in R,$$  
where $a_{i,j}: S/JS \longrightarrow e_i(S/JS)e_j, \, \, \, x \mapsto e_ixe_j,$ is the canonical projection.  
\end{itemize}  
\end{lemma}  
  
\index{Zrhoi@$\rho_i$, a representation of $R$ over a quotient of $A$,   
lifting $\rhob_i$}  
  
\begin{pf} The idempotent $f_l$ corresponding to $\{i\}$  
is then $e_i$. Note that $e_i(S/JS)(1-e_i)(S/JS)e_i$ is a two-sided ideal of $e_i(S/JS)e_i  
\simeq M_{d_i}(A/J)$ whose trace is $0$ by assumption and formula (\ref{meaningIP}),  
which shows (i). As a consequence,  
$a_{i,i}$ coincides with the map in formula (\ref{formuletl}) (with of course $S$ replaced by  
$S/JS$), which proves (ii). The last assertion is immediate from the fact that $e_ixye_j-(e_ix(e_i+e_j)ye_j)$ lies in $$ e_i(S/JS)(1-(e_i+e_j))(S/JS)e_j=\sum_{k  
\neq i, j} e_i(S/JS)e_k(S/JS)e_j. $$  
\end{pf}  
  
\subsubsection{An explicit construction of extensions betweens the $\rho_i$'s.} \label{extensions}  
  
We keep the assumptions of \S\ref{assredloc}, and {\it we fix a Cayley-Hamilton  
quotient $S$ of $(R,T)$.} We fix a data $\cal{E}$ on $S$, using Theorem \ref{structure} (i), such that $(S,\cal{E})$ is a GMA  
 and set $$\cal{A}'_{i,j}= \sum_{k \neq i,j} \cal{A}_{i,k}\cal{A}_{k,j}.$$ We have  
by definition $\cal{A}'_{i,j} \subset \cal{A}_{i,j}$.\ps  
  
\index{Aijprime@$\cal{A}'_{i,j}$, the submodule   
$\sum_{k \neq i,j} \cal{A}_{i,k} \cal{A}_{k,j}$ of $\cal{A}_{i,j}$}  
  
Fix $i\neq j \in \{1,\dots,r\}$. Let $\PP$   
be any partition of $\{1,\dots,r\}$ such that the   
singletons $\{i\}$ and $\{j\}$ belong to $\PP$, and $J$ an ideal containing  
$I_\PP$. By Definition \ref{lesrhoi}, for $k=i, j$, we have a representation   
$\rho_k : R/JR \longrightarrow M_{d_k}(A/J)$. By an extension of $\rho_j$ by $\rho_i$   
we mean a representation $R/JR \longrightarrow \End_{A/J}(V)$ together with an exact sequence of $R/JR$-module $0 \longrightarrow \rho_i \longrightarrow V \longrightarrow  
\rho_j \longrightarrow 0$. Hence $V$ is in particular a free $A/J$-module of rank   
$d_1+d_2$. Such an extension defines an element in the module   
$\Ext^1_{R/JR}(\rho_j,\rho_i)$.  
  
\begin{theorem} \label{extension1} There exists a natural injective map of   
$A/J$-modules   
$$\iota_{i,j} : \Hom_A(\cal{A}_{i,j}/\cal{A}'_{i,j},A/J) \hookrightarrow   
\Ext^1_{R/JR}(\rho_j,\rho_i).$$  
\end{theorem}  
\index{Ziota@$\iota_{i,j}$, a morphism $\Hom(\cal{A}_{i,j}/\cal{A}'_{i,j},A/J) \longrightarrow \Ext^1_{R/JR}(\rho_j,\rho_i)$}  
\begin{pf}  The map $\iota_{i,j}$ is constructed as follows. Pick an   
$f \in \Hom_A(\cal{A}_{i,j}/\cal{A}'_{i,j}, A/J)$.   
We see it as a linear form $f : \cal{A}_{i,j} \longrightarrow A/J$, trivial on $\cal{A}'_{i,j}$.   
It induces a linear application, still denoted by $f$ :   
$M_{d_i,d_j}(\cal{A}_{i,j}) \longrightarrow M_{d_i,d_j}(A/J)$. We consider the following linear application  
$R \longrightarrow S \longrightarrow M_{d_i+d_j}(A/J)$,  
\begin{eqnarray} \label{explicite}  
x &\mapsto &   
\left(\begin{array}{cc} a_{i,i}(x) \pmod{J} & f(a_{i,j}(x))   
\\ 0 & a_{j,j}(x) \pmod{J} \end{array} \right)   
\end{eqnarray}   
By assumption, $f$ is trivial on $\sum_{k\neq i, j}  
M_{d_i,d_k}(\cal{A}_{i,k})M_{d_k,d_j}(\cal{A}_{k,j})\subset  
M_{d_i,d_j}(\cal{A}'_{i,j})$, hence Lemma \ref{conseq} (ii) and (iii) show that the map  
(\ref{explicite}) is an $A/J$-algebra homomorphism which is an extension of  
$\rho_j$ by $\rho_i$. As a consequence, it defines an   
element $\iota_{i,j}(f)$ in $\Ext^1_{R/JR}(\rho_j,\rho_i)$. \par  
It is clear by the Yoneda interpretation of the addition in $\Ext^1$ that the map $\iota_{i,j}(f)$ is linear.   
Let us prove that $\iota_{i,j}$ is injective. Assume $\iota_{i,j}(f)=0$. This means that the extension is split. As it factors by construction through $S/JS$, it is certainly split when restricted to any subalgebra of $S/JS$. Let us restrict it to the subalgebra $e_i S/JS e_j$ (without unit, but we can add $A/J(e_i+e_j)$ if we like).   
The restricted extension is   
\begin{eqnarray*} x &\mapsto &   
\left(\begin{array}{cc} 0 & f(a_{i,j}(x))   
\\ 0 & 0 \end{array} \right)   
\end{eqnarray*}     
and such an extension is split if and only if $f = 0$.   
\end{pf}   
\ps  
The construction above is a generalization of the one of Mazur and Wiles, directly giving the matrices of the searched extensions. We will give a second construction in the next subsection, more in the spirit of Ribet's one, which will   
realize the extensions constructed before as subquotient of some explicit $R$-modules.   
Our second aim is to characterize the image of $\iota_{i,j}$ and to verify that this image is the biggest possible subset   
of the above $\Ext$-group seen by $S$.   
  
\subsubsection{The projective modules $M_i$ and a characterization of the  
image of $\iota_{i,j}$.}

\label{Mi}  
  
We keep the assumptions and notations  of \S\ref{extensions}. For each $i$, we define the $A$-modules   
$$M_i:=SE_i=\oplus_{j=1}^r e_jSE_i.$$  
\index{Mi@$M_i$, an $R$-module}  
%As $A$-modules, we have $M_i=N_i \oplus e_iSE_i$.  
\noindent Note that $M_i$ is a left ideal of $S$, hence an $S$-module. It is even a  
projective $S$-module as $S=M_i\oplus S(1-E_i)$.   
  
\begin{theorem}\label{extension2} \label{extensions2}  
Let $j \in \{1,\dots,r\}$, $\PP$ a partition containing $\{j\}$ and $J$ an ideal containing $I_\PP$, then  
\begin{itemize}  
\item[(0)] there is a surjective map of $S$-modules   
$M_j/JM_j \longrightarrow \rho_j$ whose kernel has the property that   
any of its simple $S$-subquotients   
is isomorphic to $\rhob_k$ for some $k \neq j$.   
Moreover $M_j$ is the projective hull of $\rho_j$ (and of $\rhob_j$) in the category of $S$-modules.  
\end{itemize}  
Let $i \neq j \in \{1,\dots,r\}$, $\PP$ a partition   
containing $\{i\}$ and $\{j\}$, and $J$ is an ideal containing $I_\PP$. Then moreover: \begin{itemize}  
\item[(1)] the image of the map   
$\iota_{i,j}$ of Theorem \ref{extension1} is  
exactly $\Ext^1_{S/JS}(\rho_j,\rho_i)\subset \Ext^1_{R/JR}(\rho_j,\rho_i)$,  
\item[(2)] any $S/JS$-extension of $\rho_j$ by  
$\rho_i$ is a quotient of $M_j/JM_j\oplus \rho_i$ by an $S$-submodule whose  
every simple $S$-subquotient is isomorphic to some $\rhob_k$ for   
$k \neq j$.  
\end{itemize}  
\end{theorem}  
  
\begin{pf} First note that we may replace $A$ by $A/J$ and $S$ by $S/JS$, that is we may assume that $J=0$ in $A$ (which   
simplifies the notations). Indeed,   
$(S/JS)E_j \simeq M_j \otimes_A A/J = M_j \otimes_S S/JS = M_j/JM_j$. Hence assertions (1) and (2) are automatically proved for $A$ once they are proved   
for $A/J$. As for assertion (0), if we know the corresponding assertion for   
$A/J$, namely "The $S/JS$-module $M_j/JM_j$ is the projective hull of   
$\rho_j$",  then (0) follows, because the map of $S$-modules $M_j \longrightarrow     
M_j/JM_j \longrightarrow \rho_j$ is essential as $JS \subset \m S \subset \rad(S)$, and because $M_j$ is projective.\ps   
Assume that $\PP$ contains $\{j\}$ and that $J \supset I_{\PP}$.  
Let us consider the natural exact (split) sequence of $A$-modules  
\begin{equation}\label{exact} 0 \longrightarrow N_j:=\oplus_{i\neq j}  
e_iSE_j \longrightarrow M_j \longrightarrow e_jSE_j \longrightarrow 0.  
\end{equation}  
We claim that $N_j$ is an $S$-submodule of $M_j$, and that $M_j/N_j\simeq  
\rho_j$. It suffices to show that  
for $k\neq j$, $e_jSe_kN_j \subset N_j$. But this follows from Lemma  
\ref{conseq} (i), as $e_jSe_kSe_j=0$. As a consequence, $M_j/N_j \simeq  
e_jSE_j$ is an $S$-module, which is isomorphic to $\rho_j$ by Lemma \ref{conseq}  
(ii). \ps  
Let us prove the first assumption in (0).   
Recall that by lemma \ref{radical}, we have  
$$S/\rad(S) \simeq \prod_{i=1}^r\End_k(\rhob_i)$$  
(see the  
formula (\ref{isom}) in the proof of lemma \ref{idempotents}). So if $U$ is a simple $S$-subquotient of $N_j$, then $U \simeq \rhob_k$ for some $k \in \{1,\cdots,r\}$. But by construction, $e_j N_j=0$, hence $e_jU=0$, and $\rho_j(e_j)=1$ by Lemma \ref{idempotents}, so $k \neq j$ and we are done. \par  
We prove now that $M_j \longrightarrow \rho_j$ is a projective hull. We just have to show that this surjection is essential. If $Q \subset M_j$ is a $S$-submodule which maps  
surjectively to $M_j/N_j=e_j(M_j/N_j)$, then $e_jQ \subset e_jSE_j$ maps  
also surjectively to $M_j/N_j$, hence $e_jQ=e_jSE_j$. But then $E_j \in Q$,  
hence $Q=M_j$, and we are done.   
\ps  
  
Now we suppose that $\PP$ contains $\{i\}$ and $\{j\}$.  
Let us apply $\Hom_S(-,\rho_i)$ to the exact sequence (\ref{exact}).   
As $M_j$ is a projective $S$-module, it takes the form:   
$$ 0 \longrightarrow \Hom_S(\rho_j,\rho_i) \longrightarrow   
\Hom_S(M_j,\rho_i) \longrightarrow  
\Hom_S(N_j,\rho_i) \overset{\delta}{\longrightarrow} \Ext_S(\rho_j,\rho_i) \longrightarrow 0 $$  
  
We claim first that $\delta$ is an isomorphism.   
We have to show that any $S$-morphism $M_j \longrightarrow \rho_i$   
vanishes on $N_j$.   
But by Lemma \ref{conseq} (ii), if $k\neq j$ we have $e_k \rho_j=0$. We are done as $N_j=\sum_{k \neq j} e_kN_j$ by definition. \ps   
  
It is well known that if $f \in \Hom_S(N_j,\rho_i)$, we have the following   
commutative diagram defining $\delta(f)$:  
  
\begin{eqnarray} \label{diagrammext}  
\xymatrix{ & 0 \ar[r] & N_j \ar[d]_{f} \ar[r] & M_j  \ar[d]^{x \mapsto (x,0)} \ar[r] & \rho_j \ar[r] \ar[d]^{id} & 0 \\  
\delta(f): & 0 \ar[r] & \rho_i \ar[r]_{x\mapsto (0,x)} & \frac{M_j\oplus \rho_i}{Q} \ar[r] & \rho_j \ar[r] & 0 }  
\end{eqnarray}     
\noindent where $Q$ is the image of the $S$-linear map $u: N_j \longrightarrow  
M_j\oplus \rho_i, \, \, \, x \mapsto (x,0)-(0,f(x))$. This will prove (2) if  
we can show that each simple subquotient of $Q$ is isomorphic to some $\rhob_k$ with $k\neq j$. But as in the proof of (0), this  
follows from the fact that $e_jQ=u(e_jN_j)=0$. \ps  
  
We claim now that we have a sequence of isomorphisms   
\begin{eqnarray*} \Hom_S(N_j,\rho_i) & \isomo & \Hom_{e_iSe_i}(e_iSE_j/(\sum_{k\neq j, i}e_iSe_kSE_j),\rho_i) \\  
& \isomo & \Hom_{A}(\cal{A}_{i,j}/(\sum_{k\neq j, i}\cal{A}_{i,k}\cal{A}_{k,j}),A)\end{eqnarray*}  
The first one is induced by the restriction map, the fact that it is an isomorphism   
is a simple matter of orthogonal idempotents, using that $e_k\rho_i=0$ if $k\neq i$ and that   
$N_j=\oplus_{k\neq j}e_kSE_j$. The second one is induced by the Morita equivalence $A=E_iA \longrightarrow e_iSe_i=M_{d_i}(A)$. \ps  
It is now easy, using the diagram (\ref{diagrammext}) and the fact that   
$(M_j\oplus \rho_i)/Q$ is naturally isomorphic as $A$-module to $e_jSE_j\oplus e_iSE_i$,   
to check that in terms of the isomorphisms above, the map $\delta$ is exactly the map $\iota_{i,j}$ given by formula (\ref{explicite}), which proves (1).  
\end{pf}  
  
\begin{remark} By the same method, we could give an expression for the higher  
Ext-groups $\Ext_{S/JS}^n(\rho_j,\rho_i)$ in terms of the $\cal{A}_{i,j}$'s.  
For example, when $r=2$, the exact sequence  
(\ref{exact}) implies that $\Ext^2_{S/JS}(\rho_j,\rho_i) \isomo  
\Ext^1_{A/J}(\cal{A}_{i,j}\otimes A/J,A/J)$.   
However, for the usual applications of pseudocharacters, it is less interesting because when $n\geq 2$, the natural map $$\Ext^n_{S/JS}(\rho_j,\rho_i)  
\longrightarrow \Ext^n_{R/JR}(\rho_j,\rho_i)$$  
is not in general injective, and we usually only care about the   
extensions between the $\rho_i$'s in the category of   
representations of $R$, not of its auxiliary quotient $S$.   
\end{remark}  
  
%\begin{remark} We assumed that $A$ is henselian throughout the whole section  
%\ref{redlocext}, however all the results and proofs here remain true verbatim if we  
%only assume that $R/\Ker T$ is a GMA whose trace is $T$ (as in the example  
%of \S\ref{defmult1}). The only fact we have to change it to replace the  
%use of theorem \ref{structure} (i) by the following fact: if we have a  
%factorisation $T: R \longrightarrow S \longrightarrow A$ which is  
%Cayley-Hamilton, then $S$ is a GMA whose trace is $T$ if, and only  
%if, $R/\Ker T$ is. Indeed, the "only if" part is immediate. For the "if"  
%one, note that the kernel of natural surjection   
%$S \longrightarrow R/\Ker T$ is a  
%nil ideal, hence we see easily the GMA structure of $R/\Ker T$ lifts to $S$.  
%\end{remark}  
  
\label{projective}  
  
\begin{remark}[Dependence on S]  
\label{dependance}  
All the constructions of \S\ref{extensions} and \S\ref{projective}   
depend on the choice of a Cayley-Hamilton quotient $S$ of $(R,T)$. If $S_1 \longrightarrow S_2$ is a morphism  
in the category of Cayley-Hamilton quotients (cf.~\ref{CHquotient}), then it is surjective and we have obviously   
$\Ext^1_{S_2/IS_1}(\rho_j,\rho_i) \subset \Ext^1_{S_1/IS_1}(\rho_j,\rho_i)$. Thus our methods construct   
the biggest group of extensions when working with $S=S_0$, and the smallest when $S=R/\ker T$. We stress that even in the most  
favorable cases, the inclusion above may be strict : an example will be given in Remark~\ref{plusextension}  
below. However, we will not be able to get much information about the $\Ext^1_{S/JS}$ except when   
$S=R/\ker T$. On the other hand, as proved in Proposition~\ref{redloc}, the reducibility ideals do not depend on $S$.  
\end{remark}  

\begin{remark}\label{extSmorS} Assume we are under the assumptions of Theorem \ref{extension2}. We claim that for $i\neq j$, the natural inclusion
$$\Ext_{S/JS}(\rho_j,\rho_i) \longrightarrow \Ext_S(\rho_j,\rho_i)$$
is an isomorphism. Indeed, let $U$ be an $S$-extension of $\rho_j$ by $\rho_i$, we have to show that 
$JU=0$. But for $f \in J$, the multiplication by $f$ induces an $S$-linear map
	$$\rho_j \longrightarrow \rho_i,$$
which is necessarily $0$ as $\Hom_S(\rho_j,\rho_i)=\Hom_{S/JS}(\rho_j,\rho_i)=0$ by Lemma \ref{conseq} (i).
\end{remark}
   
\subsubsection{Complement: topology.} We keep the hypotheses of \S 
\ref{assredloc}. We assume moreover that $A$ is a Hausdorff topological ring such that the natural functor from  
the category of topological Hausdorff finite type $A$-modules    
to the category of $A$-modules has a section endowing $A$ with its topology.   
We fix such a section, hence every finite type $A$-module is provided with an   
Hausdorff  
$A$-module topology, and any $A$-linear morphism between two of them is continuous with closed image.   
For example, this is well known to be the case when $A$ is a complete noetherian local rings,  
and it holds also when $A$ is the local ring of a rigid analytic  
space at a closed point (see~\cite[\S2.4]{lissite}).\ps  
  
\begin{prop} \label{topoprop}  
Assume that $R$ is a topological $A$-algebra and that $T: R \longrightarrow A$ is continuous.  
\begin{itemize}  
\item[(i)] Let $I$ be an ideal containing $I_\PP$ where $\PP$ is a partition containing $\{i\}$.  
Then the representation $\rho_i : R/IR \longrightarrow M_{d_i}(A/I)$ is continuous.  
\item[(ii)] Let $I$ be an ideal containing $I_\PP$ where $\PP$ is a partition containing $\{i\}$ and  
$\{j\}$, $i\neq j$. If $A$ is reduced and $S=R/\Ker T$, then the image of $\iota_{i,j}$ of Theorem  
\ref{extension1} falls into the $A$-submodule of continuous extensions   
$\Ext^1_{R,\text{cont}} (\rho_i,\rho_j)$.  
\end{itemize}  
\end{prop}   
  
\begin{pf} By Lemma \ref{conseq} (ii), we can find $e \in R$ such that for  
all $x \in R$, $T_i(x)=T(ex)$ (any lift of the element $e_i \in S/JS$  
{\it loc. cit.} works for $e$), which proves (i). Let us show (ii). Fix $f\in  
\Hom_A(\cal{A}_{i,j}/\cal{A}'_{i,j},A/I)$. By the formula (\ref{explicite}) defining  
$\iota_{i,j}(f)$, it suffices to show that the natural maps   
$\pi_{i,j}: R \longrightarrow \cal{A}_{i,j}$, $x \mapsto E_ixE_j$, are continuous. Note  
that this makes sense because by Theorem \ref{structureGMA} (iii),    
the $\cal{A}_{i,j}$'s are finite type $A$-modules. Let us choose a family  
of $A$-generators $x_1,\dots,x_n$ of $\cal{A}_{j,i}$. As $T: S  
\longrightarrow A$ is faithful by assumption, and by Lemma \ref{lesaij}, the  
map $$\mu: \cal{A}_{i,j} \longrightarrow \prod_{s=1}^n A,\, \, \, x\mapsto  
\mu(x)=(T(xx_s))_s,$$  
is injective. By assumption on the topology of finite type $A$-modules, the  
map above is an homeomorphism onto its image. It suffices then to prove  
that $\mu \circ \pi_{i,j}$ is continuous, what we can check componentwise. But for each $s$,  
$(\mu \circ \pi_{i,j})_s$ is the map $x \mapsto T(exfg)$, where $(e,f,g)\in  
R^3$ denote any lift of $(E_i,E_j,x_s) \in S^3$. This concludes the proof. \end{pf}  
  
\subsection{Representations over $A$}  
\label{representations} We keep the notations and hypotheses of \S\ref{hypomfree}: $A$ is local henselian and $d!$ is invertible in $A$. In this subsection we are mainly concerned with the following natural question   
which is a converse of Example \S\ref{mainexample} :  
if $T: R \longrightarrow A$ is a residually multiplicity free   
pseudocharacter of dimension $d$, does $T$ arise as the trace of a true representation $R \longrightarrow M_d(A)$ ? \ps  
  
When $T$ is residually absolutely irreducible, the theorem of Nyssen and Rouquier (\cite{Nys}, \cite[corollaire 5.2]{Rou}) we   
recalled in \S\ref{mainexample} shows that the answer is yes.  
Although for a given residually multiplicity free pseudocharacter, it may be difficult to determine if it arises as the trace of a representation (see next subsection for interesting particular cases), it turns out that there is a simple sufficient and (almost) necessary condition on $A$ for this to be true for every residually multiplicity free pseudocharacter of dimension $d$ on $A$.  
  
\begin{prop} \label{factfree}Assume that $A$ is a factorial domain. Then any   
residually multiplicity-free pseudocharacter $T: R \longrightarrow A$ of dimension $d$ is the trace of a representation $R \longrightarrow M_d(A)$.  
\end{prop}  
\begin{pf} We use the notations of \S\ref{defmult1} for $T$.  
As $A$ is a domain, its total fraction ring is a field $K$. By the point (i) of Theorem~\ref{structure}, there is a data $\cal{E}$ on $R/\Ker T$ that makes it a GMA,  
and by the point (ii) of the same theorem,  
there is an adapted (to $\cal{E}$) representation $\rho : R/\ker T   
\longrightarrow M_d(K)$ whose image  
is the standard GMA (see example~\ref{gma}) attached to some fractional ideals $A_{i,j}$ of $K$, $i,j \in \{1,\dots,r\}$.  
  
Let $v$ be a (discrete) valuation of $A$. Let $v_{i,j}$ be the smallest integer of   
the form $v(x)$ for a nonzero $x \in A_{i,j}$, this makes sense since $A_{i,j}$ is a fractional ideal. Because   
$A_{i,i}=A,\ A_{i,j}A_{j,i}\subset A$ and $A_{i,j}A_{j,k} \subset   
A_{i,k}$, we have  
\begin{equation} \label{vij} v_{i,i}=0,\ \ \ v_{i,j}+v_{j,i} \geq 0,\ \ \ v_{i,j}+v_{j,k} \geq v_{i,k} \end{equation}  
  
Because $A$ is factorial, there exists for each $i$ an element $x_i \in   
K^\ast$ such that $v(x_i)=v_{i,1}$ for every valuation $v$ of $A$. Let   
$P$ be the following diagonal matrix :  
$$\left( \begin{array}{cccc} x_1 {\rm Id}_{d_1} & & &\\  
& x_2 {\rm Id}_{d_2} & &\\  
& & \ddots &\\  
& & & x_r {\rm Id}_{d_r} \end{array} \right)$$  
Let $\rho':=P^{-1} \rho P$. Then $\rho'$ is adapted to $\cal{E}$ and   
its image is the standard GMA attached to the modules $A'_{i,j}=x_j   
x_i^{-1} A_{i,j}$.  
  
If $x \in A'_{i,j}$, and $v$ is a valuation on $A$, we have $v(x) \geq   
v(x_j)-v(x_i) + v_{i,j} = v_{j,1} - v_{i,1} + v_{i,j}$ which is   
nonnegative by~(\ref{vij}). Hence $x \in A$ since $A$ is factorial,  
and $A'_{i,j} \subset A$. That is, $\rho'$ is a representation $R   
\longrightarrow M_d(A)$ of trace $T$.  
\end{pf}  
  
\begin{remark} Let $A$ be a valuation ring of fraction field $K$,   
with valuation $v : K^\ast \longrightarrow \Gamma$, where $\Gamma$ is a totally ordered groups and assume   
$v(K^\ast)=\Gamma$.  
  
Then the proof above shows that the result of Proposition~\ref{factfree} holds also for this ring $A$ if the ordered group $\Gamma$ admits infima. Indeed, it suffices to define $v_{i,j}$ to be the infimum of the $v(x)$ with $x \in \cal{A}_{i,j}$ nonzero, and to choose $x_i \in K^*$ such that $v(x_i)=v_{1,i}$, which is possible by the assumption 
$v(K^\ast)=\Gamma$.   
  
Consider for example a valuation ring $A$ as above, with $\Gamma=\R$ (such a valuation ring exists by  
\cite[chapitre VI, \S3, n${}^o$ 4, example 6]{BouAc}). Then the result of Proposition~\ref{factfree} holds for $A$,   
even if $A$ is a non factorial domain ($A$ has no irreducible elements !). Note however that $A$ is non noetherian.    
  
If on the contrary we do not assume that $\Gamma$ admits infima, the result fails  
as showed for the ring $\anneau_{\C_p}$ in \cite[remark 1.14]{BCKL}.  
\end{remark}  
  
We are now interested in the converse of Proposition~\ref{factfree}. Because of the remark above, we shall assume that $A$ is noetherian.  
  
\begin{theorem} Assume $d \geq 2$ and $A$ noetherian (in addition of being local henselian). If each residually multiplicity free pseudocharacter of dimension $d$ is the trace of a representation $R \longrightarrow M_d(A)$,  
then $A$ is factorial.  
\label{factoriel}  
\end{theorem}  
  
\begin{pf} We claim first that the hypothesis implies the following purely module-theoretical assertion on $A$: \ps  
  
{\it For every $A$-modules $B$ and $C$, and every morphisms of $A$-modules  
$\phi : B \otimes C \longrightarrow m$ such that  
\begin{equation} \label{phimn}  
\phi(b,c) b'= \phi(b',c) b, \text{ for any }b,b' \in B,\ c \in   
C,\end{equation}  
there exist two morphisms $f: B \longrightarrow A$, and $g: C \longrightarrow   
A$ such that $\phi(b \otimes c)=f(b) g(c)$ for any $b \in B$, $c \in   
C$.}\ps

Let us prove the claim. Let $B, C$ be two $A$-modules with a morphism $\phi : B \otimes C\longrightarrow A$ 
satisfying the property above. Set $\cal{A}_{1,2}:=B$, $\cal{A}_{1,2}=C$,  
$\cal{A}_{i,i}=A$ for $i=1, 2$, $\phi_{1,2,1}:=\phi$, $\phi_{2,1,2}(c   
\otimes b)=\phi(b \otimes c)$, and $\phi_{i,i,j}$, $\phi_{i,j,j}$ be   
the structural morphism. Then we check at once that these $\cal{A}_{i,j}$'s and $\phi_{i,j,k}$'s 
satisfy the properties (COM), (UNIT), and (ASSO) (see~\ref{reprGMA}), 
and thus defines a GMA $(R,\cal{E})$ whose they are the structural modules and morphisms. 
As $\phi(B \otimes C) \subset m$, we are in the case of Example~\ref{importantexample}, 
and the trace function $T: R \longrightarrow A$ of $(R,\cal{E})$ is a   
residually multiplicity free pseudocharacter of dimension $d$.

 \par   
The hypothesis of the theorem then implies that there is a   
trace representation $R \longrightarrow M_d(A)$.  
Because $A$ is local, every finite-type projective $A$-module is free and by Lemma~\ref{traceadapted},  
there is an adapted (to $\cal{E}$) representation $\rho: R \longrightarrow   
M_d(A)$, that is an element of $G(A)$ where the functor   
$G=G_{R,\cal{E}}$ is the one defined in \S\ref{reprGMA}. By   
Proposition~\ref{FetG}, $F(A)$ is not empty. If $(f_{i,j}) \in F(A)$, then by definition $(f,g):=(f_{1,2},f_{2,1})$ satisfies the claim, and we are done. \ps  
Using the assertion above, we will now prove in three steps that $A$ is a   
factorial domain. \ps

{\bf First step: } $A$ is a domain. \par \noindent  
Choose an $x \in A$, $x \neq 0$, and let $I$ be its annihilator. Set $B   
= A/xA$, $C=I$ and let $\phi : B \otimes C \longrightarrow A$ be the   
morphism induced by the multiplication in $A$. Then $\phi(B \otimes   
C)=I \subset m$ and the property (\ref{phimn}) is obvious. Thus there   
exist $f : B \longrightarrow A$ and $g : C  \longrightarrow A$ such that   
$\phi(b \otimes c)=f(b) g(c)$ for any $b \in B$, $c \in C$. As $xC=xI=0$, we have $x g(C)=0$ hence $g(C) \subset I$. 
As $x A/xA =  xB=0$ we have also $f(B) \subset I$. Hence $I=\phi(B\otimes C) = f(B)   
g(C) \subset I^2$. Because $A$ is local and noetherian,  
this implies $I=0$.  Hence $A$ is a domain. \ps

{\bf Second step:} if $A$ is a domain, then $A$ is normal. \par \noindent  
Let $K$ be the fraction field of $A$. Assume, by contradiction, that   
$A$ is not normal, and let $B \subset K$ be a finite $A$-algebra containing $A$, but different from $A$. 
Let $C=\{x \in K\ ;   
\; xB \subset A\}$. We have then:  
\begin{itemize}  
\item[i.] by definition, $C$ is a  
$B$-submodule of $K$ (hence an $A$-module too).  
\item[ii.] $C \subset A$, because $1 \in B$. Hence $C$ is an $A$-ideal.  
\item[iii.] We have $C \subset m$. Indeed, $C$ is an $A$-ideal by ii. As $A$ is local, 
we only have to see that $C \neq A$. But if  
$1 \in C$, $B \subset C \subset A$ by i. and ii. , which is absurd.  
\item [iv.] $C$ is non zero : if $(p_i/q_i)$ is a finite  
family of generators of $B$ as an $A$-module, with $p_i, q_i \in A$,   
$q_i \neq 0$, then $0 \neq \prod_i q_i \in C$.  
\end{itemize}  
Now let $\phi: B\otimes_A C \longrightarrow K$ be the map induced by the multiplication in $K$. By iii.   
$\phi(B \otimes C) \subset m$. Moreover, hypothesis ~(\ref{phimn}) is obviously satisfied. 
Thus there exist two morphisms $f: B \longrightarrow A$ and $g: C \longrightarrow   
A$ such that $\phi(b \otimes c)=f(b) g(c)$ for any $b \in B$, $c \in   
C$. Since $B \otimes_A K = K$, $f \otimes K: K \longrightarrow K$ is the   
multiplication by some element $x \in K^\ast$, and so is $f$. As $C\otimes_A K =K$ by iv. , $g$ has to be the multiplication by $x^{-1}$. We thus get   
\begin{equation}\label{deuxrelations} xB  \subset A, \ \ \ \ {\rm and}\ \ \ \ \ x^{-1} C \subset A.\end{equation}  
The first relation implies $x \in C$, so  
$1 \in x^{-1}C$. As $x^{-1}C$ is a $B$-module, $B \subset x^{-1} C$   
and by the second relation, $B \subset A$, which is absurd. 
(The reader may notice that this step does not use the noetherian hypothesis). \ps

{\bf Third step :} if $A$ is a normal domain, then $A$ is factorial.\par \noindent  
We may assume that the Krull dimension of $A$ is at least $2$, because   
a normal noetherian domain of dimension $\leq 1$ is a discrete valuation ring, hence   
factorial. Let $C$ be an invertible ideal of $A$, and set $B=m C^{-1} \subset K$.   
Let $\phi: B \otimes_A C \longrightarrow m$ be induced by the multiplication in $K$. Then reasoning as in the second step above, we see that there is an $x \in K^\ast$  
such that $x m C^{-1} \subset A$ and  
$x^{-1} C \subset A$, as in (\ref{deuxrelations}). \par

\renewcommand{\div}{\text{div\,}}  
  
Now, since $A$ is normal and noetherian, it is  completely integrally closed, and even a Krull ring   
(\cite[chapter VII, \S1, n${}^0$ 3, corollary]{BouAc}).  
Recall from \cite[chapter VII, \S1, n${}^0$ 2, Theorem 1]{BouAc} the  
 ordered group $D(A)$ of divisorial fractional ideals of $A$, and the projection $\div$ 
from the set of all fractional ideals of $A$ to $D(A)$. Since $x^{-1} C \subset A$, 
we have 
(using  \cite[chapter VII, \S1, n${}^0$ 2, formula (2)]{BouAc}) $$\div x^{-1} + \div C = \div(x^{-1} C) \geq 0,$$ that is   
$\div C \geq \div x$. From $x m C^{-1} \subset A$ we have $$\div x + \div m + \div C^{-1} \geq 0,$$  
but since $m$ has height greater than $2$, $\div m=0$ by \cite[chapter VII, \S1, n${}^0$ 6, corollary (1)]{BouAc}, and since $A$ is completely integrally closed,  
$\div C^{-1} = - \div C$ by   
\cite[chapter VII, \S1, n${}^0$ 2, corollary]{BouAc}).   
Hence $\div x \geq \div C$.  
Thus $\div x = \div C$, and if $C$ is divisorial, then $C=Ax$ is principal. 
But a Krull ring where every divisorial ideal is principal is factorial,   
cf.   \cite[chapter VII,  \S3, n${}^0$ 1]{BouAc}.  
\end{pf}  
%\begin{remark} The problem of finding a necessary and sufficient condition on a noetherian local henselian ring $A$ such that all (non necessarily residually multiplicity free) pseudocharacter comes from a representation stands beyond the scope of this article. Let us just say that "if $A$ is a regular domain of dimension $\leq 2$, $K:={\rm Frac}(A)$, and $R \subset M_d(K)$ is a finite $A$-algebra, then there exists $P \in \GL_d(K)$ such that $PRP^{-1} \subset M_d(A)$". Indeed, if $M:=R.A^d \subset K^d$, then its bidual $M^{**} \subset K^d$ is stable by $R$ and free of rank $d$ over $A$ by \cite[Ch VII ?]{BouAc} (of course, if $\dim(A)\leq 1$, $M$ already has this property). \end{remark}  
  
When a trace representation $\rho : R \longrightarrow M_d(A)$ does exist, we may ask what are its kernel and image. 
In some favorable cases, we can give a satisfactory answer :  
  
\begin{prop} \label{kert} Assume $A$ is reduced with total fraction ring $K$ a finite product of fields  
$K_s$. Let   
$T:R \longrightarrow A$ be a residually multiplicity free pseudocharacter and assume $T \otimes  
K_s$ irreducible for each $s$.   
If $\rho : R \longrightarrow M_d(A)$ is a trace representation then $\ker \rho=\ker T$  
and $\rho(R)\otimes K = K[\rho(R)]=M_d(K)$.  
\end{prop}  
\begin{pf} We have obviously $\ker \rho \subset \ker T$. Set $S:=\rho(R)  
\subset M_d(A)$, which is a   
Cayley-Hamilton quotient of $(R,T)$. To show that $T: S \longrightarrow A$  
is faithful, it suffices to show the last statement. By the irreducibility assumption and   
Theorem \ref{structure} (iii), $S \otimes K$ is (trace) isomorphic to  
$M_d(K)$. As a consequence, the injective map $\rho \otimes K: S\otimes K  
\longrightarrow M_d(K)$ is an isomorphism, which concludes the proof.  
\end{pf}  
  
\begin{remark} \label{plusextension}   
The proof above shows in particular that under the hypotheses of the proposition, the only Cayley-Hamilton quotient of $R$ that is torsion free as an $A$-module is $R/\ker T$. We stress the reader that we cannot omit the hypothesis ``torsion free''. Here is a counter-example : with the notations of the proof of Theorem~\ref{factoriel},  
take $A=\Z_p$, and set $B=\Z_p$, $C=\Z_p \oplus \Z/p\Z$ and let $\phi : B \otimes C \longrightarrow \Z_p$ be   
defined by $\phi(b \otimes (c,c'))=pbc$. As it is clear that $\phi$ satisfies~(\ref{phimn}), those data define a   
GMA $R$ of type $(1,1)$. Its trace function $T$ is a Cayley-Hamilton residually multiplicity free   
pseudocharacter. Hence $R$ is Cayley-Hamilton, we have $R=S_0$ in the notation of \S\ref{CHquotient}, but $R \neq R/\ker T$ because $\ker T \simeq \Z/p\Z$. Moreover this example provides a case where   
$\Ext^1_{S_0/pS_0}(\rhob_1,\rhob_2)$ has dimension $2$ whereas $\Ext^1_{(R/\Ker T)/p(R/\Ker T)}(\rhob_1,\rhob_2)$ has dimension $1$.  
\end{remark}  
  
\subsection{An example: the case $r=2$} ${}^{}$ \ps  
\label{exampler2}  
Let $A$ be a reduced, noetherian, henselian local ring and $T: R \longrightarrow A$ be a multiplicity free,   
$d$-dimensional, pseudocharacter. As before, $K$ is the total fraction ring of $A$, which is a finite product of fields $K_s$. In this subsection, we investigate the consequences of our general results   
in the simplest case where $\bar T$ is the sum of only two irreducible pseudocharacters   
$\tr \rhob_1$ and $\tr \rhob_2$. Note that in this case,  
the only reducibility locus is the total one, of ideal $I_\PP$ with $\PP=\{\{1\},\{2\}\}$.  
  
Let $S$ be a given Cayley-Hamilton quotient of $(R,T)$  
We are first interested in giving a lower bound on the dimension of $\Ext^1_{S/mS}(\rho_1,\rho_2)$, hence of  $\Ext^1_{R/mR}(\rho_1,\rho_2)$.  
  
\begin{prop} Let $n$ be the minimal number of generators of $I_\PP$. Then  $$(\dim_k \Ext^1_{S/\m S}(\rhob_1,\rhob_2))(\dim_k \Ext^1_{S/\m S}(\rhob_2,\rhob_1)) \geq n$$   
\end{prop}  
\begin{pf}  By Remark~\ref{dependance}, we may and do assume $S=R/\ker T$.  
Let $\rho : R/\ker T \longrightarrow M_d(K)$ be as in Theorem~\ref{structure} whose we use notations.   
If $i\neq j$, let $n_{i,j}$ be the minimal number of generators of the finite $A$-module $A_{i,j}$.   
By Theorems \ref{extension1} and \ref{extension2}, and   
by Nakayama's lemma, $n_{i,j}=\dim_k  \Ext^1_{S/\m S}(\rhob_j,\rhob_i)$. On the other hand, since $I_\PP=A_{1,2}A_{2,1}$, we have $n_{1,2} n_{2,1} \geq n$, and the proposition follows.  
\end{pf}   
  
This easy observation is one of the main theme of the paper: to produce many extensions of $\rhob_1$ by $\rhob_2$ we shall not only produce a pseudodeformation of $\tr \rhob_1 + \tr \rhob_2$, but to do it sufficiently non trivially so that the reducibility locus of the pseudodeformation has a big codimension. A very favorable case occurs of course when $I_\PP$ is the maximal ideal $m$.  
In this case, the above result writes  
%(using Nakayama's lemma for $m$)  
$$(\dim_k \Ext^1_{S/\m S}(\rhob_1,\rhob_2))(\dim_k \Ext^1_{S/\m S}(\rhob_2,\rhob_1))   
\geq \dim_k m/m^2 \geq \dim A.$$   
  
When moreover $T$ is the trace of a true representation, we can say more:  
  
\begin{prop} Assume that each $T \otimes K_s$ is irreducible, that $I_{\PP}$ is the maximal ideal and that there is a   
trace representation $R \longrightarrow M_d(A)$, then   
$$\max(\dim_k \Ext^1_{S/\m S}(\rhob_1,\rhob_2),\dim_k \Ext^1_{S/\m S}(\rhob_2,\rhob_1)) \geq \dim_k \m/\m^2$$  
\end{prop}  
\begin{pf} Again we may and do assume that $S=R/\ker T$. Moreover we also have $\rho(R)=R/\ker T=S$  
 by Proposition~\ref{kert}. By Lemma~\ref{traceadapted}, and Lemma~\ref{imageadapted} we may assume that the image of $\rho$ is a standard GMA  
attached to ideals $A_{1,2}$, $A_{2,1}$ of $A$.  
Then $A_{1,2}$ and $A_{2,1}$ are ideals of $A$ such that $A_{1,2}A_{2,1} = I_\PP =\m$. Hence   
$\m \subset A_{1,2}$ and $\m \subset A_{2,1}$, but we cannot have $A_{1,2}=A_{2,1}=A$,   
hence one of those ideals is $\m$. The proposition follows.\end{pf}  
\begin{remark} The inequality above does not hold when $T$ has no representation over $A$.   
Indeed set $A=k[[x,y,z]]/(xy-z^2)$ which is a complete noetherian normal local domain, but not factorial.  
 Let $K$ be its fraction field, and $A_{1,2} = y A + z A$, $A_{2,1}=\frac{x}{z} A + A$ in $K$, $A_{1,1}=A_{2,2}=A$.   
Let $R$ be the standard GMA of type $(1,1)$ associated to these $A_{i,j} \subset K$. As $A_{1,2}A_{2,1}=m$, the trace $T$ of $R$ is an $A$-valued residually multiplicity free pseudocharacter. Its reducibility locus   
$I_\PP=A_{1,2}A_{2,1}=(x,y,z)=\m$ is the maximal ideal of $A$, and $T \otimes K$ is obviously irreducible  
but $\m/\m^2$ has   
dimension $3$, whereas $\dim_k \Ext^1_{R/\m R}(\rhob_1,\rhob_2)=\dim_k \Ext^1_{R/\m R}(\rhob_2,\rhob_1)=2$.  
\end{remark}  
  
We now give a result relating the Ext groups and  
the existence of a  trace representation over $A$ :  
  
\begin{prop}  
Assume that each $T \otimes K_s$ is irreducible.  
The two following assertions are equivalent:  
\begin{itemize}  
\item[(i)] There is a representation $\rho : R \longrightarrow M_d(A)$ whose trace is $T$,   
and whose reduction modulo $m$ is a non split extension of $\rhob_1$ by $\rhob_2$,  
\item[(ii)] $\Ext^1_{(R/\ker T)/m (R/\ker T)}(\rhob_1,\rhob_2)$ has $k$-dimension $1$.  
\end{itemize}  
Moreover, if those properties hold, then the representation $\rho$ in (i) is unique up to isomorphism.  
%that is, up to multiplication by a scalar of $K^\ast$  
\end{prop}  
\begin{pf} Let us prove first (i) $\Rightarrow$ (ii). Fix $\rho$ as in (i).   
By reasoning as in the proof of the proposition above, we can assume that $\rho(R)$ is the standard GMA attached to some ideals $A_{1,2}$, $A_{2,1}$ of $A$ for image, and has $\ker T$ for kernel. Hence   
$$(\rho \otimes k)(R\otimes k)=\left( \begin{array}{cc} M_{d_1}(k) & M_{d_1,d_2}(\overline{A_{1,2}}) \\ M_{d_2,d_1}(\overline{A_{2,1}}) & M_{d_2}(k)   
 \end{array} \right) $$ where $\overline{A_{i,j}}$ is the image of the ideal $A_{i,j}$ in $A/\m = k$.   
The hypothesis tells us  
that $\overline{A_{1,2}}=0$ and $\overline{A_{2,1}} \neq 0$, hence $A_{1,2} \subset \m$ and $A_{2,1}=A$. But   
$$\Ext^1_{(R/\ker T)/\m (R/\ker T)}(\rhob_1,\rhob_2) \simeq \Hom_k(A_{2,1},k) = k,$$ which is (ii). \par  
Let us prove (ii) $\Rightarrow$ (i). Let $\rho : R \longrightarrow M_d(K)$ be a representation as in Theorem~\ref{structure}, (ii), whose kernel is $\ker T$ and whose image is the standard GMA of type $(d_1,d_2)$ attached to fractional ideals $A_{1,2}$, $A_{2,1}$ of $A$.  
Since $$k \simeq \Ext^1_{(R/\ker T)/\m (R/\ker T)}(\rhob_1,\rhob_2) \simeq \Hom_k(A_{2,1},k),$$ we have $A_{2,1}/\m A_{2,1} \simeq k$ hence   
by Nakayama's lemma $A_{2,1}= f A$ for some $f \in K$. By Theorem \ref{structure} (iii), $A_{2,1}K=K$, hence $f \in K^\ast$. Then, if we change the basis of $A^d$, keeping the $d_1$ first vectors and multiplying the $d_2$ last vectors by $f$, we get a new representation $\rho':R \longrightarrow \Gl_2(A)$ whose image is the standard GMA attached to $A'_{i,j}$, with $A'_{2,1}=A_{2,1}/f=A$, hence $A'_{1,2}\subset m$. It is then clear that the reduction modulo $m$ of that representation is a non split extension of $\rhob_1$ by $\rhob_2$. We leave the last assertion as an exercise to the reader.  
\end{pf}  
In the same spirit, we have  
\begin{prop}  
Assume that each $T \otimes K_s$ is irreducible.  
The two following assertions are equivalent:  
\begin{itemize}  
\item[(i)]  $\Ext^1_{(R/\ker T)/\m (R/\ker T)}(\rhob_1,\rhob_2)$ and  $\Ext^1_{(R/\ker T)/\m (R/\ker T),T}  
(\rhob_2,\rhob_1)$ have $k$-dimension $1$.  
\item[(ii)] The reducibility ideal $I_\PP$ is principal, with a non-zero divisor generator.  
\end{itemize}  
\end{prop}  
\begin{pf}  
We will use the notations $\rho$ and $A_{1,2}$, $A_{2,1}$ of the part (ii) $\Rightarrow$ (i) of the proof of the above proposition.\par  
Proof of (i) $\Rightarrow$ (ii). Reasoning   
as in the proof  of the proposition above, we see that $A_{1,2}=fA$ and $A_{2,1} =f'A$ with $f,f' \in K^\ast$.  
Hence $I_\PP=A_{1,2}A_{2,1}=ff' A$ with $ff' \in K^\ast \cap A$. Hence the ideal $I_\PP$ is generated by $ff'$ which is  not a zero divisor. \par  
Proof of (ii) $\Rightarrow$ (i). By hypothesis, $A_{1,2}A_{2,1}=fA$ with $f$ not a zero divisor.   
Hence there is a family of $a_i \in A_{1,2}$, $b_i \in A_{2,1}$ such that $\sum_{i=1}^n a_i b_i = f$. Let $x \in A_{1,2}$,   
then $x b_i \in fA$ so we can write $a b_i = f x_i$ for a unique $x_i \in A$. Hence   
$$f x = \sum_i (a_i b_i) x= \sum_i a_i (x b_i)= \sum_i a_i f x_i.$$   
Because $f$ is not a zero divisor, $x = \sum a_i x_i$. This shows that the $a_i$ generate  
$A_{1,2}$, and the morphism $A^n \surjection A_{1,2}$, $(x_1,\dots,x_n) \mapsto \sum a_i x_i$ has a section   
$x \mapsto (x_1,\dots,x_n)$. Hence $A_{1,2}$ is projective of finite type, hence free, and since $A_{1,2} \subset K$,   
it is free of rank one. The same argument holds of course for $A_{2,1}$, and we conclude by Theorems \ref{extension1} and \ref{extension2} (i) applied to $J=m$ and $S=R/\Ker T$. \end{pf}

\subsection{Pseudocharacters with a symmetry}\label{pseudosymmetry}  
  
\subsubsection{The set-up}  
\label{setup}  
  
In this section we return to the hypotheses of \S\ref{defmult1} :   
$A$ is a local henselian ring where $d!$ is invertible,   
$T:\ R \longrightarrow A$ is a $d$-dimensional pseudocharacter   
residually multiplicity free. \par  
  
Moreover, in this section, we suppose given an automorphism of $A$-module $\tau: R \longrightarrow R$, which is either a morphism or an   
anti-morphism of $A$-algebra and such that $\tau^2={\rm id_R}$. We note that in both cases   
$T \circ \tau$ is a pseudocharacter on $R$ of dimension $d$, and we assume  
\begin{eqnarray} \label{symmetry} T \circ \tau=T. \end{eqnarray}  
If $B$ is any $A$-algebra, and $\rho : R \longrightarrow M_{n}(B)$ is any   
representation, then we shall denote by $\rho^\bot$ the representation   
$\rho \circ \tau : R \longrightarrow M_n(B)$ if $\tau$ is a morphism of algebra, and   
${}^t(\rho \circ\tau)$ if $\tau$ is an anti-morphism of algebra. Note that   
$\rho^\bot$ is a representation whose trace is $(\tr \rho)\circ \tau$. If $\rho : R \longrightarrow M_d(K)$ is a semisimple  
representation of trace $T$, where $K$ is a field, then the hypothesis~(\ref{symmetry}) is equivalent to   
\begin{eqnarray} \label{symmetry2} \rho^\bot \simeq \rho \end{eqnarray}  
  
The hypothesis~(\ref{symmetry}) also implies that   
$\bar T \circ \tau = \bar T$, hence $\rhobar^\bot \simeq \rhobar$.   
Thus there is a permutation $\sigma$ of $\{1,\dots,r\}$ of order two, such   
that for each $i \in \{1,\dots,r\}$, we have  
$\bar T_i \circ \tau=T_{\sigma(i)}$, and equivalently,  
 $\rhob_i \circ \tau \simeq \rhob_{\sigma(i)}$. This implies $d_i=d_{\sigma(i)}$. 
 
\begin{remark}\label{quotsymker} \begin{itemize}
\item[(i)] We check at once that the ideal $\ker T \subset R$ is stable by $\tau$, hence $\tau$ induces an automorphism, or an anti-automorphism, on $R/\ker T$ 
which we will still denote by $\tau$.  
\item[(ii)] In the same vein, we have for each $x \in R$ an equality of characteristic polynomials  $$P_{x,T}=P_{\tau(x),T},$$
hence $\tau$ factors also through the maximal Cayley-Hamilton quotient of $R$ (see \S\ref{CHquotient}).
\end{itemize}
\end{remark}
  
\subsubsection{Lifting idempotents} In the following lemma, $A$ is a local henselian ring in which $2$ is invertible.  
\begin{lemma}   
Let $S$ be an integral $A$-algebra, $\tau$ an $A$-linear involution of $S$ which is  
either a morphism or an anti-morphism of algebra,   
and let $I \subset \rad(S)$ be a two-sided ideal of $S$ such that $\tau(I)=I$.  \par  
Let $(\epsilon_i)$, $i=1,\dots,k$, be a family of orthogonal idempotents in $S/I$, and assume that the set $\{\epsilon_i, i=1,\dots, k\} \subset R$ is   
stable by $\tau$. Then there is a family of orthogonal idempotents $(e_i)$,  
$i=1,\dots, k$, lifting $(\epsilon_i)$ and such that $\{e_i, i=1, \dots, k\}$ is stable by $\tau$.  
\end{lemma}  
\begin{pf} We prove the lemma by induction on $k$. It is obvious for $k=0$. Assume it is true  
for any $k'<k$. We will consider two cases.\ps   
  
{\bf First case :} $\tau(\epsilon_1)=\epsilon_1$. Let $x$ be any lifting of $\epsilon_1$ in $S$.  
Set $y=(x+\tau(x))/2$. Then $\tau(y)=y$. Let $S_1$ be the $A$-subalgebra of $S$   
generated by $y$. It is a commutative, finite $A$-algebra on which $\tau=\Id$.    
Set $I_1 := I \cap S_1$. Then $S_1/I_1 \subset S/I$ and $S_1/I_1$ contains the reduction   
of $y$ which is $\epsilon_1$. As $A$ is henselian, there exists $e_1 \in S_1$ an idempotent lifting $\epsilon_1$. Then $\tau(e_1)=e_1$. \par  
The $A$-subalgebra\footnote{  
Recall that if $e \in S$ is an idempotent and $I$ a two-sided ideal of $S$, then $eIe=I\cap eSe$ and $\rad(eRe)=e\,\rad(R)e$.} $S_2 :=(1-e_1)S(1-e_1)$ is stable by $\tau$, and if $I_2 := I \cap S_2$,   
then $S_2/I_2 \subset S/I$ contains the family $\epsilon_2, \dots, \epsilon_k$ that is stable by $\tau$. By induction hypothesis, this family can be lifted as an orthogonal family of idempotents $e_2,\dots, e_k$, stable by $\tau$, in $S_2$, and then $e_1,\dots, e_k$ is an orthogonal family of idempotents lifting $\epsilon_1,\dots, \epsilon_k$ in $S$ that is stable by $\tau$. The lemma is proved in this case. \ps  
   
{\bf Second case :} $\tau(\epsilon_1)\neq \epsilon_1$. Then up to renumbering, we may assume that $\tau(\epsilon_1)=\epsilon_2$. We claim that   
\begin{center} {\it there are two orthogonal idempotents $e_1$ and $e_2$ in $S$  lifting $\epsilon_1$ and $\epsilon_2$ respectively, such that $\tau(e_1)=e_2$.} \end{center}   
This claim implies the lemma since we may apply the induction hypothesis to lift the family $\epsilon_3,\dots,\epsilon_k$ in $(1-(e_1+e_2)S(1-(e_1+e_2))$ by the same reasoning as above. Moreover, in order to prove the claim, we may assume that $\varepsilon_1+\varepsilon_2=1$. Indeed, set $\epsilon=\epsilon_1+\epsilon_2$. This is an idempotent of $S/I$ stable by $\tau$. By the first case above, there is an idempotent $e$ in $S$ lifting $\epsilon$ and such that $\tau(e)=e$. Replacing $S$ with $eSe$, and $I$ with $I \cap eSe$, we have now $\epsilon_1+\epsilon_2=1$, and we are done. To prove our claim, we have to distinguish again two cases : \ps  
     
{\it First subcase :} $\tau$ is an automorphism of algebra. Let $f \in S$   
be any idempotent lifting $\epsilon_1$.  Set $f' := f(1-\tau(f))$. Then $f' \tau(f') = f (1-\tau(f)) \tau(f) (1-f)=0$ and $\tau(f')f'=\tau(f)(1-f)f (1-\tau(f))=0$. Hence the subalgebra $S_1$ of $S$ generated by  
$f'$ and $\tau(f')$ is commutative and stable by $\tau$. Moreover, the reduction of $f'$ modulo $I_1:=I \cap S_1$ is $\epsilon_1(1-\tau(\epsilon_1))=\epsilon_1(1-\epsilon_2)=\epsilon_1$  
and the reduction of $\tau(f')$ is $\tau(\epsilon_1)=\epsilon_2$. \par  
Now, let $g$ be an idempotent in $S_1$ lifting $\epsilon_1$, and again let $g'=g(1-\tau(g))$.  
The same computation as above shows that $g' \tau(g')=\tau(g') g'=0$, but now, since $S_1$ is commutative, $g'$ is an idempotent. Set $e_1:=g'$, $e_2=:\tau(g')$, and the claim is proved, hence the lemma in this subcase (we could also have concluded by using the fact that the lemma is easy if $S$ is a finite commutative $A$-algebra). \ps  
  
{\it Second subcase :} $\tau$ is an anti-automorphism. Let $f \in S$   
be any idempotent lifting $\epsilon_1$. Set $x := f\tau(f)$. Then $x \in I$ and $\tau(x)=x$.  
Let $S_1$ be the $A$-subalgebra of $S$ generated by $x$, $I_1:=I\cap S_1$. This is a finite commutative $A$-algebra stable by $\tau$. Note that $I_1 \subset \rad(S_1)$. Indeed, $I_1 \subset \rad(S)$, hence for all $y\in I_1$, $1+y$ is invertible in $S$, hence in $S_1$ as it is integral over $A$. We conclude as $I_1$ is a two-sided ideal of $S_1$. In particular, $x \in \rad(S_1)$. Since $A$ is henselian and $2$ is invertible in $A$, there exists a unique element   
$u \in 1+\rad(S_1)$ such that $u^2=1-x$. Such an element $u$ is invertible in $S_1$ and satisfies $\tau(u)=u$.  
Set $g=u^{-1} f u$. Then $g$ is an idempotent lifting $\epsilon_1$ and from $u\tau(u)=u^2=1-f\tau(f)$ we get  
\begin{eqnarray*} g \tau(g) & =   
%& u^{-1} f u \tau(u) \tau(f) \tau(u^{-1}) \\ & = & u^{-1} f u^2 \tau(f) u^{-1} \\ & = & u^{-1} f (1-x) \tau(f) u^{-1} \\ & = &   
u^{-1} f(1-f\tau(f))\tau(f) u^{-1} = 0.\end{eqnarray*}  
Finally, we set $e_1=g-\frac{1}{2}\tau(g) g$ and $e_2=\tau(e_1)=\tau(g)-\frac{1}{2}\tau(g)g$. Then $e_i$ lifts $\epsilon_1$ and we claim that $e_i^2=e_i$  
 and $e_1e_2=e_2e_1=0$. Indeed, this follows at once from the following easy fact:\ps  
\noindent {\it Let $R$ be a ring in which $2$ is invertible, and let $e, f$ be two idempotents of $R$ such that $ef=0$. If we set $e'=(1-\frac{f}{2})e$ and $f'=f(1-\frac{e}{2})$, then $e'$ and $f'$ are orthogonal idempotents.}   
\end{pf}  
  
\begin{lemma} \label{idempotentspsi} Assume that $T$ is Cayley-Hamilton. There are idempotents $e_1,\dots,e_r$ in $R$ and morphisms $\psi:\ e_i R e_i   
\longrightarrow M_{d_i}(A)$ satisfying properties (1) to (5) of Lemma~\ref{idempotents}   
of prop~\ref{idempotents} and moreover   
\begin{itemize} \item[(6)] For $i \in\{1,\dots,r\},$ $\tau(e_i) = e_{\sigma(i)}$. \end{itemize}  
\end{lemma}   
\begin{pf}  
We call $\epsilon_i$, $i=1,\dots,r$ the central   
idempotents of $\bar R / \ker \bar T$. Note that we have $\tau(\epsilon_i)=\epsilon_{\sigma(i)}$. Applying the preceding lemma to $S=R$ and   
$I:=\ker(R \longrightarrow \bar R / \ker \bar T) = \rad R$ (Lemma~\ref{radical}), there exists a family of orthogonal idempotents    
$e_1,\dots, e_r$ lifting $\epsilon_1,\dots,\epsilon_r$ that is stable by $\tau$. Hence $\tau(e_i)=  
e_{\sigma(i)}$, which is (6), and the other properties are proved exactly as in Lemma~\ref{idempotents}.  
\end{pf}   
  
\subsubsection{Notations and choices}\label{notationsandchoices} 
  
From now we let $S$ be a Cayley-Hamilton quotient of $R$ which is stable by $\tau$. For example, by Remark \ref{quotsymker} the faithful quotient $R/\ker T$ has this property. As $\sigma^2=\Id$, we may cut $I$ into three parts $$I = I_0 \coprod I_1 \coprod I_2$$ with $i \in I_0$ if and only if $\sigma(i)=i$, and with $\sigma(I_1)=I_2$. If $i \in I_1$ and $j=\sigma(i)$, we definitely choose $\rhobar_j:=\rhobar_i^{\bot}$, which is permitted since $\rhobar_j$ and $\rhobar_i^\bot$ are isomorphic. \ps
  
We now choose in a specific way a GMA datum on $S$ taking into account 
the symmetry $\tau$. First, Lemma~\ref{idempotentspsi}
provides us with a family of idempotents $e_i$ such that $$\tau(e_i)=e_{\sigma(i)}, \forall i \in I.$$ 
Moreover, by property (5) 
(actually Lemma~\ref{idempotents} (5)) we also have 
isomorphisms $\psi_i : e_i S e_i \longrightarrow M_{d_i}(A)$ for $i \in I$. We are happy with the $\psi_i$ for $i \in I_0 \cup I_1$, but for $j \in I_2$, $j=\sigma(i)$ with $i \in I_1$  
we forget about the $\psi_j$ given by (5) by setting    
\begin{eqnarray}\label{I2} \psi_j = \psi_{\sigma(i)} :=   
\psi_i^\bot. \end{eqnarray}  
Of course, we also have $\psi_i=\psi_j^{\bot}$ as $\tau^2={\rm id}$.   
From now on, we fix a choice of $e_i$'s and $\psi_i$'s on $S$ as above, and this choice makes $S$ a GMA. \ps 
  
Let $i \in I_0$. Note that the two morphisms $\psi_i$ and   
$\psi_i^\bot :\ e_i S e_i \longrightarrow M_{d_i}(A)$ have the same   
trace and are residually irreducible. Hence by Serre and Carayol's result (\cite{carayol}),   
that is also the uniqueness part of the Nyssen and Rouquier's result, there exists a matrix   
$P_i \in \Gl_{d_i}(A)$ such that $\psi_i = P_i \psi_i^\bot P_i^{-1}$.  
Note that $P_i$ is determined up to the multiplication by an element of $A^\ast$. We fix the choices of such a matrix $P_i$ for each $i \in I_0$.  
% If $d_i=1$ there is a canonical choice : we may and do choose $P_i = 1$.  
  For $i \in I_1 \coprod I_2$ we set $P_i:= \Id$. Note that obviously $P_i = P_{\sigma(i)}$. We have, for any   
$i \in I$,  
\begin{eqnarray} \label{psiibot} \psi_{\sigma(i)} = P_i \psi_{i}^\bot P_i^{-1}, \ \ \  \psi_{i} = P_i \psi_{\sigma(i)}^\bot P_i^{-1}.  
\end{eqnarray}  
  
\begin{lemma} \label{pi2}  
If $\tau$ is an automorphism (resp. an   
anti-automorphism) $P_i^2$ (resp. $P_i {}^t P_i^{-1}$)  
is a scalar matrix $x_i \Id_{d_i}$ where $x_i \in A^\ast$   
(resp. $x_i \in \{\pm 1\}$).  
\end{lemma}  
\begin{pf} Assume that $\psi$ is an anti-automorphism (we leave the other, simpler, case to the reader). Using the   
two equalities of~(\ref{psiibot}) we get   
$$\psi_i = P_i {}^t P_i^{-1} \psi_i (P_i {}^t P_i^{-1})^{-1},$$ hence $P_i {}^t P_i^{-1}$ is a scalar matrix $x_i Id$ with $x_i \in A^\ast$ and  
we have  $x_i {}^t P_i = P_i$ hence   
$x_i^2=1$. The result follows since $A$ is local and $2$ is invertible in $A$.  
\end{pf}  

\subsubsection{Definition of the morphisms $\tau_{i,j}$}  
\label{deftauij}

Recall from~\S\ref{lesgma} that the idempotent $E_i$ of $S$  
 is defined as $\psi_i^{-1}(E_{1,1})$ and that $\cal{A}_{i,j}$ is the $A$-module $E_i S E_j$. Set   
$p_i = \psi_i^{-1}(P_i) \in e_i Se_i$. This is an invertible element in the algebra $e_i S e_i$ and we denote its inverse {\it in this algebra}   
by $p_i^{-1}$.  
  
Applying~(\ref{psiibot}) to $\tau(E_i)$ we get easily     
$$\tau(E_i)=p_{\sigma(i)} E_{\sigma(i)} p_{\sigma(i)}^{-1}$$  
  
Assume first that $\tau$ is an automorphism. We have  
$$ \tau(\cal{A}_{i,j})=  
\tau(E_i)\, S\, \tau(E_j) = p_{\sigma(i)} E_{\sigma(i)}   
p_{\sigma(i)}^{-1} \, S \, p_{\sigma(j)} E_{\sigma(j)} p_{\sigma(j)}^{-1}.$$  
Hence we may define a morphism of $A$-modules   
$\tau_{i.j} : \cal{A}_{i,j} \longrightarrow \cal{A}_{\sigma(i),\sigma(j)}$ by   
setting   
$$\tau_{i,j} = p_{\sigma(i)}^{-1} \tau_{|\cal{A}_{i,j}} p_{\sigma(j)}.$$  
  
Assume now that $\tau$ is an anti-automorphism. We define similarly a morphism  
$\tau_{i,j}: \cal{A}_{i,j} \longrightarrow \cal{A}_{\sigma(j),\sigma(i)}$  
 by   
setting   
$$\tau_{i,j} = p_{\sigma(j)}^{-1} \tau_{|\cal{A}_{i,j}} p_{\sigma(i)}.$$  
  
\begin{lemma} \label{lestauij} Assume $\tau$ is an automorphism (resp. an   
anti-automorphism).  
\begin{itemize}  
\item[(i)] For all $i,j$, the $A$-linear endomorphism $\tau_{\sigma(i),\sigma(j)}\circ \tau_{i,j}$ (resp. $\tau_{\sigma(j),\sigma(i)} \circ  
\tau_{i,j}$) of $\cal{A}_{i,j}$ is the multiplication by an element of $A^\ast$.  
\item[(ii)] For all $i,j$, $\tau_{i,j}$ is an isomorphism of $A$-modules.  
\item[(iii)] For all $i,j,k$ and $x\in \cal{A}_{i,j}$, $y \in \cal{A}_{j,k}$ we have   
$\tau_{i,j}(x) \tau_{j,k}(y)=\tau_{i,k}(xy)$ in $\cal{A}_{\sigma(i),\sigma(k)}$ (resp. $\tau_{j,k}(y) \tau_{i,j}(x)=\tau_{i,k}(xy)$ in $\cal{A}_{\sigma(k),\sigma(i)}$).   
\item[(iv)] We have $\tau_{i,j}(\cal{A}'_{i,j})=\cal{A}'_{\sigma(i),\sigma(j)}$ (resp.   
$\tau_{i,j}(\cal{A}'_{i,j})=\cal{A}'_{\sigma(j),\sigma(i)}$).
\end{itemize}  
\end{lemma}  
\begin{pf}   
The assertion (i) is an easy computation using   
Lemma~\ref{pi2}. The assertion (ii) follows immediately  
from (i). The assertion (iii) is a straightforward computation and (iv)   
follows  from (iii), (ii) and the definition of   
the $\cal{A}'_{i,j}$ (see \S\ref{extensions}).  
\end{pf}  
  
\subsubsection{Definition of the morphisms $\bot_{i,j}$}  
\label{defbotij}  
  
Let $\PP$ be a partition of $\{1,\dots,r\}$ such that the singletons $\{i\}$ and $\{j\}$ belong to $\PP$. Let   
$I_\PP$ be the corresponding reducibility ideal. Note that by Lemma~\ref{lestauij}, $I_\PP=I_{\sigma(\PP)}$ so that we may assume without changing $I_\PP$ that the singletons $\{\sigma(i)\}$ and $\{\sigma(j)\}$ belong   
to $\PP$. Let $J$ be an ideal of $A$ containing $I_\PP$.  
  
Recall that we defined a representation $\rho_{i} : R/JR \longrightarrow M_d(A/J)$ in Def. \ref{lesrhoi}. By point   
(ii) of Lemma~\ref{conseq},   
$\rho_i$ is the reduction mod $J$ of the composite of the morphism $\psi_{i}$ with the surjection   
$R \longrightarrow S \longrightarrow e_iSe_i$. Hence we have  
\begin{eqnarray} \label{rhoibot} \rho_{\sigma(i)} = P_i \rho_{i}^\bot P_i^{-1}. \end{eqnarray}  
  
Let $c$ be an extension in $\Ext^1_{R/JR}(\rho_j,\rho_i)$. We can see it as a morphism of algebra  
\begin{eqnarray*}  
\rho_c : R/JR & \longrightarrow & M_{d_i+d_j}(A/J) \\  
x  & \mapsto &  \left( \begin{array}{cc} \rho_i(x) & c(x) \\ 0 & \rho_j(x) \end{array} \right),   
\end{eqnarray*}  
 where $c(x) \in M_{d_i,d_j}(A/J)$.  
Then setting $Q_{i,j} = \diag(P_i,P_j) \in M_{d_i+d_j}(A/J)$ 
%if $\tau$ is an automorphism (resp.   $Q_{i,j} = \diag(P_j,P_i) \in M_{d_i+d_j}(A/J)$ if $\tau$ is an anti-automorphism)   
we see using~(\ref{rhoibot})   
that if $\tau$ is an automorphism,   
\begin{eqnarray}\ \ \ \ \ \  Q_{i,j} \rho_c^\bot (x)  Q_{i,j}^{-1} =    
\left( \begin{array}{cc} \rho_{\sigma(i)}(x) & c'(x) \\ 0 &   
\rho_{\sigma(j)}(x) \end{array} \right), \text{ where }\label{cprime}   
c'(x)= P_i \  c(\tau(x)) \    
P_j^{-1} , \end{eqnarray}  
and that if $\tau$ is an anti-automorphism,  
 \begin{eqnarray}\ \ \ \ \ \  Q_{i,j} \rho_c^\bot (x)  Q_{i,j}^{-1} =   
 \left( \begin{array}{cc} \rho_{\sigma(i)}(x) & 0 \\ c'(x) &   
\rho_{\sigma(j)}(x) \end{array} \right), \text{ where }\label{anticprime} c'(x)= P_j\ {}^t c(\tau(x))\ P_i^{-1}.  
\end{eqnarray}  
Hence $Q_{i,j} \rho_c^\bot  Q_{i,j}^{-1}$ represents an element $c'$  in   
$$\Ext^1_{R/JR} (\rho_{\sigma(j)},\rho_{\sigma(i)}) \text{ (resp. in }\Ext^1_{R/JR} (\rho_{\sigma(i)},\rho_{\sigma(j)}) \text{ )}$$  
\noindent and we set $$\bot_{i,j}(c):=c',$$ 
\noindent thus defining a morphism  
$$\bot_{i,j}:\ \Ext^1_{R/JR}(\rho_j,\rho_i) \longrightarrow \Ext^1_{R/JR}(\rho_{\sigma(j)},\rho_{\sigma(i)})$$  
$$\text{(resp. }   
\bot_{i,j}:\ \Ext^1_{R/JR}(\rho_j,\rho_i) \longrightarrow   
\Ext^1_{R/JR}(\rho_{\sigma(i)},\rho_{\sigma(j)}) \text{ )}.$$

Note that all we have done also works when $R$ is replaced by its $\tau$-stable Cayley-Hamilton quotient $S$,  
and that the morphisms $\bot_{i,j}$ thus defined on the $\Ext^1_{S/JS}$'s are   
simply the restriction of the morphisms $\bot_{i,j}$ on $\Ext^1_{R/JR}$.  
  
\subsubsection{The main result}  
  
\begin{prop} \label{mainsymmetry}  
  
 If $\tau$ is an automorphism, the following diagram is commutative  
$$  
\xymatrix{ \Hom_A(\cal{A}_{i,j}/\cal{A}'_{i,j},A/J) \ar[rrr]^{\iota_{i,j}} \ar[d]^{(\tau_{i,j}^{-1})^\ast} & & &  
\Ext^1_{S/JS}(\rho_j,\rho_i) \ar[d]^{\bot_{i,j}} \\     
\Hom_A(\cal{A}_{\sigma(i),\sigma(j)}/\cal{A}'_{\sigma(i),\sigma(j)},A/J)   
\ar[rrr]^{\iota_{\sigma(i),\sigma(j)}} & & & \Ext^1_{S/JS}(\rho_{\sigma(j)},\rho_{\sigma(i)})} 
$$  
 If $\tau$ is an anti-automorphism, the following diagram is commutative  
$$  
\xymatrix{ \Hom_A(\cal{A}_{i,j}/\cal{A}'_{i,j},A/J)   
\ar[rrr]^{\iota_{i,j}} \ar[d]^{(\tau_{i,j}^{-1})^\ast} & & &  
\Ext^1_{S/JS}(\rho_j,\rho_i) \ar[d]^{\bot_{i,j}} \\     
\Hom_A(\cal{A}_{\sigma(j),\sigma(i)}/\cal{A}'_{\sigma(j),\sigma(i)},A/J)   
\ar[rrr]^{\iota_{\sigma(j),\sigma(i)}} & & & \Ext^1_{S/JS}(\rho_{\sigma(i)},\rho_{\sigma(j)}) } 
$$  
\end{prop}  
\begin{pf}  
This follows immediately from the definitions of the morphisms $\tau_{i,j}$   
(see~\S\ref{deftauij}), $\bot_{i,j}$ (see~\S\ref{defbotij},   
especially~(\ref{cprime}) and~(\ref{anticprime})) and $\iota_{i,j}$'s
(see~\S\ref{extensions}).     
\end{pf}  
  
\subsubsection{A special case}  \label{speccasesect1}
  
We keep the assumptions of~\S\ref{setup} and the notations above,   
but we assume for simplicity that  
\begin{itemize}  
\item[(i)] the ring $A$ is reduced, of total fraction ring a finite product   
of fields $K=\prod_{s=1}^n K_s$,
\item[(ii)] the pseudocharacters $T \otimes K_s$ are irreducible,
\item[(iii)] $\tau$ is an anti-automorphism.
\end{itemize}  
  
Let $\rho: S:=R/\ker T \longrightarrow M_d(K)$ be a representation   
as in Theorem~\ref{structure} (ii). By assumption (ii) above and Theorem~\ref{structure} (iii), $\rho$ induces an isomorphism 
\begin{equation} \label{absirrutile} S\otimes_A K \isomo M_d(K). \end{equation}
For $s \in \{1,\dots,n\}$, denote by $\rho_s$ the composite   
$S \stackrel{\rho}{\longrightarrow} M_d(K) \longrightarrow M_d(K_s)$.  
  
\begin{lemma} For each $s \in \{1,\dots,n\}$ there exists a matrix $Q_s \in \Gl_d(K_s)$ such that   
\begin{equation} \label{Qs} \rho_s^\bot=Q_s \rho_s Q_s^{-1},\end{equation}   
and there is a   
well-determined sign $\epsilon_s =\pm 1$ such that ${}^t Q_s = \epsilon_s Q_s$. 
If $d$ is odd then $\epsilon_s=1$.  
\end{lemma}  
\begin{pf}      
The representations $\rho_s$ and $\rho_s^\bot$ are irreducible by hypothesis (ii) and have   
the same trace hence are isomorphic. Moreover $\rho_s$ is absolutely   
irreducible by (\ref{absirrutile}), hence the existence of a $Q_s$   
such that $\rho_s^\bot=Q_s \rho_s Q_s^{-1}$, and its uniqueness up to the   
multiplication by an element on $K_s^\ast$. Using that $(\rho_s^\bot)^\bot =
\rho_s$,  
we see that ${}^t Q_s Q_s^{-1}$ centralizes $\rho_s$, hence is a scalar matrix. 
Thus ${}^t Q_s = \epsilon_s Q_s$ and $\epsilon_s = \pm 1$. The last assertion holds because there is no antisymmetric invertible matrix   
in odd   
dimension.  
\end{pf}  

We will now relate these signs $\epsilon_s$ to other signs, and prove that they are actually equal in many cases.
Recall that if $k \in \{1,\dots,r\}$ is such that $\sigma(k)=k$, we fixed in \S\ref{notationsandchoices} 
a $P_k \in GL_{d_k}(A)$ such that $\psi_k=P_k\psi_k^{\bot}P_k^{-1}$, and we showed that 
${}^t P_k P_k^{-1} = \pm 1 \in A^*$ is a sign in Lemma \ref{pi2}. As explained there, $P_k$ is uniquely determined up to an element of $A^*$, so this sign 
is well defined, let us call it $\epsilon(k)$. By reducing those equalities mod $m$,
$\epsilon(k)$ is also "{\it the sign}" of the residual representation $\rhob_k\simeq \rhob_k^{\bot}$ in the obvious sense. 

\begin{lemma} Assume that $\sigma(k)=k$ for some $k \in \{1,\dots,r\}$. Then for each $s$, $\epsilon_s=\epsilon(k)$ is 
the sign of $\rhob_k$.
\end{lemma}

\begin{pf} As $\tau(e_k)=e_{\sigma(k)}=e_k$, we have $\tau(e_kSe_k)=e_kSe_k$. Recall that by 
the assumptions in \S\ref{notationsandchoices}, we have
$$\rho_{|e_kSe_k}=\psi_k: e_kRe_k \isomo M_{d_k}(A)$$
with $\psi_k^{\bot}=P_k^{-1}\psi_k P_k$. For each $s \in \{ 1,\cdots,n\}$, we also have 
$$e_k = \rho_s^{\bot}(e_k)= Q_s \rho_s(e_k) Q_s^{-1} = Q_s e_k Q_s^{-1},$$ 
so $Q_s$ commutes with $e_k={}^t e_k$, and $e_kQ_s$ and $Q_s$ are both symmetric or antisymmetric. 
Since $\psi_k: e_k(S\otimes_A K)e_k \isomo M_d(K)$ is an isomorphism, we get that for some $\lambda_s \in K_s^*$, 
	$$e_k Q_s = \lambda_s P_k^{-1}.$$
In particular, the three matrices $e_k Q_s$, $Q_s$ and $P_k$ (which does not depend on $s$) are simultaneously symmetric or antisymmetric, and we are done.
\end{pf}

Let us fix now $i\neq j$ two integers in $\{1,\dots,r\}$ such that $\sigma(i)=j$.    
Under hypothesis (iii), the morphism $\bot_{i,j}$ is an endomorphism   
of the $A$-module $\Ext^1_{S/JS}(\rho_j,\rho_i)$ and is canonically defined.  
We will study it using Proposition~\ref{mainsymmetry} and in terms of the signs above. Recall that we also defined some $A$-linear isomorphism 
$\tau_{i,j}$ of $\cal A_{i,j}=\cal A_{\sigma(j),\sigma(i)}$.

\begin{lemma} The morphism $\tau_{i,j} : \cal{A}_{i.j} \longrightarrow \cal{A}_{i,j}$ is the multiplication by the element  
$(\epsilon_1,\dots,\epsilon_n)$ of $K^*$.  
\end{lemma}  

\begin{pf} Let $Q \in \Gl_d(K)$ be the matrix whose image in $\Gl_d(K_s)$ is $Q_s$ for $s=1,\dots,n$.  
The representation $\rho$ identifies $S$ with a standard GMA $\rho(S)$  
in $M_d(K)$ and it follows from~(\ref{Qs}) that the anti-automorphism $\tau$ on $\rho(S)$ is the restriction of the anti-automorphism $M \mapsto Q{}^tMQ^{-1}$.  
Remember that $\rho(E_i)$ is the diagonal matrix whose all entries are zero   
but the $(d_1+\dots+d_{i-1}+1)^\text{th}$ which is one, and similarly for   
$\rho(E_j)$. Remember also that $\rho$ identifies $\cal{A}_{i,j}=E_i S E_j$ with $A_{i,j}=\rho(E_i) \rho(S) \rho(E_j)$. Since $\tau(E_i)=E_j$ we have   
$\rho(E_j) = Q {}^t \rho(E_i) Q^{-1}=Q \rho(E_i) Q^{-1}$. Thus the  $2$ by $2$  
submatrix of $Q$, keeping only the $(d_1+\dots+d_{i-1}+1)^\text{th}$ and $(d_1+\dots+d_{j-1}+1)^\text{th}$  lines and row, is antidiagonal :  
$$\left( \begin{array}{cc} \rho(E_i) Q \rho(E_i) & \rho(E_i) Q \rho(E_j) \\     
\rho(E_j) Q \rho(E_i) & \rho(E_j) Q \rho(E_j) \end{array} \right) = \left( \begin{array}{cc} 0 & a \\ b & 0 \end{array} \right) \in M_2(K)$$  
But by the lemma we have ${}^t Q = (\epsilon_1,\dots,\epsilon_n) Q$, hence   
$$b = (\epsilon_1,\dots,\epsilon_n) a.$$  
Now $\tau_{i,j} : A_{i,j} \longrightarrow A_{i,j}$   
is by definition the restriction of $M \mapsto  Q M Q^{-1}$ 
to $A_{i,j}=\rho(E_i) S \rho(E_j)$. By the formula above, this map is the multiplication by $a b^{-1}$, that is by the element $(\epsilon_1,\dots,\epsilon_n)$ of $K^*$.   
\end{pf}    
  
Thus, by Proposition~\ref{mainsymmetry} and the lemmas above :  

\begin{prop}\label{signebotijfinal}\begin{itemize} 
\item[(i)] If all the signs $\epsilon_s$ are equal, then for each pair 
$i \neq j$ with $j=\sigma(i)$ the endomorphism $$\bot_{i,j} :\Ext^1_{S/JS}(\rho_j,\rho_i) \longrightarrow \Ext^1_{S/JS}(\rho_j,\rho_i)$$ is the multiplication by $\epsilon_1=\pm 1$. 
\item[(ii)] If $\sigma$ has a fixed point $k$, then all the $\epsilon_s$ are equal to the sign of $\rhob_k \simeq \rhob_k^{\bot}$. 
\item[(iii)] If $d$ is odd, all these signs are $+1$.
\end{itemize}
\end{prop}  

\begin{remark}  
Note that the hypothesis of the corollary holds obviously when $A$ is a   
domain. Note also that the fact that $\bot_{i,j}$ is the multiplication   
by $\pm 1$ implies (and in fact is equivalent to)  
that every extension $\rho$ in  
$\Ext^1_{S/JS}(\rho_j,\rho_i)$ is isomorphic to $\rho^\bot$ {\it as a representation} (not necessarily as an extension).  
\end{remark}

\newpage

\section{Trianguline deformations of refined crystalline representations}\label{trianguline}

\subsection{Introduction} 

The aim of this section is twofold. First, we study the $d$-dimensional trianguline
representations of $$\Gp:=\Gal(\overline{\Q}_p/\Q_p)$$\index{Gp@$G_p$, the absolute Galois group of $\Q_p$}for any $d\geq 1$ and with artinian ring coefficients,
extending some results of Colmez in \cite{colmeztri}.
Then, we use them to define and study some deformation problems of the
$d$-dimensional crystalline representations of $\Gp$.  \ps

These deformation problems are motivated by the theory of $p$-adic families of automorphic forms and the wish to understand the family of Galois representations carried by the eigenvarieties. They have been extensively studied in the special case of
 ordinary deformations ({\it e.g.} Hida families), however the general case is more subtle. When $d=2$, it was first dealt with by Kisin in \cite{kis}. He proved that the local $p$-adic Galois representation attached to any finite slope overconvergent modular 
eigenform $f$ admits a non trivial crystalline period on which the crystalline 
Frobenius acts through $a_p$ if $U_p(f)=a_pf$, and also that this period "vary analytically" on the eigencurve. These facts lead him to define and study some 
deformation problem he called $D^h$ in {\it loc. cit.} \S8. In favorable cases, he was then able to show that the Galois deformations coming from Coleman's families give examples of such ``$h$-deformations'' (see \S10, 11 {\it loc. cit}). In this section, we define and study a deformation problem for the $d$-dimensional case via the theory of $\fg$-modules. It turns out to be isomorphic to Kisin's one when $d=2$ but in a non trivial way. We postpone to sections \ref{kisin} and \ref{families} the question of showing that higher rank eigenvarieties produce such deformations. \ps
	The approach we follow to define these problems is mainly suggested by 
Colmez's interpretation of the first result of Kisin recalled above in \cite{colmeztri}. Precisely, Colmez proves that for a $2$-dimensional $p$-adic 
representation $V$ of $\Gp$, a twist of $V$ admits a non trivial 
crystalline period if, and only if, the $\fg$-module of $V$ over the Robba ring\footnote{Recall that the category of $\fg$-modules over the Robba ring $\Ro$ is strictly bigger than the category of $\Q_p$-representations of $\Gp$, which occurs as its full subcategory of \'etale objects.} is triangulable (\cite[Prop. 5.3]{colmeztri}). For instance, the $\fg$-module of a $2$-dimensional crystalline representation is always trigonalisable  (with non \'etale graded pieces in general) even if the representation is irreducible (that is non ordinary), which makes things interesting. This also led Colmez to define a {\it trianguline representation} as a representation whose $\fg$-module over $\Ro$ is a successive extension of rank $1$ $\fg$-modules. Although this has not yet been proved, it is believed (and suggested by Kisin's work) that the above triangulation should vary analytically on the eigencurve, so that the general ``finite slope families'' should look pretty much like ordinary families from this point of view\footnote{A related question is to describe the $A$-valued points, $A$ being any $\Q_p$-affinoid algebra, of the parameter space $\cal{S}$ of triangular $\fg$-modules defined by Colmez in \cite[\S0.2]{colmeztri}. 
The material of this part would be e.g. enough to answer the case where $A$ is an artinian $\Q_p$-algebra, at least for "non critical" triangulations. See also our results in section \ref{families}.}. \ps

In what follows, we define and study in details the {\it trianguline
deformation functors} of a given $d$-dimensional crystalline
representation for any $d$, establishing an ``infinitesimal version''
of the program above, that is working with artinian $\Q_p$-algebras as
coefficients (instead of general $\Q_p$-affinoids which would require
extra work). This case will be enough for the applications in the next
sections and contains already quite a number of subtleties, mainly
related to the notion ``non criticality''. We prove also a number of
results of independent interest on triangular $\fg$-modules, some
of them generalizing to the $d$-dimensional case some results of
Colmez in \cite{colmeztri}. Let us describe now more precisely what we
show. \ps

In \S\ref{hodgedef}, we collect the fundamental facts we shall use of the theory of $\fg$-modules over the Robba ring $\Ro$. 
We deduce from Kedlaya's theorem that an extension between two \'etale $\fg$-modules is itself \'etale (lemma \ref{etalekedlaya}). A useful corollary is the fact 
that it is the same to deform the $\fg$-module over $\Ro$ of a representation or to deform the representation itself (Proposition \ref{XVXD}). We prove 
also in this part some useful results on modules over the Robba ring with coefficients in an artinian $\Qp$-algebra. \ps

In \S\ref{triangularfg}, we study the triangular $\fg$-modules over $\Ro_A:=\Ro \otimes_{\Qp} A$ where $A$ is an artinian $\Qp$-algebra. They are defined as $\fg$-modules $D$ finite free over $\Ro_A$ equipped with a strictly increasing filtration (a {\it triangulation}) $$(\Fil_i(D))_{i=0,\dots,d}, \, \, \,  d:=\rk_{\Ro_A}(D),$$ of $\fg$-submodules which are free and direct summand over $\Ro_A$\footnote{It is important here 
not to restrict to the \'etale $D$, even if in some important applications this would be the case. Indeed, most of the proofs use an induction on $d$ and the $\Fil_i(D) \subset D$ will not even be isocline in general.}. When $D$ has rank $1$ over $\Ro_A$, we show that it is isomorphic to a ``basic'' one $\Ro_A(\delta)$ for some unique continuous character $\W \longrightarrow A^*$ (Proposition 
\ref{uni}), hence the graded pieces of $\Fil_i(D)/\Fil_{i-1}(D)$ have the form $\Ro_A(\delta_i)$ in general. The {\it parameter} $(\delta_i)_{i=1,\dots,d}$ of $D$ 
defined this way turns out to refine the data of the Sen polynomial of $D$ 
(Proposition \ref{senpol}). A first important result of this part is a {\it weight criterion} ensuring that such a $\fg$-module is de Rham (Proposition \ref{criteriumDR}); this criterion appears to be a generalization 
to trianguline representations of Perrin-Riou's criterion ``ordinary representations are semistable''(\cite[Expos\'e IV, Th\'eor\`eme]{PP}). In the last paragraphs, we define and study the functor of triangular deformations of a given triangular $\fg$-module $D_0$ over $\Ro$: its $A$-points are simply the triangular $\fg$-modules deforming $D_0$ and whose triangulation lifts the fixed triangulation of $D_0$. In the same vein, a {\it trianguline deformation} of a trianguline representation $V_0$ is a triangular deformation of its $\fg$-module\footnote{In this section, all the $\fg$-modules are understood with coefficients in the Robba ring $\Ro$.} $D_0$ (it depends on the triangulation of $D_0$ we choose in general). The main result here is a complete description of these functors under some explicit conditions on the parameter of the triangulation of $D_0$ (Proposition \ref{deform}). \ps

In \S\ref{refinement}, we show that crystalline representations are trianguline and 
study the different possible triangulations of the $\fg$-module of a given crystalline representation\footnote{For simplicity, we restrict there to crystalline representations with distinct Hodge-Tate weights. 
In fact, the results of this part could be extended to the representations becoming semi-stable over an abelian extension of $\Q_p$, and even to all the de Rham representations in a weaker form.} $V$. We show that they are in natural bijection with the {\it refinements} of $V$ in Mazur's sense \cite{Maz}, that is
the full $\varphi$-stable filtrations of $\Dc(V)$. 
More importantly, we introduce a notion of {\it non critical refinement} 
in \S\ref{crit} by asking that the 
$\varphi$-stable filtration is in general position compared to the weight filtration on $\Dc(V)$. 
We interpret then this condition in terms 
of the associated triangulation of the $\fg$-module (Proposition \ref{translation}), and compare it to other related definitions in the literature (remark \ref{remcrit}). This notion turns out to be the central one in all the subsequent results. 
The main ingredient for this part is Berger's paper \cite{bergerinv}. \ps

In \S\ref{defdef}, we apply all the previous parts to define and
study the trianguline deformation functor of a refined crystalline
representation. It should be understood as follows: the choice a
refinement of $V$ defines a triangulation of its $\fg$-module by the
previous part, and we can study the associated trianguline deformation
problem defined above. When the chosen refinement is 
non critical, we can explicitly describe the trianguline deformation functor (Theorem \ref{thmdefnoncrit}), and also describe the crystalline locus inside it. 
A maybe striking result is that ``a trianguline deformation of a non
critically refined crystalline representation is crystalline if and
only if it is Hodge-Tate'' (Theorem \ref{colcrit}). This fact may be
viewed as an infinitesimal local version of Coleman's ``small slope
forms are classical'' result; it will play an important role in the
applications to Selmer groups of the subsequent 
sections (see e.g. Corollary \ref{corredloc}). In the last paragraph, 
we give a criterion ensuring that a deformation satisfying some conditions in Kisin's style is in fact trianguline (Theorem \ref{poidscst2}). Combined with 
the extensions of Kisin's work studied in section \ref{kisin}, this result will be useful to prove that the Galois deformations living on eigenvarieties are trianguline in many interesting cases\footnote{However not in all cases; part of this result may be viewed as a trick allowing to circumvent the study of a theory of families of triangular $\fg$-modules alluded above.}. \ps

In a last \S\ref{globcons}, we discuss some applications of these
results to global deformation problems. Recall that a consequence of
the Bloch-Kato conjecture for adjoint pure motives is that a
geometric, irreducible, $p$-adic Galois representation $V$ (say crystalline above $p$) 
admits no non trivial crystalline deformation\footnote{As $\Q$-motives are
countable, it is certainly expected that there is non trivial
$1$-parameter $p$-adic family of motives, but the infinitesimal
assertion is stronger.}. Admitting this, we obtain that the trianguline
deformations of $V$ for a non critical refinement $\Ref$ 
(and with good reduction outside $p$ say) have Krull-dimension at most $\dim(V)$
(corollary \ref{selmerglob}). This ``explains'' for example why the eigenvarieties of reductive rank $d$
have dimension at most $d$, and in general it relates  
the dimension of the tangent space of eigenvarieties of $\GL(n)$
at classical points (about which we know very few) to an
``explicit'' Selmer group. As another good indication about the
relevance of the objects above, let us just say that when such a
$(V,\Ref)$ 
appears as a classical point $x$ on a unitary eigenvariety $X$ (say of
``minimal level outside $p$''), 
standard conjectures imply that $$R \simeq T \overset{\kappa}{\isomo}
L[[X_1,...,X_d]],$$ where $R$ prorepresents the trianguline deformation
functor of $(V,\Ref)$, $T$ is the completion of $X$ at $x$, and
$\kappa$ is the morphism of the eigenvariety to the weight space. \ps \ps

The authors are grateful to Laurent Berger and Pierre Colmez for very helpful
discussions during the preparation of this section. We started working on the infinitesimal properties of the Galois representations
on eigenvarieties in September 2003, and since we have been
faced to an increasing number of questions concerning non de Rham
$p$-adic representations which were fundamental regarding the
arithmetic applications. We
warmly thank them for taking the time to think about our questions
during this whole period. As will be clear to the reader, Colmez's paper
\cite{colmeztri} has been especially influential to us.

\subsection{Preliminaries of $p$-adic Hodge theory and $\fg$-modules} \label{hodgedef}
\subsubsection{Notations and conventions.} In all this section, 
$$\Gp=\Gal(\Qpb/\Q_p),$$ 
is equipped with its Krull topology. 
\index{Aa@$A$, a local artinian ring, finite over $\Q_p$|(} Let $A$ be a finite dimensional local
$\Q_p$-algebra equipped with its unique Banach $\Q_p$-algebra
topology, $m$ its maximal ideal, $L:=A/m$. 
\index{L@$L$, the residue field of the artinian local ring $A$, a finite extension of $\Q_p$|(}

By an {\it $A$-representation of $\Gp$}, 
we shall always mean an $A$-linear, continuous, representation of $\Gp$ on a finite 
type $A$-module. We fix an algebraic closure $\Qpb$ of $\Qp$, 
equipped with its canonical valuation $v$, and norm $|.|$, 
extending the one of $\Q_p$ (so $v(p)=1$ and $|p|=1/p$), and we denote by $\C_p$ its completion. We denote by $\Bc$, $\Bdr$, $\Dc(-)$, $\Ddr(-)$ etc... the
usual rings and functors defined by Fontaine (\cite[Expos\'es II et III]{PP}).\ps
We denote by $\Q_p(1)$ \index{Qp@$\Q_p(1)$} 
the $\Q_p$-representation of $G_p$ on $\Q_p$ defined by the cyclotomic character 
$$\chi: \Gp
\longrightarrow \Z_p^*.$$ \index{Zc@$\Gamma$, a quotient of $G_p$ 
isomorphic to $\Z_p^*$} If we set $\Hp:=\Ker(\chi)$ and 
$\Gamma:=\Gp/\Hp$, then $\chi$ induces a canonical isomorphism \index{Zk@$\chi$, the cyclotomic character, seen as a a character of $G_p$ or of $\Gamma$} $\Gamma
\isomo \Zp^*$. Our convention on the sign of the Hodge-Tate weights, and on the Sen
polynomial, is that $\Q_p(1)$ has weight $-1$ and Sen polynomial $T+1$. With
this convention, the Hodge-Tate weights (without multiplicities) of a de Rham representation $V$
are the jumps of the weight filtration on $\Ddr(V)$, that is the integers $i$ such that
$\Fil^{i+1}(\Ddr(V)) \subsetneq \Fil^i(\Ddr(V))$, and also the roots of the
Sen polynomial of $V$.
\subsubsection{$\fg$-modules over the Robba ring $\Ro_A$.} It will be
convenient for us to adopt the point of view of $\fg$-modules over the Robba
ring, for which we refer to \cite{fontainegr}, \cite{colmezbki}, 
\cite{kedlaya}, and \cite{bergereqdiff}. \ps
Let $\Ro_A$ be the Robba ring with
coefficients in $A$, {\it i.e.} the ring of power series $$f(z)=\sum_{n\in \Z} a_n
(z-1)^n,\, \, \, a_n \in A\,\,$$ 
converging on some annulus of $\C_p$ of the form $r(f)\leq |z-1|<1$, equipped with
 its natural $A$-algebra topology. If we set $$\Ro:=\Ro_{\Q_p},$$ we
have $\Ro_A=\Ro\otimes_{\Q_p}A$. Recall that $\Ro_A$ is equipped with
commuting, $A$-linear, continuous actions of $\varphi$ and
$\Gamma$ defined by 
$$\varphi(f)(z)=f(z^p),\, \, \, \,
\gamma(f)(z)=f(z^{\chi(\gamma)}).$$
To get a picture of this action, note that if $z \in \C_p$ 
satisfies $|z-1|<1$, we have $|z^n-1|=|z-1|$ for $n \in \Z_p^*$, 
whereas $|z^p-1|=|z-1|^p$ when $|z-1|>p^{-\frac{1}{p-1}}$. 

\begin{definition} A $\fg$-module over
$\Ro_A$ is
a finitely generated $\Ro_A$-module $D$ which is free over $\Ro$ and equipped with commuting, $\Ro_A$-semilinear,
continuous\footnote{It means that for any choice of a free basis
$e=(e_i)_{i=1...d}$ of $D$ as $\Ro$-module, the matrix map $\gamma \mapsto M_e(\gamma) \in \GL_d(\Ro)$,
defined by $\gamma(e_i)=M_e(\gamma)(e_i)$, is a continuous function on
$\Gamma$. If $P\in \GL_d(\Ro)$, then $M_{P(e)}(\gamma)=\gamma(P)M_e(\gamma)P^{-1}$,
hence it suffices to check it for a single basis.} 
actions of $\varphi$ and $\Gamma$, and such that $\Ro\varphi(D)=D$.
\end{definition}

\subsubsection{Some algebraic properties of $\Ro_A$.} \label{algpropr}

In the first part of this section, we assume that $A=L$ is a field, and we will now recall some algebraic properties of modules over
$\Ro_L$. \ps 

A first remark is that $\Ro_L$ is a domain. Moreover, although it is not noetherian, a key property is that $\Ro_L$ is an {\it adequate Bezout domain} 
(see \cite[prop. 4.12]{bergerinv}), hence the theory of finitely presented $\Ro_L$-modules is similar to the one for principal rings: \begin{itemize}
\item[(B1)] Finitely generated, torsion free, $\Ro_L$-modules are free.
\item[(B2)] For any finite type $\Ro_L$-submodule $M \subset \Ro_L^n$, there is a basis
$(e_i)$ of $\Ro_L^n$, and elements $(f_i)_{1\leq i \leq d} \in
(\Ro_L \backslash\{0\})^d$, such that $M=\oplus_{i=1}^d f_i\Ro_L e_i$. The $f_i$
may be chosen such that $f_i$ divides $f_{i+1}$ in $\Ro_L$ for $1\leq i \leq d-1$, and are unique
up to units of $\Ro_L$ if this is satisfied (they are called the {\it elementary divisors}
of $M$ in $\Ro_L^n$). \ps 
\end{itemize}

Let $M \subset \Ro_L^n$ be a $\Ro_L$-submodule, the {\it saturation} 
of $M$ in $\Ro_L^n$ is
 $$M^{\rm sat}:=M \cap \fr(\Ro_L)^n,$$
and we say that $M$ is {\it saturated} if $M^{\rm sat}=M$, or which is the same if $\Ro_L^n/M$ is torsion free\footnote{As 
$\fr{\Ro_L}=\fr{\Ro}\otimes_{\Q_p}L$, an $\Ro_L$-module 
is torsion free over $\Ro_L$ if, and only if, it is torsion free over
$\Ro$.}. By (B1) (resp. (B2)) such an $M$ is saturated if, and only if, it is a direct summand as $\Ro_L$-module (resp. if its elementary divisors are units).
Note also that by property (B2), if $M \subset \Ro_L^n$ is finite type over
$\Ro_L$, then so is $M^{\rm sat}$.
\ps
	It turns out that in a $\fg$-module situation, we can say much more. 
Let $t:=\log(z) \in \Ro$ be the usual " $2i\pi$-element". It satisfies $\varphi(t)=pt$
and $\gamma(t)=\chi(\gamma)t$ for all $\gamma \in \Gamma$. Note that $t$ is not an irreducible element of $\Ro$, as it is divisible by $z^{p^n}-1$ for all $n$.

\begin{prop} \label{structurer} Let $D$ be a $\fg$-module over $\Ro_L$ and $D'$ a $\fg$-submodule.
\begin{itemize}
\item[(i)] ${D'}^{\rm sat}=D'[1/t]\cap D$.
\item[(ii)] If $D'$ has rank $1$ over $\Ro_L$, then $D'=t^k {D'}^{\rm sat}$, $k\in \N$.
\end{itemize}
\end{prop}

\begin{pf} Part (ii) is \cite[rem. 4.4]{colmeztri}. To prove part (i), it suffices to show that the product of the elementary divisors of $D'$ is a power of $t$. But this follows from (ii) applied to $\Lambda^j(D') \subset \Lambda^j(D)$ with $j=\rk_{\Ro_L}(D')$.
\end{pf}

We end this section by establishing some basic but useful facts when we work with artinian rings; in what follows $A$ is not supposed to be a field any more.

\begin{lemma}\label{lib} \begin{itemize}
\item[(i)] Let $E$ be a free $A$-module and $E' \subset E$ a free submodule, then $E'$ is a direct summand of $E$.
\item[(ii)] Let $E$ be a $\Ro_A$-module (resp. $\Ro_A[1/t]$-module) which is free of finite type as
$\Ro$-module (resp. $\Ro[1/t]$-module), and free as $A$-module. Then $E$ is free of finite type over
$\Ro_A$ (resp. $\Ro_A[1/t]$).
\item[(iii)] Let $D$ be a finite free $\Ro_A$-module. 
Assume that $D$ contains a submodule $D'$
free of rank $1$ such that $D'/mD'$ is saturated in
$D/mD$ as $\Ro$-module. Then $D'$ is a direct summand as $\Ro_A$-submodule of $D$.
\end{itemize}
\end{lemma}
\begin{pf} Let $n\geq 1$ denote the smallest integer such that $m^n=0$. As $m$ is nilpotent, the following version of Nakayama's lemma holds for all $A$-modules $F$: $F$ is zero (resp. free) if, and only if, $F/mF=0$ (resp.
$\Tor_1^A(F,A/m)=0$).\par
To prove (i), consider the natural exact sequence $$0 \longrightarrow \Tor_1^A(E/E',A/m) \longrightarrow E'/mE' \longrightarrow E/mE.$$
We have to show that the $\Tor$ above is zero, {\it i.e.} that $mE\cap E'=mE'$. But this follows from the fact that for a free $A$-module $F$, $$mF=F[m^{n-1}]:=\{f \in F, m^{n-1}f=0\}.$$
Let us show (ii) now, set $\Ro'=\Ro[1/t]$ or $\Ro$. As $E$ is free over $A$ we have $mE=E[m^{n-1}]$  hence $mE$ is a saturated $\Ro'$-submodule of $E$. As a consequence, $E/mE$ is a torsion free, finite type, $\Ro'_L=\Ro'_A/m\Ro'_A$-module, so it is free over $\Ro'_L$ by property (B1). As $E$ is free over $A$, Nakayama's lemma shows that any $\Ro'_A$-lift $(\Ro'_A)^d \longrightarrow E$ of an $\Ro'_L$-isomorphism $(\Ro'_L)^d \isomo E/mE$ is itself an
isomorphism.\ps
 Before we prove (iii), let us do the following remark. 

\begin{remark} \label{remfree} Let $D$ be a $\fg$-module free over $\Ro_A$ and $D'
\subset D$ be a submodule also free over $\Ro_A$. Then part (i)
shows that $D'$ is a direct summand of $D$ as $A$-module. In particular,
for any ideal $I \subset A$, the natural map $D'/ID' \longrightarrow
D/ID$
is injective, and $D' \cap ID = ID'$. Note that
this remark gives sense to part (iii) of the proposition (take $I=m$).
\end{remark}

To prove (iii), let us argue by induction on the length of $A$ to show that
$D'$ is $\Ro$-saturated in $D$. We are done by assumption if $A$ is a field. Let $I \subset A$ be a proper ideal,
then $D'/ID' \subset D/ID$ satisfies the induction hypothesis by the
remark \ref{remfree} above,
hence $D'/ID'$ is a saturated $\Ro$-submodule of $D/ID$.
As $ID$ is a direct summand as $\Ro$-module, it is saturated in $D$,
hence we only have to show that $D' \cap ID=ID'$ is saturated in $ID$.
We may and do choose an ideal $I$ of length $1$.
But in this case, $ID' \subset ID$ is $D'/mD' \subset
D/mD$, which is saturated by assumption. We proved that $D'$ is
$\Ro$ saturated in $D$. As a consequence, $D/D'$ is $\Ro$-torsion free,
and free over $A$ by part (i), hence it is free over $\Ro_A$ by
part (ii). \ps
\end{pf}

\subsubsection{\'Etale and isocline $\varphi$-modules.} \label{pentes} Assume again that $A=L$ is a field. Let $D$ be a $\varphi$-module over $\Ro_L$, i.e. a free of finite rank  
$\Ro_L$-module with a $\Ro_L$-semilinear action of $\varphi$ 
such that $\Ro_L\varphi(D)=D$. 
Recall that Kedlaya's work (see
\cite[Theorem 6.10]{kedlaya}) associates to $D$ a 
sequence of rational numbers $s_1 \leq \dots \leq s_d$ (where $d$ is the rank of $D$) called the {\it slopes} of $D$. The $\varphi$-module $D$ is said 
{\it isocline} of slope $s$ if $s_1=\dots=s_d=s$ and {\it \'etale} if it is 
isocline of slope $0$. A $\fg$-module is \'etale (resp. isocline) if its 
underlying $\varphi$-module (forgetting the action of $\Gamma$) is. For more details see~\cite{kedlaya}, especially part 4 and 6, \cite{kedlaya2} or, for a concise review, \cite{colmeztri}, part 2.  \ps
 
\begin{lemma} \label{etalekedlaya}
Let $0 \longrightarrow D_1 \longrightarrow D \longrightarrow D_2
\longrightarrow 0$ be an exact sequence of $\varphi$-modules free of finite rank
over $\Ro_L$. If $D_1$ and $D_2$ are isocline of the same slope $s$, then $D$
is also isocline of slope $s$. 
\end{lemma}

\begin{pf}
Up to a twist (after enlarging $L$ if necessary) we may
assume that $s=0$, that is $D_1$ and $D_2$ \'etale, 
and we have to prove that $D$ is \'etale as well.

Assume that $D$ is not \'etale, so it has by Kedlaya's Theorem 
(\cite[Thm 6.10]{kedlaya}) a saturated $\varphi$-submodule $N$ which 
is isocline of slope $s < 0$. Note $n$ the rank of $N$ and consider the $\varphi$-module $\Lambda^n D$
which contains as a saturated $\varphi$-submodule the $\varphi$-module
$\Lambda^n N$, of slope $ns <0$. 
By assumption, $\Lambda^n M$ is a successive extensions of $\varphi$-modules of 
the form $$\Lambda^a(D_1)\otimes_L\Lambda^b(D_2),\, \, a+b=n,$$ 
which are all \'etale (see~\cite[Prop. 5.13]{kedlaya}). Since $\Lambda^n N$ has rank one, it is isomorphic to 
a submodule of one of those \'etale $\varphi$-modules 
$\Lambda^a(D_1)\otimes_L \Lambda^b(D_2)$. But by~\cite[Prop. 4.5.14]{kedlaya2}, an \'etale $\varphi$-module has no rank one $\varphi$-submodule of slope $<0$, a contradiction.
\end{pf}

\subsubsection{$\fg$-modules and representations of $\Gp$.} 

Works of Fontaine, Cherbonnier-Colmez, and Kedlaya, 
allow to define a $\otimes$-equivalence $\Dr$ between 
the category of $\Qp$-representations of $\Gp$ and 
\'etale (in the sense of \S\ref{pentes}) $\fg$-modules over $\Ro$ (\cite[prop. 2.7]{colmeztri}). 
By \cite[\S3.4]{bergerinv}, $\Dr(V)$ can be defined in Fontaine's style: there
exists a topological ring $B$ (denoted $B^{\dag,\rig}$ there) 
equipped with actions of $\Gp$ and $\varphi$ and such that $B^\Hp=\Ro$, and
$$\Dr(V):=(V\otimes_{\Qp} B)^\Hp.$$

As in Herr's theory, recall that there is an explicit complex of $\Ro$-modules computing the
cohomology of a $\fg$-module $D$ (see \cite[\S3.1]{colmeztri}), which we will denote by
$H^i(D)$. Some properties of these constructions are summarized in the following proposition.

\begin{prop} \label{drig}
\begin{itemize}
\item[(i)] The functor $\Dr$ induces an $\otimes$-equivalence of categories
between $A$-representations of $\Gp$ and \'etale $\fg$-modules over
$\Ro_A$. We have $\rk_{\Qp}(V)=\rk_{\Ro}(\Dr(V))$.
\item[(ii)] For an $A$-representation $V$ of $\Gp$, $\Dr(-)$ induces an isomorphism
$$\Ext_{A[\Gp],{\rm cont}}(A,V) \isomo \Ext_{\fg-{\rm mod}}(\Ro_A,\Dr(V))=H^1(\Dr(V)).$$
\end{itemize}
\end{prop}

\begin{pf} Part (i) is \cite[prop. 2.7]{colmeztri}, part (ii) follows from (i)
and Proposition \ref{etalekedlaya}. 
\end{pf}

\begin{lemma} \label{freeness} An $A$-representation $V$ of $\Gp$ is free over $A$ if, and only if,
$\Dr(V)$ is free over $\Ro_A$.
\end{lemma}

\begin{pf} Assume first that $V$ is free over $A$. Let $M$ be any finite length $A$-module $M$, and fix a
presentation $A^n \longrightarrow A^m \longrightarrow M \longrightarrow 0.$
As the functor $\Dr(-)$ is exact by Proposition \ref{drig} (i), and by left exactness of
$-\otimes_A \Dr(V)$, we deduce from this presentation that we have a
canonical $A$-linear isomorphism 
$$\Dr(V)\otimes_A M \isomo \Dr(V \otimes_A M).$$
In particular, when $M=A/m$ where $m$ is the maximal ideal of $A$, we get
that the minimal number of generators of $D(V)$ as $\Ro_A$-module is $d:=\rk_A(V)$. 
In particular, there is a $\Ro_A$-linear surjection $\Ro_A^d \longrightarrow \Dr(V)$. As $$\rk_{\Ro}D(V)=\dim_{\Qp}(V)=d\dim_{\Qp}(A)=\rk_{\Ro}(\Ro_A^d)$$ by 
Proposition \ref{drig} (i), any such surjection is an isomorphism by property (B1).
\ps
The proof is the same in the other direction using the natural inverse functor of
$\Dr$. 
\end{pf}

\subsubsection{Berger's theorem.} We will need to recover the usual Fontaine 
functors from $\Dr(V)$, which is achieved by Berger's work \cite{bergerinv} and
\cite{bergereqdiff} that we recall now. Let us introduce, for $r > 0 \in \Q$, the $\Q_p$-subalgebra
$$\Ro_r=\{f(z) \in \Ro, \, \, f{\rm \, \, converges\, \, on\, \,  
the\, \, annulus\, \, }
p^{-\frac{1}{r}} \leq |z-1| < 1\}.$$ 
Note that $\Ro_r$ is stable by $\Gamma$, and that $\varphi$ induces a map
 $\Ro_r \longrightarrow 
\Ro_{pr}$ when $r> \frac{p-1}{p}$ which is \'etale of degree $p$. 
The following lemma is \cite[thm
1.3.3]{bergereqdiff}: \ps

\begin{lemma} \label{defdr} Let $D$ be a $\fg$-module over
$\Ro$. There exists a $r(D)>\frac{p-1}{p}$ such that for each $r> r(D)$, 
there exists a unique finite free, $\Gamma$-stable,  
$\Ro_r$-submodule $D_r$ of $D$ such that
%Recall that for any commutative ring $C$, a $C$-linear surjection $C^n \longrightarrow C^n$ is bijective. In particular, it is equivalent here to ask that $\Ro D_r=D$ and $\rk_{\Ro_r}(D_r)=\rk_{\Ro}(D)$ or to ask that the natural map $\Ro \otimes_{\Ro_r} D_r \longrightarrow D$ is an isomorphism.} 
$\Ro \otimes_{\Ro_r} D_r \isomo D$ and that $\Ro_{pr}D_r$ has a $\Ro_{pr}$-basis in 
$\varphi(D_r)$. In particular, for $r>r(D)$, 
\begin{itemize}
\item[(i)] for $s\geq r$, $D_s=\Ro_sD_r \isomo \Ro_s \otimes_{\Ro_r} D_r$,
\item[(ii)] $\varphi$ induces an isomorphism 
$\Ro_{pr} \otimes_{\Ro_r, \varphi} D_r \isomo D_{pr} \isomo \Ro_{pr}\otimes_{\Ro_r}D_r$.
\end{itemize}
\end{lemma} 

\ps

If $n(r)$ is the smallest integer $n$ such that 
$p^{n-1}(p-1) \geq r$, then for $n\geq n(r)$ the primitive $p^n$-th roots of 
unity lie in the annuli $p^{-\frac{1}{r}} < |z-1| < 1$ and 
$t$ is a uniformizer at each of them so that we get 
by localization and completion at their underlying 
closed point a natural map $$\Ro_r \longrightarrow
K_n[[t]], \, \, n\geq n(r),\, \, r>r(D), $$
which is injective with $t$-adically dense image, 
where $K_n:=\Q_p(\sqrt[p^n]{1}\,)$. For any $\fg$-module over $\Ro$, we can
then form for $r>r(D)$ and $n\geq n(r)$ the space
$$D_r \otimes_{\Ro_r} K_n[[t]],$$
which is a $K_n[[t]]$-module free of rank $\rk_{\Ro}(D)$ 
equipped with a semi-linear continuous action of $\Gamma$. By 
Lemma \ref{defdr} (i), this space does not depend on the choice of $r$ such
that $n \geq n(r)$. Moreover, for a fixed $r$, $\varphi$ induces by the same
lemma part (ii)\, $-\,\otimes_{\Ro_{pr}}K_{n+1}[[t]]$ a
$\Gamma$-equivariant, $K_{n+1}[[t]]$-linear, isomorphism
$$(D_r \otimes_{\Ro_r} K_n[[t]])\otimes_{t\mapsto pt} K_{n+1}[[t]]
\longrightarrow D_r \otimes_{\Ro_r} K_{n+1}[[t]].$$
(Note that the map $\varphi: \Ro_r \longrightarrow \Ro_{pr}$ induces
the inclusion $K_n[[t]] \longrightarrow K_{n+1}[[t]]$ such that $t
\mapsto pt$.)
\ps

We use this to define functors $\Dsf(D)$ and $\Ddrf(D)$, as follows. 
Let $K_{\infty}=\bigcup_{n\geq 0} K_n$.
For $n\geq n(r)$ and $r>r(D)$, we define a 
$K_{\infty}$-vector space with a semi-linear action of $\Gamma$ by setting 
$$\Dsf(D):=(D_r \otimes_{\Ro_r}K_n)\otimes_{K_n} K_{\infty}.$$ 
By the discussion above, this space does not depend of 
the choice of $n, r$. In the same way, the $\Q_p$-vector spaces
$$\Ddrf(D):=(K_{\infty}\otimes_{K_n}K_n((t))\otimes_{\Ro_r}D_r)^{\Gamma},\, \, $$
$$\Fil^i(\Ddrf(D)):=(K_{\infty}\otimes_{K_n} t^iK_n[[t]]\otimes_{\Ro_r} D_r)^{\Gamma} \subset
\Ddrf(D), \, \, \, \forall i \in \Z,$$
are independent of $n\geq n(r)$ and $r>r(D)$. As
$K_{\infty}((t))^{\Gamma}=\Q_p$, $\Ddrf(D)$ so defined is a finite
dimensional $\Q_p$-vector-space whose dimension is 
less than $\rk_{\Ro}(D)$, and
$(\Fil^i(\Ddrf(D)))_{i\in \Z}$ is a decreasing, exhausting, and saturated,
filtration on $\Ddrf(D)$.\ps

We end by the definition of $\Dcf(D)$. Let
$$\Dcf(D):=D[1/t]^{\Gamma}.$$
It has an action of $\Q_p[\varphi]$ induced by the one on $D[1/t]$.
It has also a natural filtration defined as follows. Choose $r > r(D)$ and $n \geq n(r)$,
there is a natural inclusion
	$$\Dcf(D) \longrightarrow \Ddrf(D)$$
and we denote by $(\varphi^n(\Fil^i(\Dcf(D)))_{i\in \Z}$ the filtration induced from the one on
$\Ddrf(D)$. By the analysis above, this defines a 
unique filtration $(\Fil^i(\Dcf(V)))_{i \in \Z}$, independent of the above choices of $n$ and $r$. 
We summarize some of Berger's results (\cite{bergerinv}, \cite{bergereqdiff},
\cite[prop. 5.6]{colmezbki})) in the following proposition.
\ps
\begin{prop} \label{bergerthm} 
Let $V$ be a $\Qp$-representation of $\Gp$, and $$* \in \{crys, dR,
Sen\}.$$
Then $\cal{D}_{*}(\Dr(V))$ is canonically isomorphic to $D_{*}(V)$.
\end{prop}

\begin{definition} \label{deffg} 
We will say that a (not necessarily \'etale) 
$\fg$-module $D$ over $\Ro$ is crystalline (resp. de Rham) if
$\Dcf(D)$ (resp. $\Ddrf(D)$) has rank $\rk_{\Ro}(D)$ over 
$\Q_p$. The Sen polynomial of $D$ is the one of the semi-linear 
$\Gamma$-module $\Dsf(D)$. 
\end{definition}

Due to the lack of references we have to include the following lemma.

\begin{lemma} \label{fonctdr} Let $0 \longrightarrow D'
\longrightarrow D \longrightarrow D'' \longrightarrow 0$ be an exact
sequence of $\fg$-modules over $\Ro$. If $r \geq r(D), r(D'),
r(D'')$ is big enough, then
$D''_r=\Im(D_r \longrightarrow D'')$ and $D'_r=D_r \cap D'$.
\end{lemma}

\begin{pf} Fix $r_0 > r(D), r(D'), r(D'')$. By Lemma \ref{defdr}, 
\begin{equation} \label{recouvrement} 
X^\ast=\bigcup_{r \geq r_0} X_r^\ast,\, \,  {\rm for  X\, \,} \in \{D, D', D''\}. 
\end{equation}
We can then find $r_1 \geq r_0$ such that
$D'_{r_0} \subset D_{r_1}$ and $Im(D_{r_0} \longrightarrow D'') \subset
D''_{r_1}$. As $D \longrightarrow D''$ is surjective, we can choose
moreover some $r_2 \geq r_1$ such that $Im(D_{r_2} \longrightarrow D'')$
contains a $\Ro_{r_1}$-basis of $D''_{r_1}$. The exact sequence of the statement induces 
then for $r\geq r_2$ an exact sequence of $\Ro_r$-modules
\begin{equation} \label{exactlemme} 0 \longrightarrow K_r:=D_r\cap D' \longrightarrow D_r \longrightarrow D''_r \longrightarrow
0,
\end{equation}
with $K_r \supset D'_r$. As $D''_r$ is free, this sequence
splits, hence remains exact when base changed to $\Ro_s$, $s \geq r$. Using Lemma \ref{defdr}
(i), this implies that $\Ro_sK_r=K_s$ for $s\geq
r\geq r_2$. Moreover, $K_r$ is finite type over $\Ro_r$. Indeed, $D_r$ is
(free of) finite type over $\Ro_r$ and the sequence
(\ref{exactlemme}) splits. By formula (\ref{recouvrement}), we can then choose $r_3 \geq r_2$ such that
$K_{r_2} \subset D'_{r_3}$, and we get that $K_r=\Ro_rK_{r_2}=D'_r$ for $r\geq
r_3$. \end{pf}

\subsection{Triangular $\fg$-modules and trianguline representations over artinian $\Q_p$-algebras} ${}^{}$ \ps
\medskip
\label{triangularfg}
In all this subsection, $A$ is a finite dimensional local $\Q_p$-algebra with maximal
ideal $m$ and residue field $L:=A/m$.\ps

\subsubsection{$\fg$-modules of rank one over $\Ro_A$.} 

We begin by classifying all the $\fg$-modules which are 
free of rank $1$ over $\Ro_A$. Let $\delta: \Q_p^* \longrightarrow A^*$ be a continuous character. In the spirit of 
Colmez \cite[\S0.1]{colmeztri}, we define
the $\fg$-module $\Ro_A(\delta)$
which is $\Ro_A$ as $\Ro_A$-module but equipped with the $\Ro_A$-semi-linear actions of 
$\varphi$ and $\Gamma$ defined by 
$$\varphi(1):=\delta(p), \, \, \, \gamma(1):=\delta(\gamma), \forall \gamma\in \Gamma,$$
with the identification $\chi: \Gamma \isomo \Z_p^*$ understood. Recall that by class field theory $\chi$ extends uniquely to an isomorphism $\theta: \W^{\rm ab} \isomo \Q_p^*$ sending the geometric Frobenius to $p$, where $\W \subset \Gp$ is the Weil group of $\Q_p$. We may then view any $\delta$ as above as a continuous homomorphism $\W \longrightarrow A^*$. Such a homomorphism extends continuously to $\Gp$ if $v(\delta(p) \bmod m)$ is zero, and in this case we see that $$\Ro_A(\delta)=\Dr(\delta \circ \theta).$$ \ps 

Note that if $I\subsetneq A$ is an ideal, it is clear from the definition that $\Ro_A(\delta)\otimes_A A/I \isomo \Ro_{A/I}(\delta \bmod I)$. Moreover, if $D$ is a $\fg$-module over $\Ro_A$, we will set $D(\delta):=D\otimes_{\Ro_A}\Ro_A(\delta)$, and 
$$D^{\delta}=\{x \in D, \varphi(x)=\delta(p)x, \, \gamma(x)=\delta(\gamma)x\}\isomo
D(\delta^{-1})^{\varphi=1,\Gamma=1}=H^0(D(\delta^{-1})).$$

\begin{prop}\label{uni}
Any $\fg$-module free of rank $1$ over $\Ro_A$ is isomorphic to $\Ro_A(\delta)$ for a unique $\delta$. Such a module is isocline of slope $v(\delta(p) \bmod m)$.
\end{prop}

\begin{pf} By Lemma \ref{etalekedlaya}, a $\fg$-module $D$ of rank $1$ over $\Ro_A$ is automatically isocline of same slope as $D/mD$. 
As $\Ro_L(\delta \bmod m)$ has slope $v(\delta(p) \bmod m)$, and as $\Ro_A(\delta)(\delta')=\Ro_A(\delta\delta')$, we may assume that $D$ is \'etale. But in this case, the result 
follows from the equivalence of categories of 
Proposition~\ref{drig} (i), Lemma \ref{freeness}, and the fact that the continuous Galois characters $\Gp \longrightarrow A^*$ correspond exactly to the $\delta$ such that 
$v(\delta(p) \bmod m)=0$.
\end{pf}

\subsubsection{Definitions.}

\begin{definition} Let $D$ be a $\fg$-module which is free of rank $d$ over $\Ro_A$ and equipped with a strictly increasing filtration $(\Fil_i(D))_{i=0\dots d}$
:
$$\Fil_0(D):=\{0\} \subsetneq \Fil_1(D) \subsetneq \cdots \subsetneq \Fil_i(D) \subsetneq \cdots \subsetneq \Fil_{d-1}(D)
\subsetneq \Fil_d(D):=D,$$
of $\fg$-submodules which are free and direct summand as $\Ro_A$-modules. We call 
such a $D$ a {\it triangular} $\fg$-module 
over $\Ro_A$, and the 
filtration $\Tc := (\Fil_i(D))$ {\it a triangulation of $D$ over $\Ro_A$}.\ps

Following Colmez, we shall say that a $\fg$-module which is 
free of rank $d$ over $\Ro_A$ is {\it triangulable} if it can 
be equipped with a triangulation $\Tc$; we shall say that an $A$-representation 
$V$ of $\Gp$ which is free of rank $d$ over $A$ is {\it trianguline} if $\Dr(V)$ 
(which is free of rank $d$ over $\Ro_A$ by Lemma~\ref{freeness}) is triangulable. 
\end{definition}
\newcommand{\gr}{{\rm gr}}
\ps
Let $D$ be a triangular $\fg$-module. 
By Lemma \ref{uni}, for each 
$$\gr_i(D):=\Fil_i(D)/\Fil_{i-1}(D), \, \, \, 1\leq i \leq d,$$ 
is isomorphic to $\Ro_A(\delta_i)$ for some unique $\delta_i: \W
\longrightarrow A^*$. It makes then sense to define the 
{\it parameter of the triangulation} to be the continuous
homomorphism $$\delta:=(\delta_i)_{i=1,\cdots,d}: \Q_p^* \longrightarrow
(A^*)^d.$$ 
\index{Zdelta@$\delta=(\delta_i)_{i=1,\cdots,d}: \Q_p^* \longrightarrow
(A^*)^d$, the parameter of a trianguline representation over $A$ or of a triangular $\fg$-module over $\Ro_A$}

\subsubsection{Weights and Sen polynomial of a triangular $\fg$-module.}

As the following proposition shows, the parameter of a 
triangular $\fg$-module refines the data of its Sen
polynomial. It will be convenient to introduce, for a continuous 
character $\delta: \Q_p^* \longrightarrow A^*$, its {\it weight:} 

$$\omega(\delta):=-\left(\frac{\partial {\delta}_{|\Gamma}}{\partial 
\gamma}\right)_{\gamma = 1}=-\frac{\log(\delta(1+p^2))}{\log(1+p^2)} \, \, \in A.$$

\begin{prop} \label{senpol} Let $D$ be a triangular $\fg$-module over $\Ro_A$ and $\delta$ the parameter of a triangulation of $D$. Then the Sen polynomial of $D$ is 
$$\prod_{i=1}^d(T-\omega(\delta_i)).$$
\end{prop}

\begin{pf} Assume first that $d=1$, i.e. that 
$D=\Ro_A(\delta)$. We see that we may take $r(D)=(p-1)/p$ and that 
$D_r=A\Ro_r$ for $r>r(D)$. But then $\Dsf(D)$ has a $K_{\infty}$-basis on
which $\Gamma$ acts through $\delta_{|\Gamma}$, and
the result follows. The general case follows by induction on $d$ 
from the case $d=1$ and
Lemma \ref{fonctdr}. 
\end{pf}

\subsubsection{De Rham triangular $\fg$-modules. }

We now give a sufficient condition on a triangular $\fg$-module $D$ over $\Ro_A$ 
to be de Rham (see definition \ref{deffg}). A necessary condition is that
the $\Ro_A(\delta_i)$ are themselves de Rham, {\it i.e.} that each
$s_i:=\omega(\delta_i)$ is an integer (see the proof below).

\begin{prop}\label{criteriumDR} Let $D$ be a triangular $\fg$-module of
rank $d$ over $\Ro_A$, and let $\delta$ be its 
parameter. Assume that $s_i:=\omega(\delta_i) \in \Z$, and that 
$s_1 < s_2 < \cdots < s_d$, then $D$ is de Rham. \ps
If, moreover, $D_0:=D/mD$ is crystalline and satisfies
$\Hom(D_0,D_0(\chi^{-1}))=0$) (resp. is semi-stable),
then $D$ is crystalline (resp. semi-stable).
\end{prop}

\begin{pf} In this proof, $K_{\infty}[[t]]$ will always mean $\bigcup_{n\geq
1} K_n[[t]]$.\ps
Assume first that $d=1$ and $D=\Ro_A(\delta)$. 
If $s:=\omega(\delta) \in \Z$, then $\delta_{|\Gamma}\chi^s$ is a
finite order character of $\Gamma$. So $(t^sK_{\infty}[[t]]\otimes_{\Q_p} \delta_{|\Gamma})^{\Gamma}$
has $\Q_p$-dimension $1$ and is equal to $(t^sK_{\infty} \otimes_{\Q_p}
\delta_{|\Gamma})^{\Gamma}$. This concludes the case $d=1$.\ps
Let us show now that $D$ is de Rham by induction on $d$. By Lemma \ref{fonctdr}
we have for $r$ big enough and $i\in \Z$ an exact sequence
\begin{eqnarray*} 
0 \longrightarrow \Fil^i(\Ddrf(\Fil_{d-1}(D))) \longrightarrow
\Fil^i(\Ddrf(D))
\longrightarrow \Fil^i(\Ro_A(\delta_d)) \longrightarrow \\
\ \ \ \ \ \ \ \  H^1(\Gamma, (\Fil_{d-1}(D)_r)\otimes K_{\infty}[[t]]t^i).
\end{eqnarray*}
\noindent By the induction hypothesis applied to $\Fil_{d-1}(D)$,
$\Fil^{s_1}(\Ddrf(\Fil_{d-1}(D)))$ has dimension
$(d-1)\dim_{\Q_p}(A)$ and $\Fil^i(\Ddrf(\Fil_{d-1}(D))=0$ for $i>s_{d-1}$. By
the case $d=1$ studied above, it suffices then to show that 
$$H^1(\Gamma, \Fil_{d-1}(D)_r\otimes_{\Ro_r} K_{\infty}[[t]]t^{s_d})=0.$$
But the $\Gamma$-module $\Fil_{d-1}(D)_r\otimes_{\Ro_r} K_{\infty}[[t]]t^{s_d}$ 
is a successive
extension of terms of the form $$t^{s_d} K_{\infty}[[t]]\otimes_{\Q_p} \delta_i.$$
The $H^1(\Gamma,-)$ of each of these terms vanishes because $s_d>s_i$ if
$i<d$, and $$H^1(\Gamma,K_{\infty}[[t]]t^j)=0$$ 
for $j>0$. This concludes the proof that $D$ is de Rham. \ps
Recall that Berger's Theorem \cite[th\'eor\`eme A]{bergereqdiff} associates 
canonically to any de Rham $\fg$-module $D$ over $\Ro$
a filtered $(\varphi,N,\Gp)$-module $\WD(D)$. 
We see as in Lemma \ref{freeness} that $\WD(D)$ 
carries an action 
of $A$ and is free as $A$-module if (and only if) $D$ is free over $\Ro_A$. Moreover we have 
$$\WD(D_0)=\WD(D)/m \WD(D), \, \, \, {\text{ where }}\, \, \,D_0:=D/m D.$$ \noindent 
But if $D_0$ is semi-stable, the action of the inertia group $I_p \subset G_p$ on $\WD(D)/m\WD(D)$ is trivial. As the action of $I_p$ on $\WD(D)$
has finite image by definition,
 the $I_p$-module $\WD(D)$ is also trivial. If moreover $$\Hom(D_0,D_0(\chi^{-1}))=0,$$ then we 
get by induction on $i\geq 1$ that $$N: \WD(D) \longrightarrow \WD(D(\chi^{-1}))/m^i\WD(D(\chi^{-1}))$$ is zero, hence $N=0$ and so $D$ is crystalline. \end{pf}

\begin{remark} 
\begin{itemize}
\item[(i)] The proposition above may be viewed as a generalization of the fact that ordinary representations are
semi-stable (Perrin-Riou's Theorem \cite[Expos\'e IV]{PP}).
\item[(ii)] There exist triangulable \'etale $\fg$-modules of rank $2$ over $\Ro_L$
which are Hodge-Tate of integral weights $0< k$, 
but which are not de Rham (hence they have no triangulation whose parameter $\delta$ satisfies the
assumption of the Proposition~\ref{criteriumDR}. Instead, we have $s_1=k$ and $s_2=0$ 
with the notation of that proposition). For example, this is the case of the $\fg$-module of the restriction at $p$ 
of the Galois representation attached to any
finite slope, overconvergent, modular eigenform of integral weight 
$k>1$ and $U_p$-eigenvalue $a_p$ such that $v(a_p) > k-1$. 
\item[(iii)] It would be easy to show that a de Rham triangulable $\fg$-module
over $\Ro_L$ becomes semi-stable over a finite abelian extension of $\Q_p$, 
because it is true in rank $1$. 
Reciprocally, a $\fg$-module which becomes semi-stable over a finite abelian
extension of $\Q_p$ is triangulable over $\Ro_L$, where $L$ contains all the eigenvalues of 
$\varphi$ (to see this, mimic the proof of Proposition
\ref{ref=tri}). 
\end{itemize}
\end{remark}

\subsubsection{Deformations of triangular $\fg$-modules.}\label{deformtrifgdef} ${}^{}$ \ps
\newcommand{\DD}{{\rm Tri}}
Let $D$ be a fixed $\fg$-module free of rank $d$ over $\Ro_L$ and 
equipped with a triangulation $\Tc=(\Fil_i(D))_{i=0,\dots,d}$ with parameters $(\delta_i)$. We 
denote by $\cal{C}$ the category of local artinian $\Q_p$-algebras 
$A$ equipped with a map $A/m \isomo L$, and local homomorphisms inducing the
identity on $L$. \ps

Let $X_D: \cal{C} \longrightarrow {\rm Set}$ and $X_{D,\Tc}: \cal{C} \longrightarrow {\rm Set}$
denote the following functors. For an object $A$ of $\cal{C}$, 
$X_D(A)$ is the set of isomorphism classes of couples $(D_A,\pi)$ where $D_A$ is a $\fg$-module
free over $\Ro_A$ and $\pi: D_A \longrightarrow D$ is a $\Ro_A$-linear $\fg$-morphism inducing an 
isomorphism $D_A\otimes_A L \isomo D$; $X_{D,\Tc}(A)$ is the set 
of isomorphism classes of triples $(D_A, \pi, (\Fil_i(D_A)))$ 
where~:
\begin{itemize}
\item[(i)] $(D_A, (\Fil_i(D_A)))$ is a triangular $\fg$-module of rank $d$ over $\Ro_A$, 
\item[(ii)] $\pi: D_A \longrightarrow D$ is a $\Ro_A$-linear $\fg$-morphism inducing an 
isomorphism $D_A\otimes_A L \isomo D$ such that $\pi(\Fil_i(D_A))=\Fil_i(D)$.
\end{itemize}

There is a natural ``forgetting the triangulation'' morphism of functors 
$X_{D,\Tc} \longrightarrow X_D$ that makes in favorable cases $X_{D,\Tc}$ a subfunctor of $X_D$ 
(as in \cite{colmeztri}, we denote by $x$ the identity 
character $\Qp^* \longrightarrow \Qp^*$).
\begin{prop} \label{subfunctor} 
Assume that for all $i<j$, $\delta_i \delta_j^{-1} \not \in x^\N$. Then $X_{D,\Tc}$ 
is a subfunctor of $X_D$.
\end{prop} 
\begin{pf} We have to show that if $A$ is an object of $\cal{C}$, 
and  $(D_A,\pi) \in X_D(A)$ is a deformation of $D$, $D_A$ has at most one 
triangulation that satisfies (ii) above. That is to say, we have to prove that 
if $\Tc=(\Fil_i(D_A))$ is a triangulation of $D_A$ 
satisfying (ii), then $\Fil_1(D_A)$ is uniquely determined as a submodule 
of $D_A$, $\Fil_2(D_A)/\Fil_1(D_A)$ is uniquely determined as
a submodule of $D_A/\Fil_1(D_A)$, and so on. For this  
note that $\Fil_j(D)/\Fil_{j-1}(D) \simeq \Ro_L(\delta_j)$, and that 
$D/\Fil_j(D)$ is a successive extension of $\Ro_L(\delta_i)$ with 
$i>j$, so that
$$\Hom(\Ro_L(\delta_j), D/\Fil_j(D))=0$$ by the hypothesis  on the $\delta_i$ and Proposition~\ref{structurer} (ii). So
 we can apply 
the following lemma, and we are done. 
\end{pf}
\begin{lemma} Let $(D_A,\pi) \in X_D(A)$, $\delta:\Q_p^\ast \longrightarrow A^\ast$ be a continuous character and
$\bar \delta=\delta \pmod{m}$. Assume $D$ has a saturated, rank 1, $\fg$-submodule $D_0 \simeq \Ro_L(\bar \delta)$ such that 
$(D/D_0)^{\bar \delta}=0$. Assume moreover that $D_A$ contains a $\Ro$-saturated $\fg$-module $D'$ isomorphic to $\Ro_A(\delta)$.
Then $\delta$ is the unique character of $\Q_p^\ast$ having this property, and $D'$ the unique submodule.
\end{lemma}
\begin{pf}
 We may assume by twisting that $\delta=1$ (hence $\overline{\delta}=1$ also). 
Let $\delta': \Q_p^* \longrightarrow A^*$ 
lifting $1$, and assume that $D_A$ has some $\Ro$-saturated submodules
$D_1 \isomo \Ro_A$ and $D_2 \isomo\Ro_A(\delta')$. By assumption, $H^0(D/D_0)=0$, and $D_i/m
D_i=D_0$
for $i=1, 2$ (see Remark \ref{remfree}). 
A devissage and the left exactness of the functor $H^0(-)$ show that $$H^0(D_A/D_i)=0, \, \, \, i=1, 2.$$ This implies that the inclusion 
$H^0(D_i) \longrightarrow H^0(D_A)$ is an equality, hence $$H^0(D_1)=H^0(D_2) \subset D_2.$$ As $D_1=\Ro H^0(D_1)$, we have 
$D_1 \subset D_2$, and $D_1=D_2$ since $D_1$ and $D_2$ are saturated and have the same $\Ro$-rank. We conclude that $\delta'=1$ by Proposition \ref{uni}.
\end{pf} 

	We will give below a criterion for the relative representability of 
$X_{D,\Tc} \longrightarrow X_D$, but we need before to make some preliminary remarks. Let $F(-)$ be the functor on $\fg$-modules over $\Ro_L$ defined by
$$F(E)=\{ v \in E, \exists n\geq 1\, |\, \, \, \forall \gamma \in \Gamma, (\gamma-1)^nv=0, \, (\varphi-1)^nv=0\}.$$
This is a left-exact functor, and $F(E)$ inherits a commuting continuous action of $\varphi$ and $\Gamma$, hence of $\Q_p^*$, as well as a commuting action of 
$A$ if $E$ does. 

\begin{lemma}\label{unipabfg} ${}^{}$ \begin{itemize}
\item[(i)] For any $\fg$-module $E$ over $\Ro_L$, $F(E) \neq 0 \Leftrightarrow \Hom_{\fg}(\Ro_L,E) \neq 0$.
\item[(ii)] $F(\Ro_L(\delta))=0$ if $\delta \notin x^{-\N}$, and $F(\Ro_L)=L$.
\item[(iii)] Let $A \in \cal C$ and $\delta: \Q_p^* \longrightarrow A^*$ a continuous homomorphism such that $\overline{\delta}=1$.
The natural inclusion $A \subset \Ro_A(\delta)$ induces a $\Q_p^*$-equivariant\footnote{Note that $A \subset \Ro_A(\delta)$ has a natural $A$-linear 
action of $(\varphi,\Gamma)$, 
hence of $\Q_p^*$, namely via the character $\delta$ by definition.} isomorphism $A \isomo F(\Ro_A(\delta))$, as well as 
an isomorphism $F(\Ro_A(\delta))\otimes_{\Q_p} \Ro \isomo \Ro_A(\delta)$. 
\end{itemize}
\end{lemma}

\begin{pf} Assertion (i) is an immediate consequence of the definition and of the fact that $\varphi$ and $\Gamma$ commute. Let us 
check assertion (ii). We may assume that $\delta=1$ by (i) and Prop. \ref{structurer} (ii). Let $\gamma \in \Gamma$ be a nontorsion element. We claim that for $f \in \Ro_L$ and $n\geq 1$ 
$$(\gamma-1)^n f =0 \Rightarrow f \in L.$$ 
Indeed, we may assume that $n\leq 2$. If $n=1$, then $f(z^\gamma)=f(z)$ so $f$ 
is constant on each circle $|z-1|=r$ with $r\geq r(f)$, because $\gamma$ is nontorsion and an analytic function on a $1$-dimensional 
affinoid has only a finite number of zeros. If $n=2$, we may assume by absurdum that $(\gamma-1)(f)=1$, which means that $$f(z^{\gamma})=f(z)+1$$ 
on $r(f)\leq |z-1| <1$. But for $z$ a $p^m$-th root of unity in the annulus $r(f)\leq |z-1| <1$, we get by applying $p^m$ times the previous equation that 
$f(z)=f(z)+p^m$, a contradiction.\par
	Let us check assertion (iii). It is clear that $A \subset F(\Ro_A(\delta))$. Moreover, as $F(\Ro_L)=L$ by (ii), the 
left-exactness of $F$ shows that the lenght of $F(\Ro_A(\delta))$ is less than the lenght of $A$. In particular, 
the previous inclusion is an equality, and the last assertion of the stament holds by definition of $\Ro_A(\delta)$.
\end{pf}

\begin{prop}\label{representri}
Assume that for all $i<j$, $\delta_i \delta_j^{-1} \not \in x^\N$. Then $X_{D,\Tc} \longrightarrow X_D$ is relatively representable.
\end{prop}

\begin{pf} By Prop. \ref{subfunctor}, we already know that $X_{D,\Tc}$ is a subfunctor of $X_D$. By \cite[\S23]{Mazurfermat}, we have to check three conditions (see also \S19 of {\it loc.cit.}).\par 
	{\it First condition:} if $A \longrightarrow A'$ is a morphism in $\cal C$ and if $(D_A, \pi) \in X_{D,\Tc}(A)$, then $$(D_A\otimes_A A', \pi\otimes_A A') \in X_{D,\Tc}(A').$$ 
This is obviously satisfied as $(\Fil_i(D_A)\otimes_A A')$ is a triangulation of $D_A \otimes_A A'$ lifting $D$.\ps

	{\it Second condition:\footnote{This is actually called condition (3) in {\it loc.cit}.}} 
if $A \longrightarrow A'$ is an {\it injective} morphism in $\cal C$, and if $(D_A, \pi) \in X_D(A)$, 
then $$(D_A\otimes_A A', \pi\otimes_A A') \in X_{D,\Tc}(A') \implies (D_A,\pi) \in X_{D,\Tc}(A).$$ 
 Arguing by induction on $d=\dim_{\Ro_L} D$, it is enough to show that $D_A$ has a $\fg$-submodule $E$ which is free of rank 
$1$ over $\Ro_A$, saturated, and such that the natural map $$\pi: E \longrightarrow D$$
surjects onto ${\rm Fil}_1(D)$ (the fact that $E$ is a direct summand as $\Ro_A$-module will follow then from Lemma \ref{lib}). 
By twisting if necessary, we may assume that $\delta_1=1$.
\ps
By (ii) of Lemma \ref{unipabfg}, the left-exactness of $F$, and the assumption on the $\delta_i$, we have 
\begin{equation}\label{calcFdeD} F(D)=F(\Fil_1(D))=L.\end{equation}
	Let $D_{A'}=D_A \otimes_A A'$, $T'_1:=\Fil_1(D_{A'})$ and $T_1:=T'_1\cap D_A$. 
	Lemma \ref{unipabfg} (ii) again and a devissage show that $F(D_A/T_1) \subset F(D_{A'}/T'_1)=0$, so the natural inclusions 
\begin{equation}\label{reducttoA} F(T_1) \isomo F(D_A)\, \, \, \, and\, \, \, \, F(T'_1) \isomo F(D_{A'}) \end{equation}
are isomorphisms. Moreover, the fact that $F(D)=L$ and another devissage\footnote{For more details about this devissage, the reader can have a look at Lemma 
\ref{abstract} of the next section, in which it is studied in a more general situation.} show that for each ideal $I$ of $A$, and each finite lenght $A$-module $M$, if $l$ denotes the lenght function (so that $L$ has lenght $1$) 
then $$l(F(ID_A))\leq l(I)\, \, \, \, {\rm and}\, \, \, \,  l(F(D_A \otimes_A M)) \leq l(M).$$
The inequalities above combined with the exact sequences 
$$0 \longrightarrow F(D_A) \longrightarrow F(D_{A'}) \longrightarrow F(D_A \otimes_A A'/A), $$
$$0 \longrightarrow F(mD_A) \longrightarrow F(D_A) \longrightarrow F(D),$$
show then that 
\begin{equation}\label{cordala} l(F(D_A))=l(A), \, \, \, \, {\rm and}\, \, \, F(D_A)\otimes_A L \overset{\neq 0}{\longrightarrow} F(D)=F(\Fil_1(D))=L.\end{equation}
In particular, there is an element $v \in F(D_A) \subset D_A$ whose image is nonzero in $F(D) \subset D_A/mD_A$, thus this 
element $v$ generates a free\footnote{The claim here is that for $A \in \cal C$, any 
free $A$-module $M$ and any $v \in M$, if the image of $v$ in $M/mM$ is nonzero, then $Av \simeq A$. Indeed, the assumption implies 
that $\Tor_1^A(M/Av,k)=0$, so $M/Av$ is free by Nakayama's lemma (see Lemma. \ref{lib}) and so is $Av$.} $A$-submodule of $F(D_A)$.
By (\ref{cordala}) we get that $F(D_A)=Av$ is free of rank $1$ over $A$ and 
the nonzero map there is actually an isomorphism. Of course, the same assertion holds if we replace the $A$'s in it by $A'$, 
as $A'=F(T'_1)=F(D_{A'})$ by Lemma \ref{unipabfg} (iii) and (\ref{reducttoA}). As a consequence, the natural map 
\begin{equation}\label{isomFTAA'} F(T_1) \otimes_A A' \longrightarrow F(T'_1)\end{equation}
is an isomophism, at it is so modulo the maximal ideal.\footnote{In particular, if we write $T'_1=\Ro_{A'}(\delta)$, 
the isomorphism above and Lemma \ref{unipabfg} (iii) show that $\delta(\Q_p^*)\subset A^*$.} Set $$E:=\Ro F(T_1) \subset D_A.$$
We claim that $E$ has the required properties to conclude. Recapitulating, we have a sequence of maps
$$F(T_1) \otimes_A \Ro_A \hookrightarrow F(T_1) \otimes_A \Ro_{A'} \isomo F(T'_1) \otimes_{A'} \Ro_{A'} \isomo T'_1.$$
As $E$ is the image of the composition of all the maps above, we get that $$F(T_1)\otimes_A \Ro_A \isomo E$$ is free of rank 
$1$ over $\Ro_A$. We already showed that $\pi(E)=\Fil_1(D)$, hence it only remains to check that $E$ is saturated in $D_A$. 
But this holds as $E$ is saturated in $T'_1$, which is saturated in $D_A'$, and we are
done\footnote{Actually, using Lemma \ref{lemmeAA'} we even see that $E=T_1$.}. \ps
%if $f\neq 0 \in \Ro$ and $fv \in E$ for some $v \in D_A$, then $v \in T'_1$, and then $v \in E$. 
	{\it Third condition:\footnote{Precisely, this is what is left to check condition (2) of {\it loc. cit.} once (1) and (3) are known to hold, 
because with the notations there $A\times_C B \subset A \times_L B$ if $L$ is the residue field of $C$. This reduction is also explained in the proof of \cite[Prop. 8.7]{kis}.}} for $A$ and $A'$ in $\cal C$, if 
$(D_A,\pi) \in X_{D,\Tc}(A)$ and $(D_{A'},\pi') \in X_{D,\Tc}(A')$, then for $B=A\times_L A'$, the natural object 
	$$(D_B=D_A\times_D D_{A'},\pi_B=\pi \circ \pr_1=\pi' \circ \pr_2)$$ lies in $X_{D,\Tc}(B)$. But it is 
clear that the filtration $(\Fil_i(D_A)\times_{D} \Fil_i(D_{A'}))$ is a triangulation of $D_B$ lifting $\Tc$, and we are done.
\end{pf}

Let us consider the natural morphism $$\diag: X_{D,\Tc} \longrightarrow \prod_{i=1}^d X_{\gr_i(D)}.$$ 
Recall that $x$ is the identity character $\Qp^* \longrightarrow \Qp^*$; recall also that 
$\chi=x|x|$ is the cyclotomic character.
\begin{prop}\label{deform} Assume that for all $i< j, \, \, \delta_i\delta_j^{-1}\not\in \chi\, x^{\N}$, 
then \begin{itemize}
\item[(i)] $X_{D,\Tc}$ is formally smooth,
\item[(ii)] for each $A \in {\rm Ob}(\cal{C})$, $\diag(A)$ is surjective.
\end{itemize}
If we assume moreover that $\delta_i\delta_j^{-1} \not\in x^{-\N}$ for $i<j$, then
$$\dim_L X_{D,\Tc}(L[\varepsilon])=\frac{d(d+1)}{2}+1.$$
\end{prop}

\begin{pf} 
Recall that (i) means that for 
$A \in \cal{C}$ and $I \subset A$ an ideal such that $I^2=0$, the natural map 
$X_{D,\Tc}(A) \longrightarrow X_{D,\Tc}(A/I)$ is surjective.\ps
Assume first that $d=1$, so the assumption is empty. The maps $\diag(A)$ are bijective 
(hence (ii) is satisfied), and by Proposition \ref{uni}, $X_{D,\Tc}$ is isomorphic the functor 
$\Hom_{\rm cont}(\W,-)$ which is easily seen to be formally smooth 
(and even pro-representable by ${\rm Spf}(L[[X,Y]])$), hence (i) is satisfied. \ps
Let us show now (i) and (ii), we fix $A \in \cal{C}$ and $I 
\subset m$ an ideal of length $1$. Let $U \in X_{D,\Tc}(A/I)$, and let 
$V=(V_i) \in \prod_i(X_{\gr_i(D)}(A))$ be any lifting of $\prod_i \gr_i(U)$. We are looking for an element 
$U' \in X_{D,\Tc}(A)$ with graduation $V$ and reducing to $U$ modulo $I$. 
We argue by induction on $d$. By the paragraph 
above, we already know the result when $d=1$. By the case $d=1$ again, 
we may assume that $\gr_d(U)$ is the trivial $\fg$-module over $\Ro_{A/I}$ 
(note that the assumption on the $\delta_i$ is 
invariant under twisting), and we have 
to find a $U'$ whose $\gr_d(U')$ is also trivial. Let $\Tc'$ denote the triangulation $(\Fil_i(D))_{i=0,\dots,d-1}$ of $\Fil_{d-1}(D)$.
By induction, $X_{\Fil_{d-1}(D),\Tc'}$ 
is formally smooth and satisfies (ii), hence we can 
find an element $U'' \in X_{\Fil_{d-1}(D),\Tc'}(A)$ 
lifting $\Fil_{d-1}(U)$ and such that $\gr_i(U'')$ lifts $V_i$ for $i=1, \cdots, d-1$. 
It suffices then to show that the natural map 
$$H^1(U'') \longrightarrow H^1(\Fil_{d-1}(U))$$
is surjective. But by the cohomology exact sequence, 
its cokernel injects to $$H^2(\Fil_{d-1}(D(\delta_d^{-1}))).$$ 
But this cohomology group is $0$ by assumption and 
by Lemma~\ref{calculH2}.
\ps
Let us compute by induction on $d=\rk_{\Ro}(D)$ the dimension 
$$x_d:=\dim_L(X_{D,\Tc}(L[\varepsilon])).$$ \noindent We have seen that $x_1=2$.
We showed above that the natural maps 
$$X_{D,\Tc}(A) \longrightarrow
X_{\Fil_{d-1}(D),\Tc'}(A) \times X_{\gr_d(D)}(A)$$ are surjective. 
Applying this to $A=L[\varepsilon]$ with $\varepsilon^2=0$, we find that
$x_d=x_{d-1}+2+n$ where $n$ is the dimension of 
$H^1(\Fil_{d-1}(D)(\delta_d^{-1}))/L.[D(\delta_d^{-1})]$, where
$[D(\delta_d^{-1})]$ is the class of the extension defined by
$D(\delta_d^{-1})$. But recall that for
$\delta \not\in x^{-\N} \cup \chi x^{\N}$,
$H^i(\Ro_L(\delta))=0$ for $i\neq 1$ and
 $\dim_L(H^1(\Ro_L(\delta)))=1$ by \cite[prop. 3.1]{colmeztri}, Lemma
\ref{calculH2} and \cite[thm. 3.9]{colmeztri}.
By the assumptions and the cohomology exact sequence, it implies that 
$H^1(\Fil_{d-1}(D)(\delta_d^{-1}))$ has dimension $d-1$, hence $n=d-2$. 
\end{pf}

\begin{lemma} \label{calculH2} 
$H^2(\Ro_A(\delta))=0$ if $(\delta \bmod m) \not \in \chi
x^{\N}$. 
\end{lemma}

\begin{pf} By the cohomology exact sequence, we may assume that $A=L$. But then
$H^2(\Ro_L(\delta))=H^0(\Ro_L(\chi\delta^{-1}))=0$ by \cite[prop.
3.1]{colmeztri}. The fact that for any $\fg$ module $D$ 
over $\Ro$, we have $H^0(D)=H^2(D^*(\chi))$ should hold by mimicking 
Herr's original argument. We warn the reader that at the moment, 
there is unfortunately no written reference for that result.
\end{pf}

\subsubsection{Trianguline deformations of trianguline representations.}

\label{deftri}

The notions of the last paragraph have their counterpart in terms of trianguline representations.
Let $V$ be a trianguline representation over $L$, and suppose given a triangulation $\Tc$ on $D:=\Dr(V)$.\ps

We define the functor $$X_V : \cal{C} \longrightarrow {\rm Set}$$ as follows: for $A \in {\rm Ob}(\cal{C})$,
$X_V(A)$ is the set of equivalence classes of {\it deformations} of $V$ over $A$, 
that is, $A$-representations $V_A$ of $\Gp$ which are free over $A$ and equipped with an $A[\Gp]$-morphism $\pi : V_A \longrightarrow A$ inducing 
an isomorphism $V_A \otimes_A L \longrightarrow V$. In the same way, we define a functor $$X_{V,\Tc} : \cal{C} \longrightarrow {\rm Set}$$ such that $X_{V,\Tc}(A)$ is the set of 
equivalence classes of {\it trianguline deformations of $(V,\Tc)$}, that is couples $(V_A,\pi) \in X_V(A)$ together with a triangulation $\Fil_i$ of $D_\rig(V_A)$ which makes $(\Dr(V_A),\Dr(\pi),\Fil_i)$ an element of $X_{D,\Tc}(A)$.\ps

The main fact here is that those functors are not new :
\begin{prop} \label{XVXD} The functor $\Dr$ induces natural isomorphisms of functors 
$X_V \simeq X_D$ and $X_{V,\Tc} \simeq X_{D,\Tc}.$
\end{prop}
\begin{pf} The second assertion follows immediately from the first one, 
since $X_{V,\Tc}(A) = X_{V}(A) \times_{X_D(A)} X_{D,\Tc}(A)$ for any $A$ in $\cal{C}$. 

To see that $\Dr$ induces a bijection 
$X_V(A) \longrightarrow X_D(A)$, we note that the injectivity is obvious because of the full faithfulness of $\Dr$, and that the surjectivity follows from the fact that if $(D_A,\pi)$ is an element of $X_D(A)$, $D_A$ is a successive extension of $D$ as a $\fg$-module 
over $\Ro$, hence it is \'etale by Lemma~\ref{etalekedlaya}; so $D_A$ is $\Dr(V_A)$ 
for some representation $V_A$ over $L$ which is free over 
$A$ by Lemma~\ref{freeness}.
% and $\pi$ is $\Dr(\pi')$ by full faithfulness again.
\end{pf}
 
\subsection{Refinements of crystalline representations}\label{refinement} (See
\cite[\S3]{Maz}, \-
\cite[\S7.5]{Ch}, \cite[\S6]{BC})
\ps

\subsubsection{Definition.}
Let $V$ be finite, $d$-dimensional, continuous, $L$-representation of $\Gp$. 
We will assume that $V$ is crystalline and that the crystalline Frobenius $\varphi$ acting on 
$\Dc(V)$ has all its eigenvalues in $L^*$. \ps

By a {\it refinement} of $V$ (see \cite[\S3]{Maz}) \index{F@$\Ref$, a refinement of a crystalline representation} we mean the data of a full $\varphi$-stable
$L$-filtration $\Ref=(\Ref_i)_{i=0,\dots,d}$ of  $\Dc(V)$:
$$\Ref_0=0 \subsetneq \Ref_1 \subsetneq \cdots \subsetneq \Ref_d=\Dc(V).$$
We remark now that any refinement $\Ref$ determines {\it two orderings}: \begin{itemize}
\item[(Ref1)] It determines an ordering 
$(\varphi_1,\cdots,\varphi_d)$ of the eigenvalues of $\varphi$, defined by the formula $$\det(T-\varphi_{|\Ref_i})=\prod_{j=1}^i (T-\varphi_j).$$
Obviously, if all these eigenvalues are distinct such an ordering conversely determines $\Ref$. 
\item[(Ref2)] It determines also an ordering $(s_1,\cdots,s_d)$ on the set of Hodge-Tate weights of $V$, defined by the property that the jumps of the weight filtration of $\Dc(V)$ induced on $\Ref_i$ are $(s_1,\cdots,s_i)$\footnote{As $\Ref_i \subset \Ref_{i+1}$, the weights of $\Ref_i$ are also weights of $\Ref_{i+1}$, hence the definition of the $s_i$ makes sense.}.
\end{itemize}

More generally, the definition above still makes sense when
$D$ is any crystalline $\fg$-module over $\Ro_L$ (see definition 
\ref{deffg}), {\it i.e. not necessarily \'etale}, such that $\varphi$
acting on $\Dcf(D)=D[1/t]^{\Gamma}$ has all its eigenvalues in $L^*$. It will
be convenient for us to adopt this degree of generality. \ps

\subsubsection{Refinements and triangulations of $\fg$-modules.}

The theory of refinements has a simple interpretation in terms of 
$\fg$-modules that we now clarify. Let $D$ be a crystalline $\fg$-module as above and
let $\Ref$ be a refinement of $D$. We can construct from $\Ref$ a
filtration $(\Fil_i(D))_{i=0,\cdots,d}$ of $D$ by setting
	$$\Fil_i(D):=(\Ro[1/t]\Ref_i)\cap D,$$
which is a finite type saturated $\Ro$-submodule of $D$ by Lemma \ref{structurer}.

\begin{prop} \label{ref=tri} The map defined above $(\Ref_i) \mapsto
(\Fil_i(D))$ induces
a bijection between the set of refinements of $D$ and the set of
triangulations of $D$, whose inverse is
$\Ref_i:=\Fil_i(D)[1/t]^{\Gamma}$. \ps
	In the bijection above, for $i=1,\dots, d$, the graded piece
$\Fil_i(D)/\Fil_{i-1}(D)$ is
isomorphic to $\Ro_L(\delta_i)$ where $\delta_i(p)=\varphi_i p^{-s_i}$ and 
${\delta_i}_{|\Gamma}=\chi^{-s_i}$, where the $\varphi_i$ and $s_i$ are defined by (Ref1) and (Ref2). 
\end{prop}

\begin{pf} We have $\fr(\Ro_L)^{\Gamma}=L$, hence the natural $\fg$-map 
$$D[1/t]^{\Gamma} \otimes_{L} \Ro_L \longrightarrow
D[1/t]$$
is injective. But it is also surjective because $D$ is assumed to be 
crystalline, hence it is an isomorphism. We deduce from this that any
$\fg$-submodule $D'$ of $D[1/t]$ over $\Ro_L[1/t]$ 
writes uniquely as $\Ro_L[1/t]\otimes_L F=\Ro_L[1/t]F$, where $F=D'^{\Gamma}$ is a
$L[\varphi]$-submodule of $D[1/t]^{\Gamma}$. This proves the first part of
the proposition. \ps
Let us prove the second part. By what we have just said, the eigenvalues of
$\varphi$ on the rank one $L$-vector space, $(\Fil_i(D)/\Fil_{i-1}(D)[1/t])^{\Gamma}$ is
$\varphi_i$. As a consequence,  
the rank one $\fg$-module $\Fil_i(D)/\Fil_{i-1}(D)$, which has the
form $\Ro_L(\delta_i)$ for some $\delta_i$ by Proposition \ref{uni}, satisfies
$$\delta_i(p)=\varphi_i p^{-t_i}, {\delta_i}_{|\Gamma}=\chi^{-t_i}$$
for some $t_i \in \Z$ by Proposition \ref{structurer} (i). By Proposition \ref{senpol}, the $s_i$ are (with
multiplicities) the Hodge-Tate weights of $V$, and it remains to show that $t_i=s_i$. We need the following essential lemma. \ps
\begin{lemma}\label{mainlemma} Let $D$ be a $\fg$-module over $\Ro_A$, $\lambda \in A^*$, and $v \in \Dcf(D)^{\varphi=\lambda}$. Then $v \in \Fil^i(\Dcf(D))$ if, and only if, $v \in t^iD$.
\end{lemma}

\begin{pf} For $n\geq 1$, let $t_n:=\frac{z^{p^n}-1}{z^{p^{n-1}}-1} \in \Ro$. For $r> r(D)$ and $n\geq n(r)$, then for $i\in \Z$, Lazard's theorem shows that \begin{equation} \label{divisibilite} \Ro_r[1/t] \cap K_n[[t]]t^i=t_n^i \Ro_r, \, \, \, {\rm and}\, \, \, \Ro_r[1/t] \bigcap_{n\geq n(r)} K_n[[t]]t^i=t^i \Ro_r.\end{equation} 
Let now $D$ be $\fg$-module over $\Ro_A$. 
By definition of the filtration on $\Dcf(D)$, we get in particular, for $r$
big enough and $n\geq n(r)$, \begin{equation}\label{filtr} \varphi^n(\Fil^i(\Dcf(D)))=\Dcf(D)\cap t_n^i D_r,
\end{equation}
Let $v \in \Dcf(D)$ be such that $\varphi(v)=\lambda v, \lambda \in A^*$. Then (\ref{divisibilite}) and (\ref{filtr}) show that $v \in \Fil^i(\Dcf(D))$ if, and only if, $v \in t^iD$. 
\end{pf}

We now show that $t_i=t_i(D,\Ref)$ coincides with $s_i=s_i(D,\Ref)$ 
by induction on $d$. Let $v \neq 0 \in \Ref_1$. 
As $$v \in \Fil^{s_1}(\Dcf(D)) \backslash \Fil^{s_1+1}(\Dcf(D))$$ 
by assumption, Lemma \ref{mainlemma} above shows 
that $t^{-s_1}v \in D \backslash tD$. By Proposition 
\ref{structurer} (ii), this shows that $\Ro\Ref_1 t^{-s_1}$ 
is saturated in $D$, hence is $\Fil_1(D)$, and $s_1=t_1$. 
Let us consider now the $\fg$-module $$D'=D/\Fil_1(D).$$ 
It is crystalline with $\Dcf(D')=\Dcf(D)/\Ref_1$, with Hodge-Tate weights
(with multiplicities) the ones of $D$ deprived of $s_1$, 
and has also a natural refinement defined by $\Ref'_i=\Ref_{i+1}/\Ref_1$. The weight filtration on $\Dcf(D')$ is the quotient filtration $((\Fil^j(\Dcf(D))+\Ref_1/\Ref_1))_{j\in \Z}$. As a consequence, $s_i(D',\Ref')=s_{i+1}(D,\Ref)$ if $i=1,\cdots, d-1$. But by construction, for $i=1,\cdots, d-1$ we have also $t_i(D',\Ref')=t_{i+1}(D,\Ref)$. Hence $t_i=s_i$ for all $i$ by the induction hypothesis.
\end{pf}

\begin{remark} \label{crystri} 
In particular, Proposition \ref{ref=tri} shows that crystalline representations are trianguline, and that the set of their triangulations are in natural bijection with the set of their refinements.
\end{remark}

\begin{definition} Let $\Ref$ be a refinement of $D$ (resp. of $V$). 
The {\it parameter} of $(D,\Ref)$ (resp. $(V,\Ref)$) is the parameter of the triangulation of $D$
(resp. $\Dr(V)$) associated to $\Ref$, {\it i. e.} the continuous character
$$\delta:=(\delta_i)_{i=1,\dots,d}: \Q_p^* \longrightarrow (L^*)^d$$ defined
by Proposition \ref{ref=tri}.
\end{definition}

\subsubsection{Non critical refinements.}\label{crit}

Let $(D,\Ref)$ be a refined crystalline $\fg$-module as in \S\ref{refinement}. 
We assume that its Hodge-Tate weights are two by two distinct, and denote
them by $$k_1<\cdots<k_d.$$ 

\begin{definition} We say that $\Ref$ is {\it non critical} if $\Ref$ is in general position 
compared to the weight filtration on
$\Dcf(D)$, {\it i.e.} if for all $1 \leq i \leq d$, we have a direct sum
$$\Dcf(D)=\Ref_i \oplus \Fil^{k_i+1}(\Dcf(D)).$$
\end{definition}

\begin{remark} Assume that $D=\Dr(V)$ for a crystalline $V$ in what follows.
%The condition of non criticality we just defined is rather weak, but seems to be 
%reasonable :
\label{remcrit} \begin{itemize}
\item[(i)] If $d=1$, the unique refinement of $V$ is always non critical. If $d=2$, all refinements of $D$ are
non-critical, excepted when $V$ is a direct sum of two (crystalline) characters.
 
\item[(ii)] Another natural definition of non criticality would be the condition
\begin{eqnarray} \label{num} \ \ \ \ \ \ \ 
\forall i\in\{1,\dots,d-1\},\ \ \ v(\varphi_1)+\dots+v(\varphi_i) < k_1+\dots+k_{i-1}+k_{i+1}.\end{eqnarray}
We call a refinement satisfying this condition {\it numerically non critical}\footnote{Mazur introduced in \cite{Maz} a variant of this condition, 
namely $v(\varphi_{i-1})<k_i<v(\varphi_{i+1})$ for $i=2,\dots,d-1$ which is equivalent to ours for $d \leq 3$.}.
The weak admissibility of $D_\cris(V)$ shows that a numerically non critical refinement is non critical 
in our sense, but the converse is false. However, as the following example shows,
it may be very hard in practice to prove that a refinement is non critical, when it is not already numerically 
non critical.
 
Assume $d=2$. The numerical non criticality
condition (\ref{num}) reduces to $v(\varphi_1)<k_2$; note that this is the hypothesis 
appearing in the weak form of Coleman's {\it classicity of small 
slope $U_p$-eigenforms} result (\cite{col1}). 
Assume $(V,\cal{F})$ is the non ordinary refinement of the restriction at $p$ 
of the Galois representation attached to a classical, ordinary, modular form $f$. 
Then $\cal{F}$ satisfies $v(\varphi_1)=k_2$, so 
it is not numerically non-critical.  
On the other hand $\cal{F}$ is non critical
if and only if $V$ is not split, which is an open problem when $f$ is not CM
(see also \cite{lissite} for the study of the non ordinary 
Eisenstein points, the CM points can be handled in the same way).
Note that from the existence of overconvergent companion forms \cite[Thm. 1.1.3]{BrEm}, $\cal{F}$ is non critical if and only if $f$ is not in the image of the Theta operator, which is exactly the condition found by Coleman in his study of boundary case of his "classicity criterion".
We take this as an indication of the relevance of our definition of non criticality. 
\item[(iii)] ({\it Comparison with the $U$-non criticality condition of \cite[\S7.5]{Ch}}) If $\Ref$ is $U$-non critical in the sense of \cite[Def. \S7.5]{Ch}, it is numerically non critical\footnote{To see this, deduce from the formula of Proposition 7.3 {\it loc. cit.} that for each $j$ we have $v_j< m_j+1$, which is precisely the numerical criterion.}, hence non critical. However, the $U$-non criticality is strictly stronger in general, even for $d=3$ (see the formula in the example {\it loc. cit.}). 
\item[(iv)] If $\varphi$ is semisimple (which is conjectured to occur in the geometric situations), we see at once that $V$ always
admits a non critical refinement in our sense. 
However, all the refinements of $V$ may be numerically critical. 
Examples occurs already when $d=3$. Indeed, (\ref{num}) is equivalent in this case to 
$v(\varphi_1) < k_2 < v(\varphi_3)$ (use that $v(\varphi_1)+v(\varphi_2)+v(\varphi_3)=k_1+k_2+k_3)$). 
Thus any $V$ with weights $0, 1, 2$, semisimple 
$\varphi$, $v(\varphi_i)=1$ for $i=1,2,3$, and with generic weight filtration, 
is weakly admissible, hence gives such an example.
%\item[(v)] Of course, for a "generic" $V$ (for example, any $V$ corresponding to a point in a Zariski-open subset of the set $Z$ of a $p$-adic family of refined representations in the sense of section~\ref{families}, see remark \ref{remdeffamraf}), the $d!$ refinements of $V$ are all non critical. 
\end{itemize}
\end{remark}

The following proposition is an immediate consequence of Proposition \ref{ref=tri}.
\begin{prop} 
\label{translation} 
$\Ref$ is non critical if, and only if, the sequence of
Hodge-Tate weights $(s_i)$ associated to $\Ref$ by 
Proposition \ref{ref=tri} is increasing, {\it i.e.} if $s_i=k_i$ $\forall i$.
\end{prop}

It is easy to see that non criticality is preserved under crystalline twists and duality. However, we have to be more careful with tensor operations, for even the notion of refinement is not well behaved with respect to tensor products. We content ourselves with the following trivial results,
that we state for later use.

\begin{lemma} \label{tensornoncrit}
\begin{itemize}
\item[(i)] Assume that $(D,\Ref)$ is a refined crystalline $\fg$-module over $\Ro_L$ with distinct Hodge-Tate weights. Then the weight of $\Lambda^i(\Ref_i) \subset \Dcf(\Lambda^i(D))$ is the smallest Hodge-Tate weight of $\Lambda^i(D)$.
\item[(ii)] Let $D_1$ and $D_2$ be two $\fg$-modules over $\Ro_L$ with
integral Hodge-Tate-Sen weights, equipped with a one dimensional $L$-vector space $W_i \subset \Dcf(D_i)$. If the weight of $W_i$ is the smallest integral weight of $D_i$ for $i=1$ and $2$, then the weight of $W_1 \otimes_L W_2$ is the smallest integral weight of $D_1 \otimes_{\Ro_L} D_2$.
\end{itemize}
\end{lemma}

\subsection{Deformations of non critically refined crystalline representations} \label{defdef} \saut

An essential feature of non critically refined crystalline representations is that they admit a nicer deformation theory. \ps

Let $(V,\cal{F})$ be a refined crystalline representation. Let us call $\Tc$ the triangulation of
$\Dr(V)$ corresponding to $\cal{F}$ as in Proposition~\ref{ref=tri}.
Recall from~\S\ref{deftri} the functors 
 $X_{V}=X_{\Dr(V)}: \cal{C} 
\longrightarrow {\rm Set}$ (resp. $X_{V,\Tc}=X_{\Dr(V),\Tc}$) parameterizing the 
deformations of $V$ (resp. the trianguline deformations of $(V,\Tc)$).
We shall use also the notation $X_{V,\cal{F}}$ for $X_{V,\Tc}$ and we call a trianguline 
deformation of $(V,\Tc)$ a trianguline deformation of $(V,\cal{F})$.

\subsubsection{A local and infinitesimal version of Coleman's classicity theorem.}

\begin{theorem} \label{colcrit} Let $V$ be a crystalline $L$-representation of $\Gp$ 
with distinct Hodge-Tate weights and such that $\Hom_{\Gp}(V,V(-1))=0$. 
Let $\Ref$ be a non critical refinement of $V$ and $V_A$ a trianguline deformation of $(V,\Ref)$. 
Then $V_A$ is Hodge-Tate if, and only if, $V_A$ is crystalline.
\end{theorem}

\begin{pf} Assume that $V_A$ is Hodge-Tate, we have to show that $\Dr(V_A)$ is crystalline by Proposition 
\ref{bergerthm}. By assumption, $\Dr(V_A) \in X_{\Dr(V),\Tc}(A)$ for the triangulation $\Tc$ of $\Dr(V)$ 
induced by $\Ref$, which has strictly increasing weights $s_i$ as $\Ref$ is non critical and by 
Proposition \ref{translation}. As $V_A$ is Hodge-Tate, $\Dr(V_A)$ satisfies by Proposition \ref{senpol} 
the hypothesis of Proposition \ref{criteriumDR}, hence the conclusion.
\end{pf}

This result has interesting global consequences, some of which will be explained in \S\ref{globcons} below. 
It is most useful when combined with the main result of the following paragraph, which gives a sufficient condition, {\it 
à la Kisin}, for a deformation to be trianguline.

\subsubsection{A criterion for a deformation of a non critically refined crystalline representation to be trianguline.}

\label{criteria} ${}^{}$ \ps Let $D$ be a $\fg$-module free over $\Ro_A$, 
we first give below a criterion to produce a $\fg$-submodule of rank $1$ over
$\Ro_A$. This part may be seen as an analogue of \cite[prop. 5.3]{colmeztri}. 

\begin{lemma} \label{prelimpoidscst} 
Let $D$ be a $\fg$-module free over $\Ro_A$, $\delta: \Q_p^* 
\longrightarrow A^*$ be a continuous character, and $\overline{\delta}:=\delta \bmod m$.
\begin{itemize}
\item[(i)] Let $M \subset D^{\delta}$ be
a free $A$-module of rank $1$. Then $\Ro_A[1/t]M$ is a $\fg$-submodule of $D[1/t]$
which is free of rank $1$ over $\Ro_A[1/t]$, and a direct summand.
\item[(ii)] Same assumption as in (i), but assume moreover that
$$\Im\left(M \longrightarrow
(D/mD)^{\bdelta}\right) \not\subset t(D/mD).$$ Then
$\Ro_A M$ is a $\fg$-submodule of $D$ which is free of rank $1$ over 
$\Ro_A$ and a direct summand.
\item[(iii)] Assume that $k \in \Z$ is the smallest integral root of the Sen polynomial of $D/mD$ and let $\lambda\in A^*$. Let 
$M \subset \Dcf(D)^{\varphi=\lambda}$ 
be free of rank $1$ over $A$ such that $$\Fil^{k+1}(M/mM)=0.$$ Then $\Ro t^{-k}M$ is a $\fg$-submodule of $D$ which is free of rank $1$ over 
$\Ro_A$ and a direct summand.
\end{itemize}
\end{lemma}

\begin{pf} The natural map $$\Ro \otimes_{\Q_p} M = \Ro_A \otimes_A M \longrightarrow \Ro_AM
\subset D$$ is injective by a standard argument as $(\fr(\Ro))^{\Gamma}=\Q_p$. As a
consequence, $\Ro_AM$ is free of rank $1$ over $\Ro_A$. In particular,
$\Ro_A[1/t]M$ is free of rank $1$ over $\Ro_A[1/t]$, hence a $\Ro_A[1/t]M$ direct summand
of $D[1/t]$ as $\Ro$-module by Proposition \ref{structurer}, and we conclude (i) by Lemma \ref{lib}
(i) and (ii). To prove (ii), it suffices by Lemma \ref{lib} (iii) to show that  
$\Im (\Ro_A M \longrightarrow D/mD)$, which is $\Ro_{A/m}(M/mM)$ by the
remark \ref{remfree}, is $\Ro$-saturated in $D/mD$. But this is the
assumption. \ps
For part (iii), write $M=Av$, $\varphi(v)=\lambda v$. Lemma \ref{mainlemma} shows that $v \in t^kD$ and that $ \overline{v} \not\in t^{k+1}D/mD$. Part (ii) above applied to $M':=t^{-k}M$ concludes the proof.
\end{pf}

\begin{remark} When $A=L$ is a field, a $\fg$-module $D$ over $\Ro_L$ 
is triangulable over $\Ro_L$ if and only if $D[1/t]$ is triangulable over $\Ro_L[1/t]$
(with the obvious definition). However, this is no longer true 
for a general $A$. 
\end{remark}

An immediate consequence of Lemma \ref{prelimpoidscst} (iii) is the following proposition, 
which could also be proved without the help of $\fg$-modules 
(see \S\ref{defdef} for the definition of $\cal{C}$).

\begin{prop}[The "constant weight lemma"] \label{poidsconstant} Let $V$ be an $L$-representation of $\Gp$ and $\lambda \in L^*$. Assume that $\Dc(V)^{\varphi=\lambda}$ has $L$-dimension $1$ and that its induced weight filtration admits the smallest integral Hodge-Tate weight $k$ of $V$ as jump. Let $A \in {\rm Ob}(\cal{C})$, $\lambda' \in A^*$ a lift of $\lambda$, and $V_A$ a deformation of $V$ such that $\Dc(V_A)^{\varphi=\lambda'}$ is free of rank $1$ over $A$. \ps
	Then the weight filtration on $\Dc(V)^{\varphi=\lambda'}$ has $k$ as unique jump. In other words, $k$ is a constant Hodge-Tate weight of $V_A$.
\end{prop}

We are now able to give a criterion on a deformation $V_A$ of a refined crystalline
representation $(V,\cal{F})$ ensuring that it is a trianguline deformation. We need the following definition.

\begin{definition} \label{reg} We say that the refinement $\Ref$ of $V$ is regular if the ordering $(\varphi_1, \cdots, \varphi_d)$ of the 
eigenvalues of $\varphi$ it defines has the property
\begin{center} (REG) for all $1\leq i\leq d$, $\varphi_1\varphi_2\cdots\varphi_i$ 
is a simple eigenvalue of $\Lambda^i(\varphi)$.
\end{center}
In particular, the $\varphi_i$ are distinct, and $\Ref$ is the unique refinement such that
$\Lambda^i(\Ref_i)=(\Lambda^i(\Dc(V)))^{\varphi=\varphi_1\cdots\varphi_i}$.
\end{definition}

The next theorem is the bridge between this section and section
\ref{kisin}. 

\begin{theorem} \label{poidscst2} Assume that $\Ref=\{\varphi_1,\cdots,\varphi_d\}$ 
is a non critical, regular, refinement of $V$. 
Let $V_A$ be a deformation of $V$, and assume that we are given a continuous
homomorphism $\delta=(\delta_i)_{i=1,\cdots,d}: \Q_p^* \longrightarrow A^*$
such that for all $i$, \begin{itemize} 
\item[(i)] $\Dc(\Lambda^i(V_A)(\delta_1\cdots\delta_i)_{|\Gamma}^{-1})^{\varphi=\delta_1(p)\cdots \delta_i(p)}$
is free of rank $1$ over $A$. 
\item[(ii)] ${\delta_i}_{|\Gamma} \bmod m = \chi^{-t_i}$ for some 
$t_i \in \Z$, and
$\varphi_i=(\delta_i(p) \bmod m)p^{t_i}$.
\end{itemize}
Then $V_A$ is a trianguline deformation of $(V,\Ref)$ whose parameter is 
$(\delta_i x^{t_i-k_i})_{i=1,\dots,d}$.
\end{theorem}

\begin{pf} Note that the assumptions and conclusions do not change 
if we replace each $\delta_i$ by $\delta_i x^{m_i}$ for $m_i \in \Z$. Thus we can assume that $t_i=0$ for all $i$, {\it i.e.} ${\delta_i}_{|\Gamma} \equiv 1 \bmod m$.\ps
Fix $1\leq i \leq d$ and set 
$W_i:=\Lambda^i(V)\otimes_A A(\delta_1\cdots\delta_i)^{-1}_{|\Gamma}$. By
assumption (ii), $W_i/m W_i$ is crystalline with smallest 
Hodge-Tate weight $w_i:=k_1+\cdots+k_i$. Moreover, by the 
regularity property of $\Ref$, 
$\Lambda^i(\Ref_i)=\Dc(W_i/mW_i)^{\varphi=\varphi_1\cdots\varphi_i}$. 
As $\Ref$ is non critical, 
$\Lambda^i(\Ref_i) \cap \Fil^{w_i+1}(\Dc(W_i/mW_i))=0$. 
Lemma \ref{prelimpoidscst} (iii) shows that 
$$t^{-w_1}\Ro_L\Dcf(W_i)^{\varphi=\delta_1(p)\cdots\delta_i(p)} \subset 
\Dr(W_i)$$
is free of rank one over $\Ro_A$ and direct summand. 
As a consequence, 
$\Lambda^i(\Dr(V_A))=\Dr(W_i(\delta_1\cdots\delta_i)_{|\Gamma})$ 
contains a $\fg$-submodule free of rank $1$ and 
direct summand as $\Ro_A$-module, and isomorphic to 
$\Ro_A(\delta_1\cdots\delta_i x^{-w_i})$. 
It has then the form $\Lambda^i(D_i)$ for a unique 
$\fg$-submodule $D_i \subset \Dr(V_A)$ free of rank $i$ and direct summand as $\Ro_A$-module. 
By construction, 
\begin{equation} \label{defdi} \Lambda^i(D_i)\isomo \Ro_A(\delta_1\cdots\delta_i x^{-w_i}) 
\end{equation}
By (ii) and the regularity of $\Ref$, 
$D_i/mD_i$ is the unique a saturated $\fg$-submodule 
$X \subset \Dr(V)$ such that $\Lambda^i(\Dcf(X))=\Ref_i$, 
hence 
\begin{equation} \label{dires} D_i/mD_i= \Fil_i(\Dr(V)). \end{equation} \ps 
To conclude, it remains only to show that $D_i \subset D_{i+1}$ 
if $0 \leq i \leq d-1$ (set $D_0:=0$). Indeed, the formula (\ref{defdi}) forces then the parameter of 
the triangulation to be $(\delta_ix^{-k_i})$. \ps
Let us assume by induction on $0 \leq j \leq d$ that $D_j \subset D_i$ for all
$j < i \leq d$ (this is obvious when $j=0$). Let $j+1 < i \leq d$, 
we know that $D_j \subset
D_i$, and we want to prove that $D_{j+1}\subset D_i$. Consider the left-exact functor 
$F(X):=X^{\delta_{j+1}x^{-k_{j+1}}}$. Note that for
$i> j+1$, 
$$F(\Ro_L(\delta_ix^{-k_i}))=H^0(\Ro_L(\delta_i\delta_{j+1}^{-1}x^{k_{j+1}-k_i}))=0$$
by assumption (ii) and \cite[prop. 3.1]{colmeztri}, as 
$\varphi_i\neq \varphi_{j+1}$. As a consequence, (\ref{dires}) and a devissage 
show that for $s \geq j+1$, 
$$F((D/D_s)\otimes_A A/m)=F(\Dr(V)/\Fil_s(\Dr(V)))=0,$$ and
another devissage shows that 
$F(D/D_s)=0$. Consider now the exact sequence
\begin{equation}\label{exseq} 0 \longrightarrow D_s/D_j \longrightarrow D/D_j \longrightarrow
D/D_s \longrightarrow 0.\end{equation}
The induction hypothesis shows that
$\Ro_A(\delta_{j+1}x^{-k_{j+1}})\isomo D_{j+1}/D_j \subset D/D_j$, hence
applying $F$ to the sequence (\ref{exseq}) when $s=j+1$ we get that
the inclusion 
$$F(D_{j+1}/D_j) \longrightarrow F(D/D_j)$$ 
is an equality. Applying now $F$ to
the exact sequence (\ref{exseq}) for $s=i$, we get that inside $D/D_j$, we have 
$F(D_i/D_j)=F(D_{j+1}/D_j) \subset D_i/D_j$. But $D_{j+1}/D_j\isomo 
\Ro_A(\delta_{j+1}x^{-k_{j+1}})$ is generated as
$\Ro$-module by $F(D_{j+1}/D_j)$, hence $D_{j+1}/D_j \subset D_i/D_j$. 
\end{pf}

\subsubsection{Properties of the deformation functor $X_{V,\cal{F}}$}

In fact, we can in many cases describe quite simply $X_{V,\Ref}$ when $\Ref$ is non critical. 
The following results will not be needed in the remaining sections, 
but are interesting in their own. 
Recall that by definition we have a natural  transformation $$X_{V,\Ref} \longrightarrow X_V.$$

\begin{prop} \label{basic} Assume that the eigenvalues of $\varphi$ on $\Dc(V)$ are distinct, 
then $X_{V,\Ref}$ is a subfunctor of $X_V$ and $X_{V,\Ref}\longrightarrow X_V$ is relatively representable. 
Moreover, if $\Ref$ is non critical, the subfunctor $X_{V,\rm crys} \subset X_V$ of crystalline deformations factors through $X_{\Ref}$.
\end{prop}

\begin{pf} As the eigenvalues of $\varphi$ are distinct, 
the characters $\delta_i$ of the parameter $\delta$ of $\Dr(V)$ associated to $\cal{F}$ satisfy
$\delta_i \delta_j^{-1} \not \in x^\Z$ for $i \neq j$ (see Prop.~\ref{ref=tri}). The first sentence thus follows
from Prop.~\ref{subfunctor}, Prop.~\ref{representri} using Prop.~\ref{XVXD}.

\ps
Assume that $\Ref$ is non critical and let $V_A$ be a crystalline deformation of $V$. 
We have to show that $D_A:=\Dr(V_A)$ admits a (necessarily unique) triangulation 
lifting the one associated to $\Ref$. 
As the $\varphi_i$ are distinct, the characteristic polynomial of $\varphi$ on 
$\Dcf(D_A)$ writes uniquely as $\prod_i(T-\lambda_i) \in A[T]$ with $\lambda_i \equiv \varphi_i \bmod m$. As $V_A$ is Hodge-Tate 
with smallest Hodge-Tate weight $k_1$, and as $\Ref$ is non critical, 
Lemma \ref{prelimpoidscst} (iv) shows that 
$$\Ro t^{-k_1}\Dcf(D_A)^{\varphi=\lambda_1}$$
is a submodule of $D_A$ which is a direct summand as $\Ro_A$-module. We construct this way by induction the required triangulation of $D_A$. \ps
%If $A \subset A' \in \cal{C}$ and $D \otimes_A A'$ is a refined deformation we have to show that $D_A$ is also a refined deformation. We remark first that the Sen polynomial is in $A[T]$, hence so are the parameters ${\delta'_i}_{|\Gamma}$. By twisting we may assume that ${\delta'_1}_{|\Gamma}=1$
\end{pf}

The main theorem concerning non critical refinements is then the following, which may be viewed as a 
$d$-dimensional generalization of the results of some computations of Kisin in \cite{kis} \S7 
(giving a different proof of his results when $d=2$).

\begin{theorem} 
\label{thmdefnoncrit} 
Assume that $\cal{F}$ is non critical, 
that $\varphi_i\varphi_j^{-1} \not \in \{1,p^{-1}\}$ if $i < j$, and that 
$\Hom_{\Gp}(V,V(-1))=0$. Then $X_{V,\cal{F}}$ 
is formally smooth of dimension $\frac{d(d+1)}{2}+1$. Moreover, the
{\it parameter} 
map induces an exact sequence of $L$-vector spaces:
$$ 0 \longrightarrow X_{V,\rm crys}(L[\varepsilon]) \longrightarrow
X_{V,\cal{F}}(L[\varepsilon]) \longrightarrow \Hom(\Z_p^*,L^d) \longrightarrow
0.$$
\end{theorem}

\begin{pf} Let $(\delta_i)$ be the parameter of $(V,\Ref)$. If $i \neq j$, then 
$\delta_i\delta_j^{-1} \not\in x^{\Z}$ since $\varphi_i \neq \varphi_j$. Moreover if $i < j$, then $k_i < k_j$ and 
$\delta_i\delta_j^{-1} \not\in \chi x^{\N}$ as $\varphi_i \neq 
p^{-1}\varphi_j$ by assumption. The result follows then from Propositions \ref{deform},\ref{translation},
\ref{criteriumDR} and \ref{basic}.\end{pf}

\index{L@$L$, the residue field of the artinian local ring $A$, 
a finite extension of $\Q_p$|)}
\index{Aa@$A$, a local artinian ring, finite over $\Q_p$|)}

\subsection{Some remarks on global applications} \label{globcons} \saut

We now derive some consequences of these results 
in a global situation. \ps

Let $V$ be a finite
dimensional $L$-vector space equipped with a geometric continuous representation of
$\Gal(\overline{\Q}/\Q)$ and assume that $V_p:=V_{|\Gp}$ is crystalline, with distinct Hodge-Tate weights and distinct Frobenius 
eigenvalues. 

Let $\cal{F}$ be a refinement of $V_p$, and recall that $X_{V_p,\cal{F}}$ denote the trianguline deformation
functor of $(V_p,\cal{F})$. Let $X_{V,\cal{F}}$ 
denote the subfunctor of the full deformation functor of $V$ consisting of
the deformations whose
restriction at $p$ is in $X_{V_p,\cal{F}}$ (that is, trianguline), and whose restriction at $l \neq p$ satisfies the usual
{\it finite} condition
(for example, are unramified at $l$ if $V$ is). Then Theorems~\ref{colcrit} and~\ref{thmdefnoncrit} imply :

\begin{cor} \label{selmerglob} If $\cal{F}$ is non critical, then there is a natural exact sequence 
$$0 \longrightarrow H^1_f(\Q,{\rm ad}(V)) \longrightarrow X_{V,\cal{F}}(L[\varepsilon])
\overset{\kappa}{\longrightarrow} \Hom(\Z_p^*,L).$$
In particular, if $H^1_f(\Q,{\rm ad}(V))=0$ (which is conjectured to be
the case if $V$ is absolutely irreducible), then $\dim_L(X_{V,\cal{F}}(L[\varepsilon])\leq \dim_L(V)$.
\end{cor}

In this setting, the question of the determination of 
$\dim_L(X_{V,\cal{F}}(L[\varepsilon]))$ seems to be quite subtle, even conjecturally. 
Among many other things, it is linked to the local dimension of the eigenvarieties of $\GL_d$, which are still quite mysterious 
(see the work of Ash-Stevens \cite{ashstevens} and of M. Emerton \cite{eminterp}). 

However, there are similar questions for which the theory of $p$-adic families of automorphic forms suggests a nice answer\footnote{We consider here a simplified setting, an appropriate condition on the Mumford-Tate group of $V$ should suffice in general.}. 
As an example, let us consider now an analogous case where $V$ is an irreducible, $d$-dimensional, geometric
$L$-representation of $\Gal(\overline{E}/E)$, $E/\Q$ a quadratic imaginary
field, satisfying $V^{c,*} \simeq V(d-1)$. Assume that $p=vv'$ splits in $E$,
fix an identification $\Gp \isomo \Gal(\overline{E_v}/E_v)$, and assume
that $V_p:=(V_{|\Gp},\cal{F})$ is crystalline and provided with a refinement $\cal{F}$, with distinct
Hodge-Tate weights. Let 
$X_{V,\cal{F}}$ denote the subfunctor of the full deformation 
functor of $V$ consisting of deformations whose
restriction at $v$ is in $X_{V_p,\cal{F}}$, satisfying ${V_A}^{c,*}=V_A(d-1)$ and the $f$ condition outside $p$.  \ps \medskip

\noindent {\bf Conjecture}:  
Assume that $\cal{F}$ is non critical, then 
$X_{V,\cal{F}}$ is prorepresented by $${\rm Spf}(L[[X_1,\cdots,X_d]])$$ 
and $\kappa$ is an isomorphism. \ps \medskip

In the subsequent paragraph \S\ref{eigenatregcripts}, we will give more details about the proofs of the facts alluded here and 
we will explain how we can 
deduce this conjecture in many cases
%\footnote{At least under the " $U$-non criticality assumption", see Remark~\ref{remcrit} (iii).} 
assuming the conjectured vanishing of $H^1_f(E,\ad(V))$, and using
freely the results predicted by Langlands philosophy on the correspondence between automorphic forms for suitable 
unitary groups $G$ (in $d$ variables) 
attached to the quadratic extension $E/\Q$. As we will explain, we can even get an "$R=T$" statement for 
${\rm Spf}(R)=X_{\Ref}$ and $T$ the completion of a well chosen eigenvariety of $G$ at the point corresponding to $(V,\Ref)$. \ps

To sum up, the eigenvariety of $G$ at {\it irreducible, classical, non-critical} points should be smooth, and
neatly related to deformation theory. 
By contrast, a much more complicated (but interesting) situation is expected at 
reducible, critical points, and this is the main object of subsequent sections of this paper.    
 \ps

\newpage

\section{Generalization of a result of Kisin on crystalline periods}

\label{kisin}

\subsection{Introduction}

In this section, we solve, generalizing earlier results of Kisin,
some questions of ``Fontaine's theory in family'' concerning continuation
of crystalline periods.

Let $X$\index{X@$X$, a reduced rigid analytic space over $\Q_p$} be a reduced rigid analytic space over $\Q_p$ and 
$\MM$ \index{Mt@$\MM$, a coherent torsion-free sheaf of $\anneau_X$-modules with a continuous action of $G_p$}
 a family of $p$-adic representations of $G_p=\Gal(\overline{\Q}_p/\Q_p)$ over $X$, that is,
in this section, a coherent torsion-free sheaf of $\anneau_X$-modules
with a continuous action of the group  $G_p$\index{Oy@$\anneau_X$, or $\OO$, the structural sheaf of the rigid analytic space $X$}. Note that we do not assume that $\MM$ is locally free.\footnote{Indeed, we will apply the results of this section to modules
associated to pseudocharacters in the neighborhood of a reducible point,
and it is one of the main results of section~\ref{pseudocharacters} that they
 do not in general come from representations over free modules.}
For each point $x \in X$ of residue  field $k(x)$ \index{kx@$k(x)$, a finite 
$\Q_p$-algebra, the residue field of $X$ at $x$}
the $k(x)$-vector space $\MB_x$ \index{Mxb@$\MB_x$, a finite dimensional representation of $G_p$ (or of $G$) over $k(x)$, the stalk of $\MM$ at $x$}\index{Mx@$\MM_x$, the $\OO_x$-module of germs of sections of $\MM$ at $x$.} is then a continuous
representation of $G_p$, to which we can apply the $p$-adic Hodge theory of
Tate and Sen and all their generalizations by Fontaine.
The questions concerning ``Hodge-Tate theory in family'' were
completely solved by Sen :
in particular he shows in our context that there exist $d$ analytic functions
$\kappa_1,\dots,\kappa_d$ on $X$, where $d$ is the generic rank of $\MM$,
such that $\kappa_1(x),\dots,\kappa_d(x)$ are the Hodge-Tate-Sen weights of
$\MB_x$ for a Zariski-dense open set of $X$.
 We shall assume in this introduction, to simplify the discussion,
that $\kappa_1=0$. We shall also assume that in our family
the other weight functions $\kappa_2,\dots,\kappa_d$ move widely
(in a technical sense we do not want to make precise here,
but see~\ref{directhyp} below), as it happens for families supported by
eigenvarieties. In particular, the families we work with are quite different from the
families with constant weights studied by Berger and Colmez.

Suppose we know that $\MB_z$ is crystalline with positive Hodge-Tate weights
for a Zariski dense
subset $Z$ \index{Z@$Z$, a Zariski dense subset of $X$} of points of $X$ and that for all $z \in Z$ it has a crystalline
 period that is an eigenvector of the crystalline Frobenius $\varphi$
with eigenvalue $F(z)$, $F$ being a fixed 
rigid analytic function on $X$. 
In other words, assume that $\Dc(\MB_z)^{\varphi=F(z)}$ is non zero for
$z \in Z$. Can we deduce from this that
\begin{itemize}
\item[(1)] for each $x$ in $X$,
$\MB_x$ has a crystalline period, which is moreover an eigenvector for $\varphi$ with eigenvalue $F(x)$ ?
\end{itemize}
Or, more generally, that
\begin{itemize}
\item[(2)] for each $x \in X$, and $\spec A$ a thickening of $x$,
(i.e. $A$ an artinian quotient of the rigid analytic local ring $\anneau_x$ of $x$ at $X$) 
$\MM$ has a non-torsion crystalline period over $A$ which is
an eigenvector for $\varphi$ with eigenvalue the image $\bar F$ of $F$ in $A$
?
In other words, is that true that $\Dc(\MM \otimes A)^{\varphi=\bar F}$ has
a free of rank one $A$-submodule ?
\end{itemize}

Kisin was the first to deal with those questions
and most of his works in~\cite{kis} is an attempt to answer them in the
case where $\MM$ is a free $\anneau_X$-module. Under this freeness
assumption
(plus some mild technical hypothesis on $Z$ we will not state neither
mention further in this introduction),
he also proves many cases of question (2), although his results are
widespread along his paper and sometimes not explicitly stated.
If we collect them all, we get
that Kisin proved that question (2) has a positive answer (when $\MM$
is free) for those
$x$ that satisfy two conditions :
\begin{itemize}
\item[(a)] The representation $\MB_x$ is indecomposable,
\item[(b)] $\Dc(\MB_x^{\ses})^{\varphi=F(x)}$ has dimension $1$.
\end{itemize}

Condition (b) is probably necessary. But condition (a) is not, and
appears because of the use by Kisin of some universal deformation
arguments. In \S\ref{kisindirect}, using mostly arguments of Kisin, but
simplifying and reordering them, we prove that when $\MM$ is a locally free
module, question (2) has a positive answer for all $x$ satisfying the
condition b) above. We hope that our redaction may clarify for many
readers the beautiful and important results of Kisin.

But our main concern here is to generalize those results to the case of
an arbitrary torsion-free coherent sheaf $\MM$.
We are able to prove that question (2) (hence also question (1)) still
has a positive answer in this case provided that $x$ satisfies hypothesis
(b) above. This is done in \S\ref{kisinindirect}.

Let us now explain the idea of the proof : basically we do a reduction to
the case where $\MM$ is locally free.
To do this we use a rigid analytic version of a ``flatification'' result of Gruson-Raynaud which gives
us a blow-up $X'$ of $X$ such that the {\it
strict transform} $\MM'$ of $\MM$ on $X'$ is locally free. Hence we know the
(positive) answers to questions (1) and (2) for $\MM'$ and the problem is to
``descend'' them to $\MM$. This is the aim of \S\ref{descenteblowup}.

For this the difficulties are twofolds. The first difficulty is that
if $x'$ is a point of $X'$ above $x$ (let us say to
fix ideas with the  same residue field, since a
field extension here would not harm)
then $\MB'_{x'}$ is not isomorphic to $\MB_{x}$ but to a quotient of
it. Since the functor $\Dc(-)^{\varphi=F}$ is only left-exact, the
positive answer to question (1) for $x'$ does not imply directly the
positive answer for $x$ -- and of course, neither for question (2).
 The second difficulty arises only when dealing with
question (2) : it is not possible in general to lift the
thickening $\spec(A)$ of $x$ in $X$ to a thickening of $x'$ in $X'$,
whatever $x'$ above $x$ we may choose\footnote{The reader may convince
himself of this assertion by looking
at the case where $X=\spec L[[T^2,T^3]]$ is the cusp and
$X'=\spec(L[[T]])$ the blow up of $X$ at its maximal ideal (that is, its normalization). The principal ideal $T^2L[[T^2,T^3]]$ has not the form $L[[T^2,T^3]]\cap
T^nL[[T]]$ for $n\geq 0$, hence $A=L[[T^2,T^3]]/(T^2)$ is a counter-example. }. So the direct strategy of
descending a positive answer to question (2) from a $\spec(A)$ in $X'$
to $\spec(A)$ in $X$ can not work.
To circumvent the second difficulty, we use a lemma of Chevalley to
construct a suitable thickening $\spec(A')$ of $\spec(A)$ in $X'$,
and then some rather
involved arguments of lengths to deal with the first one as well as the
difference between $\spec(A')$ and $\spec(A)$. 
As Chevalley's lemma requires to work at the level of complete noetherian ring, and as we have to use rigid analytic local rings when dealing with interpolation arguments, we need also at some step of the proof to compare various diagrams with their completion. For all these reasons, the total argument in \S\ref{descenteblowup} is rather long.

Finally let us say that the idea of using a blow-up was already
present in Kisin's argument in the free case\footnote{There Kisin does not
use
the blow-up to make $\MM$ free, since it already is, but instead to make an
ideal of crystalline periods locally principal. He does not prove a direct descent
result as ours, using instead a comparison of universal deformation rings.},
and is still present in the locally free case in \S\ref{kisindirect}.
This is why our descent result of \S\ref{descenteblowup} is used twice, once 
in \S\ref{kisindirect} and
once in \S\ref{kisinindirect}. However, were it to be used only in the locally
free case, the descent method could be much simpler\footnote{The first
difficulty above vanishes, since in that case the strict transform of $\MM$
is simply its pull-back, and the second may be dealt with much more
easily.}.
\subsection{A formal result on descent by blow-up}

\subsubsection{Notations.} \label{notationsrigide} Let $X$ be a reduced, separated, rigid analytic space over $\Q_p$, $\OO_X$ (or simply $\OO$) its structural sheaf, and 
let $\MM$ be a coherent $\OO$-module on $X$. For $x$ a point of $X$, we shall note $\anneau_{x}$ \index{Ox@$\OO_x$, the rigid analytic local ring of $X$ at $x$}
the rigid analytic
local ring of $X$ at $x$, $m_x$ its maximal ideal, and 
$k(x)=\OO_x/m_x$ its residue field. Moreover, we denote by
$\MM_{x}$ the rigid analytic germ of $\MM$ at $x$,  that is 
$\MM_x = \MM(U) \otimes_{\anneau(U)} \anneau_x$
where $U$ is any open affinoid containing $x$, and by 
$\overline{\MM}_x:=\MM_x \otimes_{\anneau_x} k(x)$ the fiber of $\MM$ at $x$. \ps

Let $G$ be a topological group and assume that $\MM$ is equipped with a continuous $\OO$-linear action of $G$. This 
means that for each open affinoid $U \subset X$, we have a continuous morphism $G \longrightarrow \Aut_{\OO(U)}(\MM(U))$, whose formation is compatible with the restriction to any open affinoid $V \subset U$. For $x$ a point of $X$, 
$\MM_x$ (resp. $\MB_x$) is then a continuous $\OO_x[G]$-module (resp. $k(x)[G]$-module) in a natural way. \ps

\begin{remark}\label{remtorsionfreegen}({\it On torsion free modules}) In this section and the subsequent ones, we will sometimes have to work with torsion free modules. Recall that a module $M$ over a {\it reduced ring} $A$ is said to be torsion free if the natural map $M \longrightarrow M \otimes_A K$ is injective where $K={\rm Frac}(A)$ is the total fraction ring of $A$ (see \S\ref{reducedcase}). \par
	If $X$ is a reduced affinoid and $\cal M$ a coherent $\OO_X$-module, then $\MM(X)$ is torsion free over $\OO(X)$ if, and only if, $\MM_x$ is torsion free over $\OO_x$ for all $x \in X$. Indeed, this follows at once from the faithful flatness of the maps $\OO(X)_x \longrightarrow \OO_x$ and the following lemma.\end{remark}
	%Recall also that for an affinoid $X$, $\OO(X)$ is reduced if and only if $\OO_x$ is reduced for all $x \in X$ by \cite[]{BGR}, in which case we say that $X$ is reduced. 
\begin{lemma}\label{gentorsionfree} Let $A$ be a reduced noetherian ring and $M$ a 
finite type $A$-module. The following properties are equivalent:
\begin{itemize}
\item[(i)] $M$ is torsion free over $A$,
\item[(ii)] $M$ is a submodule of a $K$-module,
\item[(iii)] $M$ is a submodule of a finite free $A$-module,
\item[(iv)] $M_x$ is torsion free over $A_x$ for all $x \in {\rm Specmax}(A)$,
\item[(v)] there is a faithfully flat $A$-algebra $B$ such that $M\otimes_A B$ is a $B$-submodule of a finite free $B$-module.
\end{itemize}
\end{lemma}

\begin{pf} It is clear that (i), (ii) and (iii) are equivalent (for (ii) $\Rightarrow$ (iii) note that any $K$-module embeds into a free $K$-module as $K$ is a finite product of fields). The equivalence between (i) and (iv) follows now from the injection $M \longrightarrow \oplus_{x \in {\rm Specmax A}}M_x$, and the fact that ${\rm Frac}(A_x)$ is a factor ring of $K$: namely the product of the fraction fields of the irreducible component of ${\rm Spec}(A)$ containing $x$. \par Note that condition (iii) is equivalent to ask that the natural map $M \longrightarrow \Hom_A(\Hom_A(M,A),A)$ is injective. But this can be checked after any faithfully flat extension $B$ of $A$ as the formation of the $\Hom$'s commute with any flat base change when the source is finitely presented, thus (i) $\Leftrightarrow$ (v). 
\end{pf}

\subsubsection{The left-exact functor $D$.} Fix a point $x \in X$. Let $D$ be an additive left-exact functor from the (artinian) category of finite length, continuous, $\OO_x[G]$-modules, to the category of finite length $\OO_x$-modules. Here are some interesting examples: \begin{itemize}
\item[(i)] Let $G:=\Gp=\Gal(\Qpb/\Qp)$ and let $B$ be any topological $\Q_p$-algebra equipped with a continuous action of $\Gp$. Assume that $B$ is $\Gp$-regular in the sense of Fontaine \cite[Expos\'e III, \S1.4]{PP}. For $Q$ any $\anneau_{x}$-module of finite
length equipped with a continuous $G_p$-action (hence a finite dimensional $\Q_p$-representation of $\Gp$), we let 
$$D(Q):=(Q \otimes_{\Qp} B)^G.$$
The functor $D$ satisfies our assumptions by {\it loc. cit.} . As an example, we can take $B=\Q_p$ or $B=\C_p$ above. \ps
\item[(ii)] Fix $F \in \OO_x^*$. For any $Q$ as above, then 
$$D(Q):= \Dcp(Q)^{\varphi=F} =\{v \in (Q \otimes_{\Q_p} 
\Bcp)^{G_p}, \varphi(v)=Fv\},$$
where $\Bcp$ is the subring of $\Bc$ defined by Fontaine in \cite[Expos\'e II, \S2.3]{PP}, satisfies again our assumptions. 
\end{itemize}

In the sequel, we will be mainly interested in the case (ii) above.

\subsubsection{Statement of the result.}\label{subsubsttres} Assume that $\MM_x$ is torsion 
free over $\OO_x$ (recall that $\OO_x$ is reduced). Let $\pi: X'\longrightarrow X$ be a proper and birational 
morphism of rigid spaces with $X'$ reduced. Here {\it birational} means that for some coherent sheaf of ideals 
$H \subset \OO_X$, $U:=X-V(H)$ is Zariski dense in $X$ (where $V(H)$ is the closed subspace defined by $H$), 
$\pi$ is an isomorphism over $U$, and 
$\pi^{-1}(U) \subset X'$ is Zariski dense in $X'$. As an important example, we may take for $\pi$ the
blow-up\footnote{We refer to \cite[\S5.1]{GR} for the basics on blow-ups and to
\cite[\S2.3, \S4.1]{conradqcoh} for the notion of relative Spec and blow-ups in the context of rigid geometry.} of $H$. 
Let $\MM'$ be the strict transform of $\MM$ by this morphism (see below). 

\begin{prop} \label{descenteblowup}
Assume that for all $x' \in \pi^{-1}(x)$ and for every ideal
 $I'$ of $\anneau_{x'}$ of cofinite length, we have 
$$l(D(\MM'_{x'} \otimes \anneau_{x'}/I')) = l (\anneau_{x'}/I').$$
Assume moreover that $$l(D(\MB_x^\ses)) \leq 1$$
Then we also have, for every ideal $I$ of cofinite length of 
$\anneau_{x}$ :
$$l(D (\MM_{x} \otimes \anneau_{x}/I))=l(\anneau_{x}/I).$$
\end{prop}

\begin{remark} \begin{itemize}\label{generalitepropdescente}
\item[(i)] More precisely, we show that Proposition \ref{descenteblowup} holds when we replace the assumption of $\MB_x^{\ses}$ by the following slightly more general one: 
for any $k(x)[G]$-quotient $U$ of $\MB_x$, $l(D(U))\leq 1$ (see the proof of Lemma \ref{abstract}, which is the only place where the assumption is used).
\item[(ii)] As will be clear to the reader, the analogue of Proposition \ref{descenteblowup} in the context of schemes instead of rigid analytic spaces would hold by the same proof.
\end{itemize}
\end{remark}

This whole subsection is devoted to the proof of the proposition. Let us fix a coherent sheaf of ideals $H \subset \OO_X$ 
such that $U:=X-V(H)$ is Zariski dense in $X$, that 
$\pi$ is an isomorphism over $U$, and that $\pi^{-1}(U) \subset X'$ is Zariski dense in $X'$. Let us first recall how the strict transform $\mathcal M'$ of a coherent $\OO_X$-module is defined: 
it is a coherent $\OO_{X'}$-module which is 
locally the quotient of the coherent sheaf $\pi^\ast \mathcal M$ by its submodule of sections whose support
is in the fiber of $\pi$ over $V(H) \subset X$. In other words, if $H'$ is a coherent sheaf of ideals 
of $\OO_{X'}$ defining the closed subset $\pi^{-1}(V(H)) \subset X'$, then $\MM'$ is the quotient of $\pi^*\MM$ by its ${H'}^{\infty}$-torsion. Note that it depends on the choice of 
$H$ in general. This description makes clear that the action of $G$ on the $\OO_X$-module 
$\mathcal M$ defines an $\OO_{X'}$-linear continuous action of $G$ on $\mathcal M'$, and that the natural map 
$\pi^\ast\mathcal M \longrightarrow \mathcal M'$ is $G$-equivariant. 
A useful fact about the notion of strict transform is that the subsheaf of 
${H'}^{\infty}$-torsion of $\pi^*\MM$ is precisely the kernel of the natural 
morphism\footnote{We are grateful to Brian Conrad for pointing this to us. 
Here is the general statement: if $S$ is a rigid space, $I \subset \OO_S$ a 
coherent sheaf of ideals, 
$j: U:=S-V(I) \hookrightarrow S$ the inclusion of the complement of 
$V(I)$ and $\cal F$ a coherent $\OO_S$-module, then the $I^{\infty}$-torsion 
of $\cal F$ is the kernel of the natural 
map $\cal F \rightarrow j_*{\cal F}_{|U}$. Indeed, we may assume that $S$ is affinoid. Set $F=\cal F(S)$ and take $m \in F$ 
such that $m_s=0 \in F \otimes \OO_s$ for all $s \in U$ and we want to show that $m$ is killed by a power of $I(S)$. The faithfull 
flatness of $\OO(S)_s \rightarrow \OO_s$ shows that the closed points of the support of $m$ lie in $V(I(S))$, and we conclude as $\OO(S)$ is a Jacobson ring. } $\pi^*\MM \longrightarrow j_*(\pi^*\MM_{|\pi^{-1}(U)})$. 
As a simple application, if $\MM$ is torsion free then $\MM'$ is torsion free as well and does not 
depend on the choice of $H$ as above.

Since $\MM_x$ is torsion free over $\OO_x$,  
it can be embedded in a free of finite rank $\OO_x$-module, so we can choose an injection $$i: \MM_x \longrightarrow \anneau_x^n.$$ Fix $x' \in \pi^{-1}(x)$, $i$ induces a 
morphism $i':\, \MM_x \otimes_{\anneau_x}\anneau_{x'} \longrightarrow \anneau_{x'}^n$.
We check easily using the aforementionned useful fact that the kernel of $i'$ is the submodule of $\MM_x \otimes_{\anneau_x}\anneau_{x'}$ whose elements are killed by a power of $H'_{x'}$ so the image of $i'$ is $\MM'_{x'}$.
%\begin{pf} Set $f:=f_{x'}$. The kernel of $i'$ certainly contains the elements killed by a power of $f$ since 
%$f$ is not a divisor of $\anneau_{x'}$. On the other hand, if 
%an element $v=\sum_j m_j \otimes a_j$ is in the kernel of $i'$, then $i'(v)=\sum_j\, a_j \, \, i(m_j)$ is zero in $\anneau_{x'}^n$. 
%There is an integer $s\geq 1$ such that $f^s a_j \in \anneau_x$ for each $j$, hence $f^s v = (\sum_j f^s a_j m_j) \otimes 1$. But $$i(\sum_j f^s a_j m_j)= \sum_j (f^s a_j) i(m_j)= \sum_j f^s a_j
%i(m_j)=0$$ in $\anneau_x^n \subset \anneau_{x'}^n$. Since $i$ is injective, $f^s v=0$ and we are done. The second assertion follows from 
%the first one by definition.
%\end{pf}
%
We thus have a commutative diagram of $\anneau_x$-modules 
(and even of $\anneau_{x'}$-modules for the half right of the diagram)
\begin{eqnarray} 
\label{diag} \xymatrix{\MM_x \ar[r] \ar@{^{(}->}[d]^i & \MM_x \otimes_{\OO_x} \anneau_{x'} \ar@{->>}[r] \ar[rd]^{i'} \ar[r] & \MM'_{x'} \ar@{^{(}->}[d] \\ 
\anneau_x^n \ar[rr] & & \anneau_{x'}^n.}
\end{eqnarray}

We call $\hat{\anneau}_x$ (resp. $\hat{\anneau}_{x'}$) the completion of the local ring $\anneau_x$ (resp. of $\anneau_{x'}$) for the $\m_{x}$-adic (resp. 
$\m_{x'}$-adic) topology, and $\hat{\MM}_x = \MM_x \otimes_{\anneau_x} \hat{\anneau}_x$ (resp. $\hat{\MM}'_{x'}=\MM'_{x'} \otimes_{\anneau_{x'}} 
\hat{\anneau}_{x'}$) the completion of $\MM_x$ (resp. of $\MM'_{x'}$).

As $\anneau_x \longrightarrow \anneau_{x'}$ is a local morphism, it is continuous for the $\m_x$-adic topology at the source
and the $\m_{x'}$-adic topology at the goal. This is also true for any 
morphism form a finite type $\anneau_x$-module to a finite type 
$\anneau_{x'}$-module.
Hence such a morphism can be extended in a unique continuous way to their 
completion. We get this way morphisms 
$\hat{\anneau}_x \longrightarrow \hat{\anneau}_{x'}$ and 
$$\hat{\MM}_x \longrightarrow \hat{\MM_x\otimes_{\anneau_x}{\anneau_{x'}}}
= (\MM_x \otimes_{\anneau_x} \anneau_{x'}) 
\otimes_{\anneau_{x'}} \hat{\anneau}_{x'} = \hat{\MM}_x \otimes_{\hat{\anneau}_{x}} \hat{\anneau}_{x'}$$ the last equality being obtained by applying twice
the transitivity of the tensor product. We thus have a commutative diagram 
\begin{eqnarray} 
\label{diagcomp} \xymatrix{\hat{\MM}_x \ar[r] \ar@{^{(}->}[d] & \hat{\MM}_x \otimes_{\hat{\OO}_x} \hat{\anneau}_{x'} \ar@{->>}[r] \ar[rd] \ar[r] & \hat{\MM}'_x \ar@{^{(}->}[d] \\ 
\hat{\anneau}_x^n \ar[rr] & & \hat{\anneau}_{x'}^n.}
\end{eqnarray}
The injectivity of the vertical maps comes from the injectivity of the analogue maps in~(\ref{diag}) and the flatness of $\hat{\anneau}_x$ over $\anneau_x$ and  of $\hat{\anneau}_{x'}$ over $\anneau_{x'}$. The surjectivity of the upright horizontal map comes directly from the surjectivity  
of the analogue map from~(\ref{diag}).

To simplify notations, we shall note $A$ the local ring $\anneau_{x}$, $k$ its residue field
 and $M$ the $A$-module $\MM_x$. We set  $$\hat{A}' := \prod_{x' \in \pi^{-1}(x)} \hat{\anneau}_{x'},$$
and we will see it with the product topology.
We call $\hat{M}$ the completion of $M$, that is also $M \otimes_A \hat{A}$. By definition, it is $\hat{\MM_x}$. 
We set $$\hat{M}' := \prod_{x' \in \pi^{-1}(x)} \hat{\MM}^{'}_{x'}.$$
Note that $\hat{M}'$ is an $\hat{A}'$-module.

\begin{lemma} \label{DA} For each open (hence cofinite length) ideal $\hat{J}'$ of $\hat{A}'$, 
$$l(D(\hat{M}'/\hat{J}'\hat{M}'))=l(\hat{A}'/\hat{J}').$$
\end{lemma}

\begin{pf} Since $\hat{J'}$ is open,  
$\hat{A}'/\hat{J}'$ is a finite product of finite length rings of the form 
$\hat{\anneau}_{x'_i}/\hat{J}'_i$. For each such $i$, $\hat{J}'_i$ is open in $\hat{\anneau}_{x'_i}$ so the ideal $J'_i := \hat{J}'_i \cap 
\anneau_{x'_i}$ of $\anneau_{x'_i}$ satisfies 
$\anneau_{x'_i}/J'_i = \hat{\anneau}_{x'_i}/\hat{J}'_i$. 
 By the hypothesis of the proposition we are proving,
we thus have $l(D(\hat{M}'_{x'_i} / \hat{J}'_i \hat{M}'_{x'_i}))=l(  \hat{\anneau}_{x'_i}/\hat{J}'_i)$.
The lemma then results from the additivity of the functor $D$ and of $l$.
\end{pf}  

\begin{lemma}{\rm (\cite[Lemma 10.7]{kis})}\footnote{As stated there, the lemma assumes that $\pi$ is a blow-up, but it is only used in the proof that $\pi$ is proper and birational.}\label{injcomp} The morphism 
$\hat{A} \longrightarrow \hat{A}'$ is injective. 
\end{lemma}

By~(\ref{diagcomp}), we have a commutative diagram 
\begin{eqnarray} 
\label{diag2} \xymatrix{\hat{M}  \ar@{^{(}->}@/^5mm/[rr] \ar@{^{(}->}[r] \ar@{^{(}->}[d] 
& \hat{M} \otimes_{\hat{A}} \hat{A}'\ar@{->>}[r] \ar[rd]  & \hat{M}' \ar@{^{(}->}[d] \\ 
\hat{A}^n \ar@{^{(}->}[rr] & & \hat{A}'^n}
\end{eqnarray}
The injectivity of the vertical maps is obvious from~(\ref{diagcomp}) and the injectivity of the horizontal lower map
is Lemma~\ref{injcomp}. The injectivity of the upper horizontal map follows. \ps

The following lemma is an application of Chevalley's Theorem (cf. \cite[ex 8.6]{matsumura}) which we recall : 
let $\hat{A}$ be a complete noetherian local ring, $\hat{M}$ a finite type $\hat{A}$-module, $\hat{I}$ a cofinite length ideal of $\hat{A}$ 
and $\hat{M}_n$ a decreasing, exhaustive (that is $\cap_n \hat{M}_n =(0)$) sequence of 
submodules of $\hat{M}$. 
Then for $n$ big enough, $\hat{M}_n \subset \hat{I}\hat{M}$.  \ps

Now we go back to the proof of Proposition \ref{descenteblowup}. Let $I$ be a cofinite length ideal of $A$, and  
note $\hat{I} \subset \hat{A}$ its completion. We recall also that $\hat{M} \subset \hat{M}'$ by diagram~(\ref{diag2}).

\begin{lemma} \label{chevalley}
There exist a cofinite length ideal $\hat{J} \subset \hat{I}$ of $\hat{A}$ and an open ideal $\hat{J}'$ of  
$\hat{A}'$ such that 
\begin{itemize}
\item[(i)] $\hat{J}=\hat{J}' \cap \hat{A}$,
\item[(ii)] $(\hat{J}' \hat{M'} \cap \hat{M}) \subset \hat{I}\hat{M}$. 
\end{itemize}
\end{lemma} 
\begin{pf} 
We let $\hat{J}'_n := (\prod_{x' \in \pi^{-1}(x)} \hat{\m}_{x'}^n) \subset \hat{A}'$. 
By Krull's theorem, $\cap_n \hat{J}'_n=0$ and $\cap_n (\hat{J}'_n \hat{M}')=0$.

We set $\hat{J}_n := \hat{J}'_n \cap \hat{A}$, the intersection being in $\hat{A'}$. 
Similarly, we set $\hat{M}_n= (\hat{J}'_n \hat{M}') \cap \hat{M}$, the intersection being in $\hat{M}'$.
Then  $\cap_n \hat{J}_n=0$ and $\cap_n \hat{M}_n=0$.

By Chevalley's theorem, applied twice, once to the finite module 
$\hat{M}$ over the local complete noetherian ring $\hat{A}$, and once to 
$\hat{A}$ as a module over 
itself, we know that for $n$ big enough, 
$\hat{M}_n \subset \hat{I}\hat{M}$, and $\hat{J}_n \subset \hat{I}$.  

We fix such an $n$.  We set $\hat{J}:=\hat{J}_n$. It is clear that $\hat{J}$ is of cofinite length since it 
contains $\hat{m}_x^n$. 
We thus have by construction $\hat{J} \subset \hat{I}$, 
$\hat{J}=\hat{J}'_n \cap \hat{A}$ and $(\hat{J}'_n \hat{M'}) \cap \hat{M} \subset \hat{I}\hat{M}$.

However, $\hat{J}'_n$ is not open.
If $F$ is a finite subset of $\pi^{-1}(x)$, we let $\hat{J}'_F$ be the ideal 
$\prod_{x' \in F}  \hat{m}_{x'}^n \times \prod_{x' \in \pi^{-1}(x)-F} \hat{\anneau}_{x'}$ of $A'$.  It is clear that $\hat{J}'_n= \cap_F \hat{J}'_F$, 
and that the $\hat{J}'_F$ are open ideals of $A'$. Because $\hat{A}/\hat{J}$ and $\hat{M}/\hat{M}_n$ are artinian, 
there is a finite $F$ such that 
$\hat{J}'_F \cap \hat{A} = \hat{J}$ and $\hat{J}'_F \hat{M}' \cap \hat{M} = \hat{M}_n$. We set $\hat{J'}$ equal to this $\hat{J'}_n$ and we are done.
\end{pf}

By (i) of this lemma, the morphism of $\hat{A}[{G}]$-modules 
$\hat{M} \longrightarrow \hat{M}'$ induces a morphism of 
$(\hat{A}/\hat{J})[G]$-modules
$$ f : \hat{M}/\hat{J}\hat{M} \longrightarrow \hat{M}'/\hat{J'}\hat{M'}.$$ 
Indeed, the image of $\hat{J}\hat{M} \subset \hat{M}$ in $\hat{M}'$ is included in 
$\hat{J}\hat{M}'$ which is included in $\hat{J}'\hat{M}'$.   

We shall denote by $K$, $C$ and $Q$ the kernel, cokernel and 
image of $f$, respectively. 
Thus we have two exact sequences of $(\hat{A}/\hat{J})[{G}]$-modules :
\begin{eqnarray} \label{imagef} 0 \longrightarrow K \longrightarrow 
\hat{M}/\hat{J}\hat{M} \longrightarrow Q \longrightarrow 0 \end{eqnarray}
\begin{eqnarray} \label{conoyau} 0 \longrightarrow Q \longrightarrow 
\hat{M}'/\hat{J}'\hat{M}' \longrightarrow C \longrightarrow 0 \end{eqnarray}
Note that the five modules involved here are all 
of finite length as $\hat{A}/\hat{J}$-modules.

\begin{lemma} \label{CQuotient} As an $\hat{A}[G]$-module, 
$C$ is a 
quotient of $({\hat{M}}/\hat{J}{\hat{M}}) 
\otimes_{\hat{A}} (\hat{A}'/\hat{J}')/(\hat{A}/\hat{J})$.
\end{lemma}
\begin{pf} This is formal. Indeed,  
we have a commutative diagram, where the vertical arrows are surjective :
$$\xymatrix{
%0 \ar[d] & 0 \ar[d] & 0 \ar[d] \\ {\hat{J}}{\hat{M}} \ar[r] \ar [d] & J'({\hat{M}} \otimes_{\hat{A}} {\hat{A}}') \ar@{->>}[r] \ar[d] & {\hat{J}}'{\hat{M}}' \ar[d] \\  
{\hat{M}} \ar[r] \ar[d] & {\hat{M}} \otimes_{\hat{A}} {\hat{A}}'   \ar[r]^s \ar[d] & {\hat{M}}' \ar[d]
\\  {\hat{M}}/{\hat{J}}{\hat{M}} \ar@/_2pc/[rr]^{f}\ar[r]^g & {\hat{M}}
 \otimes_{\hat{A}} {\hat{A}}'/{\hat{J}}' \ar[r]^h & {\hat{M}}'/{\hat{J}}'{\hat{M}}' 
 %\\ 0 & 0 & 0 
}$$

This diagram makes clear that the map labeled $h$ is surjective, since the one labeled $s$ is.
Hence the cokernel $C$ of $f$ is a quotient of the cokernel of the map labeled 
$g :\, {\hat{M}}/{\hat{J}}{\hat{M}} \longrightarrow {\hat{M}} \otimes_{\hat{A}} {\hat{A}}'/\hat{J}' = 
({\hat{M}}/{\hat{J}}{\hat{M}}) \otimes_{{\hat{A}}} {\hat{A}}'/\hat{J}'$ and this cokernel is, by right-exactness of the tensor product by 
$\hat{M}/\hat{J}\hat{M}$, the module $({\hat{M}}/{\hat{J}}{\hat{M}}) \otimes_{\hat{A}} ({\hat{A}}'/\hat{J}')/({\hat{A}}/{\hat{J}})$.
\end{pf}

We now prove an abstract lemma concerning the left-exact functor $D$ and 
length of modules. 

\begin{lemma}\label{abstract} Let $V$ be an $A$-module of finite length with 
a continuous action of ${G}$, such that $$l(D((V \otimes k)^\ses) \leq 1.$$ 
 Let $N$ be an $A$-module of finite 
length\footnote{We view it as a $G$-module for the trivial action.} and $\pi: V \otimes_A N \longrightarrow Q$ a surjective $A[{G}]$-linear 
map.
\begin{itemize}
\item[(i)] $l(D(Q)) \leq l(N)$.
\item[(ii)] Assume that equality holds in (i), and that there is a surjective map of $A$-modules $N \longrightarrow N'$ such that the natural induced surjection $V \otimes_A N \longrightarrow V 
\otimes_A N'$ factors through $\pi$. Then $l(D(V \otimes_A N'))=l(N')$.
\item[(iii)] Let $J$ be a cofinite length ideal of $A$.
If $l(D(V/JV))=l(A/J)$, then for each ideal $J'\supset J$, 
$l(D(V/J'V))=l(A/J')$.
\end{itemize}
\end{lemma}
\begin{pf} 
First remark that the hypothesis $l(D((V \otimes k)^\ses) \leq 1$ implies, by left exactness of $D$, that $l(D(U)) \leq 1$ for any
subquotient $U$ of $(V \otimes k)$ (as a $k[{G}]$-module). \ps

Let us prove (i). There is a filtration $N_0 \subset \dots \subset N_i \subset \dots \subset N_{l(N)}=N$ of $N$ such that 
$N_i/N_{i-1}\simeq k$. We denote by $VN_i$ the image of $V \otimes_A N_i$ into $V \otimes_A N$ and by $Q_i$ the image of $VN_i$ in $Q$. It is clear that $VN_i / VN_{i-1}$ is a 
quotient of $V \otimes k$, and that $Q_i/Q_{i-1}$ is a quotient of $VN_i/VN_{i-1}$, hence we have $l(D(Q_i/Q_{i-1})) \leq 1$ by the remark beginning the proof. By left exactness of $D$, this proves (i).
Note also that if $l(D(Q))=l(N)$, all the inequalities above have to be equalities, so that $l(D(Q_i))=i$ for each $i$. \par

Let us prove (ii).  In the proof of (i) above, we can certainly choose the $N_i$ such that one of them, say $N_k$, is the kernel of the surjection $N \longrightarrow N'$. Then 
$k=l(N')-l(N)$. We have an exact sequence $0 \longrightarrow VN_k \longrightarrow V \otimes N \longrightarrow V\otimes N' \longrightarrow 0$, hence (using the hypothesis) an exact 
sequence $$0 \longrightarrow Q_k \longrightarrow Q \longrightarrow V \otimes N' \longrightarrow 0.$$ Because $D$ is left exact, we have 
$l(D(V \otimes N')) \geq l(D(Q)) - l(D(Q_k))$. But by hypothesis, we have $l(D(Q))=l(N)$, which implies by the remark at the end of the proof of (i) that $l(D(Q_k))=k$. Hence $$l(D(V \otimes N')) \geq l(N)- k = l(N) -(l(N)-l(N'))=l(N').$$ The other equality follows from (i), hence $l(D(V \otimes N'))=l(N')$.

The assertion (iii) is a special case of (ii) : apply (ii) to $Q=V \otimes_A N=V/IV$, $\pi=\Id$ 
and $N'=A/J'$.
\end{pf}

Going back to the proof of the Proposition \ref{descenteblowup} we get the following lemma.

\begin{lemma} \label{lastlemma}We have 
\begin{itemize} 
\item[(i)] $l(D(C)) \leq l(\hat{A}'/\hat{J}') - l(\hat{A}/\hat{J})$,  
 \item[(ii)] $l(D(Q)) = l(\hat{A}/\hat{J})$,
\item[(iii)] $l(D(\hat{M}/\hat{I}\hat{M}))=l(\hat{A}/\hat{I})$.
\end{itemize}
\end{lemma}

\begin{pf}
Lemma~\ref{CQuotient} tells us that $C$ is a quotient of the module 
$$({\hat{M}}/\hat{J}{\hat{M}}) \otimes_{\hat{A}} (\hat{A}'/\hat{J}')/(\hat{A}/\hat{J}).$$
We now apply the point (i) of Lemma~\ref{abstract} to $V=\hat{M}/\hat{J} \hat{M}$ and $N= (\hat{A}'/\hat{J}')/(\hat{A}/\hat{J})$. 
We note that $V \otimes_{\hat{A}} k$, that is $\hat{M} \otimes_{\hat{A}} k= M \otimes_A k$, satisfies the hypothesis of Lemma \ref{abstract} by hypothesis.
So $l(D(C)) \leq l(D(N))$, hence (i).

To prove (ii) note that by the exact sequence~(\ref{conoyau}), \begin{eqnarray*}
l(D(Q)) &\geq& l(D(\hat{M}'/\hat{J}'\hat{M}'))-l(D(C)) \\ &\geq& 
l(D(\hat{M}'/\hat{J}'\hat{M}'))-l(\hat{A}'/\hat{J'}) + l(\hat{A}/\hat{J}), \, \, \, {\rm by (i)}.
\end{eqnarray*}

Since $l(D(\hat{M}'/\hat{J}'\hat{M}')) = l(\hat{A}'/\hat{J}')$ by the Lemma~\ref{DA} we get 
$$l(D(Q)) \geq l(\hat{A}/\hat{J}).$$
To get the other inequality, recall that $Q$ is by construction a quotient of 
$\hat{M}/\hat{J}\hat{M}=\hat{M} \otimes_{\hat{A}} \hat{A}/\hat{J}$, so by point (i) of Lemma~\ref{abstract} we have
$l(D(Q)) \leq l(\hat{A}/\hat{J})$.

Let us prove (iii). Assertion (ii) of Lemma~\ref{chevalley} tells that $\hat{M}/\hat{J}\hat{M} \longrightarrow
\hat{M}/\hat{I}\hat{M}$ factors through
  the canonical surjection $\hat{M}/\hat{J}\hat{M} \longrightarrow Q$.
 We apply point (ii) of Lemma~\ref{abstract} to $Q$, with 
$V=\hat{M}/\hat{J}\hat{M}$, $N=\hat{A}/\hat{J}$, $N'=\hat{A}/\hat{I}$. This is possible because $l(D(Q))=l(N)$ by (ii) above, 
and that gives us
$l(D(V \otimes N'))=l(N')$, which is (iii). 
\end{pf}

Now recall that since $I$ is of cofinite length, 
$A/I \simeq \hat{A}/\hat{I}$ and $M/IM \simeq \hat{M}/\hat{I}\hat{M}$.
Hence by (iii) of Lemma \ref{lastlemma} above, 
$$l(D(M/IM))=l(A/I).$$
The proof of Proposition~\ref{descenteblowup} is complete.

\subsection{Direct generalization of a result of Kisin}

\label{kisindirect}

\subsubsection{Notations and definitions}
\label{accumulate}

We keep the general notations of paragraph \ref{notationsrigide}. We fix $p$ a prime number and set $$\Gp=\Gal(\Qpb/\Qp).$$

Recall that a subset $Z \subset X$ is said to be Zariski-dense 
if the only analytic subset of $X$ containing $Z$ is $X_{\rm red}$ 
itself. We shall need below some arguments involving the notion of irreducible component of a rigid analytic space, for which we refer to \cite{con}. \par
We will say that a subset $Z \subset X$ {\it accumulates} 
at $x \in X$ if there is a basis of affinoid neighborhoods $U$ of $x$ 
such that $U \cap Z$ is Zariski-dense in $U$.

\subsubsection{Hypotheses}
\label{directhyp}

We assume that we are given a couple of maps $(F,\kappa) \in \OO(X)^*\times 
\OO(X)$, and a Zariski-dense subset $Z \subset X$ satisfying the following conditions.\ps
\begin{itemize}
\item[(CRYS)] For $z \in Z$, $\MB_z$ is a crystalline representation of 
${G_p}$ whose smallest Hodge-Tate weight is $\kappa(z) \in \Z$, and that 
satisfies $D_{\rm crys}(\MB_z)^{\varphi=p^{\kappa(z)}F(z)}\neq 0$.
\item[(HT)] For any non-negative integer $C$, if $Z_C$ denotes the 
subset of $z \in Z$ such that the Hodge-Tate weights of $\MB_z$ other 
than $\kappa(z)$ are bigger that $\kappa(z)+C$, then $Z_C$ accumulates 
at any point of $Z$. \par \smallskip \noindent
\end{itemize}
\begin{remark}\label{remoncondzd} The assumption (HT) together with the Zariski-density of $Z$ in $X$ imply that $Z$ accumulates at each of its points. This stronger density condition on $Z$, introduced in \cite{jlpch} under the terminology "Z is {\it very Zariski-dense} in X", turns out to be rather well-behaved and allows to avoid some pathological Zariski-dense subset.\footnote{As an exercise, the reader can check that there are Zariski-dense subsets of $\AAA^2$ whose intersection with any affinoid subdomain $V \subset \AAA^2$ is not Zariski-dense in $V$. However, if $Z$ is a very Zariski dense subset of a rigid space $X$, then for any irreducible component $T$ of $X$ there is an open affinoid of $T$ in which $Z$ is Zariski-dense.}
\end{remark} 
For some technical reasons, we shall also need to know that:\ps
\begin{itemize}
\item[($\ast$)] There exists a continuous character $\Z_p^* 
\longrightarrow \OO(X)^*$ whose derivative at $1$ is the map $\kappa$ and whose evaluation at any point $z \in Z$ is the elevation to the $\kappa(z)$-th power.
\end{itemize}
\ps
Condition $(\ast)$ allows us to define by composition with the cyclotomic character $\chi$ a continuous character $$\psi: \Gp \underset{\chi}{\isomo} \Z_p^* 
\longrightarrow \OO(X)^*$$ 
whose evaluation of at any point $z \in Z$ is then the $\kappa(z)$-th power of the cyclotomic character (whence crystalline). \par
\index{Zkappa@$\kappa$ or $\kappa_i$, a Hodge-Tate-Sen weight in $\OO(X)$ or the associated character $\Zp^* \longrightarrow \OO(X)^*$}
\begin{definition} \label{defkappa} We shall often denote by $\kappa: \Gp \longrightarrow \OO(X)^*$ the character $\psi$ defined above, and if $N$ is any sheaf of $\OO[G_p]$-modules on $X$, we will also  
denote by $N(\kappa)$ the $\OO$-module $N$ whose $G_p$-action is twisted by the character $\psi$.
\end{definition}
\subsubsection{The subspace $X_{fs}$} The arguments in this part will follow closely Kisin's paper \cite[\S5]{kis}. We want first to apply Kisin's construction \cite[Prop. 5.4]{kis} to prove that $X_{fs}=X$. Precisely, this proposition determines a Zariski closed subspace of $X$ (there 
denoted by $X_{fs}$) under the assumption that $\MM$ is a free module 
(not only locally free). However, it is  formal  that under the weaker assumption 
``$\MM$ is locally free'', \cite[Prop. 5.4]{kis} still holds
 if we relax its condition (2) by asking that 
the maps $f$ considered there fall in an admissible open subset of $X$ on 
which $\MM$ is a free module.

Indeed, it suffices to apply {\it loc. cit.} to an admissible covering 
$(U_i)$ of $X$ by affinoids on which $\MM$ is free. We define $U_{i,fs}$ 
by the Proposition loc. cit. and we have to show that they glue, that 
is $U_{i,fs}\cap U_j = U_{j,fs} \cap U_i$. 
But if $U \subset V$ are open affinoids on which 
$\MM$ is free, $V_{fs}=U_{fs} \cap V$ by the last 
assertion of {\it loc. cit.}~, hence the 
intersections above both coincide with $(U_i\cap U_j)_{fs}$. \ps
 
We still denote by $X_{fs}$ the subspace of $X$ defined by the 
generalization explained above of \cite[Prop 5.4]{kis} when $\MM$ is locally 
free.
 
\begin{theorem}
\label{thmkisin}
Assume $\MM$ is locally free. 
\begin{itemize}
\item[(i)] For all $x \in X$, then 
$\Dcp(\MB_x(\kappa(x)))^{\varphi=F(x)}$ is non zero. Moreover, 
$X_{fs}=X$. 
\item[(ii)] Let $x \in X$ and assume that $\Dcp(\MB_x^{\rm 
ss}(\kappa(x)))^{\varphi=F(x)}$ has $k(x)$-dimension $1$. Then for all 
ideal $I$ of cofinite length of $\OO_x$, 
$\Dcp((\MM_x/I\MM_x)(\kappa))^{\varphi=F}$ is free of rank $1$ over $\OO_x/I$.
\end{itemize}
\end{theorem}

\begin{remark} \label{remarkdcp}
Part (i) of this theorem is a combination of results of 
Kisin in \cite{kis}. Moreover, he proved {\it loc. cit.} some cases of 
part (ii), essentially the cases where 
$\MB_x$ is an indecomposable $k(x)[{G_p}]$-module (although it is not stated explicitly, this is done during the proof of Proposition 10.6 of \cite{kis}, 
page 444 and 445). The proof we give 
here simplifies a bit some arguments of \cite[section 8]{kis} and avoids all use of universal deformation ring, using some length arguments and our lemma of descent by blow-up instead. 
It also paves the way for the proof of Theorem~\ref{kisinnonlibre} below.
\end{remark}
\newcommand{\oc}{\widehat{\otimes}}

\begin{pf}
By replacing $\MM$ by $\MM(\kappa)$, we may assume that $\kappa=0$. Let
  $$TQ(T)\in \OO(X)[T]$$ be the Sen polynomial of $\MM$, whose roots at $x 
\in X$ are the generalized Hodge-Tate weights of $\MB_x$. Let $W \subset 
X$ denote the subset consisting of the points $x \in X$ such that the 
Sen polynomial of $\MB_x$ has $0$ as unique root in the integers $\N$ 
(and which is a simple root).

\begin{lemma} \label{zarht} For each admissible open $U$ of $X$, 
$W \cap  U$ is Zariski-dense in $U$.
\end{lemma}

\begin{pf} For each $k\geq 0$, and $U \subset X$ admissible open, let 
$U_k$ denotes the (reduced) zero locus of $Q(k)$ in $U$, so $$W\cap U=U-\bigcup_{k\geq 0}U_k.$$
Let $T$ be a closed 
analytic subset of $U$ such that $U=T\cup \bigcup_{k\geq 0} U_k$. Let 
$T'$ be any irreducible component of $U$. If $T' \not \subset T$, then 
$T' \subset U_k$ for some $k$ by \cite[Lemma 5.7]{kis}. Let $T''$ be an 
irreducible component of $X$ such that $T''\cap U \supset T'$, then 
$T'' \subset X_k$, which is not possible by (HT) applied to $C=k+1$. 
Hence $T=U$, which proves the lemma.
\end{pf}

To prove that $X_{fs}=X$ it suffices to show (as Kisin does to prove his Theorem 6.3) that

\begin{lemma} \label{denseW} The set $\{x \in W,\,  
\Dcp(\MB_x)^{\varphi=F(x)}\neq 0\}$ is Zariski-dense in $X$.  \end{lemma}

Indeed, by Tate's computation of the cohomology of $\C_p(i)$ for $i \in 
\Z$, the natural map 
$$\Ddr^+(\MB_x) \longrightarrow (\MB_x\otimes_{\Q_p} \C_p)^\Gp$$ 
is an isomorphism between $k(x)$-vector-spaces of dimension $1$ when $x \in W$. In particular, if $x$ is in the subset of Lemma \ref{denseW}, the natural injection
$$\Dcp(\MB_x)^{\varphi=F(x)} \longrightarrow \Ddr^+(\MB_x)$$
is an isomorphism, hence $x \in  X_{fs}$.

\begin{pf} Let us fix first some $z \in Z$ and choose an open affinoid 
$U \subset X$ containing $z$ which is small enough so that $\MM$ is 
free over $U$, $U$ is $F$-small (\cite[(5.2)]{kis}), and such that $Z$ is Zariski-dense in $U$ (it exists by (HT)).
Assumption (HT) implies then that $Z_C \cap U$ is Zariski-dense in $U$ for
any $C$.

 We now apply \cite[Prop. 5.14]{kis} and its corollary \cite[Cor 5.15]{kis}
 to $\cal{R}:=\anneau(U)$, $M:=\MM_{|U}$, $I:=Z\cap U$, $\cal{R}_i:=k(i)$ and 
$I_k:=Z_{k+\sup_U|F|+1}$. Note that we just checked condition (3) there (that is, $I_k$ is Zariski-dense
in $U$) 
and that condition (2) follows from our assumption (ii). Moreover, 
condition (1) follows from (CRYS) and the weak admissibility of $D_{\rm 
crys}(\MB_x)$, $x \in I_k$, applied to the filtered $\varphi$-submodule 
$\Dcp(\MB_x)^{\varphi=F(x)}$. As a consequence, \cite[cor. 
5.15]{kis} tells  that for all $x \in U$, 
$\Dcp(\MB_x)^{\varphi=F(x)}\neq 0$. We conclude the proof by Lemma 
\ref{zarht}.
\end{pf}
Applying now \cite[cor. 5.16]{kis}, we first get the point (i) of our theorem~\ref{thmkisin}.
\begin{remark} \begin{itemize}
\item[(i)] We note the extreme indirectness of this method of proof (which is entirely Kisin's) :
to prove that $\Dcp(\MB_x)^{\varphi=F(x)}\neq 0$ for every $x \in X$, knowing that this is true for the points of $Z$, we use the closed set $X_{fs}$, which by definition contains the points
satisfying this properties provided that they are in the set $W$ - in particular, not in $Z$ ! 
\item[(ii)] The proof of Lemma~\ref{denseW} shows that if $X$ is an affinoid space, $F$-small, on which 
$\MM$ is free, then in the proof of point (i) of Theorem~\ref{thmkisin} 
condition (HT) may be replaced  by the weaker condition \par \smallskip
\begin{itemize}
\item[(HT')] : for every non negative $C$,
$Z_C$ is Zariski-dense in $X$.
\end{itemize}
\end{itemize}
\end{remark}

We now prove point (ii) of our Theorem~\ref{thmkisin}.
 Let us fix some $x 
\in X$ (but not necessarily in $Z$) and choose an $F$-small open affinoid neighborhood $U$ of $x$ such that $\MM$ is
free over $U$. As $U \subset X=X_{fs}$, by the corollary {\it loc. cit.} we get 
that $\Dcp(\MM(U))^{\varphi=F}$ is generically free of rank $1$ over 
$\anneau(U)$. 
More precisely, if $H \subset \anneau(U)$ denotes the smallest ideal such 
that\footnote{\label{coef}The Banach $\anneau(U)$-module $\MM(U) 
\oc_{\Q_p}B_{\rm 
crys}^+$ is ON-able, $H$ is the ideal of $\OO(U)$ generated by all the 
coefficients in a given $ON$-basis of all the elements of 
$\Dcp(\MM(U))^{\varphi=F}$. It does not depend on the choice of the ON-basis 
as the ideals of $\anneau(U)$ are all closed.}$$\Dcp(\MM(U))^{\varphi=F} \subset 
H(\MM(U)\oc_{\Q_p} B_{\rm crys}^+),$$ then $U-V(H)$ is Zariski-dense in $U$. 
Let $$\pi: U' \longrightarrow U,$$ be the blow-up of the ideal $H$ and 
$\cal M'$ the pullback of $\cal M$ on $U'$.

\begin{lemma}\label{surblow} Let $x' \in U'$ and let $V \subset U'$ 
be a sufficiently small open affinoid containing $x'$.

\begin{itemize}
\item[(i)] The ideal of $\anneau(V)$ generated by all the coefficients (see 
the footnote~\ref{coef}) of $\Dcp(\MM'(V))^{\varphi=F} \subset 
\MM'(V)\oc_{\Qp}\Bcp$ is $\anneau(V)$ itself.
\item[(ii)] If $I'$ is a cofinite length ideal of $\anneau_{x'}$ then 
$\Dcp(\MM'_{x'}/I'\MM'_{x'})^{\varphi=F}$ is free of rank 
$1$ over $\anneau_{x'}/I'$.
\end{itemize}
\end{lemma}

\begin{pf}
By the universal property of blow-ups, for $V$ sufficiently small $H \anneau(V)$ is a principal ideal
generated by a non zero divisor $f_V$ of $\anneau(V)$. 
As a consequence, the ideal of 
the statement is $\anneau(V)$ itself, as it contains $H \anneau(V)/f_V$. 
Indeed, it 
is clear that if $\Dcp(\MM'(V))^{\varphi=F}$ contains $f v$ for some non 
zero divisor $f \in \anneau(V)$ and $v \in \MM'(V)\oc_{\Qp}\Bcp$, it contains 
$v$. This proves (i). 

It follows that
the natural map $\Dcp(\MM'_{x'}/I'\MM'_{x'})^{\varphi=F} 
\longrightarrow \Dcp(\MB'_{x'})^{\varphi=F}$
is non-zero. Moreover $\Dcp(\MB'_{x'})^{\varphi=F}=\Dcp(\MB_{x})^{\varphi=F} \otimes_{k(x)} k(x')$, hence it has $k(x')$-dimension $1$ by assumption on $\MB_x$ and part (i) of Theorem \ref{thmkisin}. So the first assertion of the following lemma (applied to $D=\Dcp(-)^{\varphi=F}$, $A=\OO_{x'}$, $J=I'$, 
$V=\MM'_{x'}/I'\MM'_{x'}$) 
implies the result.
\end{pf}

The following lemma holds in the same context as Lemma~\ref{abstract}.

\begin{lemma}\label{abstract2} Let $J$ be
a cofinite length ideal of $A$,
$V$ a continuous $(A/J)[G_p]$-module that is {\it free} of finite rank over $A/J$ and such that 
$D(V \otimes_A k)$ has $k$-dimension 1.
Assume moreover that one of the following two conditions holds:
\begin{itemize}\item[(i)] $D(V) \longrightarrow D(V \otimes_A k)$ is non-zero,
\item[(ii)] $l(D(V))=l(A/J)$.
\end{itemize}
Then $D(V)$ is free of rank one over $A/J$.
\end{lemma}
\begin{pf}
Under assertion (i), the lemma is exactly \cite[Lemma 8.6]{kis}. Under assertion (ii), it can be proved using similar ideas: we prove that $D(V \otimes_A A/J')$ is free of rank one over 
$A/J'$ for any ideal $J'$ containing $J$, by induction on 
the length of $A/J'$. There is nothing to prove for $J'=m$. Assume the result known for ideals of colength $<k$, and let $J'$ be an ideal containing $J$ of colength $k$. Let $J''$ be an ideal such that $J' \subset J''$, the first inclusion being proper and of colength one. 
We have (since  $V \otimes_A A/J'$ is free over $A/J'$) 
an exact sequence :
$$0 \longrightarrow D(V\otimes_A k) \otimes_k J''/J' 
\longrightarrow D(V \otimes_A A/J') \longrightarrow D(V \otimes_A A/J'').$$
By (iii) of Lemma~\ref{abstract}, $l(D(V \otimes_A A/J'))=l(A/J')$ and 
similarly for $J''$. Hence  
the last morphism of the exact sequence above is surjective. So we have 
$D(V \otimes_A A/J') \otimes_A A/J''=D(V \otimes_A A/J'')$, hence 
$D(V \otimes_A A/J') \otimes_A k = D(V \otimes_A A/J'') \otimes_A k$. By induction,
the latter has $k$-dimension $1$. 
Hence by Nakayama's lemma, the $A/J'$-module $D(V \otimes_A A/J')$ is generated by a single element
and since its length is $l(A/J')$, it is free of rank one over $A/J'$.
\end{pf}

We can now use our ``descent result'' (Proposition~\ref{descenteblowup})
for the blow-up $\pi : U' \longrightarrow U$.
Assertion (ii) of Lemma~\ref{surblow} shows that for every 
$x' \in \pi^{-1}(x)$,
and every cofinite length ideal $I'$ of $\OO_{x'}$,  
$$l\left( \Dcp(\MM'_{x'} \otimes_{\OO_{x'}} \anneau_{x'}/I')^{\varphi=F} \right)
= l(\anneau_{x'}/I').$$
Thus by Proposition~\ref{descenteblowup}, we have for every cofinite length ideal
$I$ of $\anneau_{x}$,
$$l\left( \Dcp(\MM_x \otimes_{\OO_x} \OO_x/I)^{\varphi=F} \right)
= l(\OO_x/I).$$
To conclude that $\Dcp(\MM_x \otimes_{\OO_x} \OO_x/I)^{\varphi=F}$ is free of rank one over 
$\OO_x/I$ we simply invoke Lemma~\ref{abstract2} (ii) with $I=J,V=Q=\MM_x/J\MM_x$. 
The proof of Theorem~\ref{thmkisin} is now complete.
\end{pf}

\subsection{A generalization of Kisin's result for non-flat
modules}\label{kisinindirect}

In this subsection we keep the assumptions of \S\ref{directhyp}, but we 
do not assume that $\MM$ is locally free, but only that $\MM$ is torsion-free.

\begin{theorem} \label{kisinnonlibre}
 Let $x \in X$ and assume that\footnote{In fact, the result holds more generally under the assumption of Remark~\ref{generalitepropdescente}, but we state it as such for short.} $\Dcp(\MB_x^{\rm 
ss}(\kappa(x)))^{\varphi=F(x)}$ has $k(x)$-dimension $1$. Then for all 
ideal $I$ of cofinite length of $\OO_x$, 
$$l \left(\Dcp(\MM_x/I\MM_x(\kappa))^{\varphi=F}\right)=l(\OO_x/I).$$
\end{theorem}

We will rely on the following flatification result whose scheme theoretic analogue is an elementary case of a result of Gruson-Raynaud (\cite[Thm. 5.2.2]{GR}). Recall that $X$ is reduced and separated.

\begin{lemma}\label{grra} There exists a proper and birational morphism $\pi: X'\longrightarrow X$ (with $X'$ reduced) such that the strict transform of $\mathcal M$ by $\pi$ is a locally free coherent sheaf of modules $\mathcal M'$ on $X'$. More precisely, we may choose $\pi$ to be the blow-up along a nowhere dense closed subspace of the normalization $\tilde{X}$ of $X$. 
\end{lemma} 

\begin{pf} Let $f: \tilde{X} \longrightarrow X$ be the normalization of $X$ (see \cite[\S 2.1]{con}), then $\tilde{X}$ is reduced, $f$ is finite (hence proper), and $f$ is birational by \cite[Thm. 2.1.2]{con}. Moreover, the strict transform $\MM'$ of $\MM$ by $f$ is torsion free as $\MM$ is, hence by replacing $(X,\MM)$ by $(\tilde{X},\MM')$ we may assume that $X$ is normal. We may then assume that $X$ is connected.

We claim that there is an integer $r\geq 0$ such that for each open affinoid $U \subset X$, $\MM(U)$ is generically free of rank $r$ over $\OO(U)$. If $U$ is connected (hence irreducible), let us denote by $r_U$ this generic rank. There is an injective $\OO(U)$-linear map $\MM(U) \longrightarrow \OO(U)^{r_U}$ which is an isomorphism after inverting some $f\neq 0 \in \OO(U)$. In particular, for each $x$ in a Zariski-open subset of $U$, we have $\MM_x \isomo \OO_x^{r_U}$. As a consequence, for each open affinoid $U'\subset U$, the $\OO_x$-module $\MM_x$ is free of rank $r_U$ on a Zariski open and dense subset of $U'$, thus $r_{U'}=r_U$ if $U'$ is connected. A connectedness argument shows then that $r_U$ is independent of $U \subset X$, and the claim follows. 
In particular, for all $x \in X$ the torsion free $\OO_x$-module $\MM_x$ has also generic rank $r$. 
%Note that if $r=0$, then $\MM=0$ and we are done, so we may assume that $r>0$.

Let us recall now some facts about the Fitting ideals (see \cite[XIX,\S 2]{lang}, \cite[\S 5.4]{GR}). For each open affinoid $U \subset X$ it makes sense to consider the $r$-th Fitting ideal $F_r(\cal \MM(U))$ of the finite $\OO(U)$-module $\MM(U)$. Its formation commutes with any affinoid open immersion so those $\{F_r(\cal \MM(U))\}$ glue to a coherent sheaf of ideals $F_r(\MM) \subset \OO_X$. A point $x \in X$ lies in $V(F_r(\MM))$ if and only if $\dim_{k(x)}(\MB_x)>r$ and $X-V(F_r(\MM))$ is the biggest admissible open subset of $X$ on which $\MM$ can be locally generated (on stalks) by $r$ elements. By what we saw in the paragraph above, $X-V(F_r(\MM))$ is actually Zariski dense in $X$. 
Moreover, if $x \in X-V(F_r(\MM))$ then $\MM_x$ is free of rank $r$ over $\OO_x$. Indeed, it can be generated by $r$ elements and we saw that $$\MM_x \subset \MM_x \otimes_{\OO_x} {\rm Frac}(\OO_x) \isomo {\rm Frac}(\OO_x)^r$$ for each $x\in X$, and we are done.

Let $\pi: X' \longrightarrow X$ be the blow-up of $F_r(\MM)$, we will eventually prove that $\pi$ has all the required properties. Note that $X'$ is reduced as $X$ is and that $\pi$ is birational as $X-V(F_r(\MM))$ is Zariski dense in $X$. As a general fact, the coherent sheaf of ideals $F_r(\MM)\OO_{X'}$ coincides with the $r$-th sheaf of Fitting ideals $F_r(\pi^*\M)$ of $\pi^*\MM$, and it is an invertible sheaf by construction. Let $Q \subset \pi^*\MM$ be the coherent subsheaf of $F_r(\MM)\OO_{X'}$-torsion of $\pi^*\MM$. We claim that $(\pi^*\MM)/Q$ is locally free of rank $r$. This can be checked on the global sections on an open affinoid $U \subset X'$. But if $A$ is a reduced noetherian ring and $M$ a finite type $A$-module such that $M$ is generically free of rang $r$ and whose $r$-th Fitting ideal $F_r(M)$ is invertible, then $M/{\rm Ann}_M(F_r(M))$ is locally free of rank $r$ by \cite[Lemma 5.4.3]{GR}. This proves the claim if we take $A=\OO(U)$ and $M=\pi^*(\MM)(U)$!
 .

By definition, the strict transform $\MM'$ of $\MM$ is the quotient of $\pi^*\MM$ by its $(F_r(\MM)\OO_{X'})^{\infty}$-torsion. 
The natural surjective morphism $\pi^*\MM \rightarrow \MM'$ factors then through $(\pi^*\MM)/Q$, which is locally free of rank $r$ by what we just proved, so $(\pi^*\MM)/Q \isomo \MM'$ is locally free of rank $r$, and we are done.
\end{pf}

\begin{pf} (of Theorem \ref{kisinnonlibre}) 
%We claim first that we may assume that $X$ is quasi-compact. Indeed, let $U$ be an affinoid neighborhood of $x$ in $X$. By the Zariski-density and accumulation property of $Z$ in $X$ (by (HT)), there exists a quasi-compact admissible open subspace $U \subset U' \subset X$ such that $Z \cap U'$ is Zariski-dense in each  irreducible component of $U'$ meeting $U$. Let $X_0$ be the Zariski-closure of $Z \cap U'$ in $U'$. The immersion $X_0 \rightarrow X$ is an open immersion at $x$ (actually, over $U$), and we check at once that $(\MM_{|X_0},F_{|X_0},\kappa_{|X_0},Z\cap X_0)$ still satisfies the assumptions of \S\ref{directhyp}, which proves the claim. We may  then apply Lemma~\ref{grra}.
Let us choose a $\pi$ as in Lemma~\ref{grra}, as well as a coherent sheaf of ideals $H \subset \OO_X$ attached to $\pi$ as in \S~\ref{subsubsttres}. As $X-V(H)$ is Zariski dense in $X$, and as $Z$ accumulates at $Z$ by assumption (HT), 
$Z \cap (X-V(H))$ is Zariski-dense in
$X$. Moreover (CRYS), (HT) and (*) are still satisfied when we replace $Z$ by $Z\cap (X-V(H))$
in their statement, so we may assume that $Z \cap V(H) =\emptyset$. \par

Let us denote by $Z'$ the set
of $z' \in X'$ such that $\pi(z') \in Z$.  Since 
$X'-\pi^{-1}(V(H)) \isomo X-V(H)$
is Zariski-dense in $X'$, $Z'$ is Zariski dense in $X'$. 
Note that for $z' \in Z'$, we have $\MB'_{z'}=\MB_z$ if $z=\pi(z')$.  
Define $\kappa'$ and $F'$ on $X'$ as $\kappa \circ \pi$ and 
$F \circ \pi$. Then it is obvious that $X',Z',\MM',F',\kappa'$ satisfy the 
hypothesis (CRYS), (HT) and (*). Because $\MM'$ is locally free we may apply
to it Theorem~\ref{thmkisin} at any $x' \in X'$. This implies 
Theorem~\ref{kisinnonlibre} by our descent Proposition~\ref{descenteblowup}. 
\end{pf}

\begin{remark} \begin{itemize}
\item[(i)] In the applications of Theorem \ref{kisinnonlibre} to section \ref{families}, we will use some coherent sheaves $\MM$ on an affinoid $X$ which are in fact direct sums of coherent torsion-free $\OO$-modules of generic ranks $\leq 1$, for which Lemma \ref{grra} is obvious.
\item[(ii)] As Brian Conrad pointed out to us, there is an alternative proof of the first assertion of Lemma~\ref{grra} using rigid analytic Quot spaces (see \cite[Thm. 4.1.3]{conradqcoh}).
\end{itemize}
\end{remark} 

\newpage

\section{Rigid analytic families of refined $p$-adic representations}
\label{families}

\subsection{Introduction}

In this section, we define and study the notion of 
{\it $p$-adic families of refined Galois representations}. As explained in
the general introduction, the general framework is the data of a continuous
$d$-dimensional pseudocharacter $$T: G \longrightarrow \OO(X),$$
where $X$ is a reduced, separated, rigid analytic space. Here $G$ is a
topological group equipped with a continuous map
$\Gp=\Gal(\overline{\Q}_p/\Q_p) \longrightarrow G$, and we shall be mainly
interested in the properties of the restriction of $T$ to $\Gp$. The
presence of the group $G$ is an extra structure that will only play a role when discussing the
reducibility properties of $T$, and we invite the reader to assume 
that $G=\Gp$ at a first reading. 

We assume that for all $z$ in a 
{\it Zariski-dense} subset $Z \subset X$, the (semisimple, continuous) representation
$\rhob_z$ of $G$, whose trace is the evaluation 
$T_z$ of $T$ at $z$, has the following properties (see \S\ref{defrefi}):\begin{itemize}
\item[(i)] $\rhob_z$ is crystalline,
\item[(ii)] its Hodge-Tate weights are distinct, and if we order them by 
$\kappa_1(z) < \cdots < \kappa_d(z)$, then the maps $z \mapsto \kappa_n(z)$ extend to
analytic functions on $X$ and each difference $\kappa_{n+1}-\kappa_n$ varies a lot
on $Z$.
\item[(iii)] its crystalline eigenvalues $\varphi_1(z), \cdots, \varphi_d(z)$ are
distinct, and their {\it normalized versions}
$z \mapsto F_n(z):=\varphi_n(z)p^{-\kappa_n(z)}$ extend to analytic
functions on $X$.
\end{itemize}
These hypotheses may seem a little bit complicated, 
but this is because we want 
them to encode all the aspects of the families of Galois 
representations arising on eigenvarieties. We refer to \S\ref{defrefi} for a
detailed discussion of each assumption. Let us just mention two things here.
First, although families with "constant Hodge-Tate weights" 
have been studied by
several people, the study of the kind of families above has been
comparatively quite poor let alone works of Sen and Kisin. A reason is maybe that 
the very fact that the weights are moving implies that 
the generic member of such a family is not even a Hodge-Tate representation,
and in particular lives outside the De Rham world. Second, each $\rhob_z$ is 
equipped by assumption (iii) with a natural ordering of its crystalline Frobenius
eigenvalues, that is with a {\it refinement} $\Ref_z$ of $\rhob_z$ (hence the name of the
families). \ps

Our aim is the following: we want to give a schematic upper bound of the reducibility loci at the points $z \in Z$ by proving that the infinitesimal deformations of the $\rho_z$ inside their reducibility loci (that we defined in section \ref{pseudocharacters}) are trianguline, and in favorables cases even Hodge-Tate or crystalline. 
Let us describe now precisely our results.\par Assume first that $z \in Z$ is such that
$\Lambda^n(\rhob_z)$ is irreducible\footnote{This irreducibility assumption
applies for the $\Lambda^n \rhob_z$ viewed as $G$-representations.} for each $n=1,\cdots,d$, and that $\Ref_z$ is a
{\it non 
critical regular} refinement of $\rhob_z$. Then on each thickening
$A$ of
$z$ in $X$, we show 
that $T \otimes A$
is the trace of a unique trianguline deformation of $(\rhob_z,\Ref_z)$
to the artinian ring $A$ (Theorem~\ref{tridefirr}).
\par When $\rhob_z$ is reducible, 
the situation turns out to be much more complicated, but still rather nice
in some favorable cases\footnote{In the applications to Selmer groups, we
will "luckily" be in that case.}. Assume that $\rhob_z=\oplus_{i=1}^r
\rhob_i$ is multiplicity-free. The refinement $\Ref_z$ of $\rhob_z$ induces
refinements $\Ref_{z,i}$ of the $\rhob_i$. 
This combinatorial data is in fact controlled by a permutation $\sigma \in \got{S}_d$
that we introduce in \S\ref{permutation}. Assume again that $\Ref_z$ is regular, but not that
$\Ref_z$ is non critical. Instead, we assume only that \newcommand{\Red}{{\rm Red}}
each $\Ref_{z,i}$ is a non critical refinement of $\rhob_i$, and that each
$\Ref_{z,i}$ is a "subinterval" of $\Ref_z$ (see \S\ref{conditionint}). 
As before, we also have to make some explicit $G$-irreducibility assumption
on the $\rhob_i$ for which we refer to 
\S\ref{hypredloc}. Our main result concerns then the total reducibility locus, say
$\Red_z$, of $T$ at the point
$z$. We show that each difference of weights
$$\kappa_n-\kappa_{\sigma(n)}$$
is constant on this reducibility locus $\Red_z$. We stress here that this result is schematic, it means that the closed subscheme $\Red_z$ lies in the schematic fiber of each map $\kappa_n-\kappa_{\sigma(n)}$ at $z$.
Moreover, on each thickening
$A$ of $z$ lying in the reducibility locus $\Red_z$, we show
that $T \otimes A$ writes uniquely as the sum of traces of true representations
$\rho_i$ over $A$, each $\rho_i$ being a trianguline deformation of
$(\rhob_i,\Ref_{z,i})$ (Theorem\ref{tridefgen}). 
We end the section by 
giving another proof of the assertion above on the weights on the reducibility
locus under some slightly different kind of assumptions (Theorem
\ref{thmredloc}). \par
As an example of application of the results above, let us assume
that $\sigma$ acts transitively on $\{1,\cdots,d\}$ (in which case we say that $\Ref_z$ is an
{\it anti-ordinary refinement}), so each difference of weights
$\kappa_n-\kappa_m$ is constant on $\Red_z$. If some $\kappa_m$ is moreover constant
(what we can assume up to a twist), we get that all the
weights $\kappa_i$ are constant on the total reducibility locus at $z$,
hence are distinct integers. 
In particular the deformations $\rho_i$ above of $\rhob_i$
are Hodge-Tate representations, and our work on trianguline deformations shows
then that they are even crystalline (under some mild conditions on
the $\rhob_i$, see corollary \ref{corredloc}). This fact will be very
important in the applications to eigenvarieties and global Selmer groups of the last section, as it will allow us to prove that the scheme $\Red_z$ coincides with the {\it reduced point} $z$ there. \ps

We end this introduction by discussing some aspects of the proofs and other results. We fix $z \in Z$ as above, let $A:=\OO_z$ and we consider the composed pseudocharacter $$T: G \longrightarrow \OO(X) \longrightarrow A$$
again denoted by $T$. It is residually multiplicity free and $A$ is henselian, hence $T$ fulfills the assumptions of our work in section \ref{pseudocharacters}. Some important role is played by some specific $A[G]$-modules called $M_j$ (introduced \S\ref{Mi}) whose quite subtle properties turn out to be enough to handle the difficulties coming from the fact that $T$ may not be the trace of a representation over $A$. We extend those modules, with the action of $G$, to torsion free coherent $\OO$-modules in an affinoid neighborhood $U$ of $z$ in $X$ (\S\ref{analytic}) to which we apply the results of section \ref{kisin}. \par 
 However, this only gives us a part of the information, namely the one concerning the first eigenvalue $\varphi_1$ of the refinement. Indeed, this eigenvalue is the only one that varies analytically (if $\kappa_1$ is normalized
to zero say) and so to which we can apply section \ref{kisin}. To deal with the other eigenvalues as well, we will work not only with the family $T$, but with
all its exterior powers $\Lambda^k T$. Some inconvenience of using these exterior products however appears in the fact that our definition of a refined family
is not stable under exterior powers (see \S\ref{weakly}, and the last paragraph of this introduction). This leads us to introduce the notion of $p$-adic family of {\it weakly refined} Galois representations, which is a modification of the one given above where we only care about $\kappa_1$ and $F_1$ (see Definition \S\ref{weakreffam}). Any exterior power of a refined family is then a
{\it weakly refined family}. Let us note here that an important tool to get the trianguline assertion at the end is Theorem \ref{poidscst2} of section \ref{trianguline}.\par 

In fact, our results mentioned above have analogues in the context of weakly refined families (in which case they hold for every $x \in X$), that we prove in Theorems \ref{Drhojlibre} and \ref{redlocfaible}. Another interesting result here is the proof (Theorem \ref{crystext}) that there exists a non-torsion crystalline period attached to the eigenvalue $\varphi_1$ in the infinitesimal extensions between the $\rho_i$ constructed in section \ref{pseudocharacters} (that is, in the image of $\iota_{i,j}$). 

\ps 
Though the trick of using exterior powers is not at all unfamiliar in the context of Fontaine's theory, we have the feeling that
it is not the best thing to do here, and that the use of exterior products
is responsible for some technical hypotheses to appear later in this section (e.g. assumptions (REG) and (MF') in \S\ref{hypredloc}).
But we have not found a way to avoid it. Actually, by using only similar arguments as in \cite{kis}, it seems quite hard to argue inductively (as we would like to) by "dividing modulo the families of eigenvectors for $\varphi_1$". Among other things, a difficulty is that although the points in $Z$ belong to Kisin's $X_{fs}$, they do for quite indirect reasons (see e.g. \cite[Remark 5.5 (4)]{kis}), which makes many arguments there -- and here also -- quite delicate. As a possible issue,  our work in this section and in section \ref{trianguline} comfort Colmez's idea that the construction of $X_{fs}$ in \cite{kis} should be reworked from the point of view of $\fg$-modules over the Robba ring\footnote{E.g. , for any $x$ in a refined family $X$, $\rhob_x$ should be trianguline.} and suggests that $X_{fs}$ should directly contain the points of $Z$ which are non critically refined. As this would have led us quite away from our initial aim, we did not follow this approach. We hope however that the present work sheds lights on aspects of this interesting problem.

\subsection{Families of refined and weakly refined $p$-adic representations}\label{deffamilies}

\subsubsection{Notations} \label{notationsGGp} As in sections~\ref{trianguline} and~\ref{kisin}, we set $\Gp=\Gal(\bar \Q_p/\Q_p)$. Moreover we suppose given a topological group $G$ together with a continuous morphism $G_p \longrightarrow G$.

\begin{example} The main interesting examples\footnote{Actually, 
our result in the case (a) would implies our result in the case (b), 
were there not the technical, presumably unnecessary, irreducibility 
hypothesis (MF) in~\S\ref{hypredloc} below, that we can sometimes verify in case (b) and not 
in case (a).}  are
\begin{itemize}\label{examplegg} 

\item[(a)] $G=G_p$ and the morphism is the identity.   

\item[(b)] $G=G_{K,S}=\Gal(K_S/K)$ where $K$ is a number field, $S$ a set of places of $K$, and $K_S \subset \bar K$ the maximal extension which is unramified outside $S$; the morphism sending $G_p$ to a decomposition group of $K$ at some prime $\p$ of $K$ such that $K_\p=\Q_p$.
\end{itemize} 
\end{example}

If $\rho$ is a representation of $G$, it induces a representation of $\Gp$ that we shall denote by 
$\rho_{|G_p}$. We will replace $\rho_{|G_p}$ by $\rho$ without further comments 
when the context prevents any ambiguity, for example in assertions such as 
``$\rho$ is Hodge-Tate'', or ``$\rho$ is crystalline''.
 
\subsubsection{Rigid analytic families of $p$-adic representations}\label{sectdefrigfam}
\begin{definition}\label{deffam}
A {\it (rigid analytic) family of $p$-adic representations} is the data of a reduced and separated rigid analytic space $X/\Q_p$ and a 
continuous\footnote{We recall that for each admissible open $U \subset X$ (not necessarily affinoid, e.g $U=X$), $\OO(U)$ is equipped with the coarsest locally convex topology (see \cite{Sch}) such that the restriction maps $\OO(U) \longrightarrow \OO(V)$, $V \subset U$ an open affinoid (equipped with is Banach algebra topology), are continuous. This topology is the Banach-algebra topology when $U$ is affinoid.}
 pseudocharacter $T : G \longrightarrow \anneau(X)$.
\end{definition}

The {\it dimension} of the family is the dimension of $T$; it will usually be denoted by $d$ in the sequel. For each point $x \in X$, we call {\it evaluation} of $T$ at $x$ and note 
$$T_x : G \longrightarrow k(x),$$ the composition of $T$ with the evaluation map : $\anneau(X) \longrightarrow k(x)$ at the residue field $k(x)$ of $x$. Then $T_x$ is a continuous $k(x)$-valued pseudocharacter. By a theorem of Taylor, it is the trace of a 
(unique up to isomorphism) continuous semisimple representation $$\rhob_x: G \longrightarrow \GL_d(\overline{k(x)}),$$ which is actually defined over a finite extension of $k(x)$. \ps

In other words, a family of $p$-adic representations parameterized by the rigid space $X$ is a collection of representations $\{\rhob_x, x \in X\}$ (or even over a Zariski-dense subset of $X$) for which we assume that the traces map $T(g): x \mapsto \tr(\rho_x(g))$ are analytic functions on $X$ for each $g \in G$, and such that $g \mapsto T(g)$ is continuous. Examples are given by the continuous representations of $G$ on locally free $\OO$-modules on $X$, but our definition is more general as we showed in section \ref{representations}. In particular, the families of $p$-adic Galois representations parameterized by Eigenvarieties turn out to be families in this "weak" sense in general.

%\begin{prop} \label{lemmafamilies} Assume $G=G_{K,S}$ as in case (ii) of Example \ref{examplegg}, and that $S$ has co-density one. Let $X/\Qp$ be a reduced rigid analytic space which is nested\footnote{That is, an admissible increasing union of open affinoids $X=\bigcup_{i\in \N} U_i$ such that the restrictions $\OO(U_{i+1}) \longrightarrow \OO(U_i)$ are compact. In particular, it is separated. }. Assume given: \begin{itemize} \item[(i)] for each finite place $v$ of $K$ which is not in $S$, a map $T_v \in \OO(X)$ bounded by $1$, \item[(ii)] a Zariski-dense subset $Z \subset X$, \item[(iii)] for $z \in Z$, a continuous semisimple representations $r_z: G \longrightarrow \GL_d(\overline{k(z)})$ such that $\tr(r_z(\Frob_v))=T_v(z)$ for each $v\notin S$. \end{itemize} Then there exists a unique family of $p$-adic representations parameterized by $X$ and such that $r_z \simeq \rhob_z$ for all $z\in Z$. \end{prop} \begin{pf} This is \cite[Prop. 7.1.1]{Ch} under the assumption that  $\{T_v, v\notin S\}$ has a compact closure in $\OO(X)$. But this last fact is implied by the nested hypothesis and by the fact that the $T_v$ are bounded by $1$. \end{pf}

\subsubsection{Refined and weakly refined families of $p$-adic representations}

\begin{definition}\label{defrefi} A {\it (rigid analytic) family of refined $p$-adic representations} (or shortly, a {\it refined family}) 
of dimension $d$
is a family of $p$-adic representations $(X,T)$ of dimension $d$ together with the following data
\begin{itemize}
\item[(a)] $d$ analytic functions $\kappa_1,\dots,\kappa_d \in \anneau(X)$,
\item[(b)] $d$ analytic functions $F_1,\dots,F_d \in \anneau(X)$,
\item[(c)] a Zariski dense subset $Z$ of $X$;
\end{itemize}
subject the following requirements.
\begin{itemize}
\item[(i)] For every $x \in X$, the Hodge-Tate-Sen weights of $\rhob_x$ are, with multiplicity, $\kappa_1(x),\dots,\kappa_d(x)$.
\item[(ii)] If $z \in Z$, $\rhob_z$ is crystalline (hence its weights $\kappa_1(z),\dots,\kappa_d(z)$ are integers).
\item[(iii)] If $z \in Z$, then $\kappa_1(z)<\kappa_2(z)<\dots<\kappa_d(z)$.
\item[(iv)]  The eigenvalues of the crystalline Frobenius acting on $\Dc(\rhob_z)$ are distinct and are
$(p^{\kappa_1}(z) F_1(z),\dots, p^{\kappa_d}(z) F_d(z))$.
\item[(v)]  For $C$ a non-negative integer, let $Z_C$ be the set
$$\ \ \ \ \{z \in Z,\  
|\kappa_I(z) - \kappa_J(z)| >C\  \ \forall I,J \subset \{1,\dots,d\},\  |I|=|J| > 0,\  I \neq J\},$$ 
where $\kappa_I = \sum_{n \in I} \kappa_n$. Then $Z_C$ accumulates at any point of $Z$ for all $C$ (see~\S\ref{accumulate}).

%({\bf Pour prouver que les familles construites section 6 verifient (v), il faudra et suffira certainement
%verifier (v) quand $X=\C_p^d$, $Z=\Z^{+}$. Dans ce cas, il suffit de noter que 
%$\{z, \kappa_{n+1}(z)-\kappa_{n}(z) > d^2(\kappa_{n}(z)-\kappa_{n-1}(z)) + C\}$ est dans $Z_C$, car cet ensemble est deja Zariski-dense et s'acummule en tout point de $Z$. Pour voir que c'est dans $Z_C$,
%prendre $I$ et $J$ comme dans la definition, appeler $n+1$ le plus grand des elements de la difference symmetique de $I$ et $J$ et suppposer qu'il est dans $I$, 
%Aloirs $\kappa_I-\kappa_J=\kappa_{n+1}   $ une somme alterne de $\kappa_m$, $m \leq n$,
%de "poids total -1", et de cardinalite $2|I|-1$. On peut ecrire cette somme quitte a rajouter des 
%$-\kappa_k+\kappa_k$, comme $-\kappa_n +$ une somme de au plus $d^2$ difference $\kappa_m-\kappa_{m-1}$ avec $m \leq n$. Ces termes sont plus petits que $\kappa_n-\kappa_{n-1}$. Alors $|\kappa_I-\kappa_J| \geq |\kappa_{n+1}-\kappa_n|- d^2 |\kappa_{n}-\kappa_{n-1}| > C$.})

\item[($\ast$)] For each $n$, there exists a continuous character $\Z_p^* 
\longrightarrow \OO(X)^*$ whose derivative at $1$ is the map $\kappa_n$ and whose evaluation at any point $z \in Z$ is the elevation to the $\kappa_n(z)$-th power.
\end{itemize}
\end{definition}
The data (a) to (c) are called {\it a refinement} of the family $(X,T)$. 

\begin{definition} \label{defraffz} Fix a refined family as above and let $z \in Z$. The (distinct) eigenvalues of $\varphi$ on $\Dc(\rhob_z)$ are naturally ordered by setting $$\varphi_n(z):=p^{\kappa_n(z)}F_n(z), \, \, \, n \in \{1,\cdots,d\},$$
which defines a refinement $\Ref_z$ of the 
representation $\rhob_z$ in the sense of \S\ref{refinement}. 
\end{definition}

\begin{example} \label{examplefam}The main examples of refined families 
arise from eigenvarieties\footnote{In particular, their construction if mostly {\it global} at the moment.}. A refined family is said to be {\it ordinary} if $|F_n(x)|=1$ for each $x \in X$ and $n \in \{1,\dots,d\}$. Many ordinary families (in the context of example \ref{examplegg} (ii)) have been constructed by Hida. In this case we could show that $T_{|\Gp}$ is a sum of $1$-dimensional families. Non ordinary refined families of dimension $2$ have been first constructed by Coleman in \cite{col2} (see also \cite{CM}, \cite{Maz}), and in this case $T_{|\Gp}$ is in general irreducible. Examples of non ordinary families of any dimension $d>2$ have been constructed by one of us in \cite{Ch}.
\end{example}

Let us do some remarks about Definition \ref{defrefi}.

\begin{remark} \label{remdeffamraf}
\begin{itemize}
\item[(i)]({\it Weights}) We stress that condition (iii) is not automatic, even up to a renumbering of the $\kappa_n$. Condition (v) impose that the Hodge-Tate weights (and their successive differences) vary a lot on $Z$. Condition (*) appears for the same reason as in \S\ref{directhyp}.  
\item[(ii)] ({\it Frobenius eigenvalues}) Assumption (iv) means that the eigenvalues of the crystalline Frobenius $\varphi$ acting on $\Dc(\rho_z)$ do not vary analytically on $Z$, but rather that they do when appropriately normalized. Note that when $d>1$, even if some eigenvalue vary analytically, i.e. if some $\kappa_n$ is constant, then the others do not by assumption (v). Moreover, because of the fixed ordering on the $\kappa_n$ by assumption (iii), $(\{\kappa_n\},\{F_{\sigma(n)}\},Z)$ is not a refinement of the family $(X,T)$ when $\sigma\neq 1\in \got{S}_d$.
\item[(iii)]({\it Generic non criticality}) Let $Z_{{\rm num}} \subset Z$ be the subset consisting of points $z \in Z$ such that $\Ref_z$ is numerically non critical in the sense of Remark \ref{remcrit} formula (\ref{num}). Then $Z_{\rm num}$ is Zariski-dense in $X$ (use (v) and the fact that around each point of $X$, each $|F_n|$ is bounded). In particular, the $\Ref_z$ are "generically" non critical in the sense of \S\ref{crit}. 
%\item[(iv)] Another fact of interests is that by (iv) and (v), any $z \in Z$ is the limit of a sequence $(z_n)_{n\geq 1}$, $z_n \in Z_n$, where the Newton polygon and the Hodge polygon of $\Dc(\rhob_{z_n})$ tend to coincide when $n$ grows.
\item[(iv)]({\it Subfamilies}) If $(X,T)$ is a refined family, and if $T$ is the sum of two pseudocharacters $T_{1}$ and $T_{2}$, then under mild conditions $T_{1}$ and $T_{2}$ 
inherits the refinement of $T$. See Prop.~\ref{propsubfamilies} below.
\end{itemize}
\end{remark}

It will also be useful to introduce the notion of {\it weakly refined} families (resp. of {\it weak refinement} of a family).

\begin{definition} \label{weakreffam}A {\it weak refinement} of a family $(X,T)$ of dimension $d$ is the data of 
\begin{itemize} 
\item[(a)] $d$ analytic functions $\kappa_n \in \anneau(X)$,
\item[(b)] an analytic function $F \in \anneau(X)$,
\item[(c)] a Zariski dense subset $Z \subset X$.
\end{itemize}
subject to the following requirements
\begin{itemize}
\item[(i), (ii)] as in Definition \ref{defrefi}. 
\item[(iii)] If $z \in Z$, then $\kappa_1(z)$ is the smallest weight of $\rhob_z$.
\item[(iv)]   For $C$ a non-negative integer, let 
$Z_C=\{z \in Z,\ \ \forall n \in \{2,\dots,d\},\  
\kappa_n(x) > \kappa_1(z)+C\}$. Then $Z_C$ accumulates at any point of $Z$ for all $C$.
\item[(v)] $\varphi_1(z):=p^{\kappa_1}(z) F_1(z)$ is a multiplicity-one eigenvalue of the crystalline Frobenius acting on $\Dc(\rhob_z)$.
\item[($\ast$)] There exists a continuous character $\Z_p^* 
\longrightarrow \OO(X)^*$ whose derivative at $1$ is the map $\kappa_1$ and whose evaluation at any point $z \in Z$ is the elevation to the $\kappa_1(z)$-th power. As in Def. \ref{defkappa}, we denote also by $\kappa_1: \G_p \longrightarrow \OO(X)^*$ the associated continuous character.
\end{itemize}
\end{definition}

\begin{remark} The conditions (i) to (v) and ($\ast$) are invariant by any permutation
in the order of the weights $\kappa_2,\dots,\kappa_d$ (not $\kappa_1$). Two weak refinements differing 
only by such a permutation should be regarded as equivalent.
\end{remark}

\subsubsection{Exterior powers of a refined family are weakly refined}

\label{weakly}
Let $(X,T)$ be a family of $p$-adic representations of dimension $d$. For $k\leq d$, then $(X,\Wk T)$ is a family of $p$-adic representations of dimension $\dwk$ (see~\S\ref{tenspseudo}), and we have $(\Wk T)_x= \tr (\Wk \rhob_x)$ for any $x \in X$. \ps

Assume that $(X,T,\{\kappa_n\},\{F_n\},Z)$ is refined. The Hodge-Tate-Sen weights of $\Wk T$ are then the $\kappa_I=\sum_{j \in I} \kappa_j$ where $I$ runs among the subsets of cardinality $k$ of $\{1,\dots,d\}$. Moreover, the $\Wk \rhob_z$ are crystalline for $z \in Z$. However, there is no natural refinement on $(X,\Wk T)$ in general \footnote{For one thing, there is no natural order on the set of subsets $I$ of $\{1,\dots,d\}$ 
of cardinality $k$ that makes the application $I \longrightarrow \kappa_I(z)$ increasing for all $z \in Z$.}. We set $$F:=\prod_{j \in \{1,\dots,i\}} F_j, \, \,  \kappa'_1:=\kappa_{\{1,\dots,k\}}=\kappa_1+\cdots+\kappa_i,$$ 
and $\kappa'_2, \dots, \kappa'_{\dwk}$ any numbering of the $\kappa_I$ for $I$ any subset of $\{1,\dots,d\}$  
of cardinality $k$ which is different from $\{1,\dots,k\}$. The following lemma is clear.

\begin{lemma}\label{wegderef=wref} The data $(\kappa'_1,\dots,\kappa'_{\dwk},F,Z)$ is a {\it weak refinement} on the family $(X,\Wk T)$.
\end{lemma} 
\subsection{Existence of crystalline periods for weakly refined families}
\label{kisinfamily}

\subsubsection{Hypotheses}

\label{hypkisinfamily}

In this subsection, $(X,T,\kappa_1,\dots,\kappa_d,F,Z)$ is a family of dimension $d$ of {\it weakly} refined $p$-adic representations. \ps

Fix $x \in X$. As in section~\ref{kisin} we shall denote by $A$ the rigid analytic local
ring $\anneau_{x}$, by $m$ its maximal ideal, and by $k=A/m=k(x)$ its residue field. We still denote by $T$ the composite pseudocharacter $G \longrightarrow \anneau(X) \longrightarrow A$. Our aim in this section is to prove that the infinitesimal pseudocharacters $T: G \longrightarrow A/I$, $I \subset A$ an ideal of cofinite length, have some crystalline periods in a sense we precise below. For this, we will have to make the following three hypotheses on $x$. \ps

\begin{itemize}
\item[(ACC)] The set $Z$ accumulates at $x$\footnote{This hypothesis is probably unnecessary but to remove it would 
require quite a big amount of supplementary work, such as a global generalization of what was done on in section~\ref{pseudocharacters} (that is
on $X$ instead of $A$). Note that any $z \in Z$ satisfies (ACC). Moreover, in the applications to eigenvarieties, (ACC) will be 
satisfied for all the $x$'s corresponding to  
$p$-adic finite slope eigenforms whose weights are in $\Z_p$, which is
more than sufficient for our needs.}.
\item[(MF)] $T$ is residually multiplicity free. 
\item[(REG)] $\Dc^+(\rhob_x(\kappa_1(x)))^{\varphi=F(x)}$ has $k(x)$-dimension $1$.
\end{itemize}\ps

Recall from Definition~\ref{defmult1} that (MF) means that $$\rhob_x=\oplus_{i=1}^r \rhob_i$$ 
where the $\rhob_i$ are absolutely irreducible, defined over $k(x)$, and two by two distinct\footnote{The results of this section will apply also in the case where the the $\rhob_i$ are not defined over $k(x)$. Indeed, it suffices to apply them to the natural weakly refined family on $X \times_{\Qp} L$, $L$ any finite extension of $\Q_p$ on which the $\rhob_i$ are defined.}. In particular, this holds of course when $\rhob_x$ is irreducible and defined over $k(x)$. As in~\S\ref{hypomfree} we shall note $d_i = \dim_k \rhob_i$, so that $\sum_{i=1}^r d_i=d$. Note that $A$ is a henselian ring (\cite[\S2.1]{berk}) and a $\Q$-algebra. In particular, $d!$ is invertible in $A$, and $T: A[G] \longrightarrow A$ satisfies the hypothesis of~\S\ref{hypomfree}.\par
Note moreover that hypothesis (REG) (for ``regularity'') is, as (MF), a kind of multiplicity free hypothesis. Indeed, Theorem~\ref{kisinnonlibre} implies easily (see below) that for 
any $x$ satisfying (ACC),  
$\Dc^+(\rhob_x(\kappa_1(x)))^{\varphi=F(x)}$ 
has $k(x)$-dimension {\it at least} $1$. 

\begin{remark} \label{remj}  The assumptions above define a $j \in \{1,\dots,r\}$ as follows. By property (REG), $F(x)$ is a multiplicity one eigenvalue of $\varphi$ on 
$$\Dcp(\rhob_x(\kappa_1(x)))= \Dcp(\rhob_1(\kappa_1(x))) \oplus \dots \oplus \Dcp(\rhob_r(\kappa_1(x))).$$ Hence this is an eigenvalue of $\varphi$ on one (and only one) of the
$\Dcp(\rhob_i(\kappa_1(x)))$ say $\Dcp(\rhob_j(\kappa_1(x)))$, which defines a unique $j \in \{1,\dots,r\}$. \end{remark}

\subsubsection{The main results}\label{resultsweakfam} 
We will use below some notations and concepts introduced in section \ref{pseudocharacters}.
Let $K$ be the total fraction ring of $A$ and let $\rho: A[G] \longrightarrow M_d(K)$ be a representation whose trace is $T$ and whose kernel is $\Ker T$. It exists by Theorem 
\ref{structure} (ii) and Remark \ref{remthmstru} as $A$ is reduced and noetherian. Fix a GMA structure on $S:=A[G]/\Ker T$ given by the theorem cited above, $j$ as defined in Remark \ref{remj}, and let $M_j \subset K^d$ the "column" $S$-submodule defined in \S\ref{Mi}. It is finite type over $A$ by construction and Remark \ref{remthmstru}. \par
Let moreover $\PP$ be a partition of $\{1,\dots,r\}$. Recall that if $\PP$ contains $\{j\}$, then for every ideal $I$ containing the reducibility ideal $I_\PP$ (see \S\ref{assredloc}), there is a unique continuous representation 
$$\rho_j : G \longrightarrow \Gl_{d_j}(A/I),$$ whose reduction mod $m$ 
is $\rhob_j$ and such that $T \otimes A/I = \tr \rho_j + T'$, where
$T':G \longrightarrow A/I$ is a pseudocharacter of dimension $d-d_j$ (see Definition~\ref{lesrhoi},
Proposition~\ref{topoprop}). 

\begin{theorem} \label{Drhojlibre}
Assume that $\PP$ contains $\{j\}$ and let $I$ be a cofinite length ideal of $A$ containing $I_\PP$. Then $\Dcp(\rho_j(\kappa_1))^{\varphi=F}$ and $\Dcp(M_j/IM_j(\kappa_1))^{\varphi=F}$ are free of rank one over $A/I$. 
\end{theorem}

\begin{pf} We will prove the proposition assuming the following lemma, whose proof is postponed 
to the next subsection.
\begin{lemma} \label{fact}Let $I$ be a cofinite length ideal of $A$, then
\begin{itemize}
\item[(i)] The Sen operator of $\Ds(M_j/IM_j)$ is annihilated by $\prod_{i=1}^{d}(T-\kappa_i)$,\par
\item[(ii)] l$\left(\Dc^+(M_j/IM_j(\kappa_1))^{\varphi=F}\right)=l(A/I)$.
\end{itemize}
\end{lemma}
\noindent By Theorem~\ref{extensions2}(0), there is an exact sequence of $(A/I)[G]$-modules
$$0 \longrightarrow K \longrightarrow M_j/IM_j \longrightarrow \rho_j \longrightarrow 0$$ 
where $K$ has a Jordan-Holder sequence whose all subquotients are 
isomorphic to $\rhob_i$ for some $i \neq j$. If $X$ is a finite length $A$-module equipped with a continuous $A$-linear action of $\Gp$ we set $D(X):=\Dcp(X(\kappa_1))^{\varphi=F}$. As $D(\rhob_i)=0$ for $i \neq j$ by (REG), we have $D(K)=0$, hence applying 
the left exact functor $D$ to the above sequence, we get an 
injection $$D(M_j/IM_j) \hookrightarrow D(\rho_j).$$
Thus by Lemma~\ref{fact} (ii) we have $l(D(\rho_j)) \geq l(A/I)$.
Applying Lemma~\ref{abstract}(i) to the $A/I$-representation $\rho_j$
 gives $l(D(\rho_j))=l(A/I)$,
hence an isomorphism $$D(M_j/IM_j) \simeq D(\rho_j),$$ 
and case (2) of Lemma~\ref{abstract2}
gives that $D(\rho_j)$ is free of rank $1$ over $A/J$. Hence the result.
\end{pf}

\begin{theorem}\label{redlocfaible} Assume that $\rhob_x$ has distinct Hodge-Tate-Sen weights and that the weight $k$ of $\Dc^+(\rhob_j(\kappa_1(x)))^{\varphi=F(x)}$ 
is the smallest integral Hodge-Tate weight of $\rhob_j(\kappa_1(x))$. Then
for the unique $l$ such that $\kappa_l(x)-\kappa_1(x)=k$,
we have that $\kappa_l$ is a weight of $\rho_j$ and
$$(\kappa_l - \kappa_1) - (\kappa_l(x) - \kappa_1(x)) 
\in I_\PP \text{ \ if } \{j\} \in \PP.$$
\end{theorem}
\begin{pf} Let $I \supset I_{\PP}$ be a cofinite length ideal of $A$. By Theorem~\ref{Drhojlibre}, 
$$\Dc^+(\rho_j(\kappa_1))^{\varphi=F}$$ is free of rank one over $A/I$. Moreover, $k$ is the smallest integral Hodge-Tate weight of 
$\rhob_j(\kappa_1(x))$. 
Thus we can apply Proposition~\ref{poidsconstant} to 
$V:=\rhob_j(\kappa_1)$ which shows that $V$ has a constant weight $k$, {\it i.e.} that $(V\otimes_{\Qp}\C_p)^{\Hp}$ contains a free $A/I$-submodule of rank $1$ on which the Sen operator acts as the multiplication by $k$. By lemma \ref{fact} (i), this implies that $$\prod_{i=1}^d(k-(\kappa_i-\kappa_1))=0 \,\, \, \, \, \,  {\rm in  }\, \, \, A/I.$$
The difference of any two distinct terms of the product above is invertible in (the local ring) $A/I$ as $\kappa_i(x)\neq \kappa_{i'}(x)$ if $i\neq i'$. Hence one, and only one, of the $k-(\kappa_i-\kappa_1)$ above is zero, and reducing mod $m$ gives $i=l$. In particular, $\kappa_l$ is a Hodge-Tate-Sen weight of $\rho_j$ and
$$k=\kappa_i-\kappa_1=\kappa_l(x)-\kappa_1(x) \,\, \, \, \, \,  {\rm in  }\, \, \, A/I.$$
 We conclude the proof as $I_{\PP}$ is the intersection of the $I$ of cofinite length containing it, 
by Krull's theorem.
\end{pf}
\begin{remark} 
\begin{itemize} 
\item[(i)] The conclusion of the theorem can be rephrased as :
$\kappa_l-\kappa_1$ is constant on the reducibility locus corresponding to 
$\PP$, if $\PP$ contains $\{j\}$.
\item[(ii)] The hypothesis that $k$ is the smallest weight is satisfied 
in many cases. For one thing, it is obviously satisfied when $k$ is the only 
integral Hodge-Tate weight of $\rhob_j(\kappa_1(x))$, which is the generic situation. More interestingly, it is also satisfied for crystalline $\rhob_x$ whenever 
$v(F(x))$ is smaller than the {\it second} (in the increasing order) Hodge-Tate weight of $\rhob_j(\kappa_1)$ since, by weak admissibility, 
$k \leq v(F(x))$. This is always true when $\rhob_j$ has 
dimension $\leq 2$, since by admissibility, the second (that is, the greatest) 
weight is greater than the valuation of {\it any} eigenvalue of the Frobenius.
\item[(iii)] The assumption that $\rhob_x$ has distinct Hodge-Tate weights implies that $\rhob_x$ has no multiple factors, hence (MF) if these factors are defined over $k(x)$. 
\end{itemize}
\end{remark}

Now let $i \neq j$ be an integer in $\{1,\dots,r\}$. Recall that if $\PP$ contains $\{i\}$ and $\{j\}$, and if $I$ contains $I_\PP$,
then (see Theorem~\ref{extensions} , Theorem~\ref{extensions2}(1) and Proposition~\ref{topoprop})
there is a map $\iota_{i,j}$ whose image is 
$\Ext^1_{S/JS, {\rm cont}}(\rho_j,\rho_i)$. 
\begin{theorem}\label{crystext} Assume that $\PP$ contains $\{i\}$ and $\{j\}$ and let $I$ be a cofinite length ideal of $A$ containing $I_\PP$.
Let $\rho_c : G \longrightarrow \Gl_{d_i+d_j}(A/I)$
 be an extension of $\rho_j$ by $\rho_i$ which belongs to the image of $\iota_{i,j}$.
Then $\Dc^+(\rho_c(\kappa))^{\varphi=F}$ is free of rank one over $A/I$. 
\end{theorem}
\begin{pf} The proof is exactly the same as the proof of 
Theorem~\ref{Drhojlibre} except that we start using point (2) of 
Theorem~\ref{extensions2} instead of point (0). 
\end{pf}

\subsubsection{Analytic extension of some $A[G]$-modules, and proof of Lemma \ref{fact}}\label{analytic}

We keep the assumptions and notations of \S\ref{resultsweakfam}. Let $M \subset K^d$ by any 
$S$-submodule which is of finite type over $A$. 

\begin{lemma} \label{prolon} There is an open affinoid subset $U$ of $X$ containing 
$x$ in which $Z$ is Zariski-dense and a torsion-free coherent sheaf $\MM$ on $U$ with a 
continuous action of $G$ such that 
$\MM(U) \otimes_{\anneau(U)} A \simeq M$ as $A[G]$-modules and topological $A$-modules.

	If moreover $K.M=K^d$, we may choose $U$ and $\MM$ such that $\MM(U)\otimes_{\OO(U)}{\rm Frac}(\OO(U))$ 
is free of rank $d$ over ${\rm Frac}(\OO(U))$, and carries a semisimple representation of $G$ with trace $T \otimes_{\OO(X)} 
\OO(U)$.
\end{lemma} 
\begin{pf}
By (ACC), we may choose a basis of open affinoid neighbourhoods $(V_i)_{i \in I}$ 
of $x \in X$ such that $Z$ is Zariski-dense in $V_i$ for each $i$. We may view $I$ as a directed set if we set $j\geq i$ if $V_j \subset V_i$, 
and then $\underset{i }{\rm indlim}\OO(V_i)=A$.

By construction we have $\tr(\rho(G)) \subset \OO(X)$. As each $\OO(V_i)$ is reduced and noetherian, a standard argument 
implies that the $\OO(V_i)$-module 
$$\OO(V_i)[G]/\ker (T \otimes_{\OO(X)}\OO(V_i))$$ is of finite type (see e.g. \cite[Lemma 7.1 (i)]{BC}). As a consequence, its quotient 
$\OO(V_i)[\rho(G)] \subset M_d(K)$ is also of finite type over $\OO(V_i)$.

As $M$ is finite type over $A$, we can find an element $0 \in I$ and a finite type $\OO(V_0)$-submodule $M_0$ of $M$ such that $A M_0=M$. 
We define now $N_0$ as the smallest 
$\anneau(V_0)$-submodule of $M$ containing $M_0$ and 
stable by $G$. It is finite type over $\OO(V_0)$ as
we just showed that 
$\OO(V_0)[\rho(G)] \subset M_d(K)$ is. 
Moreover, the map 
$G \longrightarrow {\rm Aut}_{\OO(V_0)}(N_0)$ 
(resp. $G \longrightarrow {\rm Aut}_{A}(M)$) is 
continuous by \cite[Lemma 7.1 (v)]{BC} (resp. by its proof).

For $i\geq 0$, we set $N_i=\OO(V_i)N_0 \subset M$. The following abstract lemma implies that for $i$ big enough, the morphism 
$N_i \otimes_{\anneau(V_i)} A \longrightarrow M$
is an isomorphism. We fix such an $i$, set $U=V_i$ and define $\MM$ as the coherent sheaf on $V_i$ whose global sections are $N_i$. 
It is torsion free over $\OO(V_i)$ as $N_i \subset M \subset K^d$, which concludes the proof of the first assertion.

Assume moreover that $K.M=K^d$ and let $N_i \subset K^d$ the module constructed above, so $K.N_i=K^d$. 
The kernel of the natural map $$N_i \otimes_{\OO(V_i)} {\rm Frac}(\OO(V_i)) \longrightarrow N_i \otimes_{\OO(U_i)} K = K^d$$ 
is exactly supported by the minimal primes of the irreducible components of $\OO(V_i)$ 
that do not contain $x$, 
and at the other minimal primes $N_i$ is free of rank $d$ with trace $T$, and it 
is semisimple because so is its scalar extension to $K$ by construction and Lemma~\ref{lemmeutile} (i) below. Let $U' \subset V_i$ be 
the Zariski open subset of $V_i$ whose complement is the (finite) union of irreducible components of $V_i$ not containing $x$. 
Choose $j\geq i$ such that $V_j \subset U'$, then $U:=U_j$ and $\MM(U):=N_j$ have all the required properties.
\end{pf}

\begin{lemma} \label{indlim} Let $(A_i)_{i \in I}$ be a directed family of commutative rings and let $A$ be the inductive limit of $(A_i)$. Assume $A$ is noetherian. Let $M$ be a finite type $A$-module and $N_0$ a finite-type $A_0$-submodule of $M$
such that $A\, N_0=M$. For $i\geq 0$, set $N_i:=A_iN_0 \subset M$. \par Then for $i$ big enough, the natural morphism 
$N_i \otimes_{A_i} A \longrightarrow M$ is an isomorphism. 
\end{lemma}

\begin{pf} 
Define $K_i$ by the following exact sequence :
$$0 \longrightarrow K_i \longrightarrow N_i \otimes_{A_i} A \longrightarrow M \longrightarrow 0.$$ 

For $i \leq j$, we have a commutative diagram 
$$\xymatrix{
0 \ar[r] & K_i \ar[d] \ar[r] & N_i \otimes_{A_i} A \ar[r] \ar[d] & 
M \ar[r]\ar[d] & 0 \\
0 \ar[r] & K_j  \ar[r] & N_j \otimes_{A_j} A \ar[r] & 
M \ar[r] & 0 \\}$$
The horizontal lines are exact sequences, the right vertical arrow is the identity and the middle one is surjective by the 
associativity of the tensor product. Hence the left vertical arrow 
$K_i \longrightarrow K_j$ is surjective. Because $K_0$ is a finite type
 $A$-module, and $A$ is noetherian, there is an $i$ such that for each 
$j \geq i$, $K_i \longrightarrow K_j$ is an isomorphism.

Let $x \in K_i$. We may write $x = \sum_k n_k \otimes a_k$ with $n_k \in N_i$ and $a_k \in A$, and $\sum_k n_k a_k =0$ in $M$.
Take $j \geq i$ such that all the $a_k$'s are in $A_j$. Then the 
image of $x$ in $N_j \otimes_{A_j} A$ is $0$, and $x$ is $0$ in $K_j$.
But then $x=0$ in $K_i$, which proves that 
$K_i=0$ and the lemma.
\end{pf}

\begin{lemma}\label{lemmeutile} \begin{itemize}
\item[(i)] $S\otimes_A K$ is a semisimple $K$-algebra. 
\item[(ii)] There exists a finite type $S$-module $N \subset K^d$ such that
$(N \oplus M_j)K=K^d$ and that $(N\otimes_A k)^{\ses}$ is isomorphic to a sum of copies of $\rhob_i$ with $i\neq j$.
\end{itemize}
\end{lemma}

\begin{pf} Recall from \S\ref{resultsweakfam} that $S=A[G]/\ker T$. Since $K \supset A$ is a fraction ring of $A$, we have $\Ker (T\otimes_A K)=K.\Ker T$ in $K[G]$. As a consequence, the natural map $S \otimes_A K
\longrightarrow K[G]/\Ker(T\otimes_A K)$ is an isomorphism, and Lemma \ref{radical} proves (i).\par
Let us show (ii). By (i) we can chose a $K.S$-module $N' \subset K^d$ such that $K.M_j \oplus N'= K^d$. As $S$ is finite type over $A$ 
by Remark \ref{remthmstru}, we can find a $S$-submodule $N \subset N'$ such that $N$ is finite type over $A$ and $K.N=N'$. We claim that $N$ has the required property. By construction we only have to prove the assertion about $(N\otimes_A k)^{\ses}$. Arguing as in the proof of Theorem \ref{extension2} (0), it suffices to show that $e_jN=0$, where $e_j$ is as before the idempotent in the fixed GMA structure of $S$. But $e_j(K^d) = e_j (K.M_j)$ by definition of $M_j$ and Theorem \ref{structure} (ii). So $e_j(K.N)=0=e_jN$, and we are done.
\end{pf}

We are now ready to prove Lemma \ref{fact}. 
%It will turn out that part (i) holds in fact not only for $M_j$ but for all the modules $M$ as chosen above, and part (ii) for all such modules satisfying moreover that $\Dc(M\otimes(\kappa_1(x))^{\ses})^{\varphi=F(x)}$ has $k$-dimension $\leq 1$. 
\ps
\begin{pf} (of Lemma \ref{fact}). Let us show (ii) first. We set $M=N\oplus
M_j$, where $N$ is given by Lemma \ref{lemmeutile}. \par
By Proposition~\ref{extensions2}(0) and Lemma \ref{lemmeutile} (ii), $(M \otimes k)^\ses \simeq
\oplus_{i=1}^r n_i \rhob_i$ where 
the $n_i$ are natural integers and $n_j=1$.
But by (REG) $\Dc^+(\rhob_i(\kappa_1(x))^{\varphi=F(x)}$ 
has dimension $\delta_{i,j}$. In particular, 
\begin{equation} \label{dcrisbon1} \dim_k(\Dc^+((M \otimes
k)(\kappa(x))^\ses)^{\varphi=F(x)})=1.
\end{equation} 
Moreover,
$\Dcp(M/IM(\kappa_1))=\Dcp(M_j/IM_j(\kappa_1))\oplus\Dcp(N/IN(\kappa_1))$,
and
\begin{equation} \label{dcrisbon2} \Dcp(N/IN(\kappa_1))^{\varphi=F}=0 \end{equation}
by a devissage and the same argument as above.\par

We claim now that the equality follows directly from Theorem~\ref{kisinnonlibre} applied 
to the module $\MM$ over $U$ associated to $M$ given by Lemma \ref{prolon} (applied in the case $K.M=K^d$).
By formula (\ref{dcrisbon2}), we just have to verify that $\MM$ satisfies the hypotheses 
(CRYS), (HT) and ($\ast$) of \S\ref{directhyp}, and we already checked that 
$\Dc^+((M \otimes k)(\kappa(x))^\ses)^{\varphi=F(x)}$ has length one in
(\ref{dcrisbon1}).

By assumption (iv) of weakly refined families, $Z_C \cap U$ accumulates at every point of 
$Z \cap U$. As $\MM(U)$ is torsion free of generic rank $d$ and with trace $T$, and by the generic flatness theorem, 
there is a proper
Zariski closed subspace $F$ of $U$
such that for $y \in U-F$, $\MB_y^{\ses} = \rhob_y$. Recall that 
the ${\rm Frac}(\OO(U))[G]$-module $\MM \otimes_{\OO(U)}{\rm Frac}(\OO(U))$ is semisimple. So enlarging $F$ is necessary, we have that for $y \in U-F$,
$\MB_y=\MB_y^{\ses}$, hence $\MB_y = \rhob_y$. We replace $Z$ by $(Z \cap U) - (F \cap Z \cap U)$, so by (ACC) $Z$ is a Zariski dense subset of $U$ 
and still has the property that $Z_C$ accumulates at any point of $Z$. Property (CRYS) follows then from (ii) and (v) of the definition of a weak refinement, and property (HT) from (iii) and (iv). This concludes the proof. \ps

Let us show (i) now. If $E$ is a $\Qp$-Banach space, we set\footnote{All the $\Qp$-Banach spaces of this proof to which we apply the functor $-_{\C_p}$ are discretely normed. We use freely the fact that any continuous closed injection $E \longrightarrow F$ between such spaces induces an exact sequence $0 \longrightarrow E_{\C_p} \longrightarrow F_{\C_p} \longrightarrow
(F/E)_{\C_p} \longrightarrow 0$ by 
\cite[1.2]{ser2}, and also that any submodule of a finite type module over an affinoid algebra is closed.}
$E_{\C_p}:=E\widehat{\otimes}_{\Qp}\C_p$. Recall that Sen's theory \cite{sen} attaches in particular to any continuous morphism $\tau: \Gp \longrightarrow \cal{B}^*$, $\cal{B}$ any
Banach $\Q_p$-algebra, an element
$\varphi \in \cal{B}_{\C_p}$ whose
formation commutes with any continuous Banach
algebra homomorphism $\cal{B} \longrightarrow \cal{B}'$. When $\tau$
is a finite dimensional $\Qp$-representation of $\Gp$, this element is the usual Sen operator. Applying this to the Banach algebra
$$\cal{B}:=\End_{\OO(U)}(\MM(U))$$ we get such an element $\varphi$. \par 
We claim that $\varphi$ is killed by the polynomial 
$$P:=\prod_{i=1}^d(T-\kappa_i).$$
Indeed, arguing as in the proof of (ii) above me may assume that for all 
$z \in Z$ we have $\MB_z  \simeq \rhob_z$ and $\cal{B}/m_z\cal{B}\simeq
\End_{k(z)}(\MB_z)$. As a consequence, using the
evaluation homomorphism $\cal{B} \longrightarrow \cal{B}/m_z\cal{B}$ 
and assumption (i) in Definition \ref{weakreffam}, we get that
$P(\varphi) \in m_z \cal{B}_{\C_p}$. But  
$\OO(U)_{\C_p}$ is reduced by \cite[Lemma 3.2.1 (1)]{con}, so $\cal{B}_{\C_p}$ is a (finite type) 
torsion free $\OO(U)_{\C_p}$-module.
Since $Z$ is Zariski-dense in $U$, hence in $U(\C_p)$, and since affinoid algebras are Jacobson rings, we obtain that $P(\varphi)=0$ in $\cal{B}_{\C_p}$. We
conclude the proof as the operator of the statement 
of Lemma \ref{fact}(i) is the image of $\varphi$ under $\cal{B}_{\C_p} \longrightarrow
\End_{A/I}(M_j/IM_j)_{\C_p}$. \end{pf}
  
\subsection{Refined families at regular crystalline points}

\subsubsection{Hypotheses}

\label{hypredloc}
In this subsection, $(X,T,\kappa_1,\dots,\kappa_d,F_1,\dots,F_d,Z)$ 
is a family of dimension $d$ of refined $p$-adic representations. We fix $z \in Z$ (and not only in $X$). 
As in~\S\ref{hypkisinfamily} we write $A=\anneau_{z}$
and still denote by $T$ the composite pseudocharacter 
$G \longrightarrow \anneau(X) \longrightarrow A$. We assume moreover that $T$ is residually multiplicity free, and we use the same notation as before:
$$\rhob_z = \oplus_{i=1}^r \rhob_i, \, \, \, d_i = \dim \rhob_i.$$ 
Recall from Definition \ref{defraffz} that $\rho_z$ is equipped with a refinement $$\Ref_z=(\varphi_1(z),\cdots,\varphi_d(z))$$ satisfying $\varphi_n(z)=p^{\kappa_n(z)}F_n(z)$. As $\Dc(\rhob_z)=\oplus_{i=1}^r \Dc(\rhob_i)$ this refinement induces for each $i$ a refinement of $\rhob_i$ that we will denote by $\Ref_{z,i}$. We will make the following hypotheses on $z$.
\begin{itemize}
\item[(REG)] The refinement $\Ref_z$ is {\it regular} (see Example~\ref{reg}) :
for all $n \in \{1,\dots,d\}$, 
$p^{\kappa_1(z)+\dots+\kappa_n(z)}F_1(z)\dots F_n(z)$ is an eigenvalue of $\varphi$ on 
$\Dc(\Lambda^n \rhob_z)$ of multiplicity one.
\item[(NCR)] For every $i \in \{1,\dots,r\}$, $\Ref_{z,i}$ is a non-critical refinement (cf.~\S\ref{crit}) of $\rhob_i$.
\item[(MF')] For every family of integers $(a_i)_{i=1,\dots,r}$ with $1 \leq a_i \leq d_i$, the representation
$\rhob_{(a_i)} := \bigotimes_{i=1}^r \Lambda^{a_i}\rhob_i$ is absolutely irreducible. Moreover, if $(a_i)$ and $(a'_i)$ are two distinct sequences as above with $\sum_{i=1}^r a_i=\sum_{i=1}^r a'_i$, then 
$\rhob_{(a_i)} \not \simeq \rhob_{(a'_i)}.$
\end{itemize}

Note that the hypothesis (NCR) does not mean at all that the refinement of $\rhob_z$ is noncritical : if for example $d=r$, that is the $\rhob_i$ are characters, any refinement of $\rhob_z$ satisfies (NCR).

Although it does not seem possible to weaken significantly the hypotheses (REG), (NCR) in order to prove Theorem~\ref{thmredloc} below, hypothesis (MF') is probably unnecessary. It is equivalent 
to the assertion that for all $k \in \{1,\dots,d\}$, $\Lambda^k T$ is a residually multiplicity free pseudocharacter with residual irreducible component the traces of the representations $\rhob_{(a_i)}$ with $\sum_{i=1}^r a_i =k$.
 
\subsubsection{The residually irreducible case {\rm ($r=1$)}.}\label{resirrcasetridef} We keep the hypotheses above. We first deal with the simplest case for which $\rhob_z$ is irreducible and defined over $k(z)$. In this case (REG) and (NCR) mean that $\Ref_z$ is a regular non critical refinement of $\rhob_z$, and (MF') that
$\Lambda^k \rhob_z$ is irreducible for each $1\leq
k \leq d$. \ps
Recall that in this residually irreducible case, 
there exists a unique continuous representation 
$\rho: G \longrightarrow \GL_d(A)$ whose trace is $T$ by the theorem of Rouquier and Nyssen (the continuity follows from Proposition \ref{topoprop} (i)). 
We define a continuous character 
$\delta: \Qp^*\longrightarrow (A^*)^d$ by 
setting : $$\delta(p):=(F_1,\cdots,F_d), \, \, \, \, \, \delta_{|\Z_p^*}=(\kappa_1^{-1},\cdots,\kappa_d^{-1}).$$
Recall that each $\kappa_n$ may be viewed as a 
character $\Zp^* \longrightarrow A^*$ 
in the same way as in Definition \ref{defkappa}, using property (*) of Definition \ref{defrefi}.
 
\begin{theorem} \label{tridefirr}For any ideal $I \subsetneq A$ of cofinite length, $\rho \otimes A/I$ is a trianguline deformation of $(\rhob_z,\Ref_z)$ whose parameter is $\delta \otimes A/I$.
\end{theorem}

\begin{pf} Fix $I$ as in the statement and $V:=\rho \otimes A/I$. By Proposition \ref{poidscst2}, it suffices to show that for each $1\leq k \leq d$, $\Dc(\Lambda^k V (\kappa_1\cdots\kappa_k))^{\varphi=F_1\cdots F_k}$ is free of rank $1$ over $A/I$. Indeed, by definition of the characters $\kappa_i$ and of the $t_i$ {\it loc. cit.} , we have $t_i=k_i$ for each $i$. \par
Fix $1\leq k \leq d$ and consider the family $(X,\Lambda^k T)$. As seen in~\S\ref{weakly}, 
this family is naturally weakly refined, with same set $Z$,  
$$ F=\prod_{n=1}^k F_n,$$
and first weight 
$$\kappa = \sum_{n=1}^k \kappa_n.$$
This weakly refined family satisfies all the hypotheses 
of \S\ref{hypkisinfamily} 
(indeed, each of the hypothesis there is an immediate 
consequence of the corresponding hypothesis of \S\ref{hypredloc}. Namely, (ACC) 
comes from the fact that $z \in Z$, (MF) from (MF'), and (REG) from (REG)). By assumption, $(\Lambda^k T)_z$ is irreducible and $\Lambda^k T \otimes A/I = \tr (\Lambda^k \rho) \bmod I$, hence we may apply Theorem~\ref{Drhojlibre} to it, and we are done.
% (in this case, the corresponding $r$ is $1$, $\PP=\{\{1\}\}$ is the trivial partition of the set with one element, hence $j=1$ and the $\rho_1$ defined there is exactly $\Lambda^k \rho \otimes A/I$).
\end{pf}

\subsubsection{A permutation}
\label{permutation} In order to study the reducible cases we need to define a permutation $\sigma$ of $\{1,\dots,d\}$ that mixes up the combinatorial data of the refinement of $\rhob_x$ and of its decomposition $\rhob_x = \rhob_1 \oplus \dots \oplus \rhob_r$. 

%We proceed by induction :
%let $k \geq 1$ and assume that $\sigma(1),\dots,\sigma(k-1)$ 
%are already defined. 
%By (MF) and (Z), there is a unique $l \in \{1,\dots,r\}$  
%such that  $p^{\kappa_k(x) F_k(x)}$ is an eigenvalue  of 
%$\Dc(\rhob_l)$.
%Define $\sigma(k) \in \{1,\dots,d\}$ such that 
%$p^{\kappa_{\sigma(k)}}(x) F_{\sigma(k)}(x)$ is the smallest weight of 
%$\rhob_l$ which is  different from 
%$p^{\kappa_{\sigma(1)}(x)} F_{\sigma(1)}(x),\dots, 
%p^{\kappa_{\sigma(k-1)(x)}} 
%F_{\sigma(k-1)}(x)$.

The refinement $\Ref_z$ together with the induced refinements $\Ref_{z,i}$ of th
$\rhob_i$'s define a partition $R_1 \coprod \dots \coprod R_r$ of $\{1,\dots,d\}$: $R_i$ is the set of $n$ such that $p^{\kappa_{n}(z)} F_n(z)$ is a 
$\varphi$-eigenvalue on $\Dc(\rhob_i)$. In the same way, we define a partition $W_1 \coprod \dots \coprod W_r$ of 
$\{1,\dots,d\}$: $W_i$ is the set of integers $n$ such that $\kappa_n(z)$ is a Hodge-Tate weight of $\rhob_i$. This is a partition as the $\kappa_n(z)$ are two-by-two distinct.

\begin{definition} \label{definitionsigma} We define $\sigma$ as the unique bijection
that sends $R_i$ onto $W_i$ and that is increasing on each $R_i$.
\end{definition}

Note that $\sigma$ does not depend on the chosed ordering on the $\rhob_i$.

\begin{example}({\it Refined deformations of ordinary representations})
Assume that $r=d$, so $\rhob_z$ is a sum of characters $\rhob_1,\dots,\rhob_d$. Since there is an obvious bijection between this set of characters and the set of eigenvalues of $\varphi$ on $\Dc(\rhob_z)$, the refinement determines an ordering of those characters. We may assume up to renumbering that this order is $\rhob_1, \dots, \rhob_d$.
By definition of the permutation above, the weights of $\rhob_1,\dots,\rhob_d$ are respectively 
$\kappa_{\sigma(1)}(z),\dots,\kappa_{\sigma(d)}(z)$. Note that in this case, $\sigma$ determines the refinement. We refer to this situation by saying that the {\it representation} $\rhob_z$ is {\it ordinary}. \par
 Assume that $\rhob_z$ is ordinary.
We say that the point $z$ (and the refinement $\Ref_z$) is {\it ordinary} if moreover $\sigma = \Id$, that is if the valuation of the eigenvalues 
in the refinement are increasing. For example, the 
families constructed by Hida (see Example \ref{examplefam}) are ordinary in this strong sense: each $z \in Z$ is ordinary. \par
When, on the contrary, $\sigma$ is transitive on 
$\{1,\dots,d\}$ we call the corresponding refinement, and the point $z$, {\it anti-ordinary}. For $d=3$, examples of families with such $z$ have been constructed and studied in~\cite{BC}.  Intermediary cases are also interesting. For example, Urban and Skinner consider in 
\cite{SU} a refined family of dimension $d=4$ with a point $z \in Z$ where $\rhob_z$ is ordinary
and $\sigma$ is a transposition. They call such a point {\it semi-ordinary}. \par
In general, let us just say that we expect that any ordinary representation and any permutation $\sigma$ should occur as a member of a refined family in the above way. 
\end{example}

\subsubsection{The total reducibility locus}

Keep the assumptions and notations of \S\ref{permutation} and \S\ref{hypredloc}. We will use again some notations and concepts introduced in section \ref{pseudocharacters}, applied to the residually multiplicity free pseudocharacter $T: A[G] \longrightarrow A$. Let $\PP$ be the finest
partition $\{\{1\},\cdots,\{r\}\}$ of $\{1,\dots,r\}$, so $I_{\PP}$ is the total reducibility ideal of $T$. Recall that for every ideal $I \subsetneq A$ containing $I_\PP$ 
(see \S\ref{assredloc}), there is for each $i$ a unique continuous representation 
$$\rho_i : G \longrightarrow \Gl_{d_i}(A/I)$$ whose reduction mod $m$ 
is $\rhob_i$ and such that $T \otimes A/I = \sum_{i=1}^r \tr \rho_i$ 
(see Definition~\ref{lesrhoi}, Proposition~\ref{topoprop}).  \ps
Let $1\leq i \leq r$ and write $R_i=\{ j_1,\cdots,j_{d_i} \}$ with $s
\mapsto j_s$ increasing. We define a 
continuous character $\delta_i: \Qp^*\longrightarrow (A^*)^{d_i}$ 
by setting :
$$\delta_i(p):=(F_{j_1}p^{\kappa_{j_1}(z)-\kappa_{\sigma(j_1)}(z)},\dots,
F_{j_{d_i}}p^{\kappa_{j_{d_i}}(z)-\kappa_{\sigma(j_i)}(z)}),$$ 
$${\delta_i}_{|\Z_p^*}=(\kappa_{\sigma(j_1)}^{-1},\dots,\kappa_{\sigma(j_{d_i})}^{-1}).$$

We will need to consider the following further assumption on the partition $R_i$ defined in \S\ref{permutation}:
\label{conditionint}
\begin{itemize}
\item[(INT)] Each $R_i$ is a subinterval of $\{1,\dots,d\}$.
\end{itemize}

\begin{theorem} \label{tridefgen}Assume (INT) and let $I_{\PP} \subset I \subsetneq A$ be any cofinite length ideal. Then
for each $i$, 
$\rho_i$ is a trianguline deformation of $(\rhob_i,\Ref_{z,i})$ whose parameter is $\delta_i$.
\par Moreover, for each $n \in \{1,\dots,d\}$, we have 
$$\kappa_{\sigma(n)}-\kappa_{n}=\kappa_{\sigma(n)}(z)-\kappa_{n}(z)
\, \, \, \, {\rm in}\, \,  A/I_{\PP}.$$

\end{theorem}

\begin{pf} We argue as in the proof of Theorem \ref{tridefirr} taking in
account the extra difficulties coming from the reducible situation. By
(INT), we have for each $i$ that 
$R_i=\{x_i+1,x_i+2,\dots,x_i+d_i\}$ for some $x_i \in \{1,\dots,r\}$. 
Up to renumbering the $\rhob_i$, we may assume that $x_1=0$ and that 
$x_i=d_1+\cdots+d_{i-1}$ if $i>1$. \par
We fix $I$ as in the statement. We will prove below that 
each $\rho_i$ is a trianguline deformation of $(\rhob_i,\Ref_{z,i})$ whose
parameter $\delta'_i$ coincides with $\delta_i$ on $p$, but satisfies
$${\delta'_i}_{|\Z_p^*}=(\kappa_{j_1}^{-1}\chi^{\kappa_{j_1}(z)-\kappa_{\sigma(j_1)}(z)},
\dots,\kappa_{j_{d_i}}^{-1}\chi^{\kappa_{j_{d_i}}(z)-\kappa_{\sigma(j_{d_i})}(z)}).$$
As the Sen polynomial of $\rho_i$ is
$$\prod_{s=1}^{d_i}(T-\kappa_{\sigma(j_s)})$$ by Lemma \ref{fact} and by
definition of $\sigma$ (use the fact that the $\kappa_n(z)$ are distinct), Proposition \ref{senpol} will then conclude the
second part of the statement (argue as in the proof of Theorem \ref{redlocfaible} to go from $I$ to
$I_{\PP}$). \par
Let us prove now the result mentioned above. Fix $j \in
\{1,\cdots,r\}$ and if $j>1$ assume by induction that for each $i<j$, 
$\rho_i$ is a trianguline deformation of $(\rhob_i,\Ref_{z,i})$ whose 
parameter is $\delta'_i$ defined above. Note that $\Ref_{z,i}$ is regular by
(INT) and (REG) (see the proof below for more details about this point),
and non critical by (NCR). So by Proposition~\ref{poidscst2}, it suffices to prove that for $h=1,\dots,d_j$, 
\begin{equation} \label{trucaprouver}
\Dc((\Lambda^h
\rho_j)(\kappa_{x_j+1}+\cdots+\kappa_{x_j+h}))^{\varphi=F_{x_j+1}\cdots
F_{x_j+h}} \, \, \, \, {\rm is\, \, free \, \, of \, \, rank }\, \, 1 \, \, {\rm
over}\, \,  A/I,
\end{equation}
what we do now. \par
For $k=x_j+h$ any number in $R_j$, let 
$a_i(k) =  |R_i \cap \{1,\dots,k\}|$ for $i \in \{1,\dots,r\}$. 
In other words, we have $a_i(k)=d_i$ (resp. $a_i(k)=0$) for all $i \in
\{1,\dots,j-1\}$ (resp. for $i>j$), and $a_j(k)=h$. We want to apply 
Theorem~\ref{Drhojlibre} to the weakly refined families $\Lambda^k T$, $k\in
R_j$, as in the proof of theorem \ref{tridefirr}. We set again $F=\prod_{n=1}^k F_n$ and
$\kappa=\sum_{k=1}^n \kappa_n$. As already explained in the proof of Theorem
\ref{tridefirr}, the family $\Lambda^k T$ satisfies the assumption of \S
\ref{hypkisinfamily}. \par
We note first that the (unique by (REG))
irreducible subrepresentation of
$\Lambda^k \rhob_z$ that has the $\varphi$-eigenvalue
$p^{\kappa(z)} F(z)$ in its $\Dc$ is $\rhob_{(a_i(k))}$ with
the notations of (MF'). This representation is exactly 
$\Lambda^h(\rhob_j)$ twisted by each $\det(\rhob_i)$ with $i<j$ (twisted by
nothing if $j=1$). With the obvious definition for the $\rho_{(a_i)}$ when $(a_i)$ is
any sequence as in (MF'), we have a decomposition 
$$\Lambda^k T \otimes A/I = \sum_{(a_i)} \tr(\rho_{(a_i)}),$$
hence $I$ contains the total reducibility ideal of $\Lambda^k T$
($\Lambda^k T$ is multiplicity free by (MF')). Theorem \ref{Drhojlibre} implies then that 
\begin{eqnarray} \label{Dck1} 
\Dc(\rho_{(a_i(k))}(\kappa))^{\varphi=F} 
\end{eqnarray} is 
free of rank one over $A/I$. \par
By induction, we know that $\rho_i$ is a trianguline deformation of
$(\rhob_i,\Ref_{z,i})$ whose parameter is $\delta'_i$ for each $i<j$. In
particular, for any such $i$,
$$\det(\rho_i)(\kappa_{x_i+1}+\cdots+\kappa_{x_i+d_i})$$ is a crystalline
character of $\Gp$ whose Frobenius eigenvalue is $F_{x_i+1}\cdots F_{x_i+d_i}$.
As 
$$\rho_{(a_i(k))}(\kappa)=\Lambda^h(\rho_j)(\kappa_{x_j+1}+\cdots+\kappa_{x_j+h})
\bigotimes_{i<j} \det(\rho_i)(\kappa_{x_i+1}+\cdots+\kappa_{x_i+d_i}),$$  
we get from formula (\ref{Dck1}) that
$$\Dc((\Lambda^h\rho_j)(\kappa_{x_j+1}+\dots+\kappa_{x_j+h}))^{\varphi=F_{x_j+1}\cdots
F_{x_j+h}}$$
is free of rank $1$ over $A/I$ for $h=1,\cdots,d_j$, which is the assertion
(\ref{trucaprouver}) that we had to prove.
\end{pf}

Note that the theorem implies that $\kappa_n-\kappa_{m}$ is constant on 
the total reducibility locus whenever $n$ and $m$ are in the same 
$\sigma$-orbit.

\begin{cor} \label{corredloc} Assume (INT) and that the permutation $\sigma$ is
transitive.
\begin{itemize}
\item[(i)] Every difference of weights $\kappa_n - \kappa_m$ is constant on
the total reducibility locus.
\item[(ii)] Assume moreover that one weight $\kappa_m \in A/I$ is constant,
and that for some $i$ we have $\Hom_{\Gp}(\rhob_i,\rhob_i(-1))=0$. Then $\rho_i$ is
crystalline. 
\end{itemize}
\end{cor}

\begin{pf} The assertion (i) follows immediately from the second assertion of
Theorem \ref{tridefgen}. \par
As a consequence, if $\kappa_m$ is constant for some $m$, 
every $\kappa_n$ is constant on the total reducibility locus. 
By Theorem \ref{tridefgen} and Proposition \ref{senpol}, this means that each $\rho_j$, seen as a 
representation $\rho_j : G \longrightarrow \Gl_{d_j}(A/I)$, $I_{\PP}\subset
I\subsetneq A$, is Hodge-Tate. On the other hand, each $\rho_j$ is a
trianguline deformation of the non critically refined representation
$(\rhob_j,\Ref_{z,j})$ again by Theorem \ref{tridefgen}, 
hence $\rho_i$ is crystalline by Proposition~\ref{colcrit}. 
\end{pf}

It turns out that the "non-trianguline" part of Theorem
\ref{tridefgen}, namely that the $\kappa_n-\kappa_{\sigma(n)}$ are constant
on the total reducibility locus, can be also proved even if we do not assume
(INT), but instead the different kind of assumption:
\begin{itemize}
\item[(HT')] For each $k \in \{1,\cdots,d\}$, $\Lambda^k \rhob_z$ has distinct Hodge-Tate weights.
\end{itemize}

\begin{theorem} Assume (HT') (or (INT)).
\label{thmredloc} Then for all $n=1,\dots,d$, 
$$(\kappa_{\sigma(n)}-\kappa_n)-(\kappa_{\sigma(n)}(z)-\kappa_n(z)) \in 
I_\PP$$ 
In other words, $\kappa_{\sigma(n)}-\kappa_n$ is 
constant on the total reducibility locus.
\end{theorem} 

Of course, part (i) of corollary \ref{corredloc} also holds assuming (HT')
instead of (INT).

\begin{pf} 
It is obviously sufficient to prove that for all $k$ in 
$\{1,\dots,r\}$, 
we have 
\begin{eqnarray} \label{sumik}
\sum_{i=1}^k  \left(\kappa_{\sigma(i)} - \kappa_i - (\kappa_{\sigma(i)}(z) 
- \kappa_i(z))\right)  \in I_\PP. \end{eqnarray}

We consider the family $(X,\Lambda^k T)$. As seen in~\S\ref{weakly}, 
this family is naturally weakly refined, with same set $Z$,  
\begin{eqnarray} F=\prod_{n=1}^k F_n, \end{eqnarray}
and first weight 
\begin{eqnarray} \label{kappaform} \kappa = \sum_{n=1}^k \kappa_n.\end{eqnarray} 
As already explained in the proof of Theorem \ref{tridefirr} this family 
satisfies all the hypotheses 
of \S\ref{hypkisinfamily}, and we want to apply to 
it Theorem~\ref{redlocfaible}.

For this, we note first that the (unique by (REG)) 
irreducible subrepresentation of 
$\Lambda^k \rhob_z$ that has the $\varphi$-eigenvalue 
$p^{\kappa(z)} F(z)$ in its $\Dc$ is the one denoted $\rhob_{(a_i)}$ above,
with $a_i$ being, for $i=1,\dots,r$, the numbers of $n \leq k$ such that 
$p^{\kappa_n(z)} F_n(z)$ is an eigenvalue of $\Dc(\rhob_i)$. In other words, $a_i$ is the  
numbers of $n \in \{1,\dots,k\}$ such that $n \in R_i$, that is $a_i=|R_i \cap \{1,\dots,k\}|$.

It follows from (NCR) and Lemma~\ref{tensornoncrit} that 
$\Dc(\rhob_{(a_i)}(\kappa(z)))^{\varphi=F(z)}$ has weight
 $\kappa'(z)-\kappa(z)$, where $\kappa'(z)$ is the smallest weight of 
$\rhob_{(a_i)}$. Hence $\kappa'(z)$ is the sum, for $n=1,\dots,k$ 
of the sum of the $a_n$ smallest weights of $\rhob_n$. In other words, 
\begin{eqnarray} \label{kappaprime} \kappa'(z)=\sum_{n=1}^k \kappa_{\sigma(n)}(z) .\end{eqnarray}

We now are in position to apply Theorem~\ref{redlocfaible}, which tells us that 
$$\kappa'-\kappa - (\kappa'(z) - \kappa(z)) \in I$$ 
where $I$ is the total reducibility ideal for the pseudocharacter $\Lambda^k T$. 
But it follows immediately from the definition of reducibility ideals and from hypothesis (MF')
that $I\subset I_\PP$, the total irreducibility ideal of $T$. 
So  $$\kappa'-\kappa - (\kappa'(z) - \kappa(z)) \in I_\PP,$$
which, using (\ref{kappaform}) and (\ref{kappaprime}) is the formula (\ref{sumik}) we wanted to prove.
\end{pf}

\subsection{Results on other reducibility loci}
 
It would be nice, and certainly useful, to have a result analogous to 
Theorem~\ref{thmredloc} for arbitrary reducibility ideals $I_\PP$, not only the total
reducibility ideal.  
This result should probably be  that certain differences of
weights $\kappa_i-\kappa_j$, for suitable couples $(i,j)$ combinatorically
defined in terms of the permutation $\sigma$ and the partition $\PP$, 
should be  constant of the reducibility locus attached to $\PP$.

But when we try to apply the methods used above, 
we get into trouble because there does not exist in general a module 
$M_{I}$ for $I$ a subset of 
$\{1,\dots,r\}$, analogous to the module $M_j$ for $j \in \{1,\dots,r\}$, in the sense that for $J$ a cofinite-length ideal of $A$, the isotypic
component of the $\rhob_j$, $j \in I$ in $M_I/JM_I$ is free over $A/J$.
This lack of freeness prevents to apply the ''constant weight lemma'' to this
module, and more generally any of our main results of section 2. This may be
a strong motivation to extend the results of section 2 to the non-free case,
but this does not seem to be easy, and we postpone this question to subsequent works
(of us or others).

However, we can still get an interesting although much coarser result on 
arbitrary reducibility loci by the method of our Theorem 9.1 in \cite{BC}. We shall 
give a sufficient condition for the other (non-trivial) reducibility ideals at a point $z$
 to be torsion free. This is 
equivalent to saying that the pseudocharacter $T$ is generically irreducible over every irreducible
component of $X$ through $z$.

As our result is coarse, we do not need for it our hypotheses of
\S\ref{hypredloc}, so we release (NCR), and (MF'), 
and we only assume  below that $z$ is a point of $Z$ that satisfies (REG).
 In that context the definitions of the subsets $W_i$ and $R_i$ (for $i=1,\dots,r$) of $\{1,\dots,d\}$ 
in \S~\ref{permutation} still make sense. For every $P \subset \{1,\dots,r\}$ we define the subset 
$W_P := \coprod_{i \in P} W_i$ and $R_P := \coprod_{i \in P} R_i$.  
 
\begin{theorem} \label{redgen} Let $\PP=\{P,Q\}$ be a non-trivial partition of 
$\{1,\dots,r\}$. Assume that $W_P \neq R_P$. Then 
$I_\PP$ is a torsion-free ideal of $A$.
\end{theorem}
\begin{remark}\label{remtransirrloc}
In particular, if the permutation $\sigma$ of \S\ref{permutation} 
is transitive, then the hypothesis of this theorem holds for all $P$ since $\sigma(R_P)=W_P$. 
In this case, the conclusion may be rephrased as : 
$T$ is generically irreducible on each irreducible component of $X$ through $z$.

When $\rhob_z$ is ordinary, the
hypothesis of the theorem, for all $P$, is equivalent to the transitivity of $\sigma$. In general, the transitivity is a stronger assumption.
\end{remark}
\ps

\begin{pf} Let $K=\prod K_s$ be the total fraction ring of $A$. We have to prove that $I_\PP K =K$, that is 
that for all $s$, $I_\PP K_s=K_s$. Replacing $X$ by its normalization
$\widetilde{X}$, $A$ by its integral closure in $K_s$, the $F_i$ and $\kappa_i$ by their
composition with $\tilde{X} \longrightarrow X$, and $Z$ by its inverse image
in $\widetilde{X}$, we may assume that $A$ is a domain, that $X$ is normal irreductible,
and what we have to prove is now that $I_\PP \neq 0$. \par

Assume by contradiction that $I_\PP=0$. Then there are two $A$-valued pseudocharacters $T_P$ and $T_Q$ 
such that $$T=T_P + T_Q, \, \, \, \, {\rm and}\, \, \,  T_\ast \otimes k =
\sum_{i \in \ast} \tr
\rhob_i.$$
Reducing $X$, we may assume that $X$ is an affinoid neighbourghood of $z$ (note that $z \in Z$), 
that $T_P$ and $T_Q$ take values in $\anneau(X)$, that for $i\neq j$ the $\kappa_i-\kappa_j$ are invertible on $X$ (since so they are at
$z$), and that $T_P$ is the generic trace of a representation of $G$ on a
finite type torsion free $\OO(X)$-module\footnote{As $T_P$ is residually
multiplicity free, the existence of such a module follows for example from
Lemma \ref{prolon}.}, say $M(X)$. By the maximum principle, the $v(F_n)$, 
$n=1,\dots,d$ are bounded on $X$. Hence Prop.~\ref{propsubfamilies} below implies that there is a 
set $I \subset \{1,\dots,d\}$ and $Z_1 \subset Z$ such that $T_P$ is refined by the $\kappa_n$, the $F_n$
for $n \in I$ and $Z_1$. 

We now claim that the eigenvalues of the crystalline Frobenius on 
$$(\rhob_P)_z:=\oplus_{i \in P} \rhob_i$$ are the $p^{\kappa_n(z)}F_n(z)$ for $n \in I$ (in other words,
we claim that we could assume that $z \in Z_1$). Indeed, by Kisin's theorem applied to the torsion free quotient of $\Lambda^k M(X)$ (apply
Theorem~\ref{thmkisin} to a flatification of the latter module as in the
proof of Theorem \ref{kisinnonlibre}), $1 \leq k \leq |I|=\dim T_P$ and to the maximal ideal of $A$, we get 
denoting by $I_k$ the first $k$ elements in $I$, 
$$D_\cris(\Lambda^k (\rhob_P)_z)^{\varphi=\prod_{n \in I_k} p^{\kappa_{n}}(x) F_n(x)}\neq 0.$$
The claim follows from this and (REG). 

By definition, we thus have $R_P=I$. Similarly, since the weights of 
$\rhob_z$ are the $\kappa_n(z)$, $n \in I$, we have $W_P=I$. But this implies that $W_P=R_P$, 
a contradiction.
\end{pf}

\begin{prop} \label{propsubfamilies}
Let $(X,T)$ be a refined family as above. We assume that  $X$ 
is connected, that the $\kappa_i -\kappa_j \in \OO(X)^*$ for all $i\neq j$, and that the $v(F_n)$, $n =1,\dots,d$ are bounded on 
$X$. If $T=T_1 + T_2$ where $T_i$, $i=1,2$ are pseudocharacters $G \rightarrow \anneau(X)$. Then there is a subset $I$ of 
$\{1,\dots,d\}$ and a subset $Z_1$ of $Z$ such that $(X,T_1)$ is refined by\footnote{The implicit ordering on $I$ here 
is the natural induced by $\{1,\cdots,d\}$.} 
$((\kappa_n)_{n \in I}, (F_n)_{n\in I},Z_1)$.
\end{prop}
\begin{remark} As we saw in the proof of Theorem~\ref{redgen}) 
we can actually enlarge $Z_1$ to contain all the points of $Z$ that satisfy (REG).
\end{remark}
\begin{pf} We denote by $(\rhob_1)_x$ (resp. $(\rhob_2)_x$) the semi-simple representation of trace 
the evaluation of $T_1$ (resp. $T_2$) at $x$, so that we obviously have
\begin{eqnarray} \label{rho12} \rhob_x \simeq (\rhob_1)_x \oplus (\rhob_2)_x. \end{eqnarray}

We first prove that there is an 
$I \subset \{1,\dots,d\}$, with $|I|=\dim T_1$, such that for all 
$x \in X$, the Hodge-Tate-Sen weights of $(\rhob_1)_x$ are the 
$\kappa_n(x)$, $n \in I$. For this we will only use property (i) of Definition \ref{defrefi} of a refined family.
Since $X$ is connected, and the weights everywhere distinct, it is obviously sufficient
to prove it when $X$ is replaced by any connected open subset $U$ of an admissible covering of $X$. 
So we may assume that $X$ is an 
affinoid. Since $\anneau(X)$ is noetherian, and by replacing $X$ by a finite surjective covering if necessary, 
we may assume that there exists (see \cite[lemme 7.1]{BC}) a finite type torsion-free module $\MM_1$ (resp. $\MM_2$) 
on $\anneau(X)$ with a continuous Galois action whose trace (defined after tensorizing by the fraction field of $\anneau(X)$) is $T_1$ 
(resp. $T_2$). Replacing $X$ by a blow-up $X'$ as in Lemma~\ref{grra}, we may also assume that 
$\MM_1$, $\MM_2$ are locally free, and by localizing again, that $X$ is a connected affinoid and that 
$\MM_1(X)$, $\MM_2(X)$ are free modules. The Sen polynomial of the module $\MM_1 \oplus \MM_2$ is 
$\prod_{n=1}^d (T-\kappa_n)$. Since it is split and $X$ is connected, it is easy to see that the Sen polynomial of $\MM_1$ 
has the form $\prod_{n \in I} (T-\kappa_n)$, for some subset $I$ of $\{1,\dots,d\}$. This proves the first assertion.

Now choose an integer $C$ greater than $\sum_{n=1}^d \sup_X v(F_n)$.
Let $z \in Z_C$. By~(\ref{rho12}), there is a subset $J$ of $\{1,\dots d\}$, with $|J|=\dim T_1$,
 such that the Frobenius eigenvalues of $\rhob_1(z)$ are $p^{\kappa_n(x)} F_n(x)$, $n \in J$.
By admissibility of $\Dc((\rhob_1)_z)$, we have $$\sum_{n \in I} \kappa_n(z) = 
\sum_{n \in J} (v(F_n(z)) + \kappa_n(z)),$$ that is
$$\kappa_I(z)-\kappa_J(z) = \sum_{n \in J} v(F_n(z)),$$
which implies $$|\kappa_I(z)-\kappa_J(z)| \leq \sum_{n=1}^d |v(F_n(z))| = C,$$ so by definition of $Z_C$,  we have $J=I.$
Thus it is clear that $((\kappa_n)_{n \in I}, (F_n)_{n\in I},Z_C)$ is a refinement of $(X,T_1)$.
\end{pf}

\newpage 

\newpage 
\section{Selmer groups and a conjecture of Bloch-Kato}\label{Selmer}

We recall in this section the Galois cohomological version of the standard conjectures 
on the order of vanishing of arithmetic $L$-functions at integers.
 Our main references are \cite{BloKa} and \cite{FP}.

\subsection{A conjecture of Bloch-Kato}\label{geometric} 
\subsubsection{Geometric representations} Let $E$ be a number field\index{E@$E$, a number field|(}, $p$ a prime\index{p@$p$, a fixed prime} and $F$ a finite extension of $\Q_p$. Let $$\rho: G_E \longrightarrow \GL_n(F)$$
be a continuous representation of the absolute Galois group $G_E$\index{GE@$G_E$, the absolute Galois group of $E$|(} \index{Zrho@$\rho$, a continuous geometric representation of $G_E$|(} of $E$, 
which is {\it geometric} in the sense of Fontaine and Mazur 
(see \cite{FP}), that is:\ps

\begin{itemize}

\item[--] $\rho$ is unramified outside a finite number of places of $E$,\ps

\item[--] $\rho_{|G_{E_v}}$ is De Rham for each place $v$ dividing $p$.\ps
\index{v@$v$, a place of $E$ above $p$}
\end{itemize}

\medskip 

\noindent It is known that the natural Galois representation on the \'etale cohomology groups $$H^{i}_{et}(X_{\overline{E}},\Z_p)\otimes_{\Z_p}F(d),$$ where $X$ is proper smooth over $E$ and $d \in \Z$, is geometric. The Fontaine-Mazur conjecture is the statement that every irreducible geometric continuous $G_E$-representation $\rho$ is a subrepresentation of 
such a representation on an \'etale cohomology group. 
 
\subsubsection{Selmer groups} \label{defselmer} We now define the Selmer group $H^1_f(E,\rho)$ of a geometric representation $\rho$.\index{H1f@$H^1_f(E,\rho)$, the Selmer group of $\rho$} 
This is the $F$-subvector space of the continuous Galois cohomology group\footnote{For the basic properties of continuous cohomology in this context, see e.g. \cite[App. B]{rubinlivre}.} $H^1(G_E,\rho)$ that parametrizes the isomorphism classes of continuous extensions

\begin{eqnarray} \label{U}
0 \longrightarrow \rho \longrightarrow U \longrightarrow F \longrightarrow 0,\end{eqnarray}
where $F$ denotes the trivial $F[G_E]$-module, satisfying for each finite place $v$ of $E$: \ps

\begin{itemize}\ps

\item[i)] $\dim U^{I_v}=1+\dim \rho^{I_v}$ if $v$ does not divide $p$,\ps

\item[ii)] $\dim D_{\cris}(U_{|G_v})=1+\dim D_{\cris}(\rho_{|G_v})$.\ps

\end{itemize}

\smallskip 

For example, such an $U$ is unramified (resp. crystalline) at a place $v$ whenever $\rho$ is. Moreover, at places $v$ dividing $p$, condition ii) implies
$\dim D_{\text{dR}}(U_{|G_v})=1 + \dim D_{\text{dR}}(\rho_{|G_v})$ so $U$ is De Rham 
since $\rho$ is. In particular, $U$ is geometric. As a consequence 
(see e.g.~\cite[Prop. B.2.7]{rubinlivre}), $H^1_f(E,\rho)$ is a finite dimensional $F$-vector space. 

\begin{remark} \label{remdefselmer}
\begin{itemize}
\item[i)] The formation of $H^1_f(E,\rho)$ commutes with any finite extension 
of the field $F$ of coefficients of $\rho$.
\item[ii)] Note that the functors 
$V \mapsto V^{I_v}$ and $V \mapsto D_{\cris}(V)$ (on the category of continuous $F[G_{E_v}]$-modules) being left exact, 
both conditions i) and ii) may be viewed as the requirement that they transform the short
exact sequence~(\ref{U}) of $F[G_E]$-modules into a short {\it exact} sequence of vector spaces.
\item[iii)] By Grothendieck's $l$-adic monodromy theorem, condition $i)$ is 
automatic if $(U^{I_v})^{\ses}$ does not contain the cyclotomic character.
\end{itemize}
\end{remark}
\smallskip

\begin{example} \label{exampleBK}\begin{itemize}
\item[i)] Assume that $\rho=\Q_p(1)$ is the cyclotomic character. 
Kummer theory (or Hilbert 90) shows that there is a canonical 
isomorphism 
$$E^* \wt_{\Z} \Q_p \overset{\sim}{\longrightarrow}
H^1(E,\Q_p(1)).$$ Under this identification, it is well known
that\footnote{First show the local analogue with $E$ replaced by any
  $E_v$, and conclude using the finiteness of the class number of
  $E$.} $\OO_E^* \wt_{\Z} \Q_p \isomo H^1_f(E,\Q_p(1))$. If we relax the hypothesis $f$ at a finite set $S$ of places of $E$, we get $S$-units instead of units of $E$. 
\item[ii)] Assume that $A$ is an abelian variety over $E$ and take
  $\rho=T_p(A)\otimes \Q_p$. Then it known that the 
$f$ condition at a place $v$  cuts out precisely the elements of the
$H^1(G_E,T_p(A))$ coming locally at $v$ from an $E_v$-rational
point of $A$ (when $v$ divides $p$, see \cite{BloKa}). The Kummer  sequence becomes then:
			$$ 0 \longrightarrow A(E)\otimes_{\Z} \Q_p \longrightarrow H^1_f(E,T_p(A)) \longrightarrow {\rm Sha}_p(A) \otimes \Q_p \longrightarrow 0,$$
\noindent where ${\rm Sha}_p(A)$ is the dual of the Tate-Shafarevich group of $A$. Assuming the finiteness of the Tate-Shafarevich group, $H^1_f(E,T_p(A))$ appears to be a purely Galois theoretic definition for $A(E)\otimes_{\Z} \Q_p$. \end{itemize}
\end{example}

\subsubsection{The general conjecture} Let $\rho$ be as in \S\ref{geometric} and fix embeddings $\Qb \longrightarrow \Qpb$ and $\Qb \longrightarrow \C$. 

It is expected that the Artin $L$-function $L(\rho,s)$ attached to
$\rho$ and these embeddings converges on a right half plane and admits a meromorphic (even entire when $\rho$ is not a Tate twist of the trivial character) continuation to the whole of $\C$. This is known for example when $\rho$ corresponds to a cuspidal automorphic representation of $\GL_n(\AAA_E)$. The general conjecture is then the following.

\begin{conj}\label{conjectureBK} ${\rm ord}_{s=0}L(\rho,s)=\dim_F H^1_f(E,\rho^*(1))-\dim_F (\rho^*(1))^{G_E}$ \end{conj}

Note that this is a conjectural equality between two terms, the one on the left being only conjecturally defined in general ! There are more precise conjectures predicting the leading coefficient of $L(\rho,s)$ at $0$, 
but we won't care about them in this paper. In view of Examples \ref{exampleBK}, the above conjecture generalizes Dirichlet units theorem (together with his theorem on the
 finiteness of the class number) and (assuming the finiteness of the Tate-Shafarevich group) the Birch and Swinnerton-Dyer conjecture. \ps

When $\rho$ is a cylotomic twist of finite image representation, the
conjecture is a theorem of Soul\'e \cite{soule}. 
Moreover, in the case $n=1$ and $E$ totally real or imaginary quadratic, the conjecture 
follows from Iwasawa's main conjecture for those
fields, 
proved by Wiles and Rubin respectively. 
Aside some sporadic results concerning the sign conjecture (see below), only a few cases are known when $n=2$ and $E=\Q$, and then
the terms in the equality are $0$ or $1$ (Wiles, Rubin, Gross-Zagier, Kato). Needless to say, each of those particular cases is a very deep theorem. 

\begin{remark} \label{noncenter} Assume that $\rho$ is pure of motivic weight $w$.
Apart from the case where $w=-1$, the conjectural left hand side of the equality in Conjecture \ref{conjectureBK} can be defined explicitly without any mention of $L$-function. \begin{itemize}
\item[i)] ($w\leq -2)$ Indeed, if $w < -2$, then $0 > 1+w/2$ should be in the domain of convergence of the Euler product defining $L(\rho,s)$ by Weil's conjectures, thus $\ord_{s=0} L(\rho,s)$ should be $0$ (and so should be $H^1_f(E,\rho)$). If $w=-2$ then $0$ should be on the boundary of the domain of convergence, and a conjecture of Tate predicts that in this case $\ord_{s=0}L(\rho,s)$ should be $- \dim (\rho^*(1))^{G_E}$.
\item[ii)] ($w\geq 0$) Recall that we expect a functional equation 
 \begin{eqnarray} \label{fonceq1}
 \Lambda(\rho,s)=\varepsilon(\rho,s)\Lambda(\rho^*(1),-s)\end{eqnarray}
  where $\Lambda(\rho,s)$ is the {\it completed} $L$-function, a product of $L(\rho,s)$ by a finite number of 
some simple explicit {\it $\Gamma$-factors} (see \cite{tate} for the recipe). Since $\rho^*(1)$ has weight $-w-2$, and by i) above, the term $\ord_{s=0} L(\rho,s)$
is determined when $w\geq 0$ by the order of the poles of the $\Gamma$-factors. 
\end{itemize}
However, although we can predict the integer of Conjecture
\ref{conjectureBK} when $w \neq -1$, it is still completely
conjectural that $\dim_F H^1_f(E,\rho)$ is actually this number. When
$w=-1$, e.g. as in the Birch and Swinnerton-Dyer conjecture, the
situation is even much worse (and more interesting) as the integer in
question is completely mysterious, and might be any integer in principle.
\end{remark}

\subsubsection{The sign conjecture} 

Among the cases where the motivic weight of $\rho$ is $-1$, of special 
interest are the ones
where $0$ is the "center" of the functional equation of $\rho$, that is when the equation~(\ref{fonceq1}) takes the form :
\begin{eqnarray} \label{fonceq2}
 \Lambda(\rho,s)=\varepsilon(\rho,s)\Lambda(\rho,-s),\end{eqnarray}
In this case we have $\epsilon(\rho,0)=\pm 1$, and this number is called the {\it sign} of the functional equation of $\rho$ (or 
shortly the sign of $\rho$). As the $\Gamma$-factors do not vanish on the
real axis, Conjecture~\ref{conjectureBK}
leads to an important special case, that we will call the {\it sign conjecture} :
\begin{conj} \label{sign} Assume $\rho$ satisfies~(\ref{fonceq2}). 
If $\varepsilon(\rho,0)=-1$ then $H^1_f(E,\rho^*(1)) \neq {0}$. \end{conj}

\begin{remark}
\begin{itemize}
\item[(i)]  The sign conjecture for $E=\Q$ implies the sign conjecture
  for any $E$. For if $\rho$ is a geometric irreducible representation
  of $G_E$ whose  functional equation satisfies (\ref{fonceq2}) with sign $-1$,
 $\tau=\Ind_{G_\Q}^{G_E} \rho$ is a semi-simple representation of $G_\Q$ with same
 sign, isomorphic Selmer group, and satisfies 
 $\tau \simeq \tau^\ast(1)$ by lemma \ref{eq1eq2}. It follows that $\tau$ is a direct sum of a  
 subrepresentation $\tau_0 \oplus \tau_0^\ast(1)$ (whose sign is 1) and  
of irreducible subrepresentations $\tau_1,\dots,\tau_l$ such that $\tau_i \simeq \tau_i^\ast(1)$ for $i=1,\dots,l$. Since the product of the signs of the factors of a direct sum is the sign of that direct sum, if $\rho$ has sign $-1$
there must be an $i$ such that $\tau_i$ has sign $-1$. 
Thus the sign conjecture for $\Q$ asserts the existence of a non zero element in $H^1_f(\Q,\tau_i)$, hence in $H^1_f(\Q,\tau)=H^1_f(E,\rho)$.
\item[(ii)] Even if the analytic continuation at $0$ of $L(\rho,s)$ is
  not known, it is possible to give a non conjectural meaning to the
  sign $\epsilon(\rho,0)$ (which is a product of local terms), hence to the sign conjecture (see \cite[\S3]{gross}).
\end{itemize}
\end{remark}

As an exercise, let us determine when equation~(\ref{fonceq2}) holds. 
We need a notation : for
 $\sigma \in \text{Aut}(E)$, we denote by $\rho^\sigma$ 
the representation (well defined up to isomorphism)  $g \mapsto \rho(\tau g \tau^{-1})$ where  $\tau \in G_\Q$ is an element inducing $\sigma$ on $E$.

\begin{lemma} \label{eq1eq2} 
We assume~(\ref{fonceq1}). Then equation~(\ref{fonceq2}) holds if there exists a 
 $\sigma \in \text{Aut}(E)$ such that  $\rho^\ast(1) \simeq
 \rho^\sigma$. When $E$ is Galois (resp. $E=\Q$) and $\rho$ is irreducible (resp. semisimple), the 
converse holds. 
\end{lemma} 

\begin{pf} In view of equation~(\ref{fonceq1}), equation~(\ref{fonceq2}) holds 
if and only if $\rho$ and $\rho^\ast(1)$ have equal
$\Lambda$-functions. As any $\sigma \in \Aut(E)$ induces a norm-preserving bijection on primes ideal of $E$, it is clear 
that $\Lambda(\rho,s)=\Lambda(\rho^\sigma,s)$ and the first assertion follows.

For the converse, it is enough to show that when $E$ is Galois, if two
irreducible, continuous and almost everywhere unramified, 
representations $\rho$ and $\rho'$ of $G_E$ have the same $L$-function, 
there exists a $\sigma \in \Gal(E/\Q)$ such that $\rho \simeq
\rho'^\sigma$. When $E=\Q$ and $\rho$ and $\rho'$ are more generally
semisimple, 
that is true because they have equal characteristic
 polynomials of Frobenii for almost all $p$, hence $\rho \simeq \rho'$ by Cebotarev's theorem. 
 Now for $E$ any number field, if $\rho$ and $\rho'$ are semi-simple representations 
of $G_E$ having the same $L$-functions  then this still holds for  $\Ind_{G_\Q}^{G_E} \rho$ 
and $\Ind_{G_E}^{G_\Q} \rho'$ which hence are isomorphic. 
Taking the restrictions to $G_E$, we find that if $E$ is Galois,
we have 
$$\oplus_{\sigma \in \Gal(E/\Q)} \rho^\sigma \simeq \oplus_{\sigma \in
  \Gal(E/\Q)} \rho.$$
Hence, if $\rho$ is irreducible,
then it is isomorphic to a $\rho^\sigma$.
\end{pf}

\subsection{The quadratic imaginary case}

\subsubsection{Assumptions and notations}
\label{quad}
Throughout this paper, we will assume that $E$\index{E@$E$, a quadratic imaginary field|(} is an imaginary quadratic field, and 
we shall denote by $\sigma$ a complex conjugation in $\Gal(\bar E/\Q)$, and by $c$ its image in $\Gal(E/\Q)$, so that $\sigma^2=c^2=1$. For $U$ any representation 
of $G_E$, we set $U^\sigma(g)=U(\sigma g \sigma)$ and we denote by $U^\bot$ (pronounce ``U Bott'')\index{bot@$\bot$; $U^\bot$ is the contragredient of the conjugate (that is, the twist by the outer automorphism of $G_E$ induced by 
the conjugation of $E$) of a representation $U$ of $G_E$}
of the representation $$U^\bot:=(U^{\sigma})^\ast.$$
We shall fix a continuous geometric $n$-dimensional representation $\rho$ of $G_E$ over $F$, and we shall 
assume that \begin{eqnarray} \label{rhobot} \rho \simeq \rho^\bot(1). \end{eqnarray}
Hence $\rho$ should satisfy Equation~(\ref{fonceq2}) by
lemma~\ref{eq1eq2}. Note that in this case $$H^1_f(E,\rho^*(1))=H^1_f(E,\rho^\sigma)=H^1_f(E,\rho).$$
Our main objectives are, assuming some widely believed (and that might well be proved soon) 
conjectures in the theory of automorphic forms:
\begin{itemize}
\item[(1)] to prove the sign conjecture for such $\rho$ ;
\item[(2)] to give a lower bound of the Selmer groups $H^1_f(E,\rho)$ 
depending of the geometry of an explicit unitary eigenvariety, at an
explicit point.
\end{itemize}

\subsubsection{An important example} \label{examplerhof} 
Aside the case $n=1$, which is already of interest, 
an important class of examples is provided by base change to $E$ of classical modular forms. \ps
Let $k$ be an even integer, $N$ an integer prime to $p$, and $f$ a
normalized cuspidal newform for $\Gamma_1(N)$ with square nebentypus $\varepsilon^2$. 
If $F$ denotes the completion at a place dividing $p$ of the field of coefficients of $f$,
we shall denote by $\rho_f$ the representation $G_\Q \longrightarrow \Gl_2(F)$ attached to $f$ and normalized in such a way that $\rho_f^\ast(1)\simeq \rho_f$ (this uses the assumptions on $k$ and on the nebentypus). 
We note $\rho_{f,E}$ the restriction of $\rho_f$ to $G_E$, then obviously $\rho_{f,E}$ 
satisfies~(\ref{rhobot}). For suitable choices of $E$, the Selmer group of $\rho_{f,E}$ turns out not to be bigger than the Selmer group of $\rho_f$, as the following well known proposition shows.

\begin{prop} \label{propQE} Let $f$ be as above, and $S$ any finite set of primes. There is
an imaginary quadratic field $E$, split at every prime of $S$, such that $\rho_{f,E}$ is irreducible and
$$H^1_f(\Q,\rho_f)=H^1_f(E,\rho_{f,E}).$$
\end{prop}
\begin{pf} Indeed, we have $H^1_f(E,\rho_{f,E}) \simeq H^1_f(\Q,\rho_{f}) \oplus 
H^1_f(\Q,\rho_{f} \otimes \chi_E)$ where $\chi_E$ is the non trivial 
quadratic character of $G_\Q$ with kernel $G_E$. By the main result of \cite{hl}, generalizing
\cite{wald}, there is an infinite number of quadratic imaginary fields 
$E$ that split at every prime of $S$ and such that $L(\rho_f \otimes
\chi_E,0) \neq 0$. For such an $E$, \cite[Thm 14.2 (2)]{Kato} proves that 
$H^1_f(\Q,\rho_{f}\otimes \chi_E)=0$, hence the proposition.
\end{pf}

\subsubsection{Upper bounds on auxiliary Selmer groups}\label{auxselgroup}

In \cite{Bellaichethese} as well as in subsequent works using an automorphic method to
produce elements in Selmer groups (\cite{SU}, \cite{BC}, this paper), an important input is a result giving an {\it upper
  bound} (for instance, $0$) on the dimension of auxiliary 
Selmer groups. \ps

The most elementary case of such a result is the next proposition, which is of crucial importance in
both proofs of chapter 8 and 9. It would become 
false if $E$ were replaced by a CM field of degree greater than $2$, and it is actually the only point 
in the proof of the sign conjecture where the fact that $E$ is quadratic is really used. 
 \begin{prop} \label{propQ1} $H^1_f(E,\Q_p(1))=0$.
\end{prop}
\begin{pf}
By~\ref{exampleBK}(1), $H^1_f(E,\Q_p(1))$ is isomorphic to $\anneau_E^\ast \otimes_{\Z} \Q_p$, which is $0$ as $\anneau_E^*$ is finite.
\end{pf}

As said above, this vanishing result turns out to be the only one necessary to the proof of the sign conjecture. 
However, we will need quite a number of other vanishing results to get our second main result. The easiest ones are dealt with the following proposition.

\begin{prop} \label{propQ2} \begin{itemize}
\item[i)] $H^1_f(E,\Q_p)=0$,
\item[ii)] $H^1_f(E,\Q_p(-1))=0$. Moreover, for $\epsilon=\pm 1$,
the subspace of $H^1(E,\Q_p(-1))$ parametrizing extensions $U$ of 
$\Q_p(1)$ by $\Q_p$ such that $U^{\bot}(1) \simeq \epsilon U$ (as extensions)
 has dimension $\leq 1$.
\end{itemize}
\end{prop}
\begin{pf}
For an extension of $\Q_p$ by $\Q_p$, it is the same to be Hodge-Tate, crystalline, or unramified, so by class-field theory we have $$H^1_f(E,\Q_p)={\rm Cl}(\anneau_E)\otimes_\Z \Q_p=0,$$ which proves (i). It is well known that part (ii) follows from results of Soul\'e and from the invariance of Bloch-Kato conjecture under duality. For sake of completeness, we give an argument below. \ps
First, note that if $F/\Q_p$ is a finite extension, then $H^1_f(F,\Q_p(-1))=0$.
By Soul\'e's theorem \cite[thm. 1]{soule} $$H^2(\anneau_E[1/p],\Q_p(2))=0.$$ 
The version of Poitou-Tate exact sequence given in \cite[prop. 2.2.1]{FP} shows then that $H^1_f(E,\Q_p(-1))=0$ when applied to the Galois module $\Q_p(2)$. 
So we get an injection
$$H^1(E,\Q_p(-1)) \longrightarrow \bigoplus_{v \big\vert p}H^1(E_v,\Q_p(-1)),$$
which is compatible with the operation $U \mapsto U^\sigma$ on the domain
and the exchange of $v$ and $\bar v$ on the range when $p=v\bar{v}$ splits in $E$. We conclude as $H^1(E_v,\Q_p(-1))$ has dimension $1$ by Tate's theorem.
\end{pf}

The two other vanishing results we will need are somehow deeper. Since they 
are expected to be proved by completely different methods (e.g. using Euler systems) 
than those used in this paper, it would be artificial to limit ourselves to the case where those results have actually been proved, hence we take them as assumptions as follows. \ps
\medskip \index{BK@BK1($\rho$) and BK2($\rho$), hypotheses on $\rho$ that 
are consequences of the Bloch-Kato conjectures}
{\bf Hypothesis} BK1($\rho$) : $H^1_f(E,\rho(-1)) = 0$.
\ps \medskip
{\bf Hypothesis} BK2($\rho$) : Every deformation $\tilde \rho$ of 
$\rho$ over $F[\vareps]$ (the ring of dual numbers) that satisfies  
$\tilde \rho^\bot(1) \simeq \tilde \rho$ and whose corresponding cohomology class lies in\footnote{or, in an equivalent way, such that for each finite place $w$, $\tilde \rho_{|E_w}$ is geometric (automatic condition if $w$ is prime to $p$) with constant monodromy operator acting 
on $D_{\rm pst}(\tilde \rho_{|E_w})$.} $H^1_f(E,\ad \rho)$ is constant.
\ps
\bigskip
\begin{remark} \label{remhypBK} 
Hypothesis BK2($\rho$) can also be formulated in terms of Selmer groups : it says that the 
part invariant by the involution $V \mapsto V^\bot(1)$ on the group 
 $H^1_f(E,\ad \rho)$  is zero. Note that 
both hypotheses should follow from Conjecture~\ref{conjectureBK} for any 
$\rho$ that is pure of weight $-2$ (for BK2($\rho$)) or $-3$ (for BK1($\rho$)), and they are precisely in case i) of Remark \ref{noncenter}. Fortunately, those assumptions have already been proved in interesting cases. 

	This conjectural vanishing of $H^1_f(E,\ad \rho)$ is actually fundamental for understanding eigenvarieties. Empirically, 
it can be understood as follows. Let $$R: \G_E \longrightarrow \GL_n(L\langle t \rangle)$$ be 
a continuous morphism, and assume that for each $t \in \Z_p$ the evaluation
$R_t$ of $R$ at $t$ is a geometric irreducible representation. 
Then the Fontaine-Mazur conjectures implies that $R$ is conjecturally constant (up to isomorphism). Indeed, 
each $R_t$ is conjecturally cut out from an $E$-motive. 
But there is only a countable number of such motives, hence of $\tr(R_t)$, 
so $\tr(R)$ is constant. The assertion $H^1_f(E,\ad \rho)$ is actually a slightly 
stronger variant of that fact, in which we replace the Tate algebra $L\langle t \rangle$ by $L[t]/t^2$.
\end{remark}

\medskip

\begin{prop} \label{propBK1} BK1($\rho$) holds in the following two cases: \begin{itemize}
\item[i)] $n=1$ and $0$ is not a Hodge-Tate weight of $\rho$.
\item[ii)] $n=2$ and $\rho$ is of the form $\rho_{f,E}$ (using notations of \S~\ref{examplerhof}) for some eigenform $f$ of weight $k \geq 4$.
\end{itemize}
\end{prop}

\begin{pf} By the theory of CM forms, case i) follows from case ii), which in turn is a result of Kato 
\cite[Thm. 14.2 (1)]{Kato}.
\end{pf}

\begin{prop} \label{propBK2} BK2($\rho$) holds if $n=1$ or if $n=2$ and $\rho$ is of the form $\rho_{f,E}$ whenever $f$ is not CM and satisfies {\it one of} the following conditions :
\begin{itemize}
\item[(i)] At every prime $l$ dividing $N$, $f$ is either supercuspidal or Steinberg.
\item[(ii)] The semi-simplified reduction $\rhob_f$ of $\rho_f$ is absolutely reducible, and is (over
an algebraic closure) the sum of two characters that are distinct over $G_{\Q(\zeta_{p^\infty})}$.
\item[(iii)] For any quadratic extension $L/\Q$ with $L \subset \Q(\zeta_{p^3})$, $(\rhob_f)_{|G_L}$  is absolutely irreducible.
\end{itemize}
\end{prop}
\begin{pf} 
If $\rho$ is a character, then $\ad \rho$ is the trivial character
and BK2($\rho$) follows as in Proposition \ref{propQ2} i), hence we may assume that we are in the case $n=2$. \ps
By construction we have $\rho\simeq \rho^*(1)$. As for any $2$-dimensional representation over any ring, 
we have $\tilde \rho \simeq \tilde \rho^\ast 
\otimes \det \tilde \rho$. But by the case $n=1$ above, the character $\det \tilde \rho$ is constant. Thus
 we get $$\tilde \rho \simeq \tilde \rho^\ast(-1).$$ Together with the hypothesis $\tilde \rho 
\simeq \tilde \rho^\bot(-1)$, we get $\tilde \rho \simeq \tilde \rho^\sigma$. That is, there is an $A
\in \Gl_2(F[\vareps])$ such that for all $g \in G_E$, $A \tilde\rho(g) A^{-1}= \tilde \rho(\sigma g \sigma^{-1})$. Since $\sigma^2=\Id$, $A^2$ centralizes $\rho(G_E)$, so since $\rho_{f,E}$ is absolutely irreducible, $A^2=\lambda$ for some $\lambda \in F[\vareps]^*$. If $\bar A$ denotes the reduction  of $A$ modulo $\vareps$, we have for all $g \in G_E$,
$\bar A \rho(g) \bar A^{-1}= \rho_f(\sigma) \rho(g) \rho_f(\sigma^{-1})$ hence $\bar A^{-1} \rho_f(\sigma)$ centralizes $\rho(G_E)$. Thus we have $ \bar A=\mu \rho_f(\sigma)$ for some $\mu \in F^*$. In particular 
$\lambda = A^2 \equiv \mu^2 \pmod{\epsilon}$. Let $\tilde \mu$ be the square root of $\lambda$ in $F[\vareps]^*$ lifting $\mu$. Set $\tilde \rho_f(\sigma)=\tilde \mu^{-1} A$, and 
$\tilde \rho_f(g) = \rho_f(g)$ : this defines a deformation $$\tilde \rho_f : G_\Q \longrightarrow \Gl_2(F[\vareps])$$ of $\rho_f$  whose restriction to $G_E$ is $\tilde \rho$.\ps 
Since $(\tilde \rho_f)_{|\G_E} = \tilde \rho$ is geometric, so is $\tilde \rho_f$. But such a deformation of $\rho_f$ 
is trivial by \cite[Theorem, page 2]{kis2} in the cases (ii) and (iii), and by \cite[theorem 5.5]{weston}
in case (i), hence so is its restriction $\tilde \rho$.  
 \end{pf}

We shall actually use assumption BK1($\rho$) to bound the subspace $$H^1_{f'}(E,\rho(-1)) \subset H^1(E,\rho(-1))$$ parameterizing extensions which have good reduction at all primes not dividing $p$, that is satisfying only condition (i) of \S~\ref{defselmer}.

\begin{prop} \label{calulhunfrhodual} Assume that $p=vv'$ splits in $E$, and that $1, p^{-1}$ are not eigenvalues of the crystalline Frobenius on $\Dc(\rho_{|G_{E_v}})$. Then $$\dim_F H^1_{f'}(E,\rho(-1)) \leq n + \dim_F H^1_f(E,\rho(-1)).$$ In particular, if BK1($\rho$) holds we have  $\dim_F H^1_{f'}(E,\rho(-1))\leq n$.
\end{prop}

\begin{pf}
As $\rho(-1)^c \simeq (\rho(-1))^*(-1)$, we have an obvious exact sequence
$$ 0 \longrightarrow H^1_f(E,\rho(-1)) \longrightarrow H^1_{f'}(E,\rho(-1)) \longrightarrow H^1_f(E_v,W) \oplus H^1_f(E_v,W^*(-1))$$
with $W=\rho_{|G_{E_v}}(-1)$. 
For any de Rham $p$-adic representation $W$ of $\Gal(\Qb_p/\Q_p)$ such that  $\Dc(W)^{\varphi=1}=0$, Bloch-Kato's computation \cite{BloKa} gives us an isomorphism
%By the fundamental exact sequence $$0 \longrightarrow \Qp \longrightarrow \Bc^{\varphi=1} \longrightarrow \Bdr/\Bdr^+ \longrightarrow 0,$$
$$\Ddr(W)/\Fil^0(\Ddr(W)) \longrightarrow H^1_f(\Qp,W).$$ The proposition would then follow from the inequality $$\dim_F(\Ddr(W)/\Fil^0(\Ddr(W)))+ \dim_F(\Ddr(W^*(-1))/\Fil^0(\Ddr(W^*(-1)))) \leq n,$$
which is immediate.
\end{pf}

\newpage
\newcommand{\Norm}{{\rm Nm}}
\newcommand{\Spec}{\rm Spec}

\section{Automorphic forms on definite unitary groups: 
results and conjectures} 

\subsection{Introduction}  

This chapter recalls or proves all the results we shall need 
from the theory of representations of reductive 
groups and of automorphic forms. 

As explained in the general introduction, the main steps of 
our method regarding the proof of our two main theorems are, very roughly, as follows: 
starting with an $n$-dimension Galois representation $\rho$ such that 
$\varepsilon(\rho,0)=-1$,
we construct a very special, non tempered, 
automorphic representation $\pi^n$ for a unitary 
groups in $m=n+2$ variables. 
We deform it $p$-adically, in other words, we put it in an eigenvariety of the unitary group. 
We associate to this deformation of automorphic forms a deformation of Galois representations, or rather, 
a Galois pseudocharacter on the eigenvariety of the unitary group. 
This Galois pseudocharacter gives us the desired non 
trivial elements in the Selmer group of $\rho$.  

Unfortunately, some results needed to make work two of those steps in 
their natural generality have not yet been published or even written 
down: the first step, the existence of the ``very special'' 
automorphic representation $\pi^n$, has been announced, but a written proof is only available in
small dimension, namely $m \leq 3$;  
the third step relies on the existence and the basic properties of the Galois 
representations attached to (some) automorphic representations of unitary groups. 
Here again the desired results are only known for $m \leq 3$. Fortunately, 
this result is also in the process of being proved: it is one of the main goals of an 
ambitious project gathering many experts and participants of the GRFA seminar of the 
"Institut de math\'ematiques de Jussieu" in Paris, under the direction of Michael Harris. 
Their work should result in a four-volumes 
book (\cite{GRFAbook}) in the next few years that is expected to contain
a construction of the Galois representations attached to automorphic forms 
on unitary groups in many cases, and in particular in the cases we need. 
An important input in this project is the recent proof by Laumon and Ngo of the so-called fundamental lemma 
for unitary groups. 

In this chapter we formulate the two needed results as conjectures 
(namely conjecture REP($m$) on Galois representations attached to automorphic forms on 
unitary groups with $m$ variables, and conjecture AC($\pi$) 
on the existence of the automorphic representation $\pi^n$ constructed from 
the automorphic counterpart $\pi$ of $\rho$), and we shall
admit those conjectures in the proof of our main theorems in 
chapters 8 and 9. In view of the situation explained above, it would have been 
pointless to limit ourselves to the case where the needed results are 
already written down. 

The main reason for which we are able to write down some still unwriten
results and rely confidently of them is not that we are told they will be proved very soon, but because they are part of a much larger and very well 
corroborated set of conjectures called 
``the Langlands program'' (and its extension by Arthur).

We believe it will be of interest for the reader
to explain in great details how our conjectures (and much more) appear as 
consequences of the Langlands program, and in particular 
how the existence of our very special non tempered 
automorphic forms is enlighten as a special case of the 
beautiful ``multiplicity formulas'' of Arthur. This is the aim of the 
appendix to this book, that recalls the part of the
 Langlands and Arthur's program that we need, 
and where we show how our conjectures follow from their's. 
This appendix may be read independently, as an 
introduction to Langlands and Arthur's parameterizations and 
multiplicity formulas. Although logically independent of it, the rest of the 
chapter will make frequent references to this appendix 
for the sake of the reader's intuition.

Although we may  expect much more general results to be true, and even to be proved soon, that the 
modest conjectures  we state, we have make a 
great deal of effort, in this chapter and throughout this books, 
to keep our conjectural input to the theory of 
automorphic forms at the lowest possible level. 
One reason for doing this is obvious: the weaker the 
assumptions we have to admit, the stronger is our result, and the sooner 
it will become an unconjectural theorem. Another more serious reason is that part of our work (especially chapter 7) 
is also bound to be used in the book \cite{GRFAbook} for the construction or the proof of some 
properties of the Galois representations 
in some "limit" cases which are too complicated to handle via a direct comparaison of trace formulae. So the logical scheme would be as follows: 
in \cite{GRFAbook} should be proved "directly" for a quite "generic" set of automorphic representations the existence and 
properties of the associated Galois representations, which should be enough to check our conjecture 
Rep($m$). In turn, our work on eigenvarieties should complete the picture by providing existence and properties of 
the Galois representations attached to the remaining (cohomological) automorphic forms. For example, our conjecture Rep($m$) only needs the Galois representations 
for automorphic forms of regular weights. To give another example, our method (hence our conjecture Rep($m$))
makes no irreducibility hypothesis on the Galois representations,\footnote{Let us say that this is anyway a subtle point, 
as only the stable tempered automorphic representations should have irreducible associated Galois representations, and this 
property is very hard to detect in practice. This actually introduces an extra difficulty in the applications to 
the construction of nontrivial elements in the Selmer groups that we will explain how to circumvent. This feature was already 
present in \cite{BC}, but was absent of the earliest
stages of the method, like in \cite{Bellaichethese} or later in 
\cite{SU}.} but instead may be used to prove many cases of irreducibility.

Let us now explain more specifically the content of the chapter.

The subsection~\S\ref{unitarygroups} deals with some general facts 
about unitary groups, with an emphasize on the {\it definite} ones and their automorphic representations. 
We define explicitly the unitary groups $\U(m)$ we will work with. We need a group that is 
quasi-split at every finite place (otherwise, the representation $\pi^n$ can not be automorphic, as explained 
in the appendix - see Remark~\ref{nonquasisplitrelevant}), but that is also compact at infinity - so that 
we can apply the theory of eigenvarieties of \cite{Ch}.\footnote{Note that we may not use in this context the construction of $p$-adic families announced recently by Urban, since the ``virtual multiplicity'' of our $\pi^n$ might be zero.} This leads to the restriction that $m \not \equiv 2 \bmod 4$.

The subsections~\S\ref{localLanglands} to~\S\ref{repreal} are local preliminaries. The short subsection~\S\ref{localLanglands} recalls the local Langlands correspondence for $\Gl_m$, as characterized by Henniart 
and proved by Harris and Taylor. It will be used very frequently. The subsection~\S\ref{refinementslocal} develops the theory of refinements
(sometimes called {\it $p$-stabilizations}) of unramified representations of $\Gl_m(\Q_p)$ which a representation theoretic counterpart of the theory of refinements of crystalline Galois representations that we explained in chapter 2. We invite the reader to look at the introduction~\S\ref{refinementslocal} of that subsection for a more precise discussion on this concept. The subsection~\S\ref{repreal} recalls two descriptions of the continuous irreducible representations of the compact group $\U(m)(\R)$ and compare them.

Next come two other subsections of local preliminaries. They are both devoted to the crucial question of monodromy.\footnote{Let us say that monodromy is bound to play a crucial role 
in our final arguments. Indeed, it follows from the Arthur multiplicity formula
that under the hypothesis $\varepsilon(\rho,0)=1$ (not $-1$) 
there should exist an automorphic representation ${\pi'}^n$ for $\U(m)$, 
isomorphic to $\pi^n$ at 
every place except one, say $l$ with $l$ inert in the splitting 
quadratic field $E$ of $\U(m)$, and such that ${\pi'}^n_l$ has the same $L$-parameter as $\pi^n_l$ on $W_{\Q_l}$ but 
a greater monodromy. If it was possible to apply our method to $\pi'^n$, it would eventually 
lead to a construction of a non-trivial element in the Selmer group of $\rho$, element which should not exist 
when $L(\rho,0) \neq 0$ according to the conjecture of Bloch-Kato. This shows that a precise control of monodromy has to play a role in our 
argument. }
By ``monodromy'' of an admissible irreducible 
representation $\pi_l$ of $\U_m(\Q_l)$ 
we mean the conjectural notion encoded in the nilpotent element that 
appears as part of the conjectural 
morphism of the Weil-Deligne group of $\Q_l$ to ${}^L \U(m)$ attached to 
$\pi_l$. Concretely, what we need is threefold. 
We need to give a non-conjectural meaning to expressions such as 
``$\pi_l$ has no more monodromy that $\pi'_l$'' or 
``$\pi_l$ has no monodromy at all''. We need tools to 
be able to show in chapter 7 that some full irreducible components of the eigenvarieties of $\U(m)$ 
containing $\pi^n$ ``have no more monodromy at every place $l$ than
$\pi^n_l$ has''. Finally, we need to be able to 
translate this ``control on monodromy of $\pi_l$'' 
into a control on the action of the inertia subgroup
at $l$ on the Galois representations attached to $\pi$. 
The latter is a part of our conjecture Rep($m$).
The objective of~\S\ref{types} and~\S\ref{NMR} is to meet the two first needs. 

In \S\ref{types}, we deal with the monodromy of 
representations of $\U(m)(\Q_l)$ for $l$ split in $E$, 
that is for $\Gl_m(\Q_l)$. In this case the meaning of the monodromy is 
non conjectural and clear thanks to the local Langlands correspondence, 
that gives us for a representation $\pi_l$ a (class of conjugation of) 
nilpotent matrices $N(\pi_l)\in\Gl_m(\C)$ : 
we can simply say that $\pi_l$ has more monodromy than 
$\pi'_l$ if the adherence of the conjugacy class of $N(\pi_l)$ 
contains $N(\pi'_l)$. To be able to control the variation of the 
monodromy in a family of such $\pi_l$, we use then the existence of 
some particular $K$-types. As we will see, this will fullfill our second need since a general property of 
the eigenvarieties we will study is that the locus of points whose associated $\U(m)(\Q_l)$-representation contains a given 
$K$-type is a union of irreducible components (this actually holds for every $l$). Of course, the simplest example of such a 
$K$-type is the trivial representation of $\GL_m(\Z_l)$, 
which cuts out precisely as is well known the unramified constituent of the unramified principal series (that is, the non monodromic ones). 
For a general monodromy type, we use suitable $K$-types that have been constructed by Schneider and Zink (see \S\ref{types}). Note that the 
types constructed by Bushnell and Kutzko are {\it a priori} of no use for our purposes because they "do not see monodromy". However, let us stress 
that the construction of Schneider and Zink actually relies on those types.

In \S\ref{NMR}, we deal with representations of $\U(m)(\Q_l)$ for 
$l$ inert or ramified in $E$. The group $\U(m)(\Q_l)$ 
is a quasi-split group, but it is not split, and
the situation in this case is much less favorable. 
First we do not know the local Langlands correspondence 
for those groups, neither we know the base change to $\GL_m/E$ (from a conjectural point of view, see the final appendix). 
Hence there is no obvious way to define ``having less monodromy than'' or 
``having no monodromy at all'' for a representation $\pi$ of $\U(m)(\Q_l)$. 
Even worse, we were not able to come up with a plausible characterization, 
in terms of group theory, of those irreducible admissible representations 
of $\U(m)(\Q_l)$ that conjecturally have no monodromy\footnote{To convince the reader that this question is not easy, let us say that for $m=3$, there is a 
{\it supercuspidal} representation of $\U(3)(\Q_l)$, discovered by Rogawski and called
$\pi^s$, whose base change has a non trivial monodromy. See \S\ref{instructiveexample} in the final appendix.}. Second, there is no 
theory of types {\`a} la Bushnell-Kutzko for $\U(m)(\Q_l)$, $m>3$, 
not to speak of a theory \`a la Schneider-Zink. The first solution we imagined 
to solve those problems was to avoid them: that is, to assume all our 
automorphic representations to be unramified\footnote{Recall that the notion of
{\it unramified} representation 
makes sense for any quasi-split group : it  means having a non-zero fixed vector by a ``very special'' maximal compact subgroup, in the sense of Labesse.} at inert or ramified $l$. 
An unramified representation should certainly
be ``non-monodromic'', and unramifiedness is easy to control in deformation as explained in the $\GL_m$-case above. 
But the problem is: for odd $m$,
there is no representation of $\U(m)$ of the form $\pi^n$ that is unramified at ramified primes.\footnote{More precisely, for odd $m$
any discrete automorphic representation $\pi$ of $\U(m)$ whose $A$-packet lies in the image of the endoscopic transfert 
${}^L (\U(n) \times \U(2)) \longrightarrow {}^L \U(m)$ has the property that its base change $\pi_E$ to $\GL_m/E$ is ramified at 
each prime of $E$ ramified above $\Q$ (see the appendix \S\ref{descarthunit}). The reason is that for odd $m$ the aforementionned $L$-morphism 
contains in its definition a Hecke character $\mu$ of $E$ such that $\mu^{\bot}=\mu$ but which does not descend to $\U(1)$. Such a Hecke 
character is automatically ramified at the primes of $E$ ramified over $\Q$ (see \S\ref{introcarhecke}).} 
So this assumption is much too restrictive. Instead, we introduce a special class of principal series representations
of $\U(m)(\Q_l)$ that certainly should have no monodromy, and which will enable us to deal with a large number of $\rho$ 
also when $m$ is odd. We call those representations {\it Non Monodromic Strongly Regular Principal Series}. 
We show in \S\ref{NMR} that to be a NMSRPS is a constructible property in a family. 

After these local preliminaries, we turn to global questions. 
In subsection~\S\ref{subsectionrepm}, we state our assumption Rep$(m)$ on existence and simple
properties attached to (some) automorphic forms of $\U(m)$.
In subsection\S\ref{constructionpin} we construct place by place a 
representation $\pi^n$ of $\U(m)(\AAA_\Q)$ starting from a cuspidal automorphic representation $\pi$ of 
$\Gl_n(\AAA_E)$ satisfying some properties (recall that $m=n+2$). 
We then state as a conjecture AC$(\pi)$ 
(even if as we said earlier this has been announced) 
that this $\pi^n$, under the assumption 
that $\varepsilon(\pi,1/2)=-1$, is automorphic. 

\subsection{Definite unitary groups over $\Q$}

\label{unitarygroups}

\subsubsection{Unitary groups} 
Let $k$ be a field, $E/k$ an \'etale $k$-algebra of degree $2$ 
with non trivial $k$-automorphism $c$, and $\Delta$ a simple 
central $E$-algebra of rank $m^2$ equipped with a $k$-algebra 
anti-involution $x \mapsto x^*$ of the {\it second kind}, 
{\it i.e.} coinciding with $c$ on $E$. We can attach to this 
datum $(\Delta,*)$ a linear algebraic $k$-group $G$ whose points 
on a $k$-algebra $A$ are given by 
$$G(A):=\{ x \in (\Delta \otimes_{k} A)^*, xx^*=1\}.$$
The base change $G \times_k E$ is then isomorphic to the $E$-group $\Delta^*$ of invertible elements of $E$, hence $G$ is a twisted $k$-form of $\GL_m$. Actually, as is well known, every twisted $k$-form of $\GL_m$ is isomorphic to such a group. 

\begin{example}\label{exunitarygp} They are two essentially 
different cases.
\begin{itemize}
\item[i)] If $E \isomo k \times k$, then $\Delta \isomo \Delta_1 \times \Delta_2$ and $*: \Delta_1 \longrightarrow \Delta_2^{\opp}$ is an isomorphism. In this case, the choice of $i \in \{1, 2\}$ induces a $k$-group isomorphism $G \isomo \Delta_i^*$. In particular if $\Delta \isomo M_m(E)$, then the choice of $i$ determines a $k$-isomorphism $G \isomo \GL_m$ which is canonical up to {\it inner} automorphisms.
\item[ii)] If $E$ is a field, then we say that $G$ (and $G(k)$) is a {\it unitary group} attached to $E/k$. When moreover $\Delta = M_m(E)$, then $*$ is necessarily the adjunction with respect to a non degenerate $c$-Hermitian form $f$ on $E^m$, hence $G$ is the usual unitary group attached to this form. If $f$ is the standard anti-diagonal form
$$f(xe_i,ye_j)=c(x)y\delta_{j,m-i+1},$$ 
then $G$ is quasi-split, and will be referred in the sequel as the {\it $m$-variables quasi-split unitary group attached to $E/k$}.
\end{itemize}
\end{example}
\noindent 

\subsubsection{The definite unitary groups $\U(m)$}\label{definitionUm} Suppose from now that $k=\Q$, $E$ is a quadratic imaginary field, and assume that $\Delta=M_m(E)$ and $*$ is attached to some form $f$ on $E^m$ as in ii) above. Then $G$ is a unitary group over $\Q$. For each place $v$ of $\Q$, the local component $G \times_{\Q}\Q_v$ is then the $\Q_v$-group attached to the datum $(\Delta \otimes_{\Q} \Q_v,*)$, hence by Examples \ref{exunitarygp}: \begin{itemize}
\item[i)] If $p=xx'$ is a finite prime split in $E$, then $x: E \longrightarrow \Q_p$ induces an isomorphism $G(\Q_p) \isomo \GL_m(\Q_p)$,
\item[ii)] if $p$ is inert or ramified, then $G(\Q_p)$ is a unitary group attached to $E_p / \Q_p$,
\item[iii)] each embedding $E \longrightarrow \C$ gives an isomorphism between $G(\R)$ and the usual real unitary group $U(p,q)$, where $(p,q)$ is the 
signature of $f$ on $E^m\otimes_{\Q} \R$, $p+q=m$. 
\end{itemize}
\ps

We say that $G$ is {\it definite} if $G(\R)$ is compact, or which is the same if $pq=0$. We will be interested in definite unitary groups $G$ with some prescribed local properties. Their existence could be deduced from the Hasse's principle for unitary groups over number fields for which we refer to \cite[\S2]{clo}. For some computational reasons, we explicitly give them below. Let $N: E \longrightarrow \Q$, $x \mapsto x c(x)$ be the norm map, $m\geq 1$ an integer. 

\begin{definition} $\U(m)$ \index{Um@$\U(m)$, the $m$-variables 
unitary group attached to $E$ that is compact at infinity and 
quasi-split elsewhere, when $m \not \equiv 2 \pmod{4}$} 
is the $m$-variables unitary group attached to the positive definite $c$-hermitian form $q$ on $E^m$ defined by $$q((z_1,\cdots,z_m))=\sum_{i=1}^m N(z_i).$$
\end{definition}

\begin{prop}\label{propgpUm}
\begin{itemize}
\item[(i)] $\U(m)$ is a definite unitary group.
\item[(ii)] If $l$ does not split in $E$, and $m \not\equiv 2 \bmod 4$, then $\U(m)(\Q_l)$ is the quasi-split $m$-variables 
unitary group attached to $E_l/\Q_l$.
\end{itemize}
If $m \not \equiv 2 \bmod 4$, the group $\U(m)$ is the unique $m$-variables unitary group attached to $E/\Q$ that is quasi-split at every finite place and compact at infinity. If $m \equiv 2 \bmod{4}$, there is no group with those properties. 
\end{prop}

\begin{pf} (i) is obvious and (ii) is an immediate consequence of Lemma \ref{hermitian} below, since $\disc(q)=1$ 
(see \cite[chap. VI]{dieu} for the basics on hermitian forms and unitary groups). 
The other assertions (that we shall not use) follow from the Hasse's principle (\cite[\S2]{clo}).
\end{pf}

In the following lemma, we write $\disc(q) \in \Q_l^*/N(E_l^*)$ 
for the discriminant of a non degenerate $c$-hermitian form $q$ and denote by $q_0$ the {\it hyperbolic} form $q_0(x,y)=\frac{xc(y)+yc(x)}{2}$ on $E_l^2$. 
Note that $\disc(q_0)=-1$.

\begin{lemma}\label{hermitian} 
Let $q$ be a non degenerate $c$-hermitian form on $E_l^m$. \begin{itemize} 
\item[(a)] If $m$ is odd, then $q$ is equivalent to  
$$\sum_{i=1}^{\frac{m-1}{2}}q_0(z_{2i-1},z_{2i})+ (-1)^{\frac{m-1}{2}}\disc(q) N(z_m).$$ For $\lambda\in \Q_l^*$, $\disc(\lambda q) \equiv \lambda \disc(q)$, therefore there is a unique non-degenerate $c$-hermitian form up to a scalar.
\item[(b)] If $m$ is even, then $q$ is equivalent to $$\sum_{i=1}^{\frac{m-2}{2}}q_0(z_{2i-1},z_{2i})+N(z_{m-1})+ (-1)^{\frac{m}{2}-1}\disc(q) N(z_m).$$ The index of $q$ is $m/2$ if and only if $(-1)^{m/2}\disc(q) \in N(E_l^*)$. 
\end{itemize}
\end{lemma}

\begin{pf} Recall that a quadratic form on $\Q_l^s$ with $s\geq 5$ always has a zero (see e.g.  \cite[Chap. IV thm. 6]{Serreari}). We may view $E_l^m$ as a $\Q_l$-vector space of rank $2m$ and $q$ as a quadratic form on that space, so $q$ has a zero when $m \geq 3$. As a consequence, $q$ contains a hyperbolic plane and we may assume $m=2$ by induction (or $m=1$, but this case is obvious). Applying the previous remark to the form $q((z_1,z_2)) - N(z_3)$ on $F^3$, we get that $q(v)=1$ for some $v\in E_l^m$, which concludes the proof.
\end{pf}
%In particular, two unitary groups over $\Q$ which are isomorphic over $\Q_l$ for every prime $l$, and over $\R$, are isomorphic over $\Q$. Moreover,
%when $m \equiv 2 \bmod 4$, there is no twisted form $H$ of $\GL_m$ over $\Q$ such that $H(\R)$ is compact and $H(\Q_l)$ is quasisplit for each prime $l$.

\subsubsection{Automorphic forms and representations}\label{defautformG}

Let $G$ be a definite unitary group. 
We denote by $\AAA$ the $\Q$-algebra of $\Q$-ad\`eles and 
$\AAA \longrightarrow \AAA_f$ the projection to the finite ad\`eles. We have the following two important finiteness results: \begin{itemize}
\item[i)] As $G(\R)$ is compact, $G(\Q)$ is a discrete subgroup of $G(\AAA_f)$, hence for each compact open subgroup $K \subset G(\AAA_f)$, the arithmetic group $K \cap G(\Q)$ is finite. 
\item[ii)] By Borel's general result on the finiteness of the class 
number (\cite{borel}), for any $K$ as above $G(\Q)\backslash G(\AAA_f)/K$ is finite.
\end{itemize} 

The space of automorphic forms of $G$ is the representation of $G(\AAA)$ by right translations on the 
space $A(G)$ of complex functions on $X:=G(\Q)\backslash G(\AAA)$ 
which are smooth and $G(\R)$-finite. The space $X$ is compact by i) and ii). 
It admits a $G(\AAA)$-invariant finite Radon measure, so that $A(G)$ is a pre-unitary representation.

\begin{lemma}\label{admautG} The representation $A(G)$ is admissible and the 
direct sum of irreducible representations of $G(\AAA)$:
\begin{equation}\label{disc} A(G)=\bigoplus_{\pi} m(\pi)
\pi,
\end{equation}

\noindent where $\pi$ describes all the (isomorphism classes of) irreducible admissible representations of $G(\AAA)$, and 
$m(\pi)$ is the (always finite) multiplicity of $\pi$ in the above space. \ps

\end{lemma}

It will be convenient to denote by $\Irr(\R)$ 
the set (of isomorphism classes) of irreducible complex continuous representations of 
$G(\R)$ (hence finite dimensional). 
For $W \in \Irr(\R)$, we define $A(G,W)$ to be the $G(\AAA_f)$-representation by right translation on the 
space of smooth vector valued functions 
$f: G(\AAA_f) \longrightarrow W^*$ such that $f(\gamma g)=\gamma_{\infty}f(g)$ for all $g \in G(\AAA_f)$ and 
$\gamma \in G(\Q)$.

\begin{pf} As $G(\R)$ is compact the action of $G(\R)$ on $A(G)$ is completely reducible, hence 
as $G(\AAA)=G(\R) \times G(\AAA_f)$ representation we have:
$$A(G) = \bigoplus_{W \in Irr(\R)}W \otimes(A(G)\otimes W^*)^{G(\R)}.$$
But we check at once that the restriction map $f \mapsto f_{| 1 \times G(\AAA_f)}$ induces a $G(\AAA_f)$-equivariant isomorphism 
$$(A(G)\otimes W^*)^{G(\R)} \isomo A(G,W).$$ As a consequence, ii) shows that 
$A(G)$ is admissible, which together with the pre-unitariness of $A(G,W)$ 
proves the lemma. 
\end{pf}

\begin{definition}\label{defautrepG} An irreducible representation $\pi$ of $G(\AAA)$ is said to be {\it automorphic} if $m(\pi) \neq 0$. 
%If $W \in \Irr(\R)$ and $K \subset G(\AAA_f)$ is a compact open subgroup, define the space of 
%automorphic forms of weight $W$ and level $K$ as $A(G,W)^K$. 
%It is a finite dimensional complex vector space 
%equipped with an action of the Hecke algebra $\cal{C}_c(K\backslash G(\AAA_f)/K,\Q)$.
\end{definition}
\par \ps
Let $W \in \Irr(\R)$ and let us restrict it to $G(\Q) \hookrightarrow G(\R)$. As is well known (see \S\ref{repreal}), 
$W$ comes from an algebraic representation of $G$, hence the choice of an embedding $\Qb \longrightarrow \C$
equips $W$ with a $\Qb$-structure $W(\Qb)$ which is $G(\Qb)$-stable; $W(\Qb)$ is necessarily unique up to . As a consequence, 
the obviously defined space $A(G,W(\Qb))$ provides a $G(\AAA_f)$-stable $\Qb$-structure on $A(G,W)$.

\begin{cor}\label{pifdefqb} If $\pi = \pi_{\infty} \otimes \pi_f$ is an automorphic representation of $G$, then $\pi_f$ is defined over a number field.
\end{cor}

\subsection{The local Langlands correspondence for $\GL_m$}

\label{localLanglands}

Let $m \geq 1$ be any integer and $p$ a prime.
Let $F/\Q_p$ be a finite extension, $\W_F$ its Weil group, 
$I_F \subset \W_F$ the inertia group, and $|.|$ the absolute value of $F$ such that the norm of a uniformizer is the reciprocal of the number of elements of the residue field. We normalize the reciprocity isomorphism of local class-field theory $$\rec: F^* \longrightarrow \W_F ^{\rm ab}$$ \index{rec@$\rec@$, the 
reciprocity isomorphism of class field theory}
so that uniformizers correspond to geometric Frobenius elements. By an 
$m$-dimensional Weil-Deligne representation $(r,N)$ of $F$ 
we mean the data of a continuous homomorphism 
$$r: \W_F \longrightarrow \GL_m(\C)$$
such that $r(\W_F)$ consists of semi-simple elements, 
and of a nilpotent matrix  $$N \in M_m(\C),$$
satisfying $r(w)Nr(w^{-1})=|\rec^{-1}(w)|N$ for all $w \in \W_F$. \ps

Recall from \cite[Thm. A]{HT} that the Langlands correspondence is known for the group $\GL_m(F)$, and we shall use it with the normalization given {\it loc. cit}. This parameterization is a bijection 
$$\pi \mapsto L(\pi)=(r(\pi),N(\pi))$$
\index{Lpi@$L(\pi)=(r(\pi),N(\pi))$, the representation of the Weil-Deligne 
group attached to $\pi$ by local Langlands}
between the set $\Irr(\GL_m(F))$ of isomorphism classes of irreducible 
smooth complex representations $\pi$ of $\GL_m(F)$ 
and the set of isomorphism classes of $m$-dimensional 
Weil-Deligne representations of $F$. It satisfies various properties. 
For example:
\begin{itemize}
\item[-] When $m=1$, $\GL_1(F)=F^*$, we have 
$N(\chi)=0$ and $r(\chi)=\chi \circ \rec^{-1}$ for any smooth character $\chi: F^* \longrightarrow \C^*$. In general, the $L$-parameter of the central character of $\pi$ is $\det(L(\pi))$, and for any smooth character $\chi: F^* \longrightarrow \C^*$, $L(\pi \otimes \chi \circ \det)=L(\pi) \otimes L(\chi)$. 
\item[-] $\pi$ is superscuspidal (resp. ess. square integrale) if, and only if, $L(\pi)$ is irreducible (resp. indecomposable).
\item[-] If $\pi_i$ is an ess. square integrable representation of $\GL_{m_i}(F)$, and $\sum_i m_i=m$, then $\oplus_i L(\pi_i)$ is the $L$-parameter of the Langlands quotient 
$\bigotimes_i \pi_i$ (when it makes sense).
\end{itemize}

\subsection{Refinements of unramified representation of $\Gl_m$}\label{sectionref}

\label{refinementslocal}

In this subsection, we explain some aspects of 
the representation theoretic counterpart\footnote{Actually, the theory developed in this
part is comparatively much simpler than the Galois theoretic one of 
section \ref{trianguline}, as we are reduced here to see refinements as some 
orderings of some Frobenius
eigenvalues in the complex world. The relation could certainly be pushed much further, in the
style of the work of M. Emerton for $\GL_2(\Q_p)$ \cite{Emgl2}.} of 
the theory of refinements developed in section \ref{trianguline}. 
The simplest example of this notion is the well known fact that any classical modular eigenform of 
level $1$ (weight $k$, say) generates a two-dimensional vector space of $p$-old forms of
level $\Gamma_0(p)$. These old forms all have the same $T_l$-eigenvalues for
$l \neq p$, and the Atkin-Lehner $U_p$ operator preserves this two-dimensional space with characteristic
polynomial $X^2-t_pX+p^{k-1}$. \ps
From a representation theoretic point of view, 
this last computation is a purely
local statement, namely the computation of the characteristic polynomial of
$U_p$, a specific element of the Hecke-Iwahori algebra, 
on the space of Iwahori invariants of a 
given irreducible unramified smooth representation of $\GL_2(\Q_p)$. In what
follows, we explain how this theory generalizes to $\GL_m(\Q_p)$, focusing 
essentially on the unramified case. In \cite[\S 4.8]{Ch} and 
\cite[\S 6]{BC}, we 
explained how to deduce them from the Bernstein
presentation of the Hecke-Iwahori algebra. Here we use an alternative
approach based on the Borel isomorphism and the geometrical lemma.

\subsubsection{The Atkin-Lehner rings}\label{atkinlehner}

Let $F$ be a finite extension of $\Q_p$ with
uniformizer $\varpi$ and ring of integer $\anneau_F$. 
We denote by $G$ the
group $\GL_m(F)$, $B$ its upper Borel subgroup, $N$ the unipotent 
radical of $B$,
and $T$ the diagonal torus of $G$. Let $K:=\GL_m(\anneau_F)$, $T^0=K\cap T$, and let $I$ be the Iwahori subgroup of $G$ consisting of elements
which are upper triangular modulo $\varpi$. \ps
The Hecke-Iwahori algebra is the
$\Z[\frac{1}{p}]$-algebra $\Cc(I\backslash G/I,\Z[\frac{1}{p}])$ of bi-$I$-invariant 
and compactly supported functions on $G$ with values in $\Z[\frac{1}{p}]$,
for the convolution product normalized such that $I$ has mass $1$. If $g \in G$, we denote by $[IgI]$ the
characteristic function of $IgI$. We introduce now two important subrings of $\Cc(I\backslash
G/I,\Z[\frac{1}{p}])$, that we call the {\it Atkin-Lehner rings}\footnote{Following
Lazarus.}. Let $U \subset T$ be the subgroup of diagonal elements
whose entries are integral powers of $\varpi$, $U^- \subset U$ the 
submonoid whose elements have the form 
$$\diag(\varpi^{a_1},\varpi^{a_2},\cdots,\varpi^{a_m}),
\, \, \, \, a_i \in \Z, \, \, \, \, \, \, a_{i}\geq a_{i+1}\, \, \,  \forall i.$$ 
We define $\ATL^- \subset \Cc(I\backslash G/I,\Z)$ as the subring generated by the $[IuI]$, $u \in U^-$.
Recall that for each $u \in U^-$, $[IuI]$ is invertible in $\Cc(I\backslash G/I,\Z[\frac{1}{p}])$, hence it makes sense to define also
$$\ATL \subset \Cc(I\backslash G/I,\Z[\frac{1}{p}])$$ 
as the ring generated by the elements $[IuI]$, $u \in U^-$, and their inverse. 

\begin{prop}\label{atkinl} \begin{itemize} \item[(i)] The subset 
$M:=IU^-I \subset G$ is a submonoid, and the map $M \longrightarrow U$, 
$iui' \mapsto u$, is a well defined homomorphism.
\item[(ii)] The map $U^- \longrightarrow \ATL^-$, $u \mapsto [IuI]$, extends
 uniquely to ring isomorphisms $\Z[U^-] \isomo \ATL^-$ and $\Z[U] \isomo \ATL$.
\end{itemize}
\end{prop}

We warn the reader that when $m>1$, the above homomorphism does not send any $u \in U$ to $[IuI]$, but rather on $[IaI].[IbI]^{-1}$ for any $a,\, b \in U^-$ such that
$u=ab^{-1}$.

\begin{pf} 
By \cite[Lemma 4.1.5]{Casselman}, $M:=\coprod_{u \in U^-} IuI \subset G$ is a disjoint union, 
$\forall u,u' \in U^-, \, \, \, IuIu'I=Iuu'I$, hence $M$ is a submonoid of $G$, and also 
$[IuI].[Iu'I]=[Iuu'I]$, which proves (i) and the first part of (ii). The proposition follows then from the easy fact that $U^- \rightarrow U$ is the 
symmetrisation of the monoid $U^-$. 
\end{pf}

\begin{example}\label{atlcomp} As a consequence of Prop. \ref{atkinl}, 
we will systematically view $\ATL$-modules as $U$-modules. For example, let $\pi$ be a smooth representation of $G$, 
say with complex coefficients. The vector space $\pi^I$ of Iwahori invariant vectors inherits a
$\C$-linear action of $\Cc(I\backslash G/I,\Z[\frac{1}{p}])$, hence of $\ATL$, hence is a $U$-module in a natural way.
It turns out that this $U$-module structure on $\pi^I$ 
is related to the Jacquet-module of $\pi$ via the following result of Borel-Casselman. 
\end{example}

If $V$ is a representation of $G$, we denote by $V_N$ the Jacquet-module of $V$ with respect to $N$ (see e.g. \S\ref{jacquetreview}), that is the space of coinvariants of $N$, with its natural action of $T$. 

\begin{prop}\label{borelI} For any smooth complex representation $\pi$ of $G$,
the natural map $$\pi^I \longrightarrow
(\pi_N)^{T^0}\otimes \delta_B^{-1},$$
is a $\C[U]$-linear isomorphism.
\end{prop}

\begin{pf} 
As the $[IuI]$ are invertible in the Hecke-Iwahori algebra, we have $\pi^I=[IuI].\pi^I$ for each $u \in U^-$.
The result follows then from Prop. 4.1.4 and Lemma 4.1.1 of \cite{Casselman}, and from the fact that 
$[IuI]v=\delta_B^{-1}(u)\cal{P}_I(u(v))$ for each $u \in U^-$ and $v \in \pi$ by Lemma 1.5.1 of {\it loc. cit}.
\end{pf}

\subsubsection{Computation of some Jacquet modules} \label{computejacquet}
In order to use the previous result, we recall now the computation of the Jacquet module of some induced representations, following \cite{BZ}.
Fix $P \supset B$ a parabolic subgroup of $G$, $L$ its Levi component containing $T$. 
Let $\chi: L \longrightarrow \C^*$ be a smooth character, viewed also as a character on $P$ which is trivial on the unipotent radical of $P$.
Denote by $\Ind_B^G(\chi)$ the unitary smooth parabolic induction of $\chi$, that is the space of complex valued smooth functions $f$ on $G$ such that
$$f(pg)=\chi(p)\delta_P(g)^{1/2}f(g),\, \, \, \forall \, p \in P, \, \, g \in G,$$
viewed as a $G$ representation by right translations. Here $\delta_P$is the module character of $P$. \ps
Let $\coprod_{i=1}^r I_i=\{1,\cdots,m\}$ be the ordered partition associated to $P$. If $ m_i=|I_i|$, 
then $L=\prod_{i=1}^s \GL_{ m_i}(F)$. The subgroup $W=\got{S}_m$ of permutation of $\{1,...,m\}$ is a subgroup of $G$ in 
the usual way ($w=(\delta_{i,w(j)})$). Let $W(P) \subset W$ be subset of elements $w \in W$ such that $w(k)<w(l)$ whenever $k<l$ and both $k$ and $l$ belong to the same $I_i$.
The group $W$ acts on the characters of $T$ by the formula $\psi^w(t)=\psi(w^{-1}tw)$. 
Moreover, $\chi$ may be viewed as a character of $T$ by restriction $T \subset P$.

\begin{prop}\label{lemmegeom} The semi-simplification of the $\C[T]$-module $(\Ind_P^G \chi)_N$ is 
$$\oplus_{w \in W(P)} \chi^w\delta_B^{1/2}.$$
\end{prop}

\begin{pf} This is a special case of the general geometrical lemma \cite[Lemma 2.12]{BZ} (see also \cite[Theorem 1.2]{Ze}).
\end{pf}

\subsubsection{Unramified representations}\label{unramified}

An irreducible smooth representation of $G$ is said to be unramified if it has a non zero vector invariant by $K$. 
The classification of unramified representations is well known and due to Satake. \ps

Let $\chi: T \longrightarrow \C^*$ be a smooth character. It will be convenient for us to write $\chi$ as a product of smooth 
characters $\chi_i: F^* \longrightarrow \C^*$ such that
$$\chi((x_1,\cdots,x_m))=\prod_{i=1}^m \chi_i(x_i).$$

Assume that $\chi(T^0)=1$ and consider the induced representation $\Ind_B^G(\chi)$. 
As is easily seen, the space of its $K$-invariant vectors is one-dimensional hence 
this induced representation has a unique unramified sub-quotient $\pi(\chi)$. 
It turns out that : 
\begin{itemize}
\item[-] $\pi(\chi) \simeq \pi(\chi')$ if, and only if, $\chi=\chi^w$ for some $w \in \got{S}_m$ (see \S\ref{computejacquet} for this notation). 
\item[-] each unramified representation is isomorphic to some $\pi(\chi)$.
\end{itemize}\ps

The Langlands parameter of $\pi(\chi)$ is easy to describe. The isomorphism class of Weil-Deligne representations $L(\pi(\chi))=(r,N)$ associated to $\pi(\chi)$ 
satisfies $N=0$, $r(I_F)=1$ (hence the name unramified). It is uniquely determined by the conjugacy class of 
the image of a geometric Frobenius element of $W_F$, namely the class of  
$$\diag(\chi_1(\varpi),\cdots,\chi_m(\varpi)) \in \GL_m(\C).$$
Of course, this diagonal element is unique only up to permutation. We will frequently refer to this class as the semi-simple conjugacy class associated to $\pi(\chi)$.
\ps

\subsubsection{Refinements}\label{refinements} We fix $\pi$ an irreducible unramified representation of $G$. 

\begin{definition}\label{defrefinement} A refinement of $\pi$ is an ordering of the eigenvalues 
of the semi-simple conjugacy class above associated to $\pi$. In an equivalent way, 
a refinement of $\pi$ is a character $\chi: 
U \simeq T/T^0 \longrightarrow \C^*$ such that $\pi \simeq \pi(\chi)$, the dictionary being 
$$\chi \mapsto (\chi_1(\varpi),\cdots,\chi_m(\varpi)).$$
\end{definition}

Let us chose some refinement $\chi$ of $\pi$, so that $\pi \simeq \pi(\chi)$ is an 
irreducible subquotient of $\Ind_B^G \chi$. By Propositions \ref{borelI} and \ref{lemmegeom}, we get that as a $U$-module 

\begin{equation}\label{sousquotI} (\pi^I)^{\ss} \hookrightarrow (\Ind_B^G \chi)^{I,\ss} \simeq \oplus_{w \in \got{S}_m} \chi^w \delta_B^{-1/2}.
\end{equation}

As a corollary, we have the following Proposition-Definition.

\begin{definition}\label{defrefi2} If a character $\chi\delta_B^{-1/2}: U \longrightarrow \C^*$ occurs in $\pi^I$, 
or equivalently if $\chi\delta_B^{1/2}$ occurs in $\pi_N^{T_0}$, then $\chi$ 
is a refinement of $\pi$. We say that a refinement of $\pi$ is accessible if it occurs this way.
\end{definition}

For most of the representations $\pi$, all the refinements are accessible. 
Indeed, by \cite{BZ}[Theorem 4.2] and formula (\ref{sousquotI}) 
we have the following positive result. Set $$q:=|\anneau_F/\varpi|.$$
%Let $|.|$ denotes absolute value of the arithmetic norm of $F/\Q_p$.

\begin{prop}\label{casgenerique} Assume that $(\chi_i/\chi_j)(\varpi) \neq q$ for all $i \neq j$. Then $\pi(\chi)=\Ind_B^G(\chi)$ and all the 
refinements of $\pi(\chi)$ are accessible.
\end{prop}

\begin{example}\label{exrefi} \begin{itemize} \item[i)] If $\pi \simeq \pi(\chi)$ is tempered, then $\chi$ is known to be unitary hence Proposition \ref{casgenerique} applies. 
More generally, if $\pi$ is generic Proposition \ref{casgenerique} applies.
\item[ii)] On the opposite, if $\pi$ is the trivial representation then it has a unique accessible refinement, namely $\delta_B^{-1/2}$. It corresponds then to the ordering 
$$(q^{\frac{m-1}{2}},q^{\frac{m-3}{2}},\cdots,q^{-\frac{m-1}{2}}).$$ 
%$$\delta^{1/2}_B=(|\varpi|^{\frac{m-1}{2}},|\varpi|^{\frac{n-3}{2}},\cdots,|\varpi|^{-\frac{m-1}{2}}).$$
\item[iii)] Actually, by \cite{BZ}[Rem. 4.2.2], $\pi(\chi) \simeq \Ind_B^G(\chi)$ if and only if the assumption of Prop. \ref{casgenerique} 
is satisfied.
\end{itemize}
\end{example}

\subsubsection{Accessible refinements for almost tempered unramified representations}\label{almtempref}
In the applications, we will need to study the accessible 
refinements of some $\pi$ which are not tempered,
 which leads us to introduce the following class of unramified representations. \ps

Let $\pi$ be an unramified irreducible representation of $G$, and $X$ 
the set of eigenvalues (with multiplicities) of the semi-simple conjugacy class 
attached to $\pi$, $|X|=m$. Assume that $X$ has a partition $X=\coprod_{i=1}^r X_i$ such that : \ps
\begin{itemize}
\item[(AT1)] for each $i$, $X_i$ has the form $\{ x, x/q,\cdots,x/q^{ m_i-1}\}$ 
with $ m_i=|X_i|$,
\item[(AT2)] the real number $|\prod_{x \in X_i} x|^{1/ m_i}$ does not depend on $i$.
\end{itemize}

\begin{prop}\label{accessalmtemp} The accessible refinements of $\pi$ are the orderings $(x_1,\cdots,x_m)$ on $X$ such 
that whenever $x_k$ and $x_l$ are in the same $X_i$ and $x_k=q x_l$, then $k<l$.
\end{prop}

\begin{pf}  Let us choose a refinement $(x_1,\cdots,x_m)$ of $\pi$ satisfying the condition of the 
statement and such that $\{x_1,\dots,x_{ m_1}\}=X_1$, 
$\{x_{ m_1+1},\dots,x_{ m_1+ m_2}\}=X_2$ and so on. It exists by (AT1). Let 
$\chi: T \longrightarrow \C^*$ be the corresponding character, $\pi$ is then the unramified subquotient of 
$\pi(\chi)$ and we are going to identify it. \ps
Consider the standard parabolic $P$ of $G$ with Levi subgroup 
$$L=\GL_{m_1}(F) \times \GL_{ m_2}(F) \times \cdots \times \GL_{ m_r}(F).$$
One checks immediately that the character $\chi\delta_B^{1/2}(\delta_P)_{|B}^{-1/2}$
of $T$ extends uniquely to a character $\psi: L \longrightarrow \C^*$. Explicitly, 
$\psi(g_1,\cdots,g_r)=\prod_i \psi_i(g_i)$, where $\psi_i$ is the unramified character of 
$\GL_{ m_i}(F)$ obtained by composing the determinant $\GL_{ m_i}(F) \rightarrow F^*$ with 
the character of $F^*$ trivial on $\anneau_F^*$ and sending $\varpi$ to the element 
$$y_i:=x q^{-\frac{m_i-1}{2}},$$
where $x$ is the element of $X_i$ appearing in (AT1). As a consequence, we have an inclusion of $G$-representations:

$$\Ind_P^G \psi \subset \Ind_B^G \chi.$$

Up to a twist, we may assume that the real number occurring in property (AT2) is $1$. In terms of the $y_i$, it means that 
$|y_i|=1$ for all $i$, {\it i.e.} that $\psi$ is unitary. A theorem of Bernstein \cite{Bern2} shows then that $\Ind_P^G \psi$ is irreducible. As it contains 
obviously the $K$-invariant vectors of $\Ind_B^G \chi$, we conclude that $$\pi \simeq \Ind_P^G \psi.$$ The proposition is then an immediate consequence of Propositions \ref{borelI} and \ref{lemmegeom}.
\end{pf}

\begin{definition}\label{defalmtemp} Let us say that $\pi$ is {\it almost tempered} 
if $X$ admits a partition $\{X_i\}$ satisfying $(AT1)$ and $(AT2)$, or equivalently if up to a twist $\pi$ is the full parabolic induction of a unitary character.
\end{definition}

The equivalence of the two definitions above is a consequence of the proof of Proposition \ref{accessalmtemp}.

\begin{example} \label{exalmtemp}\begin{itemize}
\item[i)] If $\pi$ is tempered, it is almost tempered (and the $X_i$'s are singletons). If $\pi$ is one dimensional, it is also 
almost tempered, for the trivial partition of $X$ in one subset. 
\item[ii)] Assume that $\pi$ is the local component of a discrete automorphic representation of a unitary group (resp. of $\GL_m$). 
A consequence of Arthur's conjectures (resp. of Ramanujan conjecture and Moeglin-Waldspurger's theorem \cite{MoWa}) 
is that $\pi$ should be almost tempered. This is actually the main reason why we introduced this definition.
\item[iii)] We will need the following explicit example. Assume that $\pi$ is such that $X$ contains 
$q^{1/2}$ and $q^{-1/2}$ with multiplicity $1$, and all of whose other elements have norm $1$. Then the accessible 
refinements of $\pi$ are exactly the orderings of $X$ of the form
$$(\cdots,q^{1/2},\cdots,q^{-1/2},\cdots),$$
that is the ones such that $q^{1/2}$ precedes $q^{-1/2}$ in the ordering.
\end{itemize}
\end{example}

\subsection{$K$-types and monodromy for $\Gl_m$}
\label{types}

We keep the notations of the preceding subsections.

\subsubsection{An "ordering" relation on $\Irr(\GL_m(F))$.}
We define a relation $\prec_{I_F}$ on $\Irr(\GL_m(F))$ as follows. 
\ps

\begin{definition} \label{defordering} Let $\pi$, $\pi' \in \Irr(\GL_m(F))$, we will write $\pi \prec_{I_F} \pi'$ if $(r(\pi)_{|I_F},N(\pi))$ is in the Zariski closure\footnote{Or in the closure for the complex topology, which amounts to the same
here.} of the conjugacy class of $(r(\pi')_{|I_F},N(\pi'))$ in $M_m(\C)$. 
\end{definition}

\begin{remark}\label{remordering}
\begin{itemize}
\item[i)] The relation $r(\pi)_{|I_F} \simeq r(\pi')_{|I_F}$ on $\Irr(\GL_m(F))$ is called the "inertial
equivalence" relation, $\prec_{I_F}$
is a transitive and reflexive relation refining the inertial equivalence classes.
Moreover, in an inertial equivalence class, the minimal elements for $\prec_{I_F}$ are precisely the $\pi$ with $N(\pi)=0$, and each of them is actually a smallest element. 
\item[ii)] Assume $m>1$ and let $1$ and $St$ be the trivial and the Steinberg representation respectively. We have $r(1)(I_F)=r(St)(I_F)=1$, $N(1)=0$ and $N(St)$ has nilpotent index $m$, hence $1 \prec_{I_F} St$. As is well known, the $\pi \prec_{I_F} 1$ are exactly the unramified representations. 
Moreover by a well known result of Borel, the 
$\pi \prec_{I_F} St$ (i.e. $\pi \sim_{I_F} 1$) can also be abstractly characterized by the property that $\pi^I \neq 0$, where $I$ is a Iwahori subgroup, 
or which is the same by the property that $$\Hom_K(\tau,\pi) \neq 0,$$
with $\tau=\Ind_I^K 1_I$ and $K=\GL_m(\anneau_F)$.
\item[iii)] Take $m=2$ in the example above, the representation $\tau$ 
is then the direct sum of the trivial $1_H$ and the Steinberg $St_H$ representation of 
the finite group $H:=\GL_2(\F_q)$. Of course, $\pi \prec_{I_F} 1$ if and only if $\Hom(1_H,\pi) \neq 0$. 
As an exercise, the reader can check that the relation 
$\Hom_K(St_H,\pi)\neq 0$ cuts exactly the $\pi$ in the trivial inertial class which are not $1$-dimensional. 
\end{itemize}
\end{remark}

\subsubsection{Types}

Using works of Bernstein, Zelevinski, Bushnell-Kutzko and Schneider-Zink, it turns out that Remarks \ref{remordering} ii), iii)  are the simplest case of a general phenomenon. We are grateful to J.-F. Dat for
drawing our attention to the reference \cite{SZ}.

\begin{prop}\label{propquituedeschneider} Let $\pi \in \Irr(\GL_m(F))$. 
There exists an irreducible complex representation $\tau$ of $\GL_m(\anneau_F)$ such that \begin{itemize}
\item[i)] $\pi_{|\GL_m(\anneau_F)}$ contains $\tau$ with multiplicity $1$,
\end{itemize}
and for any $\pi' \in \Irr(\Gl_m(F))$,
\begin{itemize}
\item[ii)] ${\mathrm Hom}_{\GL_m(\anneau_F)}(\tau,\pi') \neq 0
\Rightarrow \, \, 
\pi' \prec_{I_F} \pi$,
\item[iii)] if $\pi'$ is tempered and $\pi' \prec_{I_F} \pi$, then ${\mathrm
Hom}_{\GL_m(\anneau_F)}(\tau,\pi') \neq 0$.
\end{itemize}
\end{prop}

\begin{pf} Up to the dictionary of local Langlands correspondence, 
this is exactly \cite[Prop. 6.2]{SZ}. 
For the convenience of the reader, we explain below the relevant translation. \ps 
Fix $\pi$ as in the statement, and let $\Omega$ be the (unique) Bernstein component of $\Irr(\GL_m(F))$ containing $\pi$. 
This component is uniquely determined by the cuspidal support of $\pi$. By the properties of the local 
Langlands correspondence, which is built from its restriction to the supercuspidal representations 
and Zelevinski's classification, this support is in turn uniquely determined by
$r(\pi)_{|I_F}$: for a $\pi' \in \Irr(\GL_m(F))$, we have 
$\pi' \in \Omega$ if and only if $\pi'$ is in the same inertial class as $\pi$.\par
The additional datum of the matrix $N(\pi)$ determines then the way Zelevinski realizes
$\pi$ as a Langlands quotient,
that is the "partition" $\cal{P}(\pi)$ such that $\pi$ lies in
$Q_{\cal{P}(\pi)}$ in the notations of \cite[\S2]{SZ}. Such a "partition" is by definition a family of
Young diagrams  (see below)
indexed by the cuspidal support of $\pi$, the size of a diagram being the multiplicity of the associated
supercuspidal. By Gerstenhaber's theorem (see Prop. \ref{gerstthm}) and by definition of the ordering on
partitions {\it loc. cit.} \S3 (which is the opposite of the dominance
ordering recalled in Appendix \ref{lemmasnilpo}), we have for a $\pi' \in \Omega$:
$$\pi' \prec_{I_F} \pi \, \, \, \Leftrightarrow\, \, \, \PP(\pi') \succ \PP(\pi).$$ 

They define then {\it loc. cit.} \S6 an explicit representation called 
$\sigma_{\cal P(\pi)}(\lambda)$ of a maximal compact subgroup of $\GL_m(F)$, 
here $\lambda$ is Bushnell-Kutzko's type of the Bernstein component $\Omega$. 
Up to conjugation we may assume that this maximal compact subgroup is $\GL_m(\anneau_F)$, 
and we set $\tau:=\sigma_{\cal P(\pi)}(\lambda)$. The proposition is then 
\cite[Prop. 6.2]{SZ}.
\end{pf}

\subsection{A class of non-monodromic representations for a 
quasi-split group}
 
\label{NMR}
\label{settingNMSRPS}

In this subsection, we let $l$ be a prime number, and $G$ the group of rational points of a connected reductive quasi-split group over a field $F$ which is a finite extension of $\Q_l$.
We denote by $S$ a maximal split torus in $G$, $T$ the centralizer of $S$ which is a maximal torus in $G$, and $B=TN$ a Borel containing $T$ (where $N$ is the unipotent  radical of $B$). We denote by $W$ the Weyl group of $S$ : $W=N(S)/C(S)=N(S)/T$; this groups acts on $T$ by conjugation. 

We denote by $F'$ a finite Galois extension of $F$ on which $G$ splits.
We denote with the same letter with a prime the set of points over $F'$
of the algebraic group defining one of the subgroup of $G$ defined above : 
hence $G'$, its Borel $B'=T'N'$, where $T'$ is a maximal torus of $G'$ (which is split). We denote by $W'$ the Weyl group of $T'$ : $W'=N(T')/T'$. We have a natural 
inclusion $W \subset W'$.

\subsubsection{Review of normalized induction and the Jacquet functor over a base ring}

\label{jacquetreview}

Let $A$ be a commutative ring $\Q$-algebra that contains a square root of $l$. 
We denote by $\delta_B$ the modulus character of $B$ which takes values in 
$l^{\Z}$ and we choose once and for all a square root 
$\delta_B^{1/2} : G \rightarrow A^\ast$ of $\delta_B$. 

We recall some terminology concerning a $A[G]$-module $M$: the module $M$ is {\it smooth}
if every $v \in M$ is fixed by some compact open subgroup $U$ of $G$ and it is $A$-admissible
if for every small enough compact open subgroup $U$, $M^U$ is a finite type $A$-module.

If $V$ is a smooth $A[T]$-module, we denote by $i_B^G(V)$ 
the normalized induction of $V$ from $B$ to $G$, that is the $A$-module 
of all locally constant functions $f:G\rightarrow V$ such that $f(bg)=
b\delta_B(b)^{1/2}f(g)$ for all $b \in B$, $g \in G$. The representation 
$i_B^G(V)$ is smooth and its formation commutes with every base change $A \rightarrow A'$.

If $M$ is an $A[G]$-module, we denote by $M_N$ 
the Jacquet module of $M$ relatively to $N$, that is the $N$-coinvariant quotient $M/M(N)$,
where $M(N)=\{v \in M, \int_{N_0} \pi(n) v \ dn=0$ for some compact subgroup $N_0 \subset N\}$, seen as a representation of $T=B/N$.
 
\begin{prop} \label{jacquet}
\begin{itemize}
\item[(a)] If $M$ is smooth, then so is $M_N$.
\item[(b)] The functor $M \mapsto M_N$ from the category of smooth $A[G]$-modules to the category of smooth $A[T]$-modules is exact, 
and commutes with $-\otimes_A M'$ for any $A$-module $M'$.
\item[(c)] If $M$ is flat (resp. if $A$ is reduced and $M$ is torsion free) as an $A$-module, then so is $M_N$.
\item[(d)] If $M$ is $A$-admissible and of finite type as an $A[G]$-module,
then $M_N$ is of finite type as an $A$-module.
\end{itemize}
\end{prop} 
\begin{pf} (a) is clear. The exactness in (b) is proved exactly as in the classical case (e.g. \cite[Proposition 3.2.3]{Casselman}) once noted than $N_0$ is a pro-$l$-group, 
hence of pro-order invertible in $A$. 

For $M'$ an $A$-module and $M$ an $A[G]$-module, we see $M\otimes_A M'$ as an $A[G]$-module for the trivial action of $G$ on $M'$. 
The natural map
\begin{equation}\label{commutejacquettensor} M_N \otimes_A M' \longrightarrow (M\otimes_A M')_N \end{equation}
is an isomorphism. Indeed, using a free presentation of $M'$ over $A$, the exactness of $V \mapsto V_N$ and the left exactness of tensor products, we are 
reduced to the case where $M'$ is free over $A$, which is obvious, hence (b) is proved.

The "torsion free" part of (c) is obvious from the exactness in (b). Assume that $M$ is flat over $A$ and let $X \hookrightarrow Y$ be an injective morphism of 
$A$-modules. Then $M \otimes_A X \rightarrow M \otimes_A Y$ is an injection of $A[G]$-module. 
Hence $(M \otimes_A X)_N \hookrightarrow (M \otimes_A Y)_N$ by (b), 
which is $M_N \otimes_A X \hookrightarrow M_N \otimes_A Y$ by (\ref{commutejacquettensor}). Thus $M_N$ is flat. 

Let us  prove (d). Since $M$ is of finite type over $A[G]$, 
we deduce easily from the compactness of $B\backslash G$ (see e.g the first paragraph of the proof of \cite[Thm 3.3.1]{Casselman}), that $M_N$ is finitely generated as 
an $A[T]$-module. Since $M_N$ is smooth and $T$ is abelian, there is a compact open subgroup $T_0$ of $T$ such that $M_N^{T_0}=M_N$. Up to 
replacing $T_0$ by a smaller group, we see by \cite[Prop. 1.4.4 and 
Thm 3.3.3]{Casselman} 
that there is a compact open subgroup $U_0$ with Iwahori factorization 
$U_0=N_0^{-}T_0M_0$ of $G$ such that the 
natural map $M^{U_0} \rightarrow M_N^{T_0}$ is surjective. Since $M^{U_0}$ is of finite type over $A$ by $A$-admissibility of $M$, 
then so is $M_N^{T_0}=M_N$.
\end{pf}

We recall the following easy reciprocity formula:
\begin{lemma} \label{reciprocity}
 If $M$ is a smooth $A[G]$-module and $T$ a smooth $A[T]$-module, we have a 
canonical isomorphism
$$\Hom_{A[G]}(M,i_B^G(V)) = \Hom_{A[T]}(M_N,V \otimes \delta_B^{1/2}).$$
\end{lemma}

\subsubsection{Non Monodromic Strongly Regular Principal Series} 
\label{defNMSRPS}

In this paragraph,  we keep the preceding notations but we also 
assume that $A=k$ is a field. We recall that a smooth character $\chi:
T \rightarrow k^\ast$ is {\it regular} when there is no $w \neq 1$ in 
$W$ such that $\chi^w=\chi$. We also recall the following elementary result 
(cf. \cite[Proposition 1]{Rodier1}) :
\begin{lemma} \label{indregular}
 Assume $\chi : T \rightarrow k^\ast$ is a smooth, regular, character. 
Then :
\begin{itemize}
\item[(a)] The representation $i_B^G(\chi)$ has a unique irreducible 
subrepresentation $S(\chi)$. 
\item[(b)] The Jacquet module $S(\chi)_N$ contains $\chi \delta_B^{1/2}$ as a 
$T$-subrepresentation.
\item[(c)] Any smooth $G$-representation $M$ such that $M_N$ 
contains $\chi \delta_B^{1/2}$ as a $T$-subrepresentation has a subquotient 
isomorphic to $S(\chi)$.
\end{itemize}
\end{lemma} 
\begin{pf} By the geometric lemma, the Jacquet module $i_B^G(\chi)^N$ is semi-simple as a 
$T$-representation and is the direct sum of the distinct characters 
$\chi^w \delta_B^{1/2}$ for $w \in W$.
Since the Jacquet functor is exact, and as $i_B^G(\chi)$ is of finite length 
(use Prop. \ref{jacquet} (d) and \cite[Rem. 3.12]{centrebern}), one and only one of the Jacquet modules of
its irreducible subquotients contains $\chi \delta_B^{1/2}$. 
Let us call this irreducible subquotient $S(\chi)$, which makes (b) tautologic.

On the other hand, by lemma~\ref{reciprocity}, the Jacquet module of
any sub-representation of $i_B^G(\chi)$ contains $\chi \delta_B^{1/2}$. 
it has an irreducible subrepresentation,
hence $S(\chi)$ is the unique irreducible subrepresentation 
if $i_B^G(\chi)$, which is (a). 
Finally if $M$ is as in (c), we have by lemma~\ref{reciprocity} a non-zero
morphism $M \rightarrow i_B^G(\chi)$. Its image admits $S(\chi)$ 
as a subrepresentation by (a), and  (c) follows. 
\end{pf}

Recall that the {\it base change} of a smooth character $\chi$ of $T$
is the character $\chi'$ of $T'$ defined as $\chi':=\chi \circ \Norm$, where $\Norm: T' \longrightarrow T$ is
the Galois norm.

\begin{definition} 
A smooth character $\chi$ of $T$ is said {\it strongly regular} if its base change $\chi'$ is regular as a character of $T'$.
\end{definition}
Since $W \subset W'$, a strongly regular character $\chi$ is also regular.\ps

We now recall some notations of \cite{Rodier1}, \cite{Rodier2}.
Let $\Delta$ be the root system of $G'$ with respect to $T'$. 
Let $X^\ast(T')$ be the group of rational character on $T'$ and 
$V=X^\ast(T) \otimes_\Z \R$. The chambers of $V$ are the connected components 
of $V-\bigcup_{\alpha \in \Delta} \ker \hat{\alpha}$. The Borel subgroup $B'$ determine the choice of a ``positive'' chamber $C^+$.
 
If $\alpha \in \Delta$, its associated coroot $\hat{\alpha}$ is a 
linear form on $V$. Let ${X_\ast}(T')$ be the group of $1$-parameter subgroups of $T'$. There is a canonical pairing $X^\ast(T') \times X_\ast(T') \rightarrow \Z$. Hence
each coroot $\hat{\alpha}$ 
determines canonically a $1$-parameter subgroup 
$t^\alpha : {F'}^\ast \rightarrow T'$. 

If $\chi'$ is a smooth character $T' \rightarrow k^\ast$ we define the set $\Sigma(\chi')$ 
as the set  of the  coroots $\hat{\alpha}$ such that 
$\chi' \circ t^\alpha(a)=|a| \in k^\ast$ for every $a \in {F'}^\ast$.
When $k=\C$, Rodier's theory \cite{Rodier1} shows that if $\chi'$ is regular, the set 
$\Sigma(\chi')$ determines the reducibility of 
$i_{B'}^{G'}(\chi')$ (in particular, this representation has 
length $2^{|\Sigma(\chi')|}$).

\index{NMSRPS@NMSRPS, a class of representations of a local quasi-split group}
\begin{definition} \label{NMSRPS}
An irreducible representation of $G$ is said to be a Non Monodromic 
Strongly Regular Principal Series ({\it NMSRPS}) if it is isomorphic to a 
representation $S(\chi)$ where 
\begin{itemize} \item[(a)] $\chi$ is strongly regular.
\item[(b)] For every $\hat{\alpha} \in \Sigma(\chi')$, we have
$\hat{\alpha}(C^+) < 0.$
\end{itemize}
\end{definition}

\begin{remark} \label{RemNMSRPS}
\begin{itemize}
\item[(1)] For a split group $G'$ and a regular character $\chi'$ of $T'$,
the local Langlands correspondence has been defined
by Rodier (\cite{Rodier2}) for the subquotients of $i_{B'}^{G'}(\chi')$, 
in a way that is compatible to the usual (that is, Henniart-Harris-Taylor's) local 
correspondence in the case of $\Gl_n(F)$. The representation 
$L(\pi)=(r(\pi),N(\pi))$
of the Weil-Deligne group of $F'$ corresponding to any of those 
subquotients $\pi$
has the same $r(\pi)$ namely the composition 
$$\phi_{\chi'} : W_{F'}  \longrightarrow {}^L T'
  \longrightarrow {}^L G'$$ where the first map is the $L$-parameter of 
$\chi'$ for the torus $T'$. The action of the monodromy $N$ 
depends on the chosen subquotient. 
The hypothesis (b) is equivalent to say, by \cite[5.2]{Rodier2}, that 
$S(\chi')$ has no monodromy, that is that $N(S(\chi'))=0$. 
In other words, the $L$-parameter for $S(\chi')$ is just the
map $\psi_{\chi'}$.

\item[(2)] In the case $G'=\Gl_n(F')$, hypothesis (b) simply says 
that if $T'$ is the diagonal torus and $B'$ the upper diagonal Borel, 
and $\chi'=(\chi'_1,\dots,\chi'_n)$,
then $\chi'_i = \chi'_j |\cdot |$ implies $i>j$.

\item[(3)] There should exists a base change map, sending $L$-packets of 
$G$ to $L$-packets to $G'$, and corresponding to the obvious restriction map on
the $L$-parameters. If $\chi$ is strongly regular, 
it is natural to expect that the base change to $G'$ of the $L$-packet of $G$ containing the representation 
$S(\chi)$ contains the representation $S(\chi')$. 

In the few cases the base change has been defined, this is actually true :
when $G=\Gl_n(F)$ and $F'/F$ is cyclic, this follows immediately from the
compatibility of local base change with local Langlands correspondence and 
from remark (1) above.  
In the more interesting case that $G=U(3)$ is the (unique) 
unitary group over $F$ that splits over the quadratic extension 
$F'/F$, and $G'=\Gl_3(F')$, 
this is verified for the base change map defined by Rogawski in~\cite{Rog4}.

Hence the {\it conjectural} $L$-parameter for $S(\chi)$ should be the composition $\psi_\chi : W_F \rightarrow {}^L T \longrightarrow {}^L G$. In particular, it should be non monodromic 
(that is, $N(S(\chi))=0$ with the notation of \S\ref{localLanglands}.)
\item[(4)] Rodier's theory does not seem to have been extended to any quasi-split group, even to every unitary groups. Thus we have not felt comfortable in assuming that for any regular $\chi$ satisfying the analog of condition (b) but for the root system of $G$ (which might be not reduced), $S(\chi)$ should have a nonmonodromic $L$-parameter. This is however the case for the rank one group $U(3)$ by results of Keys and Rogawski (see \cite[12.2]{Rog4}). 
Assuming that this is true, we could replace our NMSRPS condition the weaker "NMRPS" above and all our results would still hold verbatim.
\end{itemize}
\end{remark}

\subsubsection{The NMSRPS locus is constructible}\label{constrNMSRPS}
 
In this paragraph, we keep the notations of \S\ref{jacquetreview}, and we assume moreover that $A$ is a noetherian ring. 
We suppose given an $A$-admissible smooth $A[G]$-module which is of finite type over $A[G]$. 
For every $x \in X:=\Spec(A)$, we denote by $k(x)=A_x/xA_x$ the residue field of $A_x$, and we set $M_x:=M \otimes_A A_x$ and 
$\overline{M}_x := M \otimes_A k(x)$. \ps
        Let us denote by $X_0$ the subset
of $x \in X$ such that the $k(x)[G]$-module $\overline{M}_x$ contains an irreducible subquotient which is
a NMSRPS representation.

\begin{prop} \label{NMSRPSlocus} \begin{itemize}
\item[(i)] $X_0$ is a constructible subset of $X$. In particular,
$\overline{X}_0$ contains a dense open subset $U$ such that $U \subset X_0$.
\item[(ii)] Assume that $A$ is reduced and that $M$ is torsion free over $A$. Then 
$\overline{X}_0$ is a union of irreducible components of $X$. 
\end{itemize}
\end{prop}
\begin{pf}
By proposition~\ref{jacquet}, the $A[T]$-module $E:=M_N$ is of finite type
over $A$. We view $E$ as an $A[T']$-module via the map $A[T'] \longrightarrow A[T]$ induced by the norm 
$\Norm$. Let $B$ be the image of the $A$-algebra $A[T']$ in $\End_A(E)$. 
It is a finite $A$-algebra, let $Y=\Spec(B)$ and $g: Y \rightarrow X$ the structural map. \ps
Let us consider the subset $Y_0 \subset Y$ 
of points $y \in Y$ such that the induced character $T' \rightarrow (B_x/xB_x)^\ast$ satisfies both conditions (a) and (b) of Definition
\ref{NMSRPS}. By definition, we have 
$$Y_0 = \bigcap_{w \in W'\backslash\{1\}} \left( \bigcup_{t \in T'} D(t^w-t) \right) \bigcap_{\{ \hat{\alpha} \in \hat{\Delta} |\, \hat{\alpha}(C)>0\}}
\left( \bigcup_{f \in F'} D(t^{\alpha}(f)-|f|) \right),$$
hence $Y_0$ is an open subset as both intersections are finite. \ps
We claim that $X_0=g(Y_0)$. First, by Lemma \ref{indregular} (c) and \ref{jacquet} (b), observe that
$X_0$ is also the subset of points $x$ of $X$ such that $E\otimes_A k(x)$ contains a character 
$T' \rightarrow k(x)^*$ satisfying (a) and (b) of Def. \ref{NMSRPS}, i.e. such that the support of the 
$B$-module $E\otimes_A k(x)$ meets $U$. In particular, it is clear that $X_0 \subset g(Y_0)$. \ps
Let $x \in X$. As $\widehat{A_x}$ is henselian and $g$ is finite,
        $$\widehat{B_x}:=B \otimes_A \widehat{A_x} \isomo \prod_{\{y\, |\, \, g(y)=x\}} \widehat{B_y},$$
hence we can write accordingly $\widehat{E_x}:=E\otimes_A \widehat{A_x}=\oplus_y E(y)$ as a direct sum of $B\otimes_A \widehat{A_x}$-modules. 
Moreover, by flatness
of $A \rightarrow \widehat{A_x}$, 
$\widehat{B_x}$ identifies with its image in $\End_{\widehat{A_x}}(\widehat{E_x})$ hence 
$E(y) \neq 0$ for all $y \in g^{-1}(x)$. In particular, if $x=g(y)$ with $u \in U$, then 
$$E \otimes_A k(x)=\widehat{E_x}/x\widehat{E_x} \supset E(y)\otimes_{\widehat{A_x}}k(x)$$ and the latter $B$-module is non zero by Nakayama, hence has support $\{y\}$. This proves that 
$g(Y_0) \subset X_0$, hence the equality. In particular, $X_0$ is constructible as $Y_0$ is open, which proves the first part of (i).\par
	As $X_0$ is constructible, it is a finite union of (say) nonempty locally closed irreducible subsets $F_i \cap U_i$. So 
$\overline{X}_0=\bigcup_i F_i$ and if $U'_i=U_i\cap (F_i\backslash \cup_{j \neq i} F_j\cap F_i)$, then $U:=\cup_i U'_i$ satisfies the second 
part of assertion (i).\ps
	Let us prove (ii). As $Y_0$ is an open subset, its closure is a finite number of irreducible components of $Y$. 
As $g(Y_0)=X_0$ and $g$ is finite, we only have to check that each irreducible component of $Y$ maps surjectively to an irreducible component of $X$. 
Note that $E$ is torsion free over $A$ by assumption and Proposition\ref{jacquet}, hence so are $\End_A(E)$ and $B\subset \End_A(E)$.
We conclude then by Lemma \ref{irredtorsfree}. \end{pf}

The following lemma is a variant of \cite[Lemma 2.6.10]{Ch}.

\begin{lemma}\label{irredtorsfree} Assume that $A$ is a reduced notherian ring and that $B$ is a finite, torsion free, $A$-algebra. Then each irreducible 
component of $\Spec(B)$ maps surjectively to an irreducible component of $\Spec(A)$.
\end{lemma}
\begin{pf} We check at once that a finite type $A$-module $M$ is torsion free if, and only if, it has an $A$-embedding 
$M \hookrightarrow A^n$ for some $n$. In particular, if $M$ is torsion free over $A$, then for all $x \in \Spec(A)$ the $A_x$-module 
$M_x$ is torsion free. \par 
As the finite map $g: \Spec(B) \rightarrow \Spec(A)$ is closed, 
it suffices to show that the image of the generic point $x$ of an 
irreducible component of $\Spec(B)$ is the generic point of an irreducible component of $\Spec(A)$. By localizing $A$ at $g(x)$, we may assume that $A$ is local 
and that $g(x)$ is a closed point. As $g^{-1}(x)$ is a discrete closed subspace of $\Spec(B)$, $x$ is also a closed point, 
hence it is open as it is minimal as well, and $B_x$ is a direct factor of $B$. Thus we may assume that $B=B_x$ is artinian. 
As $B \subset A^n$ and $A$ is reduced, it implies easily that $A$ is 
itself artinian, which concludes the proof.
\end{pf}

\begin{remark}\label{variantNMSRPS} (A variant) Assume that we have a finite number of quasisplit groups $G_i$, possibly over local fields of different characteristics, each one being equipped with a datum 
$(G_i,B_i,T_i,G'_i,B'_i,T'_i)$ as in the beginning of \S\ref{settingNMSRPS}. Then we may form $G=\prod_i G_i$, as well as $B$, $T$, $G'$, $B'$ and $T'$, and all the propositions and lemmas of this part \ref{settingNMSRPS} apply verbatim to this case, as all the arguments are group theoretic. For example, in this context, a NMSRPS representation of $G$ is a tensor product
of NMSRPS representations $\pi_i$ of $G_i$.
\end{remark}

\subsection{Representations of the compact real unitary group}

\label{repreal}
 
Recall that the continuous, irreducible, complex representations 
of the compact group $\U(m)(\R)$ are all finite dimensional. There are two ways to describe them : 
either by their highest weight or by their Langlands parameters. We give here both descriptions, as well as the relation between 
them. \ps
\newcommand{\uk}{\underline{k}}
If $\uk:=(k_1,\dots,k_m)\in \Z^m$ satisfies $k_1 \geq k_2 \geq \cdots \geq k_m$, we denote
by $W_{\uk}$ the rational (over $\Q$), irreducible, algebraic representation of $\GL_m$ whose highest 
weight relative to the upper triangular Borel is the character\footnote{This means that the action of the diagonal torus of $\GL_m$ on the unique 
$\Q$-line stable by the upper Borel is given by the character above.} 
$$\delta_{\uk}: (z_1,\cdots,z_m) \mapsto \prod_{i=1}^m z_i^{k_i}.$$
For any field $F$ of characteristic $0$, we get also a natural irreducible algebraic representation 
$W_{\uk}(F):=W \otimes_{\Q} F$ of $\GL_m(F)$, and it turns out that they all have this form, for a 
unique $\uk$. \ps
Let us fix an embedding $E \hookrightarrow \C$, which allows us to see $\U(m)(\R)$ as a subgroup of $\Gl_m(\C)$ 
well defined up to conjugation (see \S\ref{exunitarygp}). So for $\uk$ as above, we can view $W_{\uk}(\C)$ as a continuous representation of $\U(m)(\R)$.
As is well known, the set of all
$W_{\uk}(\C)$ is a system of representants of all equivalence 
classes of irreducible continuous representations of $\U(m)(\R)$. We will say 
that $W_{\uk}$ has {\it regular weight} if $k_1 > k_2 > \dots > k_m$. \ps

On the other hand, the $L$-parameters of the irreducible representations 
of $\U(m)(\R)$, are determined by their restrictions to the Weil group of $W_\C=\C^\ast$ of $\C$, 
which are, up to conjugation, all the morphisms $\phi : \C^\ast \rightarrow \Gl_m(\C)$
of the form $$\phi(z) = \diag((z/\bar z)^{a_1},\dots,(z/\bar z)^{a_m})$$ where $a_1,\dots,a_m \in \Z + \frac{m+1}{2}$ and $a_1 > \cdots > a_m$. 
For the proof, see \cite[Prop. 4.3.2]{BeCl}. The relation between the two descriptions is given by 
$$a_i = k_i+\frac{m+1}{2}-i, \, \, \,  i=1 \dots m.$$

\subsection[The Galois representations of an automorphic representation of $\U(m)$]{The Galois representations attached to an automorphic representation of $\U(m)$} \ps

\label{subsectionrepm}

\subsubsection{Settings and notations}\label{setnotrepm} In this section, $m \geq 1$ is an integer such that $m \not \equiv 2 \pmod{4}$.
Let us fix a prime number $p$ that is split in $E$, algebraic closures $\Qb$
of $\Q$ and $\Qpb$ of $\Qp$, and some embeddings $\iota_{p}: \Qb \rightarrow \Qb_{p}$, $\iota_{\infty}: \Qb \rightarrow \C$.\index{Zi@$\iota_p$, $\iota_\infty$, embeddings of $\Qb$ into $\Qpb$ and $\C$}
As $G(\R)$ is compact, for any automorphic representation $\pi$ then $\pi_{\infty}$ is algebraic and 
the finite part $\pi_{f}$ is defined over $\Qb$ by 
Lemma \ref{pifdefqb}, so that we may view it over $\Qpb$ using
$\iota_{p}\iota_{\infty}^{-1}$. 
Let us fix such a $\pi$. \ps

If $l=x\bar{x}$ is a prime that splits in $E$, then
we will denote by $\pi_x$ the representation of $\GL_m(\Q_l)$ deduced from
$\pi_l$ and the identification $G(\Q_l) \isomo_x \GL_m(\Q_l)$ defined by $x$
as in \S\ref{definitionUm}. \ps

We recalled in \S\ref{localLanglands} Langlands-Harris-Taylor
parameterization of complex irreducible smooth representations. This
parameterization holds actually if we replace $\C$ everywhere there 
by any algebraically closed field of characteristic $0$, {\it e.g. $\Qpb$}.
As a consequence, to each $\pi_x$ as above viewed with $\Qpb$ coefficents via
$\iota_p\iota_{\infty}^{-1}$, is attached a unique $\Qpb$-valued Weil-Deligne
representation $(r(\pi_x),N(\pi_x))$. We recall also that from Grothendieck's 
$l$-adic monodromy theorem (see e.g. the Appendix \ref{grothlmon}), for any local field
$F$ over $\Q_l$ with $l\neq p$, there is a bijection between the isomorphism
classes of continuous representations $W_F \rightarrow \GL_m(\Qpb)$, and the
isomorphism classes of $m$-dimensional $\Qpb$-valued Weil-Deligne representations
of $F$. We shall use these bijections freely in the sequel. \ps

We let $v$ denote the (split) place of 
$E$ above $p$ induced by $\iota_p: E \rightarrow \Qpb$, and by $v_{\infty}$ the complex place
of $E$ induced by $\iota_{\infty}: E\rightarrow \C$. \ps

Let $l$ be a prime that does not split in $E$. If $E$ is unramified at
$l$, then $G(\Q_l)$ is an unramified unitary group hence it makes sense to
talk about its {\it unramified representations}: they are the irreducible smooth
representations having a nonzero vector invariant under a maximal
hyperspecial subgroup. For example, for the obvious model of $\U(m)$ over
$\Z$ then $G(\Z_l)$ is maximal hyperspecial. When $l$ is ramified, we will
also say, following Labesse, that $\pi_l$ is unramified if it has a nonzero
vector invariant under a {\it very special} compact subgroup in the sense of
\cite[\S 3.6]{labesselivre}. \ps

\subsubsection{Statement of the assumption $Rep(m)$} \label{sectrepm}We now formulate a conjecture about the existence and the basic properties of
the Galois representation attached to an automorphic representation of $\U(m)$.
We expect and hope that this conjecture (and much more) will be proved in the 
forthcoming book \cite{GRFAbook} on unitary groups 
written under the direction of Michael Harris in Paris. 

To make this hope likely, we have made a special effort throughout this book
to keep under control the properties on those Galois representations we need,
and in the next conjecture to formulate the weakest statement that we need.

\index{Repm@Rep($m$), an assumption about Galois representations 
attached to automorphic representations for $\U(m)$}

\begin{conj}[Rep($m$)] \label{repm} 
Let $\pi$ be an automorphic representation of $\U(m)$ such that $\pi_\infty$ has regular weight. There exists a continuous, semisimple, Galois representation:
$$\rho_{\pi}: \G_E \longrightarrow \GL_m(\Qpb),$$
\noindent such that the following properties are satisfied:\ps
\smallskip
\noindent {\rm (P0)} if $l=xx'\neq p$ is split and 
$\pi_x$ is unramified, then 
$\rho$ is unramified above $l$ and the characteristic
polynomial of a geometric Frobenius at $x$ is given by the Langlands 
conjugacy class of $\pi_x|\det|^{\frac{1-m}{2}}$.\ps \smallskip
\noindent {\rm (P1)} If $l \neq p$ is a prime and if $\pi_l$ is unramified, 
then ${\rho_{\pi}}$ is unramified at each prime above $l$.\ps\smallskip
\noindent {\rm (P2)} If $l=xx' \neq p$ is a prime that splits in 
$E$, then the nilpotent monodromy operator of the Weil-Deligne 
representation attached to ${\rho_{\pi}}_{|\W_{E_x}}$ is in the
closure of the conjugacy class of $N(\pi_x|\det|^{\frac{1-m}{2}})$ in ${\rm M_m}(\Qpb)$. \ps \smallskip
\noindent {\rm (P3)} If $l \neq p$ is a prime, $x$ a place of $E$ above $l$, 
and $\pi_l$ is an NMSRPS (see Definition~\ref{NMSRPS})
then the Weil-Deligne representation attached to ${\rho_{\pi}}_{|\W_{E_x}}$ has a trivial monodromy. \ps\smallskip
\noindent {\rm (P4)} The $p$-adic representation 
${\rho_{\pi}}_{|\G_{E_v}}$ is De Rham, and its Hodge Tate weights are the 
integers $$-a_1+\frac{m-1}{2},\dots,-a_m+\frac{m-1}{2}$$ where  $a_1,\dots,a_m$ are 
such that the restriction to $\C^\ast$ of the $L$-parameter of $\pi_\infty$ is 
$z \mapsto \diag((z/\bar z)^{a_1},\dots, (z/\bar z)^{a_m})$ (see \S~\ref{repreal}).
\ps \smallskip
\noindent {\rm (P5)} If $\pi_p$ is unramifed, then 
${\rho_{\pi}}_{|\G_{E_v}}$ is crystalline and
the characteristic polynomial of its crystalline Frobenius is 
the same as the one of $\iota_p\iota_{\infty}^{-1}L(\pi_v|\det|^{\frac{1-m}{2}})$.
\end{conj}
\medskip

\begin{remark}\label{remsectrepm}
\begin{itemize}
\item[(i)] 
By Cebotarev density theorem, and since 
the primes of $E$ which split above $\Q$ have density $1$, 
the property $(P_0)$ alone determines $\rho_{\pi}$ up to isomorphism. 
Moreover, it implies
that  $\rho_{\pi}$ is conjugate selfdual in the sense that:
$$ \rho_{\pi}^{\bot} \simeq \rho_{\pi}(m-1).$$ 

\item[(ii)] The properties (P0), (P1) and (P4) imply that $\rho_{\pi}$ is
geometric. 

\item[(iii)] The Langlands program and Arthur's conjecture predict 
that there should exist local and global base change from 
$\U(m)$ to $\Gl_m$, that a $\pi$ with regular weight as in 
the conjecture should be tempered, 
and that for such a tempered $\pi$ 
the global base change $\pi_E$ should be compatible at every place with the 
local base change (see~\ref{basechange} below).
Moreover it also predicts the existence of a Galois representation 
$\rho_\pi$ of $G_E$, and the Weil-Deligne representation attached to
the restriction at $\W_{E_x}$ for every place $x$ 
(prime to $p$) of $E$ should be isomorphic to 
$L((\pi_E|\det|^{\frac{1-m}{2}})_x)$.

The properties (P0) to (P3) are very special cases of those predictions.
This is clear for (P0), (P1), (P2); as for (P3), if $l\neq p$ 
is a prime such that $\pi_l$ is an NMSRPS and $x$ a place of $E$ above $S$,
then the base change $(\pi_E)_x$ should be a NMSRPS of $\Gl_n(E_x)$, whose 
$L$-parameter should be non-monodromic (see remark~\ref{RemNMSRPS}).

Moreover, properties (P4) and (P5) are also part of the standard
expectations 
for the Langlands correspondence at places dividing $p$.

\item[(iv)] The property (P3) for split primes $l$ is a special case
of (P2). This should be clear from the preceeding remark.

\item[(v)] In the following chapter, the property (P2) will allow us to 
work with representation $\rho_\pi$ that have arbitrary ramification at split 
primes. However, because of the weak form of (P3), 
we shall have to assume that the ramification, if any, is of a 
very special kind at non-split places, namely is a NMSRPS. \par
On the other hand, if we are willing to restrict ourselves to representations that are either unramified or
a NMRSPS at every finite place (and not only inert and ramified), then 
the property (P2) is not necessary. The property (P3), 
supposedly easier to prove, will be enough.
 
\item[(vi)] When $m \leq 3$, the properties (P0) and (P1) and (P3) to (P5) 
are known by the work of Blasius and Rogawski 
(cf. \cite{br1} and \cite{br3} and also \cite[\S3.3]{BC}
for some details). Property $(P2)$ is not completely known 
(to the best of our knowledge) but anyway is not necessary in view of the above remark, since 
all the $\pi$ to which we shall need to apply Rep($m$) 
to are NMSRPS at every ramified place if $m \leq 3$. Thus, in the sequel, we will allow ourselves to say that
Rep($m$) is known for $m \leq 3$.

\item[(vii)] Properties (P1) to (P5), except maybe (P3), are also known for any $m$ when $\pi$ admits 
a base change to a representation satisfying the assumptions of Harris-Taylor (\cite{HT}, \cite{TY}) and which is strong at the split places. 
This includes e.g. the case of a $\pi$ that is supercuspidal at two split places (\cite[\S 3]{HL}). Unfortunately, the
representations to which we shall apply Rep($m$) will never be of this type. \par
        Moreover, let us consider the slighty different setting where $\U(m)$ is replaced by a definite
unitary group $G$ such that $G(\Q_p) \simeq \GL_m(\Q_p)$, and that for some split $l \neq p$, $G(\Q_l)$ is
isomorphic to the group of invertible elements of a central division algebra over $\Q_l$. In this case and if $\pi_l$ is not one dimensional,
the existence of $\rho_{\pi}$ satisfying (P0), (P1), (P2), (P4) and (P5) is known by \cite{HT} and \cite[Thm.
3.1.3]{HL}.

%\item[(vii)] Properties (P1) to (P5) are also known for any $m$ 
%when $\pi$ comes from a representation satisfying the
%assumption of Harris-Taylor (\cite{HT}, \cite{TY}). This includes e.g.
% the case of a $\pi$ that is supercuspidal at a split place (\cite{HL}).
%Unfortunately, the representation to which we shall apply Rep($m$) will 
%never be of this type.
%
\item[(viii)] We of course 
expect that in the forthcoming book by Harris {\it et al.},
the representation $\pi$ (or some well chosen base change of it) will be cut off in the 
Galois cohomology of some explicit local system of the Shimura variety of
some inner form of $\U(m)$. (Here the hypothesis that $\pi_\infty$ 
has regular weights  might be helpful from a technical point of view). 
 Hence (P0), (P1), (P4) and (P5) should 
follow directly from the construction and a few standard arguments
(see e.g. \cite[Prop. 3.3]{BC} for (P4) and (P5) and \cite[Prop. 3.2]{BC} 
for (P1) at ramified primes). 

The properties (P2) and (P3) are special cases, concerning monodromy,
of the compatibility of the construction of the Galois representation to 
the local Langlands correspondence that might be harder to prove. However, they only ask 
for an ``upper bond'' on the monodromy of the Galois representation, which 
is the ``easy sense'', and this should follow from an accessible (maybe 
already known) study of the local geometry of the special fiber of the relevant Shimura varieties.

At any rate, (P3) would follow easily if the base change (local
and global, with compatibility) was known, by an argument completely 
similar to \cite[Prop. 3.2]{BC}).     

\item[(ix)] Note that the Langlands and Arthur's conjectures also predict 
some irreducibility results on the representation $\rho_E$ (for example, if $\pi$ is not endoscopic). 
Those results might be much harder to prove. 
However, a feature of our method, already present in \cite{BC},
is that we have absolutely no need of them. Instead, we shall be able to 
prove, as a by-product of our work, that in many cases $\rho_\pi$ is 
irreducible, or not too reducible.
\end{itemize}
\end{remark}

\subsection{Construction and automorphy of non-tempered 
representation of $\U(m)$}

\label{constructionpin}
 
In this subsection we fix an integer $n \geq 1$ that is not 
divisible by four, and we set $m:=n+2$, so that $m \not \equiv 2 \pmod{4}$ as 
above. For a representation $\pi$ of $\Gl_n(E_v)$, 
$v$ non split place of $E$ (resp. of
$\Gl_n(\A_E)$), we note $\pi^\bot$ the representation 
$g \mapsto \pi^\ast(c(g))$,
where $\pi^\ast$ is the contragredient of $\pi$ and $c$ denotes the map on $\Gl_n(E_v)$ (resp. $\Gl_n(\A_E)$) induced by the non trivial element $c \in \Gal(E/\Q)$.

\index{bot2@$\bot$, similar notion for an automorphic representation}

\subsubsection{The starting point}

\label{start}

\index{Zpi@$\pi$, an automorphic cuspidal representation of $\Gl_n(\AAA_E)$}

We start with a cuspidal tempered\footnote{This condition should be a consequence of the preceding one, according to the generalized Ramanujan's conjecture.} automorphic representation $\pi$ of $\Gl_n(\AAA_E)$.
We make the following assumptions on $\pi$ :
\begin{itemize}
\item[(i)] We have $\pi^\bot \simeq \pi$.
\item[(ii)] The $L$-parameter of $\pi_\infty$ has the form
$$z \mapsto \diag((z/\bar z)^{a_1},\dots,(z/\bar z)^{a_n})$$ where the $a_n$ are distinct half-integers (not integers), 
and are different from $\pm 1/2$. 
\item[(iii)] If $l$ is a nonsplit prime, then either
\begin{itemize}
\item[(iiia)] $\pi_l$ is unramified and its central character $\chi$ satisfies\footnote{Here $\varpi_l$ is a uniformiser of $E_l$. 
When $l$ is inert, $\chi(\varpi_l)=\chi(l)$ and the condition on $\chi$ is automatically satisfied, see Rem. \ref{remautomaticiiia}. The reason for the 
appearence here of this condition on $\chi$ basically comes from the fact that it is not equivalent for a character 
of $U(1)(\Q_l)$ to be unramified (i.e. trivial), and to have an unramified base change. However, 
the two notions coincide when $l$ is inert. It is maybe possible to circumvent this assumption by extending slightly \cite[\S 3.6]{labesselivre}.}  $\chi(\varpi_l)=(-1)^n$, or 
\item[(iiib)] the representation $\pi_l$ is a NMSRPS representation $S(\eta)$, where 
$\eta=(\eta_1,\dots,\eta_n)$ is a regular character of the standard 
maximal torus of $\Gl_n(E_l)$, and there is no (resp. exactly one) 
$i \in \{1,\dots,n\}$ if $n$ is even (resp. if $n$ is odd) such that 
$\eta_i^\bot=\eta_i$. 
\end{itemize}
\end{itemize}

The aim of this section is to describe, place by place, 
a representation of $U(m)$, called $\pi^n$ (the $n$ stands for "non-tempered", as $\pi^n$ turns out to be non tempered at every finite place) 
depending on $\pi$, and to state a conjecture AC($\pi$) 
(known in low dimensions) 
that $\pi^n$ is automorphic if (and actually only if)
$\varepsilon(\pi,0)=-1$. 
The representation  $\pi^n$ is an endoscopic 
transfer of $\pi$, 
and the conjecture we state is a 
particular case of the far reaching multiplicity formula of Arthur, 
as will be explained in Appendix A. 

\begin{remark}
\begin{itemize}
\item[(i)] When $n$ is even, properties (i) and (ii) are conjecturally 
sufficient conditions for $\pi$ to be the base change of a discrete automorphic representation of the quasisplit unitary group $\U(n)^*$ attached to $E/\Q$. When $n$ is odd, on the contrary, a representation satisfying (ii) is not a base change from $\U(n)^*$, but (i) and (ii) should rather ensure that $\pi \otimes \mu$ is a base change from $\U(n)^*$ for any Hecke character $\mu$ as in Lemma \ref{lemmamu}(iii) (see Example \ref{examplenecnodd} below).

\item[(ii)] Suppose that there is a representation $\rho$ of 
$G_E$ with the same $L$-function as $\pi$. 
Then property (i) says that $\rho^\bot \simeq \rho$, so that the 
representation $\Ind_E^\Q \rho$ of dimension $2n$ is autodual. Now condition 
(ii) tells that that representation is symplectic, as opposite 
to orthogonal. Hence $\varepsilon(\pi,1/2)=\varepsilon(\rho,1/2)$ 
may be either $1$ or $-1$ (in the orthogonal case it would always be 
$+1$ by a theorem of Deligne and Ribet).

\item[(iii)] For $l$ a non split prime, 
if $\pi_l$ is a subquotient of the (normalized)
 induced representation of a character $\eta=(\eta_1,\dots,\eta_n)$ 
of the standard maximal torus, 
then it follows easily that $\eta^\bot = \eta^\sigma$ for some 
$\sigma$ in the Weyl group $\got{S}_n$ of the standard maximal torus\footnote{Indeed, $\pi_l^\bot \simeq \pi_l$, hence the induced representations of $\eta^\bot$ and $\eta$ share an irreducible subquotient,
which is well known to imply the desired equality.}. If $\eta$ is regular, then $\sigma$ is an involution and  
$\eta_i^\bot = \eta_{\sigma(i)}$ for all $i=1 \dots n$. We thus see that
the condition 
in (iiib) amounts to say that $\sigma$ has the smallest possible number 
of fixed points for an involution in $\got{S}_n$, 
namely $0$ or $1$ depending whether $n$ is even or odd.

\item[(iv)]
Property (iii) is not really needed for the conjecture we are going to state, 
but it simplifies the exposition, allowing in particular to give a non 
conjectural description of $\pi^n$ at non split places. 

To be more precise, and conjecturally speaking, 
condition (iii) on $\pi_l$ is the 
condition needed for $\pi^n_l$ to be either unramified or an NMSRPS at $l$.
Without this condition, there should still exist a $\pi^n_l$ with suitable
properties, but it could be square integrable or even supercuspidal, and 
it is not possible in the present state of knowledge 
on the representation theory of local unitary groups to 
construct the needed representation. 

Moreover, the hypothesis that $\pi^n_l$ is unramified or an NMSRPS is
what we will need in the following sections to be able to deal with the monodromy at 
the nonsplit $l$. So it is not a big loss to assume it from now.
\end{itemize}
\end{remark}

\subsubsection{Hecke characters} \label{introcarhecke}

If $\mu$ is a Hecke character of $E$, that is a continuous morphism $$\mu:\A_E^\ast/E^\ast \longrightarrow \C^\ast,$$ 
recall that $\mu^\bot$ is the Hecke character 
$x \mapsto \mu(c(x))^{-1}$. 
We say that a Hecke character $\mu$ 
{\it descends to $U(1)$} if $\mu = \psi(x/c(x))$, for some
continuous character $\psi$ of $$\U(1)(\A_E)=\{ x \in \A_E^\ast, \ xc(x)=1\}.$$ Obviously a character $\mu$ that descends to 
$\U(1)$ satisfies $\mu^\bot=\mu$. If a character satisfies $\mu^\bot=\mu$, we have $\mu_{\infty}(z)=(z/c(z))^a$ for all 
$z \in (E\otimes {\mathbb R)}^\ast$ and some {\it weight}  $a$ which is either an integer or a half integer.

\begin{lemma} \label{lemmamu}
\begin{itemize} 
\item[(i)] The subgroup of Hecke characters of $E$ that descend to 
$U(1)$ is of index $2$ in the group of all Hecke character of $E$ 
that satisfy $\mu^\bot = \mu$.
\item[(ii)] 
For a Hecke character $\mu$ such that $\mu^\bot=\mu$ the 
following are equivalent :
\begin{itemize} 
\item[-] the character $\mu$ does not descend to $U(1)$,
\item[-] the weight $a$ of $\mu$ is not an integer,
\item[-] the restriction of $\mu$ to $\A_\Q^\ast$ is the order $2$ 
character $\omega_{E/\Q}$
corresponding via class field theory to the extension $E/\Q$.
\end{itemize}
In particular, a character that satisfies the above conditions is 
ramified at every ramified place of $E/\Q$, since so does 
$\omega_{E/\Q}$.
\item[(iii)] There exists a Hecke character $\mu$ of $E$, 
satisfying $\mu^\bot=\mu$ with weight $1/2$ 
and which is ramified only at ramified places of $E/\Q$.
\end{itemize}
\end{lemma}
\begin{pf}
Both (i) and (ii) result from the following observations : \begin{itemize}
\item[-] a Hecke character descends to $U(1)$ if and only if it is trivial on 
$\A_\Q^\ast/\Q^\ast$,
\item[-] a Hecke character $\mu$ satisfies such that $\mu^\bot=\mu$ if and only if it is trivial on the norm group $N(\A_E^\ast/E^\ast)$, 
\item[-] by class field theory, $N(\A_E^\ast/E^\ast)=\ker\, \,  \omega_{E/\Q}$ is an open subgroup of index $2$ in $\A_\Q^\ast/\Q^\ast$. 
\end{itemize}

For (iii), let $S$ be the set of rational primes that ramify in $E$.
For each $l \in S$, choose any finite order character
$\mu_l: \anneau_{E\otimes \Q_l}^{\ast} \rightarrow \C^{\ast}$ extending $\omega_{E/\Q,\, l}$.
We fix also an isomorphism $E\otimes \R \isomo \C$ for convenience, and set $\mu_{\infty}(z)=(z/\bar z)^{1/2}$ for $z \in \C$. Assume first that the cyclic group $U=\langle u \rangle$ of units in ${\mathcal O}_E$, is reduced to $\{\pm 1\}$. Then we can define $\mu$ on 
$\C^{\ast} \times {\widehat{\anneau}_E}^{\ast}$
to be $\mu_{\infty}\prod_l \mu_l$. As $\mu$ coincides with $\omega_{E/\Q}$ on 
$\R^{\ast}\times \widehat{\Z}^{\ast}$, $\mu(U)=\{1\}$.
As $E^{\ast} \cap (\C^{\ast} \times {\widehat{\anneau}_E}^{\ast})=U$, $\mu$ extends uniquely to an $E^{\ast}$-invariant continuous character of the open subgroup
$G:=E^{\ast} (\C^{\ast} \times {\widehat{\anneau}_E}^{\ast})$ of $\mathbb{A}_E^{\ast}$. Note that $G$ is open of finite index in $\mathbb{A}_E^{\ast}$ by the
finiteness of the class number of $E$, hence we can extend $\mu$ furthermore the $\mu$ above to a continuous character of $\AAA_E^{\ast}/E^*$. 
Note that $G$ contains ${\A_\Q}^{\ast}$ and that $\mu$ extends $\omega_{E/\Q}$ by construction, hence $\mu^{\bot}=\mu$, which concludes the proof. \ps
	When $|U|>2$, then $E=\Q(i)$ or  $\Q(j)$, and $S=\{l\}$ contains only one prime. In
this case, note that $U \cap \Z_l^{\ast}=\{\pm 1\}$ hence we may first extend
$\omega_{E/\Q,\, l}$ to $U\Z_l^{\ast}$ by choosing $\mu_l(u):=u^{-1}$,
and then extend it anyhow to a finite order character of ${\cal O}_{E\otimes \Q_l}^{\ast}$. Again, if
we define $\mu$ as before, we have $\mu(U)=\{1\}$ and the same proof works.
\end{pf} 
\begin{example}\label{examplenecnodd} The central character of $\pi$ is a Hecke character $\mu$ that satisfies $\mu^\bot=\mu$ by condition (i) of~\S\ref{start}, 
with weight $a=\sum_{i=1}^n a_i$, which is an integer if and only if $n$ is even by (ii) of~\S\ref{start}. Hence
 $\mu$ descends to $U(1)$ if and only if $n$ is even, which is also the conjectural condition for $\pi$ to descend to $\U(n)$.
\end{example}
\index{Zmu@$\mu$, a Hecke character of $E$}

\begin{remark}\label{remautomaticiiia} Assume that $l$ is inert in $E$ and, in the notations of \S\ref{start}, 
that $\pi_l$ is unramified. We claim that the central character $\chi$ of $\pi_l$ automatically satisfies $\chi(l)=(-1)^n$. 
Indeed, $\chi$ is trivial on $\anneau_E^\ast$ as $\pi_l$ is unramified, and it satisfies $\chi^{\bot}=\chi$. By Lemma \ref{lemmamu} and the
example above, $\chi_{|\Q_l^*}=1$ if, and only if, $n$ is even, hence the claim. 
\end{remark}

\begin{notation} \label{choixdemu}We choose a Hecke character $\mu$ of $E$ as follows : $\mu$ is a character like in  Lemma~\ref{lemmamu}(iii) if $n$ is odd, and 
$\mu=1$ if $n$ is even.
\end{notation}

We are now going to construct, place by place, a representation $\pi^n$ of $U(m)=U(n+2)$ whose conjectural
base change to $\Gl_m(E)$ has $L$-parameter 
$$L(\pi) \mu \oplus|\, |^{1/2}\mu \oplus |\, |^{-1/2} \mu.$$

\subsubsection{Construction of $\pi^n_l$, for $l$ split in $E$}\label{construpinsplit}

We denote by $P$ the upper parabolic subgroup of $\Gl_m(\Q_l)$ of type $(n,2)$, whose Levi subgroup is $M=\Gl_n(\Q_l)
\times \Gl_2(\Q_l)$. For $x$ a place of $E$ above $l$, we set
$$\pi^n_x:=\Ind_{P}^{\Gl_m(\Q_l)} (\pi_x (\mu_x \circ \det) \otimes (\mu_x \circ \det)).$$
Here $\Ind$ is the normalized induction. 
Since $\pi_x$ is unitary, $\pi^n_x$ is irreducible (see \cite{Bern2}). 

\begin{remark} \label{remlsplit}
Let $P'$ be the upper parabolic of type $(n,1,1)$.
Since $\pi_x$ is tempered by hypothesis, Langlands' classification theorem shows that 
$$\Ind_{P'}^{\Gl_m(\Q_l)} (\pi_x (\mu_x \circ \det) 
\otimes |\, |^{1/2} \mu_x 
\otimes |\, |^{-1/2} \mu_x),$$ 
has a unique irreducible quotient (that is, the Langlands quotient). As we have a natural $\GL_m(\Q_p)$-equivariant surjection from the representation above to the the irreducible representation $\pi^n_x$, this Langlands quotient is actually $\pi^n_x$. Thus, the $L$-parameter of $\pi^n_x$ is 
$L(\pi_x)\mu_x \oplus |\,|^{1/2}\mu_x \oplus |\,|^{-1/2}\mu_x.$
\end{remark}

Let us write $l=x\bar{x}$. By (i) of \S\ref{start}, $\pi_{\bar x}(\mu_{\bar x}\circ \det)$ is dual to $\pi_x (\mu_x \circ \det)$, 
so $\pi^n_{x}$ is dual to $\pi^n_{\bar x}$. The place $x$ defines, 
up to conjugation, an identification $$i_x : \U(m)(\Q_l) \rightarrow \Gl_m(\Q_l)$$ and so does the place $\bar x$, in such a way
that $i_{x}\circ i_x^{-1}$ is conjugate to $g \mapsto {}^t g^{-1}$. Hence we see that 
$i_x^\ast \pi_x^n \simeq i_{\bar x}^\ast \pi^n_{\bar x}$, 
using the well known result of Zelevinski that 
the representation $g \mapsto \tau({}^t g^{-1})$ is the contragredient 
of $\tau$ for any irreducible admissible representation $\tau$ of $\Gl_m(\Q_l)$.

We thus may set $\pi^n_l := i_x^\ast \pi^n_x$ and $\pi^n_l$ does not depend 
on the choice of the place $x$ above $l$.

\subsubsection{Construction of $\pi^n_l$, for $l$ inert or ramified in $E$}\label{construpinnonsplit}

We denote also by $l$ the place of $E$ above $l$. 
In this case $G:=\U(m)(\Q_l)$ is a quasi-split unitary group, 
and we shall use notations compatible to those of \S\ref{settingNMSRPS}. 
We may assume that $G$ is the unitary group 
defined by the following hermitian form on $E_l^m$ :
$$f(xe_i,ye_j)= c(x) y\delta_{j,m-i+1},$$ so that the group of diagonal
matrices in $G$, 
$$T=\{\diag(a_1,\dots,a_m),\ a_i \in E_l^\ast, a_{m-i+1}=c(a_i)^{-1}, 
i=1,\dots,m\},$$ 
is the centralizer of a maximal split torus in $\U(m)(\Q_l)$. Let $B$ be a the
upper triangular Borel.
The group $G':=\U(m)(E_l)$ is naturally identified with $\Gl_n(E_l)$ and 
$T'$ is the standard diagonal torus. Its Weyl group $W'$ is canonically 
identified with $\got{S}_{n+2}$.
The action of the non trivial element $c$ of 
$\Gal(E_l/\Q_l)$ on $T'$ is 
$$c (\diag(x_1,\dots,x_m)) = \diag(c(x_m)^{-1},\dots,c(x_1)^{-1}),$$ 
and $T$ is the subgroup of invariants of $c$ in $T'$.
 There is a norm map $\Norm : T' \rightarrow T$, 
$$x=(x_1,\dots,x_m) \mapsto xc(x)=(x_1 c(x_m)^{-1},x_2 c(x_{m-1})^{-1},\dots, 
x_m c(x_1)^{-1}).$$

By hypothesis (iii) of \S\ref{start}, and point (iii) of the remark therein,
$\pi_l$ is a subquotient of the normalized induction of a character 
$(\eta_1,\dots,\eta_n)$ of the standard torus of $\Gl_n$, with 
$\eta_i^\bot = \eta_{\sigma(i)}$ for all $i$ and some $\sigma \in {\got S}_{n}$. As $\pi_l$ is tempered,
 each $\eta_i$ is a unitary character.

We are going to define a character 
$\chi'=(\chi'_1,\dots,\chi'_m=\chi'_{n+2})$ of $T'$.
Up to reordering, the $\chi'_i$, $i=1,\dots,m=n+2$ are the $\eta_i \mu_l$, 
$i=1,\dots,n$ and $|\,|^{\pm 1/2} \mu_l$. The order is as follows :\begin{itemize} \ps

\item[-] First we define
$\chi'_1=|\,|^{-1/2} \mu_l$, $\chi'_m=|\,|^{1/2} \mu_l$. \ps

\item[-] Next, consider the set $I \subset \{1,\dots,n\}$ of $i$ 
such that $\eta_i \not \simeq \eta_i^\bot$. Clearly $|I|$ is even, say $2r$,
and we may define $\chi'_{2},\dots,\chi'_{r+1}$ and 
$\chi'_{m-r},\dots,\chi'_{m-1}$ in such a way that 
$\chi'_{m-j+1} \simeq {\chi'}_j^\bot$ for $j \in \{2,m-1\}$.
%, and also that if $\chi'_{k}=\chi'_{j}|\, |^{-1}$, for $j,k \in \{2,\dots,r+1,m-r\dots,m-1\}$, then $k > j$.
Finally, in case (iiib) we have $|I|=n$ if $n=2r$ is even 
(in which case we are done with the definition of $\chi'$) and $|I|=n-1$ if $n=2r+1$ is odd. In this 
case we have only left to define the ``midpoint'' character $\chi'_{r+1}$ 
for which we take (we have no other choice) $\eta_j \mu_l$, 
where $\eta_j \simeq \eta_j^\bot$ (this holds for a unique $j$).\ps

\item[-] In case (iiia), the characters $\eta_i$ for $i \not \in I$ satisfy $\eta_i =
\eta_i^\bot$, but since they are unramified, this implies $\eta_i(\varpi_l)= \pm 1$ 
(here $\varpi_l$ is a uniformizer of $E_l$). By the assumption on the central character 
of $\pi_l$, the set $\{i \not \in I,\, \, \eta_i(\varpi_l)=+1\}$ always has an even number of elements, say $2r'$.
For $r+1 \leq i \leq r+r'$, we set then $\chi'_i = \chi'_{m-i+1}=+\mu_l$ (with the obvious abuse of language), 
so for the remaining ones we have $\chi'_i=-\mu_l$.
\end{itemize}

\begin{lemma}\label{descentedescaracteres} The character $\chi'$ descends to $T$ i.e. there is a
smooth character $\chi$ of $T$ such that $\chi'=\chi \circ \Norm$.
Moreover, $\chi$ satisfies properties (a) and (b) of Definition~\ref{NMSRPS} in case (iiib). In case (iiia), 
$\chi$ is unramified if $m$ is even or if $l$ is inert in $E$.
\end{lemma}
\begin{pf}
By construction, in both cases, we have ${\chi'}_{m-i+1}^\bot = \chi'_i$ for all $i$. When $m=2q$ is even, we define 
$\chi(\diag(a_1,\dots,a_{2q}))=\chi'_1(a_1)\dots\chi'_q(a_q)$ and it is clear
 that $\chi \circ \Norm = \chi'$.
 
When $m=2q+1$, we remark that the middle character $\chi'_{q+1}$ of 
$E_{l}^\ast$ actually descends to a character $\chi_{q+1}$ of 
$U(1)(\Q_l)$. Indeed 
$$\chi'_{q+1} \prod_{i=2}^q \chi'_{i} {\chi'}_i^\bot = \det(\chi')$$ 
is the central character of $\pi_l \mu_l$. Since the central character of 
$\pi \mu$ has an integral weight (namely $\sum_{i=1}^n a_i + n/2$), it 
descends to $U(1)$ by Lemma~\ref{lemmamu}, and so does the central character
 of $\pi_l \mu_l$, hence also $\chi'_{q+1}$. \par
	Let $\psi$ be a smooth character of $U(1)(\Q_l)$ such that 
$\chi'_{q+1}(x)=\psi(x/c(x))$ for all $x \in E_l^*$. We define 
$\chi(\diag(a_1,\dots,a_{2q+1}))=\chi'_1(a_1)\dots\chi'_q(a_q) 
\psi(a_{q+1})$ and again it is obvious to see that 
$\chi \circ \Norm = \chi'$.

The other assertion is clear in case (iiib) as the $\eta_i$ are unitary, as well as in case (iiia) when $m$ is even. 
In the remaining case, the $\chi'_i$ are unramified for $i\neq q+1$ by choice of $\mu$ (i.e. Lemma \ref{lemmamu} (iii)), so 
we only have to check that $\psi$ is trivial. But 
$\chi'_{q+1}$ is trivial since it is unramified and satisfies $\chi'_{q+1}(l)=\mu_l(l)\eta_{q+1}(l)=+1$, 
hence the result follows from Hilbert 90.\end{pf} 

We now define $\pi^n_l$ as the unramified subquotient of $\Ind_B^G \chi$ 
in case (iiia) and as the unique subrepresentation $S(\chi)$ of 
$\Ind_B^G \chi$ in case (iiib) (see Def.~\ref{NMSRPS}). 

\begin{remark} \label{remlinert} 
The $L$-parameter of the conjectural base change of $\pi^n_l$ to 
$\Gl_m(E_l)$ is, by Remark~\ref{RemNMSRPS} in case (iiib) and 
by \cite[\S 3.6]{labesselivre} in case (iiia), the $L$-parameter attached to the character $\chi'$ of $T'$, which is by construction 
$$L(\pi_l) \mu_l \oplus |\, |^{1/2}\mu_l \oplus |\, |^{-1/2} \mu_l,$$ 
like in the split case. 
\end{remark}

\subsubsection{Construction of $\pi^s_\infty$}\label{construpisinfty}

Consider the morphism $\C^\ast \longrightarrow \Gl_m(\C)$ (recall that $m=n+2$)\begin{eqnarray*} z &\mapsto& \mu_\infty(z) 
\diag((z/\bar z)^{a_1},\dots,(z/\bar z)^{a_n},(z/\bar z)^{1/2},
(z/\bar z)^{-1/2}) \\  &=& \left\{ \begin{array}{c}
\diag((z/\bar z)^{a_1},\dots,(z/\bar z)^{a_n},(z/\bar z)^{1/2},
(z/\bar z)^{-1/2})\text{\ \ if $n$ is even}\\
\diag((z/\bar z)^{a_1+1/2},\dots,(z/\bar z)^{a_n+1/2},(z/\bar z),
1)\text{\ \ if $n$ is odd} \end{array}\right.\end{eqnarray*}

Since the $a_i$ are half-integers, and different from $\pm 1/2$,
we see by~\S\ref{repreal} that this morphism is always the 
restriction to $\C^\ast$ of the $L$-parameter 
of a unique irreducible representation $\pi^s_\infty$
of $\U(m)$. Here the $s$ stands for {\it square integrable}. 
The notation $\pi^n_\infty$ would be misleading since $\pi^s_\infty$ is, 
like any irreducible representation of a compact group, 
finite dimensional, square integrable, hence tempered.

\subsubsection{Assumption $AC(\pi)$}

\index{Zpin@$\pi^n$, an automorphic representation of $\U(m)$}
\index{AC@AC$(\pi)$, the assumption (a proof of which has been announced by Harris) that the representation $\pi^n$ 
attached to $\pi$ is automorphic}
\begin{conj} \label{AC} Assume that $\varepsilon(\pi,1/2)=-1$.
Then the irreducible admissible representation 
$$\pi^n := \pi^s_\infty \otimes {\bigotimes_{l}}{}' \pi^n_l$$
is automorphic. 
\end{conj}

The proof of this conjecture has recently been announced by Harris 
in the introduction of his preprint \cite{harrispreprint} (maybe under some local assumption). 
Since a written proof is not yet available, we prefer to be 
conservative and state it as a conjecture rather than as a theorem.
 
\begin{remark}\label{remsectAC}
\begin{itemize}
\item[(i)] The case $m=3$ (that is $n=1$) of this conjecture has been proved 
by Rogawski (\cite{Rog2}), using the Theta correspondence. 
In the case $m=4$, the needed local computations have
been published recently by Konno and Konno (\cite{konno}).

\item[(ii)] This conjecture is a very special case of the multiplicity formula 
of Arthur. Its derivation from those formulae is explained in great details in the following subsection. From that we shall see that the 
$\varepsilon(\pi,1/2)=-1$ should also be a 
necessary condition for the automorphy of $\pi^n$.

\item[(iii)] Although the construction of $\pi^n$ depends on the choice 
of the Hecke character $\mu$ (for odd $n$ : see Notation~\ref{choixdemu}), 
it is clear that the conjecture is independent of this choice. 
Indeed, if $\mu$ is changed into another character $\mu_1$, then 
$\mu_1=\mu \phi'$ where $\phi'$ is a
Hecke character of $\A_E^\ast$ that descends to a character $\phi$ of $\U(1)$.
By construction
the representation ${\pi_1^n}$ defined using $\mu_1$ is simply 
$\pi^n (\phi \circ \det)$ and it follows that the automorphy of $\pi^n$ is 
equivalent to the automorphy of ${\pi_1^n}$. 

Note also that the hypothesis in the conjecture is about
$\varepsilon(\pi,1/2)$, not about $\varepsilon(\pi (\mu\circ\det),1/2)$.
\end{itemize}
\end{remark}

\newpage
\renewcommand{\iff}{if, and only if,}
\renewcommand{\Sc}{\mathcal{S}}
\renewcommand{\Cc}{\mathcal{C}}
\newcommand{\CY}{\cal Y}  

\section{Eigenvarieties of definite unitary groups}

\subsection{Introduction}	In this section, we introduce and study in details the eigenvarieties of definite unitary groups and we prove the basic 
properties of the (sometimes conjectural) family of Galois representations that they carry. It furnishes a lot of 
interesting examples where all the concepts studied in this book occur, and provides also an important tool for the applications to Selmer groups in the next sections. 
As a first application, we define some purely Galois 
theoretic global deformation rings and discuss their relations to those 
eigenvarieties at some specific classical points (including $R=T$ like statements). We give a second application to 
the construction of many irreducible Galois representations. We prove also quite a number of results of independent interests 
regarding the theory of eigenvarieties that we explain in details below. The organization of this section is as follows. \ps
	In the first subsection \ref{dabpote}, we give an axiomatic definition of eigenvarieties and draw the general 
consequences of our definitions. In particular, we show that an eigenvariety is unique 
(up to unique isomorphism) if it exists (Prop. \ref{eigenunique}). One interest is that there are in principle many 
different ways two construct eigenvarieties: using coherent or Betti cohomology, a group or its inner forms (or any 
transfer suggested by Langland's philosophy), using Emerton's representation theoretic approach  etc... 
Each of those constructions has its own advantages but they sometimes should all lead to the same eigenvariety. The uniqueness statement 
alluded above will often show that they are indeed the same\footnote{It may be also useful to combine it with a Jacquet-Langlands correspondence.}. A second interest is that there is some abstract game that we 
can play to deduce the existence of a eigenvarieties form the existence of other ones. \par	
	The context is the following:\footnote{The choice of a quadratic imaginary field in this chapter rather than a general CM field 
(as well as the choice of a split $p$) is made mainly to simplify the exposition and also because 
this the only case that we shall use in the applications to Selmer groups. All the constructions actually extend to this more general 
setting by combining the arguments here (or of \cite{Ch}) and those of \cite{Buz} 
(see also Yamagami's work \cite{Yam}). Alternatively, the general {\it definite} case is now covered by Emerton's paper \cite{eminterp}. The main reason why we fix a 
split prime $p$ is Galois theoretic: at the moment, Kisin's arguments \cite{kis} and the theory of trianguline representations 
are only written in the case where the base field is $\Q_p$ rather that any finite extension of $\Q_p$.} $E/\Q$ is a quadratic imaginary field 
and $G/\Q$ is a unitary group in $m\geq 1$ variables attached to $E/\Q$. We assume that $G(\R)$ is the compact real unitary group ($G$ is {\it definite}) and we fix a prime $p$ such 
that $G(\Q_p) \simeq \GL_m(\Q_p)$, as well as embeddings $\overline{\Q} \rightarrow \Qpb$ and $\overline{\Q} \rightarrow \C$. 
An (irreducible) automorphic representation $\pi=\pi_{\infty}\otimes \pi_f$ of $G$ is automatically 
algebraic and has cohomology in degree $0$. The finite dimensional representation $\pi_{\infty}$ is determined by its {\it weight} 
which is a decreasing sequence of integers $\uk=(k_1\geq \cdots \geq k_m)$, and $\pi_f$ is defined over $\Qb$ hence 
may be viewed over $\Qpb$ via the chosen embeddings. We fix also a commutative Hecke-algebra $$\HH = \ATL \otimes \HH_{\rm ur}$$ that 
contains the Atkin-Lehner algebra $\ATL$ of $U$-operators at $p$ and a spherical part $\HH_{\rm ur}$. We are interested in $p$-adically 
interpolating the systems of Hecke eigenvalues $\psi_{\pi}: \HH \longrightarrow \Qpb$ cut out from the $\pi$ as above, and more precisely the pairs 
$(\psi_{\pi},\uk)$ where $\uk$ is the weight of $\pi$. Note that the systems of eigenvalues of $\ATL$ on the Iwahori invariants of $\pi_p$ 
(say if $\pi_p$ is unramified) are in bijection with the {\it refinements} of $\pi_p$ in the sense of \S\ref{sectionref}, so that $\psi_{\pi}$ contains the extra datum 
of a choice of refinement $\RR$ of $\pi_p$. To remind this refinement, we actually denote $\psi_{\pi}$ by 
$\psi_{(\pi,\RR)}$. Let us fix now a collection 
$$\cal Z \subset \Hom_{\rm ring}(\HH,\Qpb) \times \Z^m$$ 
of such $(\psi_{(\pi,\RR)},\uk)$. An {\it eigenvariety for $\cal Z$} is a 
$4$-uple $(X,\psi,\omega,Z)$ where \begin{itemize} 
\item[(a)] $X$ is a reduced rigid analytic space over $\Q_p$, 
\item[(b)] $\psi: \HH \longrightarrow \OO(X)$ is a ring homomorphism, 
\item[(c)] $\omega: X \longrightarrow \WW:=\Hom((\Z_p^*)^m,{\mathbb G}_m)$ is an analytic map to the weight space,
\item[(d)] $Z \subset X(\Qpb)$ is a Zariski-dense subset,
\end{itemize}
such that the evaluation of $\psi$ induces an injection $$X(\Qpb) \hookrightarrow \Hom(\HH,\Qpb) \times \WW(\Qpb)$$ which 
itself induces a bijection\footnote{This makes sense as $\WW(\Qp)$ naturally contains $\Z^m$, see \S\ref{sectdefeigen}.} 
$Z \isomo \cal{Z}$. To have the uniqueness property we need of course to impose some extra conditions on $(\psi,\omega,Z)$ 
for which we refer to Definition \ref{eigendef}. We show that for an eigenvariety $X$, the unit ball $\OO(X)^0$ is a compact 
subset of $\OO(X)$, which (together with (d)) is the basic property needed for the construction of Galois pseudocharacter on $X$. \ps
	In the second subsection \S\ref{secteigenvidempo} we recall the results of one of us on the existence of eigenvarieties (\cite{Ch}).
The statement is that for any idempotent in the Hecke-algebra $\Cc_c(G(\AAA^p),\Qb)$ commuting with $\HH$, there is an 
eigenvariety for the set $\cal Z_e$ parameterizing all the $p$-refined $\pi$ such that $e(\pi) \neq 0$. We discuss in 
Example \ref{exampleidemp} what kind of $\cal Z$ we can reach this way in representation theoretic terms 
(Bernstein components, type theory). In fact, those eigenvarieties {\it of idempotent type} have 
stronger properties than the general ones.  As their structure plays a crucial role in our main theorem on Selmer groups, as well as in some subsequent constructions in this section, we found it necessary to review 
their construction in details, 
as well as the theory of $p$-adic automorphic forms of $G$ developed in \cite{Ch}. This is the aim of \S\ref{eigenreview} to \S\ref{padIII}.
In fact, this part is essentially self-contained and slightly improves some results there ({\it e.g.} we do not restrict to the center of the weight space, or to a neat level, 
and we release the assumption that $p=2$ at some point). 
We rely on the work of Buzzard on eigenvarieties \cite{Buz}. 
Let us also mention here that if we had been only interested in the subset ${\cal Z}_{e, {\rm ord}} \subset {\cal Z}_e$
of $\pi$ which are $p$-ordinary, the existence of $X$ is due to Hida (actually in a much wider context \cite{hida}). Moreover, there is an alternative construction 
of $X$ due to Emerton in \cite{eminterp}.\ps
	In a third subsection \ref{bernsteinnonmonodrmic}, we show how to define some quasicoherent sheaves 
of admissible $G(\AAA_f^p)$-representations on an eigenvariety of idempotent type, and we prove their basic properties. As an application,
we show the existence of an eigenvariety for the subset $\cal Z_{e,NM}$ parameterizing the $p$-refined $\pi$ in $\cal Z_e$
such that $\pi_l$ is an NMSRPS representation (in the sense of \S\ref{settingNMSRPS}) for each $l$ in a fixed finite set of places such that $G(\Q_l)$ is a 
quasisplit unitary group. We don't know if those latter eigenvarieties are of idempotent types. As a consequence of all those constructions, 
we introduce in Def. \ref{minimaleigendef} the convenient notion of {\it minimal eigenvariety} containing 
a given $p$-refined automorphic representation, it is defined at the moment only when $\pi_l$ is either unramified or NMSRPS at each nonsplit 
prime $l$.\ps
	In the fourth part \S\ref{famillegaloisiennesurX}, we explain how the existence of the
expected $p$-adic Galois representations associated to (sufficiently many of) the $\pi$ parameterized by $\cal Z$ gives rise to a 
continuous Galois pseudocharacter $$T: \G_E \longrightarrow \OO(X).$$
Our results in this part are unconditional when $G$ is attached to certain division algebras, and conditional to the assumption 
(Rep(m)) when $G$ is the definite unitary group $\U(m)$ which is quasisplit at each finite place (so $m \neq 2 \bmod 4$). For each $x \in X(\Qpb)$, 
we have then a canonical semisimple representations $$\rhob_x: \G_E \longrightarrow \GL_m(\Qpb),$$ 
whose trace is the evaluation of $T$ at $x$. The game is to understand with this weak notion of families of Galois representations (namely the simple existence of $T$) how to deduce from a
property of the $\rhob_z$ for (a Zariski-dense subset of) $z \in Z$, a similar property for $\rhob_x$ 
for {\it any} $x \in X(\Qpb)$. We are typically interested in a property concerning the restriction to a decomposition group 
at a finite place $w$ of $E$.\par
	At a prime $w$ above $p$, we show that $T$ is a refined family in the sense of \S\ref{sectdefrigfam} hence we can apply to $(X,T)$ the results 
of section 4. At a prime $w$ not dividing $p$, it is convenient to introduced the {\it generic representation}
$$\rhog_x: \G_E \longrightarrow \GL_m(\Kb_x),$$ whose trace is the composition of 
$T: \G_E \longrightarrow \OO(X) \longrightarrow \K_x:={\rm Frac}(\OO_x)$ and where $\Kb_x$ is a product of algebraic closures of each factor field of $\K_x$ 
(i.e. of the fraction fields of the germs of irreducible components of $X$ at $x$). 
The representations $\rhob_x$ and $\rhog_x$ have an associated Weil-Deligne representation that we 
compare, and that we also compare with the ones of the $\rhob_z$ for $z \in Z$. For example, when $X$ is the minimal 
eigenvariety containing some $\pi$ of the type of Harris-Taylor, these Weil-Deligne representations are constant on $X$ when 
restricted to the inertia group. The proofs rely on some useful lemmas on nilpotent elements in general matrix rings or GMA. 
Those facts are proved separately in the first subsection of an appendix \S\ref{appendixladic} that we devote to the general study of 
$p$-adic families of Galois representations of $\Gal(\Qb_l/F)$ when $l\neq p$ and $F/\Q_l$ a finite extension. In this appendix, we also 
introduce the {\it dominance ordering} $\prec$ on nilpotent matrices and on Weil-Deligne representations, which is convenient to state our results. \par

	In the next subsection \S\ref{applidef} we give an application of the techniques and results of this book to study some global Galois 
deformation rings, as was announced in \S\ref{globcons} of section 2. We fix a continuous absolutely irreducible Galois representation
$$\rho: \G_{E,S} \longrightarrow \GL_m(L)$$
($L$ a finite extension of $\Q_p$) such that $\rho^\bot \simeq \rho(m-1)$, and which is crystalline with two-by-two distinct 
Hodge-Tate weights and crystalline Frobenius eigenvalues at the primes above $p$. \par
	We are interested in the deformations $\rho_A$ 
of $\rho$ such that $\rho_A^{\bot} \simeq \rho_A(m-1)$, where $A$ is a local artinian rings $A$ 
with residue field $L$. We introduce two subfunctors of the full deformation functor: the {\it fine deformation functor} $\XX_{\rho,f}$,
whose tangent space is $H^1_f(E,\ad(\rho))$, and the {\it refined deformation functor} $\XX_{\rho,\Ref}$, which depends on the choice
of a refinement $\Ref$ of $\rho_{|E_v}$ ($p=v\bar{v}$). We show that those functors are pro-representable, 
we compare them when $\Ref$ is non critical, and we formulate two basic conjectures (C1) and (C2) concerning their structure 
(see Conj. \ref{conjc1c2}). \par
	For a quite general case of $\rho$, we also introduce in 
\S\ref{specialcaserhomod} a definite unitary eigenvariety $X$ and a point $z \in X$. If 
$\mathbb T$ (resp. $R_{\rho,\Ref}$) denotes the completion of the local ring 
of $X$ at $z$ (resp. the ring pro-representing $\XX_{\rho,\Ref}$), we show the existence of a natural map 
$$R_{\rho,\Ref} \longrightarrow \mathbb T.$$
In this context, this arrow and the properties of eigenvarieties allow us to show that our conjectures (C1) and (C2) are actually consequences of 
the Bloch-Kato conjecture\footnote{More 
precisely, they are equivalent to the vanishing of $H^1_f(E,\ad(\rho))$.}, which provides a strong evidence for them. 
This leads us to conjecture that the arrow above is an isomorphism ("$R_{\rho,\Ref}=\mathbb T$"), and that 
a strong infinitesimal version of the {\it non critical slope form are classical} should hold: 
"eigenvarieties should be \'etale over the weight space (hence smooth) at non critical irreducible classical points". 
In turn, these last two conjectures imply (C1) and (C2), and we end the paragraph by a series of remarks concerning them.\par
Finally, as a simple application of the theory of refined families, we show in \S\ref{parappliirr} how we can construct many $m$-dimensional 
Galois representations of $\G_E$ which are unramified outside $p$ and crystalline, irreducible, and with generic Hodge-Tate weights at the two primes of $E$ 
dividing $p$. This application is conditional to (Rep(m))) but does not use any irreducibility assertion for the automorphic Galois representations. 
We rather 
start from the trivial representation and move in the {\it tame level 1} eigenvariety to find the Galois representations we are looking for. \ps

\subsection{Definition and basic properties of the Eigenvarieties}\label{dabpote}

\subsubsection{The setting}\label{basicsettingeigenv} Let $E$ be a 
quadratic imaginary field and $G$ be a definite unitary group in $m\geq 1$ 
variables attached to $E/\Q$, as in \S\ref{unitarygroups}, {\it e.g.} the group $\U(m)$ defined in \S\ref{definitionUm}.
Let us fix once and for all a rational prime $p$ as well as fields 
embeddings\footnote{Which means that we also fix once for all an algebraic closure $\Qb$ (resp. $\Qpb$) of $\Q$ (resp. $\Qp$).} $$\iota_p : \Qb \longrightarrow \Qpb, \, \, \, \, \, 
\iota_{\infty}: \Qb \longrightarrow \C.$$ 
We assume that $G(\Q_p)\simeq \GL_m(\Q_p)$. In particular, $p$ splits in $E$ and if we 
write $p=vv^c$ where $v: E \rightarrow \Q_p$ is defined by $\iota_p$, then $v$ induces a natural isomorphism 
$G(\Q_p) \isomo \GL_m(\Q_p)$. The embedding $E \longrightarrow \C$ given by $\iota_{\infty}$ induces an embedding 
$G(\R) \hookrightarrow_{\iota_{\infty}} \GL_m(\C)$ well defined up to conjugation. \ps
We fix a model of $G$ over 
$\Z$ and an associated product Haar measure $\mu$ on $G(\AAA_f)$. We use the standard conventions for ad\`eles: $\AAA_S$ (resp. $\AAA_f^S$) denotes 
the ring of finites ad\`eles with components in (resp. outside) the set of primes $S$. Moreover, we denote by $\widehat{\Z}_S=\prod_{l \in S} \Z_l$ the ring of integers of $\AAA_S$.\ps
        The definition of an eigenvariety for $G$ depends on the choice of a commutative Hecke algebra $\HH$ that we fix once for all as follows. We fix a 
subset $S_0$ of the primes $l$ split in $E$ such that $G(\Q_l)\simeq \GL_m(\Q_l)$ and $G(\Z_l)$ is a maximal compact subgroup\footnote{In most applications, $S_0$ will have Dirichlet 
density one.}, and set
$$\HH_{\rm ur}:=\Cc(G(\widehat{\Z}_{S_0})\backslash G(\AAA_{S_0})/G(\widehat{\Z}_{S_0}),\Z).$$
Recall that we defined in \S\ref{atkinlehner} a subring\footnote{Let us assume that $\mu(I)=1$.} 
$\ATL \subset \Cc(I\backslash G(\Q_p)/I,\Z[1/p])$ where $I \subset \G(\Q_p)\isomo_v \GL_m(\Q_p)$ is the standard Iwahori subgroup. We set
\footnote{We limit ourselves a bit the choice of $\HH$ here only for notational reasons and later use. 
We could add for example inside $\HH$ any commutative subring of 
$\Cc(K\backslash G(\AAA^{S_0 \cup\{p\}})/K,\Qpb)$ for some compact open subgroup $K$ and everything would apply verbatim.} $$\HH:=\ATL \otimes \HH_{\rm ur}.$$
\index{H@$\HH=\ATL \otimes \HH_{\rm ur}$, an Hecke algebra}

\subsubsection{$p$-refined automorphic representations.}\label{defprefautform}

\begin{definition}\label{prefaut} We say that $(\pi,\RR)$ is a $p$-refined automorphic representation of weight $\uk$ if: \begin{itemize}
\item[-] $\pi$ is an irreducible automorphic representation of $G$ (see \S\ref{defautrepG}),
\item[-] $\pi_{\infty} \isomo_{\iota_{\infty}} W_{\uk}(\C)$ (see \S\ref{repreal}),
\item[-] $\pi_p$ is unramified and $\RR$ is an accessible refinement of $\pi_p$ (see \S\ref{refinements}).
\end{itemize}
\end{definition}
\par \medskip

\begin{definition}\label{heckeeigenv} If $A$ is a ring, an {\it $A$-valued system of Hecke eigenvalues} is a ring homomorphism $\HH \longrightarrow A$.
\end{definition}
Let $(\pi,\RR)$ be as above, we can attach to it a $\Qpb$-valued system of 
eigenvalues $$\psi_{(\pi,\RR)}: \HH \longrightarrow \Qpb$$ 
as follows. By Definition \ref{defrefi2}, if $\chi: U \rightarrow \C^*$ is the character of the refinement $\RR$, $\chi\delta_B^{-1/2}$ 
occurs in $\pi_p^I$, and by 
\S\ref{repreal} there is also a character $\delta_{\uk}: U \rightarrow \C^*$ associated to $\uk$, so there is a unique 
ring homomorphism $\psi_p: \ATL \rightarrow \C$ such that
\begin{equation}\label{combiref} {\psi_p}_{|U}=\chi\delta_B^{-1/2}\delta_{\uk}.\end{equation}
Moreover, $\pi^{G(\widehat{\Z}_{S_0})}$ is one dimensional
hence it defines a ring homomorphism $\psi_{\rm ur}: \HH_{\rm ur} \rightarrow \C$. By Lemma \ref{pifdefqb}, 
the complex system of Hecke-eigenvalues $\psi_p \otimes_{\Z} \psi_{\rm ur}$ is actually $\Qb$-valued. 

\index{Zpsi@$\psi_{(\pi,\RR)}:\HH \rightarrow \bar \Q_p$, the system of Hecke eigenvalues attached to an automorphic 
representation $\pi$ with refinement $\RR$}  
\begin{definition}\label{syshecke} We call $\psi_{(\pi,\RR)}$ the $\Qpb$-valued system of Hecke eigenvalues associated to the 
$p$-refined automorphic representation $(\pi,\RR)$ of weight $\uk$ defined by $\iota_p\iota_{\infty}^{-1}(\psi_p \otimes \psi_{\rm ur})$. 
\end{definition}

\begin{remark} \label{remsyshecke} We have $\psi_{(\pi,\RR)}=\psi_{(\pi',\RR')}$ \iff\, $\pi_t \simeq \pi'_t$ for each $t \in
S_0\cup \{p\}$ and $\RR=\RR'$. 
\end{remark}

\subsubsection{Eigenvarieties as interpolations spaces of $p$-refined automorphic representations}\label{sectdefeigen} Let $\cal{Z}_0 \subset \Hom_{\rm ring}(\HH,\Qpb) \times \Z^m$ be the set of pairs 
$(\psi_{(\pi,\RR)},\uk)$ associated to all the $p$-refined 
automorphic representations $\pi$ of any weight $\uk$, and let us fix $\cal{Z} \subset \cal{Z}_0$ a subset. It will actually be convenient to give here a formal definition of what is an eigenvariety attached to $\cal{Z}$. 
We shall actually never use here the group $G$ and the set $\cal Z$ could be replaced by any subset of $\Hom_{\rm ring}(\HH,\Qpb) \times \Z^m$. \ps

The {\it weight space} is the rigid analytic space over $\Q_p$
\footnote{It is isomorphic to a finite disjoint union of unit open
$m$-dimensional balls.}
$$\WW:=\Hom_{gr-cont}(T^0,\mathbb{G}_m^{\rig})$$ 
whose points over any complete extension $k$ of $\Q_p$ parameterize the continuous, $k$-valued, characters of $T^0=(\Z_p^*)^m$. We view $\Z^m$ as embedded inside $\WW(\Q_p)$, by
mean of the map $$(k_1,\cdots,k_m) \mapsto ((x_1,\cdots,x_m) \mapsto
x_1^{k_1}\cdots x_m^{k_m}), \, \, \, \, \Z^m \hookrightarrow \WW(\Qp),$$ 
and we denote by $\Z^{m,-}$ the subset of $\Z^m$
consisting of strictly decreasing sequences. We denote also by $u_p$ the element 
$$u_p:=\diag(p^{m-1},\dots,p,1) \in U^- \subset \ATL^*.$$

Let us fix $L \subset \Qpb$ a finite extension of $\Qp$. In the definition below and in the sequel, we will always view $\WW$, $\AAA$ and ${\mathbb G}_m$ as rigid analytic spaces over $L$ 
even if we do not make it appear explicitly: for example, we will write $\WW$ for $\WW \times_{\Q_p} L$.

\begin{definition}\label{eigendef} An eigenvariety for $\cal{Z}$ is a reduced $p$-adic analytic space $X$ over $L$ equipped with: \begin{itemize}
\item[-] A ring homomorphism $\psi: \HH \longrightarrow \OO(X)^{\rig}$,\ps
\item[-] An analytic map $\omega: X \longrightarrow \WW$ over $L$, \ps
\item[-] An accumulation and Zariski-dense subset $Z \subset X(\Qpb)$, 
\end{itemize}
\par \smallskip \noindent such that the following conditions are satisfied:
:\begin{itemize}
\item[(i)]  
For all $x \in X$, the natural map $\HH \otimes_{\Z}
\OO_{\omega(x)} \longrightarrow \OO_x$ is surjective. \ps
\item[(ii)] The map $\nu:=(\omega,\psi(u_p)^{-1}): X \longrightarrow \WW \times \mathbb{G}_m$ is
finite.\ps
\item[(iii)] The natural evaluation map $X(\Qpb) \rightarrow \Hom_{\rm ring}(\HH,\Qpb)$, $$x \mapsto \psi_x:=(h\mapsto \psi(h)(x)),$$ induces a bijection $Z \isomo \cal{Z}$, $z \mapsto (\psi_z,\omega(z))$. 
\end{itemize}
\end{definition}
\begin{remark}\label{remoneigendef} \begin{itemize}
\item[i)] In other words, $(\psi,X)$ is a rigid analytic family of systems of Hecke-eigenvalues interpolating the ones in $\cal{Z}\isomo Z$. The Zariski-density of $Z$ and (i) ensures that $X$ is somewhat minimal with that property in some sense (see Prop. \ref{eigenunique}). 
\item[ii)] In some sense, we may like to think of (or define) such an eigenvariety as the "Zariski-closure" of $\cal{Z}$ in $\WG \times \Hom_{\rm ring}(\HH,\AAA)$. However, as this latter space is not a rigid space if $\HH$ is not finitely generated (which will be the case in the applications), we have to be a little careful. The requirement (ii) above is a way to circumvent this problem.
\item[iii)] The notation $X(\Qpb)$ is a shortcut for the union of the $X(K)$ for all $K \subset \Qpb$ 
finite over $\Qp$. Moreover, there is a slight abuse of language in the definition above: if $Z \subset X(\Qpb)$ is a subset, we say that $Z$ is Zariski-dense (resp. {\it accumulation}) if the underlying set of closed points $|Z|$ is (resp. if $|Z|$ accumulates at each point of itself, in the sense of 
\S\ref{accumulate}). 
\end{itemize}
\end{remark}
\ps
It turns out that such an eigenvariety, if exists, is unique.

\begin{prop}\label{eigenunique} If $(X_1/L,\psi_1,\nu_1,Z_1)$ and $(X_2/L,\psi_2,\nu_2,Z_2)$ are two eigenvarieties for $\cal{Z}$, there exists a unique $L$-isomorphism $\zeta: X_1 \rightarrow X_2$ such that $\nu_2 \cdot \zeta = \nu_1$, and $\forall h \in \HH$, $\psi_1(h)=\psi_2(h) \cdot \zeta \in \OO(X_1)$.
\end{prop} 

\begin{pf}
Fix $(X_i,\psi_i,\omega_i,Z_i)$,
$i=1, 2$ satisfying (i) to (iii), and denote again by $\nu_i=(\omega_i \times
\psi_i(u_p)): X_i \rightarrow \WW \times \mathbb{G}_m$ the finite map of (iii). As a consequence of Lemma~\ref{eigenproperties} (c) thereafter and assumption (iii), there is a unique bijection $\zeta: Z_1 \isomo Z_2$, such that
for all $z \in Z_1$, ${\psi_1}_z={\psi_2}_{\zeta(z)}$ and
$\omega_1(z)=\omega_2(\zeta(z))$. We will eventually prove that $\zeta$
extends to an isomorphism $\zeta: X_1 \isomo X_2$ as in the statement.
\noindent By Lemma\ref{eigenproperties} (a), such a map is actually unique if it exists. \ps
For any admissible open $V \subset \WG$, we set 
$X_{i,V}:=\nu_i^{-1}(V)$, and let $\AAA_V$ denote the affine line over $V$. 
For each finite set $I \subset \HH$ and such a $V$, we have a natural 
$V$-map $$f_{i,V,I}: X_{i,V} \longrightarrow \AAA_{V}^I,\, \, \, x\mapsto (h(x))_{h \in
I}$$
inducing a natural map $X_{i,V} \rightarrow \projlim_{I \subset \HH} \AAA_V^I$, and commuting with any base change by an 
open immersion $V' \subset V$. The morphism $f_{i,V,I}$ is closed by (iii). \ps
Assume $V$ is moreover affinoid. Lemma \ref{eigenproperties} (a) shows that there exists $I_V$ 
such that for $I \supset I_V$, $f_{i,V,I}$ is a closed immersion for both $i$. We claim that for $I \supset I_V$, we have an inclusion 
$f_{1,V,I}(X_{1,V}) \subset f_{2,V,I}(X_{2,V})$. 
By exchanging $1$ and $2$ and using that both $X_{i,V}$ are reduced, it will follow that as closed subspaces of $\AAA_V^I$, for $I \supset I_V$, we have $f_{1,V,I}(X_{1,V})=f_{2,V,I}(X_{2,V})$. \ps
Let $x \in X_{1,V}$. If $x \in Z_1$ then $f_{1,V,I}(x) \in f_{2,V,I}(X_{2,V})$ 
by definition of $\zeta$. In general, by the Zariski density of $Z_1$ in $X_1$ and 
Lemma \ref{eigenproperties} (b), we can find an open affinoid $V' \supset V$ 
such that some $z \in Z_1 \cap X_{1,V'}$ lies in the same 
irreducible component $T$ of $X_{1,V'}$ as $x$. 
By the accumulation property of $Z$, $Z$ is Zariski-dense in $T$, hence for 
$I' \supset I \cup I_{V'}$, $$f_{1,V',I'}(T) \subset f_{2,V',I'}(X_{2,V'}).$$ In particular, for such an $I'$ 
we have $f_{1,V,I'}(x) \in f_{2,V,I'}(X_{2,V})$ and by projecting to $\AAA_V^I$ we get that this holds also 
when $I'=I$, hence the claim. \ps We define now a $V$-isomorphism
$\zeta_V:  X_{1,V} \longrightarrow X_{2,V}$ by setting, for $I \supset I_V$,  
$$\zeta_V:=f_{2,V,I}^{-1}\cdot f_{1,V,I}.$$ 
This map does not depend on $I$ and it obviously extends the previously defined map $\zeta$ on $Z_1 \cap X_{1,V}$. The independence of $I$ implies that $\forall h \in \HH,\, \, \psi_1(h)=\psi_2(h) \cdot \zeta \in \OO(X_{1,V})$. We check at once that $\zeta_V \times_V V'=\zeta_{V'}$ for any 
$V' \subset V$ open affinoid, hence the 
$\zeta_V$ glue to a unique isomorphism $\zeta: X_1 \longrightarrow X_2$ and we are done. \end{pf}

\begin{lemma}\label{eigenproperties} Let $(X,\psi,\omega,Z)$ be an eigenvariety:\begin{itemize}
\item[(a)] For any open affinoid $V \subset \WG$, $\nu^{-1}(V)$ is affinoid and 
the natural map $\HH\otimes \OO(V) \longrightarrow \OO(\nu^{-1}(V))$ is surjective. \ps
\item[(b)] Any two closed points of $X$ lie in a same $\nu^{-1}(V)$ for some $V \subset \WG$ open affinoid.
\item[(c)] For any $x, y \in X(\Qpb)$, $x=y$ \iff\,  $\psi_x=\psi_y {\rm \, \, and}\, \,
\omega(x)=\omega(y).$
\end{itemize}
\end{lemma}

\begin{pf} The first part of (a) follows from the finiteness of $\nu$, and 
the second part form (i) and the faithful flatness of $\OO(\Omega) \longrightarrow
\oplus_{x \in \Omega} \OO_{x}$ for any affinoid $\Omega$. 
Assertion (b) follows from the fact that any two closed points of $\WG$ lie in a same affinoid subdomain.
Part (c) is then a consequence of (a) and (b).
\end{pf}

\begin{definition} \label{defnested} We say that a rigid space $X$ over $\Q_p$ is nested if it has an
 admissible covering by some open affinoids $\{X_i,i \geq 0\}$ such that $X_i \subset X_{i+1}$ and that the natural 
$\Q_p$-linear map $\OO(X_{i+1}) \longrightarrow \OO(X_i)$ is compact. 
\end{definition}

Note that any such $X$ is separated, and that any finite product
of nested spaces is nested. For example, $\AAA$, $\mathbb{G}_m$
and $\WW$ are easily checked to be nested, hence so is $\WG$. 

\begin{lemma}\label{lemmanested} Assume that $X/\Q_p$ is nested. \begin{itemize}
\item[(i)] If $Y \longrightarrow X$ is a finite morphism, then $Y$ is nested.
\item[(ii)] Assume that $X$ is reduced. Then $$\OO(X)^0:=\{f\in \OO(X),\ \forall x \in X,\, |f(x)|\leq 1 \}$$ is a compact subset of $\OO(X)$.
\end{itemize}
\end{lemma}
 
Recall that $\OO(X)$ is equipped the coarsest locally convex topology such that all the restriction maps $\OO(X) \longrightarrow \OO(U)$, $U \subset X$ an affinoid subdomain, are continuous ($\OO(U)$ being equipped with its Banach algebra topology). It is a separated topological $\Q_p$-algebra.

\begin{pf} To show (i), it suffices to check that for each $\Q_p$-affinoid $X$, each coherent $\OO_X$-module
$M$, and each affinoid subdomain $U \subset X$ such that $\OO(V) \longrightarrow \OO(U)$ is
compact, then the natural map $M(X) \longrightarrow M(U)$ is compact. Of course, $M(U)$ and $M(X)$ are equipped here with their (canonical) topology of finite module over an affinoid algebra. Let us fix an $\OO$-epimorphism $\OO^n \longrightarrow M$, and
consider the commutative diagram of continuous $\Q_p$-linear maps

$$\xymatrix{ \OO(X)^n \ar@{->}[d] \ar@{->>}[rr]
& & M(X) \ar@{->}[d] \\
 \OO(U)^n \ar@{->>}[rr] & & M(U) }$$
Let $\cal{B} \subset M(X)$ be a bounded subset. By the open mapping theorem, there is a bounded subset $\cal{B}' \subset \OO(X)^n$ whose image under the top 
surjection is $\cal{B}$. By assumption the left vertical arrow is compact hence the image of $\cal{B}'$ in $\OO(U)^n$ has compact closure, hence so has the image of $\cal{B}$ in $M(U)$.\ps
	We show (ii) now. Let us fix $X=\cup_i X_i$ a nested covering of $X$, and set $$Y_i:=\overline{{\rm Im}(\OO(X)^0 \rightarrow \OO(X_i))}.$$ It is a compact subspace of $\OO(X_i)$ as is the image of the unit ball of $\OO(X_{i+1})$ by assumption. But the injection $$\OO(X)^0 \longrightarrow \prod_i \OO(X_i)$$ 
has a closed image, and is a homeomorphism onto its image. We conclude as it lies in the compact subspace $\prod_i Y_i$.	
\end{pf}

\begin{cor}\label{eigennested} Eigenvarieties are nested. \end{cor}

For later use, let us introduce another notation. Let $(X/L,\psi,\nu,Z)$ be an eigenvariety. Let $u_i=(z_1,\cdots,z_m) \in U$ be the element such that $z_i=p$ and $z_j=1$ if $j\neq i$.

\begin{definition}\label{definitionFi} For $i\in \{1,\cdots,m\}$, $F_i:=\psi(u_i)\in \OO(X)^*$. By definition and by formula 
(\ref{combiref}), they are the unique analytic functions $X$ such that for each $z=(\psi_{(\pi,\RR)},\uk) \in Z \isomo {\cal Z}$, we have 
$$\iota_p\iota_{\infty}^{-1}(\RR \cdot |p|^{\frac{1-m}{2}})=(F_1(z)p^{-k_1},\cdots, F_i(z)p^{-k_i+i-1},\cdots,F_m(z)p^{-k_m+m-1}).$$
\end{definition}

\subsection{Eigenvarieties attached to an idempotent of the Hecke-algebra}\label{secteigenvidempo}

\subsubsection{Eigenvarieties of idempotent type}\label{eigenidem}

We keep the notations above, and we fix an idempotent 
$$e \in \Cc_c(G(\AAA_f^{p,S_0}),\Qb)\otimes {1_{\HH_{\rm ur}}} \subset \Cc_c(G(\AAA_f^p),\Qb).$$ 
Let $\cal{Z}_e \subset \cal{Z}$ be the subset of  $(\psi_{(\pi,\RR)},\uk)$ such that $e(\pi^{p}) \neq 0$. We will say that such a $\pi$ is {\it of type $e$}. Assume that $\cal{Z}_e$ is nonempty. \ps
We fix $L \subset \Qpb$ a sufficiently big finite extension of $\Q_p$ such that $\iota_p(e) \in \Cc_c(G(\AAA_f^p),L)$. Moreover, we will write also $e$ instead of $\iota_p(e)$ or for $\iota_{\infty}(e) \in \Cc_c(G(\AAA_f^p),\C)$ when there is no possible confusion.

\begin{theorem}{\rm (\cite[Thm. A]{Ch})}\label{Cheeigen} There exists a unique eigenvariety $(X/L,\psi,\nu,Z)$ for $\cal{Z}_e$. It has the following extra properties: \ps
\begin{itemize}
\item[(iv)] $X$ is nested and equidimensionnal of dimension $m$. Moreover, $\nu(X)$ is a Fredholm hypersurface of $\WG$, hence $X$ inherits a canonical admissible covering. Precisely, $X$ is admissibly covered by the affinoid subdomains $\Omega \subset X$ such that $\omega(\Omega) \subset \WW$ is an open affinoid and that $\omega: \Omega \longrightarrow \omega(\Omega)$ is finite and surjective when restricted to any irreducible component of $\Omega$. \ps
\item[(v)] Let $Z'$ be the subset of $x \in X(\Qpb)$ such that $\omega(x)=(k_1,\cdots,k_m)\in \Z^{m,-}$, that (see Def. \ref{definitionFi})
$$\forall\, i\neq j, \, \, p^{-k_i+i}F_i(x) \neq p^{-k_j+j-1}F_j(x),$$ 
and that
$$ v(\psi(u_p)(x)) < 1+ {\rm Min}_{i=1}^{m-1} (k_{i+1}-k_i).$$
Then $Z' \subset Z$, and $Z'$ is an accumulation Zariski-dense subset of $X$.\ps
\item[(vi)] $\psi(\HH_{\rm ur}) \subset \OO(X)^0$.
\end{itemize}
\end{theorem}

\begin{pf} The existence of $X/L$ satisfying (i)-(vi) is \cite[Thm. A]{Ch} (using \cite{Buz}) when $e$ has the form $1_K$ 
for a net compact open subgroup $K \subset G(\AAA_f^{S_0\cup\{p\}})$ and $\WW$ is replaced by its open subspace of analytic characters. Moreover, assertion (v) (and the accumulation property of $Z$ in (iii)) requires there 
$p>2$. We will explain in \S\ref{eigenreview} to \S\ref{padIII} below how to extend the construction to any $e$ and to the full weight space $\WW$, in the spirit 
of \cite{Buz}. This will make (iii) and (v) valid also when $p=2$ by the same proof as \cite[Prop. 6.4.7]{Ch}. The fact that $X$ 
is reduced follows from \cite[Prop. 3.9]{jlpch}. The uniqueness assertion is Prop. \ref{eigenunique}.
\end{pf}

\begin{remark}\label{remcheeigen} \begin{itemize}
\item[i)] An independent construction of $X$ has been given by M.
Emerton in \cite{eminterp}. The admissible open subspace $X^{\ord}
\subset X$ defined by $|\psi(u_p)|=1$ 
was previously constructed by H. Hida in a much more general context (see \cite{hida}). 
It is actually closed and the induced map
$\omega: X^{\ord} \longrightarrow \WW$ is finite.
\item[ii)] Assume that $e=e_1+e_2$ is the sum of two orthogonal idempotents. Also we will not use it in what follows, let us note that by 
\cite{jlpch}, the eigenvariety $X_i$ of ($e_i$,$\HH$) has a natural closed embedding into $X$ 
commuting with $(\psi,\nu)$. Moreover, $X = X_1 \cup X_2$ (the intersection might be non empty). Actually, we could even show that $X$ is precisely the abstract gluing of $X_1$ and $X_2$ "over $(\psi,\nu)$".
\item[iii)] (A variant with a {\it fixed weight}) Let $i_0 \in \{1,\cdots,m\}$ and $k \in \Z$, and consider the subset 
${\cal Z}_{e,k_{i_0}=k} \subset {\cal Z}_e$ whose elements are the $(\psi_{(\pi,\RR)},\uk)$ such that $\uk$ has its $i_0^{th}$ term 
$k_{i_0}$ equal to $k$. Then there exists also a unique eigenvariety $X'$ for ${\cal Z}_{e,k_{i_0}=k}$ satisfying all the properties of Thm. \ref{Cheeigen}, except 
that it is equidimensional of dimension $m-1$. 
This follows for example verbatim by the same proof (see below) if we replace everywhere the space $\WW$ in this proof by 
its hypersurface $\WW_{k_{i_o}=k} \subset \WW$ parameterizing the characters whose $i_0^{th}$ component is fixed 
and equal to $x_{i_0} \mapsto x_{i_0}^k$. In most cases, $X$ turns out to be isomorphic 
to $X' \times X_1$ where $X_1$ is a suitable eigenvariety of $\U(1)$. As those $X_1$ are explicit 
({\it e.g.} they are finite over $\WW_1$), it is in general virtually equivalent to study $X$ or $X'$. 
\end{itemize}
\end{remark}

We end this paragraph by a discussion on idempotents, so as to shed light on the kind of sets $\cal{Z}_e$ that we can obtain. Of course, in the applications we will mostly choose $e$ as a tensor product of idempotents $e_l \in \Cc_c(G(\Q_l),\Qb)$ ($l \neq p$) such that $e_l=1_{G(\Z_l)}$ for $l \in S_0$ or big enough.

\begin{example}\label{exampleidemp}({\it see e.g.} \cite[\S2]{BuKu}) 
\begin{itemize} \item[i)] Of course, the simplest class of idempotents of $\Cc_c(G(\AAA_S),\Q)$ are the $$e_K:=\mu(K)^{-1}1_K \in \Cc(K\backslash G(\AAA_S)/K,\Q)$$ for any compact open subgroup $K \subset \G(\AAA_{S})$.
\item[ii)] A little more generally, if $\tau$ is an irreducible smooth $\Qb$-representation of such a $K$ (hence finite dimensional), 
the element $$e_{\tau} \in \Cc_c(G(\AAA_S),\Qb)$$
which is zero outside $K$ and coincide with $\frac{\dim_{\Qb}\tau}{\mu(K)}\tr(\pi^*)$ on $K$ is an idempotent. We see at once that for 
each smooth representation $V$ of $G(\AAA^S)$, $e_{\tau}V \subset V$ is the $\tau$-isotypic component of $V$. \par 
\item[iii)] ({\it Special idempotents}) Let $k$ be a field of characteristic $0$, 
$H=\Cc_c(G(\AAA_S),k)$, $e \in H$ an idempotent, and $\Mod_{e}$ be the full 
subcategory of the category of smooth $k[G(\AAA_S)]$-representations whose objects $V$ are 
generated by $eV$. Following the terminology of \cite[\S3]{BuKu}, we say that $e$ is {\it special} if the functor $V \mapsto eV$, $\Mod_{e} \rightarrow \Mod(eHe)$ is an equivalence of categories. If $e$ is special, the induction functor 
$$W \in \Mod(eHe) \mapsto I(W):=He\otimes_{eHe}W \, \, \in \Mod_e.$$
is a quasi-inverse of $V \mapsto eV$, hence is exact, and for any $V \in \Mod_e$, the natural surjection induced an isomorphism\footnote{Indeed, as $eI(W)=W$ for all $W$, the induced map $eI(V) \rightarrow eV$ is necessarily 
an isomorphism, hence the kernel of the map (\ref{Iexact}) is killed by $e$, whence is zero as $e$ is special property. This argument shows actually that $e$ is special if, and only if, $\Mod_{e}(G(\Q_l))$ is stable by subobjects.}
\begin{equation}\label{Iexact}	I(eV) \isomo V.\end{equation}
\item[iv)] Set $S=\{l\}$ to simplify. If $e=e_{\tau_l}$ for some $K_l$-type $\tau_l$, then $e$ is special if, and only if, 
$\tau_l$ is a type in the sense of Bushnell and Kutzko \cite{BuKu}. The simplest example, due to Borel, is the case where $e=e_{K_l}$ and $K_l$ is a Iwahori subgroup of $G(\Q_l)$, in which case $\Mod_e$ is the unramified Bernstein component. Moreover, by \cite[Cor. 3.9]{centrebern}, there exist arbitrary small compact open subgroups $K_l$ of $G(\Q_l)$ such that $e_{K_l}$ is special. However, as is well known, if $K_l$ is a maximal compact subroup then $e_{K_l}$ is not special in general (e.g. when $G(\Q_l)\isomo\GL_m(\Q_l)$ for $m>1$). 
\item[v)] (Bernstein components) Set again $S=\{l\}$. By results of Bernstein (see \cite{centrebern} or \cite[3.12, 3.12]{BuKu}), if $e$ is special then there is a finite set $\Sigma_l$ of Bernstein components of $G(\Q_l)$ such that $\Mod_e$ is the direct sum of these components (\cite[Prop. 3.6]{BuKu}). Reciprocally, for any finite set $\Sigma_l$ of components we can find a special idempotent $e_{\Sigma_l} \in \Cc(G(\Q_l),\Qb)$ whose associated set of components is $\Sigma_l$. This idempotent is not unique however in general, but all the equivalent ones will give rise to the same set $\cal{Z}_e$, hence to the same eigenvariety by virtue of the uniqueness Prop. \ref{eigenunique}. \par 
	This remark allow us in particular to say that there exist eigenvarieties for the subset $\cal{Z} \subset \cal{Z}_0$ parameterizing $p$-refined automorphic representations whose local components in a finite set of primes all lie in a given Bernstein component.
\end{itemize}
\end{example}	
	
\begin{remark} ($K$-types versus general idempotents) The aim of type theory is to show that the special idempotents $e_{\Sigma_l}$ above can be chosen of the form $e_{\tau}$ for some explicit $K_l$-type $\tau$. In our context, this extra information is helpful from a computational point of view. For example, if $e=e_{\tau}$ then we will see that the space of $p$-adic automorphic form of type $e$ is $$\Sc(V,r)=\tau \otimes_L \left( F(\Cc(V,r))\otimes_L \tau^*\right)^{K^p}$$ 
which is computable in theory. In general, some $e_{\Sigma_l}$ are given abstractly by images of some idempotents in the Bernstein center of $G(\Q_l)$ and we have very few control on them.
\end{remark}

\subsubsection{Review of the construction of the eigenvariety (\cite{Ch})} \label{eigenreview}	
The eigenvariety $X$ associated to $e$ is constructed by some formal process
(\cite{col2}, \cite{CM}, \cite{Ch}, \cite{Buz}) from the action of the compact
Hecke operators on the ONable Banach family of
spaces of $p$-adic automorphic forms of $G$. For example, 
$X(\Qpb)$ turns out to parameterize exactly the $\Qpb$-valued 
systems of Hecke eigenvalues on finite slope $p$-adic eigenforms of type $e$ for
$G$. \ps As our main theorem relates some Selmer group to the smoothness of $X$ at some
point, and for sake of completeness, we give below an essentially self-contained 
overview of the construction of $X$ and of the theory of $p$-adic automorphic forms alluded above\footnote{We warn the reader 
that some of the various conventions that we use in this book differ from the ones used in \cite{Ch}.}. Actually, we shall use also some objects occurring in this construction to define the families of 
admissible $G(\AAA_f^p)$-representation on $X$ in \S\ref{famrepgql}, as well as to define 
their NMSRPS locus in \S\ref{NMlocus}. The construction proceeds in three steps.

\subsubsection{Step I. The family of the $U^-$-stable principal series of
a Iwahori subgroup.}\label{padI} The theory of $p$-adic automorphic forms of $G$ relies essentially 
on the existence and properties of the $p$-adic family of the $U^-$-stable principal series of 
the Iwahori subgroup $I$ of $\GL_m(\Q_p)$. We take here and below the notations of 
\S\ref{atkinlehner} with $F=\Q_p$ and $\varpi=p$, except that we shall write $G(\Q_p) \isomo_v \GL_m(\Q_p)$ 
for the $G$ {\it loc. cit.} which is already used here for the unitary group over $\Q$. \ps
\newcommand{\Nob}{\overline{N}_0}
Let $\Nob$ be the subgroup of lower triangular elements of $I$. The product map in $G(\Q_p)$ induces an isomorphism
$$ \Nob \times B \isomo I\,B.$$
If $u \in U^-$ then $u^{-1}\Nob u \subset \Nob$, hence (see Proposition \ref{atkinl}) $M^{-1}IB \subset IB$.
\indent Let $$\chi: T^0 \longrightarrow \OO(\WW)^*$$
denote the tautological character. If $V \subset \WW$ is either an open affinoid or a closed point we denote by $\chi_V: T^0 \longrightarrow A(V)^*$ 
the induced continuous character. Fix such a $V$. There exists an smallest 
integer $r_V \geq 0$ such that for any integer $r \geq r_V$, $\chi_V$ restricts to 
an analytic $A(V)$-valued function on the subgroup of elements of $T^0$ with coefficients in $1+p^r\Z_p$. Viewing the character $\chi$ of $T^0$ as a characters of $B$ which is 
trivial on $U N$, it makes then sense to consider for $r \geq r_V$ the space 
\[ \Cc(V,r) = \left\{ \begin{array}{ll}
f: IB \longrightarrow A(V),\, \, f(xb)=\chi_V(b)f(x)\, \, \, \, \forall\, \,  x \in IB, \, \, b \in B, \\
f_{|\Nob}\, {\rm\, is\,}\, 
r{\rm-analytic}.
\end{array} \right\} \]
Let us recall what $r$-analytic means. Let $\{n_{i,j}\}_{i>j}$ be the natural algebraic coordinate maps on $\Nob$. 
A function $f: \Nob \rightarrow A(V)$ is said to be {\it $r$-analytic} 
if for each $a \in \Nob$, the induced map $$f_a: \Nob \rightarrow A(V), \,
\, \, n \mapsto f(an),$$ lies in the Tate-algebra
$A(V)\langle\{p^rn_{i,j}\}_{i>j}\rangle$. If we endow this latter algebra with the sup norm, 
then the norm $|f|:=\sup_a |f_a|$ makes $\Cc(V,r)$ a Banach $A(V)$-module which is $A(V)$-ONable by construction
\footnote{It is isometric to $A(V)\langle\{n_{i,j}\}_{i>j}\rangle^{p^{rm(m-1)/2}}$.}. It is equipped with an 
integral $A(V)$-linear action of $M$ by left translations: $(m.f)(x)=f(m^{-1}x)$. If we set
$$U^{--}=\{u=(p^{a_1},\cdots,p^{a_m}) \in U, a_1>a_2>\cdots>a_m\},$$
then an immediate computation shows that the action of any $u \in U^{--}$ on $\Cc(V,r)$ is $A(V)$-compact.\ps
The family $\{\Cc(V,r),V, r\geq r_V\}$ of $M$-modules defined above satisfies some compatibilities. If $V' \subset V$ is another open affinoid or closed point, then the natural 
map $\Cc(V,r) \rightarrow \Cc(V',r)$ induces an $M$-equivariant isomorphism $\Cc(V,r)\widehat{\otimes}_{A(V)}A(V')$.
Moreover, the natural inclusion $\Cc(V,r) \longrightarrow \Cc(V,r+1)$ is $A(V)[M]$-equivariant and compact. If $r>0$ and $u \in U^{--}$, 
then the action of $u$ factors through the compact inclusion $\Cc(V,r-1) \longrightarrow \Cc(V,r)$ above. \ps

For any continuous character $\psi: T \longrightarrow L^*$ and $r \geq r_{\psi}:=r_{\psi_{|T^0}}$, let us consider similarly the 
$L$-Banach space $i_B^{IB}(\chi)$ of functions $f: IB \rightarrow L^*$ whose restriction to $\Nob$ is $r$-analytic and which 
satisfy $f(xb)=\psi(b)f(x)$ for all $x \in IB$ and $b \in B$. It has again an action of $M$ by left translations. 
If $\psi': T \rightarrow L^*$ is another continuous character and 
$r\geq r_{\psi},r_{\psi'}$, we define a natural map 
	$$i_B^{IB} \psi \longrightarrow i_B^{IB} (\psi'\psi), \, \, f \mapsto (x \mapsto \psi'(x)f(x)),$$
where for $x \in IB$, $\psi'(x):=\psi'(t)$ for $t \in T$ the unique element such that $x \in \Nob t N$. We check at once that this map 
is well defined and that it induces an $M$-equivariant isomorphism\footnote{Recall that by Proposition \ref{atkinl} (i), each character of $U^-$ extends uniquely to a character of $M$ trivial on $I$. 
} 
\begin{equation} \label{twistapasoublier} (i_B^{IB} \psi)\otimes {\psi'}^{-1} \longrightarrow i_B^{IB} (\psi\psi').\end{equation}

Assume now that $V=\{\uk\}$ with $\uk \in \Z^{m,-}$, in which case $r_V=0$. The choice of an highest-weight vector in 
$W_{\uk}(\Q_p)$ with respect to $B$ gives an $M$-equivariant embedding $$W_{\uk}(L)^*  \longrightarrow i_B^{IB}(\delta_{\uk}).$$ Hence we get by 
(\ref{twistapasoublier}) a canonical (up to multiplication by $L^*$) $M$-equivariant embedding\footnote{The twists appearing there comes from the fact that we chosen to extend trivially on $U$ the induced character in the definition of $\Cc(V,r)$. This choice could have been avoided by introducing the space of $p$-adic characters 
of $T$ rather than $T^0$. However, as $\Z^m$ is not Zariski-dense in that space, this would have introduced other nuisances...} 
$$W_{\uk}(L)^*\otimes \delta_{\uk}  \longrightarrow \Cc(\uk,0)=i_B^{IB}(\chi_{\uk}).$$
Actually, the subspace of $\Cc(\uk,0)$ defined above is exactly the subspace of functions on $IB$ which are restrictions of polynomial maps on $G(\Q_p)$. \ps

\medskip
\noindent 
\subsubsection{Step II. $p$-adic automorphic forms.}\label{padII} Using the $M$-modules defined above as coefficient systems, we can define various analogous Banach spaces of 
{\it $p$-adic automorphic forms} for $G$. Consider the subring $\HH^-:=\ATL^- \otimes \HH_{\rm ur} \subset \HH$. It will be 
convenient to introduce a functor $F: {\rm Mod}(L[M]) \longrightarrow {\rm Mod}(\ATL^- \otimes \Z[\GAFP])$ by 
$$\FT(E) := \left\{ \begin{array}{lll} f: G(\Q)\backslash G(\AAA_f) \longrightarrow E,\\ 
f(g(1\times k_p))=k_p^{-1}f(g),\, \, \forall g\in G(\AAA_f),\, \, k_p \in I,\, \\
f \, \, {\rm is\, \,smooth\, \,  outside\, \,} p. \end{array} \right\} $$
The group $\GAFP$ acts on this space in a smooth way by right translations, and it commutes with the natural action of $\ATL^-$. 
The direct summand $eF(E) \subset F(E)$ is then a $\HH^-$-module in a natural way. Let $K=I\times K^p \subset G(\AAA_f)$ be a 
compact open subgroup which is sufficiently small so that $e=e.e_{K^p}$, and such that $K_l=G(\Z_l)$ for each $l\in S_0$. If we choose a writing 
(see \S\ref{defautformG} ii))
$$G(\AAA_f)=\prod_{i=1}^{h_K} G(\Q)x_iK, \, \, \, \Gamma_i:=x_i^{-1}G(\Q)x_i \cap K,$$ 
then $\Gamma_i$ is a finite group, and we may even assume by reducing $K$ that $\Gamma_i$ is trivial for each $i$. The map $f \mapsto (f(x_1),\cdots,f(x_{h_K}))$ induces a $L$-linear isomorphism 
\begin{equation}\label{Fiswellbehaved}e_KF(E) \isomo E^{h_K}. \end{equation}
In particular, the functor $E \mapsto eF(E)$, $\Mod(L[M]) \longrightarrow \Mod(\HH^-)$, is an extremely well behaved functor, as a direct summand
of $e_KF$. \ps 

	Let $\uk \in \Z^{m,-}$. We check at once using $\iota_p\iota_{\infty}^{-1}$ that 
$eF(W_{\uk}(L)^*)$ is a $\HH^-$-stable $L$-structure of the space $\iota_{\infty}(e)A(G,W_{\uk}(\C))$ 
of complex automorphic forms of weight $W_{\uk}(\C)$ and type $e$.\ps
	Let $V \subset \WW$ is an affinoid subdomain or a closed point, and $r\geq r_V$. We define an $\HH^-$-module by setting
$$\Sc(V,r):=eF(\Cc(V,r)).$$
This is the space of {\it $p$-adic automorphic forms of weight in $V$, radius of convergence $r$ and type $e$}. It has a natural structure of 
Banach $A(V)$-module which is a topological direct summand of the ONable Banach module 
\footnote{For $f \in e_KF(\Cc(V,r))$, we set $|f|:=\sup_{x\in G(\AAA_f)} |f(x)|=\sup_{i=1}^{h_K} |f(x_i)|$. The isomorphism 
(\ref{Fiswellbehaved}) induces an isometry $e_KF(\Cc(V,r)) \isomo \Cc(V,r)^{h_K}$ hence $e_KF(\Cc(V,r))$ is $A(V)$-ONable. 
We give $\Sc(V,r) \subset e_KF(\Cc(V,r))$ the subspace topology, it is closed as the image of the continuous 
linear projector $e$ on $e_KF(\Cc(V,r))$.} (hence satisfies Buzzard's (Pr) condition). It is 
equipped with an $A(V)$-linear action of $\HH^-$, each $h \in \HH^-$ being bounded by $1$ and each element of 
$U^{--} \subset \HH^-$ being $A(V)$-compact. By formula (\ref{Fiswellbehaved}), the collection of spaces $\{\Sc(V,r), V, r\geq r_V\}$ satisfies 
exactly the same compatibilities as $\{ \Cc(V,r),\, V,\,  r\geq r_V\}$.
\ps
	For $\uk=(k_1,\cdots,k_m) \in \Z^{m,-}$, we have moreover a natural $\HH^-$-equivariant inclusion 
\begin{equation}\label{classicalinclusion} eF(W_{\uk}(L)^*)\otimes \delta_{\uk} \hookrightarrow \Sc(\uk,0),\end{equation}
whose image is usually referred as the subspace of {\it classical} $p$-adic automorphic forms. The {\it control theorem}
asserts then that an element $f \in \Sc(\uk,r)$ which is in the generalized eigenspace of the compact operator $u_p \in U^{--} \subset \HH$ with respect to any 
eigenvalues $\lambda \in \Qpb$ such that $v(\lambda) < 1 + {\rm Min}_{i=1}^m (k_{i+1}-k_i)$ is actually classical. \ps

\medskip

\noindent 
\subsubsection{Step III. Fredholm series and construction of the eigenvariety.}\label{padIII} As any $h \in U^{--}\HH^-$ acts compactly on the 
(Pr)-family of Banach modules $\{\Sc(V,r),V,r\geq r_V\}$, there is a unique Fredholm series $P_h(T) \in 1+T\OO(W)\{\{T\}\}$ such that for any $V \subset \WW$ open affinoid 
or closed point and $r \geq r_V$,
$${P_h(T)}_{|V}=\det(1-Th_{|\Sc(V,r)}) \in 1+TA(V)\{\{T\}\}.$$
\indent Set $P:=P_{u_p}$, and consider the Fredholm hypersurface $Z(P) \subset \WG$, that is the closed subspace defined by $P=0$. As any Fredholm 
hypersurface, $Z(P)$ is canonically admissibly covered by its affinoid subdomains $\Omega^*$ such that $\pr_1(\Omega^*)$ is an open affinoid of $\WW$ and that 
the induced map $\pr_1: \Omega^* \longrightarrow \pr_1(\Omega^*)$ is finite. Here $\pr_1$ is the first projection $\WG \longrightarrow \WW$. Let us 
denote by $\Cc^*$ this canonical covering. This covering is easily seen to be stable by finite intersections, by pullback over affinoid 
subdomains of $\WW$, hence to satisfy the following good property: \ps
\begin{center}(*) if $\Omega_1^*,\, \Omega_2^* \in \Cc^*$ then $\Omega_1^* \cap \Omega_2^*$ is a clopen subspace of $\Omega_1 \times_{V_1} (V_1 \cap V_2)$.
\end{center} \ps
The eigenvariety $X$ will then be constructed as a finite map $\nu: X \longrightarrow Z(P)$ as follows. 
Let $\Omega^* \in \Cc^*$ and $V:=\pr_1(\Omega^*)$. There is a unique factorization $P=QR$ in $A(V)\{\{T\}\}$ 
where $Q \in 1+TA(V)[T]$ has a unit leading coefficient and is such that $\Omega^*=Z(Q)$ is a closed and open subspace of $Z(P) \times_{\WW} V$. To this factorization corresponds, for $r \geq r_V$, 
a unique Banach $A(V)$-module decomposition $$\Sc(V,r)=S(\Omega^*) \oplus N(\Omega^*,r),$$
which is $\HH^-$-stable, and such that: \ps
\begin{itemize}
\item[-] $S(\Omega^*)$ is a finite projective $A(V)$-module, which is independent of $r\geq r_V$, \ps
\item[-] the characteristic polynomial of $u_p$ on $S(\Omega^*)$ is the reciprocal polynomial $Q^{\rec}(T)$ of $Q(T)$, and $Q^{\rec}(u_p)$ is invertible on $N(\Omega^*,r)$. 
\end{itemize}\ps
\indent The local piece $\Omega$ of $X$ is then by definition the maximal spectrum of the $A(V)$-algebra 
generated by the image of $\HH=\HH^{-}[u_p]^{-1}$ in $\End_{A(V)}(S(\Omega^*))$. It is equipped by construction 
with a ring 
homomorphism $\HH \longrightarrow A(\Omega)$, with a finite map 
$\nu: \Omega \longrightarrow \Omega^*$, and with a finite $A(\Omega)$-module $S(\Omega^*)$. We check then that the $\Omega$ and the maps above 
glue uniquely over $\Cc^*$ to an object $(X,\psi,\nu)$ as in the statement of Proposition \ref{Cheeigen}, which is easy using 
the property (*) mentioned above of the admissible covering $\Cc^*$. In other words, 
the coherent sheaves of $\OO$-algebras $\{\widetilde{A(\Omega)},\Omega^* \in \Cc^*\}$ glue canonically to a coherent 
$\OO_{Z(P)}$-algebra, and $\nu: X \longrightarrow Z(P)$ is its relative spectrum (see \cite[\S2.2]{conradqcoh}). In the same way, the locally free coherent sheaves 
$\{\widetilde{S(\Omega^*)},\Omega^* \in \Cc^*\}$ glue canonically to a 
coherent sheaf on $Z(P)$. This sheaf is a $\nu_{*}\OO_X$-module in a natural way, hence has the form 
$\nu_* \Sc$ for as a coherent sheaf $\Sc$ on $X$.

\begin{definition}\label{defCS} We denote by $\Cc$ the admissible covering $\nu^{-1}\Cc^*$ of $X$ and by $\Sc$ the 
coherent sheaf on $X$ defined above. If $\Omega=\nu^{-1}(\Omega^*) \in \Cc$, then $\Sc(\Omega)=S(\Omega^*)$.
\end{definition}

\begin{remark} \label{remacoh}({\it On quasicoherent and coherent sheaves on rigid spaces}) Let $X$ be a rigid analytic space over $k$. 
An $\OO_X$-module $\Fc$ is said to be {\it quasicoherent} (resp. coherent) if there exists an admissible covering $\{U_i\}$ 
of $X$ by affinoid subdomains such that $\Fc_{|U_i}$ is the sheaf $\widetilde{M_i}$ 
associated to some $\OO(U_i)$-module $M_i$ (resp. such that $M_i$ is finite type over $\OO(U_i)$) (see \cite[\S9.4.2]{BGR}, 
\cite[\S2.1]{conradqcoh}). Contrary to the case of schemes, it does not imply in general that for any affinoid subdomain $U$, $\Fc_{|U}$ is associated to 
an $\OO(U)$-module (see Gabber's counterexample \cite[Ex. 2.1.6]{conradqcoh}). This holds however when $U$ lies in some $U_i$,
when $\Fc$ is coherent (\cite[\S9.4.3]{BGR}), or when $\Fc$ is globally on $X$ a direct inductive limit 
of coherent $\OO_X$-modules (\cite[Lemma 2.1.8]{conradqcoh}). \par 
	In our applications, we will define some quasicoherent sheaves on $X$ using the covering $\Cc$, but they will all be 
direct inductive limit of coherent sheaves. \end{remark}

\subsection{The family of $G(\AAA_f^p)$-representations on an eigenvariety of idempotent type}\label{bernsteinnonmonodrmic}

In all this part, we keep the notations of \S\ref{eigenidem}. 
In particular, $X$ is the eigenvariety associated to the idempotent $e$ given by Theorem \ref{Cheeigen}, or its variant 
with a fixed weight as in Remark \ref{remcheeigen} (iii). 

\subsubsection{The family of local representations on $X$}\label{famrepgql} Let us fix 
some finite set $S$ of primes\footnote{We hope that there will be no possible confusion 
with the letter $S$ occurring in the spaces $S(\Omega^*)$ or $S(V,r)$.} such that $S \cap (S_0\cup \{p\})=\emptyset$. Assume moreover that the idempotent $e$ decomposes as a tensor product of idempotents
at $l \in S$ and outside $l$: $e=e_S \otimes e^S$ and $e_S=\otimes_{l\in S} e_l$, $e_l^2=e_l \in \Cc(G(\Q_l),\Qb).$\ps
The eigenvariety $X$ carries a natural sheaf of admissible $G(\AAA_S)$-representations that we will describe now. 
For $V \subset \WW$ an open affinoid and $r \geq r_V$, we have by definition a split inclusion
\begin{equation}\label{injlocl} \Sc(V,r)=eF(\Cc(V,r)) \subset F(\Cc(V,r)),\end{equation}
and the latter space is a smooth $G(\AAA_S)$-module as $p \notin S$. We fix now $\Omega^* \in \Cc^*$, 
set $\Omega=\nu^{-1}(\Omega^*)$, $V=pr_1(\Omega^*)$, and we consider $S(\Omega^*) \subset S(V,r)$ as in \S\ref{padIII}. 

\begin{definition}\label{defpil} We define $\Pi_S(\Omega)$ as the $\Z[G(\AAA_S)]$-submodule of $F(\Cc(V,r))$ 
generated by $S(\Omega^*)$.
\end{definition}

By definition, $\Pi_S(\Omega)$ is an $\HH^- \otimes A(V)$-submodule and the natural map 
$\HH^- \otimes A(V) \longrightarrow \End(\Pi_S(\Omega))$ factors through $A(\Omega)$ as it does on the generating subspace $\Sc(\Omega^*)$ 
of $\Pi_S(\Omega)$. As a consequence, $\Pi_S(\Omega)$ is an $A(\Omega)$-module in a natural way. It is independent of $r\geq r_V$ as 
$\Sc(\Omega^*)$ is.

\begin{prop}\label{proppil} Let $\Omega \in \Cc$. 
\begin{itemize}
\item[(i)] $\Pi_S(\Omega)$ is an $A(\Omega)$-admissible smooth representation of $G(\AAA_S)$. 
\item[(ii)] The natural inclusion $\Sc(\Omega) \longrightarrow e_S\Pi_S(\Omega)$ is an equality.\ps
\end{itemize}
\noindent Moreover,\begin{itemize}
\item[(iii)] The sheaves of $\OO$-modules $\{\widetilde{\Pi_S(\Omega)},\Omega \in \Cc\}$ glue canonically to a quasicoherent 
smooth $\OO_X[G(\AAA_S)]$-module, and (ii) glue to an isomorphism $\Sc \isomo e_S\Pi_S$.
\item[(iv)] For each compact open subgroup $J \subset G(\AAA_S)$ the 
subsheaf of $J$-invariants $\Pi_S^J \subset \Pi_S$ is a coherent $\OO_X$-module, and $\Pi_S=\bigcup_J \Pi_S^J$. 
\item[(v)] For all $x \in X$, $(\Pi_S)_x$ is torsion free over $\OO_{\omega(x)}$, hence also over $\OO_x$.
\end{itemize}
\end{prop}

\begin{pf} We check first assertion (ii). Let $Q \in A(V)[T]$ be the polynomial attached to $\Omega^*$ as in \S\ref{padIII}, so that 
$\Sc(\Omega^*)$ is the Kernel of $Q^{\rm rec}(u_p)$ on $eF(\Cc(V,r))$. As $p\notin S$, $Q^{\rm rec}(u_p)\Pi_S(\Omega^*)=0$, hence (ii) holds by definition. \par
We know that $\Pi_S(\Omega)$ is smooth as $F(\Cc(V,r))$ is, 
hence (i) follows from (ii) and lemma \ref{delbern} (i). Assertion (iii) follows easily from the properties of the admissible 
covering $\Cc^*$, the proof is similar to the gluing argument for the sheafs $\Sc$ and $\nu_*\OO_X$ 
so we leave the details to the reader. To prove (iv), note that for any 
$\Q$-algebra $A$, any $A$-linear representation $V$ of $J$, and any $A$-module $M$, the natural map 
	$$ V^J \otimes_A M \longrightarrow (V \otimes_A M)^J$$
is an isomorphism (argue as in Lemma \ref{jacquet} (b)). Part (iv) follows now from (i) and (iii). \ps
Before showing (v), let us recall 
that by construction $\OO_x$ is a $\OO_{\omega(x)}$-subalgebra of the endomorphism ring of a finite free $\OO_{\omega(x)}$-module, hence the total fraction ring of $\OO_x$ identifies with 
$\OO_x \otimes_{\OO_{\omega(x)}}{\rm Frac}(\OO_{\omega(x)})$. As a consequence, it suffices to check that 
$(\Pi_S)_x$ is torsion free over $\OO_{\omega(x)}$. But for each $\Omega^*, V$ and $\Omega$ as above, 
$\Pi_S(\Omega)$ is a subpace of $F(\Cc(V,r))$, which is clearly torsion free over $A(V)$.
\end{pf}

\begin{lemma}(Bernstein) \label{delbern} 
Let $k$ be a field of characteristic $0$, $A$ a noetherian $k$-algebra and $V$ a smooth $A[G(\AAA_S)]$-representation. Assume that for some decomposed idempotent 
$e \in \Cc_l(G(\AAA_S),k)$, $eV$ is finite type over $A$ and 
generates $V$ as an $A[G(\Q_l)]$-module. Then: \begin{itemize}
\item[(i)] $V$ is $A$-admissible,
\item[(ii)] if $A$ is moreover finite dimensional over $k$, $V$ is of finite length on $A[\G(\Q_l)]$.
\end{itemize}
\end{lemma}

\begin{pf} Let us show (i). By induction on $|S|$, we may assume that $S=\{l\}$. By \cite[Prop. 3.3]{centrebern}, and more precisely by "variante 3.3.1" and the remark following 
Corollary 3.4 {\it loc. cit.}, $V$ is $Z(G(\Q_l)) \otimes_k A$-admissible 
where $Z(G(\Q_l))$ is the center of the $k$-valued Hecke-algebra of $G(\Q_l)$. 
As $A[G(\Q_l)]eV=V$, the action of $Z(G(\Q_l))\otimes_k A$ on $V$ 
factors through its faithful quotient $A' \subset \End_A(eV)$. 
As $eV$ is finite type over $A$ by assumption, and $A$ is noetherian,
so is $A'$, hence $V$ is $A$-admissible. \ps
	The second assertion follows from (i) as $V$ is then $k$-admissible and of finite type over $k[G(\AAA_S)]$ 
(use \cite[Cor. 3.9]{centrebern}).
\end{pf}

For sake of completeness, we end this discussion by a study of the fibers of $\Pi_S$ at a point of $X$. 
We fix $x \in X$ with residue field $k(x)$, hence we get a natural system of Hecke-eigenvalues $\psi_x: \HH \longrightarrow k(x)$. 
To this system of Hecke-eigenvalues corresponds a generalized $\HH$-eigenspace $\Sc^{\psi_x} \subset \Sc(\omega(x),r)=eF(\Cc(\omega(x),r))$, for $r$ big enough.

\begin{definition}\label{defpilpsix} We denote by $\Pi_S^{\psi_x}$ the 
$\OO_x/m_{\omega(x)}\OO_x$-representation of $G(\AAA_S)$ generated by the (finite dimensional) subspace 
$\Sc^{\psi_x}$ of $F(\Cc(\omega(x),r))$. It is a finite length admissible representation of $G(\AAA_S)$ by Lemma \ref{delbern} (ii).
\end{definition}

\begin{definition}\label{clasvsad} Assume moreover that $x \in Z$, so $\psi_x(\HH)\subset \Qb$. We denote by $\Pi_x$ the $k(x)$-model\footnote{Defined by $\iota_p$ and $\iota_{\infty}$.} 
of the complex $\G(\AAA_S)$ subrepresentation of $A(G,W_{\uk})$ generated by the $\HH$-eigenspace of $e(A(G,W_{\uk}))$
for the system of eigenvalues $\iota_{\infty}\iota_p^{-1}(\psi_x)$. 
It is a semisimple $k(x)$-representation. \end{definition}

\begin{prop}\label{evalxpil} \begin{itemize}
\item[(i)] The natural map $(\Pi_S)_x/m_{\omega(x)}(\Pi_S)_x \longrightarrow \Pi_S^{\psi_x}$ is surjective and induces an 
isomorphism $\Sc_x/m_{\omega(x)}\Sc_x \isomo \Sc^{\psi_x}$.\ps
\item[(ii)] If $x \in Z$, then $\Pi_x$ is a subrepresentation of $\Pi_S^{\psi_x}$.
\item[(iii)] If $x\in Z$ is $U$-non critical then $\Pi_x \isomo \Pi_S^{\psi_x}$. 
\end{itemize}
\end{prop}

\begin{pf} Let $\Omega^* \in \cal C^*$ containing $x$ as above, we will argue as in the proof of Prop. \ref{proppil} (ii) of which we take the notations. As $\Sc(\Omega) \subset \Sc(V,r)=eF(\Cc(V,r))$ is projective and direct summand, the natural map
$$\Sc(\Omega)/m_{\omega(x)} \Sc(\Omega) \longrightarrow eF(\Cc(\omega(x)),r),$$
is injective and the Fredholm series of $u_p$ on $\Sc(\Omega)/m_{\omega(x)} \Sc(\Omega)$
is the evaluation of $Q(T)$ at $\kappa(x)$. By taking the $\psi_x$-generalized eigenspace we get that $$\Sc_x/m_{\omega(x)}\Sc_x \isomo \Sc^{\psi_x},$$
and (i) follows. The point (ii) follows from 
(\ref{classicalinclusion}) of \S\ref{padII}, and (iii) from the {\it small slope forms are classical} result of \cite[Prop. 4.7.4]{Ch}. 
\end{pf}

When $e_S$ is special (see \ref{exampleidemp} (iii)), most of our results hold in the stronger form.

\begin{cor}\label{corspecial} Assume that $e_S$ is special. \begin{itemize}
\item[(i)] The natural surjection induces an isomorphism $$\Cc_c(G(\AAA_S),L)e_S\otimes_{e_S\Cc_c(G(\AAA_S),L)e_S} \Sc \longrightarrow \Pi_S.$$
\item[(ii)] For all $x \in X$, $(\Pi_S)_x$ is flat over $\OO_{\omega(x)}$.
\item[(iii)] For all $x \in X$, the map of Prop. \ref{evalxpil} (i) is an isomorphism.
\end{itemize}
\end{cor}

\begin{pf} By Prop. \ref{proppil} (ii), $e_S\Pi_S=\Sc$, thus (i) follows from the discussion in Example \ref{exampleidemp}. Assertion (ii) follows formally from (i), the fact that for each $x \in X$, $\Sc_x$ is finite free over $\OO_{\omega(x)}$, 
and from the exactness of the functor $I$ defined in Example \ref{exampleidemp}. The map of Prop. \ref{evalxpil} (i) induces an 
isomorphism after projection to $e_S$ by {\it loc. cit.}, hence is an isomorphism as $e_S$ is special, which proves (iii).
\end{pf}

\subsubsection{The NMSRPS locus of $X$}\label{NMlocus}
We keep the assumptions and notations of \S\ref{eigenidem}. We fix a finite set $S_{NM}$ of primes $l$ such that $G(\Q_l)$ is quasisplit and assume that $e$ is a tensor product of idempotents $e_l$ with $l \in S$ by an idempotent outside $S_{NM}$. Recall that we defined some $\OO_X[G(\Q_l)]$-modules $\Pi_{S_{NM}}$ in \S\ref{famrepgql}. \ps

Let $X_0 \subset X$ be the subset of points $x$ such that for each $l \in S$, $\Pi_x \otimes_{\OO_x} k(x)$ contains a NMSRPS $G(\AAA_S)$-representation in the sense of \S\ref{defNMSRPS} (see Remark \ref{variantNMSRPS} when $|S|>1$). Let $X_{NM} \subset X$
be the Zariski-closure of $X_{0}$, we view it as a reduced closed subspace of $X$. Let also $\cal{Z}_{e,NM} \subset \cal{Z}_e$ be the subset parameterizing $p$-refined automorphic $\pi$ such that $\pi_l$ is NMSRPS for each $l \in S$.

\begin{prop}\label{existeigenNMSRPS} There exists a unique eigenvariety for $\cal{Z}_{e,NM}$, namely $$(X_{NM},\psi_{|X_{NM}},\nu_{|X_{NM}},Z\cap X_{NM}).$$ \ps
	$X_{NM}$ is a union of irreducible components of $X$, hence equidimensional of dimension $m$ if $\cal{Z}_{e,NM}\neq \emptyset$. It satisfies also properties (iv), (v) and (vi) of $X$ (see Theorem \ref{Cheeigen}).
\end{prop}

\begin{pf} For any open affinoid $V \subset X$, set $V_0=V\cap X_0$ and define $\overline{V}_0 \subset V$ to be the Zariski-closure of $V_0$ equipped with its reduced structure. By Prop. \ref{proppil} (see also remark \ref{remacoh}), ${\Pi_{S_{NM}}}_{|V}$ is the sheaf associated to the $A(V)$-admissible representation $\Pi_{S_{NM}}(V)$, to which we can apply the construction of \S\ref{constrNMSRPS}. By definition, $V_0$ is the intersection of the ${\rm Spec}(A(V))_0$ defined there with its subspace $V={\rm Specmax}(A(V))$ (and the Zariski-topology of $V$ is by definition the topology induced from ${\rm Spec}(A(V))$). As $A(V)$ is a Jacobson ring, and as $\overline{{\rm Spec}(A(V))}_0$ is constructible by Prop. \ref{NMSRPS} (i), we check easily that we also have $\overline{V}_0=V\cap \overline{{\rm Spec}(A(V))}_0$. By Prop. \ref{NMSRPS}(ii) and by Prop. \ref{proppil} (v), we know that $\overline{V}_0$ is a union of irreducible components of $V_0$.\ps
	We claim that for any two open affinoids $V,\, W \subset X$, \begin{equation}\label{gluingNM} \overline{V}_0 \cap W=\overline{(V\cap W)}_0.\end{equation} 
Note that $V\cap W$ is affinoid as $X$ is separated, and by replacing $W$ by $V\cap W$ in (\ref{gluingNM}), we may assume that $W \subset V$. Moreover, the inclusion $\supset$ above is clear as $V_0 \cap W=(V\cap W)_0$, thus it only remains to prove that $\overline{V}_0\cap W \subset \overline{W}_0$. As we know that $\overline{V}_0$ has a Zariski-dense open subset $V' \subset V_0$ by Prop. \ref{NMSRPS} (i), $V' \cap W$ is Zariski-dense in $W \cap \overline{V}_0$ by Lemma \ref{zardensecomparaison} (applied to $Y=\overline{V}_0$, $U=V'$, $\Omega=W\cap Y$), and we are done.\par
	As a consequence of (\ref{gluingNM}), all the 
$\overline{V}_0$ glue to a reduced closed subspace $T \subset X$. By construction, $X_0 \subset T$ is Zariski-dense, as it satisfies the much stronger assertion that for any open affinoid $V$, $V_0=V\cap X_0$ is Zariski-dense in $\overline{V}_0=V\cap T$. Hence $T=X_{NM}$ and the proposition follows at once.
\end{pf}

\begin{lemma}\label{zardensecomparaison} Let $Y$ be an affinoid, $U \subset Y$ a Zariski-open subset, and $\Omega \subset Y$ an affinoid subdomain. If $U$ is Zariski-dense in $Y$, then $\Omega \cap U$ is Zariski-dense in $\Omega$.
\end{lemma}
\begin{pf} Set $F:=Y\backslash U$, we have to show that 
$F\cap \Omega$ does not contain any irreducible component of $\Omega$. If it was the case, $F\cap \Omega$ would contain an affinoid subdomain of $\Omega$, hence $F$ would contain an affinoid subdomain of $Y$, 
as well as each irreducible component of $Y$ containing it, but this is a contradiction.
\end{pf}

\subsection{The family of Galois representations on eigenvarieties}${}^{}$\ps
\label{famillegaloisiennesurX}

In this part, we explain how the existence of Galois representations attached to classical automorphic representations for $G$ give rise to a
family of Galois representations on eigenvarieties. We keep the notations and assumptions of \S\ref{basicsettingeigenv}, as well as those of \S
\ref{setnotrepm}\footnote{Of course, we take the same choice of $\Qb$, $\Qpb$, $\iota_p$ and $\iota_{\infty}$ in both cases.}. 

\subsubsection{Setting} \label{notationsrepgal} So as not to multiply the statements, let us assume once and for all that 
$G$ is:\begin{itemize}
\item[(a)] either the group $\U(m)$ defined in \S\ref{definitionUm}, in which case we assume that hypothesis (Rep(m)) 
of \S\ref{sectrepm} 
holds,
\item[(b)] or a definite unitary group such that for any finite prime $l$, $G(\Q_l)$ is either quasisplit or 
isomorphic to the group of invertible elements of a central division algebra over $\Q_l$, 
this latter case occurring at least for one $l=:q$. Assume moreover that $G(\Q_p) \simeq \GL_m(\Q_p)$.
\end{itemize}
As explained in Remark \ref{remsectrepm} (vii), recall that in the second case, the obvious analog of condition (Rep(m)) is known 
if we forget property (P3). \ps

We assume moreover that the set $S_0$ defining $\HH_{\rm ur}$ 
has Dirichlet density one, and we fix a decomposed compact open subgroup $K^p\subset G(\AAA_f^p)$ as well as a finite set
$S$ of primes, such that $p \in S$ and that for each $l \notin S$ or in $S_0$, $K_l$ is a maximal hyperspecial or a very special compact subgroup. We choose the decomposed idempotent $e$ such that $ee_{K^p}=e_{K^p}$. 
If we are in case (b), we assume moreover that $e_q$ vanishes on the one dimensional representations of the division algebra\footnote{We need to make this technical condition to 
ensure that automorphic representations of $G$ admit associated Galois representations 
(see \cite[Thm. 3.1.4]{HL} and \S\ref{remsectrepm} (vii)).} $G(\Q_q)$.
We fix also a finite set of primes $S_{NM} \subset S$ that we assume to be empty in case (b), and
 define $\cal{Z} \subset \cal{Z}_e$ to be the set parameterizing $p$-refined 
automorphic representations of type $e$ which are $NMSRPS$ at primes in
$S_{NM}$ (note that $G(\AAA_{S_{NM}})$ is quasisplit), with one of the weights fixed if we like (sse Rem. \ref{remcheeigen} (iii)), 
and let $$(X,\psi,\omega,Z)$$ be the corresponding eigenvariety given by Prop. \ref{existeigenNMSRPS}, which is a 
closed subspace of the eigenvariety $X_e$ associated to $e$. \ps
	Recall that $\G_{E,S}$ is the Galois group of a maximal algebraic extension of $E$ which 
is unramified outside the primes above $S$. For each regular automorphic representation $\pi$, 
properties (P0)-(P1) assert the existence of a unique semisimple continuous representation 
$$\rho_{\pi} : \G_{E,S} \longrightarrow \GL_m(\Qpb),$$
such that for each prime $l=w\bar{w}\in S_0$, the trace of a geometric Frobenius at $x$, say $$\Frob_w \in \G_{E,S},$$ is the trace 
of the Langlands conjugacy class of $\iota_p\iota_{\infty}^{-1}(\pi_w.|\det|^{\frac{1-m}{2}})$. Let 
$$h_w \in \Cc(K_l\backslash \Cc(G(\Q_l)/K_l,\Z) \isomo_w \Cc(\GL_m(\Z_l)\backslash \GL_m(\Q_l)/\GL_m(\Z_l),\Z)$$
be the usual Satake element $[K_l(l,1,\dots,1)K_l]$, it satisfies $$\tr(\rho_{\pi}(\Frob_w))=\psi_{(\pi,\RR)}(h_w).$$

We denote by $Z_{\rm reg} \subset Z \isomo \cal{Z}_e$ the subset of points parameterizing the $p$-refined $(\RR,\pi)$ such that $\pi_{\infty}$ is regular, and such that the semisimple conjugacy class of $\pi_p$ (see \S\ref{unramified}) has $m$ distinct eigenvalues. If $z \in Z_{\rm reg}$ parameterizes the regular $p$-refined $\pi$, we will set $$\rhob_z:=\rho_{\pi}.$$

\begin{lemma}\label{zregestzdense} $Z_{\rm reg}$ is a Zariski-dense subspace of $X$ accumulating at each point of $Z$. \end{lemma}

\begin{pf} It is immediate from properties (iv) and (v) of the eigenvariety $X$ (see Theorem \ref{Cheeigen}).
\end{pf}

We end this paragraph by an important example-definition. 

\begin{example}\label{minimaleigen} {\it The minimal eigenvariety containing $\pi$}. 
Let $(\pi,\RR)$ be a $p$-refined automorphic representation of $G$. \par In case (a) (resp. in case (b)), 
assume that $\pi$ is either NMSRPS or unramified (resp. is unramified) at all the finite nonsplit primes.
Define $S$ as the finite set consisting of $p$ and of the primes $l$ such that either $\pi_l$ is ramified or $G(\Q_l)$ 
is the group of invertible elements of a division algebra (which occurs only in case (b)). Define also 
$S_{NM} \subset S$ as the subset of nonsplit primes $l$ such that $\pi_l$ is NMSRPS (this set is empty in case (b)). 
Choose $e=e_S$ such that: \begin{itemize}
\item[(i)] For $l \in S_{NM}$, $e_l=e_{\Sigma_l}$ is a special idempotent attached to the 
Bernstein component attached to the inertial class $\Sigma_l$ of $\pi_l$ (see Example 
\ref{exampleidemp} (v)).
\item[(ii)] For $l=w\bar{w}\neq p \in S\backslash S_{NM}$ such that $G(\Q_l)=\GL_m(\Q_l)$, $e_l=e_{\tau_l}$ where 
$\tau_l$ is the finite dimensional irreducible representation of $\GL_m({\cal O}_{E_w})$ attached to $\pi_l$ 
by Prop. \ref{propquituedeschneider}.
\item[(iii)] For $l=w\bar{w}$ such that $G(\Q_l)$ is the group of invertible elements of a division algebra, $e_l=e_{\tau_l}$ 
where $\tau_l$ is a Bushnell-Kutzko's type for the Bernstein component of $\pi_l$. Such a $\tau_l$ exists by \cite[]{BuKu}.
\end{itemize}
Choose a finite extension $L/\Q_p$ which is sufficiently big so that $\iota_p\iota_{\infty}^{-1}\pi_f$ and each $e_l$ is defined over $L$. 

\begin{definition}\label{minimaleigendef} Under these assumptions, 
the unique eigenvariety over $L$ for $\cal{Z}_{e,NM}$ given by Prop. \ref{existeigenNMSRPS} will be 
referred as {\it the minimal eigenvariety of $G$ containing $\pi$} (or $(\pi,\RR)$).
\end{definition}

Of course, in this context, if we are interested in the variant with the $i_0^{th}$ weight fixed as in Rem. \ref{remcheeigen} (iii), 
we shall always choose the integer $k$ to be the $i_0^{th}$ weight of $\pi_{\infty}$.

\end{example}

\subsubsection{The family of Galois representations on $X$}\label{existfamrepgal}

We adopt also from now on the notations of \S\ref{sectdefrigfam}. 
The first result is that the $\rhob_z$ with $z \in Z_{\rm reg}$ interpolates uniquely to a rigid analytic family of $p$-adic representations of $\G_{E,S}$ on $X$. It uses only properties (P0) and (P1).

\begin{prop}\label{defpseudocaronX} There exists a unique continuous pseudocharacter 
$$T: \G_{E,S} \longrightarrow \OO(X)^{0}$$ such that for all $z \in Z$, $T_z=\tr(\rhob_{z'})$. Moreover:\ps\begin{itemize}
\item[(i)] $T(cgc^{-1})=T(g)\chi(g)^{m-1}$ for each $g \in \G_{E,S}$ (see \S\ref{quad}),\ps
\item[(ii)] for each prime $l=w\bar{w}$ in $S_0$, we have $T(\Frob_w)=\psi(h_w)$. \end{itemize}
\end{prop}

In the statement above, $\chi: \G_{E,S} \longrightarrow \Z_p^*$ is the $p$-adic cyclotomic character. 
Moreover $c$ is the outer complex conjugation (see \S\ref{quad}).

\begin{pf} By property (vi), $\psi(\HH_{\rm ur}) \subset \OO(X)$ is a relatively compact subset, and by Lemma 
\ref{zregestzdense} $Z_{\rm reg}$ is Zariski-dense in $X$. 
The existence and uniqueness of $T$ follows then from \cite[Prop. 7.1.1]{Ch}. 
The equalities in (i) and (ii) hold as $X$ is reduced and as they hold on the Zariski-dense 
subspace $Z_{\rm reg}$ (see Remark \ref{remsectrepm} (i)).
\end{pf}

For $x \in X$, recall that $\OO_x$ is the rigid local ring at $x$, $k(x)$ its residue field and $\overline{k(x)}$ an algebraic closure of $k(x)$.
As $\OO_x$ is reduced and noetherian, its total fraction ring $$\K_x:={\rm Frac}(\OO_x)$$
is a finite product of fields, and we will denote by $\Kb_x$ a (finite) product of algebraic closures of each of those fields. By Taylor's theorem 
\cite[Thm. 1.2]{Tay}, we have then two canonical representations attached to $x$: \begin{itemize}
\item[(a)] $\rhob_x: \G_{E,S} \longrightarrow \GL_m(\overline{k(x)})$, which is the unique (up to isomorphism) continuous semisimple representation with trace 
$T_x: G \longrightarrow \OO_x \longrightarrow k(x)$. 
\item[(b)] $\rhog_x: \G_{E,S} \longrightarrow \GL_m(\Kb_x)$, which is the unique (up to isomorphism) semisimple representation with trace
$T \otimes \K_x: \G_{E,S} \longrightarrow \OO_x \longrightarrow \K_x$.
\end{itemize}

\begin{cor}\label{corrhobxeigen} For each $x \in X$, and for each prime $l=w\bar{w}$ in $S_0$, 
we have $\rhob_x^{\bot} \simeq \rhob_x(m-1)$ and ${\rhog_x}^{\bot} \simeq \rhog_x(m-1)$.
\end{cor}

\subsubsection{Properties of $T$ at the primes $l\neq p$ in $S$.} Let $l\neq p \in S$ and $w$ a prime of $E$ above $l$. We are interested in the restriction to $\W_{E_x} \longrightarrow \G_{E,S}$ of the family $T$. We invite the reader to read first the Appendix \ref{appendixladic} of which we will use the concepts and notations.\ps

\begin{lemma}\label{lemmebienutile} For each $x \in X$ and $s(x)$ a germ of irreducible component at $x$, there exists $z \in |Z_{\rm reg}|$ in the same irreducible component as $x$ such that $$\Ng_{s(x)} \sim_{I_{E_w}} \Nb_z.$$
\end{lemma}

\begin{pf} It follows from Prop. \ref{semicontinuiteWD} (i) and (ii), and the Zariski-density of $|Z_{\rm reg}|$ in $X$.
\end{pf}

Assumption (P3) has the following consequence.

\begin{prop}\label{nonmonatNMSRPS} Assume that 
$l \in S_{NM}$. For each $x \in X$, $\Nb_x=\Ng_x=0$.
\par
\end{prop}

\begin{pf} By assumption (P3), $\Nb_{z}=0$ for each $z \in |Z_{\rm reg}|$, hence we are done by Lemma \ref{lemmebienutile} and Prop. \ref{semicontinuiteWD} (iii).
\end{pf}

Assume now that $l=w\bar{w} \neq p$ splits in $E$ and that $G(\Q_l) \simeq \GL_m(\Q_l)$. Let us fix a $\Qpb$-valued $d$-dimensional Weil-Deligne representation $(r_0,N_0)$ of $E_w$. Assume that the idempotent $$e_l \in \Cc_c(G(\Q_l),\Qpb)\isomo_w \Cc_c(\GL_m(\Q_l),\Qpb)$$ occurring in the definition of $X$ (see \S\ref{notationsrepgal}) has the property that for all the irreducible smooth representations $\pi$ of $\Qpb[G(\Q_l)]$, we have 
\begin{equation}\label{propidemutile}
\, e(\pi) \neq 0\, \,\Rightarrow \, N(\pi) \prec_{I_{\Q_l}} N_0.
\end{equation}

Note that such idempotents exist by Prop. \ref{propquituedeschneider}.

\begin{prop}\label{monauxplacessplit} Assume that $l=w\bar{w} \neq p$ splits in $E$ and that $e_l$ is as above. For each $x \in X$, and each germ $s(x)$ of irreducible component at $x$, then $$\Nb_x \prec \Ng_{s(x)} \prec N_0.$$ 
\end{prop}

\begin{pf} By assumption (P2), we have $\Nb_z \prec N_0$ for all $z \in |Z_{\rm reg}|$. As $\prec_{I_{E_w}}$ implies $\prec$, we conclude by Lemma \ref{lemmebienutile} and Prop. \ref{semicontinuiteWD} (iii).
\end{pf}

\begin{remark}\label{remP2'} For sake of completeness, 
let us consider also the following stronger variant of condition (P2): 
{\it let $l=w\bar{w}' \neq p$ be a prime that splits in $E$, and $(r,N)$ (resp. $(r',N')$) the $\Qpb$-valued Weil-Deligne representation attached to $\pi_w|\det|^{\frac{1-m}{2}}$ (resp. ${\rho_{\pi}}_{|\W_{E_w}}$). If $N \prec_{I_{E_w}} N_0$, then $N' \prec_{I_{E_w}} N_0$}. Under this stronger assertion, the proof of Prop. \ref{monauxplacessplit} shows that we even have $\Nb_x \prec_{I_{E_w}} \Ng_{s(x)} \prec_{I_{E_w}} N$.
\end{remark}
\ps

Let us give another application in a more specific situation. Let us fix $x \in X$ and assume that $\rhob_x$ is irreducible and defined over $k(x)$. 
Let us view $T$ as a continuous pseudocharacter $$T: G_{E,S} \longrightarrow \OO_x,$$
and consider the faithful Cayley-Hamilton algebra\footnote{It actually coincides with the universal Cayley-Hamilton quotient 
(see \S\ref{CHquotient}) of $(\OO_x[\G_{E,S}],T)$ by theorem \ref{structure} (i).}
 $$S:=\OO_x[G_{E,S}]/\Ker T.$$ Then $S \simeq \M_d(\OO_x)$ by Thm. \ref{structure} (i), so 
that $T$ is the trace of a unique (continuous) representation $$\rho: \G_{E_S} \longrightarrow \GL_m(\OO_x).$$
Let $K$ be the total fraction ring of $\OO_x$, then $\rho \otimes K$ is absolutely irreducible as 
$\OO_x[G_{E,S}] \longrightarrow \M_m(\OO_x)$ is surjective. In particular $$\rhog_x \simeq \rho \otimes \overline{K}.$$ 
By Lemma \ref{prolon} and Prop. \ref{grothlmonofamilies}, $\rho$ admits an associated $\OO_x$-valued Weil-Deligne representation, 
say $(r,N)$, $N \in M_d(\OO_x)$.

\begin{cor}\label{monauxplacessplitcasirred} We keep the assumptions of Prop. \ref{monauxplacessplit}. Assume that $$\Nb_x \sim_{I_{E_w}} N_0,$$ then 
$N$ admits a Jordan normal form over $\OO_x$ and $N \sim_{I_{E_w}} N_0$.
\end{cor}

\begin{pf} Note that $N_1 \prec_{I_{E_w}} N_2$ and $N_1 \sim N_2$ imply $N_1 \sim_{I_{E_w}} N_2$. In particular, by Prop. \ref{monauxplacessplit} and the assumption, we get that for each germ of irreducible component $s(x)$ at $x$, 
$$N_0 \sim_{I_{E_w}} \Nb_x \prec_{I_{E_w}} \Ng_{s(x)} \prec N_0,$$
hence all the $\prec$ above are $\sim_{I_{E_w}}$. The corollary follows then from Lemma \ref{nilpoandmatrix} (ii).
\end{pf}

\subsubsection{Properties of $T$ at the prime $v$.}
We are interested in the restriction of the family $T$ to $$\Gal(\Qpb/\Q_p) \longrightarrow_v \G_{E,S}$$ given by the prime $v$ above $p$, as in section 4. For any representation $\rho$ of $G$, we will shortly say that $\rho$ is Hodge-Tate, crystalline etc... if its restriction by the map above is.\ps   
Let $z \in Z_{\rm reg}$ parameterizing the $p$-refined automorphic form $(\pi,\RR)$ of weight $\uk=(k_1,\cdots,k_m)$.
By properties (P4) of $(Rep(m))$, $\rhob_z$ is Hodge-Tate, with Hodge-Tate weights the 
following strictly increasing sequence of integers:
	$$-k_1, \, -k_2 + 1, \, \dots, -k_m+m-1.$$
For convenience, and also in order to fit with the notations of sections 2 and 3 of this book, this shift leads us 
to modify a little the map $\omega$ as follows.  
Let  $$\log_p: \WW \longrightarrow \Hom_{\rm gr}(T_0,\AAA^1).$$
be the map induced the usual $p$-adic logarithm ${\mathbb G}_m \rightarrow \AAA^1$ (killing $p$), and let us identify 
$$\Hom_{\rm gr}(T_0,\AAA^1) \isomo \Hom_{\Qp}({\rm Lie}(T),\AAA^1)  \isomo \AAA^m$$ 
via the diagonalisation $T \isomo (\Qp^*)^m$. Under these identifications, $\log_p: \WW \rightarrow \AAA^m$ associates to 
the character $\chi=(\chi_1,\dots,\chi_m) \in \WW(L)$ the element 
$$\log_p(\chi)=(\cdots,\left(\frac{\partial \chi_i}{\partial \gamma}\right)_{\gamma=1},\cdots) \in L^m.$$ 
In particular, the composition of the embedding $\Z^m \hookrightarrow \WW$ defined in \S\ref{sectdefeigen} with $\log_p$ 
is the natural inclusion $\Z^m \subset \AAA^m$.

\begin{definition}\label{defonkappai} The morphism $\kappa=(\kappa_1,\cdots,\kappa_m): X \longrightarrow \AAA^m$ is the composition of the map 
$\log_p \cdot \omega$ by the affine change of coordinates $$(x_1,\cdots,x_m) \mapsto (-x_1,-x_2+1,\cdots,-x_m+m-1).$$\par
	For each $z \in Z_{\rm reg}$, $\kappa_1(z), \cdots, \kappa_i(z)$ 
is the strictly increasing sequence of Hodge-Tate weights of $\rhob_z$.
\end{definition}

It turns out that this is enough to imply that for each $x \in X(\Qpb)$, 
the Sen polynomial of $\rhob_x$ is $\prod_{i=1}^m(T-\kappa_i(x))$.

\begin{lemma} \label{calculdeskappai} Let $T: \Gal(\Qpb/\Qp) \longrightarrow X$ be any $m$-dimensional 
continuous pseudocharacter on a separted rigid analytic space over $\Q_p$, 
$\kappa=(\kappa_1,\cdots,\kappa_m): X \longrightarrow \AAA^m$ an analytic map, and $Z \subset X(\Qpb)$ a 
Zariski-dense accumulation subset. Assume that for each $z \in Z$, 
the Sen polynomial of $\rhob_z$ is $\prod_{i=1}^m(T-\kappa_i(z))$. \par
Then for each $x \in X(\Qpb)$, the Sen polynomial of 
$\rhob_x$ is $\prod_{i=1}^m(T-\kappa_i(x))$. In particular, $\rhob_x$ is Hodge-Tate whenever 
the $\kappa_i(x)$ are distinct integers.
\end{lemma}

\begin{pf} By replacing $X$ by its normalization $\tilde{X}$ and $Z$ by its inverse image in $\tilde{X}$, we may assume that $X$ is normal and irreducible. 

	Let $\Omega \subset X$ be an open affinoid. Let $g$, $g'$, $Y$, $Y'$ and $\CY$ and $\MM_{\CY}$ be as in Lemma \ref{existrepweak}. For each open affinoid $V \subset \CY$, Sen's theory \cite{sen} attaches to the locally free continuous $\OO(V)$-representation of $\Gal(\Qpb/\Qp)$ on $\MM_{\CY}(V)$ a canonical element $$\varphi_V \in \End_{\OO(V)}(\MM_{\CY}(V))_{\C_p},$$ whose formation commutes with any open affinoid immersion $V' \subset V$. The characteristic polynomial $P_{\varphi,V}$ of each $\varphi_V$ lies in $\OO(V)[T]$, and all of them glue to a single polynomial $P_{\varphi} \in \OO(\CY)[T]$. 
	
Let $S \subset \Omega$ be a Zariski-dense subset. 
Then $g^{-1}(S)$ is Zariski-dense in $Y$ by \cite[Lemme 6.2.8]{Ch}, hence $g^{-1}(S)\cap Y'$ is Zariski-dense in $\CY$. 
Assuming that the conclusion of the lemma holds for all $x$ in $S$, 
then	$P_{\varphi}$ coincides then with $\prod_{i=1}(T-\kappa_i)$ on $\CY^{\rm red}$, so the conclusion of the lemma holds actually on the whole of $\Omega$ (note that $g$ and $g'$ are surjective). In particular, by a connectedness argument it is enough to show that the conclusion of the lemma holds for each $x$ in a single affinoid subdomain of $X$.

	As $Z$ is a Zariski-dense accumulation subset of $X$, there is an affinoid subdomain $\Omega$ of $X$ such that $Z$ is Zariski-dense in $\Omega$. We claim that the conclusion of the lemma holds for each $x \in \Omega$. Indeed, this follows from the previous paragraph for that specific $\Omega$ and for $S=Z \cap \Omega$, and we are done. 
\end{pf}

Let again $z \in Z_{\rm reg}$ parameterizing the $p$-refined automorphic form $(\pi,\RR)$ as above, and set 
$$\Ref_z=(F_1(z)p^{\kappa_1(z)},\dots,F_m(z)p^{\kappa_m(z)}).$$ 
By Definition \ref{definitionFi}, 
$$\Ref_z = \iota_p\iota_{\infty}^{-1}(\RR|p|^{\frac{1-m}{2}})$$ is an accessible refinement of 
$\pi_p|\det|^{\frac{1-m}{2}}$. By properties (P4) and (P5) of $(Rep(m))$, $\rhob_z$ is a crystalline representation 
and $\Ref_z$ is an ordering of the eigenvalues of its crystalline Frobenius. As $z \in Z_{\rm reg}$, all these eigenvalues 
are distinct hence $\Ref_z$ is also the ordered set of Frobenius eigenvalues of a unique refinement of $\rhob_x$ in the sense of \S\ref{refinement}, 
that we will call also $\Ref_z$.

\begin{prop}\label{Tisref} $(X,T,{\kappa_i},\{F_i\},Z_{\rm reg})$ is a refined family in the sense of \S\ref{defrefi}.
\end{prop}

\begin{pf} By Lemma \ref{calculdeskappai} and what we just explained, $(X,T,{\kappa_i},\{F_i\},Z_{\rm reg})$ satisfies properties (i) to (iv) of the 
definition \ref{defrefi} of a refined family. It also satisfies (*) of {\it loc. cit.} as the maps 
$\omega_i$ lift the $\kappa_i$ by Def. \ref{defonkappai}. \par
To prove the property (v) of the definition of a refined family, we need to prove for any $z \in Z$ and any integer $C$ that the set $Z_C$ accumulates at $z$. By property (iv) of Thm~\ref{Cheeigen}, there is an open 
affinoid $\Omega \subset X$ such that $\kappa(\Omega)$ is an open 
affinoid and $\kappa_{|\Omega}$ is finite and surjective when restricted to any irreducible component of $\Omega$. Thus $\kappa(\Omega)$ contains an open ball of center $\kappa(z)$, and there is an integer $N$ such
that $\kappa(\Omega)$ contains the set $Y \subset \Z^d$ of $m$-uples $(k_1,\dots,k_m)$ with $k_1<k_2<\dots <k_m$, and  
$(k_1,\dots,k_m)\equiv \kappa(z) \pmod{(p-1)p^N}$. Set $$Y_C:=\{(k_1,\dots,k_m)\in Y,\ |k_I - k_J|>C,\ \forall I,J \subset\{1,\dots,m\}, 
|I|=|J| \neq 0\}$$ where $k_I:=\sum_{i \in I} k_i$. Then by definition, $Z_C$ contains $\kappa^{-1}(Y_C)$. From the properties of $\kappa$ recalled above and by \cite[Lemma 6.2.8]{Ch},
it is thus enough to prove that $Y_C$ accumulates at $\kappa(z)$, and 
we are reduced to a simple question about the closed unit ball of dimension $m$. Consider the set $$Y'_C =\{(k_1,\dots,k_m) \in Y,\ k_{i+1} - k_i>
m^2(k_i-k_{i-1})+C,\ \forall i \in \{1,\dots,m-1\}.$$ It is obvious 
that $Y'_C$ accumulates at $\kappa(z)$ and is Zariski dense. It is a simple exercise\footnote{To see that,
 take $I$ and $J$ as in the definition of $Y_C$, and define $n+1$
 as the greatest integers that is in $I$ or $J$ but not both.
 We may assume that $n+1 \in I$.  
 Then $k_I-k_J=k_{n+1} + \sum_{l \leq n} \epsilon_l k_l$, with 
 $\epsilon_l \in \{-1,0,1\}$ and $\sum \epsilon_l = -1$, 
 We may write the last sum, up to adding terms of the forms
 $-k_l+k_l$, as $\sum_{l \leq n} \epsilon_l k_l = -k_n +$ a sum of at most 
$(n+1)^2 \leq m^2$  terms of the form $k_l-k_{l-1}$ with $l \leq n$. 
 None of those term is greatest than
 $k_n-k_{n-1}$. Thus if $(k_1,\dots,k_m) \in Y'_C$, we have  
 $|k_I-k_J| \geq k_{n+1}-k_n- m^2 (k_{n}-k_{n-1}) > C$, which proves
 that $(k_1,\dots,k_m) \in Y_C$.} to check that $Y'_C \subset Y_C$, which completes the proof.
\end{pf}

\subsection{The eigenvarieties at the regular, non critical, crystalline points and global refined deformation functors}
\label{applidef}
\label{eigenatregcripts} In this part, we give an application of the techniques and results of this book to study some global deformation rings, as 
we announced in \S\ref{globcons} of section 2. This has some counterparts concerning the geometry of the minimal 
eigenvarieties at the classical, non critical, crystalline points. We will show that those eigenvarieties 
should be smooth at those points and that they are very neatly related to deformation theory. 
By contrast, a much more complicated situation is expected at
reducible, critical points, and this will actually be the main theme of the last sections 8 and 9.

\subsubsection{Some global deformation functors and a general conjecture.} \label{globdefconj}
Let $m\geq 1$ be an integer, $E/\Q$ a quadratic imaginary field, $S$ a finite set of places of $E$, and 
$$\rho: \G_{E,S} \longrightarrow \GL_m(L)=\GL(V)$$
an absolutely irreducible continuous representation with coefficients in the finite extension $L/\Q_p$. We assume 
(see \S\ref{quad}) that
$$V^{\bot} \simeq V(m-1).$$
\noindent Let us assume also that the prime $p=v\bar{v}$ splits in $E$, that $V_p:=V_{|E_v}$ is a crystalline 
representation whose crystalline Frobenius has $m$ distinct eigenvalues in $L^*$, and that the Hodge-Tate weights of $V_p$ are two-by-two distinct.\ps
\smallskip
{\bf Choice:} Let us choose a refinement $\Ref$ of $V_p$ in the sense of \S\ref{refinement}.\ps
\smallskip

We will introduce below a deformation functor of $V$ depending on this choice, but let us first remind some general facts of 
deformation theory. Let $\cal C$ be the category of finite dimensional local $\Q_p$-algebras with fixed residue field 
$L$ that we introduced in \S\ref{deformtrifgdef}, $H$ a topological group and $V$ a finite dimensional continuous $L$-representation of $H$. Following Mazur~\cite{Mazurfermat}, let 
$$\XX_V: \cal C \longrightarrow {\rm Ens}$$
be the deformation functor of the $H$-representation $V$. For any $A \in \cal C$, $\XX_V(A)$ is by definition
the set of isomorphism classes of pairs $(V_A,\pi)$ where $V_A$ is a finite free $A$-module equipped with a continuous 
$A$-linear action of $H$ and $\pi: V_A \otimes_A L \isomo V$ an $H$-equivariant isomorphism. If $V$ is absolutely 
irreducible and if the continuous cohomology group $H^1(H,{\rm ad}(V))$ is finite dimensional, we know from 
\cite[\S10]{Mazurfermat} that $\XX_V$ is prorepresentable by a complete local noetherian ring, with tangent 
space isomorphic to the cohomology group above.

\begin{remark}\label{remdefknonfini} Note that in Mazur's theory \cite{Mazurfermat}, the residue field 
$k$ of the coefficient rings is a finite field. However, everything also applies verbatim when 
$k$ is a finite extension of $\Q_p$  (in which case it is actually even a bit simpler as $k$ automatically lifts as a subfield
of the coefficient rings). The adequate version in this setting of the 
{\it $p$-finiteness condition} of {\it loc. cit.} is the following : 
for any finite dimensional continuous $\Q_p$-representation $U$ of $H$, 
the continuous cohomology group $H^1(H,U)$ is a finite dimensional $\Q_p$-vector space. By 
\cite[Prop. B.2.7]{rubinlivre} and Tate's theorems, this condition holds if $H=G_{E,S}$ or the Galois group of a 
local field.
\end{remark}

Let us denote by $\XX_V$ and $\XX_{V_p}$ the deformation functors associated respectively to the $\G_{E,S}$-representation $V$ and to 
the $\Gp$-representation $V_p$. The choice of any embedding $E \longrightarrow \Qpb$ extending $v$ defines a natural transformation by restriction
	$$\XX_V \longrightarrow \XX_{V_p},$$ 
that is $(V_A,\pi) \mapsto (V_{A\, |E_v},\pi)$. Let us denote again by $\Ref$ the triangulation of 
$D_{\rm rig}(V_p)$ associated to our chosen refinement $\Ref$ by Prop. \ref{ref=tri}. Recall that we defined in \S\ref{deftri} a 
{\it refined deformation functor} 
$$\XX_{V_p,\Ref}: \cal{C} \longrightarrow {\rm Set}$$
of $(V_p,\Ref)$ equipped with a natural transformation $\XX_{V_p,\Ref} \longrightarrow \XX_{V_p}$.
By assumption on $V_p$ and Prop. \ref{subfunctor} and \ref{ref=tri}, $\XX_{V_p,\Ref}$ is actually a subfunctor $\XX_{V_p}$. 

\begin{definition}\label{defXVF} Define two subfunctors $\XX_{V,\Ref}$ and $\XX_{V,f}$ of $\XX_V$ as follows. If $A \in \cal C$, 
say that $(V_A,\pi) \in \XX_{V,\Ref}(A)$ (resp. $\XX_{V,f}(A)$) if, and only if:
\begin{itemize}
\item[(i)] $V_A^{\bot} \simeq V_A(m-1)$,
\item[(ii)] For $w \in S$ not dividing $p$, $V_{A}$ is constant when restricted to $I_{E_w}$, that is
$$V_A \simeq_{I_{E_w}} V \otimes_L A.$$
\item[(iii)] $(V_{A\,|\,E_v},\pi) \in \XX_{V_p,\Ref}(A)$ (resp. $V_{A\,|\,E_v}$ is crystalline).
\end{itemize}
We call $\XX_{V,f}$ the {\it fine deformation functor} of $V$, and $\XX_{V,\Ref}$ the {\it refined deformation functor} of $V$ associated to $\Ref$.
\end{definition}

Recall that the parameter of a triangulation define for each $A$
a morphism $$\delta=(\delta_i): \XX_{V_p,\Ref}(A) \longrightarrow \Hom(\Q_p^*,A^*)^m.$$ 
Here $\Hom$ means {\it continuous group homomorphisms}. In particular, the derivative at $1$ of such a morphism 
is an element of $A^m$, which gives us a morphism 
$$\dkappa: \XX_{V_p,\Ref} \longrightarrow \widehat{{\mathbb G}_a^m}.$$
Denote by 
$\bar{\delta}: \Q_p^* \rightarrow (L^*)^m$ the parameter of $\Ref$ (see \ref{ref=tri}).

\begin{prop}\label{basicfactsReffine} \begin{itemize}
\item[(i)] $\XX_{V,f}$ and $\XX_{V,\Ref}$ are prorepresentable by some local complete noetherian rings.
\item[(ii)] The parameter of a triangulation induces a canonical 
morphism $$\delta: \XX_{V,\Ref} \longrightarrow \widehat{\Hom(\Q_p^*,{\mathbb G}_m^m)}_{\bar{\delta}}.$$ 
\item[(iii)] There is a canonical isomorphism $\XX_{V,f}(L[\varepsilon]) \isomo H^1_f(E,{\rm ad}(V))$.
\end{itemize}
\end{prop}

\begin{pf} In order to prove (i), we have to check that each of the conditions (i), (ii) and (iii) in Definition \ref{defXVF} are 
{\it deformation conditions} in the sense of Mazur \cite[\S19, 23]{Mazurfermat}. \par
	For condition (i), note that as $V$ is absolutely irreducible, a deformation $V_A$ is uniquely determined up to isomorphism by its trace (Serre-Carayol's theorem).
It is then trivial to check conditions (1), (2) and (3) of \S23 of {\it loc.
cit.}\,  for that deformation condition.
For condition (ii), (1) is obvious, (3) follows from Prop.
\ref{proprelrepnilp}\footnote{If $(r,N)$ is the Weil-Deligne representation of
$V_{A\,|E_w}$, note that condition (ii) is equivalent to ask that $N$ is
constant on each isotypic component of the semisimple representation
$r_{|I_{E_w}}$.}, and (2) follows easily from (1) and (3) (see the proof of
Prop. \ref{representri}). \par	
	For condition (iii) in the refined case, it is Prop. \ref{representri}. In the fine case, it follows from Ramakrishna criterion 
(see \cite[\S25, Prop. 1]{Mazurfermat}) and from the fact that the category of crystalline representations is closed under passage to subobjects, quotients, and finite direct sums, by a result 
of Fontaine. That concludes the proof of part (i) of the proposition.  \par
	We already explained assertion (ii) before the statement, and assertion (iii) is
now immediate.
\end{pf}

By imitating assertion (iii) of the definition above, let us set 
$$H^1_{\Ref}(E,{\rm ad}(V)):=\XX_{V,\Ref}(L[\varepsilon]).$$
Recall that we defined in \S\ref{crit} a notion of {\it non critical refinement}. 

\begin{prop}\label{basicfactsRefnoncrit} If $\cal F$ is a non critical refinement of $V_p$, then $\XX_{V,f}$ is a subfunctor of $\XX_{V,\Ref}$. 
If moreover $\Hom_{\Gp}(V_p,V_p(-1))=0$, then \begin{itemize}
\item[(i)] $\XX_{V,f}$ is exactly the subfunctor of $\XX_{V,\Ref}$ defined by the equation $\dkappa=0$,
\item[(ii)] This inclusion induces the following exact sequence on tangent spaces:
$$0 \longrightarrow H^1_f(E,{\rm ad}(V)) \longrightarrow H^1_{\Ref}(E,{\rm ad}(V)) \overset{\partial \kappa(L[\varepsilon])}{\longrightarrow} L^m$$
\end{itemize}
\end{prop}

\begin{pf} The first assertion follows from Prop. \ref{basic}. Point (i) is Theorem \ref{colcrit}, 
and (ii) is then obvious.
\end{pf}

We believe in the following conjectures.

\begin{conj}\label{conjc1c2} \begin{itemize}
\item[(C1)] $\XX_{V,f}$ is a closed point.
\item[(C2)] If $\cal{F}$ is non critical, then $\dkappa$ is an isomorphism. In particular, 
$\XX_{V,\Ref}$ is (formally) smooth of dimension $m$. 
\end{itemize}
\end{conj}
By Prop. \ref{basicfactsReffine} (iii), Conjecture (C1) is actually equivalent to the conjecture 
$BK2(\rho)$ introduced in \S\ref{auxselgroup} (see also \ref{remhypBK}). Let us record this fact in 
the following corollary.

\begin{cor} Conjecture (C1) is equivalent to the conjecture $BK2(\rho)$ (see \ref{auxselgroup}). In particular, 
the Bloch-Kato conjecture implies (C1).
\end{cor}

As a consequence, (C1) is a very "safe" conjecture. In what follows, we will try to provide evidence for (C2) 
and we will relate it to eigenvarieties. In particular, this will shed some light on the expected structure of those 
eigenvarieties in some cases. 

\begin{remark}\label{C2impC1} Assuming that $\Hom_{\Gp}(V_p,V_p(-1))=0$, Prop. \ref{basicfactsRefnoncrit} (ii) shows that (C2) implies (C1). 
As $\rho$ is conjecturally pure, this assumption conjecturally always hold, hence {\it (C2) is conjecturally stronger that (C1)}. As we shall see, the input of eigenvarieties will show that they are actually equivalent.
\end{remark}

\subsubsection{A modular special case}\label{specialcaserhomod} We keep the assumptions on 
\S\ref{globdefconj}. As we want to give examples providing evidence for (C2), we will 
focus from now to some special cases (but still rather general, see Rem. \ref{specialcaseisgeneral}) 
coming from the theory of automorphic forms for which everything we shall need is known. Let us fix a prime $q \neq p$ that splits in $E$, 
as well as another split prime $q' \notin \{q, p\}$ if $m \equiv 0 \bmod 4$ and such that $q'=q$ else.

\begin{lemma}\label{existgrG} There exists a unique unitary group in $m$ variables $G$ attached to $E/\Q$ such that \begin{itemize}
\item[(i)]  $G(\R)$ is compact,  
\item[(ii)] if $l \notin \{q,q'\}$, $G(\Q_l)$ is quasisplit,
\item[(iii)] if $l=q$ or $q'$, $G(\Q_l)$ is the group of invertible elements of a central division algebra over $\Q_l$.
\end{itemize}
\end{lemma}
\begin{pf} This follows from Hasse's principle (see e.g. \cite[(2.2)]{clo}). There is no global obstruction when $m$ is odd 
(\cite[Lemme 2.1]{clo}), and a $\Z/2\Z$-obstruction when $m$ is even. In that case, the local invariant in $\Z/2\Z$ 
of a division algebra is always non zero (see (2.3) of {\it loc. cit.}), and the one at the real place is $(-1)^{m/2}$ by \cite[Lemme 2.2]{clo},
hence the lemma.
\end{pf}

Let $\pi$ be an automorphic representation of $G$. Assume that $\pi$ is not $1$ dimensional and that:
\begin{itemize}
\item[($\pi 1$)] $\pi$ is only ramified at primes that split in $E$,
\item[($\pi 2$)] $\pi_p$ is unramified and its Langlands conjugacy class has $m$ distinct eigenvalues,
\item[($\pi 3$)] $\pi_q$ is supercuspidal.
\end{itemize} 
It is easy to construct such representations. Let us fix some choices of $\iota_{\infty}$ and $\iota_p$ as in \S\ref{basicsettingeigenv}. 
By \cite[Thm. 3.1.3]{HL}, such a $\pi$ admits a strong base change $\pi_E$ to $\GL_m(\AAA_E)$ (it is cuspidal as $\pi$ is not one dimensional).
Moreover, this $\pi_E$ satisfies the assumptions of Harris-Taylor's theorem \cite{HT}, hence by {\it loc. cit.} 
we can attach to this $\pi$ and those embeddings a Galois representation $\rho$ with the following properties: \begin{itemize}
\item[($\rho 1$)] $\rho$ has all the properties of \S\ref{globdefconj}.
\item[($\rho 2$)] $\rho_{|E_w}$ is unramified for each nonsplit place $w$ of $E$, and compatible with $\pi_E|.|^{(m-1)/2}$ 
at all split $w$ (up to Frobenius semi-simplification).
\item[($\rho 3$)] $\rho_{|E_v}$ is crystalline and the characteristic polynomial of its crystalline Frobenius is the same as the one of 
$\iota_p\iota_{\infty}^{-1}L(\pi_p)$.
\end{itemize}

\begin{remark}\label{specialcaseisgeneral} If we believe in Langlands' extension of the Taniyama-Shimura-Weil and Artin conjectures, as well as Arthur-Langlands' yoga of parameters 
(see Appendix A),
any $\rho$ as in \S\ref{globdefconj} which is unramified at nonsplit places, irreducible at $q$, and indecomposable 
at $q'$, should occur this way.
\end{remark}

Under assumptions ($\rho3$) and ($\pi2$), $\iota_p\iota_{\infty}^{-1}$ induces a bijection $\RR \mapsto \Ref$ 
between the refinements of $\pi_p$ in the sense of \S\ref{refinements} and the refinements of $V_p$ in the sense of \S
\ref{refinement}. As $\pi_p$ is tempered by Harris-Taylor's theorem, all its refinements are accessible by Example \ref{exrefi}. 
In particular, to any choice of any refinement $\Ref$ of $V_p$ as in \S\ref{globdefconj} corresponds an accessible refinement of $\pi_p$ and 
vice-versa. For some technical reasons, let us also assume that: \begin{itemize}
\item[($\rho_4$)] $\Ref$ is non critical and regular (see Def. \ref{reg}).
\item[($\rho_5$)] For $i=1,\cdots,m$, $\Lambda^i \rho$ is irreducible.
\end{itemize}

All these assumptions being done, we can consider the minimal eigenvariety associated to $(\pi,\Ref)$ 
(see Example \ref{minimaleigen}). Let $z \in X$ be the $L$-point parameterizing $\pi$ equipped with its refinement $\Ref$, and set 
$$\mathbb T:=\widehat{\OO_z}.$$

\subsubsection{$R=T$ at the regular non critical crystalline points of minimal eigenvarieties} Assume that $\rho$ and 
$\pi$ are as in \S\ref{globdefconj} and \S\ref{specialcaserhomod}. Let $R_{\rho,\Ref}$ be the 
universal deformation ring of the refined deformation functor $\XX_{V,\Ref}$ given 
by Prop. \ref{basicfactsReffine} (i). 

\begin{prop}\label{mapRtoT}There is a canonical commutative diagram
\begin{eqnarray} 
\xymatrix{ R_{\rho,\Ref} \ar@{->>}[r]  & {\mathbb T} \\  L[[t_1,\cdots,t_m]] \ar[u]^{\dkappa^{\sharp}} \ar@{=}[r] & \widehat{\OO_{\kappa(z)}} \ar[u]^{\kappa^{\sharp}}  }
\end{eqnarray}
Moreover $\mathbb T$ is equidimensional of dimension $m$ and $\kappa^{\sharp}$ is a finite injective map. 
\end{prop}

\begin{pf} We claim first the existence of a natural map $R_{\rho,\Ref} \rightarrow {\mathbb T}$. Let $A=\OO_z$. As $\rho$ is 
absolutely irreducible, $T$ is the trace of a unique continuous representation 
	$$\rho_A: \G_{E,S} \longrightarrow \GL_m(A),$$
hence for each cofinite length proper ideal $I$ of $A$ we have a canonical element $\rho_A \otimes A/I \in X_V(A/I)$ (note that 
$A/I={\mathbb T}/I{\mathbb T}$). We have to show that this element falls in $X_{V,\Ref}(A/I)$, {\it i.e.} to check conditions (i) 
to (iii) in Def. \ref{defXVF}. Condition (i) follows at once from the fact that $T^{\bot}=T(m-1)$ and that $\rho$ is absolutely irreducible.
Condition (ii) follows from Cor. \ref{monauxplacessplitcasirred}, which applies by \cite{TY}. Last but not least, 
condition (iii) follows from Theorem~\ref{tridefirr} if we can check the assumptions of \S\ref{hypredloc} at the point $z$. They hold as
$(X,T,\kappa,\{F_i\})$ is a refined family by Prop. \ref{Tisref}, and as the assumptions (REG), 
(NCR) and (MF') of \S\ref{hypredloc} follow form ($\rho 4$) and ($\rho 5$). This concludes the claim. \par
	The existence of a commutative diagram as in the statement is moreover given by the 
identification of the parameter $\delta\otimes A/I$ in the statement of Theorem \ref{tridefirr}. More precisely, that theorem 
furnishes a commutative diagram
\begin{eqnarray}\label{commutedeltasurj}
\xymatrix{ R_{\rho,\Ref} \ar@{->}[rr]  \ar@{<-}_{\delta}[dr] & & {\mathbb T} \ar@{<-}[dl] \\  & \widehat{\OO_{\overline{\delta}}} &}
\end{eqnarray}
where $\overline{\delta} \in \Hom(\Qp^*,{\mathbb G}_m)$ is the parameter of $\Ref$, and where the map on the right is the natural map.
\par The assertion on $\mathbb T$ and 
$\kappa^{\sharp}$ follow from property (iv) of the eigenvariety $X$, thus it only remains to check that the upper map is a surjection. By property (i) of eigenvarieties (see Def. \ref{eigendef}), $\mathbb T$
is generated by $\HH$ as an $\widehat{\OO}_{\kappa(z)}$-algebra. As $$\HH=\ATL
\otimes \HH_{\rm ur}$$ we see that
$\mathbb T$ is generated over $\widehat{\OO}_{\kappa(z)}[F_1,\cdots,F_m]$ by
the $T(\Frob_w)$'s for the primes $l=w\bar{w} \in S_0$. But each $T(\Frob_w)$ obviously lifts in $R_{\rho,\Ref}$ as the trace of 
$\Frob_w$ in the universal refined deformation, and we are done by the commutative diagram 
(\ref{commutedeltasurj}).
\end{pf}

\begin{cor}\label{cormapRtoT} Assume that $\rho$ is associated to a $\pi$ as in \S\ref{specialcaserhomod}. \begin{itemize}
\item[(i)] Conjecture (C1) is equivalent to conjecture (C2). In particular, if the Bloch-Kato conjecture holds then (C2) holds.
\item[(ii)] Conjecture (C1) implies that all the maps of the diagram of Prop. \ref{mapRtoT} are isomorphisms.
\end{itemize}
\end{cor}

\begin{pf} By Harris-Taylor's theorem (especially property (P5) of \ref{repm}), $\rho$ is pure and 
$\Hom_{\Gp}(V_p,V_p(-1))=0$, hence Prop. \ref{basicfactsRefnoncrit} applies. In particular, (C2) implies (C1) by Remark \ref{C2impC1}. \par
	Assume that (C1) holds. We claim that $\partial\kappa^{\sharp}$ is an isomorphism. It would imply (C2) as well as assertion (ii) of the corollary, since the top arrow (resp. the right arrow) in Prop. \ref{mapRtoT} is surjective (resp. injective), so it is enough to prove the claim. \par 
		By Prop. \ref{mapRtoT}, $\partial\kappa^{\sharp}$ is injective. As it induces an isomorphism on the residue fields $L$ and as $R_{\rho,\Ref}$ is a complete local noetherian ring, it is enough to show that $\partial\kappa^{\sharp}$ induces an isomorphism on tangent spaces. Under conjecture (C1), the exact sequence of Prop. \ref{basicfactsRefnoncrit} (ii) shows that the tangent space of $R_{\rho,\Ref}$ has dimension $\leq m$, and that is enough to know that it has dimension exactly $m$ in order to conclude. But by Prop. \ref{mapRtoT} $R_{\rho,\Ref}$ has Krull dimension $\geq m$, as its quotient $\mathbb T$ has dimension $m$. The dimension theory of local noetherian rings shows then that the Krull dimension of $R_{\rho,\Ref}$ and the dimension of its tangent space both coincide with $m$. Thus $R_{\rho,\Ref}$ is regular of dimension $m$ and we are done.
\end{pf}

In particular, the map $R_{\rho,\Ref} \longrightarrow \mathbb T$ that we defined should always be an isomorphism.
Moreover, we also get that a far reaching infinitesimal version of the principle {\it a non critical 
slope form is classical} should hold for eigenvarieties. 

\begin{conj}\label{conjR=TetCRIT} \begin{itemize}
\item[(R=T)] The map $R_{\rho,\Ref} \longrightarrow \mathbb T$ is an isomorphism.
\item[(CRIT)] The map $\kappa^{\sharp}$ is an isomorphism, {\it i.e.} the weight map $\kappa$ is \'etale at $z$.
\end{itemize}
\end{conj}

Of course, the natural trend in the area since the work of Wiles and Taylor-Wiles is that we should first try to prove conjectures 
(R=T) and (CRIT) using the cohomology of Shimura varieties and the theory of automorphic forms. Note that (R=T) and (CRIT) together imply (C2) by Prop. \ref{mapRtoT}. Then (C1) would follow by Cor. \ref{cormapRtoT}.

\begin{cor}\label{R=TimpliesC} Conjectures (R=T) and (CRIT) imply (C1) and (C2).
\end{cor}

The deepest part there is certainly to show (R=T), but we will not say more here about that conjecture 
(see Kisin's paper \cite{kis2} in the case $m=2$, as well as the discussion in \S\ref{auxselgroup}). In the remark below, we discuss instead what is known about conjecture (CRIT).

\begin{remark} \label{remconjcrit}
\begin{itemize}
\item[(i)] Assuming that $\Ref$ is moreover $U$-non critical (see \ref{remcrit} (iii) or \cite[\S7.5.1]{Ch})), then by Prop. \ref{evalxpil} (iii) - that is essentially the {\it small slope forms are classical} result - the refined automorphic representation $(\pi,\RR)$ does not have any infinitesimal deformation in the space of $p$-adic automorphic forms. \par 
	This fails short to imply that $\kappa$ is \'etale because of a subtlety: for these general $G$ we do not have a good control of $X$ in terms of the spaces of $p$-adic automorphic forms, for instance like the pairing we have for $\GL_2/\Q$ given by the $q$-expansion (see {\it e.g.} \cite[Prop. 1 (c)]{lissite}). However, if we knew that the multiplicity one theorem holds for the automorphic representations of $G$ which are unramified at the nonsplit primes (which is expected), then (CRIT) would follow easily\footnote{Basically because when an $A$-module $M$ is free of rank one over $A$, then each commutative $A$-subalgebra of $\End_A(M)$ is equal to $A$ (hence etale over $A$)...}  from Prop. \ref{evalxpil} (iii). In particular, by results of Rogaswki we know (CRIT) in the $U$-non critical case when $m\leq 3$.
\item[(ii)] Even admitting this multiplicity one result for $G$, it would be very interesting to have a proof of conjecture (CRIT) in general. As explained in Rem. \ref{remcrit} (ii), in the 
classical case of the eigencurve and the group $\GL_2/\Q$, the full case of (CRIT) is known and is quite deep: it 
follows from Coleman's theorem \cite{col1} (including the so 
called {\it boundary case} where $v(a_p)=k-1$) and from the existence of overconvergent companion forms due to 
Breuil-Emerton \cite{BrEm}. 
\item[(iii)] Let us mention that we certainly believe that conjectures (R=T) and (CRIT) also hold without the assumption that 
$\Ref$ is regular or that the $\Lambda^i\rho$ are irreducible for $i>1$ (but of course still non critical). However, it seems to be an interesting
problem to understand the case where the crystalline Frobenius of $V_p$ are not assumed to be two-by-two distinct any more 
(or, which is the same, when ($\pi 2$) does not hold). Indeed, there seems to be no trivial way to make refinements of $\pi_p$ 
and $V_p$ correspond without this assumption. It is quite clear that it in this case our definition of a refinement of $\pi_p$ is not 
the good one (it is somehow an "{\it archimedean one}" rather than $p$-adic). From the point of view of the hypothetical
Breuil-Langlands correspondence between $m$-dimensional crystalline representations of $\G_p$ 
and $p$-adic representations of $\GL_m(\Q_p)$, we should rather define those refinements in terms of 
the closure (in the space of $p$-adic automorphic forms for $G$) of the representation $\pi_p \otimes \pi_{\infty}^*$ and show 
that they match (in the dictionary) with the refinements of $V_p$. In the case of $\GL_2/\Q$, this appears 
quite naturally in Emerton's theory of Jacquet modules and it is actually a theorem \cite{Emgl2}.
\item[(iv)] If we do not assume that $\Ref$ is critical (but say all the other hypotheses), the situation is actually very interesting but 
quite different. For example, the can be no map from $R_{\rho,\Ref}$ to $\mathbb T$ as the refined family $(T,X)$ is generically non critical.
We postpone its study to a subsequent work.
\end{itemize}
\end{remark}

\subsection{An application to irreducibility}

\label{parappliirr}

A simple application of the results of 
this section, together with the generic irreducibility theorem \ref{redgen},
is the existence of many $n$-dimensional Galois representations of $G_E$ 
that are irreducible, even after restriction to a decomposition 
group at a place $v$ above $p$. As an example, let us prove the following result :

\begin{theorem}\label{thmexampleirred} Assume Rep$(m)$ (actually for that matter we may release 
conditions (P2) and (P3)).  Then for any integer $C$, there exists an automorphic 
representation $\pi$ for $\U(m)$  such that the Galois 
representation $\rho_\pi$  \begin{itemize}
\item[-] is unramified at every places not dividing $p$, 
\item[-] is crystalline and irreducible at each of the two places dividing $p$, 
\item[-] and has Hodge-Tate weights $k_i$, $i=1,\dots,m$, such that $|k_i-k_j|>C$ for every $i \neq j$.
\end{itemize}
\end{theorem}

\begin{pf} Let $\pi$ be the trivial representation of $\U(m)$. 
It is unramified at all the finite primes hence we may consider the minimal eigenvariety containing $\pi$ as in \S\ref{minimaleigen}. 
By definition it is the eigenvariety $(X/\Q_p,\psi,\nu,Z)$ for the set $Z$ of $p$-refined automorphic representations $(\psi_{(\pi,\RR)},\uk)$ such that $\pi$ is unramified at all finite places and $\RR$ is an accessible refinement of $\pi_p$. For that reason, we may wish to 
call it the {\it unramified}, or {\it tame level 1}, eigenvariety for $\U(m)$. Note that $X$ 
is not empty since it 
contains the point $x_0$ corresponding to the trivial representation $\pi$ of $\U(m)$ with its unique refinement at $p$ (see Ex. \ref{exrefi} (ii)).\par
	By an argument similar than the one given in the proof of Prop.~\ref{Tisref}, there exists a point 
$x_1 \in Z$ arbitrary close to $x_0$ that corresponds
to an automorphic representation $\pi_1$ (together with a refinement $\RR_1$)
such that:\begin{itemize}
\item[(i)] $(\pi_1)_\infty$ is regular,
\item[(ii)]  the eigenvalues $\lambda_1,\cdots,\lambda_m$ 
of the Langlands conjugacy class of $(\pi_1)_p$ are distinct, and  no quotient of two of them 
is equal to $p$, 
\item[(iii)] if $I,J \subset \{1,\cdots,m\}$ are such that $|I|=|J|$, 
then $\prod_{i \in I}\lambda_i=\prod_{j \in J}\lambda_j \Rightarrow I=J$.
\end{itemize}
Condition (ii) ensures that $(\pi_1)_p$ is a full
irreducible unramified principal series, hence that 
all its refinements are accessible (cf. Prop.~\ref{casgenerique}). \par

We now use the refined family $T$ of Galois representations on $X$ 
constructed in this section. The representation $\rhob_{x_1}$ 
corresponding to $\pi_1$ has distinct Frobenius eigenvalues at $v$ 
and distinct Hodge-Tate weights by (i). In particular, we may write $\rhob_{x_1}$ as a sum $\rhob_1 \oplus \dots \oplus \rhob_r$ of non-isomorphic irreducible representations.
As explained in \S\ref{permutation}, there is a partition $\{1,\dots,m\} =W_1\coprod \dots \coprod W_r$, such that $|W_i|=d_i$ for $i=1,\dots,r$, defined by setting $W_i$ equal to the set of $j$'s such that $\kappa_j(x_1)$ is a weight of $\rhob_j$.
Moreover, as explained there, to a refinement $\cal F$ (that is, 
an ordering of the crystalline  Frobenius eigenvalues) of $\rhob_{x_1}$
is attached a second partition $\{1,\dots,m\}=R_1(\cal F) \coprod \dots \coprod R_r(\cal F)$ with $|R_i(\cal F)|=d_i$, 
and all partitions of this type are attached to some 
refinement $\cal F$. It is a simple combinatorial task, tackled in the next lemma, to see that 
there is always one (and in general, many) partition $(R_i)$ of this type 
which is ``orthogonal'' to the partition $(W_i)$, in the sense that
$$\forall {\cal P} \subset \{1,\dots,r\},\ 0<|\cal P|<r,\ \coprod_{i \in \cal P} W_i \neq \coprod_{i \in \cal P} R_i.$$
We choose such a partition, and a refinement $\cal F$ of $\rhob_{x_1}$ 
defining that partition. 

To $\cal F$ corresponds a refinement $\RR_2$ of $(\pi_1)_p$, necessarily 
accessible. We thus may define a point $x_2$ 
of $Z \subset X$ corresponding to $(\pi_1,\RR_2)$.
Note that $\cal F$ is regular by property (iii) above of $x_1$.
By Thm~\ref{redgen}, and the properties of $\RR_2$,
the family of Galois representation $T$ restricted to $D_v$ 
is generically (absolutely) irreducible near $x_2$. Hence for any integer $C$,
there is a point 
$x_3 \in Z_C$, such that $(\rhob_{x_3})_{|D_v}$ is crystalline
absolutely irreducible. So is $(\rhob_{x_3})_{|D_{\bar v}}$, 
since it is the dual of the preceding representation. 

It is clear that the automorphic representation $\pi_3$ corresponding 
to $x_3$ satisfies all assertions of the theorem.
\end{pf}
\begin{remark} The trivial representation of $\U(m)$ has a unique refinement,
and it turns out that this refinement does not allow us to conclude using 
Theorem~\ref{redgen} that the deformation it defines 
of the trivial representation is 
generically irreducible restricted to $D_v$ (and likely it is 
not for $m \geq 3$). 
Indeed, as explained in \S\ref{permutation}, since the trivial 
representation is ordinary in the sense of {\it loc. cit.}, 
this unique refinement is characterized by a permutation 
$\sigma$ of $\{1,\dots,m\}$. Actually for the trivial representation we 
have $\sigma(1)=m,\dots,\sigma(m)=1$ and we see that, for $m \geq 3$, 
$\sigma$ is not transitive, or in the language of {\it loc. cit.} the refinement is not {\it anti-ordinary} (it is not ordinary either).
That is why we had to process in two steps in the proof above. 
\end{remark}
We now prove the combinatorial lemma needed in the above proof.
\begin{lemma} For every partition
$\{1,\dots,m\} =W_1\coprod \dots \coprod W_r$, with $|W_i|=d_i$ for all $i$,
 there exists a partition $\{1,\dots,m\} =R_1\coprod \dots \coprod R_r$, 
with $|R_i|=d_i$ for all $i$ such that
$$\forall {\cal P} \subset \{1,\dots,r\},\ 0<|\cal P|<r,\ \coprod_{i \in \cal P} W_i \neq \coprod_{i \in \cal P} R_i.$$
\end{lemma}
\begin{pf} Pick up an element $t_i$ in each $W_i$. Choose a transitive permutation $\sigma$ of $\{1,\dots,r\}$. Put $t_{\sigma(i)}$ in $R_{i}$. 
Complete the construction of $R_i$ as you like. Then for all 
$\cal P \subset \{1,\dots,d\}$,  $\coprod_{i \in \cal P} W_i$ 
contains $t_i$ for $i \in \cal P$ and 
no other $t_j$, while  $\coprod_{i \in \cal P} R_i$ contains 
the $t_i$ for $i \in \sigma(\cal P)$ and no other $t_j$. 
Hence if those two unions are equal, $\sigma(\cal P)=\cal P$.
\end{pf} 
\begin{remark}
\label{remarkneededinsection8}
The proof above may be adapted to prove the existence of partitions $(R_i)$ satisfying further properties. For example, we shall need in a remark of chapter 8
to deal with a case where $d_1=d_r=1$, $W_r=\{k\}$, $W_1=\{k+1\}$, 
with $1\leq k<k+1 \leq m$; in this case we want a partition $(R_i)$, satisfying the above properties and moreover $R_1=\{1\}$ and $R_r=\{m\}$. 
It is certainly clear
for the reader how the above proof has to be adapted to prove the existence of such $R_i$.
\end{remark}

\newpage

\subsection{Appendix: $p$-adic families of Galois representations of $\Gal(\Qb_l/\Q_l)$ with $l\neq p$}\label{appendixladic}
\newcommand{\ut}{\underline{t}}
\subsubsection{Some preliminary lemmas on nilpotent matrices}\label{lemmasnilpo} Let $k$ be a field and $n \in M_d(k)$ a 
nilpotent matrix. By Jordan's normal form theorem, there exists a unique unordered partition of $\{1,\dots,d\}$ 
	$$\ut(n):=(t_1\geq t_2 \geq \cdots ),\, \, \,  t_i \in \mathbb{N}, \, \, \sum_i t_i=d,$$ 
such that $n$ is conjugate in $M_d(k)$ to the direct sum of Jordan's blocks\footnote{For us, the Jordan block $J_d \in M_d(A)$ for any commutative ring $A$ 
is the matrix of the endomorphism $n$ of $A^d$ defined by $n(e_1)=0$ and $n(e_i)=e_{i-1}$ if $i>1$.} $$J_{t_1}\oplus J_{t_2} \oplus \cdots J_{t_s},$$ where 
$s$ is the smallest integer such that $t_{s+1}=0$. If $k \longrightarrow k'$ is a field embedding, then $t(n\otimes_k k')=t(n)$.\ps
	Recall that the {\it dominance ordering} on the set of decreasing sequences $\ut=(t_1\geq t_2 \geq \cdots)$ of integers is 
the partial ordering $$ \ut \prec \ut'\, \, \, \Leftrightarrow \, \, \, \forall i\geq 1,\, \,  t_1+\cdots+t_i \leq t'_1+\cdots+t'_i.$$
We refer to \cite[\S I]{macdo} for its basic properties.

\begin{prop}[Gerstenhaber]\label{gerstthm} Let $n, n' \in M_d(k)$ be two nilpotent matrices. Then the following assertions are 
equivalent:
\begin{itemize}
\item[(i)] $n$ is in the Zariski-closure of the conjugacy class of $n'$ in $M_d(\bar k)$,
\item[(ii)] For all $i\geq 1$, ${\rm rank}\,  n^i \, \leq \, {\rm rank} \, {n'}^i$,
\item[(iii)] $t(n)\prec t(n')$.
\end{itemize}
\end{prop}

\begin{pf} The equivalence between (i) and (ii) is \cite[Thm. 1.7]{Ger}. Assertion (ii) is equivalent to ask that for each $i$, 
${\rm dim}({\rm ker}(n^i)) \geq {\rm dim}({\rm ker}({n'}^i))$, which is another way to say that $t^*(n) \succ t^*(n')$. Here, 
$\ut^*$ is the conjugate partition of $\ut$, and the result follows as $\ut \prec \ut' \Leftrightarrow \ut^* \succ {\ut'}^*$ by 
\cite[\S 1.11]{macdo}.
\end{pf}

\begin{definition}\label{defordonnilp} If $n, n'$ are any two nilpotent matrices\footnote{It is not necessary to ask, in this definition, 
that they have the same coefficient field (or even the same size).}, we write $n \prec n'$ (resp. $n \sim n'$) if $t(n) \prec t(n')$ (resp. $t(n)=t(n')$). If $n, n' \in M_d(k)$, $n \prec n'$ if, and only if, $n$ is in the 
Zariski-closure of the conjugacy class of $n' \in M_d(\bar k)$.
\end{definition}

\begin{cor} \label{ineqnilp} Let $V$ be a finite dimensional $k$-vector space and $n \in \End_k(V)$ a nilpotent element. Let $U \subset V$ 
be a $k[n]$-submodule and $n'$ is the endomorphism induced by $n$ on $U\oplus V/U$, then $n' \prec n$.
\end{cor}
\begin{pf} It is clear on condition (ii) of Prop. \ref{gerstthm}.
\end{pf}

Let us collect now some useful results about nilpotent matrices with coefficients in a ring.
Let $A$ be a commutative ring and $n \in M_d(A)$ a nilpotent matrix. 
We will say that $n$ {\it admits a Jordan normal form over $A$} if $n$ is 
$\GL_d(A)$-conjugate in $M_d(A)$ to a direct sum of Jordan's blocks 
$J_{t_1} \oplus \cdots J_{t_s}$ for some unordered 
partition $(t_1\geq t_2 \geq \cdots)$ of $\{1,\cdots,d\}$ as above. Again, 
we see by reducing modulo any maximal ideal of $m$ 
that if such a Jordan normal form exists, then the associated partition is
unique. The following proposition is probably well known.

\begin{prop}\label{normalformoverA} Let $A$ be a local ring and $n \in M_d(A)$ a nilpotent matrix. The following properties are equivalent:
\begin{itemize}
\item[(i)] $n$ admits a Jordan normal form over $A$,
\item[(ii)] for some faithfully flat commutative $A$-algebra $B$, the image of $n$ in
$M_d(B)$ admits a Jordan normal form over $B$,
\item[(iii)] for each integer $i\geq 1$, the submodule $n^i(A^d) \subset A^d$ is free over $A$ and direct summand.

\end{itemize} 
\end{prop}

\begin{pf} It is clear that (i) implies (ii), and also that (i)
implies (iii) even if we do not assume $A$ to be local. Note that
for $u \in \End_A(A^d)$ any element and $B$ a faithfully flat $A$-algebra, then
$$Im(u)\otimes_A B \isomo Im(u\otimes_A B),$$ 
and the latter is projective and direct summand in $B^n$ as a $B$-module if, and
only if, $\Im(u) \subset A^n$ has those properties as $A$-module. As a
consequence, (ii) implies that $n^i(A^d)$ are projective and direct summand $A$-modules, hence free as $A$
is local, hence (ii) implies (iii). \par
	It only remains to show that (iii) implies (i), for which we argue
as in the classical proof of Jordan's theorem. For $i\geq 0$ let 
$N_i:=\Im(n^i) \in A^d$. We construct by descending induction some $A[n]$-submodules $F_{i+1}$ and $Q_i$ of
$A^d$ for $i=d-1,\dots,0$, such that \begin{itemize}
\item[-] $F_d=0$,
\item[-] $F_{i+1}$ and $Q_i$ are free and direct summand as $A$-modules,
\item[-] $n_{|Q_i}$ has a Jordan normal form over $A$ with blocks of size $i+1$ (if any), 
\item[-] $F_i=F_{i+1} \oplus Q_i$ and $n^i(F_{i+1}) \oplus n^i(Q_i) = N_i$.
\end{itemize}
Assume that $F_j$ and $Q_j$ are constructed for $j>i$, we have to define $Q_i$.
Note that $F_{i+1}$ is free and direct summand as $A$-module, and that
$n_{|F_{i+1}}$ admits a Jordan normal form, so $n^i(F_{i+1})$ and 
$$K_i:=\ker n \cap n^i(F_{i+1})$$ are free and
direct summand as an $A$-module as well. In
particular, $K_i$ is a direct summand of $\ker n \cap N_i$, by Lemma \ref{lemmeAA'} (i) below. \par 
	As $A$ is local, we may then find a finite free $A$-module $Q'_i \subset \ker n \cap N_i$ which is a complement to $K_i$, it satisfies:  \begin{equation}\label{intrducdeKi} K_i \oplus Q'_i = \ker n \cap N_i, \, \, \, \, Q'_i\cap n^i(F_{i+1})=0.\end{equation} We claim that 
\begin{equation}\label{claimqprimei}n^i(F_{i+1}) \oplus Q'_i = N_i.\end{equation} Indeed, $n(N_i)=N_{i+1}=n^{i+1}(F_{i+1})$ implies that 
$N_i=n^i(F_{i+1})+\ker n \cap N_i$, which proves the claim by (\ref{intrducdeKi}). Note also that 
\begin{equation}\label{secondclaim} Q'_i \cap F_{i+1} = 0. \end{equation}
Indeed, the Jordan blocks of $Q_j$ for $j>i$ have size $> i+1$, thus $\Ker n \cap F_{i+1} \subset n^i(F_{i+1})$. As $n(Q'_i)=0$, we get that $Q'_i \cap F_{i+1} \subset Q'_i \cap n^{i}(F_{i+1})=0$ by (\ref{intrducdeKi}), and we are done. \par
We can now conclude the proof. If $Q'_i=0$, then we set $Q_i=0$ and we are done. If else,
we may choose
$v_1,\cdots,v_r$ in $A^d$ such that $n^i(v_1),\cdots,n^i(v_r)$ is an
$A$-basis of $Q'_i$. Set $$Q_i:=A[n]v_1+\cdots +A[n]v_r.$$
Note that $n^{i+1}(Q_i)=n(Q'_i)=0$. We check at once by applying $n$
several times that \begin{itemize}
\item[-] the $n^s(v_j)$ with $0\leq s \leq i$ and $1\leq j \leq
r$ are an $A$-basis of $Q_i$ (so in particular,
$n_{|Q_i}$ admits a Jordan normal form $J_{i+1} \oplus J_{i+1} \oplus \cdots \oplus
J_{i+1}$ ($r$ times)),
\item[-] $Q_i \cap F_{i+1}=0$ (note that $Q'_i\cap F_{i+1}=0$ by (\ref{secondclaim})). \end{itemize}
and we are done by (\ref{claimqprimei}) if we set $F_i:=F_{i+1}\oplus Q_i$.
\end{pf}

Criterion $(i)\Leftrightarrow (ii)$ of Prop. \ref{normalformoverA} shows that the property of
admitting a Jordan normal form (say over a local ring\footnote{In the general
case, the well-behaved definition for "admitting a Jordan normal form over
$A$" is certainly to ask directly that the $n^i(A^d)$ are
projective and direct summand.}) is invariant under faithfully flat base change. When we deal with
deformation theory, the following other kind of descent is useful. 

\begin{prop}\label{proprelrepnilp} Assume that $A \longrightarrow A'$ is a local homomorphism between 
artinian local rings inducing an isomorphism on the residue fields,
and let $n \in M_d(A)$ be a nilpotent matrix. Then $n$ admits a Jordan normal 
form over $A$ if, and only if, its image in $M_d(A')$ admits a Jordan normal form over $A'$.
\end{prop}

\begin{pf} Assume that the image of $n$ in $M_d(A')$ has a Jordan normal form (the other implication is obvious). For $i\geq 1$, let $N_i:=Im(n^i) \subset A^d$ and
$N'_i=Im((n\otimes_A A')^i) \subset A'^d$. As $A \subset A'$, we have 
	$$N_i \subset A'.N_i=N'_i.$$ 
Recall that $A$ and $A'$ have the same residue field $k:=A/m$, and let $\nb
\in M_d(k)$ be the image of $n$. We have two natural surjections with the same image 
\begin{equation}\label{twonatsurj} N_i \otimes_A k \longrightarrow
\Im(\nb^i),\, \, \, N'_i \otimes_{A'} k \isomo \Im(\nb^i),\end{equation}  
where the second map is an isomorphism as $N'_i\subset A'^d$ is a direct summand.
Let $\bar{v_1},\cdots,\bar{v_r}$ a $k$-basis of $\nb^i(k^d)$, and
$v_1,\cdots,v_r \in N_i$ some liftings of the $\bar{v_j}$. Set $$P_i:=\sum_{j=1}^r
A v_j \subset N_i.$$
(If $\nb^i=0$ then we set $P_i=0$). This is a (free)
direct summand of $A^d$ by Lemma \ref{famlibremodlibre} below.
Moreover, the $v_j$ generate $N'_i$ over $A'$ by (\ref{twonatsurj}) and
Nakayama's lemma, hence $P_iA'=N_iA'=N'_i$. By Lemma \ref{lemmeAA'} (ii) below,
this implies that $P_i=N_i$, thus $N_i$ is free and direct summand, and we
are done by Prop. \ref{normalformoverA}.
\end{pf}

\begin{lemma}\label{famlibremodlibre} Let $A$ be a local ring with residue field $k$. 
If some elements $v_1,\dots,v_p$ in $A^d$ have $k$-independent images in
$k^d$, then they are $A$-independent and $$\oplus_{i=1}^p A v_i \subset A^d$$
is a direct summand.
\end{lemma}
\begin{pf} Let $M \in M_{d,p}(A)$ be the matrix defining the $v_i$ in the canonical basis. By assumption, some 
$p \times p$-minor of $\bar{M} \in M_{d,p}(k)$ is nonzero, hence the same 
$p \times p$-minor of $M$ is in $A^*$, and the $v_i$ are $A$-independent. We conclude by completing the $\bar{v_i}$ 
in a basis of $k^d$ and by Nakayama's lemma.
\end{pf}

\begin{lemma}\label{lemmeAA'} Let $A$ be a commutative ring, $P \subset
N \subset M$ an inclusion of $A$-modules such that $P$ and $M$ are projective,
and $P$ is a direct summand of $M$. \begin{itemize}
\item[(i)] $P$ is a direct summand of $N$.
\item[(ii)] Assume that $A \longrightarrow B$ is an injective ring homomorphism, then
$M \longrightarrow M \otimes_A B$ is injective. If $B.P=B.N$ in $M\otimes_A B$,
then $P=N$.
\end{itemize}
\end{lemma}
\begin{pf} We can write $M=P \oplus P'$ for some $A$-module $P' \subset M$.
It is immediate to check that $N = P \oplus (P' \cap N)$, which proves
(i). To check (ii), we may assume that $M$ is a free $A$-module, and the first
assertion is then obvious. Assuming that $B.P=B.N$, we have to show by (i)
that $P'\cap N=0$. But 
$$(P'\cap N) \subset B.(P'\cap N) \subset B.P' \cap B.N= B.P' \cap B.P =0,$$ 
which concludes the proof.
\end{pf}

When a nilpotent element in $M_d(A)$ (or even in a GMA over $A$) does not
necessarily admit a Jordan normal form there are still some inequalities between
the generic and residual partitions that are satisfied.

\begin{prop}\label{nilpoandgma} Let $A$ be a commutative reduced local ring whose total fraction ring is a finite product of fields\footnote{See Prop. \ref{propfinitefrac} for a discussion of this assumption.} $K=\prod_s K_s$, 
and let $k$ be its residue field. Let $R \subset M_d(K)$ be a standard GMA of type $(d_1,\dots,d_r)$ (see Example \ref{gma}). Assume that the natural surjective map 
	$$R \longrightarrow \prod_i M_{d_i}(k)$$
is a ring homomorphism\footnote{It is the same to ask that the Cayley-Hamilton algebra $(R,T)$ is residually multiplicity free, as explained in 
Example \ref{importantexample}.}.  Let 
$n \in R$ be a nilpotent element that we write $n=(n_s) \in (M_d(K_s))$, and let $\nb \in \prod_i M_{d_i}(k) \subset M_d(k)$ be its projection under the map above. Then $\nb \prec n_s$ for each $s$.
\end{prop}
\begin{pf} Note that $\nb \in M_d(k)$ is nilpotent as the map of the statement is a ring homomorphism, hence the statement makes sense. \ps
	By Prop. \ref{gerstthm}, we have to show that for each $i\geq 1$ and for each $s$,  
\begin{equation}\label{inrkimpliesprop}	\dim_{K_s}(n^i(K_s^d)) \geq \dim_k({\nb}^i(k^d)).\end{equation}
By replacing $n$ by $n^i$ we may assume that $i=1$. Let us write $A^d=\oplus_{i=1}^r V_i$ according to the standard basis, 
$V_i=0 \times A^{d_i} \times 0$. For each $i=1,\dots,r$, let 
$$\nb(v_{i,1}),\cdots,\nb(v_{i,t_i})$$ be a $k$-basis of $\nb(k^{d_i})$, $t_i < d_i$ 
(choose no $v_i$ if $\ut(\nb)_i=0$). Let $w_{i,j} \in V_i$ be any lifting of $v_{i,j}$. To prove (\ref{inrkimpliesprop}), it suffices to show that 
the elements $$n(w_{i,j}) \in K^d, \, \, i=1\dots r, \, \,j=1\cdots t_i,$$
are $K$-independents. It suffices to check that for each $i$, if $p_i: K^d \rightarrow K^{d_i}$ denotes the canonical 
$K$-linear projection on $KV_i$, the elements $$p_i(n(w_{i,j})) \in A^{d_i}, \, \, j=1 \dots t_i,$$
are $A$-independent. By construction these elements reduce mod $m$ to the elements $\nb(v_{i,1}),\cdots,\nb(v_{i,t_i}) \in k^{d_i}$ which are 
$k$-independent, hence we conclude by Lemma \ref{famlibremodlibre}. 
\end{pf}

\begin{prop}\label{nilpoandmatrix} Let $A$ be as in the statement of Prop. \ref{nilpoandgma}, and let $n \in M_d(A)$ be a nilpotent element. 
\begin{itemize}
\item[(i)] for each $s$, $n_s \succ \nb$.
\item[(ii)] if (i) is an equality for each $s$, then $n$ admits 
a Jordan normal form over $A$.
\end{itemize}
\end{prop}

\begin{pf} Assertion (i) follows from Prop. \ref{nilpoandgma} in 
the special case when $R=M_d(A)$. Let us check (ii). 
We have to show that $N_i:=n^i(A^d)$ is free and direct summand as $A$-module. 
By assumption, 
\begin{equation} \label{eqsns} \forall s,\, \, \, \dim_{K_s}(N_i\otimes_A K_s) = \dim_k
(\nb^i(k^d))=:d_i,\end{equation}
and we also have a natural surjection 
\begin{equation}\label{ressns} N_i \otimes_A k \longrightarrow
\nb^i(k^d).\end{equation} 
By Nakayama's lemma and (\ref{ressns}), $N_i$ is generated over $A$ by $d_i$
elements, and those elements are necessarily $K$-independent by
(\ref{eqsns}), thus $N_i$ is free of rank $d_i$ over $A$ and the map in
(\ref{ressns}) is an isomorphism. We conclude by Lemma \ref{famlibremodlibre}.
\end{pf}
%
%\begin{lemma}\label{checkblocks}Let $A$ be a commutative ring and
%$$n=\left(\begin{array}{cc} J_p & \ast \\ 0 & J_q \end{array}\right) \, \, \, \in M_{p+q}(A)$$
%a matrix such that $n^p=0$ (hence $p\geq q$). Then there exists 
%$$P=\left(\begin{array}{cc} 1 & \ast \\ 0 & 1 \end{array}\right)$$
%in $\GL_{p+q}(A)$ such that $PMP^{-1}$ is the block diagonal matrix 
%$J_p \oplus J_q$.
%\end{lemma}
%
%\begin{pf} Let $e_i$ denote the canonical basis of $A^d$.
%We have $n^q(e_d) \in A[J_p](e_p)$. By assumption $n^p=0$, hence
%$$n^q(e_d) \in A[J_p](e_p)\cap \Ker n^{p-q}=n^qA[J_p](e_p).$$
%Let $f\in A[J_p]$ be such that $n^q(e_d)=n^qfe_p$, and set $e'_d=e_d-fe_p$. As $$n^i(e'_d) \equiv e_{d-i} \bmod A[J_p](e_p)$$ for $i<q$, the family $$e_1,\dots,e_p, n^{q-1}(e'_d),\dots,n(e'_d),(e'_d)$$ is an $A$-basis of $A^q$. But $n^q(e'_d)=0$ by construction, and we are done.
%\end{pf}

\begin{lemma}\label{nilpogen} Let $A$ be a noetherian commutative domain, $K$ its fraction field. Let $n \in M_d(K)$ be a nilpotent matrix. There 
exists a nonzero $f \in A$ such that $n \in M_d(A_f)$ and such that for each $x \in D(f) \subset {\rm Spec(A)}$, if 
$n_x$ denotes the image of $n \in M_d(A_x/xA_x)$, then we have $n \sim n_x$.
\end{lemma}

\begin{pf} We may assume that $n \in M_d(A)$. For $i=1\dots d$, let $M_i=n^i(A^d)$. As $M_i \otimes_A K$ is free and direct summand in $K^d$, we may assume, by replacing $A$ by some $A_f$ for a nonzero $f\in A$ if necessary, that all the $M_i$ are free and direct summand in $A^d$. But then for each 
$x \in {\rm Spec}(A)$, $\rk \,n_x^i = \dim_K\, M_i \otimes_A K$, and we are done by Lemma \ref{gerstthm}.
\end{pf}

\subsubsection{Preliminaries on general families of pseudocharacters} 
Even in the "specific" context of the pseudocharacter $T$ on the eigenvarieties $X$ introduced in \S\ref{existfamrepgal}, there is no reason to expect 
that $T$ should be the trace of a representation of $\G_{E,S}$ on a locally free $\OO_X$-module of rank $m$, or on a 
torsion free $\OO_X$-module of generic rank $m$, and this even locally. However, 
this holds for general reasons on an \'etale covering of the Zariski-open subspace of $X$ consisting of 
the $x$ such that $\rhob_x$ is absolutely irreducible by \cite[Cor. 7.2.6]{Ch}. 
Moreover, recall that in the first part of this book, we studied that question in details locally around any point $x$ such that $\rhob_x$ is multiplicity free. \ps

We collect here some general facts that might be used to circumvent this problem. 

\begin{lemma}\label{existrepweak} Let $T: \Gamma \longrightarrow \OO(X)$ be any continuous $m$-dimensional pseudocharacter of a topological group $\Gamma$ on a reduced rigid space $X$ over $\Q_p$. Let $\Omega \subset X$ be an open affinoid.
\begin{itemize}
\item[(i)] There 
is a normal affinoid $Y$, a finite dominant\footnote{Precisely, it is surjective, and the image of any irreducible component of $Y$ is an irreducible component of $\Omega$.} map $g: Y \longrightarrow \Omega$,
 and a finite type torsion free $\OO(Y)$-module $M(Y)$ of generic ranks $m$ equipped with a continuous 
representation $$\rho_Y: \Gamma \longrightarrow \GL_{\OO(Y)}(M(Y))$$ whose generic trace is $T$. 
\par Moreover, $\rho_Y$ is generically semisimple and the sum of absolutely irreducible representations. 
For $y$ in a dense Zariski-open subset $Y' \subset Y$, $M(Y)_y$ is free of rank $m$ over $\OO_y$, and $M(Y)_y \otimes \overline{k(y)}$ is semisimple, and isomorphic to $\rhob_{g(y)}$. 
\item[(ii)] There is a blow-up $g': \CY \rightarrow Y$ of a closed subset of $Y \backslash Y'$ such that the strict transform $\MM_{\CY}$ of the coherent sheaf on $Y$ associated to $M(Y)$ is a locally free $\OO_{\CY}$-module of rank $m$. That sheaf $\MM_{\CY}$ is equipped with a continuous $\OO_{\CY}$-representation of $\Gamma$ with trace $(g'g)^{\sharp}\circ T$, and for all $y \in \CY$, $(\MM_{\CY,y} \otimes \overline{k(y)})^{\rm ss}$ is isomorphic to $\rhob_{g'g(y)}$.
\end{itemize}
\end{lemma}
\begin{pf} Let us prove (i). By normalizing $\Omega$ if necessary, we may assume that $\Omega=X$ is irreducible. 
By Taylor's theorem \cite[thm. 1.2]{Tay}, there exists a finite extension $K'$ of $K:={\rm Frac}(\OO(X))$ 
such that $T: \Gamma \longrightarrow K'$ is the trace of a direct sum of absolutely simple representations 
of $\Gamma \rightarrow \GL_m(K')$. If we define $\OO(Y)$ as the normalization of $X$ in $K'$, the existence 
of a finite type, continuous, $\Gamma$-stable $\OO(Y)$-submodule $M(Y) \in {K'}^m$ is 
\cite[Lemme 7.1 (i), (v)]{BC}. \par It satisfies the "Moreover, ..." assertion by definition. 
By a classical result of Burnside and the generic flatness theorem, this latter fact implies that 
for $y$ in a Zariski-open subset of $Y$, $M(Y)_y=M(Y) \otimes_{\OO(Y)} \OO_y$ is free of rank 
$m$ over $\OO_y$ and that $M(Y)_y \otimes k(y)$ is a semisimple $k(y)[\Gamma]$-module. 
In particular, for those $y$ we have $M(Y)_y \otimes_{\OO_y} \overline{k(y)} \simeq \rhob_{g(y)}$ as they both 
have the same trace, which concludes the proof of (i).

Part (ii) follows then from (i) and Lemma \ref{grra}
(either reduce $Y'$ in (i) or note that the explicit blow-up of $Y$ used in that lemma is the blow-up of the $m$-th Fitting ideal of $M(Y)$, whose associated closed subset does not meet $Y'$.) 
\end{pf}

\subsubsection{Grothendieck's $l$-adic monodromy theorem in families.}\label{grothlmon}

From now and to the end of this section, $F$ denotes a finite extension of $\Q_l$ with $l\neq p$, $\W_F$ its Weil group, $I_F \subset W_F$ its inertia group and 
$\varphi \in W$ a geometric Frobenius. We fix also a nonzero continuous group 
homomorphism $t_p: I_F \rightarrow \Q_p$. The following lemma is well known.

\begin{lemma}\label{defWDingen} Let $B$ be any $\Q$-algebra and $\rho: \W_F \longrightarrow B^*$ any group homomorphism. Assume that
there exists some nilpotent element $N \in B$ such that $\rho$ coincides with 
$g \mapsto {\rm exp}(t_p(g)N)$ on some open subgroup of $I_F$. 
Then $N$ is the unique element with this property. Moreover, the map 
$$r: \W_F \longrightarrow B^*, \, \, \, 
\varphi^ng \mapsto \rho(\varphi^ng){\rm exp}(-t_p(g)N), \, \, \forall \, n\in \Z,\, g\in I_F,$$ 
is a group homomorphism, trivial on some open subgroup of $\W_F$.
\end{lemma}

\begin{definition} Let $\rho$ be as above. If $N$ exists, we say that $(r,N)$ is the Weil-Deligne representation associated to $\rho$. By Lemma \ref{defWDingen}, it is unique and determines $\rho$ entirely.
\end{definition}

We are interested in the study of $p$-adic analytic families of representations of $\Gal(\overline{F}/F)$, and actually a little more generally of $\W_F$. 
Let us give a version of Grothendieck's $l$-adic monodromy theorem adapted to this setting. We 
let $A$ be an affinoid algebra over $\Q_p$ and $B$ a finite type $A$-algebra equipped with 
its canonical $A$-Banach algebra topology. 

\begin{lemma}\label{grothlmonofamilies} Let $\rho: \W_F \longrightarrow B^*$ be a continuous morphism, then $\rho$ 
admits a Weil-Deligne representation.
\end{lemma}

\begin{pf} We fix a submultiplicative norm on $B$ and let $B^0 \subset B$ be its open unit ball. Then $\{1+p^nB^0, n\geq 1\}$ is a basis on open
neighborhoods of $1 \in B^*$ whose successive quotients are discrete and killed by $p$. As a consequence, 
the restriction $\rho'$ of $\rho$ to the 
the wild inertia subgroup of $I_F$ (which is pro-$l$) has a finite image, as its Kernel contains the open subgroup 
${\rho'}^{-1}(1+pB^0)$. Let $F'/F$ be a finite extension such that 
$\rho_{|I_{F'}}$ is tame and pro-$p$, so that it factors through a continuous 
morphism $$t_p(I_{F'}) \longrightarrow B^*.$$ The derivative at $0 \in t_p(I_{F'}) \subset \Q_p$ of the map above 
gives an canonical element $N \in B$. As $\varphi N \varphi^{-1} = \lambda N$ for $\lambda$ a
nonzero power of $l$, $N$ is nilpotent by Lemma \ref{nilplemmegen}, and we are done.
\end{pf}

\begin{lemma}\label{nilplemmegen} Let $A$ be a noetherian $\Q$-algebra and $B$ a finite type (non necessarily commutative) $A$-algebra. If $x \in B$ is $B^*$-conjugate to $\lambda x$ for some $\lambda \in A$ such that $\lambda-1 \in A^*$, then $x$ is nilpotent.
\end{lemma}

\begin{pf} By replacing $B$ by $\End_A(B)$ we may assume that $B=\End_A(M)$ for some finite type $A$-module $M$. When ${\rm Supp}(M)=\{P\}$ is a closed point, then $M$ has a finite length and using the $B$-stable filtration $P^nM$ we may assume that $M$ is a vector space over $A/P$, in which case the result is easy linear algebra. \par
	In the general case, we argue by noetherian induction on the closed subset ${\rm Supp}(M) \in X:={\rm Spec}(A)$. Let $P$ be the generic point of an irreducible component of ${\rm Supp}(M)$, and let $K$ be the kernel of the natural map $M \rightarrow M_P$, it is a $B$-stable submodule. By the previous case $x$ acts nilpotently on $M_P$, and by notherian induction $x_{|K}$ is nilpotent since $P \notin {\rm Supp}(K)$, and we are done.
\end{pf}

\subsubsection{$p$-adic families of $\W_F$-representations}
Let us now fix a topological group $G$, a continuous homomorphism $\W_F \longrightarrow G$, a rigid analytic space $X$ over $\Q_p$ and a continuous $d$-dimensional 
pseudocharacter
	$$T: G \longrightarrow \OO(X).$$
The reason for the appearance of the group $G$ here is exactly the same as in our study of the $p$-adic case (see section \ref{families}). \ps

Let us fix $x \in X$. As already explained in \S\ref{existfamrepgal}, we have two canonical semisimple $G$-representations $\rhob_x$ and $\rhog_x$ with 
respective traces $T\otimes k(x)$ and $T\otimes \Kb_x$. 
As $\rhob_x$ is continuous and defined over a finite extension of $k(x)$, its restriction to $\W_F$ has an associated Weil-Deligne representation. This holds also for
$\rhog_x$. Indeed, let us choose $\Omega$ an open affinoid neighborhood of $x$ and apply Lemma \ref{existrepweak} (i) to this $\Omega$. 
It gives us a continuous representation $\rho_Y: G \longrightarrow \GL_{\OO(Y)}(M(Y))$, which admits a Weil-Deligne representation 
$(r_Y,N_Y)$ by lemma \ref{grothlmonofamilies}. 
If we choose an $\OO(\Omega)$-morphism $\OO(Y) \rightarrow \Kb_x$, we can compose $(r_Y,N_Y)$ with the
ring homomorphism $\End_{\OO(Y)}(M(Y)) \rightarrow M_d(\Kb_x)$ to get a Weil-Deligne representation for ${\rhog_x}_{|\W_F}$, and we are done.

\begin{definition}\label{resgenWD} We call 
$(\rb_x,\Nb_x)$ (resp. $(\rg_x,\Ng_x)$) the residual (resp. generic) Weil-Deligne 
representation of $F$ attached to $T$ at $x$. 
\end{definition}

\begin{remark}\label{resfield}Let $x \in X$. As the $\Q_p$-algebra $\OO_x$ is local henselian and $k(x)$ is finite over $\Q_p$, there is a canonical embedding 
$k(x) \rightarrow \OO_x$ inducing the identity after composition with $\OO_x \rightarrow k(x)$. In particular, we can chose an embedding 
$$\iota_x: \kbx \longrightarrow \Kb_x,$$
and try to compare the two Weil-Deligne representations $(\rb_x \otimes_{\iota_x} \Kb_x,\Nb_x \otimes_{\iota_x} \Kb_x)$ and 
$(\rg_x,\Ng_x)$.
\end{remark}

\begin{lemma}\label{rhoWcst} ${\rg_x}_{|I_F}$ is isomorphic to ${\rb_x}_{|I_F} \otimes_{\iota_x} \Kb_x$.
Moreover, $T_{|I_F}$ is constant on the connected component of $x$ in $X$.
\end{lemma}

\begin{pf} The representation ${\rg_x}_{|I_F}$ has a finite image by construction, hence 
${\rg_x}_{|I_F}$ is actually a semisimple $I_F$-representation and its trace is $\kbx$-valued. 
But this trace coincides by definition with $\iota_x(T_x)$, which proves the first part of the lemma.\ps
Let us show the second assertion. Let $H \subset I_F$ be an open subgroup such that $$T(gh)=T(g) \in k(x), \, \, \forall g \in I_F.$$ Then we just showed that the equality $T(gh)=T(g)$ holds in $\OO_x$, which implies that 
it holds on each irreducible component of $X$ containing $x$, and actually on the whole 
connected component $X(x)$ of $x$ in $X$ by applying the same reasoning to all the points of $X(x)$. In particular, 
$$k_0:=\Q_p(T(I_K)) \subset \OO(X(x))$$ is a finite dimensional $\Q_p$-algebra. But ${\rm Spec}(\OO(X(x)))$ is
connected and reduced, as so is $X(x)$, hence $k_0$ is a field, which concludes the proof.
\end{pf}

Let $(r,N)$ be a $M_d(k)$-valued Weil-Deligne representation, with $k$ an algebraically closed field of characteristic $0$. 
Then the representation $r_{|I_K}$ is semisimple and commutes with $N$, so each of its isotypic component is preserved by $N$. If $\tau$ is any ($k$-valued) finite dimensional irreducible representation of $I_F$, let us denote by $N_{\tau}$ the induced nilpotent element acting on $\Hom_{I_K}(\tau,k^d)$. The following definition is a 
mild extension of Definition \ref{defordonnilp}, and was already studied in \S\ref{types} when $k=\C$.

\begin{definition}\label{deforderingWD} Let $(\rho_1,N_1)$ and $(\rho_2,N_2)$ be two Weil-Deligne representations as above. 
We will write $N_1 \prec_{I_F} N_2$ (resp. $N_1 \sim_{I_F} N_2$) if for each $\tau$, $N_{1,\tau} \prec N_{2,\tau}$ (resp. $N_{1,\tau} \sim N_{2,\tau}$).\par
If both $(\rho_1,N_1)$ and $(\rho_2,N_2)$ are $M_d(k)$-valued, $N_{1,\tau} \prec N_{2,\tau}$ if, and only if, $(\rho_1,N_1)$ is in the Zariski-closure
of the conjugacy class of $(\rho_2,N_2)$. \par Of course, $N_1 \prec_{I_F} N_2$ implies that $N_1 \prec N_2$.
\end{definition}

Let $x \in X$ and write $\Kb_x=\prod_{s(x)} \Kb_{s(x)}$ where $s(x)$ runs the finite set of irreducible components of ${\rm Spec}(\OO_x)$, {\it i.e.} the germs of irreducible components of $X$ at $x$. We can write in the same way $\rhog_x$ and $(\rg_x,\Ng_x)$ as a set of $\Kb_{s(x)}$-representations $\rhog_{s(x)}$ and $(\rg_{s(x)},\Ng_{s(x)})$.

\begin{prop}\label{semicontinuiteWD}Let $x \in X$, $s(x)$ a germ of irreducible component at $x$, and $W$ the\footnote{This component is defined as follows. Let $\Omega$ be any open affinoid of $X$ containing $x$. 
We have a natural map $\OO(\Omega)_x \rightarrow \OO_x$ from the Zariski-local ring at $x$ to the analytic one, which is known to be injective, and both are reduced if $\OO(\Omega)$ is, so we get an injective morphism 
$${\rm Frac}(\OO(\Omega)_x) \hookrightarrow {\rm Frac}(\OO_x).$$ 
The image of ${\rm Spec}(\Kb_{s(x)})$ in ${\rm Spec}(\OO(\Omega)_x)$ is the generic point of a (unique) irreducible component $W_{\Omega}$ of $\Omega$ containing $x$. The component $W$ alluded above is then the unique irreducible component of $X$ containing $\Omega$. It does not depend on the choice of $\Omega$. Indeed, if $\Omega' \subset \Omega$ is another open affinoid containing $x$, then $W_{\Omega'}$ is an irreducible component of $W_{\Omega} \cap \Omega'$, hence is Zariski-dense in $W$.} irreducible component of $x$ in $X$ containing $s(x)$.
\begin{itemize} 
\item[(i)] Let $y \in W$ and $s(y)$ a germ of irreducible component of $X$ at $y$ belonging to $W$. Then $\Ng_{s(x)} \sim_{I_F} \Ng_{s(y)}$.
\item[(ii)] For each open affinoid $\Omega \subset W$, there is a Zariski-dense and Zariski-open subset $\Omega' \subset \Omega$ such that 
$\Nb_y \sim_{I_F} \Ng_{s(x)}$ for all $y \in \Omega'$.
\item[(iii)] $\Nb_x \prec_{I_F} \Ng_{s(x)}$.
\end{itemize}
\end{prop}

\begin{pf} By normalizing $X$ if necessary, we may assume that $X=W$ is normal and irreducible. In particular, $\OO_x$ is a domain for each $x$ hence we will not have to 
specify the $s$ any more: $\rg_x=\rg_{s(x)}$. We may also assume that $X$ is affinoid, and it is enough to show (ii) when $\Omega=X$. \par
	Let $Y$ be a normal affinoid, as well as $g: Y \longrightarrow X$, $M(Y)$, $\rho_Y$ and $Y'$, be 
given by Lemma \ref{existrepweak} (i). By replacing $Y$ by a connected component, we may assume that $Y$ is irreducible.
Note that for each $x \in X$ and $y \in Y$ with $g(y)=x$, we have $$\rhog_y \simeq \rhog_x, \, \, \, \rhob_x \simeq \rhob_y,$$
so we may assume that $Y=X$, and that $T$ is the trace of a continuous representation 
$$\rho: G \longrightarrow \GL_{\OO(X)}(M)$$
on a finite type, torsion free, and generic rank $d$ $\OO(X)$-module $M$. \par
Let $K={\rm Frac}(\OO(X))$. For each $y \in X$, the natural map $\OO(X) \rightarrow \Kb_y$ 
extends to an embedding $\overline{K} \rightarrow \Kb_y$. As the $\overline{K}[G]$-module $M \otimes_{\OO(X)} \overline{K}$ is semisimple by Lemma \ref{existrepweak} (i), we have $\rho \otimes_K \Kb_y \simeq \rhog_y$. By Lemma \ref{grothlmonofamilies}, $\rho$ admits a Weil-Deligne representation $(r,N)$, so by the uniqueness of the Weil-Deligne representation we have 
\begin{equation}\label{equu1}(r,N)\otimes_K \Kb_y \isomo (\rg_y,\Ng_y),\, \, \forall y \in X,\end{equation}
which proves (i). \par
	Let us show (ii). By replacing $X$ by a finite etale covering coming from the base, we may assume that the irreducible representations of the finite group $r(I_F)$ are all defined over some local field $k_0 \subset \OO(X)$, hence we can write the following finite decomposition of $M_{|W_F}$ $$M=\bigoplus_{\tau}\tau \otimes_{k_0} M_{\tau}, \, \, \, M_{\tau}:=\Hom_{k_0[r(I_F)]}(\tau,M_{\tau}).$$  
Let us choose a nonzero $f \in \OO(X)$ such that each $(M_{\tau})_f$ is a free $\OO(X)_f$-module, and that 
$M \otimes k(y) \simeq k(y)^d$ is a semisimple $k(y)[G]$-representation for each $y$ in $X_f$ (use Burnside's theorem). In particular, 
\begin{equation}\label{equu2} M \otimes_{\OO(X)} k(y) \simeq \rhob_y, \ \ \forall y \in X_f.
\end{equation} 
Applying Lemma \ref{nilpogen} to $N_{\tau} \in \End_{\OO(X)_f}((M_{\tau})_f)$, we may assume by changing $f$ if necessary that 
\begin{equation}\label{equu3} N_{\tau,y} \sim N_{\tau} \otimes K,\, \, \forall y \in X_f,\end{equation}
hence (ii) holds by (\ref{equu1}), (\ref{equu2}) and (\ref{equu3}) if we take $\Omega':=X_f$. \ps
	It only remains to prove assertion (iii). We claim first that me may assume that the module $M$ defined in the second paragraph above of the proof is free. Indeed, take $g': \CY \rightarrow Y=X$ and $\MM_\CY$ as in Lemma \ref{existrepweak} (ii). If $V$ an affinoid subdomain of $\CY$, then $V \cap Y'$ is Zariski-dense in $V$, so $\MM_{\CY}(V)_y$ is a direct sum of absolutely irreducible $G$-representations for each $y$ in a Zariski-dense subset of $V$. Thus $\MM_{\CY}(V)\otimes_{\OO(V)}{\rm Frac}(\OO(V))$ has the same property, which proves the claim by replacing $X$ by $\CY$, and then by an affinoid subdomain as $\CY \rightarrow X$ is surjective. In particular, the Weil-Deligne representation $(r,N)$ is now $M_d(\OO(X))$-valued.
	
	As $M$ is free over $\OO(X)$, we have
$$(M \otimes_{\OO(X)} \overline{k(x)})^{\rm G-ss} \simeq \rhob_{x}, \, \, \forall x \in X,$$
so Lemma \ref{ineqnilp} (i) implies that it is enough to get (iii) to check that the image $N_x$ of $N$ in $\End_{k(x)}(M\otimes k(x))$ satisfies
$$N_x \prec_{I_F} N\otimes K.$$
As this is an assertion on the action of $\W_F$, we may decompose $M$ (again, up to enlarging the base field if necessary) as a sum of its isotypic components $M_{\tau}$ as above, and the result follows then from Prop. \ref{nilpoandmatrix} (i) applied to the Zariski-local ring $A$ of $X$ at $x$ and to $N$ acting on the free $A$-module $M_{\tau} \otimes_{\OO(X)} A$.
\end{pf}

\begin{remark}\label{remsemicontinuiteWD} Let $x \in X$. For each $s \in {\rm Spec}(\OO_x)$, say with residue field $k(s)$, there exists a unique (isomorphism class of) semisimple representation 
	$$\rho_s: G \longrightarrow \GL_d(\overline{k(s)})$$
whose trace is the composite $T(s): G \rightarrow \OO_x \rightarrow k(s)$. The same argument as we gave for $\rg_x$ shows that $\rho_s$ admits a Weil-Deligne representation $(r_s,N_s)$. As an exercise, the reader can check using a slight variant of the proof of Prop. \ref{semicontinuiteWD} that 
$$s \in \overline{\{s'\}}\, \, \,  \Rightarrow\, \, \,  N_s \prec_{I_F} N_{s'}.$$ 
\end{remark}

\newpage
\renewcommand{\m}{\got m}

\section{The sign conjecture}
%\secttoc

\subsection{Statement of the theorem}\label{statemtthmsign}

We use the notations of section~\ref{Selmer}, especially of~\S\ref{quad}:
$E$ is a quadratic imaginary field, $p$ a prime that is split in $E$,
$$\rho: \G_E \longrightarrow \GL_n(L)$$ an $n$-dimensional, geometric, semisimple,
representation of $G_E$ with coefficients in a finite extension $L/\Q_p$, satisfying $$\rho^{\bot}\simeq \rho(-1).$$
We fix also embeddings $\iota_{\infty}$ and $\iota_p$ as in \S \ref{setnotrepm}. We denote by $v$ and $\bar v$ the
two places of $E$ above $p$ as in {\it loc. cit}. We make the following assumptions on $\rho$ :

\begin{itemize}
\item[(1)] The dimension $n$ is not divisible by $4$.\ps
\item[(2)] There is a cuspidal tempered automorphic
representation $\pi$ of $\Gl_n(\A_E)$, satisfying properties (i),(ii)
and (iii) of \S\ref{start} and such that for every split place $x$ of $E$
the Weil-Deligne representation of $\iota_p\iota^{-1}_{\infty}\pi_x|\det|_x^{1/2}$ and the one
attached to ${\rho}_{|E_x}$ are isomorphic up to Frobenius semi-simplification.\ps
\item[(3)] The representation $\rho_{|E_v}$ is crystalline and the characteristic polynomial of its
crystalline Frobenius is the same as the one of $\iota_p\iota_\infty^{-1}(L(\pi_v|\det|^{1/2}))$.
\end{itemize}
\ps
Note that for the sake of generality, and because irreducibility
may be hard to check in applications, we do not assume that
$\rho$ is irreducible.

\begin{example}\label{exampleHT81} There are many known examples of such
$\rho$, in any dimension $n$ (say not divisible by $4$). Start with a cuspidal
representation $\pi$ of $\Gl_n(\AAA_E)$
satisfying the hypothesis (i), (ii), (iii) of \S\ref{start}. Assume moreover
that $\pi$ is square-integrable at some finite
place.\footnote{When $n=2$, it is well known that this assumption is not necessary.
When $n=3$, it is actually possible to remove it, but we have
to assume that $p$ is outside a density zero set of primes
depending on $\pi$. Indeed, if $\pi$ satisfies (i) and (ii) of \ref{start},
$\pi \mu$ descends by Rogawski's base change to the quasi-split unitary group
$\U(2,1)(\Q)$ attached to $E$ (or even to the form that is compact at
infinity), there is a Galois representation that we may write $\rho \mu$
attached to $\pi \mu$ satisfying (2) by~\cite{br1}. The temperedness of $
\pi$ together with the compatibility condition in (3)
are then the main result of \cite{belcomp}.} Then by the
main result of \cite{HT}, $\pi$ is tempered and there is a Galois representation $\rho$ satisfying (1), (2) and (3).
In this case, by \cite{TY} we also know that $\rho$ is irreducible.
It should be possible, in a near future,
to remove the square-integrability hypothesis using  results of \cite{GRFAbook}, but then the
irreducibility of $\rho$ might not be known.
\end{example}

Recall that we introduced previously assumptions Rep and AC (conjectures \ref{repm} and
\ref{AC}).
\begin{theorem}\label{signthm} Assume AC($\pi$) and Rep($n+2$).
Then the sign conjecture holds for $\rho$:
namely, if $\varepsilon(\rho,0)=-1$, then
$\dim_L H^1_f(E,\rho) \geq 1$
\end{theorem}

Since hypotheses AC($\rho$) for a character $\rho$ and Rep$(3)$ are known (see Remarks \ref{remsectrepm} (vi) and \ref{remsectAC} (ii)),
we deduce:
\begin{cor} \label{corsign1}
If $n=1$, for a $\rho$ as above, the sign conjecture holds.
\end{cor}
This result was the main result of
\cite{BC}, where it was proved by similar methods, and can also be deduced
of earlier results of Rubin (see the introduction of {\it loc. cit.}).

For $n=2$ we can prove a result avoiding the hypotheses at nonsplit primes.

\begin{cor} Let $f$ be a modular form of even weight
$k \geq 4$ and level $N$ prime to $p$, and $\rho=\rho_f$ as in Example~\ref{examplerhof}.
Assume AC$(\pi_{f,E})$\footnote{The needed representation $\pi_{f,E}$ is actually defined in the proof. Precisely, 
$\pi_{f,E}$ is the base change to a quadratic imaginary field $E$ as in Prop.~\ref{propQE} 
of the automorphic representation $\pi_f$ of $\GL_2(\AAA_\Q)$ attached to $f$, where $\pi_f$ is normalized so that $\pi_f^* \simeq \pi_f$.} and Rep$(4)$.
Then the sign conjecture holds for $\rho_f$, namely if
$\varepsilon(\rho_f,0)=-1$, we have $\dim_L H^1_f(\Q,\rho_f) \geq 1$.
\end{cor}
\begin{pf}
By Proposition~\ref{propQE}, there is a quadratic imaginary  field $E$
where $p$ and
all the primes dividing $N$ are split, and such that
$\varepsilon(\rho_{f,E},0)=\varepsilon(\rho_f,0)$
and $H^1_f(\Q,\rho_f)=H^1_f(E,\rho_{f,E})$.
So the corollary follows from Theorem~\ref{signthm} if we verify
 hypotheses (1) to (3) for $\rho=\rho_{f,E}$. Assumption (1) is clear.
For the automorphic representation $\pi$ needed in (2) we simply take $\pi:=\pi_{f,E}$ the Langland's
base change to $E$ of the automorphic representation $\pi_f$ of
$\Gl_2(\A_\Q)$ attached to $f$ which is normalized so as to be autodual: $\pi_f \simeq \pi_f^*$ and it has a trivial central character. 
It is clear that $\pi$ satisfies (i) and (ii)
of \ref{start} since the $L$-parameter of $\pi_{\infty}$ coincides with
$$z \mapsto \diag((z/\bar
z)^{\frac{1-k}{2}},(z/\bar z)^{\frac{k-1}{2}})$$ on $\C^*$ and $k>2$ is even, and also (iii)
since $\pi$ is unramified at non-split places of $E$ by construction of $E$.
It is well known that assumption (3) holds since $p$ does not divide $N$.
\end{pf}

\begin{remark}\label{remaprescasgl2} \begin{itemize}
\item[(i)] When $f$ is ordinary at $p$, this result, without its automorphic assumptions, was proved by Nekovar as a consequence of his 
{\it parity theorem}. A similar result was also proved later by Skinner-Urban in \cite{SU}, using
automorphic forms on the symplectic group $\text{GSp}_4$. Since the existence
of Galois representations attached to such forms is known
they do not assume any variant of hypothesis Rep$(4)$ but they
have stronger hypotheses on $f$, namely that $p$ is an ordinary prime
for $f$, and that $N=1$. 
\item[(ii)] It may be possible to remove the restriction $k>2$ from this
result (or, for that matter, the restriction on the weight in
Thm~\ref{signthm}) by actually deducing the $k=2$ case from the result above and a deformation argument. 
We postpone this to a subsequent work.
\end{itemize}
\end{remark}

\begin{remark} \label{suppressassumpt} As we explained in section 5, it was not our policy in this book to assume
the most general versions of Langlands and Arthur's conjectures on the
discrete spectrum of unitary groups, but rather to formulate a minimal set
of expected assumptions which we prove to be enough to imply the sign
conjecture in a large number of cases. Indeed, the version of Theorem \ref{signthm}
that we state is actually the stronger that we can prove under our
assumptions (Rep(n+2)) and (AC($\pi$)). A reason for that restriction is that is not clear to us
which part of those general conjectures will be proved first and in which
form (this might clarify itself during the completion of the book project \cite{GRFAbook}). 
This especially applies to the part concerning the Langland's 
parameterization for the local unitary group 
and the local-global compatibility of the base change from $\U(m)$ to $\GL_m$ at those primes (see Appendix A), of which our proof would 
need some properties ({\it e.g.} if we do not want to make assumption (iii) of \S\ref{start}). 
Let us simply say that we believe that {\it at the end}, those general conjectures should imply the full case of the sign conjecture for the $\rho$ attached to a regular algebraic 
cuspidal automorphic representation of $\GL_n(\AAA_E)$. We hope 
to go back to this extension as a part, or a sequel, 
of the aforementioned book project.
\end{remark}

The following two subsections are devoted to the proof of~\ref{signthm}.

\subsection{The minimal eigenvariety $X$ containing $\pi^n$}

\subsubsection{Definition of $\pi^n$ and $X$}\label{mineigenprepa} From now till the end of this book, we set 
$$m:=n+2.$$ Assume that $\varepsilon(\pi,0)=-1$ and let $\pi^n$ the automorphic representation of $\U(m)$ given by $AC(\pi)$. Recall that the representation $\pi^n$ depends on the choice 
of a Hecke character $\mu: \AAA_E^* \longrightarrow \C^*$ as in Definition \ref{choixdemu}. Recall that $\mu^{\bot}=\mu$, that $\mu=1$ if $m$ is even, and that $\mu$ does not descend to $U(1)$ when $m$ is odd. \ps

\index{X@$X$, the minimal eigenvariety containing $\pi^n$}

	By \S\ref{start}, for each prime $l$ that does not split in $E$, $\pi_l$ is either NMSRPS or unramified, so that it makes sense to consider the {\it minimal eigenvariety} $X$ containing $\pi^n$ as in Example \ref{minimaleigen}.\footnote{The results 
here will apply whether we choose or not the variant with a fixed weight.} We will use in the sequel the same notations as in \S\ref{famillegaloisiennesurX}. In particular, 
recall that $S$ is the finite set of primes consisting of $p$ and of the primes $l$ such that $\pi_l$ is ramified, and $L/\Q_p$ is a big enough finite extension of $\Q_p$ on which $\pi^n$ and $X$ are defined. Assume also that $L$ is big enough so that 
$\rho$ is a sum of absolutely irreducible representations defined over $L$. \ps

	In order to associate a point of $X$ to $\pi^n$ we have to specify an accessible refinement of $\pi_p^n$. They are given by the following lemma. Recall that the place $v$ fixes an isomorphism $U(m)(\Q_p)\isomo_v \GL_m(\Q_p)$. 

\begin{lemma}\label{refofpin} \begin{itemize}
\item[(i)] The representation $\pi_p^n$ is almost tempered (see Def. \ref{defalmtemp}). 
\item[(ii)] Its accessible refinements are the $n!\frac{(n+1)(n+2)}{2}$ orderings of the form 
$$\mu_w|.|^{-1/2}(p)(\dots,1,\dots,p^{-1},\dots),$$
where $1$ preceds $p^{-1}$.
\item[(iii)] The ordered set of the other eigenvalues in the dots above is any ordering of 
the eigenvalues of the Langlands conjugacy class $L(\pi_p|.|^{1/2})$. Each of those eigenvalues has complex norm $\sqrt{p}^{-1}$, and in particular is different from $1, p^{-1}$.
\end{itemize}
\end{lemma}

\begin{pf} By Remark \ref{remlsplit}, $\pi_p^n$ is the unramified representations such that 
$$L(\pi^n_p)=\mu_v|.|^{-1/2}(L(\pi_p|.|^{1/2})\oplus 1 \oplus |.|),$$
and the parameter $L(\pi_p)$ is bounded as $\pi_p$ is tempered by assumption. The lemma follows now from Prop. \ref{accessalmtemp}.
\end{pf}

\index{R@$\RR$, a refinement of $\pi^n_p$}
\index{z@$z \in Z \subset X$, the point of $X$ corresponding 
to $(\pi_p^n,\RR)$}
Let us choose any such refinement $\RR$ for the moment, which fixes an associated point $z=(\psi_{(\pi^n,\RR)},\uk)$. 

\begin{remark}\label{remcorcomprho} Note that we did not make any assumption on the compatibility between $\rho$ and $\pi|.|^{1/2}$ at the nonsplit primes.
Actually, we can prove a version of it under the running assumptions. Indeed, as 
$X$ is the minimal eigenvariety containing $\pi^n$, Prop. \ref{nonmonatNMSRPS} and Prop. \ref{monauxplacessplit}, $\rho$ is 
unramified outside $S$, and has a trivial monodromy operator at the primes $l \in S_{NM}$. In the same way, Lemma \ref{calculdeskappai} shows that the Hodge-Tate 
weights of $\rho_{|E_v}$ correspond to the highest weight $\uk$ of $\pi_{\infty}^{n}$ as in part (P4) of property (Rep(m)) of \S\ref{sectrepm}. In particular, 
those Hodge-Tate weights are two-by-two distinct.
\end{remark}

\subsubsection{Normalization of the Galois representation on $X$}\label{normalisationrepgal}	As explained in \S\ref{existfamrepgal}, by assumption $Rep(m)$ we have a continuous pseudocharacter $$T: \G_{E,S} \longrightarrow \OO(X),$$	
such that $T(g^{\bot})=\chi(g)^{m-1}T(g)$ for all $g \in\G_{E,S}$. Recall that $\chi$ is the cyclotomic character. It will be convenient to twist it by a constant character as follows. The following 
lemma is immediate (see \S\ref{introcarhecke}). \ps

\begin{lemma} The Hecke character $\mu^{-1}|\cdot|^{\frac{m}{2}}$ has an integral weight $\delta:=\frac{m}{2}$ if $m$ is even,
and $\delta:=\frac{m-1}{2}$ if $m$ is odd. In particular, it is an algebraic Hecke character.
\end{lemma} 

Enlarging a bit our base field $L/\Q_p$ if necessary, 
we may assume that it is defined over $L$ via $\iota_p\iota_{\infty}^{-1}$. 
By class field theory, there is an associated\footnote{This character is the unique 
continuous character $\nu$ such that for each finite prime $w$ of $E$ prime to $p$, 
$$\nu_{|\W_{E_w}}=\iota_p\iota_{\infty}^{-1}(\rec^{-1}(\mu_w^{-1}|\cdot|^{\frac{m}{2}}))$$ 
where the {\it local} $\rec$ map is the one discussed in \S\ref{localLanglands}.} continuous 
character $$\nu=\mu^{-1}|\cdot|^{\frac{m}{2}} \circ \rec^{-1}:\, \, \G_{E,S} \longrightarrow L^*.$$
By Cebotarev theorem and Rem. \ref{remlsplit}, the evaluation of $T$ at the point $z$ is the trace of the representation
$$\rho \nu^{-1} \oplus \nu^{-1} \oplus \chi\nu^{-1}.$$
This leads us to define $T'$, $\kappa_i'$ and $F_i'$ as follows: \begin{itemize}
\item[-] $T':=T\otimes\nu$, {\it i.e.} $T'(g)=T(g)\nu(g) \, \forall g\in G_{E,S}$, 
\item[-] $\kappa'_i:=\kappa_i - \delta$ for $i=1,...,m$, 
\item[-] $F'_i:=F_i\iota_p\iota^{-1}(\nu_v(p))$ for $i=1,...,m$.
\end{itemize}

\begin{definition}\textdbend \label{defnormbydelta}
From now on, we shall use the letters $T$, $\kappa_i$ and $F_i$ to denote respectively $T'$, $\kappa_i'$ and $F_i'$ above. With this choice of normalisation we have
\begin{equation}\label{conseqnorma} \rhob_z= 1 \oplus \chi \oplus \rho,\, \, \, \, \, T^{\bot}=T(-1),\end{equation}
and the new $(X,T,\{\kappa_i\},\{F_i\},Z)$ is obviously still a refined family.
\end{definition}

\subsubsection{The faithfull GMA at the point $z$.}\label{faithgmachoix}

	Let $A:=\OO_z$ be the rigid local ring at the closed point underlying to $z$, and $\m$ its maximal ideal, $k=k(z)=A/\m \simeq L$. 
We will focus on the $A$-valued $m$-dimensional continuous pseudocharacter induced by $T$, 
		$$T: \G_{E,S} \longrightarrow A,$$
that we denote also by $T$ (rather than $T \otimes_{\OO(X)} A$). Let $$R:=A[\G_{E,S}]/\ker T$$ be the faithful Cayley-Hamilton algebra associated to $T$. 
Recall that $$\rhob_z=\rho \oplus 1 \oplus \chi.$$ As the Hodge-Tate weights of $\rhob_z$ are two-by-two distincts by Rem. \ref{remcorcomprho}, $\rhob_z$ is multiplicity free. For later use, let us write $$\rho = \oplus_{j=1}^r \rho_j$$ 
where the $\rho_j$ are pairwise non isomorphic and absolutely irreducible. 
As each $\rho_j$ is defined over $L$ by assumption on $L$, $T$ is actually {\it residually multiplicity free} in the sense of Definition \ref{defmult1}. By Theorem \ref{structure} and Remark \ref{remthmstru} 
we get the following lemma.

\begin{lemma}\label{Ristfgma} $R$ is a GMA over $A$ and is a finite type, torsion free, $A$-module. 
\end{lemma}

We will be interested in the $\Ext$-groups between the irreducible constituents of $\rhob_z$. For this purpose we set $${\cal I}:=\{\chi,1,\rho_1,\cdots,\rho_r\},$$
which is also the set of simple $R$-modules by Lemma  \ref{radical}.

\begin{definition}\label{extT} If $i, j\in {\cal I}$ are two irreducible factors of $\rhob_z$, we set\footnote{Note that the last equality in the definition below 
is Remark \ref{extSmorS}.}
		$$\Ext_T(i,j):=\Ext_{R\otimes_A k}(i,j)=\Ext_R(i,j).$$
It is a finite dimensional $k$-vector space.\footnote{This follows for example from Theorem \ref{extension1} and Lemma \ref{Ristfgma}.}
\end{definition}

The following lemma follows from Prop. \ref{topoprop}.

\begin{lemma}\label{extTinj} The natural $k$-linear injection $$\Ext_T(i,j) \hookrightarrow \Ext_{k[\G_{E,S}]}(i,j)$$ falls inside 
the subspace of continuous extensions of $i$ by $j$ as $k[\G_{E,S}]$-representations. 
\end{lemma}

\begin{remark}\label{remextT} By definition, the image of the inclusion above is exactly the set of extensions of $i$ by $j$ that occur in some subquotients of some $R$-module $M$. 
As $R=A[G]/\ker T$ is the only natural Cayley-Hamilton quotient of $(A[G],T)$ that we can consider a priori here, 
$\Ext_T(i,j)$ should be thought as {\it the space of extensions of $i$ by $j$ that we can construct from the datum of the pseudocharacter
$T$}, hence the notation $\Ext_T$.
\end{remark}

We will study now the local conditions at each primes of the elements in $\Ext_T(i,j)$.

\subsubsection{Properties at $l$ of $\Ext_T(i,j)$}\label{proplTsign} Let us fix $l\neq p$ a prime, $w$ a prime of $E$ above $l$, as well as a decomposition group $\G_{E_w} \longrightarrow \G_{E,S}$. We begin by a general lemma.

\begin{lemma}\label{somlocext} Let $V$ be the semisimplification of the representation $\rho_{|\G_{E_w}}$.
\begin{itemize}
\item[(i)] Assume that $l$ splits in $E$. For $d \in \Z$, $\chi^d$ is 
not a subrepresentation of $V$, and  $$\Ext_{L[\G_{E_w}]}(V,\chi^d)=\Ext_{L[\G_{E_w}]}(\chi^d,V)=0.$$
\item[(ii)] $\Ext_{L[\G_{E_w}]}(\chi,1)=0$.
\end{itemize}
\end{lemma}

\begin{pf} Assume that $l$ splits in $E$. We claim first that for all $d \in \Z$, $l^d$ is not an eigenvalue of a Frobenius at $w$ in 
$\rho_{|E_w}^{I_w}$. Indeed, as $\pi_w$ is tempered, 
the eigenvalues of any geometric Frobenius element $\phi_w$ in the complex Weil-Deligne representation attached to $\pi_w|.|^{1/2}$ 
have norm $\sqrt{l}$. This proves already the first part of (i) by assumption (2) on $\rho$.\par
Let $W$ be either $V(d)$, $V^*(d)$ or $\chi^{-1}$. We need to show that $H^1(E_w,V)=0$. As $l\neq p$, we know from Tate's theorem that 
	$$\dim_L H^1(E_w,W) = \dim_L H^0(E_w,W)+\dim_L H^0(E_w,W^*(1)),$$ so the case $W=\chi^{-1}$ is clear and the 
other ones follow from the claim and assumption (2) on $\rho$. 
\end{pf}

\begin{prop}\label{ext1cyclol} For each $i \neq 1$ in $\cal I$, 
$\Ext_T(1,i)$ consists of extensions which are split when restricted to $I_{E_w}$.
\end{prop}

\begin{pf} Let $R_w \subset R$ be the image of $A[\G_{E_w}]$ in $R$ via the natural map $A[\G_{E,S}] \longrightarrow R$. 
It is of finite type over $A$ as $R$ is and as $A$ is noetherian. Let $K$ be the total fraction field of $A$, and 
set $R_K=R\otimes_A K$. As $R$ is torsion free over $A$, $R \subset R_K$. Let us choose any datum of idempotents $e_{\chi},e_1,\dots$ 
as well as a representation $R_K \rightarrow M_m(K)$ adapted to the chosen $\{e_i\}$ as in Theorem \ref{structure}, 
and consider the induced representation
	$$\rho_K: \G_{E,S} \longrightarrow \GL_m(K).$$
By Lemma \ref{lemmeutile} and Prop. \ref{specialcase}, $\rho_K$ is semisimple and the sum of absolutely irreducible representations, so $\rho_K \otimes \Kb \simeq \rhog_z$. 
In particular ${\rho_K}_{|E_w}$ has an associated Weil-Deligne representation $(r,N)$ with values in $R_w$. \par
	The argument will be different according as $l$ splits or not in $E$. As the proposition is obvious if $l \notin S$, 
we may assume that $l \in S$. Let us assume first that $l$ does not split in $E$, which implies that $l \in S_{NM}$. 
Note that for a continuous extension $U$ of $1$ by $i\neq 1 \in \cal I$ to be trivial when restricted to $I_{E_w}$, 
it is enough to check that $I_{E_w}$ acts through a finite quotient on $U$, as $\Q$-linear 
representations of finte groups are semisimple. By Theorems \ref{extension1} and \ref{extension2} (1), 
it suffices to show that the image of $N$ in $R \subset R_K$ is trivial, because then $I_{E_w} \rightarrow R^*$ 
factors through a finite quotient. But as $X$ is a minimal eigenvariety containing $\pi^n$, Prop. \ref{nonmonatNMSRPS} 
shows that $N=0$ when $w \in S_{NM}$, and we are done in this case. \ps
	Let us assume now that $l=w\bar{w}$ splits in $E$. By Lemma \ref{somlocext} (i), there is nothing to prove when $i \neq \chi$, hence 
we concentrate from now on $\Ext_T(1,\chi)$. We will need to choose a specific GMA structure of $R$. 

\begin{lemma}\label{prepacentidem} Let $A$ be local henselian noetherian commutative ring, $\m$ its maximal ideal, $S$ an $A$-algebra (non necessarily commutative) which is finite type as $A$-module.
Let $\Irr(S)$ be the (finite) set of simple $S$-modules, or which is the same of simple $S/\m S$-modules, and let $\cal P \subset \Irr(S)$ a subset with the following property:
	$$\forall M \in \cal P,\, \, \, N \in \Irr(S)\backslash {\cal P},\, \, \Ext_{S/\m S}(M,N)=\Ext_{S/\m S}(N,M)=0.$$
Then there is a unique central idempotent $e \in S$ such that for each $M \in \Irr(S)$, $e(M)=M$ if $M \in \cal P$, $0$ if else.
\end{lemma}

\begin{pf} Note that if $M$ is a simple $S$ module, it is monogenic over $S$ hence of finite over $A$, thus $\m M=0$ by Nakayama's lemma. It shows 
that $\Irr(S/\m S) \rightarrow \Irr(S)$ is bijective. Moreover, $\m S \subset \rad(S)$. \ps
Assume first that $A=k$ is a field, hence $S$ is any finite dimensional $k$-algebra. Let $M$ be any finite type $S$-module, $M$ has finite lenght. 
Define $M_{\cal P}$ (resp. $M^{\cal P}$) as the largest submodule of $M$ all of whose simple subquotients lie in $\cal P$ (resp. in $\Irr(S)\backslash 
\cal P$). Obviously, $M_{\cal P}\cap M^{\cal P}=0$. We claim that $M=M_{\cal P} \oplus M^{\cal P}$. Ab absurdum, 
as $M$ is of finite lenght, we can find a submodule $M_{\cal P} \oplus M^{\cal P} \subset M' \subset M$ such that 
$$Q:=M'/(M_{\cal P} \oplus M^{\cal P}) \in \Irr(S).$$
 By the $\Ext$-assumption and an immediate induction, note that $\Ext_S(A,B)=\Ext_S(B,A)=0$ whenever $A$ and $B$ are 
finite lenght $S$-modules such that each irreducible subquotient of $A$ (resp. of $B$) lies in $\cal P$ (resp in $\Irr(S)\backslash \cal P$). 
Assume for example that $Q \notin \cal P$. The remark above shows that there is an $S$-submodule 
$M_{\cal P} \subsetneq M_0 \subset M'$ such that $M'=M_{\cal P} \oplus M_0$. 
But $M_0$ has all its subquotients 
in $\Irr(S)\backslash \cal P$, a contradicton. The case $Q \in \cal P$ is similar, which proves the claim. \ps
	We check at once that the decomposition $M=M_{\cal P} \oplus M^{\cal P}$ is stable by $\End_S(M)$. 
In particular, we can write $S=S_{\cal P}\oplus S^{\cal P}$ and we get that both $S_{\cal P}$ and $S^{\cal P}$ are two-sided ideals of $S$. 
We check now at once that the element $e \in S_{\cal P}$ given by the decomposition above of $1=e+(e-1) \in S$ is a central idempotent with all the required properties, thus proving the case where $A$ is a field.\ps
	In general, we choose $e \in S/\m S$ as above. As $A$ is henselian and $S$ finite over $A$, 
there is an idempotent $f \in S$ lifting $e$. By reducing mod $\m$ the direct sum decomposition 
	$$S=fSf \oplus (1-f)Sf \oplus fS(1-f)\oplus (1-f)S(1-f),$$
and as $e$ is central, we get that $fS(1-f)\otimes_A A/\m=(1-f)Sf \otimes_A A/\m=0$, hence $fS(1-f)=(1-f)Sf=0$ in $S$ by Nakayama's lemma. In other words, $f$ is central, and we are done.
\end{pf}

We show now that up to $R^*$-conjugation, $R_w$ is bloc diagonal of type $(2,n)$ in $R$.

\begin{lemma}\label{gma2n} There is a datum of idempotents $\{e_i, i\in {\cal I}\}$ for the generalized matrix algebra $R$ such that $e:=e_{\chi}+e_1$ is in the center of $R_w$. \par
\end{lemma}

\begin{pf} We have $\rad(R)\cap R_w \subset \rad(R_w)$ by Lemma \ref{radical}, so the set ${\cal I}_w$ of simple $R_w/\m R_w$-modules is the set of irreducible subquotients of the $W_{|\G_{E_w}}$ with $W \in {\cal I}$. 
Let us consider the subset $$\cal P=\{\chi,1\} \subset {\cal I}_w.$$
By the second assertion of Lemma \ref{somlocext} (i), we can apply Lemma \ref{prepacentidem} to $S=R_w$ and the set $\cal P$ above, which gives us a central idempotent $e \in R_w$. \par 

By the first assertion of Lemma \ref{somlocext} (i), $T(e)={\bar T}(e)=2$, so
if we consider now the restriction $T_e$ of $T$ to $eRe$, it is a Cayley-Hamilton pseudocharacter of dimension $2$ (see Lemma \ref{dimpseudo}) which is residually multiplicity free with residual representations $1$ and $\chi$. 
By Lemma \ref{idempotents}, we can then write $$e=e_{\chi}+e_{1}$$
where $e_{\chi},\, e_1 \in eRe$ lift the residual idempotents $1=\epsilon_{\chi}+\epsilon_1$ (see the proof of Lemma \ref{idempotents}). We conclude the proof by lifting then successively the remaining residual primitive idempotents in $(1-e)R(1-e)$ and arguing as in the first part of the proof of Lemma \ref{idempotents}, or better by applying that lemma to $(1-e)R(1-e)$.
\end{pf}

We can now conclude the proof of Prop. \ref{ext1cyclol}. Let us choose a datum of idempotents $e_{\chi},e_1,\dots$ as in Lemma \ref{gma2n} as well as 
a representation $R_K \rightarrow M_m(K)$ adapted to those $\{e_i\}$ as above. 
Note that a continuous $\G_{E_w}$-extension of $1$ by $\chi$ is trivial if and only if its monodromy operator is trivial. 
By Theorems \ref{extension1} and \ref{extension2} (1), it suffices to show that the image $eN$ of $N$ in 
$eR_we=eR_w \subset eR_Ke$ is trivial. Write $N=(N_s) \in M_d(K)$, $K=\prod_s K_s$. By the minimality of $X$, assumption (2),
and by Prop. \ref{monauxplacessplit}, we have for each germ $s$ of irreducible component of $X$ at $x$, 
	$$\Nb_z \prec N_s \prec \Nb_z,$$ 
hence $\Nb_z \sim N_s$ (recall that $\Nb_z$ is the monodromy operator of $\rhob_z$). We have to show that $eN=0$. 
Let us write by abuse $(1-e)\Nb_z$ for the image of $\Nb_z$ in $\prod_{j=1,...,r}\End(\rho_j)$.
We have $\Nb_z \sim (1-e)\Nb_z$ as $1$ and $\chi$ are unramified. But $(1-e)\Nb_z \prec (1-e)N_s$ by Prop. \ref{nilpoandgma} applied to $(1-e)R(1-e)$ and $n=(1-e)N$, so we get
\begin{equation}\label{majonil} \forall s, \, \, \, (1-e)N_s \sim N_s \sim (1-e)\Nb_z \end{equation} 
and $eN=0$ by the first of the $\sim$'s above. 
\end{pf}

Let us record now a fact that we will use in section 9. We assume here, and only here, that $\rho$ is irreducible. Let $I_{\rm tot} \subset A$ be the total reducibility locus of $T$ (see \S\ref{assredloc} and Definition \ref{defredloc}). Let $J \supset I_{\rm tot}$ be a proper ideal of cofinite lenght of $A$. Recall that $T \otimes A/J$ writes uniquely as the sum of of three residually irreducible pseudocharacters $$T \otimes A/J = T_{\chi}+T_1+T_{\rho}$$ of respective dimension $1$, $1$ and $n$, 
lifting the decomposition of $T \otimes k$. 
Moreover, let $$R_{\rho}: \G_{E,S} \longrightarrow \GL_m(A/J)$$ be the unique (up to conjugation) continuous representation with trace $T_{\rho}$ (see Def. \S\ref{lesrhoi}, Prop. \ref{topoprop}). 

\begin{lemma}\label{lemmemonopourplustard} Assume that $\rho$ is irreducible. The monodromy operator of ${R_{\rho}}_{|E_w}$ admits a Jordan normal form over $A/J$.
\end{lemma}

\begin{pf} We keep the notations of the proof above. As we showed, this monodromy operator is $0$ if $w$ is not split, hence we may assume that is does. 
As $\rho$ is irreducible, $(1-e)R(1-e) \simeq M_n(A)$, and by Lemma  \ref{conseq} (ii), it suffices to show that $(1-e)N \in (1-e)R(1-e)$ admits a Jordan normal form over $A$. By (\ref{majonil}), 
	$$\forall s, \, \, (1-e)N_s \sim (1-e)\Nb_z $$ 
and $(1-e)\Nb_z$ is the monodromy operator of $\rho_{|E_w}$ 
as $\rho$ is irreducible, so the lemma follows from Prop. \ref{nilpoandmatrix} (ii).
\end{pf}

\subsubsection{Properties at $v$ and $\bar{v}$ of $\Ext_T(1,i)$.} 

\index{R@$\RR$, a refinement of $\pi^n_p$}
Let us assume from now that the accessible refinement $\RR$ of $\pi_p^n$ has been chosed of the form $$(1,\cdots,p^{-1}),$$
there are $n!$ such refinements by Lemma \ref{refofpin}. We fix 
for $\ast=v,\bar{v}$ some decomposition group map $\Gal(\overline{E}_\ast/E_\ast) \longrightarrow \G_{E,S}$.

\begin{prop}\label{fenp} For all $i \neq 1$, $\Ext_T(1,i)$ consists of extensions which are crystalline at $v$ and $\bar{v}$.
\end{prop}

\begin{pf} By Prop. \ref{Tisref}, $(X,T,\{\kappa_i\},\{F_i\},Z_{\rm reg})$ (see also Def. \ref{defnormbydelta}) is a refined family for the restriction $\Gal(\overline{E}_v/E_v) \rightarrow \G_{E,S}$. As a consequence, it induces a weakly refined family, and as $T^{\bot}=T\chi^{-1}$, the family
$$(X,T,-\kappa_m-1,F_m^{-1},Z_{\rm reg})$$  
is also a refined family for the restriction $\Gal(\overline{E}_{\bar v}/E_{\bar v}) \rightarrow \G_{E,S}$.
Note that by the choice of $\RR$, 
$$F_1(z)p^{\kappa_1(z)}=F_m(z)^{-1}p^{-\kappa_m(z)-1}=1.$$
By left exactness of the functor $D_{{\rm crys}}$, Lemma \ref{refofpin} and assumption (3), it suffices to check that for each extension $U$ in $\Ext_T(1,i)$, $$D_{{\rm crys}}(U_{|E_{\ast}})^{\varphi=1} \neq 0, \, \, \ast=v,\bar{v}.$$

But this follows from Theorem \ref{crystext} once we now that assumptions (ACC), (MF) and (REG) (for $v$ and $\bar{v}$) of \S\ref{kisinfamily} are satisfied. Assumption (ACC) follows from Lemma \ref{zregestzdense}, (MF) from \S\ref{faithgmachoix}, and (REG) from Prop. \ref{refofpin} again, and we are done.
\end{pf}

Recall the definition of $H^1_f$ from \S\ref{defselmer}.

\begin{cor}\label{corollf} For each $i \in {\cal I}$ with $i\neq 1$, $\Ext_T(1,i) \subset H^1_f(E,i)$.
\end{cor}

\begin{pf} It follows from Prop. \ref{fenp} and Prop. \ref{ext1cyclol}.
\end{pf}

\subsubsection{Symmetry properties of $T$}\label{symmgemastru} We choose now a particular GMA datum on $R$ using the symmetry of the pseudocharacter $T$ (see \S\ref{pseudosymmetry}). \ps
	Let $\tau: \OO(X)[\G_{E_S}] \rightarrow \OO(X)[\G_{E,S}]$ be the $\OO(X)$-linear map such that
		$$\tau(g):=cg^{-1}c^{-1}\chi(g).$$ 		
We have $\tau^2=1$ and $\tau(gg')=\tau(g')\tau(g)$, hence $\tau$ is an $\OO(X)$-linear anti-involution of $\OO(X)[\G_{E_S}]$. By \S\ref{normalisationrepgal}, it satisfies $$T \circ \tau = T.$$
As a consequence, $\tau$ induces an $A$-linear anti-involution on the $A$-algebra $R$ and we can apply to it the results of \S\ref{pseudosymmetry} (see Remark \ref{quotsymker}). The involution $\tau$ induces 
naturally an involution on ${\cal I}$ that we still denote by\footnote{This involution is denoted by $\sigma$ in \S\ref{pseudosymmetry}.} $\tau$, namely 
$$\forall i \in {\cal I}, \tau(i)=i^{\bot}\otimes \chi.$$
For example, we have $\tau(1)=\chi$. By Lemma \S\ref{idempotentspsi}, we can find a data of idempotents $\{e_i, i\in {\cal I}\}$ for the GMA $R$ such that 
\begin{equation} \label{equsymutile} \forall i \in {\cal I},\, \, \tau(e_i)=e_{\tau(i)}.\end{equation}
We choose now a GMA datum for $R$ of the form $\{e_i,\psi_i, i \in {\cal I}\}$ 
with the $e_i$ as above. It will be also convenient to fix an adapted representation $$R \hookrightarrow M_m(K)$$ 
associated to this datum in the sense of Theorem \ref{structure} (ii), so that $R$ identifies with the standard GMA of 
type $(1,1,d_1,\cdots,d_r)$ associated to $\{A_{i,j}, i,j \in {\cal I}\}$ where the $A_{i,j} \subset K$ are fractional ideals. In particular, 
each $A_{i,j}$ is finite type over $A$. \ps

\begin{lemma}\label{relsymijji} For all $i,j \in {\cal I}$, 
$A_{i,j}A_{j,i}=A_{\tau(i),\tau(j)}A_{\tau(j),\tau(i)}.$
\end{lemma}
\begin{pf} We have to show that $T(e_iRe_jRe_i)=T(e_{\tau(i)} R e_{\tau(j)} R e_{\tau(i)})$, which is 
immediat from the fact $T \circ \tau =T$, and that $\tau(e_{\ast})=e_{\tau(\ast)}$.
\end{pf}

 Let $i\neq j \in {\cal I}$. Recall that we defined in \S\ref{defbotij} a map $$\bot_{i,j}: \Ext_{k(z)[\G_{E,S}]}(i,j) \longrightarrow \Ext_{k(z)[\G_{E,S}]}(\tau(j),\tau(i)).$$
The following lemma is a consequence of Prop. \ref{mainsymmetry} and Lemma \ref{lestauij}.
	
\begin{lemma}\label{dualityapplication} The map $\bot_{i,j}$ induces an isomorphism 
$$\Ext_T(i,j) \isomo \Ext_T(\tau(j),\tau(i)).$$
\end{lemma}

In particular, using Prop. \ref{fenp} and \ref{ext1cyclol}, the lemma above 
has the following corollary. 

\begin{cor}\label{corfenpsym} For all $i \neq \chi$ in $\cal I$, $\Ext_T(i,\chi)$ consists of extensions which 
are crystalline at $v$ and $\bar{v}$, and split when restricted to $I_{E_w}$ for each $w$ prime to $p$.
\end{cor}

\subsection{Proof of Theorem \ref{signthm}}

\subsubsection{} Let us show Theorem \ref{signthm}. By Corollary \ref{corollf}, it suffices to show that for some irreducible subquotient $\rho_j$ of $\rho$, 
$$\Ext_T(1,\rho_j) \neq 0.$$

For that we will relate those Ext-groups to some reducibility ideal of $T$ and 
we will show that the associated reducibility locus is bounded. We first draw a key consequence of the vanishing of 
$H^1_f(E,\chi)$ and of the work above. Recall that we fixed in \S\ref{symmgemastru} a specific trace embedding $R \hookrightarrow M_m(K)$ 
identifying $R$ with the standard GMA of 
type $(1,1,d_1,\cdots,d_r)$ associated to $\{A_{i,j}, i,j \in {\cal I}\}$ where the $A_{i,j} \subset K$ are fractional ideals of $K$ (hence of finite type over $A$). 
\newcommand{\CA}{{\cal A}}
\begin{lemma}[{\it Vanishing of $H^1_f(E,\chi)$}]\label{vanih1ideal} We have $$\Ext_T(1,\chi)=0 \, \, \, \,{\rm  and}\, \, \, \, A_{\chi,1}=\sum_j A_{\chi,\rho_j}A_{\rho_j,1}.$$
\end{lemma}

\begin{pf} By Prop. \ref{propQ1}, $H^1_f(E,\chi)$ vanishes, hence so does $\Ext_T(1,\chi)$ by Cor. \ref{corollf}. Set 
$$A_{\chi,1}'=\sum_j A_{\chi,\rho_j}A_{\rho_j,1}.$$
By Theorem \ref{extension1} applied to $J=\m$, we get that $\Hom_A(A_{\chi,1}/A'_{\chi,1},k)=0$. 
But $A_{\chi,1}$ has finite type over $A$, hence we conclude by Nakayama's lemma.
\end{pf}

Let $\cal P$ be the following partition of ${\cal I}$: $\cal P=(\{\chi\},\{1\},{\cal I}\backslash \{\chi,1\}),$
and $I_{\cal P} \subset A$ its reducibility locus (see Def. \ref{defredloc}).

\begin{lemma}[{\it Reduction to the nonvanishing of $I_{\cal P}$}]\label{redlocuIP}${}^{}$\par
\begin{itemize}
\item[(i)] $I_{\cal P}=A_{1,\chi}A_{\chi,1}+\sum_{j} A_{1,\rho_j}A_{\rho_j,1}+\sum_{j} A_{\chi,\rho_j}A_{\rho_j,\chi}$.
\item[(ii)] $I_{\cal P}=\sum_j A_{1,\rho_j}A_{\rho_j,1}$.
\item[(iii)] If $\Ext_T(1,\rho_j)=0$ for all $j$, then $I_{\cal P}=0$.
\end{itemize}
\end{lemma}
\begin{pf} Assertion (i) is Prop. \ref{redloc}. As $\tau(1)=\chi$, Lemma \ref{relsymijji} shows that  
$$\sum_{j} A_{1,\rho_j}A_{\rho_j,1}=\sum_{j} A_{\chi,\rho_j}A_{\rho_j,\chi}.$$
But by Lemma \ref{vanih1ideal}, $A_{\chi,1}=\sum_j A_{\chi,\rho_j}A_{\rho_j,1}$, hence 
$$A_{\chi,1}A_{1,\chi}=\sum_j A_{\chi,\rho_j}A_{\rho_j,1}A_{1,\chi} \subset \sum_j A_{\chi,\rho_j}A_{\rho_j,\chi},$$
as $A_{\rho_j,1}A_{1,\chi} \subset A_{\rho_j,\chi}$. This proves assertion (ii).\par
Arguing as in the proof of Lemma \ref{vanih1ideal}, $\Ext_T(1,\rho_j)=0$ implies that 
$$A_{\rho_j,1}=\sum_{i \neq 1, \rho_j} A_{\rho_j,i}A_{i,1}.$$	
Assume that this holds for each $j$. By applying this identity or Lemma \ref{vanih1ideal} $s:=1+|\cal I|$ times, 
we get that for each $j$, $A_{1,\rho_j}A_{\rho_j,1}$ is a finite sum
of terms of the form 
\begin{equation} \label{productdegeulasse} A_{1,\rho_j}A_{\rho_j,i_1}A_{i_1,i_2}\cdots A_{i_s,1},\end{equation}
with all the $i_k \in \cal I$, $i_1 \neq \rho_j$, $i_s \neq 1$, and $i_{k+1} \neq i_k$. As $s>|\cal I|$, 
for each such term there exist $k < k'$ such that $i_k=i_{k'}$, which implies that 
$$A_{i_k,i_{k+1}}A_{i_{k+1},i_{k+2}} \cdots A_{i_{k'-1},i_{k'}} \subset A_{i_k,i_{k+1}}A_{i_{k+1},i_k} \subset \m,$$
$$A_{\rho_j,i_1}A_{i_1,i_2}\cdots A_{i_{k-1},i_k}A_{i_{k'},i_{k'+1}} \cdots A_{i_s,1} \subset A_{\rho_j,1},$$
hence that $$A_{1,\rho_j}A_{\rho_j,i_1}A_{i_1,i_2}\cdots A_{i_s,1} \subset \m A_{1,\rho_j}A_{\rho_j,1}.$$
This proves that $I_{\cal P} \subset \m I_{\cal P}$, hence that $I_{\cal P}=0$ by Nakayama's lemma.
\end{pf}

By Lemma \ref{redlocuIP} (iii), it only remains to show part (ii) of 
the following lemma. Recall that $\{\kappa_i(z)\}_{i=1\dots m}$ is the strictly increasing (see Def. \ref{defonkappai}) sequence of Hodge-Tate weights of $\rhob_z$ at $v$. 
Let $a \in \{1,\dots,m\}$ be the unique integer such that 
	$$\kappa_a(z)=0.$$

\begin{lemma}[{\it Non triviality of $I_{\cal P}$}]\label{nontotredu}${}^{}$\par
\begin{itemize}
\item[(i)] $(\kappa_a - \kappa_1) -(\kappa_a(z)-\kappa_1(z)) \in I_{\cal P}$.
\item[(ii)] $a\neq 1$ and $I_{\cal P}\neq 0$.
\end{itemize}
\end{lemma}

\begin{pf} As already said in the proof of Prop. \ref{fenp}, $(X,T,\kappa_1,F_1,Z_{\rm reg})$ is a weakly refined family for $\G_{|E_v} \rightarrow \G_{E,S}$, and the assumption (ASS), (MF) and (REG) of \S\ref{kisinfamily} are satisfied. Part (i) is then Theorem \ref{redlocfaible} as $$D_{{\rm crys}}(\rhob_z)^{\varphi=1}=D_{\rm crys}(1)^{\varphi=1}$$
has dimension $1$ by Lemma \ref{refofpin} and assumption (3).	\par
Let us show assertion (ii). By property (iv) of the eigenvariety $X$, the natural 
map $\OO_{\kappa(z)} \longrightarrow A$ is injective, so it suffices to check 
that $$(\kappa_a-\kappa_1)-(\kappa_a(z)-\kappa_1(z)) \neq 0,$$ {\it i.e.} that $a \neq 1$. If $a=1$, $\kappa_a(z)=0$ is the smallest Hodge-Tate weight of $\rhob_z$ at $v$. But this is absurd as $-1$ is a Hodge-Tate weight of $\rhob_z$ at $v$, namely the one of $\chi$.
\end{pf}

\subsubsection{Some remarks about the proof} The proof above of Theorem \ref{signthm} 
can actually be simplified in several different ways. We chosed to look at the full minimal eigenvariety $X$ containing $z$ and its associated Galois
pseudocharacter $T$ because this is the relevant point and space for which the analysis developped here 
can be pushed further (and for which the $\Ext_T$ have a maximal dimension) as we will explain in the next section. 
All the strenght of the results proved here (especially the ones in 
\S\ref{mineigenprepa}) will be used in section 9, and we found it convenient to directly include them here so 
as not to repeat half of the story there. \par 
	In the style of \cite{BC}, we could have replaced $X$ by the normalization of the germ of any irreducible 
curve $C \subset X$ containg $z$ such that $Z \cap C$ is infinite and that $\kappa_a-\kappa_1$ is not constant on $C$. 
The ring $A$ would have been a DVR which would have simplified some of the pseudocharacter theoretic arguments. Note that in the argument above, 
we do not really choose a stable "lattice" as in \cite{BC} but we rather work with the full ring theoretic image $R$ of 
the family of Galois representations. This is actually convenient and it illustrates the techniques developped in the previous sections of this book. 
Had we worked on the germ of a smooth curve as explained above, we could have used the choice of a good lattice 
as in \cite[Prop. 7.1]{BC} (as written, it requires $\rho$ to be irreducible). \par
	Moreover, a nice way to understand the combinatorics in (iii) of Lemma \ref{redlocuIP} is to compare it 
with the connected graph theorem \cite[Thm. 1]{bellaiche}.\footnote{Here is the statement: let $A$ be a henselian
DVR, $r: G \longrightarrow \GL_d(A)$ a generically absolutely irreductible representation, and assume that the semi-simplified 
residual representation $\bar{r}^{\rm ss}$ is multiplicity free. Then the oriented graph whose vertices are the irreducible 
constituents of $\bar{r}^{\rm ss}$, and with an edge from $i$ to $j$ if there is a nontrivial extension of $i$ by $j$ in the 
subquotients $A$-modules of $A^d$, is connected as an oriented graph.} In our case, there would be no edge $1 \rightarrow \chi$ by  Lemma \ref{vanih1ideal}, hence at least an edge $1 \rightarrow \rho_j$ for some $j$. 
Note that we do not claim that the pseudocharacter $T$ used in the proof above is generically irreducible, 
but Lemma \ref{nontotredu} rather says that it is not "too reducible", and this is actually enough to conclude. 
Actually, had we assumed that the eigenvalues of the Langland's conjugacy class of $\pi_p$ are "regular", we could 
have chosed a refinement $\RR$ (hence a $z$) leading to a generically irreducible $T$ (even on the curve $C$), as follows 
from Rem. \ref{remarkneededinsection8}.\par
	As is clear from the proof, the cornerstone of the argument is the fact that 
$$\Ext_T(1,\chi) \subset H^1_f(E,\chi),$$
(that is Lemma \ref{vanih1ideal}) which requires to control the deformation at {\it all the finite places}, 
from which we deduced that $\Ext_T(1,\chi)=0$ 
using the finiteness of $\OO_E^*$ (in terms of the graph alluded above, it is the step: "{\it there is no arrow $1 \rightarrow \chi$}"). 
This last fact fails for a general CM field $E$, and actually the whole argument breaks down in this 
generality because of that. We will discuss that issue in greater details in the 
subsequent Remark \ref{caseCMfield}.

\newpage

\section{The geometry of the eigenvariety at some Arthur's points and higher rank Selmer groups}

\subsection{Statement of the theorem}

We keep the notations of \S\ref{statemtthmsign}. In particular 
$$\rho: \Gal(\overline{E}/E) \longrightarrow \GL_n(L)$$
is a {\it modular} Galois representation attached to a cuspidal 
automorphic representation $\pi|.|^{1/2}$ of $\GL_n(\AAA_E)$ that satisfy conditions (1), (2) and (3) there. 
We now make the following new assumptions on $\rho$ :\ps
\begin{itemize}
\item[(4)] $\Lambda^i \rho$ is absolutely irreducible for $i=1,\dots,n$.\ps
\item[(5)] the crystalline representation $\rho_{|E_v}$ admits a regular non critical refinement 
(see \S\ref{crit} and Def. \ref{reg}),\ps
\item[(6)] the hypotheses $BK1(\rho)$ and $BK2(\rho)$ hold (see \S\ref{auxselgroup}).
\end{itemize}

\begin{remark}\label{remassumptionmainthm} \begin{itemize}
\item[(i)] The irreducibility assumption (4) is known if $n \leq 3$. 
\item[(ii)] Recall that the regularity assumption in (5) combined with (3) means that the Langlands conjugacy class 
$c \in \GL_n(\C)$ of the unramified representation $\pi_p|.|^{1/2}$ has distincts eigenvalues, 
and that those eigenvalues can be ordered as $$(\varphi_1,\cdots,\varphi_n)$$ 
in such a way that for each $j=1,\cdots,n$, 
$\varphi_1\cdots\varphi_j$ is a simple eigenvalue of $\Lambda^j(c)$. If $n\leq 3$, it is equivalent 
to only ask that the $\varphi_i$ are distinct, and if $n=2$ (resp. $n=1$) this is 
conjectured to always be the case (resp. it is obviously true). \par The non critical part of the assumption of (5) means 
that the refinement of $\rho_{|E_v}$ associated by property (3) to the ordering above is non critical in the sense of \S\ref{crit}. 
Again, this is automatically satisfied if $n=1$, and in most cases when $n=2$ (see Remark \ref{remcrit}).
\item[(iii)] As we saw in Propoositions~\ref{propBK1} and~\ref{propBK2}, the hypothesis (6) is know to hold if $n=1$ 
and also in the $n=2$ case for $\rho$ of the form $\rho_{f,E}$ for $f$ a modular forms of
even weight with a small explicit set of exceptions.
\end{itemize}
\end{remark}

Of course, we will also assume that $\varepsilon(\rho,0)=-1$, and that Rep($m$) and AC($\pi$) hold. As in \S\ref{constructionpin}, we denote by $\pi^n$ 
the non-tempered automorphic representation of $\U(m)$ attached to
$\rho$ by assumption AC($\pi$), for some choice of a Hecke character $\mu$ as in Def. \ref{choixdemu} that we fix once for all.
Recall that we defined in \S\ref{mineigenprepa} the minimal eigenvariety $X$ of $\U(m)$ containing $\pi^n$. We consider here the 
variant where we fix one of the $m$ weights (anyone), so that $X$ is equidimensional of dimension $m-1=n+1$. Any choice of 
an accessible refinement $\RR$ of $\pi^n_p$ defines a point $z \in X$. By Lemma \ref{refofpin} and assumption (5), we may choose a refinement of the form 
$$\mu_v|.|^{-1/2}(p)(1,\varphi_1,\varphi_2,\cdots,\varphi_n,p^{-1}),$$
where the $\varphi_i$ are chosen as in Remark \ref{remassumptionmainthm} to satisfy the regularity and non criticality 
assumption of (5). We fix once for all such a refinement $\RR$ \index{R@$\RR$, a refinement of $\pi^n_p$}, hence a point $$z \in X,$$
which is the {\it Arthur's} point that we refer to in the title, $z$ is defined over $L$. \ps

	In all this section, we will generally follow the notations of \S\ref{mineigenprepa}. 
In particular, recall that we defined in Def. \ref{extT} an $L$-subspace $\Ext_T(1,\rho) \subset \Ext_{\G_{E,S}}(1,\rho)$ 
which is the space of extensions of $1$ by $\rho$ that we can construct from the Galois 
pseudocharacter $T$ carried by $X$ (see Remark \ref{remextT}).
By Prop. \ref{fenp} and \ref{corfenpsym}, we know that 
	$$\Ext_T(1,\rho) \subset H^1_f(E,\rho).$$

\begin{theorem}\label{mainmainthm} Assume that $\rho$ satisfies (1) to (6), as well as AC($\pi$) and Rep($m$). Let $t$ be the dimension of 
the tangent space of $X$ at $z$ and $h$ the dimension of $\Ext_T(1,\rho)$, then  
$$t \, \,\,\leq\, \,\, h\, (\,n+\,\frac{h+1}{2}\,).$$
\end{theorem}
\ps
Note that both dimensions above are taken over the residue field $k=k(z) \simeq L$. As $\OO_{X,z}$ is equidimensional of dimension $n+1$, we have $$n+1 \,\,\leq\,\, h\,(n+\frac{h+1}{2})$$ 
so we recover in particular that $h \neq 0$ (i.e Theorem \ref{signthm} for $\rho$) and get the following corollary.

\begin{cor}\label{corlissitehigher} (same assumption) If $X$ is not smooth at $z$, then 
$$\dim_{L}H^1_f(E,\rho) \geq \dim_{L}\Ext_T(1,\rho) \geq 2.$$
\end{cor}
\medskip
When $n=\dim(\rho)=1$, {\i.e.} when $$\rho: \Gal(\overline{E}/E) \longrightarrow L^*$$
is any continuous character such that $\rho^{\bot}=\rho(-1)$, class field theory and Remark \ref{remassumptionmainthm} imply 
that conditions (1) to (6) are satisfied once we assume that the two Hodge-Tate weights of $\rho_{|E_v}$ are different from $0$ or $-1$ 
(because of condition (ii) on $\pi$ in \S\ref{start}). Moreover, by Remarks \ref{remsectrepm} (vi) and \ref{remsectAC} (ii), 
assumptions Rep($3$) and AC($\pi$) are also known.

\begin{cor}\label{cormainmain1dim} If $n=1$, Theorem \ref{mainmainthm} and its corollary above hold under the single assumption that $0$ and $-1$ are not 
Hodge-Tate weights of the character $\rho_{|E_v}$.
\end{cor}

\medskip

\begin{remark} \label{mainconjvsbk}As we already said, various versions of the main conjectures are known for Hecke characters from 
the work of Rubin \cite{rubinmc}. These "main-conjectures" also imply that $\dim_L H^1_f(E,\rho) \geq 1$ when expected. However, as far as we know, they do not allow to show that $\dim_L H^1_f(E,\rho) \geq 2$
when the $L$ function or rather its $p$-adic analogues vanish at a higher order (this phenomenum is sometimes called the possible 
{\it non semisimplicity of the Iwasawa module}). As a consequence, even the simplest case covered by the corollary above is of interest. 
\par In all this book, we have concentrated on the Galois side part of the study, letting aside the various 
$p$-adic $L$-functions that should enter in the picture. We hope that once this will be done, the Galois 
deformations studied here will shed some light on the $\leq$ part of the conjectures alluded above. From a more conjectural point of view, 
this remark actually applies to any $\rho$ satisfying conditions (1) and (2).
\end{remark}

\bigskip

\subsection{Outline of the proof} 
The proof of Theorem \ref{mainmainthm} is a refinement of the one of the sign conjecture consisting in 
a careful analysis of the Galois pseudocharacter $$T: \G_{E,S} \longrightarrow \OO_{X,z}$$ at the point $z$. 
As in \S\ref{faithgmachoix}, we let $A=\OO_{X,z}$, $\m$ the maximal ideal of $A$, $k=A/\m \simeq L$ its residue field, and 
$$R:=A[G]/\Ker T$$ 
the faithful Cayley-Hamilton GMA associated to $T$ at $z$. 
Recall that it is finite type and torsion free over $A$ 
(which is a reduced henselian noetherian ring). As $\rho$ is irreducible, we have now
$$\cal I=\{ \chi, \rho, 1\}.$$
As in \S\ref{symmgemastru}, we fix a data of idempotents $(e_{\chi},e_{\rho},e_1)$ for $R$ such that 
$\tau(e_{\ast})=e_{\sigma(\ast)}$. Note that $\sigma$ fixes $\rho$, and exchanges $1$ and $\chi$. Last but not least, 
$K=\prod_s K_s$ is the total fraction ring of $A$, and we fix a 
representation $$\rho_K: R \longrightarrow M_m(K)$$
associated to this data of idempotents, as in Theorem \ref{structure} (ii). Recall that this gives us a 
set of finite type $A$-modules $A_{i,j} \subset K$, $i,j\in \cal I$, such that $A_{i,i}=A$ for each $i$, 
$A_{i,j}A_{j,k} \subset A_{i,k}$ for each $i,j,k$ (the {\it Chasles relations}), and $A_{i,j}A_{j,i} \subset \m$ if $i\neq j$. Moreover,
$R=\rho_K(R) \subset M_m(K)$ is the standard GMA of type $(1,n,1)$ associated to these data (see Example \ref{gma}), that is

\begin{equation}\label{forme3} R=\rho_K(R)=\left( \begin{array}{ccc} A & A_{\chi,\rho}^n  & A_{\chi,1} \\
A_{\rho,\chi}^n & M_n(A) & A_{\rho,1}^n \\
A_{1,\chi} & A_{1,\rho}^n & A \end{array} \right) \subset M_m(K).\end{equation}

Our aim will be to elucidate as much as possible the structure of the $A$-modules $A_{i,j}$. \ps
One the one hand, those $A_{i,j}$ are related to the $\Ext_T(j,i)$ by Theorem \ref{extension1}.
In turns, by results already proved in \S\ref{mineigenprepa}, those $\Ext_T(j,i)$ are related to all the $6$
possible fine Selmer groups occuring here, namely $H^1_f(E,\ast)$ where 
$$\ast= \chi,\chi^{-1},\rho,\rho^*,\rho\chi^{-1},\rho^*\chi^{-1}.$$
Up to the underlying symmetries, and by $BK_1(\rho)$ when $\ast=\rho^*$, 
we will actually know all of them except the one of $\rho$, which is precisely the 
one we are interested in. \ps
	On the other hand, the $A_{i,j}$ are also related to the reducibility loci of $T$ by Prop. 
\ref{redloc}. A remarkable fact is that we are able to compute here all the reducibility ideals of $T$. Precisely, 
we will show that all the proper reducibility loci actually coincide schematically with the closed point $z$. 
In other words, $T$ is as irreducible as possible. The proof of this key fact will actually use all the machinery that we 
developped in sections 1 to 4 of the book. This will provide then the missing link between
the tangent space of $X$ at $z$ and the $A_{i,j}$, and then with the $\Ext_T(1,\rho)$.\par

\subsection{Computation of the reducibility loci of $T$}

Let us analyse the proper reducibility loci of $T$ (see \S\ref{assredloc}). 
Recall that each of them is attached to a non trivial partition $\cal P$ 
of $\cal I=\{\chi,1,\rho\}$, and there are $4$ such partitions. An especially interesting one is 
the the total reducibility ideal $I_{\rm tot}$, which is
attached to the finest partition $\{\{\chi\},\{1\},\{\rho\}\}$.

\begin{lemma}\label{computeIP} ${}^{}$\par
\begin{itemize}
\item[(i)] All these four reducibility ideals coincide with $I_{\rm tot}$.
\item[(ii)] $I_{\rm tot}=A_{1,\rho}A_{\rho,1}$.
\item[(iii)] $I_{\rm tot}K=K$.
\end{itemize}
\end{lemma}

\begin{pf} By Prop. \ref{redloc}, each proper reducibility ideal is a sum of terms of the form 
$A_{i,j}A_{j,i}$ with $i\neq j$, and contains $A_{\ast,\rho}A_{\rho,\ast}$ for $\ast=1$ or $\chi$. By Lemma \ref{relsymijji}, $A_{i,j}A_{j,i}=A_{\sigma(j),\sigma(i)}A_{\sigma(i),\sigma(j)}$ so 
\begin{equation}\label{relsymijmain} A_{1,\rho}A_{\rho,1}=A_{\chi,\rho}A_{\rho,\chi}. \end{equation}
By Lemma \ref{vanih1ideal}, we also have 
\begin{equation}\label{h1fchizeromain} A_{\chi,1}=A_{\chi,\rho}A_{\rho,1}. \end{equation}
But by the Chasles relation $A_{1,\chi}A_{\chi,\rho} \subset A_{1,\rho}$, so 
$$A_{1,\chi}A_{\chi,1} \subset A_{1,\rho}A_{\rho,1},$$
which proves assertions (i) and (ii). Part (iii) follows from (i) and Lemma \ref{nontotredu} (ii) as 
$\kappa_a-\kappa_1-(\kappa_a(z)-\kappa_1(z)) \neq 0$ in the domain $\OO_{\kappa(z)}$, and as the composition of the maps
$$\OO_{\kappa(z)} \longrightarrow A \longrightarrow K_s$$
is injective by propety (iv) of the eigenvariety $X$.
\end{pf}

\begin{lemma}\label{corcomputeIP} $\rho_K$ induces an isomorphism $R \otimes_A K \isomo M_m(K)$, and
$\rho_K \otimes K_s$ is absolutely irreducible for each $s$. \end{lemma}

\begin{pf} This is actually a general consequence of Prop. \ref{specialcase} and of the fact that $IK=K$ for all irreducibility ideals $I$,
but we argue directly. By Lemma \ref{computeIP} (ii) and (iii), $A_{1,\rho}K=A_{\rho,1}K=K$, 
and the same equality holds with $1$ replaced by $\chi$ by formula (\ref{relsymijmain}).
As $A_{1,\chi} \supset A_{1,\rho}A_{\rho,\chi}$ we get also that $A_{1,\chi}K=K$, as well as $A_{\chi,1}K=K$ by the same reasoning. 
This proves the first part of the lemma, of which the second part is an obvious consequence.
\end{pf}

Note that for the moment, we did not use assumptions (5) to (6), and only the irreducibility assumption of $\rho$ in (4). 
We will now use (4) and (5) by beginning a deeper study of $I_{tot}$. We show first that 
the total reducibility loci $$V(I_{\rm tot}) \subset {\rm Spec}(\OO_z)$$ lies in the schematic fiber of the weight morphism 
$\kappa: X \rightarrow \WW$ over $\kappa(z)$. 

\begin{prop}\label{Itotfiber} For each integer $j \in \{1,\dots,m\}$, $\kappa_j-\kappa_j(z) \subset
I_{{\rm tot}}$. 
\end{prop}

To prove this proposition, we will need to recall some aspects of refined deformations that we 
developped in \S\ref{hypredloc}. By definition of the chosed refinement
$\RR$, the refinement $\Ref_z$ of $\rhob_z=1\oplus\rho\oplus \chi$ is
	$$\Ref_z =\iota_p\iota_{\infty}^{-1}(1,\varphi_1,\cdots,\varphi_n,p^{-1}).$$
This makes sense as by assumption (5) and Lemma \ref{refofpin}, the $m$ Frobenius eigenvalues
above are distinct. Of course, this refinement induces also a refinement $\Ref_{z,\ast}$ of each
$\rhob_{\ast}$: $\Ref_{z,1}=(1)$, $\Ref_{z,\chi}=(p^{-1})$ and
$\Ref_{z,\rho}=\iota_p\iota_{\infty}^{-1}(\varphi_1,\cdots,\varphi_n)$. \par
	Recall that in this situation, we attached in \S\ref{permutation} a permutation 
$$\sigma \in \got{S}_m$$
that encaptures how the indices $i$ of the weights $\kappa_i(z)$ and the Frobenius eigenvalues
$F_i(z)p^{\kappa_i(z)}$ are related to the decomposition $\rhob_z=1\oplus
\rho \oplus \chi$: \begin{itemize}
\item[-] $R_{\ast}$ is the set of integers $i$ such that
$D_{\rm crys}(\rhob_{\ast})^{\varphi=F_i(z)p^{\kappa_i(z)}} \neq 0$. 
\item[-] $W_{\ast}$ is the set of integers $i$ such that $\kappa_i(z)$ is a Hodge-Tate weight
of $\rhob_{\ast}$.
\end{itemize}

\begin{lemma}\label{lemmeantiordref}${}^{}$
\begin{itemize}
\item[(i)] $\sigma$ is
a transitive permutation. 
\item[(ii)] $\Ref_z$ is a critical regular refinement of $\rhob_z$. 
However, for each $\ast$, $\Ref_{z,\ast}$ is a non critical regular refinement of
$\rhob_{\ast}$ and $R_{\ast}$ is a subintervall of $\{1,\cdots,d\}$.
\end{itemize}
\end{lemma}

\begin{pf} Let us show assertion (ii) first. As $(1,p^{-1})$ is a critical
refinement of $1 \oplus \chi$, $\Ref_z$ is a critical refinement of
$\rhob_z$. For the regularity property, let us fix $j\geq 1$ an integer. By Lemma \ref{refofpin}, the eigenvalues $\lambda$ 
of the crystalline Frobenius on $D_{\rm crys}(\Lambda^j \rhob_z)$ such that
$|\iota_{\infty}\iota_p^{-1}(\lambda)|=\sqrt{p}^{-j+1}$ are
exactly the products of $j-1$ elements of
$\iota_p\iota_{\infty}^{-1}(\{\varphi_1,\cdots,\varphi_n\})$. We conclude
then by assumption (5). The non criticality of $\Ref_{z,\ast}$ is obvious
for $\ast=1,\chi$ and is assumption (5) for $\ast=\rho$. Moreover, the
assertion on $R_{\ast}$ is clear, namely:
\begin{equation}\label{formuleRast} R_{1}=\{1\},\, \, \,
R_{\rho}=\{2,3,\dots,m-1\},\, \, \, R_{\chi}=\{m\}.\end{equation}
	We show now (i). Let $a \in \{1,\cdots,m\}$ be the unique integer
such that $\kappa_a(z)=0$. As already said, $a>1$ (see Lemma \ref{nontotredu} (ii)), and we have 
\begin{equation}\label{formuleWast} W_{1}=\{a\}, \, \, \,
W_{\rho}=\{1,2,\dots,a-2,a+1,\dots,m\}, \, \, \,
W_{\chi}=\{a-1\}.
\end{equation}
By definition of the permutation $\sigma$, we see now that 

$$\sigma(i)=\left\{ \begin{array}{llll} a,\, \, {\rm if}\, \, \, i=1,\\
i-1,\, \, {\rm if}\, \, \, i=2,\cdots,a-1,\\
i+1,\, \, {\rm if}\, \, \, i=a,\cdots,m-1,\\ 
a-1,\, \, {\rm if}\, \, i=m. \end{array} \right\},$$
which is a cycle, and we are done.
\end{pf}

\begin{pf} (of Prop. \ref{corItotfiber}) As one of the $\kappa_i$ is constant by assumption on $X$, the proposition is an immediat
consequence of Lemma \ref{lemmeantiordref} and Corollary \ref{corredloc}, once
we know that assumptions (REG), (NCR), (INT) and (MF') of \S\ref{hypredloc} are
satisfied. But the first three ones are by Lemma \ref{lemmeantiordref} (ii),
and (MF') follows from assumption (4) and from the fact that $1$ and
$\chi$ are one dimensional.
\end{pf}

\begin{cor}\label{corItotfiber} $I_{\rm tot}$ is a cofinite lenght ideal of $A$.
\end{cor}	

\begin{pf} It follows from Prop. \ref{Itotfiber} and from the fact that the natural map $\OO_{\kappa(z)} \rightarrow \OO_z=A$ is finite 
by property (iv) of the eigenvariety $X$. 
\end{pf}

Recall from Def. \ref{lesrhoi} that for each $\ast \in \cal
I$, there is a (unique up to isomorphism) continuous representation $$\rho_{\ast}: \G_{E,S} \longrightarrow
\GL_{d_{\ast}}(A/I_{\rm tot})$$
lifting $\rhob_{\ast}$. 

\begin{prop} \label{cristproperties} ${}^{}$ 
\begin{itemize}
\item[(i)] For each $\ast \in \{1,\chi,\rho\}$, $\rho_{\ast}$ is
crystalline at $v$ and $\bar{v}$. 
\item[(ii)] Moreover, the characteristic polynomial of
the crystalline Frobenius on the free $A/I_{\rm tot}$-module $D_{\rm
crys}(\rho_{\ast|E_v})$ is $$\prod_{i \in R_\ast}(T-F_ip^{\kappa_i(z)}) \in
(A/I_{\rm tot})[T].$$
\end{itemize}
\end{prop}

\begin{pf} Note that $\Hom_{G_{E_v}}(\rhob_{\ast},\rhob_{\ast}(-1))=0$ for
each $\ast$. Indeed, it is clear for $\ast=1,\chi$, and it holds for
$\ast=\rho$ as $\varphi_i \neq p\varphi_j$ for each $i,j \in \{1,\dots,n\}$
by Lemma \ref{refofpin}. The first part of the proposition for the place $v$ follows then from Corollary
\ref{corredloc} (ii) (we already checked that (INT), (REG), (NCR) and
(MF') hold in the proof of Prop. \ref{corItotfiber}). But 
$\rho_{\ast}^{\bot}\simeq \rho_{\tau(\ast)}$ for each $\ast$, as they
share the same trace and each $\rho_{\ast}$ is residually irreducible, so
the proposition also holds for the place $\bar{v}$, which proves (i). \ps
By Theorem \ref{tridefgen}, we know that the crystalline representation 
$\rho_{\ast}$ is trianguline over $A$ with parameters the $\delta_i$ with $i \in
R_\ast$ such that 
$${\delta_i}_{|\Z_p^*}=\chi^{-\kappa_{\sigma(i)}(z)}, \, \, \,
\delta_i(p)=F_ip^{\kappa_i(z)-\kappa_{\sigma(i)}(z)} \in (A/I_{\rm tot})^*.$$ 
By Berger's theorem \ref{bergerthm}, the characteristic polynomial of the statement
writes as the products over the $i\in R_{\ast}$ of the characteristic polynomial of $\varphi$ on 
$D_{\rm crys}(\Ro_{A/I_{\rm tot}}(\delta_i))$, hence the statement. 
\end{pf}

We now come to the main proposition of this subsection. We will use here
assumption $BK_2(\rho)$ of hypothesis (6).

\begin{prop}\label{itotestmax} $I_{\rm tot}$ is the maximal ideal of $A$. \end{prop}

\begin{pf} Note that the residue field $k:=k(z)$ of $A$ lifts canonically to a subfield
of $A$ by the henselian property. Let us fix a 
$$\psi: A/I_{\rm tot} \longrightarrow k[\varepsilon],$$
a $k$-linear ring homomorphism. We claim that for each $\ast$, $$\rho_{\ast,\psi}:=\rho_{\ast}
\otimes_{A/I_{\rm tot},\psi} k[\varepsilon]$$ is a trivial
deformation of $\rhob_{\ast}$, which means that we have an isomorphism
$$\rho_{\ast,\psi} \simeq \rhob_{\ast} \otimes_k
k[\varepsilon].$$ 
	Let us assume this claim and show how to
conclude. By property (i) of eigenvarieties (see Def. \ref{eigendef}), $A$
is generated by $\HH$ as an $\OO_{\kappa(z)}$-algebra. As $$\HH=\ATL
\otimes \HH_{\rm ur}$$ and by assumption (2), we see that 
$A$ is generated over $\OO_{\kappa(z)}$ by the $F_i$'s and by 
the $T(\Frob_w)$'s for the primes $l=w\bar{w} \in S_0$. Assuming that each
$\rho_{\ast,\psi}$ is constant, we get that for any such $w$,
$$\psi(T(\Frob_w))\in k \subset
k[\varepsilon]$$ is constant.
Moreover, Prop. \ref{cristproperties} (ii) implies that for each $i$,
$$\psi(F_i) \in
k \subset k[\varepsilon]$$ is also constant (use that the
$F_i(z)p^{\kappa_i(z)}$ are two-by-two distinct by Lemma
\ref{lemmeantiordref}). Last but not least, by Prop. \ref{Itotfiber} the image of 
	$$\OO_{\kappa(z)} \longrightarrow A/I_{\rm tot}
\longrightarrow_{\psi} k[\varepsilon]$$
also falls into $k$. As $A$ is generated over $\OO_{\kappa(z)}$ by the
$F_i$ and the $T(\Frob_w)$, we get that $$\psi(A/I_{\rm tot})=k.$$ 
As this holds for all $\psi$, $A/I_{\rm tot}=k$ and we are done. \par
	Let us show the claim now. By Prop. \ref{cristproperties}, we know
that $\rho_{\ast,\psi}$ is crystalline at $v$ and $\bar{v}$. Moreover,
$\rho_{\ast,\psi}$ is obviously unramified outside $S$. By Lemma
\ref{lemmemonopourplustard} (applied to $J=\ker \psi$) we know that for each prime $w$ of $E$ not
dividing $p$, the monodromy operator of $\rho_{\ast,\psi|E_w}$ admits a Jordan
normal form over $A/I_{\rm tot}$, hence is constant, when $\ast=\rho$. This
trivially also holds when $\ast=1$ or $\chi$, as any continuous $\G_{E_w}$-extension
of $1$ by $1$ is unramified for such a $w$. \par 
	If $\ast=1$ or $\chi$, the
finiteness of the class number of $E$, and more precisely Prop. \ref{propQ2}
(i), implies then that $\rho_{\ast,\psi}$ is constant. If $\ast=\rho$, we have
$\rho_{\ast,\psi}^{\bot}=\rho_{\ast,\psi}$ (see the first paragraph of the proof of
Prop. \ref{cristproperties}), hence hypothesis $BK_2(\rho)$ in assumption
(6) shows again that $\rho_{\ast,\psi}$ is constant, which completes
the proof. \end{pf}

\begin{remark}\label{remfinpreuveirrloc} We could also study the proper reducibility loci of the restriction of 
$T$ to $\G_{E_v}$. For example when $\rho_{|E_v}$ is irreducible (e.g. when $n=1$), the same proofs as above show that 
they all coincide and that they lie in the schematic fiber of $\kappa$ above
$\kappa(z)$. However, they do not necessarily coincide with the maximal ideal of $A$.
\end{remark}

\subsection{The structure of $R$ and the proof of the theorem} For $i\neq j \in \cal I$, let us consider the integers $$h_{i,j}:= \dim_L \Ext_T(i,j).$$
We first recapitulate all that we know about those $h_{i,j}$. 
\begin{lemma}\label{valeurshij} ${}^{}$
\begin{itemize}
\item[(i)] $h_{1,\rho}=h_{\rho,\chi}=h$, 
\item[(ii)] $h_{1,\chi}=0$ and $h_{\chi,1}\leq 1$, 
\item[(iii)] $h_{\rho,1}=h_{\chi,\rho}\leq n$.
\end{itemize}
\end{lemma} 

\begin{pf} The first equalities in (i) and (ii) follow from Lemma \ref{dualityapplication}, which proves (i). Assertion (iii) is then a consequence of Prop. \ref{ext1cyclol} and of Prop. \ref{calulhunfrhodual} (which assumes hypothese $BK_1(\rho)$ and whose assumptions are satisfied by Lemma \ref{refofpin}). \par 
	We already proved that $h_{1,\chi}=0$ in Lemma \ref{vanih1ideal}, so it only remains to show that $h_{\chi,1}\leq 1$. That will follow from Prop. \ref{propQ2} (ii) if we can show that $$\Ext_T(\chi,1) \subset H^1(E,L(-1))$$
falls into an eigenspace of the endomorphism $U \mapsto U^{\bot}(-1)$ of the latter space. But this follows from Lemma \ref{corcomputeIP} and Prop. \ref{signebotijfinal} as $\tau$ fixes $\rho \in \cal I$. We will actually show later that $\Ext_T(\chi,1)$ falls inside the part of sign $+1$.
\end{pf}

As $|\cal I|=3$, recall that from Theorem \ref{extension1} that for $i, j$ and $k$ two-by-two distinct in $\cal I$, we have an isomorphism
\begin{equation}\label{linkextaij} \Hom_L(A_{i,j}/A_{i,k}A_{k,j},L) \isomo \Ext_T(j,i).\end{equation}

\begin{lemma}\label{linkh} $h$ is the minimal number of generators of $A_{\rho,1}$.
\end{lemma}

\begin{pf} By (\ref{linkextaij}) above and Nakayama's lemma, we have to show that $$A_{\rho,\chi}A_{\chi,1} \subset \m A_{\rho,1}.$$
But $A_{\chi,1}=A_{\chi,\rho}A_{\rho,1}$ by Lemma \ref{vanih1ideal} (that is by $h_{1,\chi}=0$ and (\ref{linkextaij})), and we conclude as $A_{\rho,\chi}A_{\chi,\rho} \subset \m$.
\end{pf}

\begin{lemma}\label{cestquasifini}${}^{}$ \begin{itemize}
\item[(i)] There are $f_1,\cdots,f_n \in A_{1,\rho}$ such that $A_{1,\rho}=\sum_{i=1}^n Af_i + A_{1,\chi}A_{\chi,\rho}$.
\item[(ii)] There is a $g \in A_{1,\chi}$ such that $A_{1,\chi}=Ag +   A_{1,\rho}A_{\rho,\chi}$.
\item[(iii)] $A_{1,\rho}A_{\rho,1}=\sum_{i=1}^n f_iA_{\rho,1} + gA_{\chi,\rho}A_{\rho,1}$.
\item[(iv)] For some $\lambda \in K^*$, we have $A_{\chi,\rho}=\lambda A_{\rho,1}$.
\item[(v)] $\m=A_{1,\rho}A_{\rho,1}$.
\end{itemize}
\end{lemma}

\begin{pf} Assertions (i) and (ii) follow from Lemma \ref{valeurshij} (ii) and (iii), formula (\ref{linkextaij}) and Nakayama's lemma. By expanding $A_{1,\rho}A_{\rho,1}$ with the formulas of (i) and (ii), we get part (iii) as the missing term satisfies
$$A_{1,\rho}A_{\rho,\chi}A_{\chi,\rho}A_{\rho,1} \subset \m A_{1,\rho}A_{\rho,1},$$
hence may be deleted by Nakayama's lemma. \par 
Assertion (iv) holds as $A_{i,j}$ and $A_{\tau(j),\tau(i)}$ are $A$-isomorphic submodules of $K$ by Lemma \ref{lestauij} (ii). Part (v) is Prop. \ref{itotestmax} combined with Lemma \ref{computeIP} (ii). 
\end{pf}

\begin{pf}(of Theorem \ref{mainmainthm}) By computing the minimal number of generators of $\m$ with formulae (v) and (iii) of Lemma \ref{cestquasifini}, as well as Lemma \ref{linkh}, we get $$ t \, \, \leq \, \, nh+s,$$ 
where $s$ is the minimal number of generators of $A_{\chi,\rho}A_{\rho,1}$. But this $A$-module is isomorphic to $A_{\rho,1}A_{\rho,1}\subset K$ by (iv) of {\it loc. cit}, so $s\leq \frac{h(h+1)}{2}$ 
and we are done.
\end{pf}

Let us give a simple corollary of this analysis when $H^1_f(E,\rho)$ has dimension $1$ (hence $h=1$), which is somehow the generic situation.

\begin{cor}\label{cormainhegal1} Assume that $h=1$. Then $A$ is regular of dimension $n+1$, all the inequalities of Lemma~\ref{valeurshij} are
equalities, and up to a block-diagonal change of coordinates in $K^{n+2}$, we have 

$$R=\left( \begin{array}{ccc} A & A^n  & A \\
\m^n & M_n(A) & A^n \\
Ag+\m^2 & \m^n & A \end{array} \right) \subset M_{n+2}(K),$$
for some $g \in \m \backslash \m^2$.
\end{cor}

\begin{pf} The proof of Lemma \ref{cestquasifini} actually shows that 
$$t\,\leq\, h_{\rho,1}\, h \,+ \,h_{\chi,1} \,s.$$
If $h=1$, the term on the right is less that $n h+ h (h+1)/2=n+1$. As $t\geq n+1$, all these inequalities are equalities, thus
$h_{\chi,1}=1$ and $h_{\rho,1}=h_{\chi,\rho}=n$. Moreover, Lemma \ref{linkh} shows that $A_{\rho,1}$ is free of rank $1$ over $A$, 
as well as $A_{\chi,\rho}$ by an argument similar to Lemma \ref{cestquasifini} (iv). In particular, up to a block-diagonal change of coordinates we may then assume that 
$A_{\rho,1}=A_{\chi,\rho}=A$, and the corollary follows at once from Lemma \ref{cestquasifini}.
\end{pf}

\subsection{Remarks, questions, and complements}

\subsubsection{The case of a CM field $E$ and the sign of
Galois representations}\label{caseCMfield}

Throughout this paper, we have made the assumption that $E$ is a
quadratic imaginary field. Actually, most of the work we have done can be
extended to the case of a CM field $E$ (say quadratic over its
totally real subfield $E^+$, with $E^+$ of degree $d$ over $\Q$),
but the method (both for the sign conjecture and for this chapter)
ultimately fails if $E$ is not quadratic over $\Q$. Let us explain why.

We would work with a unitary group $\U(m)$ defined over $E^+$, which is compact at every archimedian places
 and quasi-split at every finite places. Such a group exists if $m$ is odd or
if $dm \not \equiv 2 \pmod{4}$. Starting with a couple $(\pi,\rho)$ of an automorphic cuspidal representation $\pi$ of $\Gl_n(\AAA_E)$
a Galois representaions $\rho$ of $\G_E$ satisfying the obvious analogs of
assumptions of~\ref{statemtthmsign}, there shoud exist a representation
$\pi^n$ of $\U(m)$ under the hypothesis that $\varepsilon(\rho,0)=(-1)^{d}$.
The results of chapter 7 would extend easily to this case. But in chapter $8$ and $9$, it is used in a crucial way that
$\Ext_T(1,\chi)=0$. This result is deduced form the fact that
$H^1_f(E,\chi)=0$, which in turns was deduced in chapter 5
from the equality
$H^1_f(E,\chi)=\anneau_E^\ast \otimes_\Z \Q_p$ (Example \ref{exampleBK}) and
the fact that $E$ is quadratic. In the general CM case, we see that instead
$\dim H^1_f(E,\chi)=d-1$, from which we cannot conclude that
$\Ext^1_T(1,\chi)=0$. Thus we are not able to make the proof of the
sign conjecture (chapter 8) work in this case and construct a non zero element of $H^1_f(E,\rho)$.
This is consistent with the fact that we haven't made any hypothesis
implying $L(\rho,0)=0$ (since $\varepsilon(\rho,0)$ may be $1$).

The discussion above is more or less the content of remark 9.1 of \cite{BC}. In the
context of this book, we can offer a much finer analysis of the situation
in the case of a CM field $E$. We denote by $c$ the non-trival element
in $\Gal(E/E^+)$ and by $\sigma$ a lifting of $c$ in $G_{E^+}$ satisfying $\sigma^2=1$. The notations $U^\sigma$ and $U^\bot$ (for
a representation $U$ of $\G_E$) are then defined as in~\S\ref{quad}.

The operation $U \mapsto U^\bot(1)$ defines
a linear involution $\tau$ on $H^1_f(E,\Q_p(1))$.
From Prop~\ref{signebotijfinal}, we see that there is a sign $\epsilon=\pm 1$,
such that the subspace $\Ext^1_T(1,\chi)$ of $H^1_f(E,\Q_p(1))$
is in the eigenspace of eigenvalue $\epsilon$ of $\tau$.
To be more precise about $\epsilon$ we need actually the following result of
independent interest :
\begin{theorem}\label{signedesrepht} Assume only the hypothesis (P0) of Rep$(m)$ (extended
to the case of a CM field $E$). Let $\pi$ be an
automorphic representation of $\U(m)(\A_{E^+})$ as in Rep$(m)$
such that the attached Galois representation $\rho_\pi$ is absolutely irreducible.
Then if $Q \in \Gl_m(\bar \Q_p)$ is such that
$\rho_\pi^\bot(g)=Q \rho_\pi(g)Q^{-1}\chi(g)^{m-1}$ for all $g \in \G_E$,  then we have ${}^t Q =Q$.
\end{theorem}
The existence of a $Q$ as in the statement follows from remark (i) after
Rep$(m)$, and it follows from the absolute irreducibility of $\rho_\pi$ that ${}^t Q = \epsilon Q$ with $\epsilon = \pm 1$ 
(see Lemma \ref{pi2}). If $m$ is odd, then it is clear that $\epsilon=1$.
We postpone the proof that this result also holds for an even
$m$ to a subsequent work. \ps

Going back to our specific situation, we can deduce
\begin{cor} $\Ext^1_T(1,\chi)$ is
a subspace of the $+1$-eigenspace of $\tau$ in
$H^1_f(E,\Q_p(1))$.
\end{cor}
\begin{pf}
By Lemma \ref{corcomputeIP}, the generic representation $\rho_K$ is absolutely irreducible, hence 
we are in the situation of Example \ref{speccasesect1}. We have in particular a collection of signs $\epsilon_s$ 
indexed by the irreducible components of $K$. By an argument already given in \S\ref{analytic}, and the accumulation of 
classical points at $z \in X$, Theorem \ref{signedesrepht} shows that each of those signs is $+1$. The corollary follows then
from Prop~\ref{signebotijfinal} (i).
\end{pf}

But it turns out, perhaps surprisingly,
that the information given by the above corollary is empty :

\begin{lemma} The involution $\tau$ is the identity of $H^1_f(E,\Q_p(1))$.
\end{lemma}
\begin{pf} We recall the Kummer isomorphism
 $$\kum: E^\ast/(E^\ast)^{p^n} \rightarrow H^1(E,\Z/p^n\Z(1)),$$
that sends $x \in E^\ast$ to the class of the cocycle $\kum(x)$ of $G_E$
defined by
$$\kum(x)(s)=s(u)/u,$$
where $u \in \bar E^\ast$ is an element such that $u^{p^n}=x$. The conjugation by $\sigma$ defines an involution of $H^1(E,\Z/p^n\Z(1))$
by sending an extension $U$ to $U^\sigma$. In terms of cocycle,
this involution sends a cocycle $j$ to $j^\sigma$, with
$j^\sigma(s) = j(\sigma s \sigma^{-1})$.                                        
Hence
$$\kum(x)^\sigma(s) = (\sigma s \sigma^{-1})(u)/u = \sigma
(s(\sigma(u))/\sigma(u)) =\sigma(u)/s(\sigma(u)) = \kum(c(x))^{-1}$$
(use that $\sigma^2=1$ and that $\sigma$ acts as the reciprocal on roots of unity). \par
Another natural involution on $H^1(E,\Z/p^n\Z(1))$ is $U \mapsto U^\ast(1)$, and it is
easy to see that this involution sends $\kum(x)$ on $\kum(x^{-1})$. Finally, the involution $\tau$ on $H^1(E,\Z/p^n\Z(1))$
defined by $U \mapsto U^\bot(1)$ is the composition of the two preceding involutions, and thus sends
$\kum(x)$ on $\kum(c(x))$. Taking the limit over $n$, tensorizing by $\Q_p$,
and restricting to $H^1_f$, we see that under the Kummer isomorphism
$$\kum : \anneau_E^{\ast} \otimes_\Z \Q_p \rightarrow H^1_f(E,\Q_p(1)),$$
the involution $\tau$ corresponds to the conjugation $c$. Hence the lemma is reduced to the assertion that $c$ acts by
the identity on a finite index subgroup of $\anneau_E^\ast$. But this
arithmetical statement is a well known consequence of Dirichlet's unit
theorem, that says that $\anneau_E^\ast$ and $\anneau_{E^+}^\ast$ have the same rank.
\end{pf}

\subsubsection{When is $T_A$ the trace of a representation over $A$?}

If $h=1$, then $T_A$ is the trace of a representation of over $A$ by Cor. \ref{cormainhegal1}. Another way to argue would be to say that $A$ 
is regular, hence a UFD, and use Prop~\ref{factfree}. Conversely :
\begin{lemma} If $T_A$ is the trace of a representation of $\G_{E,S}$ over $A$, then we have either $h=1$ or
$h = t \geq n+1$.
\end{lemma}
\begin{pf}
Since $T_A$ is the trace of a representation we may assume that $R \subset M_m(A)$, {\it i.e.} that the $A$-modules $A_{i,j}$ are actually ideals of $A$ 
(use Prop.\ref{kert}, Lemma \ref{traceadapted} and Prop. \ref{imageadapted}). From Lemma~\ref{cestquasifini}(v) we see that either 
$A_{\rho,1}=A$ and $A_{1,\rho}=\m$ or $A_{1,\rho}=A$ and $A_{\rho,1}=\m$, and we conclude by Lemma \ref{linkh}.
\end{pf}

In the conclusion of the above lemma, the case $h=t$ seems very unlikely. However, it is not possible to exclude it by a simple GMA analysis, 
since the data $A_{\chi,1}=\m^2,\ A_{\rho,1}=A_{\chi,\rho}=\m,\ A_{1,\chi}=A_{1,\rho}=A_{\rho,\chi}=A$
define a GMA satisfying all the assertions of Lemma~\ref{valeurshij} (which is even equipped with an obvious anti-involution). \ps
	Another related intruiguing question is to know whether $\Ext_T(\chi,1) \neq 0$. By Cor. \ref{cormainhegal1}, this is the case if $h=1$, 
and the example above shows that it is not formal from what we have proved.

\subsubsection{Other remarks and questions} From a philosophical point of
view, a very intruiguing open question is the following one. \ps

{\bf Question:} Should we expect that $\Ext_T(1,\rho)=H^1_f(E,\rho)$ ?

\ps

On the one hand, although $\Ext_T(1,\rho)$ is a canonical subspace of
$H^1_f(E,\rho)$, it is attached to the unitary group
$\U(n+2)$ so its arithmetic content is somehow included in the one of the
cohomology of the unitary Shimura varieties of dimension $n+2$.

On the other hand, the trend of ideas initiated by Mazur-Wiles'
proof of Iwasawa's main conjecture and
by Wiles' $R=T$ philosophy rather
suggests that we may have equality in our context too. This is
also confirmed by our results in \S\ref{eigenatregcripts}. \ps
       Note that by Corollary \ref{corlissitehigher}, we can detect
directly on the geometry of the
eigenvariety $X$ at $x$ if $\Ext_T(1,\rho)$ has rank $\geq 2$. It would be
very interesting to find
examples where it is indeed the case! As we saw, the space $X$ is built
from some rather explicit spaces of $p$-adic automorphic forms on the
definite unitary group $\U(m)$, thus we hope
that some numerical experiments could be made.\footnote{Actually, for
those questions as well as many others,
it would be very useful to have at hand a database like the one of William
Stein for classical modular forms, as well as a program computing slopes
as Buzzard's one.} The first step is actually to find a $\rho$ for which
the Bloch-Kato conjecture predicts that $\dim_L H^1_f(E,\rho)>1$. When
$n\leq 2$, this amounts to find some modular form of even weight $k\geq
4$, whose sign is $-1$, and whose archimedian $L$-function vanishes at
order $\geq 2$ at $k/2$. The authors do not know any such example at the
moment\footnote{Of course, if we could
have handled the case $k=2$ it would have sufficed to take an elliptic
curve over $\Q$ with
sign $-1$ and rank $\geq 2$, and there are plenty of them.}.\ps
        As explained in Remark \ref{mainconjvsbk}, we hope that we can go
further in the future and make the $L$-function of $\rho$ (or say a
$p$-adic version) enter into the picture, altough it is not clear how at
the moment.

\newpage

\section*{Appendix : Arthur's conjectures}

\renewcommand{\thesubsubsection}{A}
\label{Arthur}

\renewcommand{\thetheorem}{A.\arabic{subsection}.\arabic{theorem}}
\renewcommand{\thesubsection}{A.\arabic{subsection}}

\setcounter{theorem}{0}
\setcounter{subsection}{0}
\setcounter{subsubsection}{0}

 In this appendix, we offer a brief and somewhat personal
exposition of parts of Langlands' and Arthur's conjectural program.
This exposition will allow us to check that the assumptions 
Rep$(m)$ and AC$(\pi)$ 
about automorphic  forms 
on unitary group that we have made
in chapter 6 are predicted by that program. We do 
this for two reasons: first, this should make our 
assumptions more believable, and second, more importantly,
putting those assumptions in the general picture of 
Langlands' and Arthur's program is very helpful to understand our method and 
how it may or may not be generalized.  
For a more complete overview of the conjectures, we send the reader to 
\cite{Arthur}, \cite{br2} and \cite{Rog1} 
(assuming the reader already knows well the material of \cite{borel}).

Let $F$ be a number field and $G$ a 
connected reductive group\footnote{The conjectures we shall describe so
far has been mainly tested for classical groups -- it is possible that some 
minor changes will be needed in the exceptional cases.} over $F$.
An {\it automorphic representation} $\pi$ of $G$ 
is an irreducible constituent of the regular representation of 
$G(\AAA_F)$ on the Hilbert space\footnote{Precisely, let $Z$ be the center of $G$ and 
$Z_{\infty}$ be the connected component of its real points. The aforementionned space is the space of measurable 
complex functions $f$ on $G(F)\backslash G(\AAA_F)$ such that the associated map $g\mapsto f(g)\omega^{-1}(g)$ is 
$Z_{\infty}$-invariant and square-integrable on $G(F)Z_{\infty}\backslash G(\AAA_F)$ for a (finite) $G(\AAA_F)$-invariant measure on this latter space (which exists by a result of Borel and Harish-Chandra).} $$L^2(G(F)\backslash G(\AAA_F),\omega,\C),$$ for some admissible character 
$\omega$ of $G(\AAA_F)$ which is trivial on $G(F)$.
In general, this representation will have a discrete and a continuous 
part, which makes the previous definition rather unprecise. Recall that $\pi$ is {\it discrete} if it occurs discretely in (that is, as a sub-representation of) the space above, such $\pi$ are well defined. 
For example, it is known that the $L^2$ above are completely discrete if (and only if) $G$ is anisotropic modulo its center 
(this was the case of our definite unitary groups of \S\ref{unitarygroups}). In general, Langlands theory of 
Eisenstein series reduces the study of all the automorphic $\pi$ to the discrete ones, 
and we will focus on those ones in this appendix. For $\pi$ a unitary admissible irreducible representations of $G(\AAA_F)$
which is not a discrete automorphic representation, we set $m(\pi)=0$.
 
The aim of Arthur's program (an extension of Langland's program) is to compute $m(\pi)$ for all $\pi$. 
This is done (partially) by a very rich set of conjectures
that are not completely rigid (meaning that they do not always make precise sense, and that they are susceptible to change
slightly as our understanding progress) that proceed in two steps: 
describing a natural 
partition of the set of all discrete automorphic representations 
into ``packets'' and understanding $m(\pi)$ for $\pi$ within a given packet.   

We denote by $\Pi_\unit(G,F)$ (or $\Pi_\unit$ when there is no ambiguity) 
the set of all isomorphism classes of irreducible unitary 
representations of $G(\AAA_F)$ and by $\Pi_\disc(G,F)$ 
the subset of discrete automorphic representations. 
Two interesting subsets of $\Pi_\disc(G,F)$ are $\Pi_\cusp(G,F)$, the set of 
cuspidal automorphic representations, and $\Pi_\temp(G,F)$, the set of 
tempered
%\footnote{Here and everywhere in this appendix, tempered means ``essentially tempered''.} 
{\it discrete} automorphic representations.

\begin{example}\label{caseglm} ($G=\GL_m$) The Ramanujan conjecture
asserts that $$\Pi_\cusp(\Gl_m,F)=\Pi_\temp(\Gl_m,F)$$ 
(this is known to be false for other reductive groups).\footnote{Actually there are trivial counterexamples, like the one dimensional representations 
of $G=D^*$ for a division algebra $D$, but there are deeper ones with $G$ split, or like the representation $\pi^n$ we are especially interested in.}
Moreover, a theorem of Moeglin-Waldspurger \cite{MoWa} shows that the full discrete spectrum of $G$ is built from the 
cuspidal spectrum of the $\GL_d$ with $d$ {\it dividing} $m$, which might be suprising. For example, if $m$ is a prime, 
a discrete $\pi$ is either cuspidal or one dimensional. In this context, each packet of discrete representations is actually a singleton, 
and each discrete representation occurs with multiplicity one ({\it weak multiplicity one theorem}). 
All these facts are actually predicted by Arthur's philosophy, which not only predicts the $m(\pi)$ but also 
gives a general hint about how the discrete spectrum of a general $G$ is constructed from the tempered one, and even from the cuspidal one of the $\GL_m$.
\end{example}

\subsection{Failure of strong multiplicity one and global $A$-packets} \label{discretespectrum}  

\label{failure}

For two unitary irreducible representations $\pi$ and $\pi'$ of $G(\AAA_F)$, 
say $\pi \sim \pi'$ if $\pi_v \simeq \pi'_v$ 
or almost all primes $v$. 

When $G=\Gl_m$ (resp. an inner form of $\Gl_m$), 
it is known (resp. conjectured) 
that if $\pi$ and $\pi'$ are discrete automorphic then $\pi \sim \pi'$ implies that $\pi=\pi'$
 as a subrepresentation of $L^2(G(F)\backslash  G(\AAA_F),\omega,\C)$ ({\it strong multiplicity one})
but this statement is known know to be false for some other groups, including some groups very close to $\GL_m$ like 
$\SL_m$, and our unitary groups.

Following Arthur we should be able to define naturally certain disjoint subsets of $\Pi_\unit$ called 
{\it global $A$-packets}. Every global $A$-packet should be a subset of an equivalence class of $\sim$, and its intersection with $\Pi_\disc$ should be a full equivalence class of the restriction of $\sim$ to $\Pi_\disc$. 
In other words, $A$-packets are equivalence classes 
of $\sim$ in $\Pi_\disc$ enlarged with some
non discrete elements of $\Pi_\unit$, chosen in the same equivalence class of $\sim$. Any nonempty $A$-packet should contain a discrete representation.
The motive for such an enlargement is to allow the global $A$-packets to be products of {\it local 
$A$-packets}, as we soon shall see. For $\Gl_m$ (resp. an inner form), thus, every $A$-packet is (resp. should be)
a singleton. A global $A$-packet $\Pi$ is said tempered, if every $\pi \in \Pi$ (discrete or not) is tempered. 

To describe the set of $A$-packets we need to introduce the conjectural {\it Langlands group} $L_F$.

\subsection{The Langlands groups}
\newcommand{\Rep}{{\rm Rep}}
For $K$ a topological group, we define $\Rep_m(K)$ (resp. $\Irr_m(K) \subset \Rep_m(K)$) as the set of equivalence 
classes of complex $m$-dimensional continuous (resp. moreover irreducible) representations of 
$K$ whose range contains only semi-simple elements.
 
According to Langlands, there should exist for tannakian reasons a group $L_F$ 
(called by others the Langlands group), extension of the global Weil-group $W_F$ by a compact group,
with a natural bijection (the {\it global correspondence}) between $\Irr_m(L_F)$ 
and $\Pi_\cusp(\Gl_m,F)$. The collection of $L$-groups $\{L_F\}$ with $F$ varying should satisfies conditions similar 
to the collection of global Weil-groups $\{W_F\}$ (see \cite{tate}).

For $v$ a place of $F$, we define explicitly a local Langlands group 
$L_{F_v}$ as the Weil group $W_{F_v}$ if $v$ is archimedean, and 
$W_{F_v} \times \SU_2(\R)$ if $v$ is not. In this latter case 
the local Langlands group is closely related to the Weil-Deligne group, 
in the sense that there is a simple bijection between $\Rep_m(L_{F_v}$) and the set of Frobenius semi-simple Weil-Deligne 
representations $(r,N)$ of $F_v$ that we recalled in \S \ref{localLanglands}.
Langlands conjectured the existence of a natural bijection (the {\it local correspondence}) between $\Rep_m(L_{F_v})$ and the set of 
equivalence classes of irreducible admissible representations of $\Gl_m(F_v)$. He proved that conjecture when $v$ is archimedian, 
and the non-archimedian case is now also a theorem of Harris and Taylor \cite{HT}, 
relying on works on many people.

It is part of the conjectures that there exists a distinguished class of embedings 
$L_{F_v} \hookrightarrow L_F$, and the global 
correspondence should coincide after restriction
to the local Langlands group to the local correspondence.

\subsection{Parametrization of global $A$-packets}\label{globApack}

We now go back to a general reductive group $G$. We refer to \cite{borel} for the definition of the $L$-group $\LG$ of $G$. Let us simply recall that $\LG$ is a semi-direct product of $W_F$ by the
dual group $\hat{G}(\C)$ of $G(\C)$, that the product is direct if $G$ is 
split, and that the $L$-group of two inner forms are canonically isomorphic.

Following Arthur, a {\it global $A$-parameter} (for $G$) is 
a continuous homomorphism $$\psi : L_F \times \Sl_2(\C) \longrightarrow \LG$$ such that
\begin{itemize} 
\item[(o)] $\psi_{|\Sl_2(\C)}$ is holomorphic and falls into $\hat{G}(\C)$,
\item[(i)] for all $w \in L_F$, the image of $\psi(w)$ in the quotient $W_F$ of $\LG$ is the same as the image of $w$ by the map $L_F \rightarrow W_F$,
\item[(ii)] $\psi(w)$ is semi-simple\footnote{Recall that an element $g$ of ${}^L G$ is semisimple 
if its image in each representation ${}^L G \rightarrow {}^L \GL_m$ has a semisimple $\GL_m(\C)$ component.} for all $w \in L_F$,
\item[(iii)] the image $\psi(L_F)$  is bounded in ${}^L G$ modulo the center 
$Z(\hat{G}(\C))$ of $\hat{G}(\C)$,
\item[(iv)] $\psi$ is {\it relevant}, that is $\psi(L_F \times \Sl_2(\C))$ is not allowed to lie in a parabolic subgroup\footnote{Recall that it is a subgroup $P \subset \LG$ which surjects onto $W_F$ and which is 
the normalizer in $\LG$ of a parabolic subgroup of $\hat{G}$, see \cite[\S 3.3]{borel}.}
of $\LG$ unless the corresponding parabolic subgroup of $G$ is defined over $F$. 
\end{itemize}
Note that condition (iv) is automatic if $G$ is quasi-split since every parabolic subgroup of a quasi-split group is 
defined over the base field.    

\begin{remark} \label{reducedLgroup} 
In the definition of the $L$-group, note that the Weil group $W_F$ acts on $\hat{G}(\C)$ through a finite 
quotient $\Gal(E/F)$, where $E$ is a finite Galois extension of $F$ on which $G$ splits. For the sake of defining global $A$-parameters 
(the same remark will apply for local $A$-parameters and $L$-parameters, 
see below), it would not change anything if the $L$-group ${}^L G$ was 
replaced by the {\it reduced $L$-group} of $G$, 
namely the semi-direct product of $\Gal(E/F)$ by $\hat{G}(\C)$ 
(being understood that in condition (i) above, each occurrence of $W_F$ is to be replaced by $\Gal(E/F)$).
In particular, an $A$-parameter for a {\it split} group $G$ is 
simply a morphism $L_F \times \Sl_2(\C) \rightarrow \hat{G}(\C)$
satisfying condition (ii) and (iii) above - remember that (iv) is automatic.
\end{remark} 

Two parameters are said {\it equivalent} if they are conjugate by an element 
of $\hat{G}(\C)$, up to a $1$-cocycle of $L_F$ in $Z(\hat{G}(\C))$ which is locally  trivial at every place (see \cite{Rog1} \S 2.1). 
In the cases we will deal with, namely $\Gl_m$ and unitary groups, 
it turns out  that every such cocycle is trivial (\cite{Rog1} \S 2.2). 
Thus in those cases the equivalence relation for parameters is just the 
conjugacy by an element of $\hat{G}(\C)$. 

\begin{definition}\label{defdisctemp} \begin{itemize}
\item[-] A global $A$-parameter $\psi$ 
is said to be discrete if $C(\psi)^0 \subset Z(\hat{G}(\C))$, where $C(\psi)=\{g \in \hat{G}(\C), g\psi(w)=\psi(w)g \, \, \, 
\forall w \in L_F \times \Sl_2(\C)\}$ is the centralizer in 
$\hat{G}(\C)$ of the image of $\psi$, 
and $C(\psi)^0$ is its neutral component.
\item[-] A global $A$-parameter $\psi$ is said tempered\footnote{There is a common 
abuse of language here, as strictly we should say {\it essentially tempered}, that is tempered up to a twist.} 
if its restriction to $\Sl_2(\C)$ is trivial.
\end{itemize}
\end{definition}

\begin{example} \label{caseglm2} Assume again $G=\Gl_m$. 
We see at once that a global $A$-parameter $\psi$ is discrete if, and only if, 
the corresponding representation $L_F \times \SL_2(\C) \rightarrow \Gl_m(\C)$ is irreducible. In particular, there exists a divisor 
$d$ of $m$ and an irreducible tempered parameter $\psi': L_F \rightarrow \GL_d(\C)$ such that $\psi'=\psi \otimes [d]$, where $[d]$ denotes the unique 
$d$-dimensional holomorphic representation of $\Sl_2(\C)$. Moreover, using the fact that $L_F$ (as $W_F$) should be an extension of an abelian group (namely $\R$) 
by a compact group, we see easily that $\Irr_m(L_F)$ should be in bijection with the set of global discrete tempered $A$-parameters of $\GL_m$. Note that this formalism 
matches perfectly with Ex. \ref{caseglm}.
\end{example}
\ps 
{\bf The first main conjecture of Arthur is the existence of a natural correspondence 
which associates to every global discrete $A$-parameter of $G$ (up to 
equivalence) an $A$-packet of $G$, or the empty set.}
\ps

Two $A$-parameters $\psi$ and $\psi'$ should be sent to the same 
(non-empty) $A$-packet if and only if 
their restriction $\psi_v$ and $\psi'_v$ to $L_{F_v}$ for all $F_v$ 
are conjugate. This condition may not imply that $\psi$ and $\psi'$ are 
equivalent in general, but it will in our situations (see Prop~\ref{proppar} below). 
The correspondence above should send tempered $A$-parameters to (non-empty) tempered $A$-packets, 
in which case it should coincide with the former theory of Langlands. Note that this requirement, 
in the case $G=\Gl_m$, is the generalized Ramanujan conjecture. 

Note that some global $A$-parameter,  even satisfying the relevance condition, 
may very well be sent to the empty set. We shall give an example below 
(see Remark~\ref{trvialApacket}). However, this should not happen when $G$ is quasi-split, 
or for a tempered $A$-parameter.

To understand the ``naturality'' of the correspondence between discrete global $A$-parameters and global $A$-packets, we need to 
introduce the local counterpart of global $A$-packets and global $A$-parameters. It should be stressed here that contrary 
to Langlands theory of local $L$-packets which should apply to all the admissible irreducible representations, 
the introduction of Arthurs' $A$-packets is mainly motivated by global considerations and basically apply to local components 
of global automorphic representations.

\subsection{Local $A$-packets and local $A$-parameters}

Following Arthur, a {\it local} {\it $A$-para\-meter} is a continuous morphism  
$$\psi_v : L_{F_v} \times \Sl_2(\C) \longrightarrow \LG$$ such that the analogues of conditions (o) to (iii) of 
global parameters are satisfied (of course, $F$ has to be replaced by $F_v$ everywhere there, and $\LG$ is now the $L$-group of $G/F_v$).
Note that there is no relevance condition (iv) in the definition.  
As in the global case, a local $A$-parameter is said to be {\it tempered} if it is trivial on 
$\Sl_2(\C)$. The restriction $\psi_v$ of a global $A$-parameter $\psi$ 
to $L_{F_v} \times \Sl_2 (\C)$ is obviously a local $A$-parameter, and $\psi_v$ is tempered if $\psi$ is. 

Following Arthur, to every local $A$-parameter 
should correspond a finite set, possibly empty, of irreducible unitary representations of $G(F_v)$. 
This maps will not be injective in general. For example, 
many $A$-parameters could define the empty $A$-packet. 
Moreover, in contrast with global $A$-packets, 
local $A$-packets are not, in general, disjoint, 
and their reunion will not be the set of all unitary irreducible 
representations of $G(F_v)$ (but rather the set of such representations that 
appear as constituents of global automorphic representations.)

 According to Arthur, a global $A$-packet $\Pi$
defined by an $A$-parameter $\psi$ should be the set 
of restricted tensor products $$\Pi=\{ \pi=\otimes_v' \pi_v, \, \, \pi_v \in \Pi_v\},$$
for a set of local representations $\pi_v$ belonging for each $v$ to the local $A$-packet $\Pi_v$
corresponding to $\psi_v$, and such that almost all $\pi_v$ are unramified.

A tempered $A$-parameter should define a tempered $A$-packet, 
that is an $A$-packet all of whose members are tempered.

Ultimately, local $A$-packets should be constructed 
using our understanding of the trace formula and its  stabilization. One
key-property that an $A$-packet should satisfy is that a suitable 
non-trivial linear combination of the 
character-distributions of its members should be a stable  
(that is invariant by conjugation by elements of $G(\bar F_{v})$, not only of
 $G(F_v)$) distribution.

To say more on the correspondence between $A$-parameters and $A$-packets, 
we need to review the earlier notions of $L$-parameters and $L$-packets, due to Langlands.

\subsection{Local $L$-parameters and local $L$-packets}\label{landalocal}

Followings Langlands, there should be a {\it partition} 
of the set of equivalences classes of all {\it admissible}
irreducible representations of $G(F_v)$ into (local) $L$-packets. 
We stress that the local $L$-packets behave much more nicely that the
local $A$-packets, since the former are {\it disjoint} and that their 
reunion do not miss any admissible irreducible representation.

The set of $L$-packets should be in bijection with the set of conjugacy 
classes of relevant {\it $L$-parameters}. Recall that an $L$-parameter is a continuous morphism
$$\phi_v:\  L_{F_v} \longrightarrow {}^L G$$ 
that satisfies conditions (i) and (ii) (but not (iii)) of \S\ref{globApack}, it is said to be relevant if it satisfies moreover (iv) {\it loc. cit.} 
Moreover, an $L$-parameter is said to be {\it discrete} if it connected centralizer is central, as in Def. \ref{defdisctemp}.

A local $A$-parameter $\psi_v$ defines a local $L$-parameter 
(that may not be relevant) by the formula 
\begin{eqnarray}
\label{phipsi} 
\phi_{\psi_v}(w)=\psi_v\left(w,\left(\begin{matrix}|w|^{1/2} & 0 \\ 0 & |w|^{-1/2}\end{matrix}\right)\right).
\end{eqnarray}
Here $|.|: L_{F_v} \rightarrow \R^*$ is the composition of $L_{F_v} \rightarrow \W_{F_v}^{\rm ab} \isomo_{\rm rec^{-1}} F_v^*$ with the norm of $F_v^*$. The local $A$-packet corresponding to $\psi_v$ should contain the local $L$-packet corresponding to $\phi_{\psi_v}$ (if relevant) 
and could be larger in general, but not when $\psi_v$ is tempered. 

The problem with $L$-packets that motivated the introduction of $A$-packets 
is that it is not always possible to construct a non-trivial linear 
combinations of the characters of its members that is a stable distribution. 
This problem does not arise for tempered $L$-packets.

\subsection{Functoriality} 

If $G$ and $G'$ are two groups as above, any admissible morphism of $L$-groups 
(that is a holomorphic group homomorphism compatible to the projection to  
$W_F$)  $${}^L G \rightarrow {}^L G'$$ induces a 
map from $L$-parameters or (local and global) $A$-parameters for $G$ to 
similar parameters for $G'$. 
According to the conjectures described above, 
this should determine a map (rather, a correspondence) from the set of 
packets (local or global, $A$ or $L$) for $G$ to the set of packets for 
$G'$. Such a conjectural correspondence 
is an instance of Langlands' functorialities.

The most basic example is the case where 
$G$ and $G'$ are inner form of each other, which is sometimes called a {\it Jacquet-Langlands} transfer. Then ${}^L G={}^L G'$ and the identity map 
should define a correspondence between packets of $G$ and of $G'$. 
Note that even in this simplest case, this 
correspondence may not be a map  
(even for local $L$-packets) since a parameter relevant for $G$ 
may not be relevant for $G'$.

Note that in defining functorialities, it is useful to work with the 
full $L$-groups, not only the reduced ones, since there are more morphisms between full $L$-groups. 

\begin{example}\label{caseglm3} Let $D$ be a quaternion division algebra over $F$, $G=D^*$ and $G'=\GL_2$. 
The $A$-parameter of the trivial, discrete, representation $\pi$ of $G$ or $G'$ (global or local) is the discrete, non-tempered, parameter $1 \otimes [2]$ in the notations of Ex. \ref{caseglm2}, 
and the corresponding $A$-packets have a single element $\pi$. Of course, the Jacquet-Langlands 
$A$-functoriality makes those trival representations correspond. \par
	If we had tried to understand this simple transfer in the context of $L$-functoriality, 
we see that we could not have asked that the transfer of $\pi$ is both discrete and compatible at all the finite places with the local Langlands correspondence. Indeed, 
for each finite place $v$ such that $D_v$ is nonsplit, this latter correspondence would match the trivial representation of $D_v^*$ with the 
Steinberg representation of $\GL_2(F_v)$, which is infinite dimensional and that would contradict the strong approximation theorem for $\GL_2$. \par
In other words, as long as we are interested in the discrete spectrum (say), 
and with non tempered representations, the $A$-functoriality is better behaved than the $L$-one, and it is actually made for that. 
This phenomenum is not just a fancy problem related to the trivial representation, but it 
appears in all kind of functorialities. We will give later in Ex. \ref{instructiveexample} a deeper example due to Rogawski 
in the case of a base change.
\end{example}

\subsection{Base Change of parameters and packets}

\label{basechange}

In this paragraph, it will be convenient to assume that $F$ is either a global or a local field. 
Recall the notion of $L$-parameter we gave applies to the local context only, whereas the one of $A$-parameter does 
in both cases. Let $E$ be a finite extension of $F$ and set $$G_E:=G \times_F E.$$ The restriction of an $L$-parameter 
$\phi$ (resp. an $A$-parameter $\psi$) of $G$ to $L_{E}$  (resp. to $L_E \times \Sl_2(\C)$) defines an $L$-parameter $\phi_E$ 
(resp. an $A$-parameter $\psi_E$) of $G_E$. The map $\phi \rightarrow \phi_E$ (resp. $\psi \rightarrow \psi_E$) is called the {\it base change map} for parameters.\footnote{This base change may be viewed as a special case of the general functoriality by considering 
the natural map ${}^L G \rightarrow {}^L(\Res_F G_E)$. } 
In general, $\psi_E$ is tempered if $\psi$ is, but $\psi_E$ may be not be discrete although $\psi$ is. 

We shall mainly be interested in the case where $E$ splits $G$ and $\hat{G_E}=\GL_m$, what we assume now. In any of the 
three possible cases it is possible to attach to $\phi_E$ or $\psi_E$ a single admissible irreducible representation of 
$G(E)$ (in the local case) or of $G(\AAA_{E})$, as follows. If $\phi$ is a local $L$-parameter, $\phi_E$ is a local $L$-parameter for $G_E$ (it is automatically relevant) and thus defines a single admissible irreducible representation of $G(E)=\Gl_m(E)$ by the local Langlands correspondence. If $\psi$ is a local $A$-parameter, then so is $\psi_E$, and the map $\phi_{\psi_E}$ defined in formula 
(\ref{phipsi}) is an $L$-parameter of $G_E$ and thus defines again an admissible irreducible representation of $\Gl_m(E)$. Finally if $F$ is global and $\psi$ a global $A$-parameter, we associate to $\psi$ 
the restricted tensor product on all places $w$ of $E$
of the representation attached to the base change $\psi_w$ of $\psi_v$ (if $v$ is the place of $F$ below $w$).

To summarize, we have defined (assuming Langlands and Arthur's parameterization) a base change to $G_E=\Gl_m$ 
of a global {\it $A$-packet}\footnote{Note that it is clear form the definition that even if an $A$-packet is defined by two different parameters, the base change representation attached to those parameters is the same.} $\Pi$ of $G$ which is a single irreducible admissible representation of $\Gl_m(\A_E)$, and also for local $A$-packets and local $L$-packets for $G$ (and the result is then an admissible irreducible local of $\Gl_m(E)$). Note that by definition, the 
base change for global and local $A$-packets are compatible in the 
obvious sense.

\subsection{Base change of a discrete automorphic representation} 

\label{Basechange2}

We keep our assertion that $G_E = \Gl_m$.

If $F$ is a global field and $\pi$ is a global discrete automorphic representation of $G$, it belongs to a unique global $A$-packet $\Pi$ which has a well defined base change as we saw above. 
We define the {\it base change} of $\pi$ as the base change $\pi_E$ of its $A$-packet. If $\psi$ is an $A$-parameter corresponding to $\Pi$ and if 
$\psi_E$ is discrete, then $\pi_E$ should be a discrete automorphic representation of $\GL_m/E$.

If $F_v$ is a local field and $\pi_v$ is an irreducible representation of $G(F_v)$,
it belongs to a single $L$-packet $\Pi_v$. Hence we may define
its base change as the base change of this local $L$-packet. 
Note that we can not use local $A$-packets to define local base change 
inambiguously since a representation may belong to several local 
$A$-packets that have different base change (see Ex. \ref{instructiveexample} below).

The inconvenient of defining local and global base change for representations 
using different types ($A$ and $L$) of packets is that in general
there is {\it no compatibility} between local and global base change 
for representations.  However, if $\pi$ is a global discrete automorphic representation which is {\it tempered}, it belongs to a 
tempered $A$-packet $\Pi$, product of tempered local 
$A$-packets $\Pi_v$ that are also $L$-packets and contain $\pi_v$. 
Hence it is clear that the $v$-component of the global base change of $\Pi$ is in this case 
the base change of the local components $\pi_v$ of $\pi$.

\subsection{Parameters for unitary groups and Arthur's conjectural description of the discrete spetrum}\label{descarthunit}

In this paragraph we specify to unitary groups the formalism developped above, including 
parameters ($L$ or $A$, local or global) and the base change to $\GL_m$. We fix a unitary group 
$G:=\U(m)$ (quasi-split or not) in $m$ variables attached to a CM extension $E/F$ of numberfields. The reduced $L$-group of $\U(m)$ is the semi-direct product 
$${}^L\U(m)=\Gl_m(\C) \sd \Gal(E/F)$$ 
where $\Gal(E/F)=\Z/2\Z=\langle c \rangle$ acts on $\GL_m(\C)$ by 
\begin{eqnarray} 
\label{lgroup} 
c M c^{-1} := \phi_m {}^tM^{-1} \phi_m^{-1}, \, \, \, M \in \GL_m(\C),
\end{eqnarray} where,
\begin{eqnarray} \label{phim} \phi_m:= 
\left( \begin{array}{cccc}  & & & (-1)^{m+1} 
\\ & & \dots &  \\ & -1 & &  \\ 1 & & &  \end{array} \right).
\end{eqnarray}
Note that
\begin{eqnarray}  
\label{propphim} {}^t \phi_m = \phi_m^{-1} = (-1)^{m+1} \phi_m, \, \, \, c\phi_m=\phi_mc.
\end{eqnarray}

Note that $\U(m)_E$ is an inner form of $\Gl_m$, hence the theory of base change to $\GL_m/E$ explained in \S \ref{basechange}, \ref{Basechange2} applies. 

An $A$-packet $\Pi$ of $G$ will be a tensor products of local $A$-packets $\Pi_v$, 
where $\Pi_v$ will have one element when $v$ splits in $E$, but more than one in general for the other places. In particular, 
$\Pi$ may be infinite in general. We now review what Arthur's theory of parameters implies for the structure of the discrete spectrum of $G$ (and for these packets),
following Rogawski's analysis \cite[\S 2.2]{Rog1}. If 
$$\psi: L_F \times \Sl_2(\C) \rightarrow {}^L \U(m)$$ is a discrete $A$-parameter, 
Rogawski shows that $\psi_E: L_E \times \Sl_2(\C) \rightarrow \GL_m(\C)$ is a 
direct sum of $r$ pairwise nonisomorphic irreducible representations $\rho_j$ (see the proof of Lemma 2.2.1 {\it loc. cit.}) 
that satisfy $\rho_j^{\bot} \simeq \rho_j$.
%Moreover, the group of connected components of $C(\psi)/Z(\hat{G})$ is then isomorphic to $(\Z/2\Z)^{r-1}$.
He defines $\psi$ to be {\it stable} if $\psi_E$ is irreducible. As explained in Ex. \ref{caseglm2}, $\psi$ is stable if and only if
$\psi_E$ is discrete, in which case there is a discrete automorphic representation $\pi_E$ of $\GL_m/E$ which is the base change of $\psi$ as 
explained in \S \ref{Basechange2}. That representation $\pi_E$ is cuspidal if and only if $\psi$ is tempered.

In general, for any (unordered) partition $$m=m_1+\cdots+m_r,$$
Rogawski defines an admissible map\footnote{This map is a special case of the so-called endoscopic functoriality as $\U(m_1)\times \cdots \times \U(m_r)$ is not a Levi subgroup of $\U(m)$
when $r>1$ (see \cite[\S 1.2]{Rog1}). When some $m_i$ has not the same parity as $m$, $\xi$ actually not defined at the level of the reduced $L$-groups, 
as a character $\mu$ as in \S\ref{introcarhecke} occurs in its definition, and $\mu$ is not of finite order. It is an 
important fact that $\xi$ is also not canonical at all in this case, as it depends on this choice of $\mu$. The uniqueness 
assertion in the Rogawski's description also assumes that such an $\mu$ has been fixed once and for all.}
	$$\xi: {}^L ( \U(m_1)\times \cdots \times \U(m_r) ) \rightarrow {}^L  \U(m).$$
He shows then that any discrete $A$-parameter $\psi$ of $G$ writes uniquely as 
$$\xi \circ (\psi_1 \times \cdots \times \psi_r)$$ where the $\psi_j$ are distinct stable parameters 
of the quasisplit group $\U(m_i)$ and for a unique unordered partition $m=m_1+\cdots+m_r$ as above (\cite[Lemma 2.2.2]{Rog1}). 
We say that $\psi$ is endoscopic if it is not stable, {\it i.e.} if $r>1$. 

If $G$ is split, this reduces conjecturally the study of the discrete spectrum of $G$ 
to the stable parameters (compare with Ex. \ref{caseglm2}), and the general case is 
then a matter of relevance.\footnote{The situation is actually not as simple as it may seem in the nontempered case, 
as a relevant parameter may exceptionally lead to an empty packet. 
Moreover, the multiplicity formula in more complicated in the nonsplit case.} This {\it structure} of the discrete 
spectrum of $G$, as well as other predictions of Arthur, 
have been verified by Rogawski when $m\leq 3$ \cite{Rog4}. 

We will refine slighty this study in \S\ref{descentum} by giving some sufficient conditions on 
an $A$-parameter $\psi'$ of $\GL_m/E$ to {\it descend to $\U(m)$}, {\it i.e.} ensuring that $\psi'=\psi_E$ for some 
discret $A$-parameter $\psi$ of $\U(m)$. As an exercise,\footnote{Note that 
$\phi_m$ is antisymetric if and only if $[m]$ is symplectic.} the reader can already check that the parameter $1 \otimes [m]$ 
descend to a stable nontempered parameter of $\U(m)$. Its associated $A$-packet has a single element, namely 
the trivial representation. 

\begin{remark} Let $\pi$ be a discrete automorphic representation of $G$. If $\pi$ is nontempered, 
the presence of a non trivial representation of the $\SL_2(\C)$ in the $A$-parameter of (the $A$-packet containing) $\pi$ 
imposes strong restrictions on the $\pi_v$. For example, if $G(F_v)$ is compact, 
$\pi_v$ cannot be regular in the sense of \S\ref{sectrepm}. In particular, 
the regularity assumption there should imply that the $\pi$ is tempered, hence that its local and global base change are compatible (see \S\ref{Basechange2}).

By contrast, there are very few conditions ensuring that $\pi$ belongs to a stable $A$-packet. A standard sufficient condition (as in the works of Kottwitz, 
Clozel and Harris-Taylor) is that $\pi_v$ is square-integrable at a split place $v$ (this follows easily from the Arthur's formalism), but 
this condition pushes aside a lot of very interesting stable packets.
\end{remark}

\subsection{An instructive example, following Rogawski.}
\label{instructiveexample} We give now a very instructive example of a nontempered $A$-packet for the group $\U(3)$ which 
illustrates most of the subtleties that appeared till now. It was found by Rogawski (\cite{Rog4},\cite{Rog2}), and it
is probably the simplest of such examples. We stress that it should not be thought as exotic, but rather as an important 
intuition for the general situation. Moreover, it is exactly the kind of packets that we use in the arithmetic 
applications of chapters 8 and 9. 

We keep here the notations of \S\ref{descarthunit} and take $m=3$. We are interested in the nontempered $A$-parameter 
associated to the partition $3=2+1$. These parameters have actually a nonconjectural meaning as they factor through 
$W_F$. To fix ideas we fix $\eta$ an automorphic character of $\AAA_E^*$ such that $\eta^{\bot}=\eta$ and that 
$\eta$ does not descend to $\U(1)$ (see \S\ref{introcarhecke}). By class field theory, we may view it as a continuous 
character of $W_E$. We may actually use this $\eta$ to define $\xi$ and we are interested in the
parameter $\xi \circ (1 \times 1 \otimes [2])$. More explicitely, let us simply say\footnote{This suffices for the discussion here. 
We will say more about those parameters 
and their extension to $L_F$ in a more general context in \S\ref{apretpackpin}.} that there is a unique parameter $\psi$ 
whose base change $$\psi_E: W_E \times \Sl_2(\C) \longrightarrow \GL_3(\C)$$ fixes the vector $e_2$ 
of the canonical basis $(e_1,e_2,e_3)$ of $\C^3$, and which acts as $$(w\times g) \mapsto \eta(w)g$$ on $\C e_1 \oplus \C e_3 = \C^2$.

We are going to describe completely the $A$-packet $\Pi$
associated to $\psi$ following Rogawski. As predicted, it is a product of local $A$-packets $\Pi_v$, 
so we are reduced to describe each of the $\Pi_v$. 
When $v$ splits in $E$, $\Pi_v$ is a singleton and coincide with its associated (non-tempered) $L$-packet defined in \S\ref{landalocal}, so we concentrate on the nonsplit case.

Consider first the local quasi-split unitary group $U_3(\Q_l)$ attached 
to quadratic extension $E_v$ of $\Q_l$. Rogawski has defined
a non-tempered representation $\pi^n(\eta_v)$, where 
$\eta_v : E_v^* \longrightarrow \C^\ast$ is the restriction of $\eta$ at $v$. Recall that $\eta_v = \eta_v^\bot$ 
but that $\eta_v$ does not come by base change from $\U(1)(\Q_l)$.
That representation is actually a twist of the one we constructed in \S\ref{construpinnonsplit} 
as the unique subrepresentation $\pi^n_l$ of a principal series of $\U(3)(\Q_l)$. This principal series has in this case two components. 
The other is called $\pi^2(\eta_v)$. According to Rogawski it is square integrable.

The representation $\pi^n(\eta_v)$ forms an $L$-packet on its own. 
This $L$-packet is not tempered, and $\pi^n(\eta_v)$ is not stable.
The $L$-packet containing $\pi^2(\eta_v)$ is $$\{\pi^2(\eta_v),\pi^s(\eta_v)\},$$ where $\pi^s(\eta_v)$ is a 
supercuspidal representation  that Rogawski constructs 
using global considerations involving the trace formula. 
Since this $L$-packet is tempered, it is also an $A$-packet. 
There is one $A$-packet containing $\pi^n(\eta_v)$, namely $\Pi_v$, and we have\footnote{Note that $\pi^n(\eta_v)$ is a NMSRPS representation, but not $\pi^s(\eta_v)$.}
$$\Pi_v=\{\pi^n(\eta_v),\pi^s(\eta_v)\}.$$ In particular,
$\pi^s(\eta_v)$ belong to {\it two} $A$-packets, 
and actually those representations are the only ones (up to a twist) 
that belongs to several $A$-packets.  

The base change of the $A$-packet $\{\pi^n(\eta_v), \pi^s(\eta_v)\}$, 
and of the $L$-packet $\{\pi^n(\eta_v)\}$ is the irreducible admissible representation of $\Gl_3(E_v)$
whose $L$-parameter is $$\eta_v |\,|^{1/2} \oplus \eta_v |\,|^{-1/2} \oplus 1.$$ The base change of the $L$ and $A$-packet 
$\{\pi^n(\eta_v),\pi^s(\eta_v)\}$ coincide with this latter parameter on the Weil-group, but is nontrivial on the $\SU_2(\R)$-factor 
(mixing $\eta_v|.|^{1/2}$ and $\eta_v|.|^{-1/2}$). Hence the two $A$-packets containing $\pi^s$ have different base changes.

Assume till the end of this subsection that $G(F_v)$ is compact for each archimedian $v$, which is our main case of interest in this book.
For $v$ archimedian, $\Pi_v$ is empty if $\eta_v$ has weight $\pm 1/2$ (see \S\ref{introcarhecke}), a singleton else: 
namely the one we defined in \S\ref{construpisinfty} (up to a twist). This ends the description of $\Pi$.

As a consequence of all of that, we first see that $\Pi$ is infinite. Moreover, Rogawski computes then $m(\pi)$ for each $\pi \in \Pi$, hence 
$\Pi \cap \Pi_\disc(G,F)$. He shows that $m(\pi)$ is always $0$ or $1$.\footnote{Hence strong multiplicity one holds for the packet $\Pi$, 
actually Rogawski shows that it holds for the full discrete spectrum of $G$.} Precisely, he assigns a sign $\varepsilon(\pi_v)=\pm 1$ to each 
$\pi_v \in \Pi_v$ as follows: $\varepsilon(\pi_v)=1$ except when $v$ is archimedian, or when $v$ is a finite nonsplit place and 
$\pi_v = \pi^s$. The final result \cite{Rog2} is that $m(\pi)=1$ if, and only if, 
$$\prod_v \varepsilon(\pi_v)=\varepsilon(\eta,1/2),$$ where 
$\varepsilon(\eta,1/2)$ is the sign of the global functional equation of $\eta$. In particular, {\it one half} of $\Pi$ is actually automorphic. 

The formula above is a special case of Arthur's multiplicity formula. We will discuss in more details this multiplicity formula in a more general case in \S\ref{formulemultart}.

\subsection{Descent from $\GL_m$ to $\U(m)$}\label{descentum}

In this paragraph we explain the algebraic formalism
relating the parameters ($L$ or $A$, local or global)
of a unitary group $G=\U(m)$ as in \S\ref{descarthunit}
and their base change to $\Gl_m$. This formalism
also applies to Galois representations instead of parameters.

We will use systematically the notation $\U(m)$ for $G$,
which frees the letter $G$ (and $G'$) for other notational purposes.

We consider a group $G$ and a subgroup $H$ of index $2$.
In this \S, we call {\it parameter} any morphism
$\psi : G \longrightarrow {}^L\U(m)$
such that $H$ is the kernel of the composition of $\psi$ and the
projection ${}^L\U(m) \rightarrow \Z/2\Z$.

We denote by $d$ a fixed element of $G-H$.
If $\rho : H \rightarrow \Gl_n(\C)$ is a morphism,
define $\rho^{\bot}(h)={}^t \rho(d^{-1} h d)^{-1}$. The
representation $\rho^\bot$ depends on $d$ only up to isomorphism.

In the applications, $E/F$ may either be an extension of global or of local
fields, and $G$ and $H$ may be respectively $L_F$ and $L_E$, their Arthur's
variant $L_F \times \Sl_2(\C)$ and
$L_E \times \Sl_2(\C)$, the Weil groups $W_F$ and $W_E$, or the absolute
Galois groups $G_F$ and $G_E$. Note that in general, $G$ is not a semi-direct
product of $\Z/2\Z$ by $H$.

If $\psi$ is a parameter, we may write $$\psi(d)=Ac.$$
where $A \in \Gl_m(\C)$,
since the image of $d$ in $\Z/2\Z$ is non trivial. For the same
reason, we may write $\psi(d^{-1})=Bc$. We thus have
$$1=\psi(d)^{-1} \psi(d)= Bc(Ac) = B \phi_m {}^t A^{-1} \phi_m^{-1}$$
so $B = \phi_m {}^t A \phi_m^{-1}$ and
\begin{eqnarray} \label{psidinv}
\psi(d^{-1}) = \phi_m {}^t A \phi_m^{-1} c.
\end{eqnarray}

From this we  deduce, calling $\rho$ the restriction of $\psi$ to $H$:
\begin{eqnarray}
\label{autod}
\forall h \in H,\ \ \rho(dhd^{-1}) =  C {}^t \rho(h)^{-1} C^{-1},
\end{eqnarray}
where $C=A \phi_m$. Indeed,
\begin{eqnarray*} \psi(dhd^{-1}) &=& A c \psi(h) \phi_m {}^t A \phi_m^{-1} c
\text{\ \ by (\ref{psidinv})}\\
&=& A \phi_m {}^t\psi(h)^{-1} {}^t \phi_m^{-1} A^{-1} {}^t\phi_m
\phi_m^{-1} \text{\ \ by (\ref{lgroup})} \\
&=& A \phi_m {}^t \psi(h)^{-1} (A \phi_m)^{-1}
\text{\ \ using (\ref{propphim})}.
\end{eqnarray*}

Note that in particular, we have $$\rho \simeq \rho^\bot.$$

We also have $\rho(d^2)=\psi(d^2)=\psi(d)^2= A \phi_m {}^t A^{-1} \phi_m^{-1}$ hence
\begin{eqnarray}
\label{d2}
\rho(d^2) = (-1)^{m+1} C {}^t C^{-1}.
\end{eqnarray}

\begin{lemma}\label{lemmepar} A morphism $\rho : H \longrightarrow
\Gl_m(\C)$ is the
restriction to $H$ of a parameter of $G$ if and only if there exists a
matrix $C \in \Gl_m(\C)$ that satisfies conditions
(\ref{autod}) and (\ref{d2}).
\end{lemma}
\begin{pf} We have already seen that those conditions were necessary. To
prove they are sufficient, assume they are satisfied for some
$C \in \Gl_m(\C)$ and define a map $\psi : G \rightarrow{}^L\U(m)$
by setting for all $h \in H$, $\psi(h):=\rho(h)$ and $\psi(hd):=\psi(h) Ac$
where $A:= C \phi_m^{-1}=(-1)^{m+1} C \phi_m$.
We only have to check that $\psi$ is a group homomorphism.
Let $g,g' \in G$. If $g \in H$ then it is clear by definition that
$\psi(gg')=\psi(g) \psi(g')$. So suppose $g=hd$. We distinguish two
cases : if $g'=h' \in H$, then we have
\begin{eqnarray*}
\psi(g) \psi(g')&=& \psi(hd) \psi(h') \\
 &=& \rho(h) Ac \psi(h') \\
&=& \rho(h) A c \rho(h') c^{-1} A^{-1} A c \\
&=& \rho(h) A \phi_m {}^t \rho(h')^{-1} \phi_m^{-1} A^{-1} A c
\text{\ \ using (\ref{lgroup})} \\
&=&  \rho(h) C {}^t \rho(h')^{-1}  C^{-1} A c \\
&=&  \rho(h) \rho(d h'd^{-1}) Ac \text{ using (\ref{autod})} \\
&=& \rho(hdh'd^{-1}) Ac \\
&=& \psi(hdh')=\psi(gg')
\end{eqnarray*}
Similarly, if $g' = h'd$, we have
\begin{eqnarray*}
\psi(g) \psi(g')&=& \psi(hd) \psi(h'd) = \psi(hd) \psi(h') A c\\
&=& \rho(hdh'd^{-1}) (Ac)^2 \text{ like in  the first six lines of the above computation} \\
&=& \rho(hdh'd^{-1}d^2)  = \psi(gg'). \end{eqnarray*}
\end{pf}

\begin{remark} \label{rempar} \begin{itemize}
\item[(i)] The lemma above gives a criterion for a representation
$\rho : H \rightarrow \Gl_m(\C)$
to be the restriction of a parameter of $G$.
Note that the criterion depends on a choice of an element $d$ in $G-H$.
In each particular context, a clever choice of $d$ may simplify the
computations.

\item[(ii)] As an exercise, let us consider the special case where $G$ is a
semi-direct product of $\Z/2\Z$ by $H$. This
case\footnote{A similar study is done in this context in the first pages of \cite{CHT}.}
occurs for example when $E/F$ is a CM extension of numberfields, $G=G_F$ and $H=G_E$.
Then we may and do choose a $d$ such that $d^2=1$.

Let $\rho : H \rightarrow \Gl_m(\C)$ such that
$\rho^\bot \simeq \rho$. Thus there is a $C$ such that (\ref{autod}) is satisfied.
By applying this relation twice, and using $d^2=1$, we see that
$C {}^t C^{-1}$ centralizes $\Im \rho$. Assume that $\rho$ is irreducible.
Then the latter means ${}^t C = \lambda C$, from which $\lambda^2=1$, and we
see that $C$ is either antisymmetric, or symmetric. If $m$ is odd, $C$ has
to be symmetric, but if $m$ is even,
it is clear that both situations may happen.
On the other hand, (\ref{d2}) reads $1=\rho(d^2)=(-1)^{m+1} C {}^t C^{-1}$ which simply
says that $C$ is symmetric if $m$ is odd, and antisymmetric if $m$ is even.

To summarize, an irreducible $\rho$ such that $\rho^\bot \simeq \rho$ comes from a parameter of $G$
always if $m$ is odd, and `half the time'' if $m$ is even.
\end{itemize}
\end{remark}

\begin{prop}\label{propparunicity}
Let $\rho : H \rightarrow \Gl_m(\C)$ be a semisimple representation such that $\rho\simeq \rho^\bot$.
We assume that 
\begin{itemize}
\item[(i)] either $\rho$ is a sum of pairwise distinct irreducible representations,
\item[(ii)] or $\rho(H)$ is abelian.
\end{itemize}
Then a parameter $\psi : H \rightarrow {}^LU(m)$ extending $\rho$ is unique up to
conjugation by an element of $\GL_m(\C)$ (if it exists).
\end{prop}
\begin{pf}
As $\rho$ is semisimple, we may and do assume, up to changing $\rho$ by a conjugate, that the algebra generated
by $\rho(H)$ is stable by transposition.

Let $\psi$ (resp. $\psi'$) be a parameter of $G$ whose restriction to $H$ is $\rho$,
and let $C$ (resp. $C'$) be as above
the matrix such that $\psi(d)=C \phi_m^{-1}
c$ (resp. $\psi'(d)=C' \phi_m^{-1} c$). The matrix $C$ (resp. $C'$) satisfies
(\ref{autod}) and (\ref{d2}). Hence:
\begin{itemize}
\item[(a)] $C C'^{-1}$ is in the centralizer of $\rho(H)$.
\item[(b)] $C{}^t C^{-1}=C' {}^t {C'}^{-1}=(-1)^{m+1} \rho(d^2)$.
\end{itemize}

We want to find a matrix $B \in \GL_m(\C)$ in the centralizer of $\rho(H)$ such that
$$B C {} {}^tB = C'.$$
Indeed, it is clear that for such a $B$, $\psi'=B\psi B^{-1}$.

Assume first that we are in {\bf case (i)}. In this case, we may write
$\rho=\rho_1 \oplus \dots \oplus \rho_r$, with $\rho_1,\dots,\rho_r$
irreducible of dimension $d_1,\dots,d_r$, and choose a basis in which
$\rho(h)=\diag(\rho_1(h),\dots,\rho_r(h))$ for every $h \in H$.
Since the $\rho_i$ are distinct, and since $C$ satisfies (\ref{autod}),
it may be written as $$C=C_0 \sigma$$
where $C_0$ is in the standard Levi $L \subset \GL_m(\C)$ of type $(d_1,\dots,d_r)$ and $\sigma$ is the
 permutation of $\{1,\dots,r\}$ satisfying
$\rho_{\sigma(i)}^\bot\simeq \rho_i$ (hence $\sigma^2=1$)
seen as a permutation matrix (in $\Gl_m(\C)$) by blocks of type $(d_1,\dots,d_r)$. By (a),
we may write $$C'=CD$$ where $D$ is in the centralizer of $\rho(H)$, that is of the Levi $L$,
hence it is of the form
$D=\diag(a_1,\dots,a_1,a_2,\dots,a_2,\dots,a_r,\dots,a_r)$ where each $a_i$ is repeated $d_i$ times.
Now we see that
\begin{eqnarray*} C {}^t C^{-1} &=& C'{}^t{C'}^{-1} \text{ by (b)}
\\ &=&(C_0 \sigma D) ({}^t C_0^{-1} \sigma D^{-1})
\text{ using that ${}^t \sigma = \sigma^{-1}=\sigma$ and that ${}^t D=D$} \\
&=& (C_0 \sigma) ({}^t C_0^{-1} \sigma) (\sigma D \sigma D^{-1}) \text{ using that $DC_0=C_0D$}\\
&=& C {}^t C^{-1} (\sigma D \sigma D^{-1}). \end{eqnarray*}
Hence $$\sigma D \sigma = D,$$
and $a_i=a_{\sigma(i)}$ for all $i$.
We thus may choose complex numbers $b_i$, $i=1,\dots,r$ such that $b_i b_{\sigma(i)} = a_i$ for all $i$,
and set
$B=\diag(b_1,\dots,b_1,b_2,\dots,b_2,\dots,b_r,\dots,b_r)$ where each $b_i$ is repeated $d_i$ times.
Then $\sigma B \sigma {}^t B=D$, $B$ is in the centralizer of $\rho(H)$ and
$$B C {}^tB = B C_0 \sigma {}^t B=C_0 \sigma \sigma B \sigma {}^t B = C D =C',$$
and we are done in case (i).

Assume now that we are in {\bf case (ii)}. Then $\rho$ is a sum of distinct characters
$\chi_1,\dots,\chi_r$, each of them with multiplicity $m_1,\dots,m_r$.
So we have $m=m_1+\dots+m_r$.
We may assume that $\rho$ acts by $\chi_1$ on the $m_1$ first vectors of the basis,
then by $\chi_2$ on the $m_2$ next vectors, and so on.
Hence $\rho(H)$ is made of diagonal matrices of the form
$\diag(a_1,\dots,a_1,a_2,\dots,a_2,\dots,a_r,\dots,a_r)$ where each $a_i$ is repeated $m_i$ times. In particular, 
$(-1)^{m+1} \rho(d^2)=\diag(d_1,\dots,d_1,d_2,\dots,d_2,\dots,d_r,\dots,d_r)$ is of that form by (b),
and the centralizer of $\rho(H)$ is the standard Levi $L$ of type $(m_1,\dots,m_r)$.

Since $C$ (resp. $C'$) satisfies (\ref{autod}), it may be written
as $$C=C_0 \sigma \text{ (resp. }C'=C'_0 \sigma\text{)}$$
where $C_0=\diag(C_1,\dots,C_r)$ (resp. $C'_0=\diag(C'_1,\dots,C'_r)$) is
in the Levi $L$  and $\sigma$ is the
 permutation of $\{1,\dots,r\}$ satisfying
$\chi_{\sigma(i)}^\bot\simeq \chi_i$ (hence $\sigma^2=1$)
seen as a permutation matrix (in $\Gl_m(\C)$) by blocks of type $(m_1,\dots,m_r)$.

Fix an $i \in \{1,\dots,r\}$. By (b) we see that
$$C_i {}^t C_{\sigma(i)}^{-1}=C'_i {}^t{C'_{\sigma(i)}}^{-1} = d_i.$$
If $\sigma(i)=i$, then $d_i \in \C^*$ satisfies $d_i^2=1$, which implies that $C_i$ and $C'_i$ are either both symmetric or both antisymmetric,
and in any case that there exists a $B_i \in \Gl_{m_i}(\C)$ such that
$B_i C_i {}^t B_i = C'_i$. In the other case we have $\sigma(i)=j \neq i$, 
and we assume that $i<j$ to fix the ideas. Then $m_i=m_j$ and
$C_j {}^t C_i^{-1} =C'_j {}^t {C'_i}^{-1}$. From 
that equality, it follows that if we set $B_i=C'_i C_i^{-1}$ and $B_j=\Id_{m_i}$, then
$$\diag(B_i,B_j) \diag(C_i,C_j) \tau {}^t\diag(B_i,B_j) = \diag(C'_i,C'_j)$$
where $\tau$ is the restriction of $\sigma$ to the set $\{i,j\}$ 
(and is the only non-trivial permutation of that set) seen as a matrix 
by blocks in $\Gl_{2m_i}(\C)$. 

Finally $B:=\diag(B_1,\dots,B_r) \in \GL_m(\C)$ is in the centralizer of
$\rho(H)$ and clearly satisfies $B C {}^t B=C'$, and we are done.
\end{pf}

\begin{example}\label{parametreseriesprincipales}Here is an example where the case (ii) of the proposition above may be used. 
Assume that $E/F$ is a quadratic extension of local fields and that $\U(m)$ is the quasisplit 
unitay group in $m$ variables attached to $E/F$. Let $T$ be a maximal torus of $\U(m)(F)$ 
and $\chi: T(F) \longrightarrow \C^*$ an admissible character of $T(F)$. By the duality for tori 
this character $\chi$ defines an $L$-parameter 
$$\phi(\chi): W_F \longrightarrow {}^L T \subset {}^L \U(m).$$ In this book, the only representations of $\U(m)(F)$ that we consider for $F$ nonarchimedian (actually 
$F=\Q_l$) are in $L$-packets of this type. Precisely, they will be either unramified, 
which means that $\chi$ is trivial on the maximal compact subgroup of $T(F)$, 
or of NMSRPS type, in which case $\chi$ is as in Def. \ref{NMSRPS}. 

The norm map $\Norm: T(E) \longrightarrow T(F)$ (see \S\ref{construpinnonsplit}) defines a {\it base change}
$$\chi_E:=\chi \circ \Norm: T(E) \longrightarrow \C^*$$ and $\phi(\chi)_E$ is 
then the $L$-parameter of $\GL_m/E$ attached in the same way to $\chi_E$. 
Case (ii) of Proposition~\ref{propparunicity} shows that $\phi(\chi)$ is actually (up to conjugation) the unique $L$-parameter of $\U(m)$ whose 
base change is $\phi(\chi)_E$. Conversely, we may start from a principal series $L$-parameter $\phi_E: W_E \longrightarrow \GL_m(\C)$, or 
which is equivalent from a character $\chi_E : T(E) \longrightarrow \C^*$, and ask whether it descends to an 
$L$-parameter of $\U(m)$ of the form $\phi(\chi)$ for some $\chi$ as above. This is requires strong conditions on $\chi_E$. In the two cases (unramified or NMSPRS) we are interested in, 
this analysis is precisely the work done in \S\ref{construpinnonsplit}: Lemma \ref{descentedescaracteres} shows that the 
conditions (iiia) or (iiib) of \S\ref{start} on the (unique) representation $\pi_l$ of $\GL_m(E)$ whose $L$-parameter is $\phi_E$ are sufficient.
\end{example}

This quite general uniqueness result is completed by the following more restrictive,
but still very useful, existence
result. We suppose given a subgroup $G'$ of $G$ which is not a subgroup of
$H$. Hence $H':=G' \cap H$ has index two in $G'$. We choose the element $d$ in $G'-H'$.

\begin{prop}\label{proppar} Let $\rho:H \rightarrow \Gl_m(\C)$ be a representation such
that $\rho^\bot \simeq \rho$. Let $\rho'$ be its restriction to $H'$.
We assume that $\rho'$ is a sum of distinct irreducible representations $\rho'_i$
such that $\rho'_i \simeq {\rho'_i}^\bot$. Then $\rho$ extends to a unique parameter $\psi$ of
$G$ whenever $\rho'$ extends to a parameter $\psi'$ of $G'$.

Moreover, if this holds, then the centralizer in $\Gl_m(\C)$ of the image of $\psi$ 
is finite.
\end{prop}
\begin{pf} Let $C'$ be the matrix attached to the
parameter $\psi'$ of $G'$. Arguing as in the proof of the above proposition (case (i)) 
applied to $\rho'$, we may assume that $\rho'(H')$ lies in the standard Levi $L$ of type
$(d_1,\dots,d_r)$ (here $d_i = \dim \rho'_i$),
that the centralizer of $\rho'(H')=\rho(H')$ is the centralizer (and the center)
$Z(L)$ of $L$, and that $C'$ is in $L$ (note that $\rho'_i \simeq {\rho'_j}^{\bot}$ 
if and only if $i=j$, so that $\sigma=\Id$).

Now let $C$ be a matrix that satisfies (\ref{autod}) for $\rho$.
Then $CC'^{-1}$ centralizes $\rho(H')$, hence it lies in the centralizer $Z(L)$ of $L$, and
commutes with $C'$. From that we deduce
that $$C{}^tC^{-1}=C'{}^tC'^{-1}$$
and condition~(\ref{d2}) holds for $C$ since by hypothesis it holds for $C'$.
Hence the existence of a parameter $\psi$ whose restriction to $H$ is $\rho$
follows from Lemma~\ref{lemmepar}. 
The uniqueness follows immediately from Proposition~\ref{propparunicity} (case (i)).

Finally, the centralizer $C(\psi) \subset \GL_m(\C)$ of $\psi(G)$ is a subgroup of the analogous 
centralizer $C(\psi')$ of $\psi(G')=\psi'(G')$. This is the subgroup of
 the centralizer of $\psi'(H')=\rho'(H')$ that is fixed by the map $g \mapsto C{}^tg^{-1}C^{-1}$.
Since the
centralizer of $\rho'(H')$ is $Z(L)$, and since $C \in L$, this map is $g \mapsto g^{-1}$ and
$C(\psi')\simeq (\Z/2\Z)^r$. Hence $C(\psi)$ is finite.
\end{pf}

\begin{remark} If $\psi$ is a {\it discrete} $A$-parameter for ${}^L U(m)$, then Rogawski 
(\cite[Lemma 2.2.1, 2.2.2]{Rog1}) has shown that its
restriction $\rho$ to $L_E \times \Sl_2(\C)$ is a sum of irreducible, pairwise non isomorphic representations $\rho_i$
satisfying $\rho_i^\bot=\rho_i$. The above proposition, in the case $G=G'$, provides a converse to this result.
\end{remark}

\begin{cor}\label{corpar} Let $E/F$ be a CM extension of number fields, and
$$\rho: L_E  \longrightarrow \GL_m(\C)$$ a tempered $A$-parameter for $\GL_m/E$ that
satisfies $\rho^\bot \simeq \rho$. Assume that there is an infinite place
of $F$ such that, for the corresponding inclusion $W_\C \hookrightarrow L_E$
(well defined up to conjugation in $L_E$), the restriction of $\rho$ to $W_\C$
extends to a {\rm discrete} $L$-parameter of $W_\R$ (see Remark~\ref{remdiscreteinf} below). Then $\rho$ extends
to a discrete tempered $A$-parameter $L_F \longrightarrow {}^L U(m)$ of the quasisplit unitary group $\U(m)$.
\end{cor}

\begin{pf} This is the proposition for $G=L_F$, $H=L_E$, $G'=W_\R$,
$H'=W_\C$. The hypothesis of the proposition on the restriction
$\rho'$ of $\rho$ to $H'$ follows immediately from the Remark \ref{remdiscreteinf} below.
The obtained parameter is discrete as is its restriction to $W_\R$.
\end{pf}

Via the philosophy of parameters, this result shows that a sufficient condition
for a cuspidal automorphic representation $\pi$ of $\Gl_m(\AAA_E)$ that satisfies
$\pi^{\bot}\simeq \pi$ to come by base change from the quasi-split unitary group in
$m$ variables attached to $E/F$ is that there is a place $v$ at infinity such that $\pi_v$ comes
 by base change from a representation of the {\it compact} unitary group.
Note that this is exactly the assumption at infinity on the representations $\pi$
studied by Clozel \cite{clo} and Harris-Taylor \cite{HT}.

\begin{remark}\label{remdiscreteinf}
Recall that the discrete $L$-parameters $\phi: W_\R \longrightarrow {}^L \U(m)$
are exactly the ones whose restriction $\phi_{\C}$ to $\C^*$ is conjugate to
$$z \mapsto ((z/\bar z)^{a_1},\dots,(z/\bar z)^{a_n}),$$ where the
$a_i$ are in $\frac{m+1}{2}+\Z$ and strictly decreasing (see e.g. \cite[Prop. 4.3.2]{BeCl}). For each such sequence
$(a_i)$ there is a unique such parameter.\footnote{This existence and uniqueness follows for 
example easily from Lemma \ref{lemmepar}, and Prop. \ref{proppar}, \ref{propparunicity}:
choose $d$ to be the usual element $j \in W_\R$ such that $j^2=-1$, $C=1$ and note that
$\rho(j^2)=\rho(-1)=(-1)^{m+1}$. These parameters $\phi$ satisfy $\phi(j)=\phi_m^{-1}c$. They are relevant because,
if they lie in a parabolic subgroup $P$ of ${}^L \U(m)$, then $P=\langle P_0, \phi(j)\rangle$
for some parabolic $P_0$ of
$\GL_m(\C)$ normalised by $\phi(j)$, and we see that ${}^t P_0 = P_0$, hence $P={}^L \U(m)$.}
As they are discrete, they are relevant for each inner form of $\U(m)_{\R}$,
and in particular for the compact one.
\end{remark}

\begin{prop}\label{locequivimplequiv} For unitary groups $\U(m)$, two everywhere locally
equivalent discrete $A$-parameters are actually equivalent.
\end{prop}

\begin{pf} Let $\psi_1$ and $\psi_2$ be two such discrete $A$-parameters, and set $G=G'=L_F \times \SL_2(\C)$,
$H=H'=L_E \times \SL_2(\C)$, and $\rho_j={\psi_j}_{|H}$.
By Rogawski's classification, the $\rho_j$ are semisimple and satisfy the assumption of Prop.~\ref{proppar}.
By assumption, the $\rho_j$ are also everywhere locally equivalent. By `Cebotarev's theorem
for Langlands' group'', the reunion of the conjugates of $L_{E_w}$ is dense in $L_E$, hence a trace consideration
implies that the $\rho_j$ are actually equivalent. By Prop.~\ref{proppar}, the same thing holds then for the $\psi_j$.
\end{pf}

The following corollary is an immediate consequence of the proposition above
and of the simplest case of Arthur's mutliplicity formula.

\begin{cor} If $\Pi$ is a stable $A$-packet for $\U(m)$, then for each $\pi \in \Pi$ we should have $m(\pi)=1$.
\end{cor}

\begin{pf} Indeed, there is a unique $A$-parameter $\psi$ of $\U(m)$ giving rise to $\Pi$ by
Prop. \ref{locequivimplequiv}. As $\Pi$ is stable,
Arthur's group $\bs_{\psi}$ is trivial by \cite[\S 2.2]{Rog1},
hence Arthur's multiplicity formula (\cite[(8.5)]{Arthur}) reduces to $m_{\psi}(\pi)=m(\pi)=1$.
\end{pf}

\subsection{Parameter and packet of the representation $\pi^n$}\label{apretpackpin}
${ }^{ }$
\par \medskip
{\it In this \S\  and the next one, we take as granted all the formalism of 
Langlands and Arthur as described above and in \cite{Arthur}. 
All the lemmas, propositions and theorems we state are thus conditional to 
this formalism. Our aim is to study the $A$-packet of the representation $\pi^n$ that we introduced in \S\ref{constructionpin}.}

\par \medskip
 
\label{parpin}

We use from now the notations of~\S\ref{constructionpin}. In particular $m=n+2$, 
$n$ is not a multiple of $4$, $\U(m)$ is definite and quasisplit at all finite places, and 
we fix an embedding $E \rightarrow \C$. Moreover, $\mu=\mu^{\bot}$ 
is a Hecke character of $\AAA_E/E^\ast$ as in Notation \ref{choixdemu}. 
Remember that $\mu$ is trivial if $m$ (or $n$) is even, and that $\mu(z)=(z/\bar z)^{1/2}$ for $z \in \C^\ast \subset \AAA_E^\ast$ if $n$ is odd. 
We will see $\mu$ as a character of $W_E$, hence of $L_E$, when needed.

To the representation $\pi$ of that subsection 
should correspond a tempered, irreducible, 
$A$-parameter $$\rho : L_E \longrightarrow \Gl_n(\C).$$ 
The hypothesis (i) on $\pi$ there translates to\footnote{Strictly, this uses the Cebotarev theorem for $L_E$, or (which is related) the weak mutliplicity one theorem for the discrete spectrum of $\GL_m/E$.} $$\rho^\bot \simeq \rho.$$
 
We now define an $A$-parameter for $\Gl_m(E)$ denoted
$\psi_{E} : L_E \times \Sl_2(\C) \longrightarrow \Gl_m(\C)$
by $$\psi_{E} (w \times g) = \left( \begin{array}{cc} 
\rho(w) \mu(w) & 0 \\ 0 & g \mu(w) \end{array} \right), \, \, \, w \in L_E,\, \, g \in \SL_2(\C).$$
 
We fix an embedding $L_\R = W_\R \hookrightarrow L_\Q$, giving an embedding $W_\C=\C^\ast \hookrightarrow L_\Q$. 
The image of the element $j$ of $W_\R$ is the unique non trivial 
element $c$ in $L_\Q/L_E = W_\R/W_\C=\Gal(E/\Q)$. 
By hypothesis (ii) of \S\ref{start}, we may and do assume
(up to changing $\rho$ to a conjugate) that for $z \in W_\C = \C^\ast$,
$$\psi_E(z)=\diag((z/\bar z)^{a_1},\dots,(z/\bar z)^{a_n})$$ where the 
$a_i$ are in $\frac{1}{2}\Z-\Z$, strictly decreasing, and different from $\pm 1/2$.

\begin{lemma}\label{lemmapsiextends} The $A$-parameter $\psi_E$ extends (uniquely
up to isomorphism, that is up to conjugation) to a discrete (relevant) $A$-parameter 
$\psi : L_\Q \times \SL_2(\C) \rightarrow {}^L \U(m)$ 
of the group $\U(m)$. We may choose  
$$\psi(j)=\left( \begin{array}{ccc} 1_n & 0 & 0 \\ 0 & 0 & 1 \\ 0 & -1 & 0 \end{array} \right) \phi_m^{-1} c.$$
\end{lemma}
\begin{pf} The restriction of $\psi_E$ to $W_\C \times \Sl_2(\C)$ is 
$$\psi_{E}(z,g)=\diag((z/ \bar z)^{a_1} \mu_\infty(z), \cdots,
(z/\bar z)^{a_n} \mu_\infty(z),g \mu_\infty(z)).$$
We thus see that for $$C:=\left( \begin{array}{ccc} 1_n & 0 & 0 
\\ 0 & 0 & 1 \\ 0 & -1 & 0 \end{array} \right)$$ we have for all $z,g$ :
$$\psi_{E}(\bar z,g)=\psi_{E}(jzj^{-1},g)=
C{}^t\psi(z,g)^{-1}C^{-1}.$$
In other words, 
the relation (\ref{psidinv}) holds for that $C$. We compute 
$$\psi_{E}(j^2,1)=\diag((-1)^{m+1},\dots,(-1)^{m+1},(-1)^m,(-1)^m),$$
using that the $a_i$ are half-integers and that $\mu_\infty(-1)=(-1)^m$. On the 
other hand, $C{}^tC^{-1}=\diag(1,\dots,1,-1,-1)$ so the 
relation~(\ref{d2}) holds for the restriction of $\psi_{E}$ to $W_\C \times \Sl_2(\C)$. 
This means by Lemma~\ref{lemmepar},  that the restriction of 
$\psi_{E}$ to $W_\C \times \Sl_2(\C)$ extends to a parameter 
$\psi_{\infty}$ of 
$W_\R \times \Sl_2(\C)$, that moreover sends $j$ to $C\phi_m^{-1}$ 
(see~(\ref{autod})). By Proposition~\ref{proppar}, $\psi_E$ 
extends to a unique, discrete, parameter $\psi$ of $L_\Q \times \Sl_2(\C)$, and we may even choose $\psi$ such that 
$\psi(j)=C\phi_m^{-1}$ by the last paragraph of the proof of that proposition.

It remains to show that $\psi$ is relevant, which means, since $\U(m)$ 
has no proper parabolic defined over $\Q$, that the only parabolic subgroup $P$ of ${}^L \U(m)$ 
containing $\Im(\psi)$ is ${}^L \U(m)$ itself. Let $P$ be such a subgroup, by definition 
$P=\langle P_0, \psi(j)\rangle$ for some parabolic $P_0$ of $\GL_m(\C)$ normalized by $\psi(j)$. We see then that
$C {}^t P_0 C^{-1}= P_0$. But $P_0$ contains $1_n \times \SL_2(\C)=\psi(\SL_2(\C))$, hence $C \in P_0$. As a consequence, 
$P_0={}^t P_0$, which implies that $P_0=\GL_m(\C)$, and we are done. 
\end{pf} 

\begin{remark}\label{remautremethodeunp2} Here is another way to view the parameter $\psi$ in terms of Rogawski's classification recalled in 
\S\ref{descarthunit}. First, using Cor.~\ref{corpar}, we see 
that $\rho\mu$ extends to a discrete stable tempered parameter $\psi_0$ of the quasisplit $\U(n)$. Then, using again $\mu$ to 
define an admissible morphism $$\xi: {}^L ( \U(n) \times \U(2) ) \longrightarrow {}^L \U(m),$$
we have actually $\psi=\xi \circ (\psi_0 \times (1\otimes [2]))$ (see Ex.\ref{caseglm2}). In particular, $\psi$ is nontempered, and 
endoscopic of type $(n,2)$. When $n=1$, it is exactly the $A$-parameter that we studied in details in \S\ref{instructiveexample}.
\end{remark}

Let us denote by $$\Pi=\prod'_v \Pi_v$$ 
the $A$-packet corresponding to $\psi$. Our aim is now to check that the representation 
$\pi^n$ defined in \S\ref{constructionpin} belongs to $\Pi$. By definition, this amounts to check that for each place $v$ the representation 
$\pi^n_v$ defined there lies in $\Pi_v$. Recall that for some reasons we have called $\pi_{\infty}^s$ the archimedian component of $\pi^n$ (see \S\ref{construpisinfty}).
%Let us choose for each prime $l$ a distinguished embedding $L_{\Q_l} \hookrightarrow L_\Q$.

\begin{lemma} \label{lemmaapacket} The global representation $\pi^n$ belongs to the global $A$-packet defined by $\psi$. 
Moreover, $\Pi_{\infty}=\{\pi^s_{\infty}\}$ and for each prime $l$, $\pi^n_l$ is in the local $L$-packet 
$\Pi_{\varphi_l} \subset  \Pi_l$ (see \S\ref{landalocal}). 
\end{lemma}

\begin{pf}
By definition, for each place $v$ of $\Q$ the $L$-parameter $\phi_v:=\phi_{\psi_v}$ associated to $\psi_v$ satisfies for all $w \in L_E$
 $$\phi_v(w) = 
\diag( \rho_v(w) \mu_v(w), |w|^{1/2} \mu_v(w), |w|^{-1/2}\mu_v(w)).$$ 
For any place $v$ of $\Q$, the $L$-parameter $\phi_v$ defines an $L$-packet $\Pi_{\phi_v}$ of representations of $\U(m)(\Q_v)$, which is a 
subset of the $A$-packet $\Pi_v$.

Assume first that $v=l$ is a prime. When $l$ splits in $E$, Remark~\ref{remlsplit} shows 
that $$\Pi_{\varphi_l}=\{\pi^n_l\} \subset \Pi_l.$$ When $l$ does not split, we defined in 
Lemma~\ref{descentedescaracteres} a smooth character $\chi$ of the maximal torus $T(\Q_l)$ of 
$\U(m)(\Q_l)$. There were two cases. If $\chi$ satisfies conditions (a) and (b) of Def. \ref{NMSRPS}, then 
$\pi^n_l$ is the NMSRPS representation $S(\chi)$. As suggested by Rodier's work \cite{Rodier2}, the $L$-parameter of 
$S(\chi)$ is conjecturally the $L$-parameter $\phi(\chi)$ defined in Remark \ref{parametreseriesprincipales}. But by the same remark 
and by construction, 
the base change $\phi(\chi)_{E_l}$ is isomorphic to $(\phi_l)_{E_l}$, hence $$\phi(\chi) \simeq \phi_l$$ by Prop. \ref{propparunicity} and 
we are done. In the other case $\chi$ is unramified and $\pi_l^n$ is by definition a constituent of the full induced 
representation defined by $\chi$ having a nonzero vector invariant by a maximal hyperspecial (resp. very special) compact open subgroup of $\U(m)(\Q_l)$ 
if $l$ is inert (resp. ramified) in $E$. Thus $\pi^n_l$ conjecturally belongs again to the L-packet defined by $\phi(\chi)$ 
(this is a standard expectation when $l$ is inert, and it is indicated by Labesse's work \cite[\S 3.6]{labesselivre} in general), 
and we conclude as above that $\phi(\chi) \simeq \phi_l$. 

The end of the proof is now devoted to the subtler case where $v=\infty$ is archimedian. 
Note that the $L$-parameter $\phi_\infty$ is not relevant for the compact group 
$\U(m)(\R)$. Indeed, it contains the non unitary characters $\mu_\infty|\, |^{\pm 1/2}$, thus it cannot be discrete 
(see \S\ref{remdiscreteinf}). As a consequence, $$\Pi_{\phi_{\infty}} =\emptyset,$$ 
but the $A$-packet $\Pi_{\infty}$ may be larger. However, note that $\Pi_{\infty}$ is a 
singleton if nonempty, since every representation of a compact real reductive group is stable (cf. \cite{AJ}).
We shall review below the description of $\Pi_{\infty}$ given in section 5 of \cite{Arthur}, following \cite{AJ}. 

For this, we see $\U(m)(\R)$ as the unitary group for the standard diagonal positive definite hermitian form,
and we consider its diagonal maximal torus $T(\R)=\U(1)(\R)^m$. We denote by $$L(\R) = U(1)(\R)^n \times U(2)(\R) \subset \U(m)(\R)$$
the subgroup of matrices which are diagonal by blocks of size $(1,\dots,1,2)$, so that $T \subset L$ 
and $L(\C)$ is a Levi subgroup of $\U(m)(\C)=\Gl_m(\C)$. In $\hat{G}(\C)\simeq \Gl_m(C)$, $\hat{T}(\C)$ is the diagonal torus and $\hat{L}(\C)$ the standard Levi of type $(1,\dots,1,2)$.
We thus have $$Z(\hat{L}(\C)) \subset \hat{T}(\C) \subset \hat{L}(\C) \subset \hat{G}(\C).$$
It turns out that those inclusions extend naturally to inclusions $${}^L 
Z(L) \subset {}^L T \stackrel{\xi_{L,T}}{\hookrightarrow} 
{}^L L \stackrel{\xi_{G,L}} \hookrightarrow {}^L G.$$ While the first inclusion is obvious, the others two need a construction, 
which is recalled in \cite{Arthur}. From this construction, we shall only need the following description of the restriction 
of the embedding $\xi_{L,T}$ to $W_\C=\C^\ast$ (see \cite[page 31]{Arthur}) :
\begin{eqnarray*} \xi_{L,T}(z) = \diag(1,\dots,1,(z/\bar z)^{1/2},(z/\bar z)^{-1/2}). \end{eqnarray*}
Similarly $\xi_{G,L}(z) \in \hat{T}(\C)$ is a diagonal matrix that we do not need to compute explicitly because it will 
cancel out in the following computations. 

Now let us consider the unique $L$-parameter $\phi_\tau$ for the group $T$ such that for $z \in W_\C = \C^\ast$, 
$$\phi_\tau(z) = \xi_{G,L}(z)^{-1} \diag((z/\bar z)^{a_1} \mu_\infty(z),\dots,(z/\bar z)^{a_n} \mu_\infty(z),\mu_\infty(z),
\mu_\infty(z)).$$ 
It is clear that $\phi_\tau$ maps $W_\R$ to ${}^L \hat{Z(L)} \subset {}^L T$.
To such an $L$-parameter is attached in \cite[page 30 second \S, page 31 first \S]{Arthur}
an $A$-parameter for the group $L$, called $$\psi_L:\,W_R \times \Sl_2(\C) \longrightarrow {}^L L.$$ By 
definition, $\psi_L$ is $\phi_\tau$ on $W_\R$ and sends 
$\left( \begin{matrix} 1 & 1 \\ 0 & 1 \end{matrix} \right) \in \Sl_2(\C)$ to a principal unipotent element in $\hat{L}$. 
Thus it is clear that (up to conjugation by an element of $\hat{L}(\C)$) the $A$-parameter 
$\xi_{G,L} \circ \psi_L$ coincides with our $A$-parameter $\psi_{\infty}$ on $W_\C$ and $\Sl_2(\C)$. Thus $\psi_{\infty} \simeq \psi_L$ 
by Lemma~\ref{propparunicity}.

Arthur's conjecture provides a description (cf. \cite[page 33]{Arthur}) 
of the $A$-packet attached to $\xi_{G,L} \circ \psi_{L}=\psi_{\infty}$. He defines for that some $L$-packets parameterized
by the set $$\Sigma:=W(L,T) \backslash W(G,T) / W_\R(G,T),$$ and for each element of this set a specific representation in the associated $L$-packet. 
Here $\Sigma=\{1\}$, and the unique $L$-parameter he defines is $\phi_1:=\xi_{G,L} \circ \xi_{L,T} \circ \phi_\tau$ 
as $L$ is anisotropic \cite[page30-31]{Arthur}. On $W_\C$, we thus have 
\begin{eqnarray*} \phi_1(z)&=&\xi_{G,L}(z) \xi_{L,T}(z) \phi_\tau(z)\\ 
&=& \diag((z/\bar z)^{a_1} \mu_\infty(z),
\dots,(z/\bar z)^{a_n} \mu_\infty(z), (z/\bar z)^{1/2} \mu_\infty(z),(z/\bar z)^{-1/2} \mu_\infty(z)).\end{eqnarray*}
Note that $\phi_1$ is relevant since the $a_i,1/2,-1/2$ are distinct half-integers.
Actually, $\phi_1$ is exactly by definition the $L$-parameter of $\pi^s_\infty$, and its associated $L$-packet is a singleton. 
According to Arthur, we thus have $\Pi_{\psi_\infty}= \Pi_{\phi_1}=\{\pi_\infty^s\}$.
\end{pf}
\begin{remark} \begin{itemize} \label{trvialApacket}
\item[(i)] In Lemma \ref{lemmapsiextends}, and especially in the proof that $\psi$ is relevant,
the fact that the $a_i$ are distinct from $\pm 1/2$ is actually not needed. However, as the proof above shows,
this latter assumption is necessary to ensure that $\Pi_{\infty}$ (hence $\Pi$) is not empty.
In particular, if one of the $a_i$ is equal to $\pm 1/2$, we get an example of a parameter $\psi$ which is relevant 
and whose associated $A$-packet is empty.
\item[(ii)] \label{nonquasisplitrelevant} As an exercise, the reader can check that the $A$-parameter $\psi$ is not 
relevant for an inner form of $\U(m)$ that is not quasi-split at every 
finite place.
\item[(iii)] We have $\Pi_{\infty}=\{\pi^s_{\infty}\}$ and $\Pi_l=\{\pi^n_l\}$ when $l$ splits in $E$. When $l$ does not split, $\Pi_l$ (and even $\Pi_{\phi_l}$) will have more that one element in general, but it does not seem possible at the moment to describe the full packet $\Pi_l$ for a general $m$ as we have done for $m=3$ in \S\ref{instructiveexample}.
\end{itemize}
\end{remark}

\begin{remark}\label{autrepossibiliteun} In \S\ref{constructionpin}, our point of view for defining the 
representation $\pi^n$ was to start from a cuspidal automorphic representation $\pi$ of $\GL_n(\AAA_E)$ such that 
$\pi^{\bot} \simeq \pi$ and that satifies \S~\ref{start} (ii). We defined then {\it by hand} all the local components 
of the representation $\pi^n$ of $\U(m)(\AAA_{\Q})$. An interest of this presentation is that it avoids admitting that $\pi\mu$ comes 
by base change from an automorphic representation $\pi_0$ of the quasisplit unitary group $\U(n)$, and more generally any appeal to a unitary group in $n$ variables.
As noticed in Remark \ref{remautremethodeunp2}, such a $\pi_0$ should however always exist. More precisely, 
arguing as in Rem.\ref{parametreseriesprincipales}, conditions (iiia) and (iiib) of \S~\ref{start} imply 
that we should be able to find such a $\pi_0$ which is either unramified or NMSRPS at all nonsplit finite places. 
Conversely, if $\pi_0$ is such a representation which satisfies also (ii) of 
\S\ref{start} and which has a cuspidal base change $\pi$ to $\GL_n/E$
(this is known to hold for example when $\pi_0$ is supercuspidal at two split places (\cite[Thm.2.1.1,3.1.3]{HL}), then $\pi$ satisfies 
our conditions. This gives a way to produce such examples (maybe using some inner forms of $\U(n)$ as well).
\end{remark}

\subsection{Arthur's multiplicity formula for $\pi^n$}\label{formulemultart}

We have checked that $\pi^n \in \Pi$ and we ask now if $\pi^n \in \Pi_\disc(\U(m),\Q)$. For that we will actually compute $m(\pi^n)$ using Arhtur's multiplicity formula.

Following \cite[page 52]{Arthur}, let us consider the subgroup 
$S_{\psi} \subset \hat{G}(\C)=\Gl_m(\C)$. In our situation\footnote{Precisely, using the triviality of $\ker^1(L_\Q,Z(\hat{G}(\C)))$ for $G=\U(m)$, as remarked in \cite[\S 2.2]{Rog1}.}, 
this group is actually $Z(\hat{G}(\C))\cdot C_{\psi}$ where $C_{\psi}$ is the centraliser in $\hat{G}(\C)$ of the image of $\psi$. We set also 
$s_{\psi}=\psi(1,\diag(-1,-1)) \in S_{\psi}$ (see \cite[page 26]{Arthur}) and 
 $$\bs_{\psi}=S_{\psi}/S_{\psi}^0Z(\hat{G}(\C)).$$

\begin{lemma} \label{lemmaSpsi} $$S_{\psi}=\{(\diag(a,\dots,a,\epsilon a,\epsilon a),\ \  a \in \C^{\ast}, \epsilon=\pm 1\} \simeq \C^\ast \times \{\pm 1\}.$$ The character $S_{\psi} \rightarrow \{\pm 1\}$ sending $\diag(a,\dots,a,\epsilon a,\epsilon a)$ to $\epsilon$ factors through $\bs_{\psi}$ and induces an isomorphism $\epsilon: \bs_{\psi} \isomo \{\pm 1\}$. Moreover, $\bs_{\psi}$ is generated by the image of $s_\psi$. 
\end{lemma} 
\begin{pf} As $\rho$ is irreducible, the centralizer of $\psi(L_E \times \SL_2(\C))$ in $\GL_m(\C)$ is 
$$\{\diag(a,\dots,a,b,b)\} = \C^\ast \times \C^\ast.$$ Among those elements,
those who commute with $\psi(j)$ are the one of order two modulo $Z(\GL_m(\C))$, hence
$S_{\psi}=\{\diag(a,\dots,a,\epsilon a, \epsilon a)\}$ and 
$\bs_{\psi} = \pi^0(S_{\psi}/Z(\hat{G})) = \{\pm 1\}$. Finally, the element $s_{\psi}$ clearly generates $\bs_\psi$. 
\end{pf}

We now introduce following Arthur (\cite[page 54]{Arthur})
the representation 
$$\tau : S_{\psi} \times L_\Q \times \Sl_2(\C) \longrightarrow \Gl(M_m(\C))$$
defined by $\tau(s,w,g)=\Ad(s\psi(w,g))$. As the kernel of $S_{\psi} \rightarrow \bs_{\psi}$ is the center of $\GL_m(\C)$, the action of $S_{\psi}$ actually factors through $\bs_{\psi}$.
We note $\tau_E$ the restriction of $\tau$ to $S_{\psi} \times L_E \times \Sl_2(\C)$.

Recall that the representation $\psi_E$ on $\C^m$ 
is the direct sum of two representations of 
$L_E \times \Sl_2(\C)$ : the representation $\rho \mu \otimes 1$ on the $n$-dimensional space $V_1$ (generated by the first $n$ vectors of the 
standard basis of $\C^m$) and the 
representation $\mu \otimes [2]$ on the $2$-dimensional space $V_2$ (generated by the last two vectors), 
where $[2]$ is the standard representation of 
$\Sl_2(\C)$. Hence $M_m(\C)$ may be written as $$(V_1 \otimes V_1^\ast) \oplus (V_2 \otimes 
V_2^\ast) \oplus (V_1 \otimes V_2^\ast) \oplus (V_2 \otimes V_1^\ast),$$ all four spaces in the decomposition 
being stable by the adjoint action of $L_E \times \Sl_2(\C)$. 
Moreover, the adjoint action of $S_{\psi}$ (that is, of $\bs_{\psi}$) preserves also this decomposition and is trivial on
$V_1 \otimes V_1^\ast \oplus V_2 \otimes V_2^\ast$, and given by the non-trivial character $\epsilon$ on $V_1 \otimes V_2^\ast \oplus V_2 \otimes V_1^\ast$.  

\begin{lemma}\label{calcadjact} The spaces $V_1 \otimes V_1^\ast$, $V_2 \otimes V_2^\ast$ and $V_1 \otimes V_2^\ast \oplus V_2 \otimes V_1^\ast$ are stable by $\tau$. The last one is isomorphic to 
$\epsilon \otimes \Ind_{L_E}^{L_\Q} \rho \otimes [2]$.
\end{lemma}
\begin{pf} The adjoint action of $\tau(j)$ on $M_\m(\C)$ is given by 
$M \mapsto -C{}^t M C^{-1}$. In particular, it stabilizes 
$V_1 \otimes V_1^\ast$ and 
$V_2 \otimes V_2^\ast$ which thus are stable by $\tau$, and it 
interchanges $V_1 \otimes V_2^\ast$ and $V_2 \otimes V_1^\ast$, from which the lemma follows easily.
\end{pf}

Arthur defines then in \cite[(8.4)]{Arthur} a character 
$$\epsilon_{\psi}: \bs_{\psi} \longrightarrow \{\pm 1\}.$$ 

\begin{prop}\label{carArthur} The character $\epsilon_\psi$ in 
 is the trivial character if $\varepsilon(\pi,1/2)=1$ and the character $\epsilon$ if $\varepsilon (\pi,1/2)=-1$. In other words, $\epsilon_{\psi}(s_{\psi})=\varepsilon (\pi,1/2)$.
\end{prop}
\begin{pf} According to Arthur's recipe, to compute $\epsilon_{\psi}$ we have to decompose the (semisimple) representation $\tau$ into its irreducible components $\tau_k=\lambda_k \otimes \rho_k \otimes \nu_k$ as a representation of $\bs_{\psi} \times L_\Q \times \SL_2(\C)$. By definition, $\epsilon_\psi=\prod_{\tau_k \text{ special}} \lambda_k$ 
where {\it $\tau_k$ special} means that $\tau_k \simeq \tau_k^\ast$, which in our 
context is equivalent to $\rho_k \simeq \rho_k^\ast$, and that $\varepsilon(\rho_k,1/2)=-1$.

As we already saw, we may ignore the $\tau_k$'s 
arising as components of either $V_1 \otimes V_1^\ast$ or 
$V_2 \otimes V_2^\ast$ since the corresponding $\lambda_k$ are trivial.
By Lemma~\ref{calcadjact}, the remaining $\tau_k$'s are the constituents of $\epsilon \otimes \Ind_{L_E}^{L_\Q} \rho \otimes [2]$. Note that $\Ind_{L_E}^{L_\Q} \rho$ is autodual as $\rho^{\bot} \simeq \rho$. 
Let us decompose $\Ind_{L_E}^{L_\Q} \rho$ as a sum of its $r$ irreducible constituents, that we may note $\rho_1,\dots,\rho_r$ up to renumbering. Since $\rho$ is irreducible, we have $r=1$ or $2$. 

If $r=1$, $\tau_1=\tau_1^\ast$ so $\tau_1$ is special if and only if 
$\varepsilon(\rho_1,1/2)=-1$, but $\varepsilon(\rho_1,1/2)=\varepsilon(\Ind_{L_E}^{L_\Q} \rho,1/2)=\varepsilon (\rho,1/2)=\varepsilon(\pi,1/2)$ and the proposition follows. 

The second case $r=2$ occurs exactly when $\rho$ is self-conjugate, hence when $\rho$ is autodual since $\rho\simeq \rho^\bot$. In this case $\rho$ extends to a representation $\rho_1$ of $L_\Q$, and we have $\Ind_{L_E}^{L_\Q} \rho = \rho_1 \oplus \rho_2 = \rho_1 \oplus \rho_1\omega_{E/\Q}$ and $\rho_1 \not\simeq \rho_2$. We have
$$\varepsilon(\pi,1/2)=\varepsilon(\Ind_{L_E}^{L_\Q} \rho,1/2)= \varepsilon(\rho_1,1/2) \varepsilon(\rho_2,1/2).$$
If $\rho_1$ and $\rho_2$ are autodual, we see that there is exactly one (resp. $0$ or $2$) 
$\rho_i$ that is special if $\varepsilon(\pi,1/2)=-1$ (resp $+1$) and the proposition follows. If $\rho_1^*\simeq \rho_2$, then the functional equations of the $L$-functions of $\rho_1$ and $\rho_2$ show that $\varepsilon(\rho_1,s)\varepsilon(\rho_2,1-s)=1$, so $\varepsilon(\pi,1/2)=+1$ and  there are no special $\tau_k$, which concludes this case as well.
\end{pf}

The last ingredient in the multiplicity formula is a conjectural 
canonical 
pairing (\cite[page 54]{Arthur}) $$\bs_\psi \times \Pi \rightarrow \R,\text{ denoted }\langle s,\pi \rangle.$$
 However, this ingredient is certainly the most 
difficult one in 
Arthur's expostion of its multiplicity formula.

Let us recall some features of this pairing in a general context, 
for a reductive group $G/F$, and a global $A$-packet
$\Pi$ with $A$-parameter $\psi$. Together with the global pairing 
should be defined for each place $v$ of $F$ a local pairing $$\bs_{\psi_v} \times \Pi_v \rightarrow \R.$$
However, this local pairing should not be canonical, but rather depends on the 
choice of a basis representation in the local $A$-packet $\Pi_v$.
Still there should be a way, after a global choice $\nu$ (see below), to choose the 
local pairing such that the product formula holds
\begin{eqnarray} \label{product} \langle s , \pi \rangle = 
\prod_v \langle s, \pi_v \rangle_{v,\nu}. \end{eqnarray}
Above $\pi = \otimes'_v \pi_v$ is in $\Pi$, the pairings on the right hand side should be the chosen local pairings depending of the global choice $\nu$, 
$s$ in the left hand side should be any element of $\bs_{\psi}$ and $s$ in the 
right hand side denotes the image of $s$
by the injective natural morphism $\bs_{\psi} \hookrightarrow \bs_{\psi_v}$.
Moreover, almost all the terms in the product should be $1$.

When $G$ is a quasi-split group $G^\ast$, the global choice $\nu$ may be a nondegenerate character of the unipotent radical of a Borel subgroup of $G^*$ defined with the help of a non trivial admissible character $F\backslash \AAA_F \longrightarrow \C^*$. Then for each $v$,
the $A$-packet $\Pi_v$ should contain one and only one representation $\pi_v^{\nu_v-\text{gen}}$
that is $\nu_v$-generic in the sense explained\footnote{Be careful that this representation is $\nu_v$-generic in the usual sense only 
for a temepered $\Pi_v$.} in \cite[4.4]{br2}. This representation should actually belong to the $L$-packet $\Pi_{\psi_v}$. When that representation is 
chosen as the base point to define the local pairing (which are then denoted 
$\langle\ ,\ \rangle_{v,\nu}$) the product formla (\ref{product}) should hold.

From now on, we work with our group $G=\U(m)/\Q$ defined in chapter 6.
Assuming the choice of $\nu$ is made as above 
for its quasi-split form $G^\ast$, we may use the local pairing 
$\langle\ ,\ \rangle_{v,\nu}$ already chosen for $G^\ast_v$ 
for every finite place $v$, since $G^\ast_v \simeq G_v$. Moreover, since there is only one infinite place $\infty$, there is a unique choice of the local
pairing at infinity which makes the formula (\ref{product}) true. We 
still denote it  as $\langle\ , \ \rangle_{\infty,\nu}$.

The pairings above have much nicer features when restricted to the subgroup
(always of order $1$ or $2$) of ${\bs}_\psi$ generated by the canonical element $s_\psi$.
So we are very lucky in our case, because that subgroup is the full 
${\bs}_\psi$ (Lemma~\ref{lemmaSpsi}). In particular, the following assertions should 
hold\footnote{at least in all cases we will use them. We are not completely sure of their generality.}. Below $v$ is a place, 
$\pi_v$ any representation of the local $A$-packet $\Pi_v$.

\begin{itemize}
\item[(a)] 

There should be a sign $e(\pi_v,\nu_v)=\pm 1$ such that 
$$\langle s_\psi^a, \pi_v \rangle_{v,\nu} = e(\pi_v,\nu_v)^a$$ for every 
$a \in \Z$ (or for that matter, for $a=0,1$).

\item[(b)] The sign $e(\pi_v,\nu_v)$ should depend on $\pi_v$ 
only through the $L$-parameter of $\pi_v$, and that sign is $+1$ if this $L$-parameter is $\phi_{\psi_v}$ and if $G_v$ is quasi-split. 
This should be understood in the strong sense that if two $\pi_v$'s, even for two different inner forms of $G_v$, have the same $L$-parameter, then they have the same sign.

\item[(c)] The sign $e(\pi_v,\nu_v)$, hence the pairing
$\langle s_{\psi}, \cdot \ \rangle_{v,\nu}$ is independent of $\nu$.  
 (Hence we may and will drop $\nu$ 
from the notation.) Moreover, any local pairing $\langle\ ,\ \rangle'_v$ such that $\langle s_{\psi}^a, \pi_v \rangle'_v$ is $1$ for $a=0, 1$ and 
a given representation $\pi_v$ of $L$-parameter $\phi_{\psi_v}$ 
is actually equal to that pairing on the subgroup generated by $s_\psi$: 
$$\langle s_{\psi}^a ,\cdot \rangle'_v = \langle s_{\psi}^a ,\cdot \rangle_v.$$
\end{itemize}

Hence the sign $e(\pi_v,\nu_v)$ is simply denoted $e(\pi_v)$, and sometimes even $e(\phi)$ where $\phi$ is the $L$-parameter of $\pi_v$.

Indeed, (a) follows from \cite[Conjecture 6.1(iii)]{Arthur}, the first assertion of (b) is clear
if $v=\infty$ by the description of the pairing given \cite[page 33]{Arthur},
and seems implicitly assumed in the general case. Anyway, we will only use it 
for representations in the canonical $L$-packet $\Pi_{\phi_v}$ of $\Pi_v$ for which it follows from \cite[Conjecture 6.1(iv)]{Arthur}. The second assertion in (b) is clear since $\pi_v^{\nu_v-\text{gen}}$ belongs to that $L$-packet and has sign $+1$ by definition. The last assertion on (b) is not explicitly written down in \cite{Arthur} but is quite natural and we believe in it (it holds for example for the inner forms of $\U(3)$ by Rogawski's work). The point (c) follows from (b) together 
with \cite[Conjecture 6.1(iii)]{Arthur} 
since the $\pi_v^{\nu_v-\text{gen}}$ belongs to the same 
$L$-packet, independently of $\nu$. 

\begin{lemma} \label{pairingArthur} The map $\bs_\psi \rightarrow \R, \ s \mapsto \langle s,\pi^n \rangle$ is the non-trivial 
character $\epsilon$.
\end{lemma}
\begin{pf} According to the above,
$$\langle s_{\psi}, \pi^n \rangle= \prod_v e_{v}((\pi^n)_v)$$
and $e_v(\pi^n_v)=1$ for every finite place $v$ since 
$\pi^n_v \in \Pi_{\phi_v}$, so we are reduce to show that 
$e((\pi^n)_\infty)=e(\pi^s_\infty)=-1$.

Let $\phi_s$ be the $L$-parameter of $\pi^s_\infty$ and remember 
from the proof of Lemma~\ref{lemmaapacket} that this is {\it not the same} as  
$\phi_\infty=\phi_{\psi_\infty}$. Note that all those $L$-parameters 
of $G(\R)=\U(m)(\R)$, as well as the $A$-parameter $\psi_\infty$, may be seen as
parameters of $G^\ast(\R)$ since those groups have the same $L$-group ${}^L G$. By (b) above we may work with the group $G^\ast(\R)$, and the aim is to show that $e(\phi_s)=-1$.

Arthur, following Adams and Johnson, describes an algorithm to compute 
the $L$-parameter $\phi$'s 
of the representations belonging to $\Pi_{\psi_\infty}$ (we already used it
 for the group $G(\R)$ in the proof of Lemma~\ref{lemmaapacket}) 
and to compute the local pairing. This algorithm, as well as the resulting pairing, 
depend on a choice of a conjugacy class of a Levi subgroup $L^\ast$
of $G^\ast$ whose associated $L$-group is the ${}^L L$ defined in the proof of Lemma~\ref{lemmaapacket}. Here we choose $L^\ast$ to be quasi-split. We denote by $\langle\ , \rangle_{L^\ast}$ the pairing described by Arthur 
using $L^\ast$.

The elements $\pi_w$ in $\Pi_{\psi_\infty}$ are parametrized by the elements 
$w$ of the set 
$$\Sigma^\ast= W(L^\ast,T)\backslash W(G^\ast,T)/W_\R(G^\ast,T)$$
where $T$ is a compact torus of $G^\ast$ contained in $L^\ast$. 
Contrarily to  the case of the corresponding set $\Sigma$ for the compact 
group $G(\R)$ used in the proof of Lemma~\ref{lemmaapacket}, 
this set $\Sigma^\ast$ is not a singleton, corresponding to the fact that the $A$-packet $\Pi_{\psi_\infty}$ for $G^\ast(\R)$ is not a singleton. 
By construction,
the $L$-parameter $\phi_w$ of the representation $\pi_w$ and the values
 $\langle s,\pi_w \rangle_{L^\ast}$ for $s \in {\bs}_\psi$
 depend only on $w$ through the Levi subgroup 
$L_w := w L^\ast w^{-1}$. This Levi subgroup is an inner form of $L^\ast$ defined over $\R$, but is not in general conjugate to $L^\ast$ in $G^*(\R)$.

The parameter $\phi_1$ for $w=1$ is, using the fact that $L^\ast$ is 
quasi-split,
$$\phi_1=\phi_{\psi_\infty}$$
after \cite[last sentence of the first paragraph page 32]{Arthur}.
This ensures from (c) above that the pairing defined by Arthur
using $L^\ast$ is the canonical pairing:
$$\langle s_{\psi} ,\cdot  \rangle_{L^\ast}=\langle s_{\psi} ,\cdot  \rangle_v.$$

Let $w$ be an element of $\Sigma^*$ such that $L:=L_w$ is the compact
inner form of $L^\ast$. Then we have $\phi_w = \phi_s$. The needed computation to check that was actually done during the proof of Lemma~\ref{lemmaapacket}
since the only ingredient used there was that $L$ was a compact 
Levi subgroup.

We thus are reduced to compute 
$$e(\pi_s)=\langle s_\psi,\pi_w\rangle_\infty$$ for $w$ as above. 

For this we have to be a little bit more explicit. We can take for $T$, 
compatibly with the choice already done, the diagonal torus in 
an orthogonal basis $e_1,\dots,e_m$ (in the complex hermitian vector space 
$(V,q)$ 
used to define $G^\ast(\R)$), such that $q(e_{m-1})=q(e_{m})=1$ but $q(e_1)=-1$ (this is always possible since $m \geq 3$ and $G^\ast$ is quasi-split.).
We may define $L$ as the Levi subgroup of matrices stabilizing 
the plane generated by $e_{m-1}$ and $e_m$, and the lines generated by $e_1,\dots,e_{m-2}$: 
it is a compact group. And we may take for $L^\ast$ 
the Levi subgroup of matrices stabilizing the plane generated by $e_1$ 
and $e_m$, and the lines generated by $e_2,\dots,e_{m-1}$: 
it is a quasi-split group. 
Now it is clear that if 
$$w \in W(G^\ast,T)=W(G^\ast(\C),T(\C)) \simeq \got{S}_m$$ is the 
transposition $(1,m-1)$, then $w L^\ast w^{-1}=L$. But that $w$ is the
reflexion $w_\beta$ (cf. \cite[page 33]{Arthur}) 
about the simple root $\beta$ of $G^\ast(\C)$ such that $\beta(\diag(x_1,\dots,x_m))=x_1/x_{m-1}$. Since this root is non compact,
we have by \cite[(5.6) \& (5.7)]{Arthur} :
$e(\pi_s)=\langle s_\psi,\pi_w \rangle = \beta^\vee(s_\psi) =-1$ using 
that $s_\psi=\diag(1,\dots,1,-1,-1)$.
\end{pf}

\begin{theorem} The multiplicity $m(\pi^n)$ of the representation $\pi^n$ in the discrete spectrum of $\U(m)$ is $1$ if $\varepsilon(\pi,1/2)=-1$ 
and $0$ otherwise.
\end{theorem}
\begin{pf} By Prop.~\ref{locequivimplequiv}, $\psi$ is the only $A$-parameter defining the $A$-packet $\Pi$ of $\pi^n$, so we have $m(\pi^n)=m_\psi(\pi^n)$ according to Arthur's definitions. 
By \cite[(8.5)]{Arthur}, 
$$m_\psi(\pi^n)=\frac{1}{|\bs_{\psi}|} \sum_{s \in \bs_\psi} \epsilon_\psi(s) \langle
s,\pi^n \rangle = \frac{1}{2}(1 - \epsilon_\psi(s_\psi))$$ 
using Lemma~\ref{pairingArthur}. 
The theorem then follows from Prop.~\ref{carArthur}.
\end{pf}

\index{A@$A$, a local henselian ring|)}
\index{A@$A$, a commutative ring|)}
\index{E@$E$, a quadratic imaginary field|)}
\index{E@$E$, a number field|)}
\index{GE@$G_E$, the absolute Galois group of $E$|)}
\index{Zrho@$\rho$, a continuous geometric representation of $G_E$|)}

\newpage
\printindex

\end{document}